\documentclass[11pt,a4paper]{myarticle}
\newif\ifarxiv
\arxivtrue
\ifarxiv\else
\usepackage[width=15truecm,height=23truecm,center]{crop}
\def\cropactivated{}
\fi

\let\Spar\S
\errorstopmode
\usepackage{graphicx,color}
\definecolor{brown}{cmyk}{0, 0.87, 0.68, 0.32}

\graphicspath{ {images/} }   
\usepackage{amsmath,amsfonts,amssymb}
\usepackage{mathrsfs}

\usepackage[T1]{fontenc}
\usepackage[applemac]{inputenc}  
\usepackage[french,english]{babel}
\usepackage[all]{xy}
\def\initxy{\catcode`;=12\catcode`:=12\xy }
\SelectTips {cm}{}
\def\xymatrixc#1#2{\vcenter{\hbox{$\displaystyle\xymatrix#1{#2}$}}}
\def\xylbl[#1]#2{\ar @{} [#1] |{\hbox{#2}}}

\input mymacros.sty
\usepackage{makeidx}\makeindex\makeglossary
\input habille.sty
\input tableaux.sty

\def\udr#1 #2 #3 {\rlap{\kern#1pt\vtop{\kern#2pt\hrule width#3pt}}}

\def\mod{\ \mathord{\rm mod\,}}

\def\rmc{{\rm c}}

\def\Hc{H_{\rmc}} 
\def\bm{{\raise0.5pt\hbox{\scriptsize \scshape bm}}}
\def\Bm{{\raise0pt\hbox{\scriptsize \scshape bm}}}
\def\Hbm{H_{\kern-0.5pt\bm}} 
\def\vH{\check\Hr}

\def\cv{{\pi!}}

\def\Hpi{H_{\cv}} 
\def\cHpi{\H_{\cv}}
\def\cHr{\H}
\def\cHHcv{\cl\H_{\cv}}
\def\Hr{H} 
\def\Pc{\P_{\rmc}} 
 
\def\Pcv{\P_{\cv}} 
\def\Hi{H_{!}}
\def\Ind{\mathop{\hbox{\bfseries\slshape Ind}}\nolimits}
\def\Res{\mathop{\rm Res}\nolimits}

\def\Pc{\mathop{\rm \P_c}\nolimits}

\def\Pgoth{{\mathfrak P}}
\def\pgoth{{\mathfrak p}}
\def\qgoth{{\mathfrak q}}
\let\chiorg\chi
\def\chi{{\mathord{\setbox0=\hbox{$\chiorg$}\raise\dp0\hbox{$\chiorg$}}}}
\def\chic{{\mathord{\setbox0=\hbox{$\chiorg$}\raise\dp0\hbox{$\chiorg$}_{\rmc}}}}
\def\chibm{{\mathord{\setbox0=\hbox{$\chiorg$}\raise\dp0\hbox{$\chiorg$}_{\scriptscriptstyle\rm BM}}}}

\def\uacyclique{$\cup$-acyclique}
\def\uacycliques{$\cup$-acycliques}

\def\gr{\mathop{\bf Gr}\nolimits}

\def\Mod{\mathop{\rm Mod}\nolimits}

\let\Sigmaorg\Sigma
\def\largeSigmaorg{\raise-1pt\hbox{\large$\Sigmaorg$}}
\def\Sigmafragile{\mathop{\mathchoice
{{\largeSigmaorg}}
{{\largeSigmaorg}}
{\Sigmaorg}{\Sigmaorg}}}
\def\Sigma{\protect\Sigmafragile}

\def\Gr{\mathop{\bf Gr\mkern1mu}\nolimits}

\def\Modgr{\mathop{\rm Mod\sp{\gr}}\nolimits}

\def\Pr#1{{}\,\widetilde{\!#1}{}}
\def\proj#1/{{}\,\widetilde{\!#1}{}}

\def\Proj#1/{\widetilde{#1}{}}

\def\rest{{\,\vrule height8pt depth3pt\,}}
\def\varrest#1{{\,\smashbot{\vrule height8pt depth#1pt}\vrule depth3pt width0pt\,}}
\def\srest{{\,\vrule height5pt depth2pt\,}}

\newdimen\tooowd \tooowd=1cm
\newdimen\toowd \toowd=0.7cm
\newdimen\towd \towd=5mm

\def\tooht{0.6}
\def\toodp{0.5}

\def\onto{\mathrel{\to\mkern-15mu\to}}

\def\too{\mathop{\goodsmash{\tooht}{\toodp}{\hbox to\toowd{\hss\rightarrowfill\hss}}}\limits}
\def\ontoo{\mathop{\goodsmash{\tooht}{\toodp}{\hbox to\toowd{\hss\rightarrowfill\hss}\hskip-16pt\to}}\limits}
\def\tooo{\mathop{\goodsmash{\tooht}{\toodp}{\hbox to\tooowd{\rightarrowfill}}}\limits}
\def\vartoo#1{\mathop{\goodsmash{\tooht}{\toodp}{\hbox to#1{\rightarrowfill}}}\limits}
\def\longonto{\mathop{{\vrule height5pt depth5pt width0pt\hfonto{}{}{\toowd}}}\limits}

\def\tot{\mathop{\bf tot}\nolimits}
\def\_{\raise-3.5pt\hbox{$\scriptstyle-$}}

\def\aronto{\ar@{->>}}

\def\arhook{\ar@{>->}}
\def\arinto{\ar@{^(->}}
\def\ariinto{\ar@{_(->}}

\def\arhookl{\ar@{_{(}->}}
\def\ardash{\ar@{-->}}\def\ardots{\ar@{..>}}
\def\aregal{\ar@{=}}
\def\arequiv{\ar@{<->}^{\cong}}
\def\arsim{\ar@{<->}^{\sim}}

\def\put#1#2#3{\vtop to0pt{\kern-#2\hbox to0pt{\kern#1#3\hss}\vss}}

\let\Borg\B\def\B{\mathop{\Borg{}}\nolimits}
\let\Zorg\Z\def\Z{\mathop{\Zorg{}}\nolimits}

\def\cone#1{{\hbox{\boldmath$\hat c$}}(#1)} 
\def\coneh#1{{\hbox{\boldmath$\hat c_{h}$}}(#1)}

\def\Rggfragile{\ifmmode R(\ggoth)\else$\Rgg$\fi}
\def\Rgg{\protect\Rggfragile}

\def\Sgbin{\mathbfsl{S}}
\def\Sgfragile{\ifmmode \Sgbin\else$\Sg$\fi}

\def\mdgr{\hbox{\ifboldmath\bfseries\fi\slshape mdg}}
\def\mgr{\hbox{\ifboldmath\bfseries\fi\slshape mg}}

\def\myR{R}
\def\Rmdgrfragile{\ifmmode \myR\tiret\mdgr\else$\Rmdgr$\fi}
\def\Rmdgr{\protect\Rmdgrfragile}
\def\Rmgrfragile{\ifmmode \myR\tiret\mgr\else$\Rmgr$\fi}
\def\Rmgr{\protect\Rmgrfragile}

\def\Diffgrfragile{\ifmmode\mathop{\rm Diff\sp{\gr}}\nolimits\else$\Diffgrfragile$\fi}
\def\Diffgr{\protect\Diffgrfragile}

\def\RDiffgrfragile{\ifmmode \Diffgr(\myR)\else$\RDiffgrfragile$\fi}

\def\RMgrfragile{\ifmmode\Modgr(\myR)\else$\RMgrfragile$\fi}

\def\dual{^\vee}

\def\underlinee#1#2#3{\vtop{\setbox0=\hbox{$#1$}\dp0=#2\copy0\kern1.5pt\hrule\ifnum#3>1\kern1pt\hrule\fi
}}

\def\fs#1{\underlinee{#1}{0pt}{2}}

\def\cxymatrix#1#2{\vcenter{\hbox{$\xymatrix#1{#2}$}}}

\def\gui#1'{{\rm#1'}}

\def\stackover#1#2{\vbox{\mathsurround0pt\offinterlineskip
\setbox1=\hbox{$\scriptstyle#1$}
\setbox2=\hbox{$#2$}
\dimen1=\wd1
\ifdim\dimen1<\wd2 \dimen1=\wd2\fi
\hbox to\dimen1{\hss$\scriptstyle#1$\hss}
\hbox to\dimen1{\hss$#2$\hss}
}}

\def\point{\mathchoice
{\raise2pt\hbox{$\scriptstyle\bullet$}}
{\raise2pt\hbox{$\scriptstyle\bullet$}}
{\raise1pt\hbox{$\scriptscriptstyle\bullet$}}
{\raise1pt\hbox{$\scriptscriptstyle\bullet$}}
}

\newdimen\eccadre
\def\mcadre#1#2{\vbox{\eccadre=#1pt\kern\eccadre\vtop{\hbox{\kern\eccadre $\scriptstyle#2$\kern\eccadre}\kern\eccadre}}}

\def\sep#1 {\vskip1ex
\hbox{\color{red}\footnotesize\tt\count100=#1\relax\loop
\ifnum\count100>0 \%\advance\count100 by-1\repeat
}\vskip1ex}
\newif\ifendpoint
\endpointtrue
\def\nopoint{\global\endpointfalse}
\def\endpoint{\ifendpoint.\else\global\endpointtrue\fi}

\def\IB{{I\mkern-6mu B}}
\def\tilda#1{{\,\nobreak\,\nobreak\widetilde{\!\!{#1}\!}\,}}

\def\supnl{\mathop{\rm sup}\nolimits}

\newenvironment{lcases}{\def\cr{\,,\qquad}}{}

\def\iivdots{\vbox{\baselineskip4pt \lineskiplimit0pt
    \kern6pt\hbox{.}\hbox{.}}}

\mathchardef\Gamma="0100
\mathchardef\Omegai="010A

\def\fsk{\fs k}
\def\Vec{\mathop{\rm Vec}\nolimits}
\def\Gammac{\Gamma_{\rmc}}
\def\AS{\fs{\mathord{\AAA}}}

\def\ttt#1#2{\vrule height#1pt depth#2pt width0pt}

\def\mput#1#2#3{\rlap{\hskip#1\vtop to0pt{\vskip-#2
\hbox{$#3$}\vss}}}

\def\mvcenter#1{\vcenter{\hbox{$#1$}}}
\def\car{\mathop{\rm car}\nolimits}
\def\pt{{\rm pt}}

\def\[#1\]{{[\mkern-3mu[#1]\mkern-3mu]}}
\def\<#1>{\langle#1\rangle}

\def\siiiorg[#1]{\hbox{\scriptsize$\mathsurround0pt\thinmuskip=1mu\medmuskip=1mu {[\mkern-2mu[}{#1}{]\mkern-2mu]}$}}

\def\ssiiiorg[#1]{\hbox{\scriptsize$\scriptstyle\mathsurround0pt\thinmuskip=1mu\medmuskip=1mu {[\mkern-2mu[}{#1}{]\mkern-2mu]}$}}

\def\iiiorg[#1]{\hbox{$\mathsurround0pt\thinmuskip=1mu\medmuskip=1mu {[\mkern-2mu[}{#1}{]\mkern-2mu]}$}}

\def\iii[#1]{\mathchoice
{\iiiorg[#1]}
{\iiiorg[#1]}
{\siiiorg[#1]}
{\ssiiiorg[#1]}
}

\def\iijorg[#1]#2{\hbox{$#2\mathsurround0pt\thinmuskip=1mu\medmuskip=1mu {[\mkern-2mu[}{#1}\mskip1mu{[\mkern-2mu[}$}}

\def\iij[#1]{\mathchoice
{\iijorg[#1]\textstyle}
{\iijorg[#1]\displaystyle}
{\iijorg[#1]\scriptstyle}
{\iijorg[#1]\scriptscriptstyle}
}

\def\bxx#1{{\color{blue}{\{\!}#1{\!\}}}}
\def\bxx#1{{\setbox100=\hbox{$#1$}\smash{\setminus}\mkern1mu
\ifdim\wd100>1em(#1)\else#1\fi}}

\def\fillwith#1 #2 {{\count111=#1\loop#2\hskip0pt plus1pt\penalty0\advance\count111 by-1\ifnum\count111>0\repeat}}

\def\Pr{\P}
\def\Pc{\P_{\rmc}}

\let\Bigorg\Big
\newbox\mybox
\newdimen\mydimen
\def\binomePars#1#2#3#4#5#6{\hbox{#4}$#6\vcenter{\offinterlineskip
\def\un{$\scriptstyle#1$}\def\deux{$\scriptstyle#3$}
\setbox\mybox=\hbox{\un}\mydimen=\wd\mybox\setbox\mybox=\hbox{\deux}
\ifdim\mydimen<\wd\mybox \mydimen=\wd\mybox\fi
\hbox to\mydimen{\hss\un\hss}
\kern#2
\hbox to\mydimen{\hss\deux\hss}
}#6$\hbox{#5}}

\def\binome#1#2{\hbox{\mathsurround0pt
\binomePars
{#1}{8pt}{#2}
{\small\Bigorg(}
{\small\Bigorg)}\,}}

\def\sbParties#1#2{\hbox{\mathsurround-0.1pt
\let\scriptstyle\scriptscriptstyle
\binomePars
{#1}{4pt}{#2}
{\rm\small$\{$}
{\rm\small$\}$}
{}}}

\def\Parties#1#2{\hbox{\mathsurround-0.1pt
\binomePars
{#1}{10pt}{#2}
{\small$\Bigorg\{$}
{\small$\Bigorg\}$}
{}}}

\def\Cycles#1#2{\hbox{\mathsurround0.2pt
\binomePars
{#1}{10pt}{#2}
{\small$\Bigorg[$}
{\small$\Bigorg]$}
{}}}

\def\sbCycles#1#2{\hbox{\mathsurround-0.1pt
\let\scriptstyle\scriptscriptstyle
\binomePars
{#1}{4pt}{#2}
{\rm\small$[$}
{\rm\small$]$}
{}}}

\def\mParties#1#2{\hbox{$|\Pgoth_{#2}(#1)|$}}

\def\Ga{\Gg_a}
\def\smalldelim#1{\mathchoice
{\hbox{\mathsurround0pt\footnotesize$\displaystyle#1$}}
{\hbox{\mathsurround0pt\footnotesize$\textstyle#1$}}
{\hbox{\mathsurround0pt\footnotesize$\scriptstyle#1$}}
{\hbox{\mathsurround0pt\footnotesize$\scriptscriptstyle#1$}}
}
\def\Ga#1{{#1}^{\mkern-2mu{\smalldelim[}a{\smalldelim]}}}
\def\GaDelta_#1{{\def\fin{\end}\doGaDelta#1,\end,{}//}}

\def\doGaDelta#1,#2,#3//{\def\donne{#2}\ifx\donne\fin
\def\nextGa{\Ga\Delta_{#1}}\else
\def\nextGa{\Delta_{#1}^{\mkern-2mu #2{\smalldelim[}a{\smalldelim]}}}\fi
\nextGa}

\def\BF{\Bg\Fg}\def\Fl{\Fg_{\ell}}
\def\BFm{\Bg\Fg_{m}}
\def\Fm{\Fg_{m}}
\def\Fn{\Fg_{n}}
\def\Fmm{\Fg_{\mm}}
\def\Fpp{\Fg_{\pp}}
\def\Fmmo{\Fg_{m{-}1}}
\def\BFmm{\Bg\Fmmo}
\def\Fa{\Fg_a}
\def\Fb{\Fg_{b}}
\def\Fba{\Fg_{b+a} }
\def\pia{\pi_{a} }
\def\comment{\begingroup\color{green}}
\def\endcomment{\endgroup}
\long\def\commentt#1\endcommentt{\relax}
\long\def\comment#1\endcomment{\relax}

\def\Fix{\mathop{\rm Fix}}

\mathchardef\varPhi="7108

\def\mb{{\scriptscriptstyle\bullet}}

\def\oo#1{\mathchoice
{\dooaccent{\displaystyle}{\scriptstyle}{#1}}
{\dooaccent{\textstyle}{\scriptstyle}{#1}}
{\dooaccent{\scriptstyle}{\scriptscriptstyle}{#1}}
{\dooaccent{\scriptscriptstyle}{\scriptscriptstyle}{#1}}
}
\def\dooaccent#1#2#3{{\,\vbox{\offinterlineskip\setbox0=\hbox{$#1\!#3$}\hbox to\wd0{\hss$#2\circ$\hss}\kern1pt\copy0}}}

\def\mmoins{\mathbin{\raise1.5pt\hbox{\footnotesize\mathsurround0pt$\smallsetminus$}}}

\let\SSSorg\SSS
\def\SSSfrag{{\mkern1.5mu\nobreak\SSSorg}}
\def\S{\protect\SSSfrag}
\def\Sm{\S_{m} }
\def\Sa{\S_{a}}
\def\Sb{\S_{b}}
\def\St{\S_{t}}
\def\SL{\S_{L}}
\def\SLp{\S_{L'}}
\def\SA{\S_{A}}
\def\SAp{\S_{A'}}
\def\mm{{m{+}1}}
\def\mt{{m{-}t}}
\def\mmt{{\mm{-}t}}
\def\mmo{{m{-}1}}
\def\pp{p{+}1}
\def\Smm{\S_{\mm}}

\def\Smmo{\S_{\mmo}}
\def\ellmo{{\ell{-}1}}
\let\timesorg\times
\def\stimes{{\timesorg}}
\let\times\stimes
\def\soustitre#1{\vskip\parskip
\noindent{\bfseries\slshape\boldmath
 #1\par\nobreak\vskip3pt \mou\nobreak\parskip0pt\mou\noindent }\ignorespaces}

\def\soustitreline#1{\vskip\parskip\penalty-100\noindent{\bfseries\slshape\boldmath #1 }}

\def\soustitreline#1{\ifhmode\par\fi
\skip0=\lastskip\vskip-\lastskip\penalty-100\advance\skip0 by1ex\mou\vskip\skip0\noindent{\bfseries\slshape\boldmath #1 }}

\let\leq\leqslant
\let\geq\geqslant
\def\circleplustimes{\mathchoice
{\mrlap{{\oplus}}{\otimes}}
{\mrlap{{\oplus}}{\otimes}}
{\mrlap{\scriptstyle{\oplus}}{\otimes}}
{\mrlap{\scriptscriptstyle{\oplus}}{\otimes}}
}
\def\tr{\mathop{\rm tr}\nolimits}
\def\relp{\gooddownstackrel{0pt}\sim\pgoth}
\def\relq{\gooddownstackrel{0pt}\sim\qgoth}
\def\varrel#1{\gooddownstackrel{0pt}\sim#1}

\def\decale#1pt#2{\hskip-#1pt{}{#2}{}\hskip#1pt}
\def\hdecale#1#2{\hskip-#1{}{#2}{}\hskip#1}
\def\notdiv{\mrlap{\hskip2pt\diagup}{{}\div{}}}

\def\div{{\mkern2mu\mid\mkern2mu}}
\def\notdiv{{\mkern2mu\nmid\mkern2mu}}

\def\uline#1{\vtop to0pt{\mathsurround0pt\hbox{$\scriptstyle#1$}
\kern-0\prevdepth\kern1.5pt\hrule height1pt\vss}}
\def\oline#1{\vbox to0pt{\vss\mathsurround0pt\hrule height1pt\kern1.5pt
\hbox{$\scriptstyle#1$}}}

\def\Oline#1#2#3{\kern#1\vbox{\mathsurround0pt\hrule height1pt\kern1.5pt
\hbox{$\kern-#1\relax#3$\kern#2}}}
\def\Uline#1#2#3{\vtop{\mathsurround0pt\hbox{\kern#1$#3$\kern-#2}\kern1.5pt\hrule height1pt}\kern#2}

\def\dofd#1#2#3{{\muskip0=#1
\mkern-\muskip0\uline{\mkern\muskip0{#3}}#2}}

\def\usp#1{^{\dofd{1mu}{}{#1}}}

\def\uspp#1[#2]{^{\setbox0=\hbox{$#2$}
\ifdim\wd0>2ex\dofd{1mu}{\raise-0.25pt\hbox{$*$}(#2)}{#1}\else
\dofd{1mu}{\raise-0.25pt\hbox{$*$}#2}{#1}\fi
}}
\def\FD#1{^{\dofd{0mu}{}{#1}}}

\def\dofa#1#2{^{\muskip0=#1
\oline{{#2}\mkern\muskip0}\mkern-\muskip0}}
\def\osp{\dofa{1mu}}
\def\FA{\dofa{0mu}}

\def\ostirling#1#2{\Oline{1pt}{0.3pt} s({#1},{#2})\mkern1mu}

\def\ustirling#1#2{\Uline{0.5pt}{1pt} s({#1},{#2})\mkern.75mu}
\def\oStirling #1#2{\Oline{2pt}{0.5pt} S({#1},{#2})\mkern1mu}
\def\uStirling#1#2{\Uline{0.5pt}{1pt} S({#1},{#2})\mkern.75mu}

\def\diag{\mathop{\rm diag}\nolimits}

\let\trans\pitchfork
\def\actson{\mathbin{\raise0.5pt\hbox{\rotatebox[origin = c]{-90}{$\circlearrowright$}}}}

\def\sumnl{\sum\nolimits}
\def\prodnl{\prod\nolimits}
\def\ind{\mathop{\rm ind}\nolimits}
\def\card{\mathop{\rm Card}\nolimits}
\def\coprodnl{\coprod\nolimits}

\def\dobouquet#1#2{\vtop{\mathsurround-0.5pt\hbox{$#1\vee$}\hrule height#2}}
\def\bouquet{\mathbin{\mathchoice
{\dobouquet\textstyle{0.4pt}}
{\dobouquet\textstyle{0.4pt}}
{\dobouquet\scriptstyle{0.35pt}}
{\dobouquet\scriptscriptstyle{0.3pt}}
}}
\def\amalgam{+}

\def\vscalesum#1#2{\mathop{\setbox0=\hbox{$\sum$}\dimen0=#1\dp0\vrule depth2.5\dimen0 width0pt
\vtop to0pt{\hbox{\scalebox{1.2}[#1]{\raise#2\hbox{$\displaystyle\sum$}}}}}}

\let\boxtimesorg\boxtimes
\def\boxtimes{\mathbin{\mathchoice
{\hbox{\footnotesize$\boxtimesorg$}}
{\hbox{\footnotesize$\boxtimesorg$}}
{\hbox{\footnotesize$\scriptstyle\boxtimesorg$}}
{\hbox{\footnotesize$\scriptscriptstyle\boxtimesorg$}}}
}

\def\sldash#1{\noindent{\sl$\circ$ #1\unskip.\,---\kern1ex }\ignorespaces}
\def\bulldash#1{\noindent{\sl$\bullet$ #1\unskip.\,---\kern1ex }\ignorespaces}
\def\bfsldash#1{\noindent{\slshape$\bullet$ #1 }}

\long\def\proclaim#1{
{\sl\leavevmode
\relax{\boldmath$\vdash$ }\ignorespaces #1}}
\long\def\subproclaim#1{{\sl\leavevmode\llap{$\Vdash$ }\ignorespaces #1}}

\def\Bettic{\mathop{\rm Betti}\nolimits_{\rm c}}
\def\Bettibm{\mathop{\rm Betti}\nolimits_{\bm}}
\def\Betti{\mathop{\rm Betti}\nolimits}
\def\FI{\mathop{\Fg\Ig}\nolimits}
\def\FB{\mathop{\Fg\Bg}\nolimits}
\def\ModFI#1{\mathop{\Mod(#1[{\bf FI}])}\nolimits}
\def\ModFIsc#1#2{\mathop{\Mod_{#1}(#2[{\bf FI}])}\nolimits}
\def\ModFB#1{\mathop{\Mod(#1[{\bf FB}])}\nolimits}
\def\ev{\mathop{\rm ev}\nolimits}
\def\miseengarde{\bfseries\boldmath
\sffamily}
\def\CCc{k_{\rm c}}

\def\Gl{\mathop{\rm Gl}\nolimits}
\def\sgn{\mathop{\rm sgn}\nolimits}

\def\rk{\mathop{\rm rg}\nolimits}
\def\rks{\mathop{\rm rg_{\rm s}}\nolimits}
\def\rkm{\mathop{\rm rg_{\rm m}}\nolimits}
\def\rkms{\mathop{\rm rg_{\rm ms}}\nolimits}
\def\rkmse{\mathop{\rm rg_{\rm ms}^{e}}\nolimits}
\def\rkme{\mathop{\rm rg_{\rm m}^{e}}\nolimits}
\def\rke{\mathop{\rm rg^{e}}\nolimits}

\def\degstab{\mathop{\rm deg\tiret stab}\nolimits}
\def\deginj{\mathop{\rm deg\tiret inj}\nolimits}

\let\poids\PPP

\def\poids{\mathop{{}\PPP}\nolimits}
\def\poidsinf{\mathop{{}\PPP}\nolimits_{\infty}}

\def\bigbsl{\raise-3pt\hbox{\Large$\backslash$}}
\def\bigsl{\raise-3pt\hbox{\Large$/$}}

\def\biquo#1#2#3{{\scriptspace0pt{\;}\raise-7pt\hbox{$\!#1\!$}
\bigbsl\raise2pt\hbox{$#2$}\bigsl
\kern-2pt\raise-7pt\hbox{$#3$}{\,}}}

\def\lquo#1#2{\hbox{\small\big[{\,}\raise-4pt\hbox{$\!#1\kern-6pt$}
\raise-2.5pt\hbox{\Large$\backslash$}
\raise0pt\hbox{\normalsize\kern-4.5pt$#2$}\big]}}

\def\bigbslsb{\raise-1pt\hbox{$\scriptstyle\backslash$}}
\def\bigslsb{\raise-1pt\hbox{$\scriptstyle/$}}

\def\biquosb#1#2#3{{\scriptspace0pt
{\;}\raise-3pt\hbox{$\scriptstyle\!#1\!$}
\bigbslsb
\raise2pt\hbox{$\scriptstyle#2$}
\bigslsb
\raise-3pt\hbox{$\scriptstyle#3$}{\,}}}
\def\txtquo#1#2#3{#1\backslash #2/#3}

\def\ulambda{{\underline\lambda}}
\def\utau{{\underline\tau}}
\def\uell{{\underline\ell}}

\def\Thetabf{\hbox{\boldmath$\Theta$}}

\def\Ibf {\mathop{{}\Ig}\nolimits}
\def\Slambda{\SSS_{\lambda}}
\def\Plambda{\PPP_{\lambda}}
\def\Glambda{\myG_{\lambda}}

\def\Pulambda{\PPP_{\ulambda}}

\def\myG{G}

\def\Ya{\agoth}
\def\Yb{\bgoth}
\def\Yc{\cgoth}
\def\Yl{\lgoth}

\def\rep{{\mkern1mu:\mkern1mu}}
\def\Rep{{\mkern2mu:\mkern2mu}}
\let\actson\rep

\def\killline{\vskip-1ex\nobreak}

\def\XXfragile#1{\mathchoice
{\hbox{\footnotesize\sffamily\slshape #1}}
{\hbox{\footnotesize\sffamily\slshape #1}}
{\hbox{\tiny\sffamily\slshape #1}}
{\hbox{\tiny\sffamily\slshape #1}}
}
\def\XX{\protect\XXfragile X}
\def\YY{\protect\XXfragile Y}
\def\Qsf{\protect\XXfragile Q}

\def\Zsf{\protect\XXfragile Z}

\def\Times{\mathop{\raise-1pt\hbox{\Large$\times$}}}
\def\vCech{\v Cech}
\def\vC{\check C}
\def\vCm{\check C_{<}}
\def\StS#1#2{\S_{#1}\times\S_{#2}}
\def\Omegac{\Omega_{\rm c}}

\def\fsOmega{\fs\Omega}
\def\sssearrow{\,{\scriptscriptstyle\searrow}\,}

\def\alt{\hbox{\boldmath$\sigma$}}

\def\IEs{\IE_{\sigma}}
\def\xysimeq{\vrule height7pt depth4pt width0pt\simeq}
\def\Mgg{\Mg_{>0}}\def\Mgge{\Mg_{\geq0}}

\let\cL\L
\let\cSS\SS
\def\tthref#1#2{{\tt#1}}
\def\varhref{\ifx\href\undefined
\let\nexthref\tthref\else
\let\nexthref\href\fi\nexthref}
\ifarxiv
\else
\usepackage{hyperref}
\hypersetup{
    colorlinks=true,
    linkcolor=blue,
    filecolor=magenta,      
    urlcolor=cyan,
}
\fi
\def\arcs/{\ggggg\smash{\rlap{\color{red}\vrule height-1pt depth4pt width1cm}}arcs}

\def\dimch{\mathop{\rm dim}\nolimits_{\rm ch}}

\def\binspace#1{\medmuskip#1mu}
\def\relspace#1{\thickmuskip#1mu}
\def\nobinspace{\medmuskip0mu}

\def\mycaption#1{\hbox to\hsize{\hss#1\hss}}

\def\ee#1 {\hbox to5mm{\rm\hss#1)\ }}
\def\Otimes_#1{\mathrel{\displaystyle\mathop{\textstyle\bigotimes}\limits_{#1}}}
\def\varOtimes#1#2_#3{\mathrel{\displaystyle
\hskip-#1\mathop{\textstyle\bigotimes}\limits_{\vrule height#2 width0pt#3}\hskip-#1}}

\def\bigoplusnl{\bigoplus\nolimits}

\def\bfit{\bfseries\slshape}

\def\ddeg{\mathop{\underline{\smash{\rm deg}}}\nolimits}

\def\leftbrace#1#2{\smash{\raise#1\hbox{$\left\{\vrule height#2 width0pt \right.$}}}

\def\dMg{{d_{\Mg}}}
\def\dXg{{d_{\Xg}}}

\def\mobius#1#2{\mu\big(\myfrac{#1}{#2}\big)}

\def\fgg{\mathbfsl  f}

\def\minicroch[#1]{{\scriptscriptstyle[}#1{\scriptscriptstyle]}}
\def\Qdeg#1{\Q^{\minicroch[#1]}}
\def\chideg#1{\chi^{\minicroch[#1]}}
\def\minideg#1{^{\minicroch[#1]}}

\def\UV{\V}

\let\ssubsubsection\subsubsubsection
\let\ssubsubsectionline\subsubsubsectionline

\begin{document}
\overfullrule1cm
\hfuzz7pt

\def\expression#1{{\ifnum0=\fontdimen1\font\sl\else\rm\fi #1\/}\index{#1}}
\def\expressiong#1{{{\rm\og}\ifnum0=\fontdimen1\font\sl\else\rm\fi #1\/{\rm\fg}\index{\relax#1}}}
\let\expression\expressiong

\def\xpress#1{{{\rm`}\ifnum0=\fontdimen1\font\sl\else\rm\fi #1{\rm'}}}

\def\UU{U}
\def\WW{W}

\newcount\shwky
\newbox\skbox
\def\showkeys#1{\relax}
\def\showkeys#1{\global\shwky=#1\relax
\ifnum#1=0\else\def\shwkydc{\ifodd\shwky\kern8.5pt\else\kern-7pt\fi}
\def\shwkyodd{\ifodd\shwky$\downarrow$\else$\uparrow$\fi}
\let\labelorg\label
\def\label##1{\global\advance\shwky1
\setbox\skbox=\hbox to0pt{\hss\vbox to0pt{\vss\tiny\hbox{\color{red}\tt\shwkyodd<##1>\shwkyodd}\shwkydc}\hss}\ifhmode\copy\skbox\else\fi\labelorg{##1}}
\let\reforg\ref
\def\ref##1{\global\advance\shwky1\leavevmode\hbox to0pt{\hss\vbox to0pt{\vss\tiny\hbox{\color{blue}\tt\shwkyodd\{##1\}\shwkyodd}\shwkydc}\hss}\reforg{##1}}
\let\citeorg\cite
\def\cite##1{\global\advance\shwky1\leavevmode\hbox to0pt{\hss\vbox to0pt{\vss\tiny\hbox{\color{cyan}\tt\shwkyodd[##1]\shwkyodd}\shwkydc}\hss}\citeorg{##1}}
\fi}

\def\sk{\ifx\labelorg\undefinded\showkeys1\else\global\advance\shwky1\fi}
\def\skk{\ifx\labelorg\undefinded\showkeys2\else\global\advance\shwky1\fi}


\mysettings
\ifx\cropactivated\undefined\else\advance\voffset-0.5cm\fi
\parskip2pt\mou
\selectlanguage{french}
\mylistskips

\begin{center}\renewcommand{\thefootnote}{\hbox{\scriptspace=0pt$(^*)$}}
\null\vskip-2em\vskip-3em\Large\bold Espaces de configuration généralisés\\Espaces topologiques $i$-acycliques\\Suites spectrales basiques

\medskip\large\bf Alberto Arabia\rlap{ \thefootnote}
\catcode`\@=11
\def\@makefnmark{\thefootnote\,}\footnotetext{Université Paris Diderot-Paris 7, IMJ-PRG, CNRS, Bâtiment Sophie Germain, bureau 608, Case 7012, 75205. Paris Cedex 13, France. Contact: {\tt alberto.arabia@imj-prg.fr}.}

\ifarxiv\else
{\normalsize décembre 2016}\fi
\end{center}

{\renewcommand{\thefootnote}{\hbox{}}\footnotetext{\hskip-2em\sl MSC-class: 55-XX (Algebraic topology), 55R80 (Configuration spaces), 
20-XX (Group theory), 20C30 (Representations of finite symmetric groups), 18G40 (Spectral sequences)}}

\begingroup
\selectlanguage{english}
\ifarxiv\vskip1.5em
\parskip2pt\mou
\smallskipamount=0.5\smallskipamount\mou
\medskipamount=0.5\medskipamount\mou
\else\vskip1.5em
\parskip2pt\mou
\smallskipamount=0.8\smallskipamount\mou
\medskipamount=0.8\medskipamount\mou
\fi
\small\mydisplayskips

\noindent {\bf Abstract. }The \emph{generalized\/} configuration spaces of a topological space $\Xg$ are the subspaces $\Delta_{?\ell}\Xg^{m}\dans\Xg^{m}$, defined, for all $0\leq\ell\leq m\in\NN$, as:
$$
\begin{casesalign}\noalign{\kern-2pt}
\Delta_{\leq\ell}\Xg^{m}&:=&
\bigset(x_1,\ldots,x_{m})\in\Xg^{m}\mid\card{\set x_1,\ldots,x_{m}/}\leq \ell/\,,\\\noalign{\kern2pt}
\Delta_{\ell}\Xg^{m}&:=&\Delta_{\leq\ell}\Xg^{m}\mmoins\Delta_{\leq\ell-1}\Xg^{m}\,,\quad\Fm(\Xg):=\Delta_{m}\Xg^{m}\,.\hfill\\\noalign{\kern-1pt}
\end{casesalign}\eqno(*)
$$ 
They are endowed with the action of the symmetric group $\S_m$  permuting coordinates. 
Our first motivation in this work was to transpose standard problems on classical configuration spaces $ \Fm (\Xg) $ to the generalized ones $(*)$ and to try to solve them for large families of spaces  using  uniform methods.
Among the questions we considered, the following had quite complete answers.
\begin{list}{--}{\mylistskips\parskip0pt
\leftmargin\parindent\itemsep2pt
\itemindent0em\topsep4pt
\parsep1pt
\labelsep1ex
}
\item Compute the character of the representation $\S_m\rep\Hr(\Delta_{? \ell}\Xg^{m})$. 

\item Compute the Poincaré polynomial of the quotients of $\Delta_{? m-a}\Xg^{m}$ by finite subgroups of $\S_m$. Show that, for fixed $a\in\NN$ the Betti numbers of $\Delta_{? m-a}\Xg^{m}$ are given by  \xpress{universal\/} polynomials on $m$ and the Betti numbers of $\Xg$.
\item Prove the degeneracy of the Leray spectral sequences for usual maps between generalized configuration spaces, e.g. the canonical projections $\Fg_{b+a}(\Xg)\onto\Fg_{a}(\Xg)$.

\item Given $a\in\NN$, estimate the ranks for representation stability and character polynomiality, in the Church-Farb sense (\cite{chu-far}), of the family $\set\S_m\rep\Hr(\Delta_{?m-a}\Xg^{m})/_{m}$. 

\end{list}

In our approach to these problems, we used what we call \xpress{the fundamental complex of $\Xg$ for $\Delta_{\leq\ell}\Xg^{m}\,$}. This is a particular complex of graded $\S_m$-modules
$$\def\term#1{\Delta_{#1}\Xg^{m}}
0\to\Hc^{*-\ell+1}(\term1)\to\cdots\to
\Hc^{*-1}(\term{\ell-1})\to
\Hc^{*}(\term\ell)\to
\Hc^{*}(\Delta_{\leq\ell}\Xg^{m})\to0\,,
$$
having the remarkable property of being exact whenever the \xpress{interior} cohomology of $\Xg$, \idest  the
image of the natural map $\Hc(\Xg)\to\Hr(\Xg)$, vanishes. Moreover, these  are equivalent properties when $\Xg$ is an oriented topological manifold (th. \ref{theo-complexe-exact}).

We call \xpress{$i$-acyclic} any space with no interior cohomology. First examples of such spaces are non-compact acyclic spaces, non-compact connected real Lie groups, and, more generally, if $\Xg$ is $i$-acyclic, by every open subset $U\dans\Xg$, every quotient $\Xg/W$ by a finite group $W$, and every cartesian product $\Xg\times \Yg$ for arbitrary $\Yg$.

An important part of this work is about what follows the fact that the fundamental complex is an $\S _{m}$-module resolution for $\Hc(\Delta_{\leq\ell}\Xg^{m})$, and thereafter, that questions about Betti numbers,
character formulas and representation stability,
may be tracked through recursive methods. 

Another important part, is devoted to the design of a spectral sequence that creates a bridge from configuration spaces of $i$-acyclic spaces to configuration spaces of general spaces. The device allows the generalization of some of our main results for configuration spaces of $i$-acyclic spaces to the general case.

\medskip
\noindent -- For the character formula problem, 
we generalize Macdonald's well-known 
formula of the character of the cohomology of the cartesian products $\Xg^{m}$ (\cite{mac}), to the case of the configuration spaces $\Fm(\Xg)$.

\smallskip
\noindent{\slshape Theorem (\ref{theo-trace-gen}). Let $\Xg$ be an $i$-acyclic space. If $\alpha\in\S _{m}$, we have
$$\relax{{\chic(\Fm (\Xg))(\alpha,T)\over T^{m}}=\prodnl_{d=1}^{m}d^{\XX_d}
\Big(\sumnl_{e\div d}\mobius de 
{\Pc(\Xg)(-T^{e})\over dT^{e}}\Big)\usp {\XX_d}}\,,
$$
where $\chic(\Fm )(\alpha,T):=\sum_{i\in\ZZ}\tr(\alpha\actson\Hc^{i}(\Fm ))\,(-T)^{i}$, 
and $(1^{\XX_1},2^{\XX_{2}},\ldots,m^{\XX_{m}})\vdash m$ is the type of the permutation $\alpha$, 
$\mu(\_)$ is the Möbius function, $\Pc(\_)$ is the Poincaré polynomial of $\Hc(\_)$, and $(\_)\usp r$ denotes the falling factorial.}

\medskip
\noindent -- For the Poincaré polynomial problem, we settled the case of $\Fm(\Xg)$ by a simple closed formula, almost immediate consequence of the $i$-acyclicity property, while for the quotients of $\Fm(\Xg)$ by finite subgroups of $\S_m$, we use the previous character formula (\ref{theo-trace-gen}). As examples, we worked through the cases of the \xpress{cyclic} configuration space $\Cg\Fm (\Xg):=\Fm (\Xg)/\Cg_{m}$, where $\Cg_{m}:=\langle(1,\ldots,m)\rangle\dans\S_m$, and of the \xpress{unordered}
configuration space $\Bg\Fm (\Xg):=\Fm (\Xg)/\S _{m}$. 

\smallskip
Let $\Pc(\_)$ denote the Poincaré polynomial of $\Hc(\_)$, 
$\phi(\_)$ the Euler $\phi$ function,
$\mu(\_)$ the Möbius function, and
$(\_)\usp r$ the falling factorial. 
The following equalities hold whenever $\Xg$ is an $i$-acyclic space.

\begingroup\parindent0pt
\slshape
\vskip1em

\noindent\rlap{Theorem (\ref{pol-poincare}):}\hfill
$\displaystyle
{\Pc(\Fm (\Xg))(-T)\over T^{m} }= \Big({\Pc(\Xg)(-T)\over T}\Big)\usp m\,.$\hfill\null

\medskip 
Theorem (\ref{theo-conf-cycliques}):
$$\relax{{\Pc(\Cg\Fm )(-T)
\over T^{m}}=\relax{1\over m}\sum_{d\div m}\phi(d)\,
d^{m/d}
\Big(\sum_{e\div d}\mu\big(\myfrac de \big)
{\Pc(\Xg)(-T^{e})\over dT^{e}}\Big)\usp {m/d}}
\,.\postskip0pt$$

Theorem (\ref{theo-conf-symetriques})
$$\mathrigid1mu
\hss\relax{{\Pc(\Bg\Fm (\Xg))(-T)\over T^{m}}=\decale-5pt{1\over m!}
\mkern-25mu
\vscalesum{1.75}{-1.3pt}_{\vrule height0pt width0pt\lambda:=(1^{\XX_1},\ldots,m^{\XX_{m}})\vdash m}
\mkern-30mu 
\decale3pt{h_{\lambda}}
\prodnl_{d=1}^{m}d^{\XX_d}
\bigg(\sum_{e\div d}\mu\big(\myfrac de \big)
{\Pc(\Xg)(-T^{e})\over dT^{e}}\bigg)\usp {\XX_{d}}}
\,,\hss$$
where  $h_{\lambda}$ is the cardinal of the set of permutations of $\S _{m}$ whose cycle decomposition
is of type $\lambda:=(1^{\XX_1},\ldots,m^{\XX_{m}})\vdash m$.

\endgroup

\medskip 
\noindent -- The particularly simple form of the Poincaré polynomial of $\Fg_{m}(\Xg)$ suggested a sort of cohomological triviality for the projections $\pi_a:\Fba (\Xg)\to\Fa (\Xg)$ and  consequently,  the degeneration of the associated Leray spectral sequences. We show that this is indeed the case when $\Xg$ is $i$-acyclic and locally connected (Th.~\ref{degen}).
\medskip

\noindent -- For the representation stability problem, we prove the  following theorems.

\nobreak\smallskip\noindent
{\slshape
 Theorem (\ref{theo-stabilite-BM-pseudo}).
Let $\Mg$ be a connected oriented pseudomanifold of dimension $\geq 2$. For $a,\,i\in\NN$, the family of representations 
$\set\S_m\rep\Hbm ^{i}(\Delta_{?m-a}\Mg^{m})/_{m}$
is monotone and stationary for $m\geq4i+4a$, if $d_{\Mg}=2$, and for  $m\geq2i+4a$, if~$d_{\Mg}\geq3$. The corresponding families 
of characters and  Betti numbers are (hence) polynomial and the family $\set \Bettibm^{i}({\Delta_{?m-a}\Mg^{m}/\S_m})/_{m}$ is constant
 within the same ranges of $m$. \parskip4pt

\noindent Proposition \ref{prop-stabilite-betti-dim3-BFm} states moreover that the family  $\set \Bettibm^{i}({\Bg\Fm(\Mg)})/_{m}$ is constant for $m\geq 2i$, if $d_{\Mg}=2$, and for $m\geq i$, if $d_{\Mg}\geq 3$.

}

\smallskip
These theorems were proved by Church (\cite{chu}, 2012) for $\Mg$ smooth and for the family $\set{\Sm\rep\Hbm(\Fm(\Mg))}/_{m}$. 
We succeeded in generalizing Church's theorems following two directions. First by removing the regularity assumption in the space $\Mg$, and second by incorporating
the families of generalized configuration spaces $\set\Delta_{\leq m-a}\Mg^{m}/_m$ (singular even if $\Mg$ is smooth), for which there was no previous conjecture. 
Our methods are completely different from those of Church, who relies on the work of Totaro (\cite{tota}) on the Leray spectral sequence associated to the embedding $\Fm(\Mg)\hook\Mg^{m}$ when $\Mg$ is smooth.

Our strategy was to prove first the theorems when $\Mg$ is $i$-acyclic (\ref{theo-stabilite-BM-i-acyclique}) using a combinatorial argument based on the exactness of the fundamental complexes that 
allows the computation of the stability and monotonicity ranks from those of the spaces $\Delta_{\leq m}\Mg^{m}=\Mg^{m}$ for which the answer is quite simple. The combinatorics make use of two induction functors in the category of $\FI$-modules, the functors
$\Ibf^{a},\Thetabf^{a}:\ModFI k\to\ModFI k$ (\ref{foncteur-I-a}) that shift  $\S_{m-a}$-modules to $\Sm$-modules and for which we can control the way they modify the stability and monotonicity ranks (thm. \ref{theo-Ind-FI}). 

The statement for a general pseudomanifold $\Mg$ is afterwards handled through the fact that $\Mg$ is the difference $\Mg=\Mgge{\mmoins\,}\Mgg$ where both $\Mgge:=\Mg\times\RR_{\geq0}$ and $\Mgg:=\Mg\times\RR_{>0}$ are $i$-acyclic. One is then naturally lead to construct a spectral sequence $(\IE_{\sigma}(\U^{m})_{r},d_r)$ (\ref{theo-suite-spectrale-basique}) converging to $\Hbm(\Fm(\Mg))$ and such that the $\IE_1$ page only concerns configuration spaces for $i$-acyclic spaces. We have (\ref{prop-tableau-normaux-ss-basiques}-(\ref{prop-tableau-normaux-ss-basiques-b})): 
$$\def\tt{{\!\vrule height0pt depth4pt width0pt}}
\hss\relax{\IEs(\U^{m})_{1}^{p,q}=
\bigoplusnl_{\tau\in\TTT(p+1,m)}
\ind\tt^{\Sm}_{\HHH_{\tau}}\ 
\alt\otimes\Hbm^{Q}(\Fg_{p+1}(\Mgg))
\ \Rightarrow\ 
\Hbm^{i}(\Fm(\Mg))}\hss\,,$$
with $\medmuskip0mu
Q:=i-(m-(\pp))\,(d_{\Mg}-1)$. We denote by
 $\TTT(p{+}1,m)$ the set of Young tableaux of $m$ boxes with first column $(m{-}p,\ldots,m)$, and by $\HHH_\tau$ the stabilizer of $\tau$ in $\S_{m{-}(p{+}1)}\times\S_{\pp}$ acting on $\Hbm(\Fg_{\pp}(\Mg_{>0}))$ through the alternate representation $\alt$ of $\S_{m{-}(p{+}1)}$.
 We call this spectral sequence \xpress{basic\/}. The construction makes it compatible with pull-backs and allows the estimation of
the stability and monotonicity ranks of the family $\set\Hbm^{i}(\Fm(\Mg))/_{m}$
from those already known  in the $i$-acyclic case of $\set\Hbm^{i}(\Fm(\Mgg))/_{m}$. The induction functor $\Ibf^{a}$ is then used again to settle the case of  $\set \Hbm^{i}(\Delta_{m-a}\Mg^{m})/$, and a decreasing induction based on long exact sequences handles the remaining case of $\set\Hbm^{i}(\Delta_{\leq m-a}\Mg^{m})/$.

The study of the families  $\set \Bettibm^{i}({\Bg\Fm(\Mg)})/_{m}$ follows the same approach. When $\Mg$ is $i$-acyclic,
the explicit formula for the character $\chic(\Fm (\Xg))$ in theorem \ref{theo-trace-gen} gives us quite precise information on the multiplicity of the trivial representation of $\Sm$ in $\Hbm(\Fm(\Mg))$. We are then able to prove the following fact.

\smallskip\noindent
{\slshape
Proposition \ref{theo-stabilite-i-acyclique-betti-BFm}.
Let $\Mg$ be a connected oriented $i$-acyclic pseudomanifold. The family $\set\Bettibm^{i}(\BFm(\Mg;\QQ))/_{m}$ is constant for each $i\in\NN$ and all $m\geq i$.
}

\smallskip
Based on this result, the use of basic spectral sequences settles
the case already mentioned (prop. \ref{prop-stabilite-betti-dim3-BFm}) where $\Mg$ is a general pseudomanifold.
\smallskip
\penalty-500
We believe that even if there are many questions still to be settled, the usefulness of fundamental complexes as a combinatorial tool in the case of $i$-acyclic spaces, and of the basic spectral sequences to pervade the broader category of pseudomanifolds, should be clear following this work.


\endgroup
\begingroup
\selectlanguage{french}
\ifarxiv
\nobreak\vskip0.1em
\hbox to\hsize{\hss \vrule height4pt depth-3.6pt width2cm \hss}
\parskip3.pt\mou
\smallskipamount=0.9\smallskipamount\mou
\medskipamount=0.9\medskipamount\mou
\small\mydisplayskips
\noindent
{\bf Résumé. }Les espaces de configuration \expression{généralisés} d'un espace topologique $\Xg$ sont les sous-espaces $\Delta_{?\ell}\Xg^{m}\dans\Xg^{m}$ définis, pour $0\leq\ell\leq m\in\NN$, par 
$$
\begin{casesalign}\noalign{\kern-2pt}
\Delta_{\leq\ell}\Xg^{m}&:=&
\bigset(z_1,\ldots,z_{m})\in\Xg^{m}\mid\card{\set z_1,\ldots,z_{m}/}\leq \ell/\,,\\\noalign{\kern2pt}
\Delta_{\ell}\Xg^{m}&:=&\Delta_{\leq\ell}\Xg^{m}\mmoins\Delta_{\leq\ellmo}\Xg^{m}\,,\quad\Fm(\Xg):=\Delta_{m}\Xg^{m}\,.\hfill\\\noalign{\kern-1pt}
\end{casesalign}
$$ 
Ils sont munis de l'action du groupe symétrique $\Sm$ 
par permutation de coordonnées. Notre première motivation dans ce travail a été de transposer certaines questions standard sur les espaces de configuration classiques $\Fm(\Xg)$ aux espaces de configuration généralisés, pour ensuite tenter de les résoudre pour des larges familles d'espaces à l'aide de méthodes uniformes. Parmi ces questions, les suivantes ont trouvé des réponses très complètes pour les coefficients de cohomologie dans un corps de caractéristique nulle.
\begin{list}{--}{\mylistskips
\parskip2pt
\leftmargin\parindent\itemsep0pt
\itemindent0em\topsep4pt
\parsep2pt
\labelsep1ex
}
\item Calculer le caractère de la représentations $\Sm \rep\Hr(\Delta_{? \ell}\Xg^{m})$. 

\item Calculer le polynôme de Poincaré des quotients de $\Delta_{? m-a}\Xg^{m}$ pas des sous-groupes de $\Sm $. Montrer que, pour $a\in\NN$ fixé, les nombres de Betti de $\Delta_{? m-a}\Xg^{m}$ sont donnés un polynôme \expression{universel} en $m$ et en les nombres de Betti de $\Xg$.

\item Montrer que la suite spectrale de Leray des fibrations entre espaces de configuration dégénèrent, notamment celles associées aux projections $\Fg_{b+a}(\Xg)\onto\Fg_{a}(\Xg)$.

\item Pour $a\in\NN$ fixé, estimer les rangs de stabilité de représentation et de polynomialité de caractère, au sens de Church-Farb (\cite{chu-far}), de la famille  $\set\Hr(\Delta_{?m-a}\Xg^{m})/_{m}$. 
\end{list}

Ces questions sont abordées au moyen de ce que nous appelons le \expression{complexe fondamental de $\Xg$ pour $\Delta_{\leq\ell}\Xg^{m}$}. Il s'agit d'un complexe de $\Sm$-modules gradués:
$$\def\term#1{\Delta_{#1}\Xg^{m}}
0\to\Hc^{*-\ell+1}(\term1)\to\cdots\to
\Hc^{*-1}(\term{\ellmo})\to
\Hc^{*}(\term\ell)\to
\Hc^{*}(\Delta_{\leq\ell}\Xg^{m})\to0\,,
$$
{\spaceskip0.8ex minus0.5ex
\looseness-1
\tolerance8000
qui a la propriété remarquable d'être exact lorsque la cohomologie \expression{intérieure} de $\Xg$, \idest  de
l'image de l'application naturelle $\Hc(\Xg)\to\Hr(\Xg)$, est nulle. Il y a même équivalence de ces propriétés lorsque $\Xg$ est une variété topologique orientable (th. \ref{theo-complexe-exact}).

}
\looseness-1
Nous appelons \expression{$i$-acyclique} tout espace sans cohomologie intérieure. Les premiers exemples de tels espaces sont les espaces acycliques connexes non compacts, les  groupes de Lie réels connexes non compacts, et, plus généralement, si $\Xg$ est $i$-acyclique,
tout ouvert $U\dans\Xg$, tout quotient $\Xg/W$ par un groupe fini $W$, et tout produit cartésien  $\Xg\times\Yg$, où $\Yg$ est quelconque, le sont aussi. 
\comment La somme amalgamée de variétés topologiques orientables $i$-acycliques est aussi $i$-acyclique.\endcomment


\smallskip
Une large partie de ce travail s'intéresse aux conséquences du fait que le complexe fondamental est une résolution par $\S _{m}$-modules de $\Hc(\Delta_{\leq\ell}\Xg^{m})$. Les questions de stabilité de représentations, de formules de caractères, de nombres de Betti, sont alors susceptibles d'approches récursives. Une autre partie importante est consacrée à la conception d'une suite spectrale particulière qui établit un pont entre les espaces de configuration des espaces $i$-acycliques et ceux des espaces généraux, ce qui permet d'entreprendre la généralisation des résultats connus du cadre $i$-acyclique au cadre général.

\smallskip
\noindent -- Pour le problème de la formule des caractères, nous étendons la formule bien connue de Macdonald  du caractère de la cohomologie des produits cartésiens $\Xg^{m}$ (\cite{mac}), au cas des espaces de configuration $\Fm (\Xg)$.

\smallskip
\noindent{\slshape Théorème. Soit $\Xg$ un espace $i$-acyclique. Si $\alpha\in\S _{m}$, on a
$$\displayindent-0\hsize
\relax{{\chic(\Fm (\Xg))(\alpha,T)\over T^{m}}=\prod\nolimits_{d=1}^{m}d^{\XX_d}
\Big(\sum\nolimits_{e\div d}\mu\big(\myfrac de \big)
{\Pc(\Xg)(-T^{e})\over dT^{e}}\Big)\usp {\XX_d}}\,,
$$
où $\chic(\Fm )(\alpha,T):=\sum_{i\in\ZZ}\tr(\alpha\actson\Hc^{i}(\Fm ))\,(-T)^{i}$, 
et où $(1^{\XX_1},2^{\XX_{2}},\ldots,m^{\XX_{m}})\vdash m$ est le type de la permutation $\alpha$, 
$\mu(\_)$ est la fonction de Möbius, $\Pc(\_)$ est le polynôme de Poincaré de $\Hc(\_)$, et $(\_)\usp r$ est la factorielle décroissante.}


\smallskip
\noindent -- 
Pour le problème des polynômes de Poincaré, le cas de $\Fm(\Xg)$ relève d'une formule fermée très simple conséquence de la $i$-acyclicité de $\Xg$, tandis que les cas plus généraux des quotients de $\Fm(\Xg)$ par des sous-groupes finis de $\Sm $ sont traités via la formule des caractères \ref{theo-trace-gen}. C'est ainsi que nous procédons pour l'espace des configurations \expression{cycliques} $\Cg\Fm (\Xg):=\Fm (\Xg)/\Cg_{m}$, où $\Cg_{m}:=\langle(1,\ldots,m)\rangle\dans\Sm $,
 et celui des configurations \expression{non ordonnées} $\Bg\Fm (\Xg):=\Fm (\Xg)/\S _{m}$.
 
\smallskip
Notons, $\Pc(\_)$ le polynôme de Poincaré de $\Hc(\_)$, 
$\phi(\_)$ la fonction indicatrice d'Euler,
$\mu(\_)$ la fonction de Möbius, 
et $(\_)\usp r$ la factorielle décroissante. 
Pour tout espace $i$-acyclique $\Xg$, nous prouvons les égalités suivantes.

\begingroup
\parindent0pt\parskip0pt
\slshape
\mynobreak\nobreak
\vskip1em

\noindent\rlap{Théorème (\ref{pol-poincare}):}\hfill
$\displaystyle
\smashbot{{\Pc(\Fm (\Xg))(-T)\over T^{m} }= \Big({\Pc(\Xg)(-T)\over T}\Big)\usp m\,.}$
\hfill\null

\bigskip 
Théorème (\ref{theo-conf-cycliques}):
$$\relax{{\Pc(\Cg\Fm )(-T)\over T^{m} }=
\relax{\relax{1\over m}\sum_{d\div m}\phi(d)\,
d^{m/d}
\Big(\sum_{e\div d}\mu\big(\myfrac de \big)
{\Pc(\Xg)(-T^{e})\over dT^{e}}\Big)\usp {m/d}}
\,.}\postskip0pt$$

Théorème (\ref{theo-conf-symetriques})
$$\mathrigid1mu
\hss\relax{{\Pc(\Bg\Fm (\Xg))(-T)\over T^{m}}=
\decale-5pt{1\over m!}
\mkern-25mu
\vscalesum{1.75}{-1.3pt}_{\vrule height0pt width0pt\lambda:=(1^{\XX_1},\ldots,m^{\XX_{m}})\vdash m}
\mkern-30mu 
\decale3pt{h_{\lambda}}
\prodnl_{d=1}^{m}d^{\XX_d}
\bigg(\sum_{e\div d}\mu\big(\myfrac de \big)
{\Pc(\Xg)(-T^{e})\over dT^{e}}\bigg)\usp {\XX_{d}}}
\,,\hss$$
où  $h_{\lambda}$ est le cardinal de l'ensemble des permutations de $\S _{m}$ dont la décomposition en cycles disjoints est de type $\lambda:=(1^{\XX_1},\ldots,m^{\XX_{m}})\vdash m$.

\endgroup

\medskip
\noindent -- 
L'expression simple du polynôme de Poincaré de $\Fg_{m}(\Xg)$ suggérait l'existence d'une forme de
trivialité cohomologique des projections $\pi_a:\Fba (\Xg)\to\Fa (\Xg)$, et donc la dégénérescence des suites spectrales de Leray associées. 
Nous montrons que c'est en effet le cas lorsque $\Xg$ est $i$-acyclique et localement connexe (Th.~\ref{degen}).

\medskip

\noindent -- Pour la stabilité de représentations, nous prouvons les théorèmes suivants.

\noindent
{\slshape\binspace0
Théorème (\ref{theo-stabilite-BM-pseudo}).
Soit $\Mg$ une pseudovariété connexe orientable de dimension \hbox{$\geq2$}. Pour  $a,i\in\NN$,
la famille de représentations 
$\set\Sm \rep\Hbm ^{i}(\Delta_{?m-a}\Mg^{m};k)/_{m}$
est monotone et stationnaire pour $m\geq4i+4a$, si $d_{\Mg}=2$, et pour $m\geq2i+4a$, si $d_{\Mg}\geq3$. 
Les familles des caractères et des nombres de Betti correspondantes sont (donc) polynomiales
et la famille $\set \Bettibm^{i}({\Delta_{?m-a}\Mg^{m}/\Sm })/_{m}$ est constante sur ces mêmes intervalles de $m$.\parskip2pt

\noindent La proposition \ref{prop-stabilite-betti-dim3-BFm} va encore plus loin pour la famille  $\set \Bettibm^{i}({\Bg\Fm(\Mg)})/_{m}$ en montrant qu'elle est constante pour $m\geq 2i$, si $d_{\Mg}=2$, et pour $m\geq i$, si $d_{\Mg}\geq 3$.

}

\smallskip
Ces théorèmes sont dus à Church  (\cite{chu}, 2012) dans le cas où $\Mg$ est lisse et pour la famille $\set\Sm\rep\Hbm(\Fm(\Mg))/_{m}$. 
Nos énoncés généralisent ceux de Church dans deux directions. Premièrement, en s'affranchissant de l'hypothèse de lissité de $\Mg$, et, deuxièmement, en incorporant les familles des espaces de configuration généralisés $\set\Delta_{?m-a}\Mg^{m}/_{m}$ (singuliers, même si $\Mg$ ne l'est pas) pour lesquelles il n'y avait pas de conjecture.
Notre approche est totalement différente de celle de Church qui s'appuie sur les travaux de Totaro (\cite{tota}) sur la suite spectrale de Leray associée au plongement  $\Fm(\Mg)\hook\Mg^{m}$ lorsque $\Mg$ est lisse.

La stratégie de la preuve a consisté à se restreindre dans un premier temps aux pseudovariétés $\Mg$ qui sont $i$-acycliques (\ref{theo-stabilite-BM-i-acyclique}) et pour lesquelles on peut élaborer une combinatoire basée sur l'exactitude des complexes fondamentaux qui réduit les questions de stabilité et monotonie au cas des familles d'espaces produit $\set\Mg^{m}/_m$ où la réponse est assez simple. La combinatoire repose sur deux foncteurs d'induction dans la catégorie des $\FI$-modules, les foncteurs
$\Ibf^{a},\Thetabf^{a}:\ModFI k\to\ModFI k$ (\ref{foncteur-I-a}) 
qui décalent les $\S_{m-a}$-modules vers des $\Sm$-modules 
et pour lesquels nous avons un contrôle assez précis de la manière dont ils perturbent les rangs de stabilité et monotonie  (thm. \ref{theo-Ind-FI}). 

Le cas où $\Mg$ est une pseudovariété générale est ensuite abordé moyennant la remarque que $\Mg$ se réalise comme différence $\Mg=\Mgge\mmoins\Mgg$, où $\Mgge:=\Mg\times\RR_{\geq0}$ et $\Mgg:=\Mg\times\RR_{>0}$ sont  des pseudovariétés $i$-acycliques. Cette observation simple conduit naturellement à la construction de la suite spectrale $(\IE_{\sigma}(\U^{m})_{r},d_r)$ (\ref{theo-suite-spectrale-basique}) qui converge vers $\Hbm(\Fm(\Mg))$ et dont la page $\IE_1$ concerne uniquement des espaces de configuration pour des espaces $i$-acycliques. Nous avons (\ref{prop-tableau-normaux-ss-basiques}-(\ref{prop-tableau-normaux-ss-basiques-b})):
$$\def\tt{\!{\vrule height0pt depth4pt width0pt}}
\hss\relax{\IEs(\U^{m})_{1}^{p,q}=
\bigoplusnl_{\tau\in\TTT(p+1,m)}
\ind^{\Sm}_{\HHH_{\tau}}
\alt\otimes\Hbm^{Q}(\Fg_{p+1}(\Mgg))
\ \Rightarrow\ 
\Hbm^{i}(\Fm(\Mg))}\hss\,,$$
où $\medmuskip0mu
Q:=i-(m{-}(\pp))\,(d_{\Mg}{-}1)$. Nous y avons noté 
$\TTT(p{+}1,m)$ l'ensemble des tableaux de Young  à $m$ boites et première colonne $(m{-}p,\ldots,m)$, et $\HHH_\tau$ le stabilisateur de $\tau$ dans $\S_{m{-}(p{+}1)}\times\S_{\pp}$ et dont l'action sur $\Hbm(\Fg_{\pp}(\Mg_{>0}))$ est tordue par
la signature $\alt$ de $\S_{m{-}(p{+}1)}$. Nous appelons ces suites spectrales \expression{basiques}. 
Leur construction est compatible aux images-inverses et permet 
d'estimer les rangs de monotonie et stabilité de 
$\set\Hbm^{i}(\Fm(\Mg))/_{m}$ à partir de ceux déjà connus dans le cas $i$-acyclique de $\set\Hbm^{i}(\Fm(\Mgg))/_{m}$. Le foncteur d'induction $\Ibf^{a}$ est ensuite de nouveau utilisé pour atteindre $\set \Hbm^{i}(\Delta_{m-a}\Mg^{m})/_{m}$ et une récurrence descendante basée sur des suites longues de cohomologie fixe le cas de $\set\Hbm^{i}(\Delta_{\leq m-a}\Mg^{m})/$.

L'étude des familles
$\set \Bettibm^{i}({\Bg\Fm(\Mg)})/_{m}$ suit la même approche. Lorsque $\Mg$ est $i$-acyclique,
la formule explicite du caractère $\chic(\Fm (\Xg))$ donnée par le théorème \ref{theo-trace-gen} nous permet d'obtenir une information assez précise de la multiplicité de la représentation triviale de $\Sm$
dans $\Hbm(\Fm(\Mg))$. Nous sommes alors en mesure de prouver le résultat suivant.

\smallskip\noindent
{\slshape
Proposition \ref{theo-stabilite-i-acyclique-betti-BFm}.
Soit $\Mg$ une pseudovariété connexe orientable et $i$-acyclique. La famille
$\set\Bettibm^{i}(\BFm(\Mg;\QQ))/_{m}$ est constante pour chaque $i\in\NN$ et tout $m\geq i$.
}

\smallskip
\`A partir de là, les suites spectrales basiques permettent d'examiner le cas où $\Mg$ est une pseudovariété générale et d'établir
la proposition \ref{prop-stabilite-betti-dim3-BFm} déjà mentionnée.

\medskip
Nous pensons que même s'il reste beaucoup de questions ouvertes, l'utilité des complexes fondamentaux en tant qu'outil d'approche 
combinatoire dans le cas des espaces $i$-acycliques, et l'utilité des suites spectrales basiques pour s'étendre dans le cadre plus large des pseudovariétés, devrait être clair d'après ce travail.
\vskip2em

\fi
\endgroup
\newpage
{\small
\setcounter{tocdepth}{2}
\tableofcontents
}



\setcounter{section}{-1}\begingroup

\def\soustitreni#1{\vskip-\lastskip\vskip1em
\noindent\refstepcounter{subsection}{\bfseries
\thesubsection.
\slshape\boldmath #1.\ }\ignorespaces}
\secskip=1.ex\mou

\def\defaulttheotitlestyle{\slshape}
\let\theotitlestyle\defaulttheotitlestyle
\section{Introduction}
\noindent
Les motivations à l'origine de ce travail ont été les suivantes.
\def\varitemizeseps{\itemsep0pt\mou\parsep0pt\parskip0pt\topsep2pt
\leftmargin1em
\itemindent0em
}
\begin{itemize}
\item La recherche d'une large classe d'espaces  $\Xg$ pour lesquels le polynôme de Poincaré des espaces de configuration $\Fm (\Xg)$ et le caractère de la représentation par permutation de coordonnées du groupe symétrique $\Sm $ sur $H(\Fm(\Xg))$ sont donnés par une formule fermée ne dépendant que de $m$ et du polynôme de Poincaré de $\Xg$. 
\item La généralisation des théorèmes de Church (2012 \cite{chu}) sur la polynomialité de la famille de représentations $\set \Sm \rep H(\Fm(\Xg))/_{m}$ du cas où $\Xg$ est une variété différentiable au cas où $\Xg$ est une variété algébrique complexe.
\item L'étude de la dégénérescence des suites spectrales de Leray associées aux projections $\Fba(\Xg)\to\Fa(\Xg)$.
\item L'étude de toutes ces questions pour les espaces de configuration généralisés $\Delta_{?\ell}\Xg^{m}$ (\cf\ref{intro-ecg}).
\end{itemize}

\penalty-500
\soustitreni{Polynômes de Poincaré}On suppose $m>0$. Pour tout espace topologique $\Xg$, le produit cartésien $\Xg\stimes\Fm (\Xg)$ contient $\Fmm (\Xg)$ comme partie ouverte. Son complémentaire est l'ensemble:
$$\Delta_{m}(\Xg\stimes\Fm (\Xg))=\set (x_0,x_1,\ldots,x_{m})\in\Xg\stimes\Fm (\Xg)\mid x_0\in{\set x_1,\ldots,x_{m}/}/\,,$$
réunion disjointe des sous-espaces fermés:
$$\Delta_{(0,i)}(\Xg\stimes\Fm (\Xg))=\set (x_0,x_1,\ldots,x_{m})\in\Xg\stimes\Fm (\Xg)\mid x_0=x_i/\,,$$
où $i=1,\ldots,m$, clairement  homéomorphes à $\Fm (\Xg)$.
On a donc
$$\Hc(\Delta_{m}(\Xg\stimes\Fm (\Xg)))\sim\Hc(\Fm (\Xg))^{m}\,,
$$
et la décomposition en parties respectivement ouverte et fermée
$$\vadjust{\kern-4pt}\Xg\stimes\Fm (\Xg)=\Fmm (\Xg)\ \sqcup\ \Delta_{m}(\Xg\stimes\Fm (\Xg))\,,\eqno(\diamond)
$$ 
donne lieu à la suite exacte longue de cohomologie (\footnote{\`A coefficients dans un corps de caractéristique nulle.}) à support compact
$$\def\too^#1{\,\mathop{\to}\limits^{#1\; }\,}
\to
\Hc(\Fmm (\Xg))
\too^{\iota}
\Hc(\Xg\times\Fm(\Xg))
\too^{\rho}
\Hc(\Delta_{m}(\Xg\stimes\Fm (\Xg)))
\too^{c}\eqno(\ast)
$$
dont on déduit la relation de récurrence entre caractéristiques d'Euler des cohomologies à support compact
$\chic(\Fmm )=(\chic(\Xg)-m)\cdot\chic(\Fm)$ à l'origine de l'égalité
$$\chic(\Fmm (\Xg))=\chic(\Xg)\usp{\mm}\,,$$
où $(\_)\usp{\mm}$ désigne la factorielle décroissante (\ref{factorielles}). Mais cette égalité ne renseigne pas sur la valeur de chaque nombre de Betti de $\Fmm (\Xg)$ séparément, ce pour quoi il faudrait une certaine forme de scindage de la suite 
$(\ast)$. Cela arrive, par exemple, lorsque $\Xg$ est un groupe de Lie $\Gg$ réel connexe et \emph{non compact}.
Dans ce cas,
le morphisme de restriction 
$$\rho:\Hc(\Gg\stimes\Fm (\Gg))\to\Hc(\Delta_{m}(\Gg\stimes\Fm (\Gg)))
$$
est nul et la suite courte 
$$0\to\Hc(\Fm (\Xg))^{m}[-1]
\to
\Hc(\Fmm (\Xg))
\to
\Hc(\Xg)\otimes\Hc(\Fm (\Xg))\to0\,,\eqno(\ddagger)$$
extraite de $(\ast)$ est exacte.
En effet, le morphisme $\rho$ est somme des restrictions
$$\rho_i:\Hc(\Gg\stimes\Fm (\Gg))\to\Hc(\Delta_{(0,i)}(\Gg\stimes\Fm (\Gg)))\,,$$
et, d'autre part, on dispose de l'application 
$\mathalign{\varphi_i:&\Gg\stimes\Fm (\Gg)&\to&\Gg\stimes\Fm (\Gg)
}$, $(g,\cl x)\mapsto(g^{-1}\cdot x_i,\cl x)$,
qui est un difféomorphisme qui échange les inclusions:
{\small$$
\left(\empileccc{\Gg\stimes\Fm (\Gg)}{\rotatebox{90}{$\dans$}}{\Delta_{(0,i)}(\Gg\stimes\Fm (\Gg))}\right)
\mathop{\longleftrightarrow}_{\varphi_i}
\left(\empileccc{\Gg\stimes\Fm (\Gg)}{\rotatebox{90}{$\dans$}}{\set e/\stimes\Fm (\Gg)}\right)
\,,
$$}\relax
de sorte que l'annulation de $\rho_i$ équivaut à l'annulation du morphisme de restriction $\Hc(\Gg\times\Fm (\Gg))=\Hc(\Gg)\otimes\Hc(\Fm (\Gg))\to
\Hc(\set e/)\otimes\Fm (\Gg)$, elle-même équivalente à l'annulation de
la restriction $\Hc(\Gg)\to\Hc(\set e/)$, autrement dit à la non compacité de $\Gg$.

\comment
Or, par Künneth, on a le diagramme commutatif
$$\xymatrix@R=6mm{
\kern-5mm\Hc\big(\Gg\times\Fm (\Gg)\big)\ar[r]\ar@<-1.2em>@{<-}[d]_(0.53){\simeq}&\Hc\big(\set e/\times\Fm (\Gg)\big)\kern0.45cm\ar@<-0.8em>@{<-}[d]_(0.53){\simeq}\\
\Hc(\Gg)\otimes\Hc(\Fm (\Gg))\ar[r]^(0.48){r\otimes\id}&\Hc(\set e/)\otimes\Hc(\Fm (\Gg))
}$$
où $r:\Hc(\Gg)\to\Hc(\set e/)$ est le morphisme de restriction. Ce morphisme est évidemment nul lorsque $\Gg$ est connexe et non compact.
\endcomment

\smallskip
L'exactitude de la suite courte $(\ddagger)$ donne alors la relation de récurrence
(\footnote{On note 
$\Pc(\Mg):=\sum_{i\in\NN}\dim_{\QQ}\Hc^{i}(\Mg)\,T^{i}$
 (resp. $\Pr(\Mg):=\sum_{i\in\NN}\dim_{\QQ}\Hr^{i}(\Mg)\,T^{i}$) le polynôme de Poincaré pour la cohomologie à support compact
  (resp. ordinaire) de $\Mg$.})
$
\Pc(\Fmm (\Gg))(T)=(\Pc(\Gg)(T)+mT)\cdot\Pc(\Fm (\Gg))(T)\,,
$
dont on déduit
$$
\Pc(\Fmm (\Gg))(T)=
\relax{\prodnl_{i=0,\ldots,m}(\Pc(\Gg)(T)+i\,T)}\,,\eqno(\dagger)$$
et, par dualité de Poincaré,
$$
\Pr(\Fmm (\Gg))(T)=\prodnl_{i=0,\ldots,m}(\Pr(\Gg)(T)+i\,T^{\,\dim\Gg-1})\,,\eqno(\dagger\dagger)$$ 
ce qui est le type de formule fermée que nous avions en vue.

\soustitreni{Les espaces de configuration généralisés}
Dans\label{intro-ecg} nos recherches nous avons aussi été guidés par une autre nécessité: celle de trouver un cadre aussi symétrique que possible pour préserver l'action du groupe des permutations $\Sm$ sur les coordonnées de $\Fm (\Xg)$, ce qui n'est pas le cas de la décomposition $(\diamond)$. Cela nous a conduit à nous intéresser aux espaces de configuration \expression{généralisés} suivants. 
Pour $0\leq\ell\leq m\in\NN$, on pose
$$\begin{casesalign}
\Delta_{\leq \ell}\Xg^{m}&:=&\set(x_1,\ldots,x_{m})\in\Xg^{m}\mid\card{\set x_1,\ldots,x_{m}/}\leq \ell/\,,\\\noalign{\kern1pt}
\hfill \Delta_{\ell}\Xg^{m}&:=&\Delta_{\leq\ell}\Xg^{m}\mmoins\Delta_{\leq\ellmo}\Xg^{m}\,,\quad\Fm(\Xg):=\Delta_{m}\Xg^{m}\,.\hfill
\end{casesalign}$$
La décomposition en parties $\S _{m}$-stables respectivement ouverte et fermée :
$$\Delta_{\leq \ell}\Xg^{m}=
\Delta_{\ell}\Xg^{m}\ \sqcup\ \Delta_{\leq \ellmo}\Xg^{m}\eqno(\diamond\diamond)$$
donne alors lieu à la suite courte de $\Sm$-modules (à priori non exacte) 
$$0\to
\Hc(\Delta_{\leq\ellmo}\Xg^{m})[-1]
\to\Hc(\Delta_{\ell}\Xg^{m})\to
\Hc(\Delta_{\leq\ell}\Xg^{m})\to0\,,
\eqno(\ddagger\ddagger)$$
extraite de la suite  longue de cohomologie à support compact associée à  $(\diamond\diamond)$.

\smallskip
En cherchant à montrer l'exactitude de ($\ddagger\ddagger$) lorsque $\Xg$ est un groupe de Lie non compact, nous avons réalisé que ce n'était pas tant le fait que $\Xg$ possède une structure de groupe, mais plutôt que sa cohomologie intérieure est nulle, qui est à l'origine du phénomène de scindage. Le théorème \expression{de scindage} \ref{theo-scindage} énonce cette observation sous l'a forme de l'implication $(A)\Rightarrow(B)$ des propriétés suivantes.
\def\varitemizeseps{\itemsep1pt\mou\parsep0pt\parskip0pt\topsep4pt\leftmargin1cm}\begin{itemize}\mynobreak\nobreak
\item[$(A)$]
La cohomologie \expression{intérieure}  de $\Xg$, \idest
l'image de l'application naturelle  $\Hc(\Xg)\to\Hr(\Xg)$,  est nulle.

\item[$(B)$] Pour tous $0\leq\ell\leq m$, les suites courtes $(\ddagger)$ et $(\ddagger\ddagger)$ sont exactes.
\end{itemize}
Des propriétés qui sont même équivalentes lorsque $\Xg$ est une variété topologique orientable.

\soustitreni{Les espaces $i$-acycliques}\`A partir du théorème de scindage, une partie importante de notre travail va se concentrer dans les espaces vérifiant la propriété~$(A)$. Ce sont les espaces que nous appelons \expression{$i$-acycliques}. 

Il existe des familles assez larges d'exemples de tels espaces:
\def\varitemizeseps{\itemsep0.5pt\mou\parsep0pt\parskip0pt\topsep0pt}\begin{itemize}
\item Les espaces acycliques (p.e. contractiles) non compacts.
\item Les ouverts des espaces $i$-acycliques.
\item Les groupes de Lie $\Gg$ tels que \smash{$\Hc^{0}(\Gg)=0$.}
\item Tout produit $\Xg\stimes\Yg$, où $\Xg$ est $i$-acyclique et $\Yg$ est quelconque.
\comment\item Les sommes amalgamées de variétés topologiques $i$-acycliques.\endcomment
\end{itemize}

\soustitreni{Le complexe fondamental de $\Xg$ pour $\Delta_{\leq\ell}\Xg^{m}$}
Dans \ref{complexe-fondamental}, on concatène les suites $(\ddagger\ddagger)$ pour construire le complexe de $\S _{m}$-modules gradués
$$
0\to\Hc(1)[-\ell+1]\to\cdots\to
\Hc(\ellmo)[-1]\to
\Hc(\ell)\to
\Hc(\Delta_{\leq\ell}\Xg^{m})\to0$$
avec $\Hc(a):=\Hc(\Delta_{a}\Xg^{m})$, que nous appelons
\expression{complexe fondamental de $\Xg$ pour $\Delta_{\leq\ell}\Xg^{m}$}. Son principal intérêt réside dans l'assertion (\ref{theo-complexe-exact-a}) suivante.

\nopoint\begin{theo*}[(\ref{theo-complexe-exact})]
\def\varlistskips{\topsep0cm\itemsep0.5pt\mou\parskip0pt\mou}
\begin{enumerate}
\item Les complexes fondamentaux d'un espace $i$-acyclique sont  exacts.
\item Une variété topologique orientable est $i$-acyclique, si et seulement si, ses complexes fondamentaux sont exacts.
\end{enumerate}\end{theo*}

\comment
La $i$-acyclicité apparaît donc comme caractéristique du fait que les complexes fondamentaux sont des résolutions de $\S _{m}$-module de $\Hc(\Delta_{\leq\ell}\Xg^{m})$. 
\endcomment

Les complexes fondamentaux sont donc des résolutions de $\Sm$-modules 
de $\Hc(\Delta_{\leq\ell}\Xg^{m})$ et permettent, lorsque $\Hc(\Xg)<\infty$, de réduire la détermination du caractère 
$\Hc(\Delta_{\leq\ell}\Xg^{m})$
à ceux de $\Hc(\Delta_{a}\Xg^{m})$ pour $a\leq\ell$. 
\comment
Ces résolutions constituent alors un outil commode pour l'étude récursive de propriétés communes aux espaces de configuration généralisés (on en verra plusieurs exemples), et aussi pour l'élaboration d'algorithmes de calcul. 
\endcomment

\soustitreni{Les polynômes de Poincaré de $\Delta_{\ell}\Xg^{m}$ et de $\Delta_{\leq\ell}\Xg^{m}$}Les remarques des premiers paragraphes, étant basées sur le scindage de la suite longue $(\ddagger)$, s'appliquent par conséquent lorsque $\Xg$ est $i$-acyclique. On a donc:

\begin{prop*}[(\ref{pol-poincare})]Si\label{intro-pol-poincare} $\Xg$ est $i$-acyclique et que $\Hc(\Xg)<+\infty$, le polynôme de Poincaré $\Pc(\Fm (\Xg))$ est le polynôme
$$\Pc(\Fm (\Xg))(T)=
\prodnl_{i=0}^{\mmo}\big(\Pc(\Xg)(T)+i\cdot T\big)\,.$$
\end{prop*}

\medskip
\`A partir de là, le passage de $\Fm (\Xg)$ à $\Delta_{\ell}\Xg^{m}$ est assez simple dans la mesure où $\Delta_{\ell}\Xg^{m}$ admet une décomposition en parties ouvertes homéomorphes à $\Fl(\Xg)$ indexée par les partitions de l'intervalle $\iii[1,m]$  en $\ell$ parties non vides. La proposition suivante établit alors que $\Pc(\Delta_{\ell}\Xg^{m})$ seul dépend de $\Pc(\Xg)$, le lien étant donné par un certain polynôme \expression{universel} de $\ZZ[P,T]$. 

\begin{prop*}[(\ref{=l-m-poly-univ})]Si $\Xg$ est $i$-acyclique et tel que $\Hc(\Xg)<+\infty$, le polynôme de Poincaré $\Pc(\Delta_{\ell}\Xg^{m})$ s'obtient en évaluant en $P:=\Pc(\Xg)$ le
polynôme homogène de degré $\ell$ de $\ZZ[P,T]$ :
$$\Qg_{\ell}^{m}(P,T)=\mParties m\ell\cdot\prodnl_{i=0,\ldots,\ellmo}(P+i\,T)\,.$$
Ici, $\Pgoth_{\ell}(m)$ désigne l'ensemble des partitions de $\iii[1,m]$ en $\ell$ parties non vides. 
\comment Le nombre $|\Pgoth_{\ell}(m)|$ est donné par le nombre de Stirling de seconde espèce:
$$|\Pgoth_{\ell}(m)|=\Parties m\ell:={1\over \ell! }\sumnl_{j=0}^{\ell}(-1)^{\ell-j}\binome{\ell}{j}j^{m}\,.
$$\endcomment
\end{prop*}

\bigskip
Le polynôme de Poincaré de $\Pc(\Delta_{\leq\ell}\Xg^{m})$ est ensuite obtenu comme somme alternée des $\P(\Hc(\Delta_{a}\Xg^{m})[a-\ell])$, grâce aux complexes fondamentaux.
\begin{prop*}[(\ref{pol-poin-gen})]Soit $\Xg$ un espace $i$-acyclique.
Le polynôme de Poincaré $\Pc(\Delta_{\leq \ell}\Xg^{m})$ est le polynôme homogène de $\ZZ[\Pc(\Xg),T]$, de degré $\ell$, donné par la somme alternée {\rm(\cf\ref{connexes})}
$$
\Pc(\Delta_{\leq\ell}\Xg^{m})=\sumnl_{0\leq a<\ell}(-1)^{a}\cdot\mParties m{\ell-a}\cdot\Pc(\Fg_{\ell-a}(\Xg))\cdot T^{a}\,.
$$
\end{prop*}

\soustitreni{Le caractère de $\Hc(\Delta_{?\ell}\Xg^{m})$}
Pour $\Zg\dans\Xg^{m}$ stable sous l'action de $\S _{m}$ et tel que $\dim\Hc(\Zg)<\infty$, notons $\chic(\Zg;i)$ le caractère du $\S _{m}$-module $\Hc^{i}(\Zg,\QQ)$. Le théorème suivant est conséquence immédiate de l'exactitude des complexes fondamentaux.
\begin{theo*}[(\ref{theo-caracteres})]Soit\label{intro-theo-caracteres} $\Xg$ un espace $i$-acyclique. \halfdisplayskips
\begin{enumerate}
\item[\rm\ref{theo-caracteres-a})] Le caractère $\chic(\Fm (\Xg);i)$ du $\Sm$-module $\Hc^{i}(\Fm (\Xg))$ vérifie 
$$\chic(\Fm (\Xg);i)=\chic(\Xg^{m};i)+\chic(\Delta_{\leq \mmo}\Xg^{m};i-1)\,.$$
\item[\rm\ref{theo-caracteres-b})] Le caractère $\chic(\Delta_{\leq\ell}\Xg^{m};i)$ du $\Sm$-module $\Hc^{i}(\Delta_{\leq\ell}\Xg^{m})$ vérifie 
$$\def\tt{\vrule height7pt width0pt }
\chic(\Delta_{\leq\ell}\Xg^{m};i)=
\sum_{\tt0\leq a<\ell}\hskip-1.25ex
\sum_{\decale-13pt{\tt\lambda\in\Y_{\ell-a}(m)}}\hskip-1ex(-1)^{a}\,
\Ind_{\Glambda }^{\Sm}\chic(\Fg_{\ell-a}(\Xg);i-a)\,.
$$
où $\Y_{\ell-a}(m)$ est l'ensemble des décompositions de $m$ en $\ell{-}a$ entiers positifs, et où,
si $\lambda=(\lambda_1,\ldots,\lambda_{\ell-a})=(\lambda_1^{\mu_1},\ldots,\lambda_r^{\mu_r})$, on a noté
$$\Glambda :=N_{\Sm}(\Plambda)/\Plambda= \S_{\mu_1}\stimes\cdots\stimes\S_{\mu_r}$$ 
où $\Plambda =\S_{\lambda_{1}}\stimes\cdots\stimes\S_{\lambda_{\ell-a}}$.
\end{enumerate}
\end{theo*}

\smallskip
Le caractère $\chic(\Delta_{\leq\ell}\Xg^{m})$ apparaît ainsi comme combinaison des caractères $\chic(\Fa (\Xg))$ pour $a\leq\ell$. D'autre part, grâce à l'exactitude des suites $(\ddagger\ddagger)$, le caractère  $\chic(\Fa (\Xg))$ est somme de $\chic(\Delta_{\leq a-1}\Xg^{a})$, avec (donc) $a-1<\ell$, et de $\chic(\Xg^{a})=\chi_{\Sa }(\Hc(\Xg)^{\otimes a})$, caractère bien connu d'après Macdonald \cite{mac} (\cf \ref{theo-trace-gen-Xm}). On a donc tous les ingrédients pour un algorithme de calcul des caractères de $\Hc(\Fm (\Xg))$ et $\Hc(\Delta_{?\ell}\Xg^{m})$ à partir uniquement de la connaissance de $\chic(\Xg^{\ell})$ pour $\ell\leq m$.
Cette observation est développée dans la section \ref{inductions-iterees} qui introduit certains opérateurs d'induction dans les groupes de Grothendieck des catégories de représentations des groupes symétriques
$$
\Thetabf_{\ell}^{m},\Ibf ^{m}_{\ell}:K_0(\Mod(k[\S_{\ell}]))\fonct K_0(\Mod(k[\Sm ]))
$$ 
qui permettent dans le théorème suivant de relier les caractères des cohomologies des espaces de configuration généralisés à ceux des produits cartésiens. 

\begin{theo*}[\ref{theo-caracteres-devissage}-(\ref{theo-caracteres-devissage-a})]Soit $\Xg$ un espace $i$-acyclique tel que $\dim\Hc(\Xg)<\infty$. 
\begin{enumerate}
\mynobreak\nobreak\item Pour tout $m\geq\ell >0$ et tout $i\in\NN$, on a
$$\mathalign{
\ \llap{\rm i)\ }\hfill\chic(\Fm(\Xg);i)&=&\sumnl_{0\leq a<m}\Thetabf^{m}_{m-a}\big(\chic(\Xg^{m-a};i-a)\big)\hfill\\\noalign{\kern4pt}
\ \llap{\rm ii)\ }\hfill\chic(\Delta_{\ell}\Xg^{m};i)&=&
\Ibf ^{m}_{\ell}\Big(\sumnl_{0\leq a<\ell}\Thetabf^{\ell}_{\ell-a}\big(\chic(\Xg^{\ell-a};i-a)\big)\Big)\hfill\\\noalign{\kern4pt}
\mathrigid0mu
\ \llap{\rm iii)\ }\chic(\Delta_{\leq\ell}\Xg^{m};i)&=&
\mathrigid0mu
\sum_{0\leq b<\ell}(-1)^{b}\ 
\Ibf ^{m}_{\ell-b}\Big(\sumnl_{\vrule height18pt width0pt\hskip-0.5cm 0\leq a<\ell-b}\hskip-1.7em\Thetabf^{\ell-b}_{\ell-b-a}\big(\chic(\Xg^{\ell-b-a};i-b-a)\big)\Big)}
$$
\end{enumerate}
\end{theo*}
A partir de là nous avons suivi deux voies de recherche, celle de l'étude de la polynomialité des familles de caractères $\set\chi(\Delta_{m-a}\Xg^{m},i)/_{m}$ et celle de la détermination explicite des caractères $\chic(\Fm(\Xg))$. Les sections \ref{stabilite} à \ref{stabilite-des-familles} sont consacrées
à la première question et la section \ref{caractere} à la seconde.

\soustitreni{Les familles de représentations
$\set\Sm \rep\Hr^{i}(\Delta_{? m-a}\Xg^{m})/_{m}$}Bien de recherches sur les espaces de configuration concernent le comportement asymptotique d'invariants cohomologiques. Par exemple, le résultat pionnier d'Arnold (1970~\cite{arnold}) qui établit que pour $i\in\NN$, la suite $\mathrigid0mu
\set{\Betti^i(\Fm(\CC)/\S _{m}})/_{m} $ est stationnaire, ou, plus récemment, celui de Church (2012 \cite{chu}) qui montre que si $\Xg$ est une variété différentielle connexe orientable et que $d_{\Xg}\geq2$, les multiplicités des facteurs irréductibles des représentations de $\Sm$ dans $\Hr^{i}(\Fm(\Xg))$ sont stationnaires (dans un sens à préciser, \cf\ref{monotonie-et-stabilite}).

Les complexes fondamentaux s'avèrent particulièrement commodes pour aborder ces questions de nature qualitative. Ils indiquent aussi que la direction à suivre pour les généraliser doit concerner les familles  $\set\Delta_{?m-a}\Xg^{m}/_{m}$ pour $a\in\NN$ fixe.
Les sections \ref{stabilite} à \ref{stabilite-des-familles} sont  consacrées à ces questions. Nous y rappelons la théorie des $\FI$-modules de Church et Frab~(\cite{cef}). Un $\FI$-module est un foncteur covariant de la catégorie $\FI$ des ensembles finis et des applications injectives vers la catégorie des $k$-espaces vectoriels, il peut être représenté par la donnée d'une famille dénombrable d'applications $k$-linéaires $\V=\set \phi_m:V_m\to V_\mm/_{m\in\NN}$ où $V_{m}$ est un $\Sm$-module et $\phi_m$ est compatible aux actions de $\Sm$ et $\Smm$. Nous rappelons les concepts importants de  rang de \expression{monotonie}, de \expression{stabilité} et de \expression{polynomialité} (\cf\ref{comment-poids-degre}) de $\FI$-modules et nous prouvons la généralisation suivante des résultats de Church.

\begin{theo*}[(\ref{theo-stabilite-BM-i-acyclique})]\binspace2
Soit\label{intro-theo-stabilite-BM-i-acyclique} $\Xg$ une pseudovariété $i$-acy\-clique, connexe orientable et de dimension $d_{\Xg}\geq2$. Pour  $a,i\in\NN$,
la famille 
$\set\Sm \rep\Hbm ^{i}(\Delta_{?m-a}\Xg^{m})/_{m}$
est monotone pour $m\geq i+a$ et est monotone et stable pour $m\geq4i+4a$, si $d_{\Xg}=2$, et pour $m\geq2i+4a$, si $d_{\Xg}\geq3$. 
\end{theo*}

La démonstration est basée sur le théorème \ref{theo-caracteres-devissage}-(\ref{theo-caracteres-devissage-a}) déjà mentionné.
Dans le cas du $\FI$-module $\V(i):=\set p_m^{*}:\Hbm ^{i}(\Fm(\Xg))\to\Hbm ^{i}(\Fmm(\Xg))/_{m}$ où $p_m^{*}$ est dual de l'intégration sur les fibres de la projection  sur les $m$ premières coordonnées $p_m:\Fmm(\Xg)\to\Fm(\Xg)$,  l'idée est de montrer que les foncteurs $\Thetabf^{m}_{m-a}:\Mod(k[\S_{m-a}])\fonct\Mod(k[\Sm])$ 
définissent un foncteur  dans la catégorie de $\FI$-modules $\Thetabf^{a}:\ModFI k\fonct\ModFI k_{\geq 2a}$ donnant lieu à une égalité de $\FI$-modules virtuels
$\binspace1
\V(i)=\sumnl_{0\leq a<m}(-1)^{a}\,
\Thetabf^{a}\big(\W(\dXg(m-a)-i+a)\big)\,,$
où $\W(j)$ est le $\FI$-module $\set p_m^{*}:\Hbm^{j}(\Xg^{m})\to \Hbm^{j}(\Xg^{\mm})/$ dont les rangs sont faciles à déterminer. Il est alors essentiel de comprendre comment $\Thetabf^{a}$ perturbe les rangs de monotonie et stabilité $\rkm$ et $\rks$.
Le théorème \ref{theo-stabilite-BM-i-acyclique} ci-dessus est alors corollaire du résultat suivant où $\rkms:=\sup\set\rkm,\rks/$.

\noendpoint\begin{theo*}[(\ref{theo-Ind-FI})]\def\Rg{\Ibf}
\def\varlistskips{\topsep1pt\itemsep0.pt\mou\parskip0pt\mou}
\begin{enumerate}\let\Ibf\Thetabf
\item Le foncteur
$\mathrigid2mu
\Rg^{a}:\ModFI k\fonct\ModFI k_{\geq 2a}$
est covariant, additif, exact.
\item  Si $\V$ est (de type fini) engendré en degrés $\leq d$, le $\FI$-module $\Rg^{a}(\V)$ est (de type fini) engendré en degrés $\leq \sup(d+a,2a)$.
\item On a 
$\rkms(\Rg^{a}\V)\leq \rkms(\V)+4a
\text{ et }
\rkm(\Rg^{a}\V)\leq \rkm(\V)+a$.
\end{enumerate}
\end{theo*}


\soustitreni{Suites spectrales basiques}Dans le but d'étendre la portée du  théorème \ref{theo-stabilite-BM-i-acyclique} 
aux espaces non $i$-acycliques, nous avons introduit dans la section 
\ref{cohomologie-Borel-Moore} la \expression{suite spectrale basique}. Associée à un espace localement compact $\Mg$ de dimension cohomologique finie, elle converge vers $\Hbm(\Fm(\Mg))$ et a la propriété remarquable de ne faire intervenir que des espaces de configuration associés à l'espace $i$-acyclique
 $\Mgg:=\Mg\times\RR_{>0}$.

\begin{theo*}[(\ref{theo-suite-spectrale-basique}, \ref{prop-tableau-normaux-ss-basiques}-(\ref{prop-tableau-normaux-ss-basiques-b}))]Soit $\Mg$ une pseudovariété  orientée de dimension $\dMg$. La suite spectrale $\IEs(\U^{m})$ converge en tant que suite spectrale de complexes de $\Sm $-modules vers le $\Sm $-module bi-gradué associé 
au $\Sm $-module gradué filtré $\Hbm(\Fm(\Mg))[1{-}m]$.
Pour  $i\in\ZZ$, on a:
$$\halfdisplayskips
\def\tt{{\vrule height0pt depth8pt width0pt}}
\hss\relax{\IEs(\U^{m})_{1}^{p,q}=
\bigoplus_{\tau\in\TTT(p+1,m)}
\ind^{\Sm}_{\HHH_{\tau}}
\
\alt\otimes\Hbm^{Q}(\Fg_{p+1}(\Mgg))
\ \Rightarrow\ 
\Hbm^{i}(\Fm(\Mg))}\hss$$
où $\mathrigid1mu
q=i+(m{-}(\pp))$ et  $\medmuskip0mu
Q:=i-(m{-}(\pp))\,(d_{\Mg}{-}1)$, où
$\TTT(p{+}1,m)$ est l'ensemble des tableaux de Young  à $m$ boites et première colonne $(m{-}p,\ldots,m)$, et où $\HHH_\tau$ est le stabilisateur de $\tau$ dans $\S_{m{-}(p{+}1)}\times\S_{\pp}$ dont l'action sur $\Hc(\Fg_{\pp}(\Mg_{>0}))$ est tordue par
la signature $\alt$ de $\S_{m{-}(p{+}1)}$.\end{theo*}

Une propriété remarquable des suites spectrales basiques est leur compatibilité aux morphismes  $p_m^{*}:\Hbm ^{i}(\Fm(\Mg))\to\Hbm ^{i}(\Fmm(\Mg))$. Cela nous a permis de généraliser le théorème \ref{theo-stabilite-BM-i-acyclique} (p.~\pageref{intro-theo-stabilite-BM-i-acyclique}) aux pseudovariétés, en particulier, aux variétés algébriques complexes. L'énoncé suivant ne diffère de \ref{theo-stabilite-BM-i-acyclique} que par le fait que $\Mg$ n'est plus supposée $i$-acyclique et aussi par la perte de  l'estimation du rang de monotonie.

\begin{theo*}[(\ref{theo-stabilite-BM-pseudo})]
Soit $\Mg$ une pseudovariété connexe orientable de dimension  $\dMg\geq2$. Pour  $a,i\in\NN$ fixés, la famille
$\set\Sm \rep\Hbm ^{i}(\Delta_{?m-a}\Mg^{m})/_{m}$
est monotone et stable pour $m\geq4i+4a$, si $\dMg=2$, et pour $m\geq2i+4a$, si $\dMg\geq3$. Les familles des caractères et des nombres de Betti correspondantes sont (donc) polynomiales
et la famille $\set \Bettibm^{i}({\Delta_{?m-a}\Mg^{m}/\Sm })/_{m}$ est constante, sur les mêmes intervalles.\end{theo*}

\soustitreni{Le calcul explicite du caractère de $\Hc(\Fm (\Xg))$}Dans la section~\ref{caractere}, on revient sur les espaces $i$-acycliques. On y considère suivant Macdonald (\cite{mac}), la \expression{série des caractères} de $\Zg\dans\Xg^{m}$ qui vaut
$$\chic(\Zg)(\alpha,T):=\sumnl_{i\in\ZZ}(-1)^{i}\tr(\alpha\actson\Hc^{i}(\Zg))\,T^{i}\,,\qquad\forall \alpha\in\S _{m}\,,$$
et on mène à terme le calcul des séries $\chic(\Fm (\Xg))(\alpha,T)$ 
grâce notamment à l'exactitude des complexes fondamentaux de $\Xg$. On prouve:

\begin{theo*}[(\ref{theo-trace-gen})]Soit\label{intro-theo-trace-gen} $\Xg$ un espace $i$-acyclique. Pour  $\alpha\in\S _{m}$, on a
$$\relax{{\chic(\Fm (\Xg))(\alpha,T)\over T^{m}}=\,\prod_{d=1}^{m}d^{\XX_d}
\Big(\sum_{e\div d}\mobius de
{\chic(\Xg)(\1,T^{e})\over dT^{e}}\Big)\usp {\XX_d}}\,,
$$
où $(1^{\XX_1},2^{\XX_{2}},\ldots,m^{\XX_{m}})\vdash m$ est le type de la permutation $\alpha$, où $\mu(\_)$ est la fonction de Möbius et où $(\_)\usp r$ est la factorielle décroissante (cf.~\ref{factorielles}).
\end{theo*}

\soustitreni{Les polynômes de Poincaré des espaces quotients $\Fm (\Xg)/\Hg$}Le théorème  précédent s'applique aussitôt pour donner les dimensions des sous-espaces invariants $\Hc^{i}(\Fm (\Xg))^{\Hg}$, quel que soit le sous-groupe $\Hg\dans\Sm$. Les polynômes de Poincaré des quotients $\Fm (\Xg)/\Hg$ en résultent. 

\penalty-1000
La section~\ref{quotients} illustre le procédé en déterminant le polynôme de Poincaré de l'espace des configurations \expression{cycliques} $\Cg\Fm (\Xg):=\Fm (\Xg)/\Cg_{m}$, où $\Cg_{m}:=\langle(1,\ldots,m)\rangle\dans\Sm $, et celui  de l'espace des configurations \expression{non ordonnées} $\BFm (\Xg):=\Fm (\Xg)/\S _{m}$. 
Pour $\Xg$ $i$-acyclique, on obtient :
\def\varitemizeseps{\leftmargin1.5ex\labelwidth1em \labelsep0em}\begin{itemize}\slshape\halfdisplayskips
\vskip1.5pt
%
\item Théorème (\ref{theo-conf-cycliques}):
$$\relax{{\Pc(\Cg\Fm )(-T)\over T^{m}}=\relax{1\over m}\sum_{d\div m}\phi(d)\,
d^{m/d}
\Big(\sum_{e\div d}\mu({d/ e})
{\Pc(\Xg)(-T^{e})\over dT^{e}}\Big)\usp {m/d}}
\,.\postskip-1em$$

\item Théorème (\ref{theo-conf-symetriques})
$$\mathrigid1mu
\hss\relax{{\Pc(\BFm (\Xg))(-T)\over T^{m}}=\decale-5pt{1\over m!}
\mkern-30mu
\vscalesum{1.75}{-1.3pt}_{\vrule height0pt width0pt\lambda:=(1^{\XX_1},\ldots,m^{\XX_{m}})\vdash m}
\mkern-35mu 
\decale3pt{h_{\lambda}}
\prodnl_{d=1}^{m}d^{\XX_d}
\bigg(\sum_{e\div d}\mobius de
{\Pc(\Xg)(-T^{e})\over dT^{e}}\bigg)\usp {\XX_{d}}}
\,,\hss$$
où  $h_{\lambda}$ est le cardinal de l'ensemble des permutations de $\S _{m}$ dont la décomposition en cycles disjoints est de type $\lambda:=(1^{\XX_1},\ldots,m^{\XX_{m}})\vdash m$.
\end{itemize}

\soustitreni{Rangs de stabilité des familles $\set\Bettibm^{i}(\Bg\Fm(\Mg))/_{m}$}Les formules explicites de caractères \ref{theo-trace-gen} (p.~\pageref{intro-theo-trace-gen}) permettent aussi de raffiner la détermination des multiplicités de la représentation triviale dans $\Hbm(\Fm(\Xg))$. On a:

\begin{prop*}[(\ref{theo-stabilite-i-acyclique-betti-BFm})]
Soit $\Xg$ une pseudovariété $i$-acyclique de type fini telle que $\dim\Hc^{\dXg}(\Xg;\QQ)\leq1$. 
Alors, pour $i\in\NN$, la famille $\set\Betti^{i}(\BFm(\Xg;\QQ))/_{m}$ est constante pour tout $m\geq i$. 
\end{prop*}

Ce résultat s'étend ensuite aux pseudovariétés connexes orientées à l'aide des suites spectrales basiques. La proposition suivante termine  alors la section sur les question de stabilité. (\footnote{La section \ref{stabilite-Betti-Delta-cas-general}, rajoutée tardivement, améliore aussi les bornes pour les rangs de stabilité des familles $\set \Bettibm^{i}(\Delta_{?m-a}\Mg^{m}{/}\Sm)/_m$ données dans le théorème \ref{theo-stabilite-BM-pseudo}.})

\begin{prop*}[(\ref{prop-stabilite-betti-dim3-BFm})]Soit $\Mg$ est une pseudovariété connexe orientée. La famille $\set\Bettibm^{i}(\Bg\Fm(\Mg))/_{m}$ est constante pour $m\geq 2i$, si $d_{\Mg}=2$, et pour $m\geq i$ si $d_{\Mg}\geq 3$.
\end{prop*}

\soustitreni{La dégénérescence des suites spectrales de Leray}Pour $a,b\in\NN$, notons $\pi_a:\Xg^{b+a}\to\Xg^{a}$ la projection sur les $a$ dernières coordonnées. Sa restriction à $\Fba(\Xg)$ est la fibration $\pi_a:\Fba (\Xg)\to\Fa (\Xg)$ (\footnote{Généralement non localement triviale, sauf si $\Xg$ est une variété topologique.}) de fibres de la forme $\Xg\mmoins\, \ag$ où $\ag$ désigne un sous-ensemble de $\Xg$ de cardinal $a$. 

Lorsque $\Xg$ est $i$-acyclique et que $\Hc(\Xg)<+\infty$, l'expression du polynôme de Poincaré de la proposition \ref{pol-poincare} (p.~\pageref{intro-pol-poincare}) montre que l'on a 
$$\Pc(\Fba (\Xg))=\Pc(\Fa (\Xg))\cdot\Pc(\Fg_{b}(\Xg\mmoins\, \ag))\,,$$
ce qui suggère une certaine forme de trivialité cohomologique pour la fibration $\pi_a$. C'est en effet le cas et c'est le sujet de la section \ref{Leray}.

\smallskip
Notons, plus généralement $\GaDelta_{?\ell}\Xg^{m}$, pour $0<a\leq\ell\leq m\in\NN$, l'ouvert des
$m$-uplets de $\Delta_{?\ell}\Xg^{m}$ dont les $a$ dernières coordonnées sont deux à deux distinctes, soit:
$$\GaDelta_{?\ell}\Xg^{m}:=\Delta_{?\ell}\Xg^{m}\cap\big(\Xg^{m-a}\stimes\Fa (\Xg)\big)\,.$$
La section \ref{degenerescence} est consacrée à l'étude des suites spectrales de Leray associées à
l'application $\pi_a:\GaDelta_{?\ell}\Xg^{m}\to\Fa(\Xg)$ dont la fibre au-dessus de $\cl x\in\Fa(\Xg)$ est $\Delta_{?\ell}(\Xg^{m-a}\times\cl x)$. 
Cela nous a emmené à nous intéresser également à la cohomologie à support $\pia$-propre  que nous notons
$\Hpi(\GaDelta_{\ell}\Xg^{m})$, et aussi aux faisceaux de cohomologie à support $\pia$-propre:
$$\cHpi^{i}(\GaDelta_{?\ell}\Xg^{m}):=\IR^{i}\pi_{a!}(\fs k_{\GaDelta_{?\ell}\Xg^{m}})\,,$$
Le principal résultat concernant ces faisceaux est le suivant.

\begin{theo*}[(\ref{autre-scindage}-(\ref{autre-scindage-c}))]Si $\Xg$ $i$-acyclique et localement connexe, les faisceaux  $\cHpi^{i}(\GaDelta_{?\ell}\Xg^{m})$ sont constants sur les composantes connexes de $\Fa(\Xg)$.
\end{theo*}

\`A partir de là, l'étude des suites spectrales de Leray pour les cohomologies $\Hc(\GaDelta_{?\ell}\Xg^{m})$ et $\Hpi(\GaDelta_{?\ell}\Xg^{m})$ se simplifie et nous montrons le théorème suivant.

\begin{theo*}[(\ref{degen})]
Soient $a\leq\ell\leq m\in\NN$. Soit $\Xg$ un espace $i$-acyclique localement connexe. 
\comment Les applications $\pi_a:\GaDelta_{?\ell}\Xg^{m}\to\Fa(\Xg)$ sont des fibrations localement triviales. \endcomment
Notons
$$
(\IE_r(\GaDelta_{?\ell}\Xg^{m})_{\rmc},d_r)
\text{\quad et\quad }
(\IE_r(\GaDelta_{?\ell}\Xg^{m})_{\cv},d_r)
\eqno(\IE_r)\,.$$
les suites spectrales de Leray associées à $\pi_a:\GaDelta_{?\ell}\Xg^{m}\to\Fa(\Xg)$ qui convergent respectivement vers $\Hc(\GaDelta_{?\ell}\Xg^{m})$ et $\Hpi(\GaDelta_{?\ell}\Xg^{m})$.
On a:
$$
\begin{cases}\noalign{\kern-2pt}
\IE_2(\GaDelta_{?\ell}\Xg^{m})_{\rmc}\sim\Hc^{p}(\Fa(\Xg))\otimes\Hc^{q}(F)\Rightarrow \Hc^{p+q}(\GaDelta_{?\ell}\Xg^{m})\\\noalign{\kern4pt}
\IE_2(\GaDelta_{?\ell}\Xg^{m})_{\cv}\sim\Hr^{p}(\Fa(\Xg))\otimes\Hc^{q}(F)\Rightarrow \Hpi^{p+q}(\GaDelta_{?\ell}\Xg^{m})\\\noalign{\kern-2pt}
\end{cases}
$$
où 
$F$ est une fibre quelconque de $\pi_a$. Les suites spectrales $(\IE_{r})$ sont dégénérées, \idest $d_{r}=0$ pour $r\geq2$.
\end{theo*}

\soustitreni{Exemples et contre-exemples d'espaces $i$-acycliques}
La section~\ref{appendice-exemples} est un appendice destiné à donner
des exemples d'espaces $i$-acycliques, ainsi que des contre-exemples à certaines propriétés les concernant.

\soustitreni{D'autres approches à la stabilité}Lors de l'écriture de cet article nous avons eu connaissance du travail de Dan Petersen \cite{petersen} (2016) où il associe à un espace topologique stratifié
une suite spectrale qui calcule l'homologie de Borel-Moore de l'espace en termes de celle de ses strates fermées. Dans le cas des espaces de configuration $\Fm(\Mg)$, il retrouve la suite spectrale de Getztler \cite{getzler} qui généralise celle de Totaro (\cite{tota}) pour le plongement $\Fm(\Mg)\hook\Mg^{m}$, du cas où $\Mg$ est une variété différentielle au cas où $\Mg$ est générale. Petersen est alors en mesure de montrer, lorsque $\Mg$ une variété algébrique complexe,  que les $\FI$-modules $\set\Fm(\Mg)/_{m}$ sont de type fini. La stabilité de $\set\Fm(\Mg)/_{m}$ en découle mais sans pour autant permettre d'obtenir des indications sur les rangs de monotonie et stabilité concernés.
\endgroup
\let\theotitlestyle\defaulttheotitlestyle

\section{Espaces $i$-acycliques}
\addtoglossary{\vskip-1.5ex}
\glossarytitle{Généralités sur les espaces $i$-acycliques}
\subsection{Généralités sur les espaces topologiques considérés}Dans ce travail on entend par \emph{espace (topologique)} tout espace $\Xg,\Yg,\Zg,\dots$, métrisable,  localement compact et dénombrable à l'infini. De tels espaces sont alors à base dénombrable, séparables et (totalement) paracompacts. Tout fermé et tout ouvert de $\Xg$ est alors également métrisable, localement compact et dénombrable à l'infini. Les \emph{pseudovariétés} (dénombrables à l'infini), en particulier les variétés algébriques complexes, constituent une large famille d'exemples de tels espaces (\footnote{Dans \cite{borel-sem}: A.~Haefliger, ``{\slshape Introduction to piecewise linear intersection homology}'' (p.~1); et A.~Borel, 
``{\slshape Sheaf theoretic intersection cohomology\/}'' 
(V-\Spar 2.1, p.~60).}).

\subsubsectionline{Cohomologies.}On désignera par $k$ un corps de caractéristique arbitraire sauf mention explicite du contraire. On note $\fs k_{\Xg}$ le faisceau constant sur $\Xg$ de fibre $k$, puis $\Mod(\fsk_{\Xg})$ la catégorie des faisceaux de $k$-espaces vectoriels sur $\Xg$. Les espaces de cohomologie respectivement: ordinaire, à support compact et à support dans une partie localement fermée $\Zg\dans\Xg$, notés 
 $\Hr(\Xg;k)$,
 $\Hc(\Xg;k)$, et 
 $H_{\Zg}(\Xg,k)$, 
sont les foncteurs dérivés des foncteurs de sections globales
$$\Gamma(\Xg;\_), \Gammac(\Xg;\_), \Gamma_{\Zg}(\Xg;\_):\Mod(\fsk_{\Xg})\to\Vec(k)\,.$$
On a donc
$$
 \Hr^{i}(\Xg;k):=\IR^{i}\Gamma(\Xg;\fs k_{\Xg})\,,\quad
 \Hc^{i}(\Xg;k):=\IR^{i}\Gammac(\Xg;\fs k_{\Xg})\,,
$$
$$
 \Hr_{\Zg}^{i}(\Xg;k):=\IR^{i}\Gamma_{\Zg}(\Xg;\fs k_{\Xg})\,.
$$

\subsubsectionline{Résolution $\varPhi$-molle canonique.}
La\label{resolution-c-molle} résolution du faisceau constant $\fs k_{\Xg}$ de référence sera celle des faisceaux de \expression{germes de cochaînes d'Alexander-Spanier de $\Xg$}\index{Alexander-Spanier}\index{germes de cochaînes d'Alexander-Spanier} (\footnote{\cite{god}, \Spar2.5, exemple 2.5.2, p.~134, puis \Spar3.7, exemple 3.7.1, p. 157.}), noté $(\AS^{\bullet}(\Xg;k),d_*)$\glossary{$(\AS^{\bullet}(\Xg;k),d_*)$:complexe des faisceaux de germes de cochaînes d'Alexander-Spanier de $\Xg$}. Les faisceaux 
$\AS^{i}(\Xg;k)$ sont $\varPhi$-mous pour toute famille paracompactifiante $\varPhi$ (\loccit), ils sont donc $\Gamma(\Xg,\_)$, $\Gammac(\Xg;\_)$ et $\Gamma_{\Zg}(\Xg,\_)$-acycliques.
Le complexe des cochaînes d'Alexander-Spanier
$$\0\to\fs k_{\Xg}\to\AS^{0}(\Xg;k)\too^{d_0}\AS^{1}(\Xg;k)\too^{d_1}\cdots$$
est une résolution $\varPhi$-molle de $\fs k_{\Xg}$ et l'on a
 $$
 \Hr^{i}(\Xg;k):=h^{i}(\Gamma(\Xg;\AS^{\bullet}(\Xg;k)))\,,\quad
 \Hc^{i}(\Xg;k):=h^{i}(\Gammac(\Xg;\AS^{\bullet}(\Xg;k)))\,,$$
 $$
\Hr_{\Zg}^{i}(\Xg;k):= h^{i}(\Gamma_{\Zg}(\Xg;\AS^{\bullet}(\Xg;k)))\,.
 $$
\comment
$$\begin{cases}
\Hr^{i}(\Xg;k)=h^{i}(\Gamma(\Xg;\AS^{\bullet}(\Xg;k)))\\
\Hc^{i}(\Xg;k)=h^{i}(\Gammac(\Xg;\AS^{\bullet}(\Xg;k)))\\
\Hr_{\Zg}^{i}(\Xg;k)=h^{i}(\Gamma_{\Zg}(\Xg;\AS^{\bullet}(\Xg;k)))
\end{cases}$$
\endcomment
\comment
Lorsque $\Xg$ sera une variété différentiable et que $k$ sera une extension de $\RR$ nous utiliserons 
comme résolution de référence pour $\fs k_{\Xg}$
éventuellement, et sans nécessairement le préciser, le complexe de faisceau de formes différentielles 
$$\def\AS{\fs\Omega}\0\to\fs k_{\Xg}\to\AS^{0}(\Xg;k)\too^{d_0}\AS^{1}(\Xg;k)\too^{d_1}\cdots$$
avec $\fs\Omega^{i}(\Xg;k):=\fs\Omega^{i}_{\Xg}\otimes_{\RR}k$.
 
Une classe de cohomologie à support compact sera donc représentée par une cochaîne ou une forme différentielle à support compact dans $\Xg$, le contexte devant permettre de faire la différence si nécessaire.

\endcomment

\begingroup\halfdisplayskips
\subsubsectionnumber\label{limind-limproj}Notons
$(\KKK(\Xg),\dans)$\glossary{$\KKK(\Xg)$: famille des parties compactes $K\dans\Xg$} la famille des parties compactes $K\dans\Xg$ munie de l'ordre d'inclusion. Le morphisme naturel en homologie singulière
$$\limind{\vrule depth3pt width0pt}_{K\in\K(\Xg)}H_{*}(K;k)\too^{\simeq} H_{*}(\Xg,k)$$ est un isomorphisme et induit (par dualité vectorielle) l'isomorphisme de cohomologies ordinaires
$$\Hr^{*}(\Xg,k)\downstackrel\too\simeq\limproj{\vrule depth3pt width0pt}_{K\in\KKK(\Xg)}\Hr^{*}(K,k)\,,$$
qui fait correspondre à une classe de cohomologie $\omega\in\Hr^{*}(\Xg,k)$ la famille de ses restrictions $\set\omega\rest_{K}\in\Hr^{*}(K,k)/_{K\in\K(\Xg)}$.

De manière duale, nous disposons pour chaque $K\in\K(\Xg)$ d'un morphisme naturel en cohomologie $\Hr_{K}^{*}(\Xg,k)\to\Hc^{*}(\Xg,k)$. La limite inductive de ces morphismes
$$\limind{\vrule depth3pt width0pt}_{K\in\K(\Xg,k)}\Hr_{K}^{*}(\Xg)\to\Hc^{*}(\Xg,k)$$
est alors aussi un isomorphisme.
\endgroup

\subsubsectionline{Suite exacte longue de cohomologies à support compact.}\label{sex-coh-compacte}
Si $\Xg=\Ug\sqcup \Zg$\glossary{$\Xg,\Yg,\Zg,\dots$: pseudovariétés} 
est une partition en parties respectivement ouverte et fermée, dont on note $j:\Ug\to\Xg$ et $i:\Zg\to\Xg$ les inclusions,
on dispose de la suite longue de cohomologie à support compact
$$\smash{\Hc(\Ug)\too^{j_!}\Hc(\Xg)\too^{i^{*}}\Hc(\Zg)\smash{\too}}$$
où $j_!$ est le \expression{prolongement par zéro} et $i^{*}$ est la \expression{restriction}.

\subsubsectionline{Isomorphismes de Künneth.}\label{iso-kunneth}Si $p:\Xg\times\Yg\to\Xg$ et $q:\Xg\times\Yg\to \Yg$ sont les projections canoniques, les morphismes 
$$\boxtimes:\AS^{a}(\Xg;k)\otimes\AS^{b}(\Yg;k)\to
\AS^{a+b}(\Xg\times\Yg;k)$$
induits par le cup-produit de cochaînes 
$\boxtimes:\alpha\otimes\beta\mapsto \alpha\boxtimes\beta:=p^{*}(\alpha)\cup q^{*}(\beta)$
\comment
(\footnote{\emph{Mutatis mutandis} dans le cas des variétés différentiables en remplaçant $\AS$ par le faisceau des formes différentielles $\fs\Omega$ et le cup-produit $\cup$ par le produit extérieur de formes différentielles $\wedge$.})
\endcomment
 induisent des isomorphismes gradués (\footnote{L'isomorphisme de Künneth en cohomologie ordinaire exige une hypothèse de finitude de cohomologie sur l'un des deux espaces.})
$$\Hr(\Xg\times\Yg)\simeq\Hr(\Xg)\otimes\Hr(\Yg)
\ \text{et}\ 
\Hc(\Xg\times\Yg)\simeq\Hc(\Xg)\otimes\Hc(\Yg)\,.$$

Le cup-produit en cohomologie résulte alors du morphisme de restriction à la diagonale, soit des composées
$$\halfdisplayskips\arraycolsep2pt\begin{array}{rccccc}
\cup:&\Hr(\Xg)\otimes \Hr(\Xg)&\too^{\boxtimes}& \Hr(\Xg\times\Xg)&\too^{\delta^{*}_{\Xg}}&\Hr(\Delta_{\Xg})\\[4pt]
\cup:&\Hc(\Xg)\otimes \Hc(\Xg)&\too^{\boxtimes}& \Hc(\Xg\times\Xg)&\too^{\delta^{*}_{\Xg}}&\Hc(\Delta_{\Xg})
\end{array}$$

\subsection{Espaces $i$-acycliques\ et \uacycliques}\label{defis-0}

\subsubsectionline{$i$-acyclicité.}Pour tout espace $\Xg$, on notera \glossary{${\epsilon_{\Xg}:\Hc(\Xg)\to\Hr(\Xg)}$:morphisme induit par l'inclusion $\Omega_{\rm c}(\Xg)\dans\Omega(\Xg)$}
$\epsilon_{\Xg}:\Hc(\Xg)\to\Hr(\Xg)$
le morphisme induit par l'inclusion  $\Gammac(\Xg;\_)\dans\Gamma(\Xg;\_)$.
L'image de $\epsilon_{\Xg}$ dans $\Hr(\Xg)$,  notée $\Hi(\Xg)$, est  \expression{la cohomologie intérieure de $\Xg$}\index{intérieure (cohomologie)}\glossary{${\Hi(\Xg):=\im(\epsilon_{\Xg})}$: cohomologie \expression{intérieure} de $\Xg$}\index{cohomologie intérieure}, on pose donc 
$$\halfdisplayskips\Hi(\Xg):=\im\big(\epsilon_{\Xg}:\Hc(\Xg)\to\Hr(\Xg)\big)\,.$$
On dira que $\Xg$ est \expression{$i$-acyclique}\index{iacyclique@$i$-acyclique} lorsque $\Hi(\Xg)=0\,.$

\subsubsectionline{(Totale) $\cup$-acyclicité.}Le plongement diagonal
$\delta_\Xg:\Delta_{\Xg}\hook\Xg\times\Xg$ induit des morphismes de restriction rendant commutatif le diagramme:\glossary{${\delta_{\Xg}^{*}:\Hc(\Xg\times\Xg)\to\Hc(\Xg)}$:restriction à la diagonale}
$$\preskip2pt\xymatrix@R=10mm{
\Hc(\Xg)\otimes\Hc(\Xg)\ar[r]^(0.54){\boxtimes}_(0.54){\simeq}\ar[d]|{\ttt84\epsilon_\Xg\otimes\epsilon_{\Xg}}&\Hc(\Xg\times\Xg)\ar[r]^(0.55){\delta_{\Xg}^{*}}
\ar[d]|{\ttt76\epsilon_{\Xg\times\Xg}}&\Hc(\Delta_\Xg)\ar[d]|{\ttt77\epsilon_{\Delta_{\Xg}}}\\
\Hr(\Xg)\otimes\Hr(\Xg)\ar[r]^(0.54){\boxtimes}_(0.54){\simeq}
&\Hr(\Xg\times\Xg)\ar[r]^(0.55){\delta_{\Xg}^{*}}&\Hc(\Delta_\Xg)
}
$$

On dira que $\Xg$ est \expression{\uacyclique}\index{\uacyclique} lorsque 
$\delta_{\Xg}^{*}:\Hc(\Xg\times\Xg)\to\Hc(\Delta_{\Xg})$ est nul, autrement dit, lorsque le cup-produit $\cup:\Hc(\Xg)\otimes\Hc(\Xg)\to\Hc(\Xg)$ est nul.
Plus généralement, $\Xg$ sera dit \expression{totalement \uacyclique}\index{totalement-\uacyclique} lorsque le cup-produit
$\cup:\Hc(\Xg)\otimes\Hc(\Ug)\to\Hc(\Ug)$
est nul pour tout ouvert $\Ug\dans\Xg$.

\subsubsection{Premières propriétés d'espaces $i$-acycliques}
La condition d'$i$-acyclicité est une condition de non compacité, aussi, son étude nous place d'emblée en opposition aux espaces compacts.

\noendpoint\begin{prop}\label{prop-acycliques}
\begin{enumerate}
\item\leavevmode\label{prop-acycliques-0}Un espace contractile non compact est $i$-acyclique.
\item\leavevmode\label{prop-acycliques-a}$i$-acyclique$\,\Rightarrow\,$ totalement \uacyclique.
\item\leavevmode\label{prop-acycliques-b}Sur une variété topologique orientable, $i$-acyclique$\,\Leftrightarrow\,$\uacyclique\ {\rm (\footnote{Voir \ref{U-non-I}, p. \pageref{U-non-I}, pour un contre-exemple si la variété n'est pas orientable.})}.

\comment
\item\leavevmode\label{prop-acycliques-bb}Sur une variété topologique orientable dimension $n$, la dualité de Poincaré induit une dualité entre $\Hi^{k}(\Xg)$ et $\Hi^{n-k}(\Xg)$.
\endcomment

\item\leavevmode\label{prop-acycliques-c}Un ouvert $\Ug$ d'un espace $i$-acyclique  $\Xg$, est  $i$-acyclique. Et de même en remplaçant \emph{$i$-acyclique} par \emph{totalement \uacyclique}.
\item\leavevmode\label{prop-acycliques-d}Un produit cartésien $\Xg\times\Yg$ est $i$-acyclique  si et seulement si l'un des facteurs l'est.
\end{enumerate}
\end{prop}
\demo \parskip1ex
(\ref{prop-acycliques-0}) Comme $\Xg$ est contractile, $\Hr^{i}(\Xg)=0$ pour $i>0$ et comme $\Xg$ est connexe non-compacte
$\medmuskip0mu
\Hc^{0}(\Xg)=0$, dans tous les cas 
$\medmuskip0mu\Hi^{i}(\Xg)=\nobreak0$.

(\ref{prop-acycliques-a})
Si $\omega$ et $\varpi$ sont des cocycles à support compact d'un espace $i$-acyclique, on a $\omega=d\alpha$ pour une certaine cochaîne $\alpha$ (à support fermé), mais alors $\omega\cup \varpi=d\alpha\cup\varpi=d(\alpha\cup\varpi)$
où $\alpha\cup\varpi$ est clairement à support compact.

(\ref{prop-acycliques-b}) Supposons que pour $\omega\in\Hc(\Xg)$ on ait $\epsilon_{\Xg}(\omega)\not=0$. Il existe alors, par dualité de Poincaré, une classe $\varpi\in\Hc(\Xg)$ telle que $\int_{\Xg}\omega\cup\varpi=1$, mais alors $\omega\cup\varpi\not=0$ dans $\Hc(\Xg)$ et $\Xg$ n'est pas \uacyclique.

\comment
(\ref{prop-acycliques-bb}) Les classes de $\Hi^{k}(\Xg)$ sont celles de $\Hr^{k}(\Xg)$ qui peuvent être représentées par un cocycle à support compact dans $\Xg$. Maintenant, si $\epsilon_{\Xg}\omega\in\Hi^{k}(\Xg)$ est non nulle, il existe $\varpi\in\Hc^{n-k}(\Xg)$ telle que $\int_{\Xg}\epsilon_{\Xg}(\omega)\cup\varpi=1$. 
L'accouplement $\int_{\Xg}:\Hi^{k}(\Xg)\times\Hi^{n-k}(\Xg)\to k$ est donc parfait.
\endcomment

(\ref{prop-acycliques-c}) Résulte du fait que  $\epsilon_\Ug:\Hc(\Ug)\to H(\Ug)$ est la composée de
$$
\xymatrix@C=5mm@R5mm{
\Hc(\Ug)\ar@{=}[d]\ar[rrr]|{\,\epsilon_\Ug\,}&&&\Hr(\Ug)\ar@{=}[d]\\
\Hc(\Ug)\ar[r]^{\jmath}&
\Hc(\Xg)\ar[r]^{\epsilon_\Xg}&
\Hr(\Xg)\ar[r]^{r}&
\Hr(\Ug)
}$$
où $\jmath$ désigne le prolongement par zéro et $r$ est le morphisme de restriction.

(\ref{prop-acycliques-d}) La condition est suffisante d'après le factorisation de Künneth:
$$\xymatrix@R5mm@C1.2cm{
\Hc(\Xg\times\Yg)\ar@{=}[d]\ar[r]|{\,\epsilon_{\Xg\times\Yg}\,}&\Hr(\Xg\times\Yg)\\
\Hc(\Xg)\otimes\Hc(\Yg)\ar[r]^{\epsilon_\Xg\otimes\epsilon_\Yg}&
\Hr(\Xg)\otimes\Hr(\Yg)\ar[u]_{\kappa}}
\postdisplaypenalty10000$$
où \smash{$\kappa(\alpha\otimes\beta)=p_{\Xg}^{*}(\alpha)\cup p_{\Yg}^{*}(\beta)$}. La nécessité résulte de l'injectivité de $\kappa$.
\enddemobox
\enddemo

\subsection{Caractérisations de la $i$-acyclicité}

\begin{theo}\label{caracterisation}Pour un espace $\Xg$, il y a équivalence entre
\def\varlistskips{\topsep2pt\itemsep2pt\mou\parskip0pt\mou}\begin{enumerate}
\item $\Xg$\label{caracterisation-a} est $i$-acyclique.
\item\leavevmode\label{caracterisation-a'}Pour tout espace $\Yg$ et toute application continue $f:\Yg\to\Xg$, le morphisme image-inverse 
$f^{*}:\Hc(\Xg)\to\Hr(\Yg)$ est nul.
\item\leavevmode\label{caracterisation-b}Pour tout espace  compact $\Kg$ et toute application continue $f:\Kg\to\Xg$, le morphisme image-inverse 
$f^{*}:\Hc(\Xg)\to\Hr(\Kg)$ est nul.
\item\leavevmode\label{caracterisation-clef}\'Etant donnés $\Zg$ et $\Yg$ et une application continue $f:\Zg\to\Xg\times\Yg$, notons $p_2:\Xg\times\Yg\to\Yg$, $(x,m)\mapsto m$, puis $f_2:=p_2\circ f$. Soit $j:V\hook\Yg$ un plongement ouvert, notons $f'_{2}:f^{-1}(V)\to V$ la restriction de $f_{2}$, et considérons le diagramme commutatif suivant
$$\preskip2pt\postskip4pt
\xymatrix@R=5mm{
f_{2}^{-1}(V)\ar[rd]|{f'_2}\ar@{^(->}[r]&\Zg\ar[rd]|{f_2}\ar[r]|(0.4){\,f\,}&\Xg\times\Yg\ar[d]^{p_2}\\
&V\ar@{^{(}->}[r]|(0.45){\,j\,}&\Yg\mrlap{\,.}
}$$
Alors, si $f$ et $f'_2$ sont propres et si $\mathrigid0mu
j_{!}:\Hc(V)\to\Hc(\Yg)$ est surjective, on a
$$\preskip4pt\big(f^{*}:\Hc(\Xg\times\Yg)\to\Hc(\Zg)\big)=0\,.$$
\end{enumerate}
\end{theo}
\demo\displayskips8/10
 (\ref{caracterisation-a}$\Leftrightarrow$\ref{caracterisation-a'}) \'Evident.
(\ref{caracterisation-a}$\Leftrightarrow$\ref{caracterisation-b}) Le morphisme $f^{*}:\Hc(\Xg)\to\Hr(\Kg)$ se factorise à travers $\epsilon_{\Xg}$ puisque l'on a le diagramme commutatif
$$\xymatrix@R=5mm{
\Hc(\Xg)\ar[r]|{\,f^{*}}\ar[d]_{\epsilon_{\Xg}}&\Hc(\Kg)\ar@{=}[d]^{\epsilon_{\Kg}}\cr
\Hr(\Xg)\ar[r]|{\,f^{*}}&\Hr(\Kg)
}$$ 
L'implication \ref{caracterisation-a}$\Rightarrow$\ref{caracterisation-b} en résulte. Réciproquement, soit $(\K(\Xg),\dans)$ l'ensemble des parties compactes de $\Xg$, ordonné par inclusion. La famille des morphismes 
$$
\xymatrix@C=1.2cm@R=5mm{
\Hc(\Xg)\ar[r]^{\epsilon_{\Xg}}&\Hr(\Xg)
\mrlap{\hf{(\_)\srest_{\Kg}}{}{1.3cm}\Hr(\Kg)}
&&
}$$
 avec $\Kg\in\K(\Xg)$ induit alors un morphisme sur la limite projective
$$\xymatrix@C=1.2cm{
\Hc(\Xg)\ar[r]^{\epsilon_{\Xg}}&\Hr(\Xg)
\mrlap{\hf{\limproj(\_)\srest_{\Kg}}{}{2cm}
\smash{\limproj_{\Kg\in\K(\Xg)}\Hr(\Kg)}}
&&
}$$
dont la composée est nulle lorsque (b) est vérifié. L'espace $\Xg$ est donc $i$-acyclique puisque le morphisme $\limproj(\_)\srest_{\Kg}$ est bijectif.

(\ref{caracterisation-a}$\Rightarrow$\ref{caracterisation-clef}) Soient $p_1:\Xg\times\Yg\to\Xg$ la projection canonique et $f_{1}:=p_{1}\circ f$. Pour prouver (\ref{caracterisation-clef}), il suffit, par Künneth, de montrer que pour tous cocycles à support compact $\omega\in\Z_{\rmc}(\Xg)$ et $\varpi\in\Z_{\rmc}(\Yg)$, le cocycle
$$f_{1}^{*}\omega\cup f_{2}^{*}\varpi\in\Z_{\rmc}(\Zg)\,.\eqno(\ast)$$
est la différentielle d'une cochaîne à support compact de $\Zg$.

Or, il existe par hypothèse $\tau\in\Z_{\rmc}(V)$ qui représente $\varpi$ dans $\Hc(\Yg)$. On peut donc remplacer dans $(\ast)$
$f_{2}^{*}\varpi$ par $f_{2}^{*}\tau$, et affirmer que $f_{2}^{*}\tau$ est à support compact dans $f_{2}^{-1}(V)$ donc dans $\Zg$.

Cela étant, comme $\Xg$ est $i$-acyclique, on a $f_{1}^{*}\omega=d\alpha$ pour une certaine cochaîne $\alpha$ de $\Zg$ (à support non nécessairement compact), et alors
$$f_{1}^{*}\omega\cup f_{2}^{*}\tau=d(\alpha\cup f_{2}^{*}\tau)\,,$$ 
où $\alpha\cup f_{2}^{*}\tau$ est une cochaîne à support compact de $\Zg$ puisqu'il en est ainsi de $f_{2}^{*}\tau$. Par conséquent, $f^{*}(\omega\otimes\varpi)=0\in\Hc(\Zg)$.

\nobreak (\ref{caracterisation-clef}$\Rightarrow$\ref{caracterisation-b}) Dans (\ref{caracterisation-clef}), on prend $\Zg$ compact, $\Yg:=\set\pt/$ et $V:=\Yg$. 
\enddemo
\subsubsectionline{Complémentaires des parties finies d'un espace $i$-acyclique.}
Pour\label{complementaire-partie-finie} $a\in\NN$, la notation \og$\Xg\mmoins a$\fg\ sera un raccourci pour \og $\Xg$ privé de $a$ points\fg\glossary{$\Xg\mmoins \ag$: ``complémentaire'' dans $\Xg$ d'une partie finie $\ag$ de cardinal $a$}. On rappelle qu'en dehors du cas où $\Xg$ est une variété topologique connexe de dimension $>1$, le type d'homéomorphie de $\Xg\mmoins a$ dépend des points qu'on enlève. Cependant, lorsque $\Hc^{0}(\Xg)=0$, p.e. si $\Xg$ est $i$-acyclique, la suite
$0\to
\Hc(F)[-1]\to
\Hc(\Xg\mmoins F)\to
\Hc(\Xg)\to
0
$
est exacte pour tout $F\dans\Xg$ fini, et la famille
\smashtop{$\set\Bettic^{i}(\Xg\mmoins F)/_{i\in\NN}$}, et donc le polynôme de Poincaré $\Pc(\Xg\mmoins F)$ pour la cohomologie à support compact, seul dépendent de $\# F$.

\subsubsectionline{Trivialité de la monodromie.}Lorsque\label{monodromie} $\Xg$ est un espace $i$-acyclique \emph{localement connexe}, la remarque précédente peut être raffinée pour montrer l'existence d'une action par monodromie de $\Pi_{1}(\Fa(\Xg),\cl x)$ sur $\Hc(\pia^{-1}(\cl x))$ où $\pia:\Fg_{1{+}a}(\Xg)\to\Fa(\Xg)$ est la projection sur les $a$ dernières coordonnées, et ce, \emph{même lorsque cette projection n'est pas localement triviale}. Pour le voir, considérons dans un premier temps le cas où $a=1$.

Notons $\KKK(\Cg)$ la famille des parties compactes et connexes contenues dans une composante connexe $\Cg$ de $\Xg$.
\smash{\'Etant} donnés $z\in Z\in\KKK(\Cg)$, on dispose du morphisme naturel de suites exactes courtes
$$
\xymatrixc{@R=6mm}{
0\ar[r]&
\Hc(Z)[-1]\ar[r]\ar[d]_{\rho(Z)_z}&
\Hc(\Xg\mmoins Z)\ar[r]\ar[d]_{\rho(Z)}&
\Hc(\Xg)\ar[r]\aregal[d]&
0\\
0\ar[r]&
\Hc(z)[-1]\ar[r]&
\Hc(\Xg\mmoins z)\ar[r]&
\Hc(\Xg)\ar[r]&
0\\
}\eqno(\diamond)$$
où le prolongement par zéro $\rho(Z)$ est bijectif en degrés cohomologiques $i\leq1$ et induit pour $i\geq2$ l'isomorphisme canonique:
$$\halfdisplayskips
\tilde\rho(Z)_{z}:{\Hc^{i}(\Xg\mmoins Z)\over\Hc^{i-1}(Z)}\too^{\cong}\Hc^{i}(\Xg\mmoins z)\,,\quad\forall z\in Z\,.$$
Cette construction est naturelle par rapport à l'inclusion de compacts  et font de $\KKK(\Cg)$ un système inductif filtrant dont la réunion recouvre $\Xg$ tout entier. On en déduit un isomorphisme \emph{canonique}
$$\halfdisplayskips
\limproj_{Z\in\KKK(\Cg)}{\Hc^{i}(\Xg\mmoins Z)\over\Hc^{i-1}(Z)}\too^{\cong}\Hc^{i}(\Xg\mmoins x)\,,\quad\forall x\in\Cg\,,$$
dont on conclut que l'isomorphisme
$$\phi(Z)_{y,x}:=(\tilde\rho(Z)_{y}\circ \tilde\rho(Z)_{x}^{-1}):
\Hc(\Xg\mmoins x)\to\Hc(\Xg\mmoins y)$$
est indépendant du compact connexe $K\cont\set x,y/\dans\Cg$. 
En particulier, une application continue $\gamma:[0,1]\to\Xg$ détermine l'isomorphisme
$$\phi_{\gamma}:=\phi(\gamma([0,1]))_{\gamma(1),\gamma(0)}=\Hc(\Xg\mmoins\gamma(0))\to
\Hc(\Xg\mmoins\gamma(1))\eqno(\diamond\diamond)$$
qui seul dépend de extrémités de $\gamma$. C'est l'action par\expression{monodromie le long du chemin $\gamma$}, et elle est clairement triviale lorsque $\gamma(0)=\gamma(1)$.

\noindent Ces remarques se généralisent à tout $a\in\NN$.\killline\killline

\begin{prop}Pour\label{prop-monodromie} $a\in\NN$\label{independance-Hc-X-a}, soit $\pia:\Fg_{1{+}a}(\Xg)\to\Fa(\Xg)$ la projection sur les dernières $a$ coordonnées. Lorsque $\Xg$ est $i$-acyclique, il existe pour tout $\cl x\in\Fa(\Xg)$ une action par monodromie de $\Pi_{1}(\Fa(\Xg),\cl x)$ sur $\Hc(\pia^{-1}(\cl x))$. Cette action est triviale. 
\end{prop}
\def\Demonstration{Indication}\demo\displayskips8/10
Commençons par justifier l'existence de l'action par monodromie. Soit $\cl\gamma:=(\gamma_1,\ldots,\gamma_a):[0,1]\to\Fa(\Xg)$ une application continue. Par compacité, il existe une famille finie $\relspace1
\TTT:=\set 0=t_0<t_1,\cdots<t_r=1/$ telle que 
pour chaque $0\leq s<r$ la famille 
$\set Z(s)_{k}:=\gamma_{k}([t_s,t_{s{+}1}])/_{1\leq k\leq a}$, 
est constituée de parties compactes connexes deux à deux disjointes. Notons $Z(s):=\bigcup_{k}Z(s)_{k}$ et reprenons le diagramme $(\diamond)$ pour $\cl z:=\set z_k\in Z(s)_{k}/_{k}$, on a
$$
\xymatrixc{@R=6mm}{
0\ar[r]&
\Hc(Z(s))[-1]\ar[r]\ar[d]_{\rho(s)_{\cl z}}&
\Hc(\Xg\mmoins Z(s))\ar[r]\ar[d]_{\rho(s)}&
\Hc(\Xg)\ar[r]\aregal[d]&
0\\
0\ar[r]&
\bigoplus_{k}\Hc(z_k)[-1]\ar[r]&
\Hc(\Xg\mmoins \cl z)\ar[r]&
\Hc(\Xg)\ar[r]&
0\\
}\eqno(\diamond_{s})$$
où le prolongement par zéro $\rho(s)$ est un isomorphisme en degrés $\leq 1$ et induit pour $i\geq2$
l'isomorphisme canonique:
$$\tilde\rho(s)_{\cl z}:{\Hc^{i}(\Xg\mmoins Z(s))\over\Hc^{i-1}(Z(s))}\too^{\cong}\Hc^{i}(\Xg\mmoins \cl z)\,,\quad\forall  \cl z\in Z(s)_{1}\times\cdots\times Z(s)_{a}\,,$$
d'où l'isomorphisme
$$\phi_{s}:=\rho(s)_{\cl\gamma(t_{s+1})}\circ\rho(s)_{\cl\gamma(t_{s})}^{-1}:\Hc(\Xg\mmoins\cl\gamma(t_s))\cong\Hc(\Xg\mmoins\cl\gamma(t_{s{+}1}))\,.$$
La composée des $\phi_{s}$ est \expression{l'action par monodromie le long du chemin $\cl\gamma$}\glossary{$\hbox{\mathrigid3mu $\phi_{\cl\gamma}:\Hc(\Xg{\mmoins}\cl\gamma(0))\to\Hc(\Xg\mmoins\cl\gamma(1))$}$:action par monodromie le long du chemin $\cl\gamma$}
$$\phi_{\gamma}:\Hc(\Xg\mmoins\cl\gamma(0))\to\Hc(\Xg\mmoins\cl\gamma(1))\,,
\postdisplaypenalty10000$$
dont on vérifie aisément qu'il ne dépend que de $\cl\gamma$ et non pas de la famille $\TTT$.

Montrons maintenant, par induction sur $a$, que $\phi_{\cl\gamma}$ seul dépend des extrémités de $\cl\gamma$. Lorsque $a=1$, c'est l'indépendance $(\diamond\diamond)$. Dans le cas où $a\geq2$, notons $\cl\gamma':=(\gamma_1,\ldots,\gamma_{a-1})$ le chemin dans $\Fg_{a-1}(\Xg)$. La naturalité de la construction de l'action par monodromie donne lieu à un morphisme des suites courtes, extraites de suites longues de Mayer-Vietoris, 
$$\def\dec{40pt}
\xymatrix@C=4mm@R=5mm
{
0\ar[r]&\Hc(\Xg\mmoins\cl\gamma(0))\ar[d]_{\phi_{\cl\gamma}}\ar[r]&
\Hc(\Xg\mmoins\cl\gamma'(0))\oplus\Hc(\Xg\mmoins\gamma_a(0))\ar@<-\dec>[d]_{\phi_{\cl\gamma'}}\ar@<\dec>[d]^{\phi_{\gamma_a}}\ar[r]&
\Hc(\Xg)\aregal[d]\ar[r]&0\\
0\ar[r]&\Hc(\Xg\mmoins\cl\gamma(1))\ar[r]&
\Hc(\Xg\mmoins\cl\gamma'(1))\oplus\Hc(\Xg\mmoins\gamma_a(1))\ar[r]&
\Hc(\Xg)\ar[r]&0
}
$$
où les lignes sont exactes puisque $\Xg$ est $i$-acyclique. 
Comme $\phi_{\cl\gamma'}$ et $\phi_{\gamma_a}$ seul dépendent de leurs extrémités par hypothèse de récurrence, on conclut qu'il en est de même pour $\phi_{\cl\gamma}$.
\enddemo

\def\Remarque{Commentaire}\begin{rema}Nous verrons lors de l'étude des suites spectrales de Leray associées à l'application $\pia:\Fba(\Xg)\to\Fa(\Xg)$,  que l'action par monodromie de $\Pi_1(\Fa(\Xg),\cl x)$ sur $\Hc(\pia^{-1}(\cl x))$ existe et 
est triviale\label{comm-independance-Hc-X-a} 
aussi pour tout $b>1$  (\cf\ref{ameliore-independance-Hc-X-a}).
\end{rema}
\section{Espaces de configuration généralisés}

\glossarytitle{Espaces de configuration généralisés}
\subsectionline{Définitions et notations.}
\noindent Nous\label{notas} rappelons quelques notations habituelles et introduisons d'autres nouvelles.
\begin{list}{\theenumi)}{\def\theenumi{N-\arabic{enumi}}\usecounter{enumi}\labelsep1ex\labelwidth3mm}
\item\leavevmode\label{nota-Delta-}Pour tout $\ell\in\NN$ et pour tout sous-ensemble $\Zg\dans\Xg^{m}$, \expression{l'ensemble de configuration généralisé dans $\Zg$ pour le cardinal $\leq\ell$ (resp. $=\ell$)}, est\glossary{$\Delta_{?\ell}\Xg^{m}$, $\Delta_{?\ell}^{m}$:espace de configurations (ordonnées) généralisé}
$$\begin{casesalign}
\Delta_{\leq \ell}\Zg&:=&\set(z_1,\ldots,z_{m})\in\Zg\mid\card{\set z_1,\ldots,z_{m}/}\leq \ell/\,,\\
\Delta_{\ell}\Zg&:=&\set(z_1,\ldots,z_{m})\in\Zg\mid\card{\set z_1,\ldots,z_{m}/}=\ell/\,.\\
\end{casesalign}$$

\noindent$\triangleright$ La notation $\Delta_{?\ell}$ désignera l'une quelconque des notations $\Delta_{\leq\ell}$ ou $\Delta_{\ell}$.

\noindent$\triangleright$ Si $m>0$, on a $\Delta_{?0}\Xg^{m}=\emptyset$. Si $\nobinspace
 0=\ell=m$, on a $\Delta_{?0}\Xg^{0}=\set\pt/$.
 \comment et si $m\in\NN$ et $\ell<0$ on a $\Delta_{?\ell}\Xg^{m}=\emptyset$.\endcomment

\item\leavevmode\label{nota-Fg}L'\expression{espace des configurations (ordonnées) de $n$ éléments de $\Xg$}\glossary{$\Fm (\Xg)$:espace de configurations (ordonnées) classique, noté aussi $\Delta_{m}\Xg^{m}$}, traditionnellement noté $\Fm (\Xg)$ est:
$$
\Fm (\Xg):=\Delta_{m}\Xg^{m}\,,\text{ si $m>0$}\,,\quad
\Fg_0(\Xg)=\set\pt/\,.
$$
\item\leavevmode\label{nota-FG}\label{nota-pi-a} Pour $\Zg\dans\Xg^{m}$, et $0<a\in\NN$, on note \glossary{$\Ga \Zg$: intersection $\Zg\cap(\Xg^{m-a}\stimes\Fa(\Xg))$}
$$\mathalign{
\Ga \Zg&:=&\Zg\cap(\Xg^{m-a}{\times}\Fa(\Xg))\hfill\\
&=&\set(z_1,\ldots,z_{m})\in\Zg\mid\card{\set z_{m-a+1},\ldots,z_{m}/}=i/\,.
}
$$
\item\leavevmode\label{pia}Pour $0<a\leq m\in\NN$, on note\glossary{$\pia=\Xg^{m}\to\Xg^{a}$:projection sur les $a$ dernières coordonnées}\glossary{$p_a=\Xg^{m}\to\Xg^{a}$:projection sur les $a$ premières coordonnées}
$$\pia=\Xg^{m}\to\Xg^{a}\,,\quad (\hbox{resp. $p_a=\Xg^{m}\to\Xg^{a}$})$$
la projection sur les $a$ dernières coordonnées $\pia (\cl x):=(x_{m-a+1},\ldots, x_{m})$ (resp. les $a$ premières coordonnées $\pia (\cl x):=(x_{1},\ldots, x_{a})$)
. Les restrictions de $\pia$ seront notées par abus de la même manière, par exemple dans l'écriture, pour tout $\Zg\dans\Xg^{m}$,
$$\pia :\Ga\Zg\to\Fa(\Xg)\,.$$
\end{list}

\def\vartitle{Convention}
\begin{var}Lorsque la référence à $\Xg$ sera superflue elle sera parfois omise et les notations  $\Delta_{?\ell}\Xg^{m}$ et $\Fm(\Xg)$ seront abrégées en $\Delta_{?\ell}^{m}$ et $\Fm$.
\end{var}

\subsubsectionline{Topologie et caractéristique d'Euler.}Lorsque\label{topologie-espace-configurations}\label{suite-longue-de-base} $\Xg$ est un espace topologique, $\Xg^{m}$ est muni de la topologie produit et tout $\Zg\dans\Xg^{m}$ de la topologie induite. 
Pour $\ell,m\in\NN$, et $\Zg\dans\Xg^{m}$, l'espace $\Delta_{\leq\ell}\Zg$ admet la décomposition en parties respectivement ouverte et fermée
$\Delta_{\leq\ell}\Zg=\Delta_{\ell}\Zg\ \sqcup\ \Delta_{\leq\ellmo}\Zg\,,$
d'où la suite exacte longue de cohomologie à support compact que l'on rencontrera fréquemment dans ce travail
$$
\cdots\to\Hc^{*}(\Delta_{\ell}\Zg)\to
\Hc^{*}(\Delta_{\leq\ell}\Zg)\to
\Hc^{*}(\Delta_{\leq\ellmo}\Zg)\to\cdots\,.
$$
Dans le cas où $\Zg=\Xg\times\Fmmo(\Xg)$ et $\ell=m$, on obtient 
$$
\cdots\to\Hc^{*}(\Fm(\Xg))\to
\Hc^{*}(\Xg\times\Fmmo(\Xg))\to
\Hc^{*}(\Delta_{m-1}(\Xg\times\Fmmo(\Xg)))\to\cdots\,.
$$
où $\Delta_{m-1}(\Xg\times\Fmmo(\Xg))$ est la réunion disjointe de $\mmo$ copies de $\Fmmo(\Xg)$, ce qui permet une récurrence pour la détermination de la caractéristique d'Euler de $\Fm(\Xg)$. La proposition suivante en découle aussitôt.

\begin{prop*}[(\ref{pol-poincare})]Supposons que $\Hc(\Xg)<+\infty$ et notons $\chic(\Fm(\Xg))$ la caractéristique d'Euler de la cohomologie à support compact de $\Fm(\Xg)$. On a
$$\chic(\Fm (\Xg))=\prodnl_{i=0}^{\mmo}\chi(\Xg)- i.$$
En particulier, on a la série génératrice 
$$\sumnl_{m\geq0}\chic(\Fm(\Xg))\;{t^{m}\over m!}=(1+t)^{\chic(\Xg)}\,.$$
\end{prop*}

\begin{rema}On trouve cette série génératrice déjà dans Félix-Thomas (\cite{FT}, 2000) pour la cohomologie ordinaire et lorsque $\Xg$ est une variété topologique orientable de dimension paire (cas auquel $\chic(\Xg)=\chi(\Xg)$).
\end{rema}

\subsection{Fibrations des espaces de configuration généralisés}
La\label{fibrations} proposition suivante est une extension assez légère du théorème de trivialité locale de la projection $\pi_a:\Fba(\Xg)\to\Fa(\Xg)$ de Fadell et Neuwirth (\cite{FN}), elle nous a servi dans l'élaboration de la section \ref{localisation} où l'hypothèse de lissité sur $\Xg$ s'est averée à postériori superflue.

\begin{prop}Si\label{fibrations-pi} $\Xg$ une variété topologique et $a\leq\ell\leq m$,
l'application
$$\pia :\GaDelta_{?\ell}\Xg^{m}\to\Fa(\Xg)$$
est une fibration localement triviale. 
\comment 
Si de plus $\Xg$ est connexe de dimension $>1$, les fibres ont le type d'homéomorphie de l'espace $\Xg$ privé de $a$ points, ce que l'on notera  $\Xg\mmoins a$.\endcomment
\end{prop}
\demo 
 Soit $\Gg$ le groupe des homéomorphismes $\phi:\Xg\to\Xg$ dont le support $|\phi|:=\cl{\set x\in\Xg\mid\phi(x)\not=x/}$ est compact. Soit $d_\Xg:\Xg\times\Xg\to\RR$ une distance. Munissons $\Gg$ de la distance $d_{\Gg}(\phi,\phi'):=\sup_{x\in\Xg} d_{\Xg}(\phi(x),\phi'(x))$. Le groupe $\Gg$, muni de la topologie associée à $d_{\Gg}$, est un groupe topologique et pour chaque $x\in\Xg$, l'application d'évaluation 
$$\mathop{\rm ev}\nolimits_{x}:\Gg\to\Xg\,,\quad \phi\mapsto\phi(x)\postskip0pt$$
est continue.

Pour $\cl x\in\Fa$, notons $\IB(\cl x,\epsilon)$
le produit des boules ouvertes $\prod_{i=1}^{a}\IB(x_i,\epsilon)$.
Pour $\epsilon>0$ assez petit, on a $\IB(\cl x,\epsilon)\dans\Fa$ et l'application $\prod \mathop{\rm ev}\nolimits_{x_i}:\Gg\to\Fa$ admet des sections locales continues $\sigma:\IB(\cl x,\epsilon)\to\Gg$, soit
$$\preskip0em\mathalign{
\sigma(\cl y)(\cl x)=\cl y\,,\qquad&
\relax{\xymatrix@C=1.5cm{
\IB(\cl x,\epsilon)\ar@/_18pt/[rr]|(.54){\ \dans\ }
\ar[r]^(0.6){\sigma}&\Gg\ar[r]^{\prod \mathop{\rm ev}\nolimits_{x_i}}&\Fa
}}}
$$
L'application
$$\mathalign{
\IB(\cl x,\epsilon)\times \pia^{-1}(\cl x)&\too&\pia^{-1}(\IB(\cl x,\epsilon))&\quad\dans\quad&\GaDelta_{\leq\ell}\Xg^{m}\,,\\
\cl y\times \cl w&\longmapsto&\sigma(\cl y)(\cl w)
}
$$
avec $\cl w=(\cl x,x_{a+1},\ldots,x_{m})$ et
$\sigma(\cl y)(\cl w)=(\cl y,\sigma(\cl y)(x_{a+1}),\ldots,\sigma(\cl y)(x_{m}))$
 est alors un homéomorphisme et une trivialisation de $\pia$ au-dessus de $\IB(\cl x,\epsilon)$. 
 
De manière entièrement analogue, l'application
$$\mathalign{
\IB(\cl x,\epsilon)\times \pia ^{-1}(\cl x)&\too&\pia ^{-1}(\IB(\cl x,\epsilon))&\quad\dans\quad&\Fba(\Xg)\,,\\
\cl y\times \cl w&\longmapsto&\sigma(\cl y)(\cl w)
}
$$
est une trivialisation locale de $\pia$. La fibre $\pia^{-1}(\cl x)$ est canoniquement homéomorphe au sous-espace $\Fg_{b}(\Xg\mmoins \cl x)\dans\Xg^{b}$. 
\enddemo

\begin{rema}[à propos des fibres de $\pia$]Si\label{fibres-pia} $\Xg$ est une variété topologique connexe de dimension $>1$, il existe des homéomorphismes $\phi:\Xg\to\Xg$ tels que $\phi(\cl x)=\cl y$, et les fibres de $\pia:\GaDelta_{\leq\ell}\Xg^{m}\to\Fa$ sont deux à deux homéomorphes (\cf\ref{complementaire-partie-finie}).
Par contre, si $\Xg$ n'est pas connexe, cette propriété peut être en défaut. Par exemple, si $\Xg=\Ug\sqcup\Vg$ est une réunion disjointe d'ouverts non vides, les fibres de $\pi_2:\Fg_{3}(\Xg)\to\Fg_2(\Xg)$ en $(x,y)\in\Fg_2(\Xg)$ sont de la forme:
$$\def\hbox#1{}
\Ug\minus\set x,y/\sqcup\Vg\hbox{\ si $x,y\in\Ug$}\,,\quad
\Ug\sqcup\Vg\minus\set x,y/\hbox{\ si $x,y\in\Vg$}\,,\quad
\Ug\minus\set x/\sqcup\Vg\minus\set y/\hbox{\ si $x\in\Ug$ et $y\in\Vg$}
$$
suivant que $(x,y\in\Ug)$, que $(x,y\in\Vg)$, ou que ($x\in\Ug$ et $y\in\Vg$). 

Il est facile de trouver des exemples pour $\Ug$ et $\Vg$ tels que, non seulement il n'y a pas d'homéomorphisme entre ces trois types de fibres, mais il n'y a pas, non plus, d'isomorphisme entre leurs cohomologies. 
Lorsque $\Xg$ est $i$-acyclique de type fini (\footnote{\label{tf}On dira que $\Xg$ est de type fini lorsque $\Hc(\Xg)$ et $\Hr(\Xg)$ sont de dimension finie.})\glossary{\expression{$\Xg$ est de type fini}: lorsque $\Hc(\Xg)$ et $\Hr(\Xg)$ sont de dimension finie}, la situation change puisque les nombres de Betti des fibres sont constants (\ref{complementaire-partie-finie}). Par exemple, si $\Ug=\RR^{2}$  et $\Vg=\RR\times\cSS^{1}$, on a 
$\Hr_{*}(\pia^{-1}(u,v))=k(0)^{2}\oplus k(1)^{3}$
quel que soit $(u,v)$. 

Dans ce cas, les algèbres de cohomologie à support compact $\Hc(\pia^{-1}(u,v))$ sont isomorphes, mais pas les algèbres de cohomologie ordinaire qui valent
$$\preskip0.75em\begin{cases}
\displaystyle\Hr(\Ug\mmoins\set u,v/)\oplus\Hr(\Vg)=\smash{k[X,Y]\over (X,Y)^{2}}\oplus \smash{k[Z]\over (Z)^2}\\\noalign{\kern15pt}
\displaystyle\Hr(\Ug)\oplus\Hr(\Vg\mmoins\set u,v/)=k\oplus \smash{k[X,Y,Z]\over (X,Y,Z)^{2}}
\end{cases}\postskip1em$$
et où l'on remarque que  dans la première $\Annul_{\Hr^{0}}(\Hr^{1})=0$, tandis que dans la seconde $\Annul_{\Hr^{0}}(\Hr^{1})=k\cdot(1,0)$.

\smallskip
Tout ceci indique que même dans le cas où $\Xg$ est une variété topologique $i$-acyclique, les monodromies des systèmes locaux :
$$\left\{\vcenter to27pt{}\right.\mathalign{
\cHpi^{i}(\GaDelta_{?\ell}\Xg^{m})&:=&\IR^{i}\pi_{a!}\,\fs k_{\GaDelta_{?\ell}\Xg^{m}}\hfill
\\\noalign{\kern0pt}
\smashbot{\cHr^{i}(\GaDelta_{?\ell}\Xg^{m})}&:=&\smash{\IR^{i}\pi_{a*}\,\fs k_{\GaDelta_{?\ell}\Xg^{m}}}
}$$
de fibres respectives $\Hc(\pia^{-1}\cl x)$ et $\Hr(\pia^{-1}\cl x)$ (en dualité), ont peu de chances d'être triviales sur $\Fa(\Xg)$ tout entier. On verra dans la section \ref{constance-de-faisceau}, consacrée à leur étude, que ces faisceaux sont cependant constants sur toute composante connexe de $\Fa(\Xg)$ (\cf comm.~\ref{comm-independance-Hc-X-a} et th.~\ref{autre-scindage}-(\ref{autre-scindage-c})).
\end{rema}

\subsection{Sous-espaces $\Fg_{\pgoth}(\Xg)$ et décomposition de $\Delta_{\ell}\Xg^{m}$}
\subsubsectionline{Partitions d'un ensemble.}Pour\label{nota-partitions} tout ensemble $E$ et tout $\ell\in\NN$, on note $\Pgoth_{\ell}(E)$\glossary{${\Pgoth_{\ell}(E)}$:ensemble des partitions de l'ensemble $E$ en $\ell$ parties non vides} l'ensemble des partitions de $E$ en $\ell$ parties non vides. On pose ensuite $\Pgoth(E)=\bigcup_{\ell\in\NN}\Pgoth_{\ell}(E)$. 

Une partition $\pgoth\in\Pgoth$ définit une relation d'équivalence `$\relp$'\glossary{$(i\relp j)$:Si $\pgoth\in\Pgoth(E)$, on pose $(i\relp j)\Leftrightarrow_{\rm def}(\exists\, I\in\pgoth)$ t.q. $\set i,j/\dans I$} sur $E$ par 
$$(x\relp y)\Leftrightarrow(\exists I\in\pgoth)(\set x,y/\dans I)\,.$$
Pour tout $m\in\NN$, on note $\Pgoth_{\ell}(m):=\Pgoth_{\ell}(\iii[1,m])$\glossary{$\iii[a,b]$ (resp. $\iij[a,b]$):entiers naturels $m$ vérifiant $a\leq m\leq b$ (resp. $a\leq m< b$)}.

\subsubsectionline{Les sous-espaces $\Fg_{\pgoth}(\Xg)$.}\label{nota-pgoth}Si $\pgoth\in\Pgoth_{\ell}(m)$, on notera $\Fg_{\pgoth}(\Xg)$ l'ensemble des $m$-uplets $(x_1,\ldots,x_{m})$ tels que $(x_i=x_j)\Leftrightarrow(i\relp j)$.
\comment On a $\Fg_{\pgoth}(\Xg)\sim\Fg_{\ell}(\Xg)$.\endcomment
\glossary{$\Fg_{\pgoth}(\Xg)$: Si $\pgoth\in\Pgoth(m)$, on a \smashbot{$\cl x\in\Fg_{\pgoth}(\Xg)\Leftrightarrow_{\rm def}((x_i=x_j)
\Leftrightarrow(i\relp j))$}}

\begin{prop}
Pour\label{connexes}  $0\leq\ell\leq m\in\NN$, on a la décomposition ouverte:
$$
\Delta_{\ell}\Xg^{m}=\coprodnl_{\pgoth\in\Pgoth_{\ell}(m)}\Fg_{\pgoth}(\Xg)\eqno(\ast)$$
{\rm(\cf notation \ref{notas}-(\ref{nota-Fg}))} où 
$\Fg_{\pgoth}(\Xg)\simeq\Fg_{\ell}(\Xg)$.
En particulier,
$$\Hc(\Delta_{\ell}\Xg^{m})=\Hc(\Fg_{\ell}(\Xg))^{|\Pgoth_{\ell}(m)|}\,,$$
où le nombre $|\Pgoth_{\ell}(m)|$ est le nombre de Stirling de seconde espèce {\rm(\cf\ref{Stirling})}:
$$\big|\Pgoth_{\ell}(m)\big|=\Parties m\ell:={1\over \ell! }\,\sumnl_{j=0}^{\ell}\,(-1)^{\ell-j}\binome{\ell}{j}j^{m}\,.\postskip0pt
$$
De plus,
\begin{enumerate}
\item\leavevmode\label{connexes-b}Si $\set\Cg_\alpha/_{\alpha\in\Agoth}$ est la famille des composantes connexes de $\Xg$. On a la décomposition en parties ouvertes
$$
\Fg_{\ell}(\Xg)=\!\!\coprod_{\epsilon:\iii[1,\ell]\to\Agoth}\!\!
\Fg(\Cg_{\epsilon(1)}\times\cdots\times\Cg_{\epsilon(\ell)})=\!\!
\coprod_{\epsilon:\iii[1,\ell]\to\Agoth}\
\prod_{\alpha\in\im(\epsilon)}\!\!
\Fg_{|\epsilon^{-1}(\alpha)|}(\Cg_{\alpha})\,.
$$
\item\leavevmode\label{connexes-a}Si $\Xg$ est une variété topologique connexe de dimension $>1$, l'espace $\Fg_{\ell}(\Xg)$ est connexe et la décomposition $(\ast)$ est la décomposition de $\Delta_{\ell}\Xg^{m}$ en composantes connexes.
\end{enumerate}
\end{prop}

\subsection{Dimension et finitude cohomologique de $\Delta_{?\ell}\Xg^{m}$}\label{finitude-Hc}

\begin{defi}Soit $\Xg$\label{dimch} un espace topologique.
Notons $\D_k(\Xg)$ l'ensemble des $d\in\NN$ tels que pour tout ouvert $U\dans\Xg$, on ait $\Hc^{i}(U;k)=0$,  $\forall i>d$.
Si $\D_k(\Xg)\not=\emptyset$, on dit que $\Xg$ est de
\relax{dimension cohomologique finie (sur $k$) de dimension
$d_{\Xg}:=\dimch(\Xg)=\inf\D_k(\Xg)$}\glossary{$\dimch(\Xg)$:dimension cohomologique de $\Xg$}. Autrement, $\Xg$ est dit de dimension cohomologique infinie et l'on pose $d_{\Xg}=+\infty$ 
 (\cf\cite{borel-sem} V, def.~1.15, p.~55).
\end{defi}

\begin{prop}Soient\label{prop-finitude-Hc} $0\leq\ell\leq m\in\NN$. Si $\Xg$ est respectivement localement compact, localement connexe (par arcs), localement cohomologiquement trivial, 
il en est de même de \smashbot{$\Delta_{?\ell}\Xg^{m}$}. De plus,
\def\varlistskips{\topsep0pt\itemsep4pt}
\begin{enumerate}
\item Le\label{prop-finitude-Hc-0}\label{faisceau-Omega} sous-complexe $(\fs\Omega^{\bullet}_{\Xg},d_{\bullet}):=\tau_{\leq d_{\Xg}}(\AS^{\bullet}(\Xg;k),d_*)$\glossary{${(\fs\Omega^{\bullet}_{\Xg},d_{\bullet}):=\tau_{\leq d_{\Xg}}(\AS^{\bullet}(\Xg;k),d_*)}$:résolution $c$-molle d'un espace $\Xg$ de dimension $d_{\Xg}$} du complexe des germes de cochaînes d'Alexander-Spanier, à savoir 
$$\fs\Omega^{i<d_{\Xg}}_{\Xg}=
\AS^{i}(\Xg;k)\,,\qquad
\fs\Omega^{d_{\Xg}}_{\Xg}:=\fs\ker(d_{d_{\Xg}})\,,\qquad
\fs\Omega^{i>d_{\Xg}}_{\Xg}=0\,,
$$\comment
$$\def\quad{\ }
\fs\Omega^{i}_{\Xg}=\begin{cases}
\AS^{i}(\Xg;k)\,,&\text{si $i<d_{\Xg}$,}\\
\fs\ker(d_{d_{\Xg}})\,,&\text{si $i=d_{\Xg}$,}\\
\hfil0\,,&\text{si $i>d_{\Xg}$,}
\end{cases}$$\endcomment
est une résolution $c$-molle de $\fs k_{\Xg}$. Pour toute partie $\Yg\dans\Xg$  localement fermée, les complexes
$\Gamma(\Yg;(\fs\Omega^{\bullet}_{\Xg},d_{\bullet}))$
et $\Gammac(\Yg;(\fs\Omega^{\bullet}_{\Xg},d_{\bullet}))$
 calculent respectivement $\Hr(\Yg;k)$ et $\Hc(\Yg;k)$.

\item On\label{prop-finitude-Hc-a} a $\dimch(\Delta_{?\ell}\Xg^{m})=\ell\, d_{\Xg}\,.$

\item Si\label{prop-finitude-Hc-b} $\dim_{k}\Hc(\Xg;k)<+\infty$, on a
$
\dim_{k}\Hc(\Delta_{?\ell}\Xg^{m};k)<+\infty\,.$
\end{enumerate}
\end{prop}
\def\Demonstration{Indications}
\demo 
\halfdisplayskips
Les assertions préliminaires à propos de la structure locale sont immédiates pour $\Xg^{m}$, puis pour $\Fm(\Xg)$ car ouvert dans $\Xg^{m}$, et alors pour $\Delta_{\ell}\Xg^{m}$ d'après
la décomposition $(\ast)$ de \ref{connexes}. Le cas de  $\Delta_{\leq\ell}\Xg^{m}$ résulte aussi de la décomposition $(\ast)$ qui montre qu'il est réunion finie des fermés $\cl{\Fg_\pgoth(\Xg)}$, tous homéomorphes à $\Xg^{\ell}$.

(\ref{prop-finitude-Hc-0}) Voir \cite{borel-sem} V, \Spar1.14 (p.~55) ou \cite{KS} exercice II.9 (p.~132).

(\ref{prop-finitude-Hc-a}) L'égalité classique $\dimch(\Xg^{m})=md_{\Xg}$ donne  $\dimch(\Fm(\Xg))\leq m d_{\Xg}$, d'où $\dimch(\Delta_{\ell}\Xg^{m})\leq \ell d_\Xg$.
Ensuite, pour tout ouvert $U\dans\Delta_{\leq \ell}\Xg^{m}$, la suite longue de cohomologie à support compact
$$
\dots\to\Hc^{i-1}(U\cap\Delta_{\leq \ellmo}\Xg^{m})\to
\Hc^{i}(U\cap\Delta_{\ell}\Xg^{m})\to
\Hc^{i}(U)\to\cdots
\eqno(\ast)$$
et un raisonnement inductif sur $\ell\geq0$ permettent de voir
que l'on a
$$\dimch(\Delta_{\leq \ell}\Xg^{m})\leq\ell d_{\Xg}\,.\eqno(\diamond)$$

L'assertion (\ref{prop-finitude-Hc-a}) en découle aussitôt lorsque $d_\Xg=0$. Lorsque $d_\Xg\geq1$, on fait une récurrence sur $m$. Les cas où $m\leq1$ sont triviaux. Lorsque $m\geq2$, on fait une récurrence sur $\ell$. Les cas où $\ell\leq1$ sont triviaux. Si $1<\nobreak\ell<\nobreak m$, on  a $\dimch(\Delta_{\ell}\Xg^{m})=\dimch(\Fg_{\ell}(\Xg))=\ell d_{\Xg}$ par l'induction en $m$, et comme $\dimch(\Delta_{\ell}\Xg^{m})\leq\dimch(\Delta_{\leq\ell}\Xg^{m})$, l'assertion
(\ref{prop-finitude-Hc-a}) pour $\ell<m$ résulte d'appliquer la majoration 
 $(\diamond)$.
Lorsque $\ell=m$, il existe un ouvert $U\dans\Xg^{m}$ tel que  $\smash{\Hc^{md_{\Xg}}}(U)\not=0$ et comme nous venons de prouver que $\dimch(\Delta_{\leq m-1}\Xg^{m})=(m{-}1)d_{\Xg}$ et que $(m{-}1)d_{\Xg}<md_{\Xg}$ puisque $d_{\Xg}\geq1$, l'exactitude de $(\ast)$ montre l'existence d'une surjection \smash{$\Hc^{m d_{\Xg}}(U\cap\Fm(\Xg))\onto\Hc^{md_{\Xg}}(U)\not=0$} et l'on conclut que $md_{\Xg}\leq \dimch(\Fm(\Xg))$. L'assertion
(\ref{prop-finitude-Hc-a}) pour $\ell=m$ résulte alors d'appliquer, une fois de plus, la majoration $(\diamond)$.

(\ref{prop-finitude-Hc-b}) 
On raisonne par induction sur $\ell$. Si $\ell\leq1$, l'assertion est claire, quel que soit $m\in\NN$. Dans le cas général, on suppose la proposition établie pour $\ellmo$ et tout $m\in\NN$. Par la suite exacte longue (\ref{suite-longue-de-base})
$$\def\Zg{\Xg^{m}}
\to\Hc^{*}(\Delta_{\ell}\Zg)\to
\Hc^{*}(\Delta_{\leq\ell}\Zg)\to
\Hc^{*}(\Delta_{\leq\ellmo}\Zg)\to\,,
$$
on a $\dim\Hc(\Delta_{\leq\ell}\Xg^{m})<+\infty$, si et seulement si, $\dim\Hc(\Delta_{\ell}\Xg^{m})<+\infty$, donc, si et seulement si $\dim\Hc(\Fg_{\ell}(\Xg))<+\infty$ (\ref{connexes}). 
Or, l'utilisation de la suite exacte longue
$$\def\Zg{\Xg^{\ell}}
\to\Hc^{*}(\Fg_{\ell}(\Xg))\to
\Hc^{*}(\Xg)^{\otimes\ell}\to
\Hc^{*}(\Delta_{\leq\ellmo}\Zg)\to\,,
$$
associée à la décomposition 
$\Xg^{\ell}=\Fg_{\ell}(\Xg)\sqcup\Delta_{\leq\ellmo}\Xg^{\ell}$, 
montre qu'une condition suffisante pour 
 $\dim\Hc(\Fg_{\ell}(\Xg))<+\infty$ est que $\dim\Hc(\Delta_{\leq\ellmo}\Xg^{\ell})<+\infty$, ce qui fait partie de l'hypothèse inductive.
\enddemo

\section{Théorèmes de scindage et complexes fondamentaux}\label{thms-de-scindage}
\glossarytitle{Théorèmes de scindage et complexes fondamentaux}
\subsection{Théorème de scindage pour $\Delta_{\leq\ell}\Xg^{m}$}
Dans ce travail, en parlant d'une suite exacte longue, on dira qu'elle est \expression{scindée} lorsque un morphisme sur trois est nul. La suite exacte longue se scinde alors en une famille dénombrable de suites exactes courtes.
$$\def\c#1{C^{#1}}
\xymatrix@C=6mm@R=3mm{
\cdots\ar[r]&\c{-1}\ar[r]^{(0)}&\c0\ar[r]&\c1\ar[r]&\c2\ar[r]^{(0)}&\c3\ar[r]&\c4\ar[r]&\c5\ar[r]^{(0)}&\c6\ar[r]&\cdots\\
\cdots\ar@{->>}[r]&\c{-1}&\c0\ar@{>->}[r]&\c1\ar@{->>}[r]&\c2&\c3\ar@{>->}[r]&\c4\ar@{->>}[r]&\c5&\c6\ar@{>->}[r]&\cdots\\
}$$

Le théorème suivant donne les deux résultats de scindage de suites longues de cohomologie à support compact qui ont motivé notre intérêt pour les espaces $i$-acycliques. Nous verrons que lorsque $\Xg$ est localement connexe, un tel scindage est intimement lié à la dégénérescence des suites spectrales de Leray (\cf\ref{autre-scindage}) associées aux fibrations $\pia:\GaDelta_{?\ell}\Xg^{m}\to\Fa(\Xg)$ de \ref{fibrations}.

\noendpoint\begin{theo}[de scindage pour $\Delta_{\leq\ell}\Xg^{m}$]\label{theo-scindage}\mynobreak
\def\varlistskips{\topsep0pt\itemsep2pt}
\begin{enumerate}
\nobreak
\item\leavevmode\label{theo-scindage-a}Pour  $m\in\NN$, si $\Xg$ est $i$-acyclique,  le morphisme de restriction
$$\Hc(\Xg\times\Fm (\Xg))\to
\Hc(\Delta_{\leq m}(\Xg\times\Fm (\Xg)))\eqno(\ddagger)$$
est nul. La suite exacte longue de cohomologie associée à la décomposition
$\medmuskip1.mu
\Xg\times\Fm =(\Fmm (\Xg))\sqcup(\Delta_{\leq m}(\Xg\times\Fm ))$
est scindée et la suite extraite:
$$0\to
\Hc(\Fm (\Xg))[-1]^{m}\to
\Hc(\Fmm (\Xg))\to
\Hc(\Xg\times\Fm (\Xg))\to0\,,
$$
est exacte.

\item\leavevmode\label{theo-scindage-b}Pour  $0<\ell\leq m\in\NN$, si $\Xg$ est $i$-acyclique,  le morphisme de restriction
$$\Hc(\Delta_{\leq\ell}\Xg^{m})\to\Hc(\Delta_{\leq\ellmo}\Xg^{m})\eqno(\ddagger\ddagger)$$
est nul. La suite exacte longue de cohomologie associée à la décomposition
$(\Delta_{\leq \ell}\Xg^{m})=((\Delta_{\ell}\Xg^{m}))\sqcup(\Delta_{\leq \ellmo}\Xg^{m})$ est scindée et la suite extraite:
$$0\to
\Hc(\Delta_{\leq \ellmo}\Xg^{m})[-1]\to
\Hc(\Delta_{\ell}\Xg^{m})\to
\Hc(\Delta_{\leq \ell}\Xg^{m})\to
0\,,
$$
est exacte.

\item\leavevmode\label{theo-scindage-c}
Soit $\Xg$ une variété topologique orientable. L'annulation
des morphismes $(\ddagger)$ pour tout $m$ (resp. 
des morphismes $(\ddagger\ddagger)$ pour tous $\ell\leq m$)
équivaut à la $i$-acyclicité de $\Xg$.
\end{enumerate}
\end{theo}
\demo (\ref{theo-scindage-a}) L'ensemble $\Delta_{\leq m}(\Xg\times\Fm )$ est celui des $(\mm)$-uplets d'éléments de $\Xg$ de la forme $(y,x_1,\ldots,x_{m})$ avec $x_i\not=x_j$ si $i\not=j$ et où $y\in\set x_1,\ldots,x_{m}/$. On comprend donc que dans le diagramme commutatif\label{demo-scindage}
$$\xymatrix@R=5mm{
\Delta_{\leq m}(\Xg\times\Fm )\ar[rd]|(0.55){\,f_2\,}\ar@{^(->}[r]|(0.58){\,f\,}&\Xg\times\Fm \ar[d]^{p_2}\\
&\Fm 
}$$
où $f$ est l'inclusion fermée, l'application $f_2:=p_{2}\circ f$ est un revêtement trivial à $m$ nappes. Les hypothèses de \ref{caracterisation}-(\ref{caracterisation-clef}) sont clairement vérifiées et le morphisme $(\ddagger)$ est bien nul.

(\ref{theo-scindage-b}) Si $\ell=1$, on a $\Delta_{\leq1}\Xg^{m}\simeq\Xg$ et $\Delta_{\leq0}\Xg^{m}=\set\pt/$, donc $(\ddagger\ddagger)$ est immédiate. Nous sommes ainsi réduits aux cas où $\ell>1$. 
On raisonne par récurrence pour $m\geq2$. Supposons l'assertion établie pour  $(m{-}1)$, et soit $1<\ell\leq m$. Par  la factorisation:
$\Delta_{\leq\ellmo}^{m}\dans
\Xg\times\Delta_{\leq\ellmo}^{\mmo}
\dans\Delta_{\leq\ell}^{m}\,,$
l'assertion résulte de montrer que le morphisme de restriction
$\Hc(\Xg\times\Delta_{\leq\ellmo}^{\mmo})\to\Hc(\Delta_{\leq\ellmo}^{m})$
est nul.
On considère pour cela le diagramme commutatif
$$\halfdisplayskips
\xymatrix@R=5mm{
f_2^{-1}(\Delta_{\ellmo}^{\mmo})\ar@{^(->}[r]\ar[rd]|{f'_{2}}&\Delta_{\leq\ellmo}^{m}\ar[rd]|{f_2}\ar@{^(->}[r]|(0.4){\,f\,}&\Xg\times\Delta_{\leq\ellmo}^{\mmo}\ar[d]^{p_2}\\
&\Delta_{\ellmo}^{\mmo}\ar@{^(->}[r]|{\;j\;}&\Delta_{\leq\ellmo}^{\mmo}
}$$
où $f$ est l'inclusion (fermée) et $j$ est l'inclusion (ouverte). 

La restriction $f'_{2}$ de $f_2:=p_{2}\circ f$ est propre. En effet, $f_{2}^{-1}(\Delta_{\ellmo}^{\mmo})$ est l'ensemble des $m$-uplets
$(y,x_1,\ldots,x_{\mmo})\in\Delta_{\ellmo}^{m}$ tels que $\#\set x_1,\ldots,x_{\mmo}/=\ellmo$ et  $y\in\set x_1,\ldots,x_{\mmo}/$, de sorte que $f'_{2}$ est une fibration trivial à $\ell{-}1$ nappes au-dessus de chaque composante ouverte $\Fg_{\pgoth}$ de $\Delta_{\ellmo}^{\mmo}$ (\ref{connexes}-(\ref{connexes-a})). Comme $\ell{-}1\geq1$ le morphisme $j_{!}:\Hc(\Delta_{\ellmo}^{\mmo})\to\Hc(\Delta_{\leq\ellmo}^{\mmo})$ est surjectif par hypothèse de récurrence. Ainsi, toutes les conditions d'application de \ref{caracterisation}-(\ref{caracterisation-clef}) sont vérifiées et  $f^{*}=0$. Le morphisme $(\ddagger\ddagger)$ est donc bien nul.

(\ref{theo-scindage-c}) En prenant $m=1$ dans (\ref{theo-scindage-a}), ou $m=\ell=2$ dans (\ref{theo-scindage-b}), l'espace $\Xg$ est \uacyclique\ donc $i$-acyclique d'après \ref{prop-acycliques}-(\ref{prop-acycliques-b}).\enddemo

\noendpoint
\begin{remas}[sur le théorème de scindage]\label{remas-theo-scindage}\begin{enumerate}\mynobreak\nobreak
\item Les\label{remas-theo-scindage-a} preuves des assertions \ref{theo-scindage}-(\ref{theo-scindage-a},\ref{theo-scindage-b}) montrent aussi l'exactitude des suites courtes
$$
\mathalign{0\to
\Hc(\Delta_{m}\Zg)[-1]^{m}\to
\Hc(\Delta_{\mm}(\Xg\times\Delta_{m}\Zg))\to
\Hc(\Xg\times\Delta_{m}\Zg)\to0\\\noalign{\kern4pt}
0\to
\Hc(\Delta_{\leq \ellmo}\Zg)[-1]\too
\Hc(\Delta_{\ell}\Zg)\too
\Hc(\Delta_{\leq \ell}\Zg)\to
0\,,}
$$
avec $\Zg\dans\Xg^{m}$ de la forme $\Yg_{1}\times\cdots\times\Yg_{m}$ où les $\Yg_{i}$ sont des sous-espaces $i$-acycliques de $\Xg$ tels que
l'ensemble $\set \Yg_{1},\ldots,\Yg_{m}/$ est totalement ordonné par la relation d'inclusion.

\item Dans\label{lambda>m-a}  \ref{theo-scindage}-(\ref{theo-scindage-b}) l'hypothèse $\ell\leq m$ est indispensable. En effet, lorsque $\ell>m$, on a  $\Delta_{\leq \ellmo}\Xg^{m}=\Delta_{\leq \ell}\Xg^{m}=\Xg^{m}$ et le morphisme $(\ddagger\ddagger)$ est l'identité. Dans un tel cas, on a aussi $\Delta_{\ellmo}^{\mmo}=\emptyset$ et l'argument à la fin de la démonstration précédente, basé sur la surjectivité de $j_{!}$, ne s'applique pas. (Cette remarque est un rapport avec la remarque \ref{rema-theo-scindage-F-lambda}.)
\end{enumerate}
\end{remas}

\begin{rema}[et corollaire]\label{rema-coro}Notons $\oo\Xg$\glossary{$\bullet$:point d'un espace}\glossary{$\Xg^{\circ}$:l'espace $\Xg$ privé d'un point $\bullet$} le complémentaire dans $\Xg$ d'un singleton $\set\bullet/$. Lorsque $\Xg$ est $i$-acyclique, il en est de même de $\oo\Xg$ et nous avons grâce à \ref{theo-scindage}-(\ref{theo-scindage-a}) un morphisme de suites exactes courtes
$$\xymatrix@C=3mm@R=5mm{
0\ar[r]&
\Hc(\Fm (\oo\Xg))[-1]^{m}\ar[r]\ar[d]&
\Hc(\Fmm (\oo\Xg))\ar[r]\ar[d]&
\Hc(\oo\Xg)\otimes\Hc(\Fm (\oo\Xg))\ar[r]\ar[d]&
0\\
0\ar[r]&
\Hc(\Fm (\Xg))[-1]^{m}\ar[r]&
\Hc(\Fmm (\Xg))\ar[r]&
\Hc(\Xg)\otimes\Hc(\Fm (\Xg))\ar[r]&
0
}$$
qui permet de montrer, par induction sur $m\geq1$, que le prolongement par zéro
$\Hc(\Fm (\oo\Xg))\to\Hc(\Fm(\Xg))$ est un morphisme surjectif.

Notons aussi que
$\halfdisplayskips\Fm (\Xg)\mmoins \Fm (\oo\Xg):=\coprod\nolimits_{i=1,\ldots,m}
\Fm ^{x_i=\bullet}(\Xg)=:\Fm ^{\bullet}(\Xg)\,,
$
où $\Fm ^{x_i=\bullet}(\Xg)$\glossary{$\Fm ^{x_i=\bullet}(\Xg)$:ensemble des $\cl x\in\Fm (\Xg)$ avec $x_i=\bullet$}\glossary{$\Fm ^{\bullet}(\Xg)$:ensemble des $\cl x\in\Fm (\Xg)$ avec l'un des $x_i$ égal à $\bullet$} désigne l'ensemble des $\cl x\in\Fm (\Xg)$ avec $x_i=\bullet$. On a une identification évidente $\Fm ^{x_i=\bullet}(\Xg)=\Fg_{\mmo}(\oo\Xg)$
et $\Fm (\Xg)\mmoins \Fm (\oo\Xg)$ se voit comme réunion disjointe fermée de $m$ copies de $\Fg_{\mmo}(\oo\Xg)$.
On a ainsi montré le corollaire suivant (qui sera utilisé dans \ref{moins-un-point-comparaison}).
\begin{coro*}Si\label{moins-un-point} $\Xg$ est $i$-acyclique, la suite de $\Sm$-modules gradués
$$\halfdisplayskips
0\to\Hc\big(\Fm ^{\bullet}(\Xg)\big)[-1]\to\Hc(\Fm (\oo\Xg))\to
\Hc(\Fm (\Xg))\to0$$
est exacte.
\end{coro*}
\end{rema}

\subsubsectionline{\`A propos des différentes hypothèses d'acyclicité.}Dans le théorème de scindage \ref{theo-scindage}, l'hypothèse d'$i$-acyclicité n'est pas optimale même lorsque $\Xg$ est une variété topologique. En effet, nous verrons dans dans \ref{U-non-I}
 que la bouteille de Klein épointée vérifie les scindages (\ref{theo-scindage-a}) et (\ref{theo-scindage-b}) alors qu'elle
 n'est n'est pas $i$-acyclique. On y montre aussi qu'elle n'est pas totalement $\cup$-acyclique (\ref{defis-0}), mais qu'elle est $\cup$-acyclique. En somme, la $i$-acyclicité est une condition suffisante et la $\cup$-acyclicité est une condition nécessaire, mais en dehors du cas où $\Xg$ est une variété topologique orientable, nous ne connaissons pas de condition sur $\Xg$ équivalente aux scindages. Nous verrons cependant, plus loin, que la condition de $i$-acyclicit\'{e} est bien nécessaire et suffisante pour la  localisation des scindages  (\ref{theo-scindage-basee}), ce qui conforte l'idée que la $i$-acyclicité est la bonne hypothèse de travail.

\comment
En fait, on peut prouver ces assertions sous cette hypothèse plus faible (\cf\ref{prop-acycliques}-(\ref{prop-acycliques-a})). 

Voici à titre d'illustration les indications de la preuve de (\ref{theo-scindage-a}) pour $\Xg$ totalement $\cup$-acyclique. En suivant la preuve en page \pageref{demo-scindage}, il nous faut justifier que les applications
$$\Hc(\Xg)\otimes \Hc(\Fm (\Xg))\to \Hc(\Fm (\Xg))\,,\quad\alpha\otimes\omega\mapsto p_{i}^{*}\alpha\cup\beta\,,$$
où $p_i:\Fm \Xg\to\Xg$, $p_i(\cl x):=x_i$, $i=1,\ldots,m$, sont nulles. Montrons que c'est plus généralement le cas de 
$$\Hc(\Xg)\otimes \Hc(\Wg)\to \Hc(\Wg)\,,\quad\alpha\otimes\omega\mapsto p_{i}^{*}\alpha\cup\beta\,,$$
pour tout ouvert $\Wg\dans\Xg^{m}$. L'idée est d'utiliser le produit des faisceaux de complexes des chaînes d'Alexander-Spanier $\Gammac(\Wg;\AS(\Xg)^{\otimes m})$ pour $\Hc(\Wg)$
\endcomment
\subsection{Complexe fondamental de $\Xg$ pour $\Delta_{\leq\ell}\Xg^{m}$}\label{complexe-DeltaXm}
\subsubsectionline{Action du groupe symétrique.}L'espace\label{action-Sm} $\Xg^m$ est muni de l'action du groupe symétrique $\Sm$\glossary{$\Sm$:groupe des permutations de l'ensemble fini $\set1,2,\ldots,m/$} par permutation des coordonnées. On pourra noter cette action par `$\Sm \rep\Xg^{m}$'\glossary{$\Sm \rep\Xg^{m}$, $\Sm \rep\Hc^{i}(\Zg)$, $\Sm \rep\Hr^{i}(\Zg)$,~\dots:action de $\Sm$ sur $\Xg^{m}$ par permutation des coordonnées et actions induites, pour $\Zg$ partie $\Sm$ stable de $\Xg^{m}$}, et de même pour celles induites sur des sous-espaces $\Sm $-stables $\Zg\dans\Xg^{m}$,  leurs cohomologies $\Sm \rep\Hc^{i}(\Zg)$, $\Sm \rep\Hr^{i}(\Zg)$,~etc\dots\ 
Les espaces de configuration généralisés $\Delta_{? \ell}\Xg^{m}$ et leurs cohomologies sont considérés munis de cette action de $\S _{m}$. 
\comment
Les suites exactes longues de cohomologie à associées à la décomposition en parties respectivement ouverte et fermée
$$\Delta_{\leq\ell}\Xg^{m}=\Delta_{\ell}\Xg^{m}\ \sqcup\ \Delta_{\leq\ellmo}\Xg^{m}\,,$$
sont des suites exactes longues de $\Sm$-modules gradués. 
\endcomment

\begin{defi}\label{complexe-fondamental}Pour $\ell\leq m\in\NN$, le \emph{complexe  fondamental de $\Xg$ pour
$\Delta_{\leq\ell}\Xg^{m}$} est, par définition, la suite des morphismes de $\Sm$-modules gradués
$$\mathrigid2mu
0\to\Hc(1)[-\ell+1]\to\cdots\to
\Hc(\ellmo)[-1]\to
\Hc(\ell)\to
\Hc(\Delta_{\leq\ell}\Xg^{m})\to0
\eqno(\star)$$
avec $\Hc(a):=\Hc(\Delta_{a}\Xg^{m})$,
où le morphisme $\Hc(a-1)[-1]\to\Hc(a)$ est la composée des morphismes
$$\Hc(\Delta_{a-1}^{m})[-1]\hf{\iota_{a-1}[-1]}{}{15mm}\Hc(\Delta_{\leq a-1}^{m})[-1]\hf{c_a}{}{7mm}\Hc(\Delta_{a}^{m})\,,$$
où $\iota_{a-1}$ est le prolongement par zéro et $c_a$ est le morphisme de liaison de la suite longue de cohomologie à support compact associée à la décomposition $\Sm$-stable $\Delta_{\leq a}^{m}=\Delta_{a}^{m}\sqcup\Delta_{\leq a-1}^{m}$ (\ref{suite-longue-de-base}). 

\nobreak La suite $(\star)$ est un complexe puisque, par construction, $\iota_a\circ c_a=0$.
\end{defi}

\nopoint\begin{theo}\label{theo-complexe-exact}
\def\varlistskips{\topsep1pt\itemsep1pt\parskip0pt\parsep0pt}\begin{enumerate}
\item\leavevmode\label{theo-complexe-exact-a}Les complexes fondamentaux d'un espace $i$-acyclique sont  exacts.
\item\leavevmode\label{theo-complexe-exact-b}Une variété topologique orientable est $i$-acyclique, si et seulement si, ses complexes fondamentaux sont exacts.
\end{enumerate}\end{theo}
\demo (\ref{theo-complexe-exact-a}) Immédiat d'après  \ref{theo-scindage}-(\ref{theo-scindage-b}), et  (\ref{theo-complexe-exact-b}) s'ensuit par \ref{theo-scindage}-(\ref{theo-scindage-c}).
\enddemo


\subsection{Le théorème de scindage pour $\Delta_{\leq\ell}\Fg^{\lambda}(\Xg)$}Nous\label{complexe-relatif}\label{generalisation} étendons dans cette section
le théorème de scindage \ref{theo-scindage}-(\ref{theo-scindage-b})
à d'autres espaces que $\Xg^{m}$, à savoir: aux produits $\smash{\Fg^{\lambda}}:=\Fg_{\lambda_1}\times\cdots\times\Fg_{\lambda_{r}}$, et aux espaces $\Fg^{\qgoth}\dans\Xg^{m}$ où $\qgoth$ est une partition de $\iii[1,m]$ (\cf\ref{FexpP}).

Ces généralisations seront utilisées dans la section \ref{caractere} pour la détermination du caractère de la représentation de $\S _{m}\rep\Hc(\Fm (\Xg))$, notamment pour expliciter la trace de l'action sur $\Hc(\Fm (\Xg))$ d'une permutation $\alpha\in\Sm$  dont la décomposition en produit cycles disjoints est de type $\lambda$ (\cf\ref{theo-trace-gen}).

\medskip
On considère les suites (peut être vides) de nombres entiers $\lambda=(\lambda_1,\ldots,\lambda_r)$ avec $\lambda_i>0$ et $r\geq0$. Lorsque $r>0$, on note 
$|\lambda|:=\sum_i\lambda_i$ et on pose
$$\Fg^{\lambda}(\Xg):=\Fg_{\lambda_1}(\Xg)\times\cdots\times\Fg_{\lambda_{r}}(\Xg)\,.$$
Pour la suite vide $()$, on conviendra que $|()|=0$ et que $\Fg_{()}(\Xg)=\set\pt/$. 

La \expression{concaténation} de deux suites
$\lambda=(\lambda_1,\ldots,\lambda_r)$ et $\lambda'=(\lambda'_1,\ldots,\lambda'_s)$ est, par définition,  la suite $\lambda\vee\lambda:=(\lambda_1,\ldots,\lambda_r,\lambda'_1,\ldots,\lambda'_s)$. On a donc
$$\Fg_{\lambda\vee\lambda'}=\Fg_{\lambda}\times\Fg_{\lambda'}\,,$$
ainsi que les égalités évidentes 
$$\Fg^{(1,\ldots,1)}(\Xg)=\Xg^{m}\,,\quad
\Fg^{(m)}(\Xg)=\Fm (\Xg)=\Delta_{m}\Fg^{(1,\ldots,1)}(\Xg)\,.$$

\def\aa{{a}}
\begin{theo}[de scindage]\label{theo-scindage-F-lambda}Soit $\Xg$ un espace $i$-acyclique. 
\begin{enumerate}
\mynobreak\nobreak\item\leavevmode\label{theo-scindage-F-lambda-a}Pour tout $\aa \geq 0$ et tout $\ell\in\NN$, la décomposition en parties respectivement ouverte et fermée
$$\Delta_{\ell}\big(\Xg\times\Fg_{\aa }\times\Fg^{\lambda}\big)=
\Delta_{\ell}\big(\Fg_{1+\aa }\times\Fg^{\lambda}\big)\sqcup
\Delta_{\ell}\big((\Delta_{\aa }(\Xg\times\Fg_{\aa }))\times\Fg^{\lambda}\big)
$$
donne lieu à une suite exacte longue de cohomologie dont la suite courte extraire
$$\skip0=3cm
\mathalign{0\to\Hc(\Delta_{\ell}(\Fg_{1+\aa }\times\Fg^{\lambda}))\to
\Hc(\Delta_{\ell}(\Xg\times\Fg_{\aa }\times\Fg^{\lambda}))\to\hfill\\
\hskip2\skip0
\to
\Hc(\Delta_{\ell}(\Delta_{\aa }(\Xg\times\Fg_{\aa })\times\Fg^{\lambda}))
\to0\hfill}\postskip-.5ex$$
est exacte.
\item\leavevmode\label{theo-scindage-F-lambda-b}Pour tout $\ell\leq|\lambda|$, le morphisme de restriction
$$\Hc(\Delta_{\leq\ell}(\Fg^{\lambda}))\to\Hc(\Delta_{\leq\ellmo}(\Fg^{\lambda}))$$
est nul. 
La suite exacte longue de cohomologie associée à la décomposition en parties respectivement ouverte et fermée  
$$\Delta_{\leq\ell}(\Fg^{\lambda})=\Delta_{\ell}(\Fg^{\lambda})\sqcup
\Delta_{\leq\ellmo}(\Fg^{\lambda})$$ est scindée et les suites courtes extraites:
$$\def\Xg{\Fg^{\lambda}}
0\to
\Hc(\Delta_{\leq \ellmo}(\Xg))[-1]\to
\Hc(\Delta_{\ell}(\Xg))\to
\Hc(\Delta_{\leq \ell}(\Xg))\to
0\,,
$$
sont exactes. 
\end{enumerate}
\end{theo}
\demo \def\relation{\gooddownstackrel{0pt}\sim\pgoth}L'assertion (\ref{theo-scindage-F-lambda-a}) est évidente si $\aa =0$. Pour $\aa>0$, notons $m:=1+\aa +|\lambda|$ et considérons l'inclusion
$$\Xg\times\Fg_{\aa }\times\Fg^{\lambda}
\dans
\Xg\times\Xg^{\aa}\times\Xg^{|\lambda|}=\Xg^{m}\,.$$
Comme $m>0$, le cas $\ell=0$ est immédiat, on peut donc supposer que $\ell>0$.
Nous avons montré dans \ref{connexes} que $\Delta_{\ell}\Xg^{m}$ se décompose en réunion disjointe des sous-espaces ouverts $\Fg_{\pgoth}(\Xg)$ où $\pgoth:=\set I_1,I_2,\ldots,I_\ell/$ désigne une partition de  $\iii[1,m]$ en $\ell$ parties non vides. 
L'ensemble $\Fg_{\pgoth}(\Xg)$ est l'ensemble des $m$-upltes $(x_1,\ldots,x_{m})$ tels que $(x_i{=}x_j)\Leftrightarrow (i\relation j)$ où l'on note $(i\relation j)\Leftrightarrow_{\rm def} (\exists k)(\set i,j/\dans I_k)$ (\cf\ref{nota-partitions}). Il s'ensuit que $\Xg\times\Fg_{\aa }\times\Fg_{\lambda_1}\times\cdots\times\Fg_{\lambda_r}$ est  la réunion des $\Fg_{\pgoth}$ tels que la relation \smash{$(i\gooddownstackrel{0pt}\sim\pgoth j)$} \emph{n'est pas vérifiée} lorsque $i$ et $j$ sont des coordonnées correspondantes à un même facteur $\Fg_{\aa}$ ou $\Fg_{\lambda_i}$. On obtient ainsi une partition de $\Delta_{\ell}(\Xg\times\Fg_{\aa }\times\Fg^{\lambda})$ en sous-espaces ouverts et fermés de $\Delta_{\ell}\Xg^{m}$. Ces mêmes remarques s'appliquent à $\Delta_{\ell}(\Fg_{1+\aa }\times\Fg^{\lambda})$ qui est alors aussi un sous-espace ouvert et \emph{fermé} de $\Delta_\ell\Xg^{m}$, de même donc qu'à son complémentaire
$\Delta_{\ell}\big(\Delta_{\aa}(\Xg\times\Fg_{\aa })\times\Fg^{\lambda}\big)$.
On donc une décomposition en somme directe
$$\Hc\big(\Delta_{\ell}(\Xg\times\Fg_{\aa }\times\Fg^{\lambda})\big)=
\Hc\big(\Delta_{\ell}(\Fg_{1+\aa }\times\Fg^{\lambda})\big)\oplus
\Hc\big(\Delta_{\ell}\big(\Delta_{\aa }(\Xg\times\Fg_{\aa })\times\Fg^{\lambda}\big)\big)\,,$$
dont on déduit aussitôt  (\ref{theo-scindage-F-lambda-a}).

\smallskip
\noindent(\ref{theo-scindage-F-lambda-b}) {\sl Lemme préliminaire. Si $\Zg \dans\Xg^{m}$ est telle que
$$\halfdisplayskips\big(\Hc(\Delta_{\leq\ell}\Zg )\to \Hc(\Delta_{\leq\ellmo}\Zg )\big)=0\,,\mrlap{\quad\forall
\ell\leq m\,,}\postskip0ex\hskip5mm
\eqno(\ast)$$
alors, on a aussi}
$$\preskip1ex
\hskip-6mm\big(\Hc(\Delta_{\leq\ell}(\Xg\times\Zg ))\to \Hc(\Delta_{\leq\ellmo}(\Xg\times\Zg ))\big)=0\,,\relax{\ \ \forall
\ell\leq m+1\,.}\hskip-8mm\eqno(**)$$

\noindent{\slshape Preuve du lemme. }On peut supposer $\ell{-}1>0$. L'inclusion $\Delta_{\leq\ellmo}(\Xg\times\Zg )\dans\Delta_{\leq\ell}(\Xg\times\Zg )$ se factorise en
$$
\def\dans{\ \subset\ }
\Delta_{\leq\ellmo}(\Xg\times\Zg )\xymatrix@C=6mm{\ar@{^(->}[r]|(0.45){\,f\,}&}
\Xg\times\Delta_{\leq\ellmo}\Zg 
\xymatrix@C=6mm{\ar@{^(->}[r]&}\Delta_{\leq\ell}(\Xg\times\Zg )\,,$$
et l'on a le diagramme commutatif
$$\preskip2pt
\xymatrix@R=4mm{
f_2^{-1}(\Delta_{\ellmo}\Zg )\ar@{^(->}[r]\ar[rd]|{f'_{2}}&
\Delta_{\leq\ellmo}(\Xg\times\Zg )
\ar[rd]|{f_2}\ar@{^(->}[r]|(0.5){\,f\,}&
\Xg\times\Delta_{\leq\ellmo}\Zg 
\ar[d]^{p_2}\\
&
\Delta_{\ellmo}\Zg 
\ar@{^(->}[r]|{\;j\;}&
\Delta_{\leq\ellmo}\Zg 
}$$
où $f_{2}:=p_2\circ f$. On raisonne alors comme dans la preuve de \ref{theo-scindage}-(\ref{theo-scindage-b}). L'application $j_{!}:\Hc(\Delta_{\ellmo}\Zg )\to\Hc(\Delta_{\leq\ellmo}\Zg )$ est surjective par l'hypothèse $(\ast)$. Si $\Delta_{\ellmo}\Zg=\emptyset$, on déduit aussitôt que $\Hc(\Xg\times\Delta_{\leq\ellmo}\Zg)=0$ et alors $f^{*}=0$. Autrement, la restriction $f'_{2}$ de $f_2$ est propre ($\ell{-}1>0$), les hypothèses de \ref{caracterisation}-(\ref{caracterisation-clef}) sont donc vérifiées et $f^{*}=0$. Tout ceci implique  $(**)$.\hfill$\boxminus$

\smallskip
Prouvons (\ref{theo-scindage-F-lambda-b}) par induction sur le nombre $r$ des termes de $\lambda=(\lambda_1,\ldots,\lambda_r)$.

\noindent {\bf Cas \boldmath$r\leq1$. }Si $r=0$, on a $\ell\leq0$ et donc $\Delta_{\ellmo}(\_)=\emptyset$. Si $r=1$, $\Fg^{\lambda}$ est de la forme $\Fg_{\aa}$, mais alors $\ell\leq\aa$ et $\Delta_{\leq\ell}\Fg_{\aa}\not=\emptyset$ seulement si $\ell=\aa$, auquel cas $\Delta_{\ellmo}\Fg_{\aa}=\emptyset$.

\noindent {\bf Cas \boldmath$r>1$. }Supposons (\ref{theo-scindage-F-lambda-b}) établie pour toute suite $\lambda$ comportant $r\geq1$ termes. Soit $\lambda=(\lambda_1,\ldots,\lambda_r)$ et montrons,  par induction sur l'entier $\aa$, que $\Fg_{\aa}\times\Fg^{\lambda}$ vérifie également (\ref{theo-scindage-F-lambda-b}).

Le cas $\aa=1$ était l'objet du lemme préliminaire. Supposons maintenant que $\Fg_{\aa}\times\Fg^{\lambda}$ vérifie (\ref{theo-scindage-F-lambda-b}) et montrons qu'il en est de même de $\Fg_{1+\aa}\times\Fg^{\lambda}$.

Pour $\ell\leq1+\aa+|\lambda|$, on considère le diagramme commutatif $(\DDD_{\ell})$
suivant où l'on a omis d'écrire le symbole $\Hc(\_)$ pour gagner de la place.
\begin{figure}[!]
$$\preskip-0ex\hss
\def\to{\ar[r]}
\def\d{\ar[d]}
\def\ds{\ar@{->>}[d]}
\def\di{\ar@{>->}[d]}\def\ui{\ar@{<-<}[u]}
\def\xFg{\times\Fg^{\lambda}}
\def\tt{\vrule depth7pt width0pt }
\xymatrix@C=4mm@R=6.5mm{
&\d|(0.35){\ominus}&\d|(0.35){\circledcirc}&\d|(0.35){\otimes}\\
\to|(0.11){\circleplustimes}&
\Delta_{\leq\ellmo}(\Fg_{1+\aa}\xFg)[-1]\to^(0.58){\gamma}\d&
\Delta_{\ell}(\Fg_{1+\aa}\xFg)\to&
\Delta_{\leq\ell}(\Fg_{1+\aa}\xFg)\to|(0.84){\circleplustimes}\d&\\
\to|(0.11){\circledcirc}&
\Delta_{\leq\ellmo}(\Xg\times\Fg_{\aa}\xFg)[-1]\ar@{>->}[r]\d&
\Delta_{\ell}(\Xg\times\Fg_{\aa}\xFg)\ar@{->>}[r]\ui+<0pt,-11pt>\ds&
\Delta_{\leq\ell}(\Xg\times\Fg_{\aa}\xFg)\to|(0.84){\circledcirc}\d^{\alpha}&\\
\to|(0.11){\odot}&
\Delta_{\leq\ellmo}(\aa\Fg_{\aa}\xFg)[-1]\to\d|(0.55){\ominus}&
\Delta_{\ell}(\aa\Fg_{\aa}\xFg)\to^{\beta}\d|(0.55){\circledcirc}&
\Delta_{\leq\ell}(\aa\Fg_{\aa}\xFg)\to|(0.84){\odot}\d|(0.55){\otimes}&\\
&&&
 }\hss\postskip-0.25ex$$
\mycaption{\small Diagramme $\DDD_{\ell}$}
\end{figure}
On y a noté~$\aa\Fg_{\aa}$, l'espace qui consiste en $\aa$ copies de l'espace $\Fg_{\aa}$, ce qui correspond très exactement à $\Delta_{\aa}(\Xg\times\Fg_{\aa})$. Le diagramme $(\DDD_{\ell})$ se prolonge indéfiniment par ses quatre côtés via les suites longues de cohomologie à support compact. Dans ce prolongement, le nombre $\ell$ reste bien sûr constant et seuls les degrés cohomologiques changent.
Les quatre flèches marquées `$\circledcirc$' y sont nulles. Cela résulte, pour la colonne et la ligne centrales, respectivement par l'assertion~(\ref{theo-scindage-F-lambda-a}) et par le lemme préliminaire.

Lorsque $\ell\leq\aa+|\lambda|$, l'hypothèse inductive s'applique et les flèches `$\odot$' sont nulles. Le morphisme $\beta$ et surjectif et on en déduit la surjectivité de 
$\alpha$. Les flèches  `$\otimes$' sont donc nulles. Comme ces propriétés sont également vérifiées sur le diagramme $(\DDD_{\ellmo})$ et que les flèches `$\otimes$' de $(\DDD_{\ellmo})$ sont les flèches `$\ominus$' de $(\DDD_{\ell})$, on conclut que ces dernières sont nulles. \`A partir de là, une chasse au diagramme élémentaire montre que le flèches `$\circleplustimes$' sont également nulles, ce qui prouve l'assertion (\ref{theo-scindage-F-lambda-b}) pour l'espace $\Fg_{1+\aa}\times\Fg^{\lambda}$ et pout $\ell\leq\aa+|\lambda|$.

Soit maintenant $\ell=1+\aa+|\lambda|$. Les flèches de la première colonne de $(\DDD_{1+\aa+|\lambda|})$ marquées `$\ominus$' sont nulles puisqu'elles coïncident avec  les flèches `$\otimes$' du diagramme $(\DDD_{\aa+|\lambda|})$ dont la nullité a déjà été établie. On en déduit l'injectivité de $\gamma$, et donc la nullité des flèches `$\circleplustimes$'. Ceci prouve l'assertion (\ref{theo-scindage-F-lambda-b}) pour l'espace $\Fg_{1+\aa}\times\Fg^{\lambda}$ et pout $\ell=1+\aa+|\lambda|$.

La récurrence par rapport à $\aa$ est terminée et l'assertion (\ref{theo-scindage-F-lambda-b}) est vérifiée par $\Fg_{\aa}\times\Fg^{\lambda}$ pour tout $\aa\in\NN$ et pour tout $\ell\leq\aa+|\lambda|$. Ceci à son tour termine la récurrence par rapport au nombre $r$ dans $\lambda=(\lambda_1,\ldots,\lambda_r)$. L'assertion (\ref{theo-scindage-F-lambda-b}) est par conséquent vérifiée par $\Fg^{\lambda}$ pour tout $\lambda=(\lambda_1,\ldots,\lambda_r)$, pour tout $r\in\NN$ et pour tout $\ell\leq|\lambda|$. \cqfd
\enddemo

\begin{rema}On\label{rema-theo-scindage-F-lambda} prendra garde du fait que l'analogue de \ref{theo-scindage-F-lambda}-(\ref{theo-scindage-F-lambda-a}) pour l'opérateur $\Delta_{\leq\ell}$, \idest l'exactitude des suites 
$$\halfdisplayskips
\skip0=3cm
\mathalign{0\to\Hc(\Delta_{\leq\ell}(\Fg_{1+\aa }\times\Fg^{\lambda}))\to
\Hc(\Delta_{\leq\ell}(\Xg\times\Fg_{\aa }\times\Fg^{\lambda}))\to\hfill\\\noalign{\kern0pt}
\hskip2\skip0
\to
\Hc(\Delta_{\leq\ell}(\Delta_{\aa }(\Xg\times\Fg_{\aa })\times\Fg^{\lambda}))
\to0\,,\hfill}$$
n'est pas vrai pour $\ell=1+\aa+|\lambda|$, cas auquel
on peut effacer  $\Delta_{\leq\ell}(\_)$. Dans ce cas, le morphisme de droite est la restriction
$\Hc(\Xg\times\Fg_{\aa}\times\Fg^{\lambda})\to\Hc(\aa\Fg_{\aa}\times\Fg^{\lambda})$ qui est nulle d'après 
\ref{theo-scindage}-(\ref{theo-scindage-a}), alors que le morphisme de gauche n'est généralement pas injectif. On remarquera que la démonstration du théorème  prouve cependant que la suite en question est bien exacte pour $\ell\leq\aa+|\lambda|$.
\end{rema}

\subsection{Sous-espaces $\Fg^{\qgoth}(\Xg)$ et sous-groupes $\S ^{\qgoth}\dans\S _{m}$}\label{FexpP}
\subsubsectionline{Les sous-espaces $\Fg^{\qgoth}(\Xg)$.}Soit $\qgoth=\set I_1,\ldots,I_{r}/$ une partition de $\iii[1,m]$ en parties non vides, notons $\lambda_{i}:=|I_i|$ et  $\lambda:=(\lambda_1,\ldots,\lambda_r)$.\label{nota-qgoth}
Pour chaque $i=1,\ldots,r$, fixons arbitrairement une bijection $\varphi_i:I_i\to\iii[1,\lambda_i]$. Notons ensuite $\varphi:\iii[1,m]\to\iii[1,m]$ la bijection 
$$t\in I_{i}\mapsto\varphi(t)=\sumnl_{j<i}\lambda_j+\varphi_{i}(t)\,,$$  et soit
$\varPhi:\Xg^{m}\to\Xg^{m}$
l'homéomorphisme $(x_1,\ldots,x_{m})\mapsto(x_{\varphi(1)},\ldots,x_{\varphi(m)})$.
L'ensemble\glossary{$\Fg^{\qgoth}(\Xg)$: Si $\qgoth\in\Pgoth(m)$, $\cl x\in\Fg^{\qgoth}(\Xg)\Leftrightarrow_{\rm def}(\forall i\not=j)((i\relq j)\Rightarrow(x_i\not=x_j))$}
$$\halfdisplayskips\Fg^{\qgoth}(\Xg):=\varPhi^{-1}(\Fg^{\lambda}(\Xg))=\varPhi^{-1}(\Fg_{\lambda_1}\times\cdots\times\Fg_{\lambda_{m}})\,,$$
est indépendant de l'indexation des parties $I_i\in\qgoth$ et des choix des bijections $\varphi_{i}$, il dépend uniquement de la partition $\qgoth$.

\def\Remarque{Commentaire à propos des notations}
\begin{rema}\label{notation-partition}Il convient de souligner la différence entre les notation $\Fg_{\pgoth}(\Xg)$ et $\Fg^{\qgoth}(\Xg)$. Dans les deux cas, $\pgoth$ et $\qgoth$ désignent des partitions de $\iii[1,m]$ et nous avons (\footnote{On rappelle que si $\pgoth$ est une partition de $\iii[1,m]$, on écrit
 $(i\relp j)$, si et seulement si, il existe $I\in\pgoth$ tel que $\set i,j/\dans I$.})
$$
\mathalign{
(x_1,\ldots,x_{m})\in\Fg_{\pgoth}(\Xg)&\Leftrightarrow_{\rm def}&
(\forall i\not=j)\big((i\relp j)\Leftrightarrow(x_i=x_j)\big)\\
(x_1,\ldots,x_{m})\in\Fg^{\qgoth}(\Xg)&\Leftrightarrow_{\rm def}&
(\forall i\not=j)\big((i\relq j)\Rightarrow(x_i\not=x_j)\big)
}$$
\end{rema}

\subsubsectionline{Sous-groupe $\S ^{\qgoth}\dans\S _{m}$.}Pour\label{Sqgoth} toute partie $I\dans\iii[1,m]$, notons 
$$\S _{I}:=\Fix\nolimits_{\S _{m}}(\iii[1,m]\mmoins I)
=\set \alpha\in\S _{m}\mid \alpha(j)=j\quad\forall j\not\in I/\,,$$
puis\label{SgI}, si $\qgoth=\set I_1,\ldots,I_{r}/$ est une partition de $\iii[1,m]$, posons
$$\S ^{\qgoth}:=\S _{I_1}\times\cdots\times\S _{I_r}\,.$$
On a $\S _{I}\sim\S _{|I|}$ et 
$\S ^{\qgoth}\sim\S _{|I_1|}\times\cdots\times\S _{|I_r|}$. 

L'action de $\S ^{\qgoth}$ sur $\Xg^{m}$ laisse clairement stables les sous-espaces $\Delta_{?\ell}(\Fg^{\qgoth})$, quel que soit $\ell$. La section suivante \ref{complexe-F-goth}
 étend les résultats de la section \ref{complexe-DeltaXm} qui concernaient les $\S _{m}$-espaces $\Delta_{?\ell}\Xg^{m}$, au cas des $\S ^{\qgoth}$-espaces $\Delta_{?\ell}\Fg^{\qgoth}(\Xg)$.

\begin{defi}\label{def-transversalite}Deux partitions $\pgoth,\qgoth\in\Pgoth(\iii[1,m])$ sont dites \expression{transverses}, et l'on note, $\pgoth\trans\qgoth$, si l'on a $(\forall i\not=j)\big(i\relp j\Rightarrow i\not\relq j\big)$ (\cf\ref{notation-partition}). On note 
$$\qgoth^{\trans}:=\set \pgoth\in\Pgoth(\iii[1,m])\mid\pgoth\trans\qgoth/\,.$$
La relation de transversalité est une relation symétrique.
\end{defi}

\nopoint\begin{prop}\label{prop-transversalite}
\begin{enumerate}
\item Pour\label{prop-transversalite-a}  $\qgoth\in\Pgoth(m)$, on a
\smash{$\Fg^{\qgoth}(\Xg):=\coprod\nolimits_{\let\stacklhss\relax
\relax{
\pgoth\trans\qgoth}}\Fg_{\pgoth}(\Xg)$.}
\item 
Le\label{prop-transversalite-b} sous-groupe $\S ^{\qgoth}\dans\S _{m}$ (\ref{Sqgoth}) fixe $\qgoth$ et laisse stable $\qgoth^{\trans}$.
\end{enumerate}
\end{prop}
\demo (\ref{prop-transversalite-a}) On rappelle que l'on a, par définition (\ref{notation-partition}),
$$\let\Big\big\begin{casesalign}
\Big((x_1,\ldots,x_{m})\in\Fg_{\pgoth}\Big)
&\Longleftrightarrow&
(\forall i\not=j)\Big((i\relp j)
\Leftrightarrow
(x_i=x_j)\Big)\,,\\
\Big((x_1,\ldots,x_{m})\in\Fg^{\qgoth}\Big)
&\Longleftrightarrow&
(\forall i\not=j)\Big(
(i\relq j)\Rightarrow
(x_i\not=x_j)\Big)\,.
\end{casesalign}$$
L'inclusion $\Fg_{\pgoth}\dans\Fg^{\qgoth}$ est alors immédiate si $\pgoth\trans\qgoth$. Réciproquement, un élément $\cl x\in\Fg^{\qgoth}$ définit la partition $\pgoth_{\cl x}$ qui regroupe dans une même partie les indices des coordonnées identiques, \idest $(i\varrel{\pgoth_{\!\scriptscriptstyle \cl x}}j)\Leftrightarrow (x_i=x_j)$, or $(x_i=x_j)\Rightarrow (i\not\relq j)$. Par conséquent, $\pgoth_{\cl x}\trans\qgoth$ et $\cl x\in\Fg_{\pgoth_{\scriptscriptstyle \cl x}}\dans\Fg^{\qgoth}$.

 (\ref{prop-transversalite-b}) Pour tout $\alpha\in\S _{m}$, 
  si $\pgoth\trans\qgoth$, on a $\alpha\cdot\pgoth\trans\alpha\cdot\qgoth$. D'autre part, si $\alpha\in\S ^{\qgoth}$, l'égalité $\alpha\cdot\qgoth=\qgoth$ est évidente.
\enddemo

\subsection{Complexe fondamental de $\Xg$ pour $\Delta_{\leq\ell}\Fg^{\qgoth}(\Xg)$}\label{complexe-F-goth}

\subsubsectionline{Le théorème de scindage pour $\Delta_{\leq\ell}\Fg^{\qgoth}(\Xg)$.}Compte tenu de l'homéomorphisme
$\Fg^{\lambda}\sim\Fg^{\qgoth}$ de \ref{nota-qgoth}, le corollaire suivant de \ref{theo-scindage-F-lambda} est immédiat.

\begin{coro}\label{coro-scindage-F-pgoth}Soit $\Xg$ un espace $i$-acyclique. Pour toute partition $\qgoth$ de $\iii[1,m]$ et pour tout $\ell\leq m$, le morphisme de restriction
\halfdisplayskips
$$\Hc(\Delta_{\leq\ell}(\Fg^{\qgoth}))\to\Hc(\Delta_{\leq\ellmo}(\Fg^{\qgoth}))$$
est nul. 
La suite exacte longue de cohomologie associée à la décomposition en parties respectivement ouverte et fermée  
$$\Delta_{\leq\ell}(\Fg^{\qgoth})=\Delta_{\ell}(\Fg^{\qgoth})\sqcup
\Delta_{\leq\ellmo}(\Fg^{\qgoth})$$ est scindée et les suites courtes extraites:
$$
\def\Xg{\Fg^{\qgoth}}
0\to
\Hc(\Delta_{\leq \ellmo}(\Xg))[-1]\to
\Hc(\Delta_{\ell}(\Xg))\to
\Hc(\Delta_{\leq \ell}(\Xg))\to
0\,,
$$
sont exactes. 
\end{coro}

\begin{defi}\label{complexe-canonique-pgoth}Soit $\qgoth$ une partition de $\iii[1,m]$ et soit $\ell\leq m$.
Le \emph{complexe  fondamental de $\Xg$ associé à
$\Delta_{\leq\ell}(\Fg^{\qgoth})$} est le complexe de $\S^{\qgoth}$-modules gradués
$$\mathrigid2mu
0\to\Hc(1)[-\ell+1]\to\cdots\to
\Hc(\ellmo)[-1]\to
\Hc(\ell)\to
\Hc(\Delta_{\leq\ell}(\Fg^{\qgoth}(\Xg)))\to0
$$
avec $\Hc(a):=\Hc(\Delta_{a}(\Fg^{\qgoth}))$, obtenu par concaténation des suites courtes du corollaire \ref{coro-scindage-F-pgoth} (\cf\ref{complexe-fondamental}).\killline
\end{defi}

\begin{theo}\label{theo-complexe-pgoth-exact}Soit $\Xg$ un espace $i$-acyclique. Le complexe fondamental de $\Delta_{\leq\ell}(\Fg^{\qgoth})$ est un complexe de $\S^{\qgoth}$-modules gradués exact.\killline
\end{theo}
\demo Conséquence du corollaire \ref{coro-scindage-F-pgoth}. 
\enddemo

\section{Cohomologie des espaces de configuration, cas $i$-acyclique}\label{coh-cas-i-acyclique}
\glossarytitle{Cohomologie des espaces de configuration, cas $i$-acyclique}
\subsection{Généralités sur le polynôme de Poincaré}
\subsubsection{Polynôme de Poincaré d'un espace vectoriel gradué}Le \emph{polynôme de Poincaré} d'un
$k$-espace vectoriel gradué $V:=\bigoplus_{i\in\ZZ}V^i$ de dimension  finie est 
le polynôme $\P(V)\in\ZZ[T,1/T]$\glossary{$\P(V)\in\ZZ[T]$:polynôme de Poincaré d'un espace vectoriel gradué $V$}
$$\P(V^*)(T):=\sumnl_{i\in\ZZ}\dim_k(V^i)\,T^{i}\,.$$

\nopoint\begin{lemm}\label{poincare-elem}
\begin{enumerate}
\item\leavevmode\label{poincare-elem-a}Si $0\to V_1\to V_2\to V_3\to0$ est une suite exacte courte d'espaces vectoriels gradués, on a $\P(V_2)=\P(V_1)+\P(V_3)$.
\item\leavevmode\label{poincare-elem-b}$\P(V[-1])(T)=T\cdot\P(V)$ {\rm (\footnote{On rappelle que $V[-1]^{i}:=V^{i-1}$, par convention.})}.
\item\leavevmode\label{poincare-elem-c}$\P(V_1\otimes_k V_2)=\P(V_1)\cdot\P(V_2)$.
\item\leavevmode\label{poincare-elem-d}$\P(\mathop{\rm Homgr}_{k}(V_1,V_2))(T)=\P(V_1)(1/T)\cdot\P(V_2)(T)$
\end{enumerate}
\end{lemm}

\subsubsection{Polynômes de Poincaré d'un espace topologique}
Les\label{Poincare-r-c} polynômes de Poincaré pour la cohomologie ordinaire et à support compact d'un espace topologique $\Xg$ de type fini\index{type fini} ($^{\ref{tf}}$) sont notés\glossary{${\Pr(\Xg)},{\Pc(\Xg)}$:polynômes de Poincaré de $\Hr(\Xg)$ et $\Hc(\Xg)$ respectivement}
$$
\Pr(\Xg,k)(t):=\P(\Hr(\Xg,k))\text{\quad et\quad }
\Pc(\Xg,k)(t):=\P(\Hc(\Xg,k))\,.
$$
Le corps $k$ sera omis de ces notations lorsque son indication sera superflue.

\begin{rema}[et notation]Comme nous l'avons déjà signalé dans la remarque
\ref{complementaire-partie-finie}, lorsque\label{X-i}\label{X-a}  $\Hc^{0}(\Xg)=0$, le polynôme de Poincaré pour la cohomologie à support compacte de $\Xg\mmoins F$ seul dépend du cardinal $a:=\#F$ de la partie finie $F\dans\Xg$. C'est le sens de la notation $\Pc(\Xg\mmoins a)$\glossary{$\Pc(\Xg\mmoins a)$: Polynôme de Poincaré compact du complémentaire dans $\Xg$ d'une partie finie de cardinal $a$}, p.e. dans l'assertion  (\ref{poincare-espaces-c}) du lemme suivant.
\end{rema}

\begin{lemm}\label{poincare-espaces}Si $\Xg$ et $\Yg$ sont des espaces topologiques de type fini, on a 
\begin{enumerate}
\item\leavevmode\label{poincare-espaces-a}$\Pr(\Xg\times\Yg)=\Pr(\Xg)\cdot\Pr(\Yg)$ 
et 
$\Pc(\Xg\times\Yg)=\Pc(\Xg)\cdot\Pc(\Yg)$.
\item\leavevmode\label{poincare-espaces-b}(Dualité de Poincaré) Si de plus $\Xg$ est une variété topologique orientable, de dimension $d_{\Xg}$ et de type fini, on a 
$$\P(\Xg)(T)=\Pc(\Xg)(1/T)\cdot T^{d_{\Xg}}\,.$$
\item Si\label{poincare-espaces-c} $\Hc^{0}(\Xg)=0$ (p.e. si $\Xg$ est $i$-acyclique), on a 
$$\Pc(\Xg\mmoins a)(T)=\Pc(\Xg)(T)+a\cdot T\,,$$
et donc
$$
{\prodnl_{i=0}^{\mmo}\Pc(\Xg\mmoins i)(T)\over T^{m}}=\Big({\Pc(\Xg)(T)\over T}\Big)\osp m\,.$$
\end{enumerate}
\end{lemm}
\comment\demo(\ref{poincare-espaces-c}) Comme $\Hc^{0}(\Xg)=0$, la suite exacte longue de cohomologie à support compact
$\cdots\to\Hc(\Xg\mmoins F_i)\to\Hc(\Xg)\mathop{\to}\limits^{\rho}\Hc(F)\to\cdots$
est scindée ($\rho=0$) pour tout $F\dans\Xg$ de cardinal $a$.
On conclut par \ref{poincare-elem}.
\enddemo
\endcomment
\subsection{Polynômes de Poincaré de $\Fm (\Xg)$ et de $\Delta_{\leq \mmo}\Xg^{m}$}La proposition suivante est une application immédiate du théorème  de scindage \ref{theo-scindage}-(\ref{theo-scindage-a}). Elle donne une formule fermée  pour les polynômes de Poincaré de deux espaces de configuration. 

\def\Remarque{\miseengarde Avertissement}
\begin{rema*}{\miseengarde La donnée d'un
espace $i$-acyclique $\Xg$ présupposera que $\dim\Hc(\Xg)<\infty$.}
\end{rema*}

\begin{prop}\label{pol-poincare}Soit $\Xg$ un espace $i$-acyclique.
\mynobreak\begin{enumerate}\nobreak
\item\leavevmode\label{pol-poincare-a}Le polynôme de Poincaré $\Pc(\Fm(\Xg))$
est le polynôme
$$\Pc(\Fm (\Xg))(T)=\prodnl_{i=0}^{\mmo}\Pc(\Xg\mmoins i)(T)=\prodnl_{i=0}^{\mmo}\Big(\Pc(\Xg)(T)+i\cdot T\Big)\,,$$
soit,
$$
\smashtop{\Pc(\Fm (\Xg))(T)= T^{m}\,\Big({\Pc(\Xg)(T)\over T}\Big)\osp m\,.}$$
\item\leavevmode\label{pol-poincare-b}L'espace $\Delta_{\leq \mmo}\Xg^{m}$ est \expression{la diagonale épaisse de $\Xg^{m}$}\index{diagonale épaisse}, son polynôme de Poincaré est le polynôme 
$\Pc(\Delta_{\leq \mmo}\Xg^{m})$ homogène  de degré $\mmo$ dans l'anneaux $\ZZ[\Pc(\Xg),T]$, vérifiant:
$${\Pc(\Delta_{\leq \mmo}\Xg^{m})(T)\over T^{\mmo}}=\Big({\Pc(\Xg)(T)\over T}\Big)\osp m-\Big({\Pc(\Xg)(T)\over T}\Big)^{m}\,\cdot$$
\end{enumerate}\end{prop}
\demo (\ref{pol-poincare-a}) La suite exacte courte de \ref{theo-scindage}-(\ref{theo-scindage-a}) et les lemmes \ref{poincare-elem} et \ref{poincare-espaces} donnent aussitôt la relation de récurrence:
$$\mathalign{\Pc(\Fmm )&=&\Pc(\Xg)\cdot\Pc(\Fm )+m\cdot T\cdot \Pc(\Fm )\hfill\\\noalign{\kern2pt}
&=&(\Pc(\Xg)+m\cdot T)\cdot\Pc(\Fm )=\Pc(\Xg\mmoins m)\cdot\Pc(\Fm )\hfill}$$
qui permet de conclure.

(\ref{pol-poincare-b}) Montrons la surjectivité du morphisme de prolongement par zéro
$$\Hc(\Fm )\onto\Hc(\Xg^m)\,.\eqno(\ast)$$ 
Immédiat si $m=1$, on raisonne par induction. Si  $\Hc(\Fg_{\mmo})\to\Hc(\Xg^{\mmo})$ est surjectif, le morphisme
$\Hc(\Xg\times\Fg_{\mmo})\onto\Hg(\Xg\times\Xg^{\mmo})\,,$
l'est aussi (par Künneth), et, composé à $\Hc(\Fm )\onto\Hc(\Xg\times\Fg_{\mmo})$, surjectif d'après \ref{theo-scindage}-(\ref{theo-scindage-a}), on conclut que ($\ast$) l'est aussi.

Maintenant, le complémentaire de l'ouvert $\Fm $ dans $\Xg^m$ est la diagonale épaisse $\Delta_{\leq \mmo}\Xg^m$, et la suite longue de cohomologie à support compact associée à la décomposition $\Xg^m=\Fm \sqcup\Delta_{\leq \mmo}^{m}$ est scindée au niveau du prolongement par zéro $(\ast)$. On a donc la suite exacte courte:
$$0\to\Hc(\Delta_{<m}^{m})[-1]\to\Hc(\Fm )\to\Hc(\Xg^m)\to0\,,$$
dont résulte l'assertion (\ref{pol-poincare-b}).
\enddemo

\begin{rema}Dans\label{rema-Fg-orientable} cette proposition, si $\Xg$ est en plus une variété topologique orientable et de dimension $d_{\Xg}$, l'ouvert $\Fm \dans\Xg^{m}$ (de dimension $md_{\Xg}$) l'est également et la dualité de Poincaré (\ref{poincare-espaces}-(\ref{poincare-espaces-b})) s'applique pour donner l'égalité analogue à \ref{pol-poincare}-(\ref{pol-poincare-a}):
$$\mathalign{
\Pr(\Fm (\Xg))(T)=\smashtop{\prodnl_{i=0}^{\mmo}}\Pr(\Xg\mmoins i)(T)&=&\smashtop{\prodnl_{i=0}^{\mmo}}\P(\Xg)(T)+i\cdot T^{d_{\Xg}{-}1}\\\noalign{\kern4pt}&=&
T^{(d_{\Xg}{-}1)m}\Big({\Pr(\Xg)(T)\over T^{d_{\Xg}{-}1}}\Big)\osp m\,.\hfill
}$$

\comment
Dans ces cas, si, de plus, $\Xg$ est simplement connexe et $d\not=2$, on a $\Hr^{1}(\Fm (\Xg))=\Ext^{1}(\fs k_{\Fm },\fs k_{\Fm })=0$, et toute extension de faisceaux entre faisceaux localement constants sur $\Fm $ sera scindée. Cette remarque est en liaison avec le chapitre \ref{Leray} sur les suites spectrales de Leray,  notamment l'assertion (\ref{autre-scindage-c}) du théorème \ref{autre-scindage}.\endcomment
\end{rema}

\begin{rema}\label{PcUFm}Soient $1\leq a\leq m\in\NN$ et $\pi_a:\Fm\to\Fa$ la projection sur les $a$ dernières coordonnées
\comment
(\ref{fibrations-pi})
\endcomment
. Pour tout ouvert $\Ug\dans\Fg_a$, notons $\Ug\Fm :=\pia ^{-1}(\Ug)$. La même preuve de \ref{theo-scindage}-(\ref{theo-scindage-a}) (\cf rem. \ref{remas-theo-scindage}-(\ref{remas-theo-scindage-a}))
montre que la suite 
$$
0\to
\Hc(\Ug\Fm )[-1]^{m}\to
\Hc(\Ug\Fmm )\to
\Hc(\Xg\times\Ug\Fm )\to0\,,
\eqno(\diamond)$$
est exacte. On en déduit, comme pour \ref{pol-poincare}-(\ref{pol-poincare-a}), que l'on a 
$$
\Pc(\Ug\Fba )=\Pc(\Ug\Fg_{b-1+a})\cdot(\Pc(\Xg)+(b-1+a)T)$$
et par itération (\cf\ref{PcvUFm}):
$$
\Pc(\Ug\Fba (\Xg))=\Pc(\Ug)\cdot\Pc(\Fg_{b}(\Xg\mmoins\aa))\,.
$$

Il est intéressant d'observer qu'une condition nécessaire pour l'exactitude des suites $(\diamond)$ pour tout $\Ug\dans\Xg$, est que $\Xg$ soit totalement $\cup$-acyclique. En effet, dans un tel cas et si $m=1$, le morphisme  $\Hc(\Xg)\otimes\Hc(\Ug)\to\Hc(\Ug)$ doit être nul. On ignore si la totale $\cup$-acyclicité  suffit à l'exactitude des suites~$(\diamond)$ lorsque, par exemple, $\Xg$ est une variété topologique non orientable.
\end{rema}
\begin{rema}
Les\label{trivialite-cohomologique} fibres des applications $\pi_a:\Fba (\Xg)\to\Fa(\Xg)$ sont de la forme $\Fg_{b}(\Xg\mmoins a)$ et, d'autre part, ce qui précède montre que l'on a
$$\Pc(\pia ^{-1}\Cg)=\Pc(\Cg)\cdot\Pc(\Fg_{b}(\Xg\mmoins a))\,,\eqno(\ast)$$
pour toute composante connexe $\Cg$ de $\Fa$. Cette observation qui suggère l'existence d'une forme de trivialité cohomologique pour les applications $\pi_a$, est à l'origine de l'étude de la section \ref{degenerescence} destinée à prouver la dégénérescence des suites spectrales de Leray associées à $\pi_a$ lorsque $\Xg$ est $i$-acyclique et localement connexe (mais pas forcément de type fini). Plus précisément, on y montre que si $\Cg$ est une composante connexe de $\Fa(\Xg)$ et si $\cl x\in\Cg$, le terme $\IE_{2}$ de la suite spectrale de Leray pour la cohomologie à support compact de la fibration $\pi_a:\pia ^{-1}\Cg\to\Cg$ est
$$\IE_{2}=\Hc(\Cg)\otimes\Hc(\pia ^{-1}\cl x)\Rightarrow\Hc(\pia ^{-1}\Cg)\,.$$
On observera déjà à ce niveau que lorsque $\Xg$ est de type fini,
 l'égalité $(\ast)$ force l'annulation des différentielles de la suite spectrale en question $(\IE_{r},d_r)$ pour $r\geq2$, car autrement $\dim(\IE_\infty)<\dim(\IE_2)$, ce qui n'est pas le cas.
\end{rema}

\begin{coro}Soit\label{=l-m-poly-univ} $\Xg$ un espace $i$-acyclique. Le polynôme de Poincaré $\Pc(\Delta_{\ell}\Xg^{m})$ s'obtient en évaluant en $P:=\Pc(\Xg)$ le
polynôme homogène $\Qg_{\ell}^{m}(P,T)\in\ZZ[P,T]$, de degré $\ell$ :
$$\halfdisplayskips
\Qg_{\ell}^{m}(P,T)=T^{\ell}\cdot \mParties m\ell\cdot \smash{\Big({P\over T}\Big)\osp \ell}\,.\postskip-4pt
$$
\end{coro}
\demo Clair d'après \ref{connexes} et \ref{pol-poincare}-(\ref{pol-poincare-a}).
\enddemo

\def\Remarque{Commentaire}\begin{rema}Notre approche des polynômes de Poincaré de $\Fm(\Xg)$ pour variété topologique orientable $\Xg$ repose sur l'annulation du morphisme de restriction $\delta^{*}_{\Xg}:\Hc(\Xg\times\Xg)\to\Hc(\Xg)$ (\ref{prop-acycliques}-(\ref{prop-acycliques-b})). Cela restreint considérablement la portée de la méthode, mais n'exclut pas qu'elle puisse être appliquée à d'autres cas, et même pour $\Xg$ compact. 
Par exemple, lorsque $\Xg$ est un groupe de Lie compact connexe $\Kg$ de dimension $d_{\Kg}$, il est avantageux de profiter de l'action diagonale libre de $\Kg$ sur $\Fm (\Xg)$. Notons $\oo\Kg=\Kg\mmoins\set e/$. On a la bijection
$$\let\Gg\Kg
\Psi:\Fg_{\mmo}(\oo\Gg)\times\Gg\to\Fg(n,\Gg)\,,\quad(\cl x,g)\mapsto(x_1 g,\ldots,x_{n-1}g,g)\,,$$
et ceci donne aussitôt l'égalité
$\Hc(\Fm (\Kg))=\Hc(\Fg_{\mmo}(\oo\Kg))\otimes\Hc(\Kg)$,
où $\oo\Kg$ n'est plus compact. On peut alors se demander si
$\oo\Kg$ est \uacyclique.

Comme $\Kg$ est compact, le prolongement par zéro $\Hc(\oo\Kg)\to H(\Kg)$ identifie $\Hc(\oo\Kg)$ à $H^{+}(\Kg)$ et nous avons le diagramme commutatif
$$\xymatrix@C=1.cm@R=5mm{
\Hc(\oo\Kg)\otimes\Hc(\oo\Kg)\ar[r]|(0.61){\,\rho_{\oo\Kg}\,}\arinto[d]&\Hc(\Delta_{\oo\Kg})\arinto[d]\\
H(\Kg)\otimes H(\Kg)\ar[r]|(0.61){\,\rho_{\Kg}\,}& H(\Delta_\Kg)\\
}$$
où $\rho_{\Kg}$ est le cap-produit. 
Si $\rho_{\oo\Kg}=0$, nous devons avoir
$$H^{+}(\Kg)\wedge H^{+}(\Kg)=0\,,$$
mais ceci n'est possible, par dualité de Poincaré, que si le groupe de Lie $\Kg$ à la même homologie qu'une sphère. Lorsque le corps des coefficients est $\QQ$, cela arrive seulement dans trois cas $\Kg=\cSS^0$, $\Kg=\cSS^1$ et $\Kg=\cSS^{3}\simeq\mathop{\rm SU}(2)=\mathop{\rm Sp}(1)$.
Dans ces cas la proposition \ref{pol-poincare}-(\ref{pol-poincare-a}) s'applique et on trouve:
$$
\mathrigid2mu
\mathalign{
\Pc(\Fm (\Kg))=\Pc(\Kg)\cdot\Pc(\Fg_{\mmo}(\oo\Kg))
&=&\Pc(\Kg)\cdot T^{\mmo}\cdot \smashtop{\big({\Pc(\oo\Kg)/ T}\big){}\osp{\mmo}}\hfill\\\noalign{\kern4pt}
&=&(T^{d_\Kg}+1)\cdot T^{\mmo}\cdot \big(T^{d_{\Kg}-1}\big){}\osp{\mmo}\,.
\hfill}\postskip0pt$$
\end{rema}

\subsection{Polynômes de Poincaré de $\Delta_{\leq\ell}\Xg^{m}$}Une application immédiate de l'exactitude du complexe fondamental est la détermination du polynôme de Poincaré $\Pc(\Delta_{\leq \ell}\Xg^{m})$ lorsque $\Xg$  est de type fini. L'assertion suivante généralise la proposition \ref{pol-poincare}.

\begin{coro}\label{pol-poin-gen}Soit $\Xg$ un espace $i$-acyclique.
Le polynôme de Poincaré $\Pc(\Delta_{\leq \ell}\Xg^{m})$ est le polynôme homogène de $\ZZ[\Pc(\Xg),T]$, de degré $\ell$, donné par la somme alternée {\rm(\cf\ref{connexes})}
$$\halfdisplayskips
\Pc(\Delta_{\leq\ell}\Xg^{m})=\sumnl_{0\leq a<\ell}(-1)^{a}\cdot\mParties m{\ell-a}\cdot\Pc(\Fg_{\ell-a}(\Xg))\cdot T^{a}\,,
$$
avec comme termes de plus bas et plus haut degrés en $T$:
$$\setbox111=\hbox{$\hbox{\mParties m\ell}$}
\setbox110=\hbox{$\hbox{\mParties m{\ell-a}}$}
\mathalign{\Pc(\Delta_{\leq\ell}\Xg^{m})&=&\box111\cdot\Pc(\Xg)^{\ell}\cdot T^{0}+\cdots\hfill\\\noalign{\kern4pt}
&&\hskip-5mm{}\cdots+\sumnl_{0\leq a<\ell}(-1)^{a}\cdot\box110\cdot(\ell-a-1)!\cdot \Pc(\Xg)\cdot T^{\ellmo}\,.}\postskip0pt$$
\end{coro}
\demo Résulte d'appliquer le théorème \ref{theo-complexe-exact}, l'explicitation du corollaire \ref{=l-m-poly-univ} et le fait que le terme de plus haut degré en $T$ de $\Pc(\Fg_{\ell-a}(\Xg))$ est $(\ell-a-1)!\cdot\Pc(\Xg)\cdot T^{\ell-a-1}$ d'après \ref{pol-poincare}-(\ref{pol-poincare-a}).
\enddemo

\subsection{Polynôme universel pour $\Delta_{\leq\ell}\Xg^{m}$}
Le\label{<=l-m-poly-univ} corollaire \ref{pol-poin-gen} montre que le polynôme de Poincaré de $\Hc(\Delta_{\leq\ell}\Xg^{m})$ s'obtient en évaluant
un certain polynôme homogène de l'anneau $\ZZ[P,T]$, degré $\ell$, en $P=\Pc(\Xg)$.
Ce polynôme est indépendant de $\Xg$ et il est unique. 
\begin{defi*}
Le \expression{polynôme universel pour la cohomologie à support compact des espaces $\Delta_{\leq\ell}^{m}$}, noté $\Qg_{\leq\ell}^{m}(P,T)$\glossary{$\Qg_{\leq\ell}^{m}\in\ZZ[P,T]$:polynôme universel pour $\Pc(\Delta_{\leq\ell}\Xg^{m})$}, est le polynôme de $\ZZ[P,T]$, homogène de degré total $\ell$, donné par
$$\Qg_{\leq\ell}^{m}(P,T):=
T^{\ell}
\sumnl_{a=1,\ldots,\ell}(-1)^{\ell-a}\cdot
\mParties ma
\cdot
\smash{\Big({P\over T}\Big)\osp a}
\,,
$$
où 
\goodsmash{0.8}{0.5}{$\displaystyle\preskip-1em\mParties ma:={1\over \ell! }\ \sumnl_{j=0}^{\ell}(-1)^{\ell-j}\,\binome{\ell}{j}\;j^{m}\,.$ (\Cf\ref{connexes}.)}
\end{defi*}
\begin{rema}Il peut y avoir beaucoup de simplifications dans l'expression donnant ces polynômes. Par exemple, sachant que $\Delta_{\leq m}\Xg^{m}=\Xg^{m}$, on peut anticiper l'égalité
$\Qg_{\leq m}^{m}(P,T)=P^{m}\,,$
ce qui est loin d'être une évidence uniquement à partir de la définition de $\Qg_{\leq \ell}^{m}$.
\end{rema}

\subsubsection{Un exemple de polynômes universels}
Voici les six polynômes universels pour $m=6$.

\begingroup
\def\S(#1){\Qg_{\leq#1}(P,T)}
\def\PX{&P}\def\*{\;}\def\T{\hfill&\,T}
\let\,\relax
\def\\{\cr}
\halign{\hskip\parindent\quad${}#{}$&&${}#{}$\cr
\S(1, 6) &=& \hfill\PX\T^{0}\\
\S(2, 6) &=& 31\*\PX^2 \T^{0}&+&\hfill 30\*\PX\T\\
\S(3, 6) &=& 90\*\PX^3 \T^{0}&+&\hfill 239\*\PX^2\T&+&150\*\PX\T^2\\
\S(4, 6) &=& 65\*\PX^4 \T^{0}&+&\hfill 300\*\PX^3\T&+&476\*\PX^2\T^2&+&240\*\PX\T^3\\
\S(5, 6) &=& 15\*\PX^5 \T^{0}&+&\hfill 85\*\PX^4\T&+&225\*\PX^3\T^2&+&274\*\PX^2\T^3&+&120\*\PX\T^4\\
\S(6, 6) &=&\hfill \PX^6\T^{0}&.&\\
}
\endgroup

\section{Représentations du groupe symétrique}
\glossarytitle{Représentations du groupe symétrique}
Dans\label{representations} \ref{complexe-fondamental} nous avons muni les espaces $\Delta_{?\ell}\Xg^{m}$ 
de l'action de $\Sm $. Nous allons maintenant explorer en détail la structure de $\S _{m}$-espace de $\Delta_{?\ell}\Xg^{m}$. 
\comment
Mais avant cela, il nous faut plus de notations.
\endcomment

\subsection{Notations pour le décompositions et les diagrammes de Young}\label{nota-young}
\def\varitemizeseps{\topsep4pt\itemsep6pt\mou\parsep6pt\mou}
\begin{itemize}
\mynobreak\nobreak\item Une \expression{décomposition}\index{decomposition@décomposition} $\lambda$ d'un entier $m\geq 0$, notée $\lambda\vdash m$\glossary{$\lambda\vdash m$:\emph{alias} pour: $\lambda$ est une décomposition de $m$}, est la donnée d'une suite décroissante d'entiers positifs $\lambda:=(\lambda_1\geq\lambda_2\geq\cdots\geq\lambda_\ell>0)$ telle que $m=\lambda_1+\lambda_2+\cdots+\lambda_\ell$. On note  $|\lambda|=m$\glossary{${|\lambda|}$:nombre décomposé par $\lambda$, donc $|\lambda|:=\sum_i\lambda_i$} et $\ell(\lambda):=\ell$.\glossary{${\ell(\lambda)}$:nombre de termes (non nuls) d'une décomposition}

\item \'Etant données une famille d'entiers naturels $\set\XX_1,\ldots,\XX_m/$, telle que $m=\sum_i \XX_i\, i$, on notera de manière équivalente
$$\let\underbrace\underparenthesis
(
\mathop{\smashbot{\underbrace{m,\ldots,m}}{\vrule depth7pt width0pt}}\limits_{\XX_m},
\mathop{\smashbot{\underbrace{\mmo,\ldots,\mmo}}{\vrule depth7pt width0pt}}\limits_{\XX_\mmo},
\ldots
\mathop{\smashbot{\underbrace{1,\ldots,1}}{\vrule depth7pt width0pt}}\limits_{\XX_1}
)=
(1^{\XX_1},2^{\XX_2},\ldots,m^{\XX_m})$$

\item\leavevmode \vtop{\advance\hsize-\leftmargin\parindent0pt\addhabille-1\Habillage{\hbox{\vrule\ $
\def\m#1{\rlap{\kern2cm$(\lambda_{#1})$}}
\tableauc{}{
\m1&&&&&\\
\m2&&&\\
\m3&&&\\
\m4
}$\kern0.8cm}}{1}{-4pt}
Un \expression{diagramme de Young}\index{Young!diagramme}\index{diagramme de Young} est un empilement vertical de juxtapositions horizontales de boites dont le nombre, $\lambda_i$, décroît. 
La suite $\lambda:=(\lambda_1,\ldots,\lambda_\ell)$ est une décomposition de $m:=|\lambda|$. On identifie ainsi décompositions et diagrammes.
On note $\Y_{\ell}(m)$\glossary{${\Y_{\ell}(m):=\set \lambda\vdash m\mid\ell(\lambda)=\ell/}$:} l'ensemble des diagrammes à $m$ boites et à $\ell$ lignes, soit donc
$$\normaldisplay\preskip1.25ex\Y_{\ell}(m):=\bigset (\lambda\vdash m)\ \&\ (\ell(\lambda)=\ell)/\,.$$
\endHabillage}

\item\leavevmode \vtop{\advance\hsize-\leftmargin\parindent0pt\Habillage{\hbox{\vrule\ $\def\everyygcell{\mathrigid0mu \relax}
\tau(\lambda):=\left(\tableauc{\yght=1em\ygwd=1em}{
1&2&3&4&5&6\\
7&8&9&10\\
11&12&13&14\\
15
}
\right)$}}{1}{-3pt}
 Un \expression{tableau de Young standard}\index{tableau de Young}\index{Young!tableau} est un diagramme de Young $\lambda$ dont on rempli les boites par les entiers $1,2,\ldots,|\lambda|$. On note $\tau(\lambda)$\glossary{$\tau(\lambda)$:tableau de Young standard associé à une décomposition $\lambda$} le tableau obtenu par une numérotation successive, de gauche à droite et de haut en bas, tel qu'indiqué dans la figure ci-contre. 
\endHabillage}

\item Si $\lambda\in\Y_{\ell}(m)$, on note $\Plambda $\glossary{$\Plambda =\S_{\lambda_1}\times\cdots\times \S_{\lambda_{\ell}}$: sous-groupe des permutations de $\Sm$ qui conservent les lignes du tableau  $\tau(\lambda_1,\ldots,\lambda_{\ell})$} le sous-groupe de $\Sm$ des permutations qui conservent les lignes de  $\tau(\lambda)$. On a 
$$
\Plambda =\S_{\lambda_1}\times\cdots\times \S_{\lambda_{\ell}}\,.$$
\item Un diagramme $\lambda\in\Y_{\ell}(m)$ décompose l'intervalle $\iii[1,m]$ suivant les $\ell$ sous-intervalles définis par les lignes de $\tau(\lambda)$. On note  $\pgoth_{\lambda}\in\Pgoth_{\ell}(|\lambda|)$ la partition ainsi déterminée.
Le sous-groupe
 $$\S_{\lambda}:=N_{\Sm}\Plambda$$
où $N_{\Sm}(\Plambda )$\glossary{${\Slambda:=N_{\Sm}(\Plambda )}$:normalisateur  de $\Plambda $ dans $\Sm$} désigne le normalisateur de $\Plambda $ dans $\Sm$,
est le sous-groupe des éléments de $\Sm $ qui fixent la partition $\pgoth_{\lambda}$.
\end{itemize}

\begin{lemm}\label{lemme-young}On fait agir $\Sm$ sur $\Pgoth_{\ell}(m)$ par son action naturelle. Pour toute partition $\pgoth\in\Pgoth$, on note $\S_{\pgoth}$ le sous-groupe de $\Sm$ qui fixe $\pgoth$. Chaque orbite de $\Sm$ rencontre une unique partition de la forme $\pgoth_{\lambda}$. On a 
$$\Pgoth_{r}(m)=\coprodnl_{\lambda\in\Y_{\ell}(m)} \Sm\cdot \pgoth_{\lambda}
\simeq\coprodnl_{\lambda\in\Y_{\ell}(m)} \Sm/\Slambda\,.
$$
où $\lambda$ parcourt l'ensemble $\Y_{\ell}(m)$ des décompositions de $m$ en $\ell$ entiers $\not=0$.
\begin{enumerate}
\item Si\label{lemme-young-a} $\lambda=(\lambda_1,\ldots,\lambda_{\ell})=((\lambda_1){}^{\XX_{\lambda_1}},\ldots,(\lambda_\ell){}^{\XX_{\lambda_\ell}})$,  on a \glossary{${\Glambda :=\Slambda/\Plambda =\S_{\XX_1}\times\cdots \times \S_{\XX_r}}$:si $\lambda=(\lambda_1,\ldots,\lambda_{\ell})=(d_1^{\XX_1},\ldots,d_r^{\XX_r})$}
$$\scriptspace0pt
\Slambda=\S_{\pgoth_{\lambda}}=N_{\Sm}(\Plambda )
\text{\quad et\quad}
\Glambda :=\Slambda/\Plambda \simeq\S_{\XX_{\lambda_1}}\times\cdots \times \S_{\XX_{\lambda_\ell}}\,,$$
\item Le\label{lemme-young-b} groupe $\Glambda$ s'identifie au sous-groupe $\S_{\XX_{\lambda_1}}\times\cdots \times \S_{\XX_{\lambda_\ell}}\dans\Slambda $ et l'on a
$$\halfdisplayskips\Slambda=\big(\S_{\XX_{\lambda_1}}\times\cdots \times \S_{\XX_{\lambda_\ell}}\big)\ltimes \Plambda \,.$$
\end{enumerate}\end{lemm}

\subsection{Décomposition ouverte $\Sm$-stable de 
$\Delta_{\ell}\Xg^{m}$}\label{dec-symetrique-Delta}

Par le lemme \ref{lemme-young} la décomposition ouverte de $\Delta_{\ell}^{m}$ de \ref{connexes} s'écrit
$$\Delta_{\ell}\Xg^{m}=\coprod\nolimits_{\lambda\in\Y_{\ell}(m)}\Sm\cdot
\Fg_{\pgoth_{\lambda}}(\Xg)\,.$$
La proposition suivante est alors immédiate.

\begin{prop}\label{dec-Delta-l-m}Soit $\lambda=(\lambda_1,\ldots,\lambda_\ell)=(1^{\XX_1},\ldots,m^{\XX_m})\in\Y_{\ell}(m)$.
\begin{enumerate}\mynobreak\itemsep2pt\displayskips8/10
\item\leavevmode\label{dec-Delta-l-m-a}Le groupe $\Plambda =\S_{\lambda_1}\times\cdots\times \S_{\lambda_\ell}$ est le sous-groupe de $\Sm$ des permutations qui agissent comme l'identité sur $\Fg_{\pgoth_{\lambda}}(\Xg)$. Son normalisateur $\Slambda $ est le sous-groupe de $\Sm$ qui laisse stable $\Fg_{\pgoth_{\lambda}}(\Xg)$. 

\begingroup\medskip\leftskip-0.95\leftmargin
\noindent Soit maintenant la surjection canonique de groupes
$$\nu_{\lambda}:\Slambda \onto 
\Glambda :=\Slambda /\Plambda \simeq\S_{\XX_{\lambda_1}}\times\cdots\times \S_{\XX_{\lambda_\ell}}\dans \S_{\ell}
$$
 et considérons
$\Fg_{\pgoth_{\lambda}}(\Xg)$ muni de sa structure de $\Glambda $-espace.

\endgroup
\item\leavevmode\label{dec-Delta-l-m-b}L'application
$\phi_{\lambda}:\Fg_{\pgoth_{\lambda}}(\Xg)\to\Fg_{\ell}(\Xg)\,,$
 $\phi_{\lambda}:(x_1,\ldots,x_{m})\mapsto (y_1,\ldots,y_{\ell})$ avec
 $y_k:=x_{\lambda_1+\cdots+\lambda_{k}}$,
est un isomorphisme de $\Glambda $-espaces.
\item L'application\label{dec-Delta-l-m-c}
$$\Psi_{\lambda}:
\Sm\times_{\Slambda }\Fg_{\ell}(\Xg)\to
\Sm\cdot\Fg_{\pgoth_{\lambda}}(\Xg)\dans\Xg^{m}\,,
\quad \cl{(\alpha,\cl x)}\mapsto
\alpha\cdot\phi^{-1}(\cl x)\,,
$$
est un isomorphisme de $\Sm$-espaces.
\end{enumerate}
\end{prop}

\subsubsection{Caractères de $\S _{m}$ associés à 
\smashbot{$\Delta_{\leq\ell}\Xg^{m}$}}
Pour\label{char-i} toute partie $\Zg\dans\Xg^m$, de type fini et stable sous l'action de $\Sm$, les espaces $\Hc^{i}(\Zg,k)$ et de $\Hr^{i}(\Zg,k)$\glossary{$\chic(\Zg ;i)\,,\chi(\Zg ;i)$:caractères de $\Sm$-module de $\Hc^{i}(\Zg ,k)$ et de $\Hr^{i}(\Zg ,k)$} sont les $\S _{m}$-modules de dimension finie, leurs les caractères seront  respectivement notés: 
$$\begin{casesalign}
\chic(\Zg;i)&:\S _{m}\to k\,,&\quad&\chic(\Zg;i)(\alpha)&:=&\tr(\alpha\actson\Hc^{i}(\Zg;k))\hfill\\
\hfill\chi(\Zg;i)&:\S _{m}\to k\,,&\quad&\chi(\Zg;i)(\alpha)&:=&\tr(\alpha\actson\Hr^{i}(\Zg;k))\,. 
\end{casesalign}
$$

\comment
L'action de $\S _{m}$ sur $\Delta_{?\ell}\Xg^{m}$ induit
les espaces de cohomologie $\Hr^{i}(\Delta_{?\ell}\Xg^{m},k)$ et $\Hc^{i}(\Delta_{?\ell}\Xg^{m},k))$ de structures de $\Sm$-modules. Nous nous intéresserons maintenant à leurs caractères.
\endcomment

\def\Remarque{Commentaire}\begin{rema}Si\label{singularites} $\Xg$ est une variété topologique orientable de dimension $d_{\Xg}$, l'ouvert $\Fm \dans\Xg^{m}$, de dimension $md_{\Xg}$, est 
une variété topologique orientable, $\Hr^{i}(\Fm )$ et $\Hc^{md_{\Xg}-i}(\Fm )$ sont des $\Sm$-modules en
dualité et $\chic(\Fm ,*)$ détermine $\chi(\Fm ,*)$ (\cf\ref{def-action-image-inverse}-(\ref{def-action-image-inverse-b})). En dehors de ce cas, nos méthodes ne s'appliquent pas à l'étude de la cohomologie ordinaire $\Hr(\Delta_{?\ell}\Xg^{m})$, ni en tant qu'espace vectoriel ni, à fortiori, en tant que $\Sm$-module. Elles s'appliqueront par contre à la cohomologie de Borel-Moore $\Hbm(\Delta_{?\ell}\Xg^{m})$ (\ref{coh-BM}).
\end{rema}

\def\Remarque{\miseengarde Avertissement}
\begin{rema*}{\miseengarde Dans toutes les sections concernant les représentations des groupes symétriques on suppose
$\car(k)=0$~(\cf\ref{car(k)=0}).
\comment\footnote{\color{red}When $k$ is a field of characteristic $0$, every $\S _{n}$-representation decomposes as a direct sum of irreducible representations. In characteristic $0$ every irreducible representation of $\S _{n}$ is defined over $\QQ$. As a result this decomposition does not depend on the field $k$.})\endcomment
\killline}
\end{rema*}

\subsection{Foncteurs d'induction $\Ind_{\Glambda }^{\S _{m}}$ et $\Ibf ^{m}_{\ell}$}Pour $0<\ell\leq m$,\label{operateurs-inductions}
 et tout $\lambda\in\Y_{\ell}(m)$, on note
$$\displayboxit{\Ind_{\Glambda }^{\S _{|\lambda|}}:\Mod(k[\S_{\ell}])\fonct\Mod(k[\S_{|\lambda|}])}\eqno(\Ind)$$
le foncteur \smashtop{$\ind_{\Slambda }^{\Sm}\circ\Res^{\S_{\ell}}_{\Glambda }$}\glossary{${\Ind_{\Glambda }^{\S _{|\lambda|}}:=\ind^{\S_{|\lambda|}}_{\Slambda }\circ\Res^{\S_{\ell}}_{\Glambda }}$:où $\Slambda $ agit à travers de la surjection $\nu_{\lambda}:\Slambda \onto \Glambda $}
où $\Slambda $ agit à travers de la surjection $\nu_{\lambda}:\Slambda \onto \Glambda $ de \ref{dec-Delta-l-m}-(\ref{dec-Delta-l-m-a}).
On considère ensuite le foncteur 
$$\displayboxit{
\Ibf ^{m}_{\ell}:=\sumnl_{\lambda\in\Y_{\ell}(m)}\Ind^{\S _{m}}_{\Glambda }:\Mod(k[\S_{\ell}])\fonct\Mod(k[\Sm ])
}\eqno(\Ibf )$$ 
(Remarquer que l'on a $\Ibf ^{m}_{m}=\id$.)
On notera par la même notation l'opérateur linéaire d'induction
$\Ibf ^{m}_{\ell}:\CCc[\S _{\ell}]\to \CCc[\S _{m}]$\glossary{{$\Ibf ^{m}_{\ell}:\CCc[\S _{\ell}]\to \CCc[\S _{m}]$}:opérateur sur les fonctions centrales} défini sur les fonctions centrales.


\begin{prop}Pour\label{coro-car-l-m}
 $0<\ell\leq m$, on a un isomorphisme de $\Sm$ espaces:
\displayskips5/10
$$\Delta_{\ell}\Xg^{m}\cong\coprod\nolimits_{\lambda\in\Y_{\ell}(m)}\Sm\times_{\Slambda }\Fg_{\ell}(\Xg)\,.\postskip0pt$$
En particulier, on a 
$$
\chic(\Delta_{\ell}\Xg^{m};i)=
\Ibf ^{m}_{\ell}\big(\chic(\Fg_{\ell}(\Xg);i)\big)\text{\quad et\quad }
\chi(\Delta_{\ell}\Xg^{m};i)=
\Ibf ^{m}_{\ell}\big(\chi(\Fg_{\ell}(\Xg);i)\big)\,.
 $$
\end{prop}
\demo Corollaire immédiat de \ref{lemme-young} et \ref{dec-Delta-l-m}-(\ref{dec-Delta-l-m-c}).
\enddemo

\subsubsectionline{Présentation de $\S _{m}$-module de $\Hc^{i}(\Delta_{\leq\ell}\Xg^{m})$.}Nous pouvons à présent rassembler les résultats précédents pour donner une présentation de la représentation de $\Sm$ sur $\Hc^{i}(\Delta_{\leq\ell}\Xg^{m})$ lorsque $\Xg$ est $i$-acyclique.

\begin{theo}\label{theo-caracteres}Soit $\Xg$ un espace $i$-acyclique tel que $\dim\Hc(\Xg)<\infty$. 
\begin{enumerate}\displayskips5/10
\mynobreak\nobreak
\item\leavevmode\label{theo-caracteres-a}Le caractère du $\Sm$-module $\Hc^{i}(\Fm (\Xg))$ vérifie 
$$\chic(\Fm(\Xg) ;i)=\chic(\Xg^{m};i)+\chic(\Delta_{\leq \mmo}\Xg^{m};i-1)\,.$$
\item Pour\label{theo-caracteres-b}
$0<\ell\leq m$, le caractère du $\Sm$-module $\Hc^{i}(\Delta_{\leq\ell}\Xg^{m})$ vérifie 
$$\def\tt{\vrule height7pt width0pt }
\chic(\Delta_{\leq\ell}\Xg^{m};i)=
\sumnl_{\tt0\leq a<\ell}\,
(-1)^{a}\,
\Ibf ^{m}_{\ell-a}\big(\chic(\Fg_{\ell-a}(\Xg);i-a)\big)\,.
\postskip0pt$$
\end{enumerate}
\end{theo}
\demo
Par \ref{theo-complexe-exact}, la suite de $\Sm$-modules
$$\mathrigid1.5mu
0\to\Hc^{i-\ell+1}(\Delta_{1}^{m})\to\cdots\to
\Hc^{i-1}(\Delta_{\ellmo}^{m})\to
\Hc^{i}(\Delta_{\ell}^{m})\to
\Hc^{i}(\Delta_{\leq\ell}^{m})\to0
$$
est exacte. On a donc l'égalité
$\chic(\Delta_{\leq\ell}^{m};i)=\sumnl_{0\leq a\leq\ell}(-1)^{a}\,\chic(\Delta_{\ell-a}^{m};i-a)\postdisplaypenalty10000$
à laquelle on applique la proposition \ref{coro-car-l-m}.
\enddemo

\subsubsectionline{Expression de $\chic(\Delta_{?\ell}\Xg^{m};*)$ en termes de $\chic(\Xg^{\ell};*)$.}Le\label{relation-caracteres}\label{algo-caracteres} théorème \ref{theo-caracteres} est la base d'un algorithme de calcul pour $\chic(\Delta_{\leq\ell}\Xg^{m};i)$. En effet, l'égalité (\ref{theo-caracteres-b}) l'exprime comme combinaison des induits des $\chic(\Fg_{\ell'}(\Xg),i')$ pour $\ell'\leq\ell$ et $i'\leq i$, et l'égalité (\ref{theo-caracteres-a}) exprime chaque $\chic(\Fg_{\ell'} (\Xg),i')$ comme la somme de
$\chic(\Xg^{\ell'},i')$, caractère calculé par Macdonald \cite{mac} (\cf aussi \ref{theo-trace-gen-Xm}), et de
$\chic(\Delta_{\leq \ell'-1}\Xg^{\ell'},i')$, avec (donc) $\ell'-1<\ell$.

La section suivante précisera davantage cette idée.

\subsection{Opérateurs d'inductions itérées $\Ibf(\sigma)$ et $\Thetabf^{m}_{\ell}$}

Pour\label{inductions-iterees} toute suite  d'entiers positifs 
$\sigma=(m_0>m_1>\cdots>m_{t-1}>m_t)$
strictement décroissante, on pose $|\sigma|:=t$ et\glossary{${\Ibf (\sigma):=\Ibf ^{m_0}_{m_1}\circ\cdots\circ\Ibf ^{m_{t-1}}_{m_t}
:\CCc[\S _{m_t}]\to \CCc[\S _{m_0}]}$:opérateur d'inductions itérées}
$$\begin{cases}
\hbox{si $|\sigma|=0$,}&\Ibf ((m_0)):=\id:\CCc[\S _{m_0}]\to\CCc[\S _{m_0}]\\[4pt]
\hbox{si $|\sigma|>0$,}&\Ibf (\sigma):=\Ibf ^{m_0}_{m_1}\circ\cdots\circ\Ibf ^{m_{t-1}}_{m_t}
:\CCc[\S _{m_t}]\to \CCc[\S _{m_0}]
\end{cases}\eqno(\Ibf (\sigma))$$
et l'on définit pour $\ell\leq m$, l'opérateur d'inductions itérées\glossary{${\Thetabf^{m}_{\ell}:=(-1)^{m-\ell}\sumnl_{\sigma:m\mathbin{\sssearrow} \ell}(-1)^{|\sigma|-1}\Ibf (\sigma):\CCc[\S _{\ell}]\to \CCc[\S _{m}]}$:opérateur d'inductions itérées}:
$$\displayboxit{\Thetabf^{m}_{\ell}:=(-1)^{m-\ell}\sumnl_{\sigma:m\mathbin{\sssearrow} \ell}(-1)^{|\sigma|-1}\Ibf (\sigma):\CCc[\S _{\ell}]\to \CCc[\S _{m}]}\eqno(\Thetabf)$$
où la sommation est indexée par l'ensemble des suites strictement décroissantes qui partent de $m$ et aboutissent à $\ell$.  

On remarquera que l'on a $\Thetabf^{m}_{m}=\id$. 

\noendpoint\begin{theo}\label{theo-caracteres-devissage}Soit $\Xg$ un espace $i$-acyclique tel que $\dim\Hc(\Xg)<\infty$. 
\begin{enumerate}
\mynobreak\nobreak\item\leavevmode\label{theo-caracteres-devissage-a}Pour tous $0<\ell\leq m$ et tout $i\in\NN$, on a
$$\mathalign{
\ \llap{\rm i)\ }\hfill\chic(\Fm(\Xg);i)&=&\sumnl_{0\leq a<m}\Thetabf^{m}_{m-a}\big(\chic(\Xg^{m-a};i-a)\big)\hfill\\\noalign{\kern4pt}
\ \llap{\rm ii)\ }\hfill\chic(\Delta_{\ell}\Xg^{m};i)&=&
\Ibf ^{m}_{\ell}\Big(\sumnl_{0\leq a<\ell}\Thetabf^{\ell}_{\ell-a}\big(\chic(\Xg^{\ell-a};i-a)\big)\Big)\hfill\\\noalign{\kern4pt}
\mathrigid0mu
\ \llap{\rm iii)\ }\chic(\Delta_{\leq\ell}\Xg^{m};i)&=&
\mathrigid0mu
\sum_{0\leq b<\ell}(-1)^{b}\ 
\Ibf ^{m}_{\ell-b}\Big(\sumnl_{\vrule height18pt width0pt\hskip-0.5cm 0\leq a<\ell-b}\hskip-1.7em\Thetabf^{\ell-b}_{\ell-b-a}\big(\chic(\Xg^{\ell-b-a};i-b-a)\big)\Big)}
$$
\item\leavevmode\label{theo-caracteres-devissage-b}\label{theo-caracteres-Poincare} Soient $\Xg$ et $\Yg$ des espaces $i$-acycliques à cohomologies à support compact de dimensions finies. On a 
$$
\let\Big\big
\Big(\chic(\Delta_{?\ell}\Xg^{m},*)=\chic(\Delta_{?\ell}\Yg^{m},*)\Big)\ \Longleftrightarrow\ 
\Big(\Pc(\Xg)=\Pc(\Yg)\Big)\,.$$
\end{enumerate}
\end{theo}
\demo
(\ref{theo-caracteres-devissage-a}-i) 
On procède par récurrence sur $m$. Lorsque $m=1$, on a 
$\Fg_{1}(\Xg)=\Xg$, $a=0$, la somme est réduite à un seul terme et l'égalité est immédiate. Dans le cas général, l'égalité du théorème \ref{theo-caracteres}-(\ref{theo-caracteres-b}) pour $\ell=m$ donne l'égalité:
$$\chic(\Fm;i)=\chic(\Xg^{m};i)-
\hskip-0.5cm
\sum_{m>m-a>0}
\hskip-1em
(-1)^{a}\ \Ibf ^{m}_{m-a}\big(\chic(\Fg_{m-a};i-a)\big)
$$
où l'on peut remplacer, par hypothèse inductive,
$$\binspace0
\mathalign{
\chic(\Fg_{m-a};i-a)&=&\chic(\Xg^{m-a};i-a)+{}
\hskip-0.5cm\sum_{m-a>m-a-b>0}\hskip-0.5cm\Thetabf^{m-a}_{m-a-b}\big(\chic(\Xg^{m-a-b};i-a-b)\big)\,.}
$$
On exprime ainsi $\chic(\Fm;i)$ comme somme de deux termes. 
$$\begin{casesalign}
\displaystyle
A&:=&\chic(\Xg^{m};i)\ +\hskip-1ex\sum_{m>m-a>0}\hskip-1em(-1)^{a+1}\,\Ibf ^{m}_{m-a}\big(\chic(\Xg^{m-a};i-a)\big)\hfill\\\noalign{\kern4pt}
B&:=&\sumnl_{\vrule height18ptwidth0pt\hskip-0.7cm m>m-a>0}\hskip-2.em(-1)^{a+1}\,\Ibf ^{m}_{m-a}
\hskip-2em \sum_{\vrule height8ptwidth0pt
m-a>m-a-b>0}\hskip-1.5em\Thetabf^{m-a}_{m-a-b}\big(\chic(\Xg^{m-a-b};i-a-b)\big)
\end{casesalign}$$
Le terme $A$ contient  $\Thetabf^{m}_{m}=\Ibf ((m))$ et les opérateurs $(-1)^{a+1}\,\Ibf ((m,m-a))$ qui interviennent dans le développement de $\Thetabf^{m}_{m-a}$ pour $a>0$.  Le terme $B$, quant à lui, contient exactement tous les opérateurs qui manquent encore pour reconstruire le second membre de (\ref{theo-caracteres-devissage-a}-i). En effet,
$$\mathalign{\sumnl_{\vrule height18ptwidth0pt\hskip-0.5cm m>m-a>0}\hskip-2.5em(-1)^{a+1}\,\Ibf ^{m}_{m-a}
\hskip-2em \sum_{m-a>m-a-b>0}\hskip-2em\Thetabf^{m-a}_{m-a-b}&=&
\sumnl_{\vrule height18ptwidth0pt\hskip-0.5cm 
a>0{\rm\ \&\ }\sigma:m-a\mathbin{\sssearrow} m-a-b>0}\hskip-3.2cm(-1)^{a+1}\,\Ibf ^{m}_{m-a}
 (-1)^{b} (-1)^{|\sigma|-1}\Ibf (\sigma)\hfill\\\noalign{\kern8pt}
 &=&
(-1)^{a+b}\ 
\sumnl_{\vrule height18ptwidth0pt\hskip-2cm 
\sigma:m\mathbin{\sssearrow} m-(a+b){\rm\ \&\ }|\sigma|\geq2}\hskip-1.2cm(-1)^{|\sigma|}\Ibf (\sigma)\,.
\hfill}
$$
Ceci termine la preuve de l'égalité (\ref{theo-caracteres-devissage-a}-i). La formule (\ref{theo-caracteres-devissage-a}-ii) en découle aussitôt puisque
$
\chic(\Delta_{\ell}\Xg^{m};i)=
\Ibf ^{m}_{\ell}\big(\chic(\Fg_{\ell}(\Xg);i)\big)
$ (\ref{coro-car-l-m}), et
 (\ref{theo-caracteres-devissage-a}-iii) en résulte par application directe de \ref{theo-caracteres}-(\ref{theo-caracteres-b}).

(\ref{theo-caracteres-devissage-b}) L'implication $\Rightarrow$ est immédiate en prenant $m=\ell=1$. Pour la réciproque, il suffit, grâce à (\ref{theo-caracteres-a}), de montrer que les caractères $\chic(\Xg^{m};i)$ sont déterminés par $\Pc(\Xg)$. Or, ceci est clair d'après le travail de Macdonald \cite{mac}   (eq.~4.5) (\cf aussi \ref{theo-trace-gen-Xm} pour la formule explicite).
\enddemo

\section{Cohomologie des espaces de configuration, cas général}\label{cohomologie-Borel-Moore}
\glossarytitle{Cohomologie des espaces de configuration, cas général}

\subsectionline{Cohomologie de Borel-Moore.}Comme\label{coh-BM} nous l'avons déjà indiqué dans \ref{singularites}, l'utilisation des complexes fondamentaux limite la portée de nos méthodes sur au moins deux aspects: les espaces $\Xg$ sont $i$-acycliques et la cohomologie est à support compact. Ce sont des limitations assez contraignantes en particulier lors de l'étude du comportement asymptotique de la cohomologie de $\Fm:=\Fm(\Xg)$ suivant les tours de projections
$$\halfdisplayskips
\def\tt{\vrule depth0pt width0pt}
\def\rien{}\def\T#1#2{\def\donne{#2}\Fg_{#1}
\ifx\donne\rien
\def\nnext{\ar@{<<-}[r]}
\else 
\def\nnext{\ar@{<<-}[r]^(#2){\ \tt p_{#1}}}\fi
\nnext}
\xymatrix@C=12mm{
\cdots\ar@{<<-}[r]&
\T {\mmo}{0.58}&
\T m{0.5}&
\T{\mm}{}&\cdots
}$$
où $p_m:\Fmm\onto\Fm$\glossary{${p_m:\Fmm\onto\Fm}$:projection sur les $m$ premières coordonnées} désigne la projection sur les $m$ premières coordonnées, autrement dit, pour l'étude des tours de représentations des groupes symétriques: les \expression{$\FI$-modules} (\cf\ref{stabilite-des-familles-i-acycliques})
$$\halfdisplayskips
\def\tt{\vrule depth0pt width0pt}
\def\rien{}\def\T#1#2{\def\donne{#2}
\Hr(\Fg_{#1})
\ifx\donne\rien
\def\nnext{\ar[r]}
\else 
\def\nnext{\ar[r]^(#2){\tt p_{#1}^{*}\ }}\fi
\nnext}
\xymatrix@C=12mm{
\cdots\ar[r]&
\T {\mmo}{0.53}&
\T m{0.46}&
\T{\mm}{}&\cdots
}\eqno(\ast)$$

\medskip
Lorsque $\Xg$ est une variété topologique orientée la dualité de Poincaré établit un isomorphisme canonique entre la suite $(\ast)$ et la suite
$$\preskip1ex
\def\tt{\vrule depth3pt width0pt}
\def\rien{}\def\T#1#2{\def\donne{#2}
\Hc(\Fg_{#1})\dual
\ifx\donne\rien
\def\nnext{\ar[r]}
\else 
\def\nnext{\ar[r]^(#2){\tt p_{#1!}{}\dual\ }}\fi
\nnext}
\xymatrix@C=10mm{
\cdots\ar[r]&
\T {\mmo}{0.57}&
\T m{0.5}&
\T{\mm}{}&\cdots
}\eqno(**)$$
où $p_{m!}:\Hc(\Fmm)\to\Hc(\Fm)[-d_{\Xg}]$ est l'intégration sur les fibres, de sorte que si $\Xg$ est en plus $i$-acyclique nos méthodes pourront s'appliquer. Or, la suite $(**)$ a encore un sens dans le cas plus général où $\Xg$ est une pseudovariété orientée, par exemple une variété algébrique complexe.
Ainsi, dans le but d'inclure ces espaces dans nos énoncés, nous sommes conduits à remplacer dans $(\ast)$ la cohomologie ordinaire par la \expression{{\bf co}homologie de Borel-Moore}.

\smallskip
\subsubsectionnumber L'\expression{{\bf ho}mologie de Borel-Moore}\label{def-cohomologie-bm} $\Hr^{\bm}_{\ast}(\Mg)$ d'un espace localement compact  $\Mg$, est le dual de sa cohomologie à support compact. Lorsque $\Mg$ est de dimension cohomologique finie $\dMg$ (\ref{dimch}), nous définissons sa \expression{{\bf co}homologie de Borel-Moore}
par l'égalité\glossary{${\Hbm ^{i}(\Mg):=\Hc^{\dMg-i}(\Mg;k)^{\vee}}$:{\bf co}homologie de Borel-Moore d'une pseudovariété $\Mg$ dimension $d_\Mg$}
$$\displayboxit{\Hbm ^{i}(\Mg):=\Hc^{\dMg-i}(\Mg;k)^{\vee}}$$

\subsubsectionline{Fonctorialité de la cohomologie de Borel-Moore.}La fonctorialité de la cohomologie de Borel-Moore est liée à celle de la cohomologie à support compact et cela impose certaines limitations. Par exemple, si $f:\Ng\to\Mg$ est une application continue et \emph{propre}, l'image-inverse des cochaînes d'Alexander-Spanier induit bien un morphisme en cohomologie à support compact $f^{*}:\Hc(\Mg)\to\Hc(\Ng)$ et, par dualité, un morphisme 
$$f_{!}:\Hbm(\Ng)\to\Hbm(\Mg)[\dMg{-}d_{\Ng}]$$
qui étend le morphisme de Thom-Gysin en cohomologie ordinaire. La cohomologie de Borel-Moore est alors fonctorielle sur la catégorie des espaces localement compacts de dimension finie et des application continues et propres.

 En dehors de ce cas, les applications continues n'induisent pas toujours de morphisme en cohomologie de Borel-Moore, ce pour quoi il faut une approche au cas pas cas qui est le but de la section suivante.


\subsection{Quelques morphismes en cohomologie de Borel-Moore}
On\label{integration-sur-fibres}\label{comm-coh-BM-pull-back} donne des analogues à l'image-inverse de la cohomologie ordinaire dans le contexte de la cohomologie de Borel-Moore dans des cas qui intéressent dans les espaces de configuration, notamment: l'action de $\Sm $ sur $\Hbm(\Fm(\Mg))$ et la projection $p_m:\Fmm(\Mg)\onto\Fm(\Mg)$.

\subsubsectionline{Rappel du cas des variétés topologiques.}
Supposons\label{ad} $\Mg$ et $\Ng$ des variétés topologiques \emph{orientées} de dimensions respectives $\dMg$ et $d_{\Ng}$. Notons $\langle\_,\_\rangle_{\Mg}$\glossary{$\langle\_,\_\rangle_{\Mg}$:accouplement de la dualité de Poincaré sur la variété topologique orientée $\Mg$} l'accouplement de la dualité de Poincaré sur $\Mg$ (resp. $\Ng$), à savoir
$$
\langle\_,\_\rangle_{\Mg}:\Hr(\Mg)\times\Hc(\Mg)\to k\,,\quad
\langle\nu,\mu\rangle_{\Mg}:=\smashtop{\int_{\Mg}\nu\wedge\mu}\,.
$$
Pour toute application $f:\Mg\to\Ng$ continue, notons
$$f_{!}:\Hc(\Mg)\to\Hc(\Ng)[-\dMg+d_{\Ng}]$$
l'adjoint pour la dualité de Poincaré de l'image-inverse $f^{*}:\Hr(\Ng)\to\Hr(\Mg)$, il est
caractérisée par l'égalité 
$$\langle f^{*}(\nu),\mu\rangle_{\Mg}=
\langle\nu, f_{!}(\mu)\rangle_{\Ng}\,,\quad\forall
\nu\in\Hr(\Ng)\,,\ 
\forall
\mu\in\Hc(\Mg)\,.
$$
On dira alors que le couple $(f^{*},f_{!})$ est un  \expression{couple adjoint}\comment (sous-entendu pour la dualité de Poincaré)\endcomment.
\comment où $(\_)^{\vee}_{\Mg}$\glossary{$(\_)^{\vee}_{\Mg}$:adjonction par dualité de Poincaré sur $\Mg$} désigne l'adjonction pour la dualité de Poincaré sur $\Mg$.\endcomment

\medskip
\noindent On a quatre cas essentiels pour la suite.
\begin{list}{{\theenumi}}{\usecounter{enumi}\def\theenumi{[ad-\arabic{enumi}]}\mylistskips
\leftmargin1.\parindent
\itemindent0.6cm\labelwidth0.8cm}
\item L'application\label{ad0} $\iota:\Ng\dans\Mg$ est l'inclusion d'une sous-variété fermée orientée. Alors, si  $\iota_{!}:\Hr(\Ng)\to\Hr(\Mg)[\dMg{-}d_{\Ng}]$ est la multiplication par la classe de Thom de $\Ng\dans\Mg$, le couple
$$(\iota_!,\iota^{*})\hbox{\it\ est un couple adjoint.}$$

\item
L'application\label{ad1} $j:\Ug\dans\Mg$ est une inclusion ouverte et $\Ug$ est muni de l'orientation induite. Alors, $j_{!}:\Hc(\Ug)\to\Hc(\Mg)$ est le morphisme de prolongement par zéro et l'on a
$$(j^{*},j_!)\hbox{\it\ est un couple adjoint.}$$

\item
L'application\label{ad2} $f:\Mg\to\Ng$ est localement triviale de fibre~$F$, et $\Mg$ est munie d'une orientation compatible à celles de $\Ng$ et de $F$. L'opération $f_{!}$ coïncide alors à l'\expression{intégration sur les fibres} (\cite{bott-tu} p.~61) et l'on a toujours:
$$(f^{*},f_!)\hbox{\it\ est un couple adjoint.}$$

\item
L'application\label{ad3} $f:\Mg\to\Mg$ est un homéomorphisme. 
Alors
$$f_!:\Hc(\Mg)\to\Hc(\Mg)$$
est induit par l'image-directe des cochaînes à support compact. 

Soit  $\Mg=\coprod\nolimits_{\agoth  \in \Pi_0(\Mg)}\Mg_\agoth  $ la décomposition en composantes connexes. On note par $f:\Pi_0(\Mg)\to \Pi_0(\Mg)$ la bijection induite et $f_{\agoth  }:\Mg_{\agoth  }\to\Mg_{f(\agoth  )}$ la restriction de $f$.
Chaque $f_\agoth  $ est un homéomorphisme de pseudovariétés connexes et orientées. On note $\sigma_{\Mg}(f_\agoth  )$\glossary{$\sigma_{\Mg}(f_\agoth  )$:
action de $f_{\agoth\,!}$ sur la classe fondamentale de $\Mg_\agoth  $, \idest
$f_{!}([\Mg_\agoth  ])=\sigma_{\Mg_\agoth  }(f)\cdot[\Mg_{f(\agoth  )}]$}
le scalaire défini par l'action de $f_{\agoth\,!}$ sur la classe fondamentale $[\Mg_\agoth  ]$ de $\Mg_\agoth  $, \idest tel que:
$$
f_{!}([\Mg_\agoth  ])=
\sigma_{\Mg_\agoth  }(f)\cdot[\Mg_{f(\agoth  )}]\,.\postskip-2pt$$
On a alors
$$\preskip4pt
\mathalign{
\big\langle f^{*}(\nu),\mu\big\rangle{}_{\Mg}&=&
\sumnl_{\agoth  }\int_{\Mg_{\agoth  }}f^{*}(\nu_{f(\agoth  )}\wedge f_{!}(\mu_{\agoth  }))\,d([\Mg_{\agoth  }])\hfill\\
&=&
\sumnl_{\agoth  }\int_{\Mg_{\agoth  }}f^{*}(\nu_{f(\agoth  )}\wedge f_{!}(\mu_{\agoth  }))\,
\sigma_{\Mg}(f_\agoth  )\,d (f_{\agoth\,!}^{-1}[\Mg_{f(\agoth  )}])\\
&=&
\sumnl_{\agoth  }
\big\langle\nu_{f(\agoth  )}, \sigma_{\Mg}(f_\agoth  )\,f_{!}(\mu_\agoth  )\big\rangle{}_{\Mg_{f(\agoth  )}}
\hfill}
\postskip0pt
$$
Par conséquent,
$$\halfdisplayskips
\Big(f^{*},
\sum\nolimits_\agoth  \sigma_{\Mg}(f_\agoth  )\cdot
f_{\agoth\,!}\Big)\hbox{\it\ est un couple adjoint.}\eqno(\diamond)$$
Lorsque $\Mg$ est connexe la somme est réduite à un seul terme et l'adjoint à droite de $f^{*}$ est juste  $\sigma_{\Mg}(f)\cdot f_{!}$.
\end{list}

\subsubsectionline{Suite exacte longue de cohomologie de Borel Moore.}
Soit\label{sex-bm} $\Mg$ un espace localement compact de dimension cohomologique finie $\dMg$. Soit $j:\Ug\dans\Mg$ une inclusion ouverte où $d_{\Ug}=\dMg$ et notons $i:\Ng\dans\Mg$ l'inclusion du fermé complémentaire $\Ng:=\Mg\mmoins \Ug$. En dualisant la suite exacte longue de cohomologie à support compact associée à la décomposition $\Mg=\Ug\sqcup \Ng$, on obtient \expression{la suite exacte longue de cohomologie de Borel-Moore}:
$$\halfdisplayskips
\too\Hbm(\Ng)[d_{\Ng}{-}\dMg]\too^{\iota_{!}}\Hbm(\Mg)\too^{j^{*}}\Hbm(\Ug)\too\,.$$

Lorsque $\Mg$ et $\Ng$ sont des variétés topologiques orientées, $j^{*}$ est la restriction et $\iota_{!}$ est multiplication par la classe de Thom de $\Ng\dans\Mg$.

\subsubsectionline{Le cas des pseudovariétés orientées.}
Pour un espace localement compact $\Mg$ de dimension cohomologique finie, la substitution de la cohomologie ordinaire $\Hr(\Mg)$ par la cohomologie de Borel-Moore $\Hbm(\Mg)$ conduit à substituer l'accouplement $\langle\_,\_\rangle_{\Mg}$ par l'accouplement  $(\_,\_):V^{\vee}\times V\to k$, $(\alpha,v):=\alpha(v)$. Les adjonctions (ad\,2,3,4) servent alors à définir les analogues de l'image-inverse en cohomologie de Borel-Moore de telle sorte qu'ils coïncident avec l'image-inverse en cohomologie ordinaire pour les variétés topologiques orientées. Dans le cas particulier des espaces de configuration, on est conduit aux définitions suivantes.\killline

\begin{defis}Soit\label{def-action-image-inverse} $\Mg$ 
une pseudovariété orientée
(\footnote{C'est-à-dire,
munie d'une section globale nulle part nulle,
sur la partie régulière de $\Mg$, du faisceau $\H^{-\dMg}(\fs\ID_{\Mg}^{\scriptscriptstyle\bullet}(k))$, où $\fs\ID_{\Mg}^{\scriptscriptstyle\bullet}(k)$ est le complexe dualisant de $\Mg$.}).

\begin{enumerate}
\mynobreak\nobreak\item {\bold Image-inverse associée aux projections $p_{b}:\Fba(\Mg)\to\Fb(\Mg)$}\\[2pt]
L'application\label{def-action-image-inverse-a} $p_{b}:\Fba(\Mg)\to\Fb(\Mg)$, $p_b(x_1,\ldots,x_{b+a})= (x_1,\ldots,x_{b})$, est composée de l'inclusion ouverte $\iota:\Fba(\Mg)\dans \Fb(\Mg)\times \Mg^{a}$
et\vadjust{\penalty-50} de la projection $p(\vec x,\vec y):=\vec x$
$$\preskip0pt
\xymatrix@C=1.2cm@R=6mm{
\Fba(\Mg)\ar[rd]_{p_b}\ar[r]^(0.45){\iota}&\Fb(\Mg)\times \Mg^{a}\ar[d]^{p}\\
&\Fb(\Mg)
}$$
L'analogue de l'opération \expression{d'intégration sur les fibres}\glossary{${p_{b!}:\Hc(\Fba(\Mg)\to\Hc(\Fb(\Mg))[-a\,\dMg]}$:intégration sur les fibres}  de \ref{ad2} (\footnote{\label{foot-orientable}C'est pour cette unique raison que nous avons été contraints de nous limiter aux pseudovariétés \emph{orientables}.})
$$p_{b!}:\Hc(\Fba(\Mg)\to\Hc(\Fb(\Mg))[-a\,\dMg]
\postskip4pt
$$
vérifie l'égalité
$$\preskip-1ex
p_{b!}:=p_{!}\circ\iota_{!}\,,$$
où
$
\iota_{!}:\Hc(\Fba(\Mg))\to\Hc(\Fb(\Mg)\times\Mg^{a})$ est 
le prolongement par zéro, et 
$p_{!}:\Hc(\Fb(\Mg))\otimes\Hc(\Mg^{a})\to\Hc(\Fb(\Mg))[-a\,\dMg]$, $\omega\otimes\varpi\mapsto\omega\int_{\Mg}\varpi\,.
$

\smallskip
L'opérateur \expression{image-inverse} pour la cohomologie de Borel-Moore et alors défini, suite à \ref{ad2}, comme le dual vectoriel de $p_{b!}$, donc par:
$$p_b^{*}:=p_{b!}^{\vee}:\Hbm(\Fb(\Mg))\to\Hbm(\Fba(\Mg))\,.$$

\item {\bf\boldmath Action de $\Sm $ sur $\Hbm(\Delta_{?m}\Mg^{m})$. }L'action\label{def-action-image-inverse-b} de $g\in\Sm $ par image-inverse sur $\Hbm(\Fm(\Mg))$ est donnée par l'égalité $(\diamond)$ dans \ref{ad}-\ref{ad3}. L'espace $\Fm(\Mg)$ est un ouvert de $\Mg^{m}$ 
et l'action de $\Sm$ est la restriction de son action sur $\Mg^{m}$,
le scalaire \smash{$\sigma_{\Fm(\Mg)}(g_{\agoth\,!})$}, qui coïncide avec
\smash{$\sigma_{\Mg^{m}}(g_{\bgoth\,!})$} pour une certain $\bgoth\in\Pi_0(\Mg^{m})$, est indépendant de $\bgoth$. En effet, on a $\sigma_{\Mg^{m}}(g_{\bgoth\,!})=\smash{\sgn(g)^{\dim_{\Mg}}}\,,$
où $\sgn(\_)$\glossary{$\sgn(g)$:signature de $g\in\Sm$} est la signature. On pose alors, \glossary{${\sigma_{\Mg^{m}}(g):=\sgn(g)^{\dim_{\Mg}}}$:action de $g\in\Sm$ sur l'orientation de $\Mg^{m}$}
$$\sigma_{\Mg^{m}}(g):=\sgn(g)^{\dim_{\Mg}}
$$ 
de sorte que $g^{*}\rep\Hbm(\Fm(\Mg))\to\Hbm(\Fm(\Mg))$ est donnée par
$$g^{*}=\sigma_{\Mg^{m}}(g)\cdot (g_!)^{\vee}\,,$$
où $g_{!}:\Hc(\Fm(\Mg))\to\Hc(\Fm(\Mg))$ est l'image-directe et où 
$(\_)^{\vee}$\glossary{$(\_)^{\vee}$:dualité vectorielle} désigne l'adjoint pour la dualité vectorielle. Le lemme suivant est immédiat.\killline

\begin{lemm}L'application\label{caractere-image-inverse} $\Sm\ni g\mapsto\sigma_{\Mg^{m}}(g)$ est un caractère multiplicatif. Pour tout $g\in\Sm$, on a
$$\preskip1ex\sigma_{(\Mg\times\Xg)^{m}}(g)=\sigma_{\Mg^{m}}(g)\cdot \sigma_{\Xg^{m}}(g)\,.$$
\end{lemm}
\end{enumerate}
\end{defis}

\def\Remarque{Commentaires}
\noendpoint\begin{rema}\label{comm-coh-BM}
\def\varlistskips{\topsep1pt\itemsep1pt\parskip0pt\parsep0pt}
\begin{enumerate}
\item 
Si\label{comm-coh-BM-etend} $\Mg$ est lisse et orientée, on a $\Hr^{i}(\Fm(\Mg)\simeq\Hbm^{i}(\Fm(\Mg))$ par dualité de Poincaré, et les  définitions d'image-inverse concordent.

\item L'opérateur\label{comm-coh-BM-a} $p_{b!}$ est défini même si la projection $p_b:\Fba(\Mg)\to\Fb(\Mg)$ n'est pas localement triviale, donc même lorsque $\Mg$ n'est pas lisse. Par contre, l'hypothèse d'orientabilité sur $\Mg$ est indispensable.

\item Lorsque\label{comm-coh-BM-Hbm0} $\Mg$ est une pseudovariété \emph{connexe} orientée de dimension $\dMg$, on a $\Hc^{\dMg}(\Mg;\ZZ)=\ZZ$ et donc $\Hbm^{0}(\Mg;\ZZ)\simeq\ZZ\simeq\Hr^{0}(\Mg;\ZZ)$. De plus, si $f:\Mg\to\Mg$ est un homéomorphisme, l'image-inverse $f^{*}$ opère comme l'identité sur $\Hbm^{0}(\Mg;\ZZ)$.
\end{enumerate}
\end{rema}

\subsubsectionline{Ingérence de la signature dans $\Hbm(\Fm)$.}Dans \ref{def-action-image-inverse}-(\ref{def-action-image-inverse-b}) nous avons défini l'action de $\Sm$ sur $\Hbm(\Mg^{m})$ de sorte qu'elle coïncide avec l'action par image-inverse en cohomologie ordinaire $\Hr(\Mg^{m})$ lorsque $\Mg$ est lisse. Cette action {\bf n'est pas} l'action duale de l'action de $\Sm$ sur $\Hc(\Mg^{m})$. La proposition suivante précise la différence entre les deux action dans une situation importante pour la suite de cette section \ref{cohomologie-Borel-Moore}.

\begin{prop}[et notation]Soient\label{prop-signature} $\Mg$ et $\Xg$ deux pseudovariétés orientées et connexes. Pour $x\in\Xg$, posons $\Mg_{\Xg}:=\Mg\times\Xg$ et $\Mg_{x}:=\Mg\times\set x/$. 
Le morphisme dual de la restriction 
$\rho:\Hc(\Fm(\Mg_{\Xg}))\to\Hc(\Fm(\Mg_{x}))$
est le morphisme de $\Sm$-modules
$$\rho^{\vee}(\_)\otimes 1:\Hbm(\Fm(\Mg_{x}))\to\Hbm(\Fm(\Mg_{\Xg}))\otimes(\alt_{m})^{\otimes\dim\Xg}\,,\postdisplaypenalty10000$$
où $\alt_{m}$ est la \expression{représentation par signature} de $\Sm $\glossary{$\alt_{m}$:la représentation par \expression{signature} de $\Sm$} et où $\Hbm(\Mg_{x})$ et $\Hbm(\Mg_{\Xg})$ sont munis des structures de $\Sm$-modules de \ref{def-action-image-inverse}-(\ref{def-action-image-inverse-b}).
\end{prop}
\demo Par définition, la représentation de $\Sm$ sur $\Hbm(\Fm(\Mg_{\Xg}))$ est la représentation duale de celle de $\Hc(\Fm(\Mg_{\Xg}))$ tordue par le caractère $\sigma_{(\Mg\times\Xg)^{m}}=\sigma_{\Mg^{m}}\cdot\sigma_{\Xg^{m}}$ (\ref{caractere-image-inverse}) tandis que pour celle de $\Hbm(\Fm(\Mg_{x}))$, il faut tordre par $\sigma_{\Mg^{m}}$. La différence est donc le caractère
$\sigma_{\Xg^{m}}$ qui est trivial si $\dim\Xg$ est paire est qui est le caractère signature autrement.
\enddemo

\subsection{Approche de $\Hbm(\Fm(\Mg))$ à l'aide d'espaces $i$-acycliques}

Dans les\label{suite-spectrale-i-acyclique} sections \ref{thms-de-scindage} et \ref{coh-cas-i-acyclique} (resp. la section \ref{caractere}), les méthodes pour la détermination du polynôme de Poincaré (resp. du caractère de $\Sm$-module) de $\Hc(\Fm(\Mg))$ s'appliquent lorsque $\Mg$ est $i$-acyclique (\ref{pol-poincare},\,\ref{rema-Fg-orientable}), mais pas lorsque $\Mg$ est général, ce pour quoi il faut une nouvelle idée, comme par exemple, celle que nous donnons à continuation qui permet d'approcher $\Hbm(\Fm(\Mg))$, où $\Mg$ est un espace localement compact \emph{quelconque}, à l'aide d'une suite spectrale dont les termes sont de la forme $\Hbm(\Fl(\Xg))$ avec $\ell\leq m$, et, surtout, où $\Xg$ est   \emph{$i$-acyclique}, suite qu'on appellera \expression{la suite spectrale basique pour $\Hbm(\Fm(\Mg))$} (\cf thm. \ref{theo-suite-spectrale-basique}).

\subsubsectionnumber L'idée est basée sur le fait que:

\proclaim{Tout espace localement compact $\Mg$ peut être réalisé comme fermé dans un espace $i$-acyclique, de complémentaire (donc) également $i$-acyclique.}

\smallskip
En effet, si $\Xg$ est $i$-acyclique et si $x\in\Xg$, l'espace $\Mg$ s'identifie au fermé $\Mg_x:=\Mg\times\set x/$ de l'espace $i$-acyclique $\Mg_{\Xg}:=\Mg\times\Xg$ (\ref{prop-acycliques}-(\ref{prop-acycliques-d})). 
On note alors $V:=\Xg-x$ et $\Mg_{V}:=\Mg\times V$ et l'on remarque que dans la suite longue de cohomologie à support compact:
$$\halfdisplayskips
\to\Hc(\Delta_m(\Mg_{\Xg}^{m}\mmoins \Mg_{x}^{m}))\too^{\iota}
\Hc(\Fm(\Mg_{\Xg}))\too^{\rho}
\Hc(\Fm(\Mg_{x}))\to\,,
$$
on a $\rho=0$. En effet, le prolongement par zéro
$\Hc(\Mg_ V)\to\Hc(\Mg_\Xg)$
est surjectif, les arguments de la remarque \ref{moins-un-point} s'appliquent et établissent la surjectivité du prolongement par zéro
$\tilde\iota:\Hc(\Fm(\Mg_V))\onto\Hc(\Fm(\Mg_{\Xg}))\,.$
La surjectivité de $\iota$ en résulte, car $\im(\iota)\cont\im(\iota')$, et la nullité de $\rho$ s'ensuit.

\smallskip
\noindent Ces arguments et la proposition \ref{prop-signature} prouvent la proposition suivante.\vskip-0.3ex\killline

\begin{prop}
Soient\label{approche-i-acyclique} $\Mg$ et $\Xg$ des espaces localement compacts. On suppose que $\Xg$ est $i$-acyclique et l'on fixe $x\in\Xg$. 
\mynobreak\begin{enumerate}
\nobreak\item 
La\label{approche-i-acyclique-a} suite courte de $\Sm $-modules
$$\halfdisplayskips
\let\Hbm\Hc
\mathalign{
0\to
\Hbm(\Fm(\Mg_{x}))[-1]\to
\Hbm(\Delta_m(\Mg_{\Xg}^{m}\mmoins \Mg_{x}^{m}))\to
\Hbm(\Fm(\Mg_{\Xg}))
\to0}
$$
extraite de la suite longue de cohomologie à support compact, est exacte.
\item
On suppose\label{approche-i-acyclique-b} $\Mg$ et $\Xg$ des pseudovariétés orientées de dimensions cohomologiques finies $\dMg$ et $d_{\Xg}$. 
La suite courte de $\Sm $-modules {\rm(\footnote{L'énoncé\label{note-sans-representation} est valable plus généralement pour les espaces localement compacts de dimensions cohomologiques finies (\ref{dimch}), sauf  pour ce qui est des structures de $\Sm$-modules que nous n'avons pas définies dans cette généralité.})}
$$\skip0=2.5cm
\mathalign{
0\to
\Hbm(\Fm(\Mg_{\Xg}))[md_{\Xg}{-}1]\to\hfill\\
\hskip\skip0\Hbm(\Delta_m(\Mg_{\Xg}^{m}\mmoins \Mg_{x}^{m}))[md_{\Xg}{-}1]\to\hfill\\
\hskip2\skip0\Hbm(\Fm(\Mg_{x}))\otimes(\alt_{m})^{\otimes\dim\Xg}\to0}
$$
extraite de la suite longue de cohomologie de Borel-Moore, est exacte.
\end{enumerate}
\end{prop}

Le polynôme de Poincaré de $\Hbm(\Fm(\Mg))$ est donc déterminé par celui de $\Hbm(\Fm(\Mg_{\Xg}))$, déjà connu, et par celui de $\Hbm(\Delta_m(\Mg_{\Xg}^{m}\mmoins \Mg_{x}^{m}))$,
qui fera l'objet d'étude des sections suivantes.

\subsection{Une suite spectrale pour $\Hbm(\Delta_m(\Mg_{\Xg}^{m}\mmoins \Mg_{x}^{m}))$}

\subsubsectionline{Rappel: cochaînes simpliciales, ordonnées et alternées.}\'Etant\label{cochaines-de-cech} donné une famille finie d'ouverts $\U=\set \UU_1,\ldots,\UU_m/$ d'un espace topologique $\Xg$, on rappelle que l'on dispose classiquement de trois notions de $p$-cochaînes de \vCech\ pour le foncteur des section locales. A savoir,
\def\varitemizeseps{\leftmargin0.5cm\itemsep2pt\parsep0pt\parskip0pt}
\begin{itemize}
\item Le groupe des \expression{$p$-cochaînes simpliciales (non ordonnées)}\glossary{$(\vC^{p}(\U,\_),d)$:complexe de $p$-cochaînes simpliciales (non ordonnées)}
$$\halfdisplayskips
\vC^{p}(\U,\_):=\bigoplusnl_{(i_0,\ldots,i_p)}\Gamma(\UU_{i_{0},\ldots,i_p},\_)\postdisplaypenalty10000$$
où $(i_0,\ldots,i_p)$ est un suite d'éléments deux à deux distincts de $\iii[1,m]$.
\item Le groupe des \expression{$p$-cochaînes (simpliciales) ordonnées}\glossary{$(\vC^{p}_{<}(\U,\_),d)$:complexe de $p$-cochaînes (simpliciales) ordonnées}
$$\vC^{p}_{<}(\U,\_):=\bigoplusnl_{1\leq i_0<\cdots<i_p\leq m}\Gamma(\UU_{i_{0},\ldots,i_p},\_)$$
\item Le groupe des \expression{$p$-cochaînes (simpliciales) alternées}\glossary{$(\vC^{p}_{\varepsilon}(\U,\_),d)$:complexe de $p$-cochaînes  (simpliciales) alternées}
$\vC^{p}_{\varepsilon}(\U,\_)$. C'est le sous-groupe des $p$-cochaînes
$\omega\in\vC^{p}(\U,\_)$ vérifiant pour $\alpha\in\S_{\iii[0,p]}$\glossary{$\sgn(\alpha)$:signature d'une permutation $\alpha$}:
$$\omega_{i_0,\ldots,i_p}=\sgn(\alpha)\,\omega_{i_{\alpha(0)},\ldots,i_{\alpha(p)}}\,.\eqno(\ast)$$

\noindent Remarquons en passant que l'application \expression{d'antisymétrisation}\glossary{${\varepsilon_p:(\vC^{p}_{<}(\U,\_),d)\to 
(\vC^{p}_{\varepsilon}(\U,\_),d)}$:quasi-isomorphisme d'antisymétrisation} 
$$\varepsilon_p:\vC^{p}_{<}(\U,\_)\to \vC^{p}_{\varepsilon}(\U,\_)\eqno(\varepsilon)$$
définie, suivant la même égalité $(\ast)$, par
$
\varepsilon_p(\omega)_{i_{\alpha(0)},\ldots,i_{\alpha(p)}}
:=\sgn(\alpha)\,\omega_{i_0,\ldots,i_p}\,,\postdisplaypenalty10000$
pout tout $\alpha\in\S_{\iii[0,p]}$ et tout $1\leq i_0<\cdots<i_p\leq m$, est  bijective.
 \end{itemize}

\smallskip
Dans les trois cas, l'opérateur cobord $\delta_p:\vC^{p}_{?}(\U,\_)\to \vC^{\pp}_{?}(\U,\_)$
$$(\delta\omega)_{i_0,\ldots,i_{\pp}}=\sumnl_{k=0}^{\pp}
(-1)^{k}\,\omega_{i_0,\ldots,\widehat{i_{k}},\ldots,i_{\pp}}\varrest6_{\UU_{i_{0},\ldots,i_{\pp}}}
$$
a un sens et respecte chaque type de cochaîne. L'antisymétrisation
$$\varepsilon_*:(\vC^{*}_{<}(\U,\_),\delta_*)\to 
(\vC^{*}_{\varepsilon}(\U,\_),\delta_*)\,,$$
est alors un isomorphisme de complexes. 
\noindent L'assertion suivante est classique (\cf\cite{god}, \Spar I.3.8, p.~58.).\killline
 
\begin{prop}
Les inclusions de complexes
$$
(\vC^{*}_{<}(\U,\_),\delta_*)
\dans
(\vC^{*}(\U,\_),\delta_*)
\cont
(\vC^{*}_{\epsilon}(\U,\_),\delta_*)\,.
\eqno(\diamond)$$
sont des quasi-isomorphismes.
\end{prop}

Les complexes $(\diamond)$ sont donc interchangeables pour les besoins du calcul de la cohomologie de \vCech.

\comment\begin{rema}On prendra garde du fait que bien qu'on parle de permutation d'indices des cochaînes de \vCech, il n'y a en général pas d'action du groupe symétrique tout simplement parce qu'il n'y a en général aucun rapport entre $\Gamma(\UU_{i_0,\ldots,i_p},\_)$ et $\Gamma(\UU_{i_{\alpha(0)},\ldots,i_{\alpha(p)}},\_)$. Le paragraphe suivant traite le cas où une telle action existe.
\end{rema}
\endcomment

\subsubsection{Faisceaux et complexes de \vCech\ $\Sm$-équivariants}Revenons\label{symetries-complexes-cech} sur le cas de l'espace $\Delta_m(\Mg_{\Xg}^{m}\mmoins \Mg_{x}^{m})$ de \ref{approche-i-acyclique}. Munissons-le du recouvrement $\U^{m}=\set \UU^{m}_{1},\ldots, \UU^{m}_{m}/$, où:
$$\UU^{m}_{i}:=\Delta_m(\Mg_{\Xg}\times\cdots\times
\smashtop{\aboveparentesis{i}{\,\Mg_{V}}\times\cdots\times\Mg_{\Xg}})\dans \Fm(\Mg_{\Xg})\,,$$
avec $V:=\Xg-x$. Notons ensuite 
$$\UU^{m}:=\UU^{m}_1\cup\cdots\cup\UU^{m}_m\,,\qquad\UU^{m}_{i_0,\ldots,i_p}:=\UU^{m}_{i_0}\cap\cdots\cap \UU^{m}_{i_p},
\postdisplaypenalty10000$$
et même $\UU_{(i_0,\ldots,1_p)}$: une copie de $\UU_{i_0,\ldots,i_p}$ paramétrée par l'uplet $(i_0,\ldots,i_p)$.\glossary{$\UU_{(i_0,\ldots,1_p)}$: une copie de $U_{i_0,\ldots,i_p}$ paramétrée par l'uplet $(i_0,\ldots,i_p)$}

\nobreak Le groupe $\Sm$ agit sur $\UU^{m}$ par permutation des coordonnées, nous avons donc $g\cdot\UU^{m}_{i_0,\ldots,i_p}=\UU^{m}_{g(i_0),\ldots,g(i_p)}$, pour tout $g\in\Sm$. 

\medskip
\soustitreline{\bf Définition.}Un faisceau $\G$ sur $\UU^{m}$ est dit \expression{$\Sm$-équivariant}\glossary{$\G_{\UU^{m}}$:faisceaux $\Sm $-équivariant sur $\UU^{m}$} s'il est muni d'une famille d'isomorphismes 
$\set \phi_{g,V}:\Gamma(gV;\G)\to\Gamma(V;\G)/_{(g,V)}\,,$
indexée par les couples $(g,V)$ où $g\in\Sm$ et $V$ est un ouvert de $\UU^{m}$, telle que la relation cocyclique 
$\phi_{h,gV}\circ\phi_{g,V}=\phi_{hg,V}$ 
est satisfaite pour $h,g\in\Sm$ et tout ouvert $V$ de $\UU^{m}$.
Les définitions de morphisme $\Sm$-équivariant entre faisceaux $\Sm$-équivariants et de catégorie des faisceaux $\Sm$-équivariants s'ensuivent.

\medskip
Les faisceaux des germes de cochaînes de Borel-Moore $\fs\Omega^{i}_{\bm,U^{m}}$ sur $U^{m}$, que nous allons introduire dans \ref{faisceau-cochaines-Borel-Moore}, constituent le principal exemple dans ce travail de faisceau $\Sm$-équivariant. Définis comme duaux des cofaisceaux des cochaînes d'Alexander-Spanier à support compact \smash{$\fs\Omega^{\dMg{-}i}_{{\rm c},U^{m}}$}, les faisceaux \smashtop{$\fs\Omega^{i}_{\bm,U^{m}}$} héritent naturellement de l'action duale de l'action de $\Sm$ sur $\fs\Omega_{{\rm c},U^{m}}$.

\begin{prop}
Munissons\label{symetries-cech} 
$\UU^{m}_{m-p,\ldots,m}=\Delta_m(\Mg_{\Xg}^{m{-}(\pp)}\times\Mg_{V}^{\pp})$  de l'action
de 
$\StS{m{-}(\pp)}{\pp}$
par permutation de coordonnées. 
\begin{enumerate}
\mynobreak\nobreak\item L'application\label{symetries-cech-a}\
$$\preskip0pt
\mathalign{\varPhi_p:&\Sm 
\hskip-4mm
\Times\limits_{\vrule height2mm width0pt
\S_{m{-}(\pp)}\times\1_{\pp}}
\hskip-4mm
\UU^{m}_{m-p,\ldots,m}
&\too&
\coprod\nolimits_{(i_0,\ldots,i_p)}\UU^{m}_{(i_{0},\ldots,i_p)}\hfill\\\noalign{\kern0pt}
&(g,x)\rlap{ $\,\hfto{}{}{2.5cm}$}\qquad&&
\quad g(x)\in \UU^{m}_{(g(m-p),\ldots,g(m))}
}$$
où $(i_0,\ldots,i_p)$ est une suite d'éléments deux à deux distincts de $\iii[1,m]$,
est un homéomorphisme. 

\medskip

\noindent\hskip-\leftmargin Soit maintenant $\G$ un faisceau $\Sm$-équivariant sur $\UU^{m}$. On note $\sigma\mapsto g\star\sigma$ l'action de $g\in\Sm$ sur une section locale $\sigma\in\G$.

\item Le\label{symetries-cech-b} morphisme image-directe défini par l'homéomorphisme $\varPhi_p$, à savoir
$$\mathalign{
\varPhi_{p\,!}:\ind^{\Sm }_{\hskip-3mm\vrule height4mm width0pt
\S_{m{-}(\pp)}\times\1_{\pp}}
\hskip-1cm
\Gamma(\UU^{m}_{m-p,\ldots,m};\G)&\too&\vC^{p}(\U^{m};\G)=
\bigoplus_{(i_0,\cdots,i_{p})}
\Gamma(\UU^{m}_{(i_{0},\ldots,i_{p})};\G)\\
(g,\sigma)\rlap{ $\,\hfto{}{}{2.5cm}$}&&g\star \sigma\in\Gamma(\UU^{m}_{(g(m-p),\cdots, g(m))};\G)\,,\hfill
}$$
est un isomorphisme. Il induit sur $\vC^{p}(\U^{m};\G)$ l'action $\omega\mapsto g\diamond\omega$ de $\Sm$:
$$\displayskips8/10
(h\diamond \omega)_{i_0,\ldots,i_p}:=h\star(\omega_{h^{-1}(i_0),\ldots,h^{-1}(i_p)})\,.$$
 Cette action est compatible au cobord des cochaînes simpliciales et le complexe de \vCech\ augmenté:
$$\displayskips8/10
0\to\Gamma(\UU^{m};\G)\too^{\epsilon}
\vC^{0}(\U;\G)\too^{\delta_0}
\vC^{1}(\U;\G)\too^{\delta_1}\cdots
$$
est un complexe de $\Sm $-modules.
\item Le\label{symetries-cech-c} complexe 
des cochaînes alternées 
$(
\vC^{*}_{\varepsilon}(\U^{m};\G),\delta_*)
$
est un sous-complexe de $\Sm $-modules
du complexe des cochaînes simpliciales $(
\vC^{*}(\U^{m};\G),\delta_*)
$. 
L'antisymétrisation
$\varepsilon_{*}:(\vC^{*}_{<}(\U^{m};\G),\delta_*)\to(\vC^{*}_{\varepsilon}(\U^{m};\G),\delta_*)$ transfère cette structure et munit chaque groupe $\vC^{p}_{<}(\U^{m};\G)$ de l'action $\omega\mapsto g\star \omega$ de $\Sm $.
En particulier, le complexe de cochaînes ordonnées de \vCech\ augmenté:
$$0\to\Gamma(\UU^{m};\G)\too^{\epsilon}
\vCm^{0}(\U;\G)\too^{\delta_0}
\vCm^{1}(\U;\G)\too^{\delta_1}\cdots
$$
est un complexe de $\Sm $-modules.
De plus, l'application
$$\mathalign{\mathrigid0mu
\varPsi_{p}\kern2pt:\kern2pt\ind^{\Sm }_{\vrule height4.5mm width0pt\hskip-4mm
\StS{m{-}(\pp)}{\pp}}
\hskip-14mm 
\Gamma(\UU^{m}_{m-p,\ldots,m};\G)
\otimes \alt_{\pp}&\to
&\vC^{p}_{<}(\U^{m};\G)=
\hskip-7.5mm
\bigoplus_{1\leq i_0<\cdots<i_{p}\leq m}
\hskip-7.5mm
\Gamma(\UU^{m}_{i_{0},\ldots,i_{p}};\G)\\\noalign{\kern4pt}
(g,\sigma)\rlap{ $\,\hfto{}{}{2.9cm}$}&&g\diamond \sigma\in\Gamma(g(\UU^{m}_{m-p,\ldots,m});\G)
}$$
où {\rm`$\diamond$'} désigne l'action {\rm`$\star$'} 
de $\S_{m{-}(\pp)}\times\S_{\pp}$ tordue par 
le caractère signature $\alt_{\pp}$
de $\S_{\pp}$, 
est un isomorphisme de $\Sm $-modules. 
\end{enumerate}\end{prop}

\demo (\ref{symetries-cech-a}) 
L'application $\varPhi_p$ est définie sur la réunion disjointe de  copies de  $\UU^{m}_{m-p,\ldots, m}$ indexées par les éléments $\cl g\in \Sm /\S_{m{-}(\pp)}$. Pour chaque $g\in\Sm $, la restriction de  $\varPhi_p$ à $(\cl g, \UU^{m}_{m-p,\ldots,m})$ est un homéomorphisme sur $\UU^{m}_{(g(m-p),\ldots,g(m))}$. On conclut que $\varPhi_p$ est bijective en remarquant que le cardinal $|\Sm /\S_{m{-}(\pp)}|$ est précisément celui de l'ensemble des $(\pp)$-uplets $(i_0,\ldots,i_p)$ d'éléments deux à deux distincts de $\iii[1,m]$.  

\noindent (\ref{symetries-cech-b}) résulte de (\ref{symetries-cech-a}) et des identités:
$$\def\sep{\noalign{\kern4pt}}
\mathalign{
\delta(h\diamond \omega)_{i_0,\ldots,i_{\pp}}
&=&\sumnl_{j=0}^{\pp}
(h\diamond\omega)_{i_0,\ldots,\widehat{i_j},\ldots,i_{\pp}}\varrest6_{\UU^{m}_{(i_{0},\ldots,i_{\pp})}}
\hfill\\\sep
&=&\sumnl_{j=0}^{\pp}
\Big(h\star(\omega_{h^{-1}(i_0,\ldots,\widehat{i_j},\ldots,i_{\pp})})\Big)\varrest6_{\UU^{m}_{(i_{0},\ldots,i_{\pp})}}\hfill
\\\sep
&=&h\star\Big(
\sumnl_{j=0}^{\pp}
\omega_{h^{-1}(i_0,\ldots,\widehat{i_j},\ldots,i_{\pp})}\varrest6_{\UU^{m}_{h^{-1}(i_{0},\ldots,i_{\pp})}}
\Big)\\\sep
&=&(h\diamond \delta\omega)_{i_0,\ldots,i_{\pp}}\,.\hfill
}
$$

\noindent (\ref{symetries-cech-c}) Soit $\omega\in\vC^{p}_{\varepsilon}(\U^{m};\G)$. Pour \smashbot{$\alpha\in\S_{\iii[0,p]}$} et $h\in\Sm $, on a:
$$\def\sep{\noalign{\kern4pt}}
\mathalign{
(h\diamond\omega)_{i_{\alpha(0)},\ldots,i_{\alpha(p)}}
&=&
h\star(\omega_{h^{-1}(i_{\alpha(0)}),\ldots,h^{-1}(i_{\alpha(p)})}
)\hfill\\\sep
&=&
\sgn(\alpha)h\star(\omega_{h^{-1}(i_{0}),\ldots,h^{-1}(i_{p})}
)\\\sep
&=&\sgn(\alpha)(h\diamond\omega)_{i_{0},\ldots,i_{p}}\,,\hfill
}
\postdisplaypenalty10000
$$
et 
$\vC^{p}_{\varepsilon}(\U^{m};\G)$ est bien un sous-$\Sm $-module de
$\vC^{p}_{\varepsilon}(\U^{m};\G)$.

\smallskip
Le sous-espace $\varepsilon_p(\Gamma(\UU^{m}_{m-p,\ldots,m};\G))$ est
stable sous $\StS{m{-}(\pp)}{\pp}$ dont l'action est tordue par le caractère signature de $\S_{\pp}$. Le morphisme de $\Sm $-modules $\varPsi_{p}$ est donc bien défini et il est surjectif puisque $\varepsilon_p(\Gamma(\UU^{m}_{m-p,\ldots,m};\G))$ engendre clairement $\vC^{p}_{\varepsilon}(\U^{m};\G)$ en tant que $\Sm $-module. Le fait que $\varPsi_{p}$ est bijectif résulte alors du fait que $|\Sm /(\StS{m{-}(\pp)}{\pp})|$ est également le cardinal de l'ensemble des parties  $I\dans\iii[1,m]$ telles que $|I|=\pp$.
\enddemo

\subsubsectionline{Suite spectrale de Borel-Moore d'un $\Gg$-espace.}
Nous\label{faisceau-cochaines-Borel-Moore}\label{ss-coh-supp-comp-recouvrement} rappelons maintenant, pour un $\Gg$-espace donné, les bases théoriques de la construction de la suite spectrale de $\Gg$-modules associée à un recouvrement $\Gg$-stable et
 pour la cohomologie de Borel-Moore.

\soustitre{Faisceaux de germes de cochaînes de Borel-Moore}
Soit $\Mg$ un espace localement compact de dimension cohomologique $\dMg$ muni de l'action d'un groupe fini $\Gg$. Dans \ref{prop-finitude-Hc}-(\ref{prop-finitude-Hc-0}) nous avons introduit la résolution $c$-molle
$$\preskip0pt
0\to\fs k_{\Mg}\hoook
\fs\Omega^{0}_{\Mg}
\too^{d_{0}}
\fs\Omega^{1}_{\Mg}
\too^{d_{1}}
\cdots
\hf{d_{\dMg{-}1}}{}{1.cm}
\fs\Omega^{\dMg}_{\Mg}
\to0
$$
qui est une résolution dans la catégorie des faisceaux $\Gg$-équivariants (\ref{symetries-complexes-cech}).
Pour tout ouvert $U\dans\Mg$, le complexe des \expression{cochaînes à support compact}\glossary{${(\Omega_{\rm c}^{*}(U),d):=\Gamma(U;
(\AS^{\bullet}(\Xg;k),d_{*}))}$:le complexe des cochaînes  à support compact de $\U$}
$$(\Omega_{\rm c}^{*}(U),d_{}):=\Gammac(U;
(\fsOmega^{*}_{\Mg},d_{*}))$$
calcule la cohomologie à support compact $\Hc^{*}(U)$, et si $U$ est en plus stable sous l'action de $\Gg$, c'est aussi un complexe de $\Gg$-modules pour l'action d'image-directe topologique.

Si
 $\iota_{V\dans U}:V\dans U$ est une inclusion ouverte, le prolongement par zéro
$$\iota_{V\dans U!}:(\Omega_{\rm c}^{*}(V),d_{*})\to(\Omega_{\rm c}^{*}(U),d_{*})$$
est une inclusion de complexes et la correspondance
$$\def\tt{\vrule height 12pt width0pt}
\xymatrix@R=6mm{
U\ar@{~>}[r]&(\Omega_{\rm c}^{*}(U),d_{*})\\
\tt V\arinto[u]\ar@{~>}[r]&
\tt (\Omega_{\rm c}^{*}(V),d_{*})\arinto[u]_{\iota_{V\dans U}}
}$$
définit un complexe $(\fs\Omega_{\Mg,\rm c}^{*},d_{*})$\glossary{$(\fs\Omega_{\Mg,\rm c},d)$:complexe de cofaisceaux (flasques) des cochaînes à support compact sur $\Mg$} de \emph{pré-cofaisceaux flasques} sur $\Mg$ qui sont en fait des \emph{cofaisceau} puisque les $\fsOmega^{i}_{\Mg}$ sont $c$-mous (\cf\cite{bredon} V.1.6, p.~282).
\medskip

Pour tout ouvert $U\dans\Mg$, le \expression{complexe des cochaînes de Borel-Moore sur $U$}\glossary{${(\Omega_{\bm}^{*}(U),d_{*}):=((\Omega_{\rm c}^{*}(U),d_{*})^{\vee})[-d_\Mg]}$:le complexe des cochaînes de Borel-Moore sur  $U\dans\Mg$} est défini par dualité et décalage: (\footnote{On remarquera l'abus de notation qui consiste à noter de la même manière les  différentielles des deux complexes 
$(\Omega_{\rm c}^{*}(U),d_{*})$ et $(\Omega_{\bm}^{*}(U),d_{*})$.})
$$(\Omega_{\bm}^{*}(U),d_{*}):=((\Omega_{\rm c}^{*}(U),d_{*})^{\vee})[-d_\Mg]\,.\eqno(\ddagger)$$
Lorsque $U$ est $\Gg$-stable, chaque $\Omega_{\bm}^{i}(U)$ est à priori muni de l'action de $\Gg$ duale de son action par image-directe sur $\Omegac^{i}(U)$, mais cette action ne correspond pas à l'action de $\Gg$ par image-inverse lorsque $\Mg$ est une variété topologique orientée, ce pour quoi il faut tordre l'action duale par le caractère $\sigma_{\Mg}$ tel qu'expliqué dans 
\ref{def-action-image-inverse}-(\ref{def-action-image-inverse-b}).
Plus généralement, si $\Mg$ est une pseudovariété orientée, on notera~`$\star$' l'action duale de $\Gg$ sur $\Omegac^{i}(\_)$ tordue par $\sigma_{\Mg}$. Dans tous les cas, le complexe $(\ddagger)$ ci-dessus est un complexe de $\Gg$-modules.

Le complexe $(\Omega_{\bm}^{*}(U),d_*)$ calcule la cohomologie de Borel-Moore $\Hbm^{*}(U)$ et
la correspondance
$$\def\tt{\vrule height 12pt width0pt}
\xymatrix@R=5mm{
U\ar@{~>}[r]&(\Omega_{\bm}^{*}(U),d_{*})\aronto[d]^(0.45){\,\iota_{V\dans U}^{*}}\\
\tt V\arinto[u]\ar@{~>}[r]&
\tt (\Omega_{\bm}^{*}(V),d_{*})
}$$
où $\iota_{V\dans U}^{*}$ est le dual de $\iota_{V\dans U!}$, est le complexe $(\fs\Omega_{\Mg,\Bm}^{*},d_{*})$\glossary{$(\fs\Omega_{\Mg,\Bm},d_*)$:complexe de faisceaux (flasques) de germes de cochaînes de Borel-Moore de $\Mg$} des \expression{faisceaux (flasques) de germes de cochaînes de Borel-Moore} sur $\Mg$. C'est un complexe de faisceaux $\Gg$-équivariants.

\smallskip
\soustitre{Bicomplexe de cochaînes \vCech-Borel-Moore}
Soit $\U^{m}=\set\UU^{m}_{1},\ldots,\UU^{m}_{m}/$ 
une famille $\Gg$-stable d'ouverts de $\Mg$, c'est à dire telle qu'il existe une action de $\Gg$ sur $\iii[1,m]$ vérifiant $g(\UU^{m}_{i})=\UU^{m}_{g(i)}$. On pose
$$\UU^{m}:=\UU^{m}_1\cup\cdots\cup\UU^{m}_m\,,\qquad\UU^{m}_{i_0,\ldots,i_p}:=\UU^{m}_{i_0}\cap\cdots\cap \UU^{m}_{i_p},
$$
et l'on considère le bicomplexe $(\vC(\U^{m})^{*}_{\rm c,\bullet},\partial_\bullet,d_{*})$\glossary{${\vC(\U^{m})^{*}_{\rm c,\bullet}:=\big(\vC^{\bullet}(\U^{m};\fs\Omega_{\Mg,\rm c}^{*}(\_)),\partial_{\bullet},d_*\big)}$:bicomplexe de chaînes de \vCech\ à support compact de $\U^{m}$}
:
$$\preskip0.5ex\postskip0.5ex\let\delta\partial
\def\ligne#1#2{
0\ar[r]&\vC_{#1}(\U^{m},\fs\Omega_{\Mg,\rm c}^{0})\ar[r]^{d_0}\ar[d]#2&
\vC_{#1}(\U^{m},\fs\Omega_{\Mg,\rm c}^{1})\ar[r]^{d_1}\ar[d]#2&
\vC_{#1}(\U^{m},\fs\Omega_{\Mg,\rm c}^{2})\ar[r]^(0.75){d_2}\ar[d]#2&}
\xymatrixc{@R=5mm@C=5mm}{
&\ar[d]&\ar[d]&\ar[d]\\
\ligne1{_{\delta_{0}}}\\
\ligne0{_\epsilon}\\
0\ar[r]&
\fs\Omega_{\Mg,c}^{0}(\UU^{m})\ar[r]^{d_0}\ar[d]&
\fs\Omega_{\Mg,c}^{1}(\UU^{m})\ar[r]^{d_1}\ar[d]&
\fs\Omega_{\Mg,c}^{2}(\UU^{m})\ar[r]^(0.75){d_2}\ar[d]&
\\
&0&0&0
}
\eqno\raise23pt\hbox to2pt{\hss$(\vC_{\rm c})$}$$
dont les colonnes sont les complexes de \expression{chaînes de \vCech\  à valeurs dans un cofaisceau} et
 sont exactes puisque les cofaisceaux $\fsOmega^{i}_{\Mg,\rm c}$ sont flasques (\footnote{\Cf Bredon \cite{bredon} chap.~VI.--- Cosheaves and \vCech\ Homology, corollary VI.4.5, p.~426.}).
\comment
En particulier, le complexe simple associé au bicomplexe  $\big(\vC^{\bullet}(\U^{m};\fs\Omega_{\Mg,\rm c}^{*}),\partial_{\bullet},d_*\big)$ calcule la cohomologie $\Hc(\UU^{m})$.
\endcomment
En dualisant, on obtient un bicomplexe de colonnes exactes
$$\preskip0.5ex\postskip0.5ex
\def\ligne#1#2{
0\ar[r]&\vC^{#1}(\U^{m},\fs\Omega_{\Mg,\Bm}^{0})\ar[r]^{d_0}\ar@{<-}[d]#2&
\vC^{#1}(\U^{m},\fs\Omega_{\Mg,\Bm}^{1})\ar[r]^{d_1}\ar@{<-}[d]#2&
\vC^{#1}(\U^{m},\fs\Omega_{\Mg,\Bm}^{2})\ar[r]^(0.75){d_2}\ar@{<-}[d]#2&}
\xymatrixc{@R=5mm@C=5mm}{
&\ar@{<-}[d]&\ar@{<-}[d]&\ar@{<-}[d]\\
\ligne1{_{\delta_{0}}}\\
\ligne0{_(0.55){\epsilon}}\\
0\ar[r]&
\fs\Omega_{\Mg,\Bm}^{0}(\UU^{m})\ar[r]^{d_0}\ar@{<-}[d]&
\fs\Omega_{\Mg,\Bm}^{1}(\UU^{m})\ar[r]^{d_1}\ar@{<-}[d]&
\fs\Omega_{\Mg,\Bm}^{2}(\UU^{m})\ar[r]^(0.75){d_2}\ar@{<-}[d]&
\\
&0&0&0
}\eqno\raise25pt\hbox to2pt{\hss$(\vC_{\bm})$}$$
qui est un bicomplexe de $\Gg$-modules puisque l'analogue du théorème \ref{symetries-cech} est vérifié pour les données en cours.

On appellera \expression{bicomplexe de cochaînes  de \vCech-Borel-Moore  de $\U^{m}$}\glossary{${\vC(\U^{m})^{\bullet,*}_{\Bm}:=\big(\vC^{\bullet}(\U^{m};\fs\Omega_{\Mg,\Bm}^{*}),\delta_{\bullet},d_*\big)}$:bicomplexe de cochaînes de \vCech-Borel-Moore de $\U^{m}$}, le bicomplexe du premier quadrant 
$$\vC(\U^{m})^{\bullet,*}_{\Bm}:=\big(\vC^{\bullet}(\U^{m};\fs\Omega_{\Mg,\Bm}^{*}),\delta_{\bullet},d_*\big)\,.$$

\smallskip
\soustitre{Filtration régulière de $\Hbm^{*}(\UU^{m})$ et suite spectrale}
En raison de l'exactitude des colonnes de $(\vC_{\bm})$, le morphisme d'augmentation $\epsilon$ induit un quasi-isomorphisme de complexes de $\Gg$-modules
 $$\epsilon_*:(\Omega_{\bm}^{*}(\UU^{m}),d_*)\to\tot^{*}(\vC(\U^{m})^{\bullet,*}_{\Bm})
$$ 
où `$\tot$' désigne le complexe simple associé. On a donc un isomorphisme de $\Gg$-modules
$$h^{i}(\epsilon_{*}):\Hbm^{i}(\UU^{m})\simeq h^{i}(\tot^{*}(\vC(\U^{m})^{\bullet,*}_{\Bm})$$
pour tout $i\in\ZZ$.

Pour $k\in\NN$, notons $\vC(\U^{m})^{\bullet,*}_{\Bm,k}$ le sous-bicomplexe de $\vC(\U^{m})^{\bullet,*}_{\bm}$ défini par
$$\begin{cases}
\vC(\U^{m})^{i,*}_{\Bm,k}=0\,,&\hbox{ si $i<k$,}\\[6pt]
\vC(\U^{m})^{i,*}_{\Bm,k}=\vC(\U^{m})^{i,*}_{\bm}
\,,&\hbox{ si $k\leq i\leq m{-}1$.}\\
\end{cases}$$
On a la filtration décroissante de bicomplexes de $\Gg$-modules
$$\vC(\U^{m})^{\bullet,*}_{\bm}=\vC(\U^{m})^{\bullet,*}_{\Bm,0}\cont
\vC(\U^{m})^{\bullet,*}_{\Bm,1}\cont\vC(\U^{m})^{\bullet,*}_{\Bm,2}\cont\cdots$$
induisant une filtration positivement graduée décroissante et régulière du complexe $\tot^{*}(\vC(\U^{m})^{\bullet,*}_{\bm})$ de même donc que pour $\Hbm^{*}(\UU^{m})$ pour lequel on pose
$$\Hbm^{*}(\UU^{m})_{k}:=h(\epsilon_*)^{-1}\big(h^{*}(\tot^{*}(\vC(\U^{m})^{\bullet,*}_{\Bm,k})))\big)\,.$$

On notera $\IF\Hbm^{*}(\UU^{m})$ et $\IF\tot^{*}(\vC(\U^{m})^{\bullet,*}_{\Bm})$\glossary{${\IF\Hbm^{*}(\U^{m}),\ \IF\tot^{*}(\vC(\U^{m})^{\bullet,*}_{\Bm})}$:filtrés par degrés de cochaînes de \vCech\ des objets concernés} ces objets filtrés.

\begin{rema}
Dans ce qui précède nous aurions tout aussi bien pu considérer les cochaînes de \vCech\ ordonnées ou alternées, ce qui aurait fournit les bicomplexes de $\Gg$-modules $\vCm(\U^{m})^{\bullet,*}_{\bm}$ et $\vC_{\epsilon}(\U^{m})^{\bullet,*}_{\bm}$.
\end{rema}

\noindent La proposition suivante est bien connue  (\footnote{\Cf Godement \cite{god} chap. I.4, thm. 4.2.2.}).\killline

\begin{prop}La\label{comm-coh-BM-suite-spectrale} suite spectrale de $\Gg$-modules $\IE(\U^{m}):=(\IE(\U^{m})_r,d_r)$ associée au complexe de $\Gg$-modules gradué filtré $\IF\tot^{*}(\vC_{?}(\U^{m})^{\bullet,*}_{\Bm}))$
converge vers le bigradué de $\IF\Hbm^{*}(\UU^{m})$. On a
$$\IE(\U^{m})_{1}^{p,q}:=\vC_{?}^{p}(\U^{m},\H_{\bm}^{q}(\_))
\postskip-1ex
$$
et
$$\preskip2pt\IE(\U^{m})_{2}^{p,q}:=\vH^{p}(\U^{m},\H_{\bm}^{q}(\_))\Longrightarrow
\Gr^{p,q}(\IF\Hbm^{*}(\U^{m}))\,.\eqno(\IE)$$
\end{prop}

\medskip\noindent
Le théorème suivant est maintenant corollaire immédiat de \ref{comm-coh-BM-suite-spectrale} et de \ref{symetries-cech}.\killline

\begin{theo}Soient $\Mg$ et $\Xg$ des espaces localement compacts de dimensions cohomologiques finies.
Le\label{theo-approche-i-acyclique} $\Sm $-module $\Hbm(\Delta_m(\Mg_{\Xg}^{m}\mmoins \Mg_{x}^{m}))$ est l'aboutissement de la suite spectrale $\IE(\U^{m})$  de {\rm \ref{comm-coh-BM-suite-spectrale}} pour le recouvrement
$\U^{m}:=\set \UU^{m}_{1},\ldots, \UU^{m}_{m}/$, où
$$\UU^{m}_{i}:=\Delta_m(\Mg_{\Xg}\times\cdots\times
\smashtop{\aboveparentesis{i}{\,\Mg_{\Xg- x}}\times\cdots\times\Mg_{\Xg}})\dans \Fm(\Mg_{\Xg})\,.$$
L'identification de \ref{symetries-cech}-(\ref{symetries-cech-c}) pour le groupe des $p$-cochaînes ordonnées
$$\vC^{p}_{<}(\U^{m};\fs\Omega_{\UU^{m},\bm}^{*})=
\ind^{\Sm }_{\vrule height4mm width0pt\hskip-4mm
\StS{m{-}(\pp)}{\pp}}
\hskip-10mm
\Gamma(\UU^{m}_{m-p,\ldots,m};\fs\Omega_{\UU^{m},\bm}^{*})\otimes\alt_{\pp}\,,
$$
munit
$(\vC^{*}_{<}(\U^{m},\_),\delta_*)$ d'une structure de complexe de $\Sm $-modules. Les ter\-mes $(\IE(\U^{m})_r,d_r)$
de la suite spectrale $\IE(\U^{m})$ héritent  d'une structure de complexe de $\Sm $-modules et la suite spectrale converge au sens de suite spectrale de complexes de $\Sm $-modules vers le $\Sm $-module bi-gradué associé au $\Sm $-module gradué $\Hbm(\Delta_m(\Mg_{\Xg}^{m}\mmoins \Mg_{x}^{m}))$ muni de filtration régulière de la proposition \ref{comm-coh-BM-suite-spectrale}. Dans le cas particulier où $r=1$, on a 
$$\IE(\U^{m})_1^{p,q}=
\ind^{\Sm }_{\vrule height4mm width0pt\hskip-4mm
\StS{m{-}(\pp)}{\pp}}
\hskip-10mm
\Hbm^{q}(\UU^{m}_{m-p,\ldots,m}
)\otimes\alt_{\pp}\,,
$$
et $d_1:\IE(\U^{m})_{1}^{p-1,q}\to \IE(\U^{m})_{1}^{p,q}$ un morphisme de $\Sm $-modules. 

\end{theo}

\begin{rema}Le même énoncé est valable pour le complexe des cochaînes non ordonnées $(\vC^{p} (\U^{m},\_),\delta_*)$, auquel cas on a
$$\IE(\U^{m})_1^{p,q}=
\ind^{\Sm }_{\vrule height4mm width0pt\hskip-4mm
\S_{m{-}(\pp)}\times\1_{\pp}}
\hskip-10mm
\Hbm^{q}(
\UU^{m}_{m-p,\ldots,m}
)\,,
$$
en raison de l'égalité
\smash{$\vC^{p}(\U^{m},\_)=
\ind^{\Sm }_{\vrule height4mm width0pt\hskip-4mm
\S_{m{-}(\pp)}\times\1_{\pp}}
\hskip-10mm
\Gamma(\UU_{m-p,\ldots,m}^{m},\_)
$}
de \ref{symetries-cech}-(\ref{symetries-cech-b}).
\end{rema}

\subsection{La suite spectrale \og{\slshape basique\/}\fg\ pour $\Hbm(\Fm(\Mg))$}\label{suite-spectrale-basique}

Nous nous restreignons maintenant à une situation qui simplifie remarquablement les considérations précédentes.
C'est le cas où l'espace $\Xg$ est l'espace $\RR_{\geq0}:=[0,+\infty[$, puis $x:=0$ et $V:=\RR_{>0}$. 
Dans la suite, $\Mg$ sera localement compact de dimension cohomologique finie $\dMg$ (\ref{dimch}), et l'on notera
$$\let\qquad\quad
\Mg_{\geq0}:=\Mg\times\RR_{\geq0}\,,\qquad
\Mg_{0}:=\Mg\times\set0/\,,\qquad
\Mg_{>0}:=\Mg\times\RR_{>0}\,.
$$
Les espaces $\Mg_{\geq0}$ et $\Mg_{>0}$ sont $i$-acycliques et  l'on a 
$$\Mg\simeq\Mg_{0}=\big(\Mg_{\geq0}\moins\Mg_{>0}\big)\,.$$

On aura remarqué que ces choix renferment l'égalité $\Hc(\RR_{\geq0})=0$, et donc le fait que l'on aura (\ref{pol-poincare}-(\ref{pol-poincare-a}))
$$\Hc(\Fm(\Mg_{\geq0}))=0\,,\quad\forall m\geq1\,.$$

\subsubsectionnumber La\label{iso-tordu} proposition \ref{approche-i-acyclique} se simplifie et donne des isomorphismes de $\Sm$-modules 
$$
\Hc(\Fm(\Mg))\simeq\Hc(\Delta_{m}(\Mg_{\geq0}^{m}\moins\Mg_{0}^{m}))[1],
\eqno(\diamond)$$
et lorsque $\Mg$ est une pseudovariété orientée (\footnote{Dans le cas général où $\Mg$ n'est pas une pseudovariété, on a toujours un isomorphisme d'espaces vectoriels, mais on perd l'aspect représentation de $\Sm$ (\cf note ($^{\ref{note-sans-representation}}$)).})
$$
\displayboxit{\Hbm(\Fm(\Mg))
\simeq
\Hbm(\Delta_{m}(\Mg_{\geq0}^{m}\moins\Mg_{0}^{m}))[\mmo]\otimes\alt_{m}}
\eqno(\diamond\diamond)$$
où l'apparition du caractère signature a été justifiée dans \ref{prop-signature}.

Avec ces données, la suite spectrale de \ref{theo-approche-i-acyclique} convenablement modifiée par le caractère $\alt_{m}$, converge vers $\Hbm(\Fm(\Mg))[1{-}m]$, on l'appellera
\expression{la suite spectrale basique pour $\Hbm(\Fm(\Mg))$}.

\begin{theo}[(des suites spectrales basiques)]Soit\label{theo-suite-spectrale-basique} $\Mg$ une pseudovariété  orientée de dimension $\dMg$ {\rm(\footnote{Comme dans la note précédente, si l'on néglige l'aspect représentations, l'énoncé est valable plus généralement pour les espaces localement compacts de dimension finie.})}. La suite $\IEs(\U^{m}):=\IE(\U^{m})\otimes\alt_{m}$\glossary{${\IEs(\U^{m}):=\IE(\U^{m})\otimes\alt_{m}}$:suite spectrale pour le complexe de faisceaux $\Sm$-équivariants $\fs\Omega_{\UU^{m},\bm}^{*}\otimes\alt_{m}$}  converge, en tant que suite spectrale de complexes de $\Sm $-modules, vers le $\Sm $-module bi-gradué associé 
au $\Sm $-module gradué $\Hbm(\Fm(\Mg))[1{-}m]$ muni de filtration régulière induite par l'isomorphisme $(\diamond\diamond)$ ci-dessus. De plus,
\begin{enumerate}\itemsep8pt
\item
Pour\label{theo-suite-spectrale-basique-a} tout $i\in\ZZ$, on a
$$\hss\displayboxit{\IEs(\U^{m})_{1}^{p,q}=
\ind^{\Sm }_{\vrule height4.5mm width0pt\hskip-4mm
\StS{m{-}(\pp)}{\pp}}
\hskip-13mm
\alt_{m{-}(\pp)}\otimes\Hbm^{q}(\UU^{m}_{m-p,\ldots,m})\ \Rightarrow\ 
\Hbm^{i}(\Fm(\Mg))}\hss$$
avec $q=i+(m{-}(\pp))$. Dans cette écriture, le groupe $\StS{m{-}(\pp)}{\pp}$ opère sur $\Hbm(\UU^{m}_{m-p},\ldots,m)$ par image-inverse {\rm(\ref{def-action-image-inverse}-(\ref{def-action-image-inverse-b}))} et le caractère $\alt_{m{-}(\pp)}$ affecte uniquement l'action du sous-groupe $\S_{m{-}(\pp)}\times\1$.

\item On\label{theo-suite-spectrale-basique-b} a une décomposition canonique d'espaces vectoriels\glossary{$\varXi^{m}_{\pp}$:décomposition canonique de $\Hbm(\UU^{m}_{m-p,\ldots,m})$}
$$\displayboxit{\varXi^{m}_{\pp}:\Hbm^{q}(\UU^{m}_{m-p,\ldots,m})\simeq
\bigoplusnl_{f\in\FFF(\pp,m)}
\Hbm^{q-(m{-}(\pp))\,\dMg}(\Fg_{\pgoth(f)}(\Mg_{>0}))}$$
où $\FFF(\pp,m)$\glossary{$\FFF(\pp,m)$:ensemble des applications 
$f:\iii[1,m]$ strictement croissantes sur $\iij[1,m-p]$ qui fixent $\iii[m{-}p,m]$}  est l'ensemble des applications $\mathrigid2mu
f:\iii[1,m]\to\iii[1,m]$ telles que 
$x<f(x)$ si $x<m-p$, et $f(x)=x$ sinon.
On a $\displaystyle
|\FFF(\pp,m)|={(\mmo)!/ p!}$.
La paramétrisation associe à $f\in\FFF(\pp,m)$ le sous-espace $\Fg_{\pgoth(f)}(\Mg_{>0})\dans\smash{\Delta_{\pp}(\Mg_{\geq0}^{m-(\pp)}\times\Fpp(\Mg_{>0}))}$ où $\pgoth(f)$ est la partition de $\iii[1,m]$ définie par les fibres de $f$ {\rm(\cf \ref{nota-pgoth})}.

\item
Pour\label{theo-suite-spectrale-basique-c} avoir $\IEs(\U^{m})_{1}^{p,q}\not=0$, il faut que
$(m{-}(\pp))\dMg\leq q\leq m\,\dMg\,.$

\item Pour\label{theo-suite-spectrale-basique-d} $i\in\NN$ donné, les termes de $\IEs(\U^{m})_{r}^{p,q}$, pour $r\geq1$, qui contribuent à $\Hbm^{i}(\Fm(\Mg))$ sont ceux pour lesquels on a
$(m{-}(\pp))(d_{\Mg}{-}1)\leq i\,.$
\end{enumerate}
\end{theo}
\def\Demonstration{Indications}\demo
(\ref{theo-suite-spectrale-basique-a}) 
Le foncteur $(\_)\otimes\alt_{m}$ étant exact, la convergence de $\IEs(\U^{m})$ vers $\Hbm(\UU^{m})\otimes\alt_{m}$ est assurée. Ensuite, 
l'égalité classique (\footnote{\label{note-ind-res}\'Etant donnée une inclusion de groupes  $H\dans G$, un $H$-module $V$ et un $G$-module W, on a $\ind_{H}^{G} V\otimes_k W\cong\ind^{G}_{H}(V\otimes_k\mathop{\rm res}^{G}_{H} W)
$.})
$$\mathalign{
\Big(\ind^{\Sm }_{\vrule height4mm width0pt\hskip-4mm
\StS{m{-}(\pp)}{\pp}}
\hskip-10mm
\Hbm^{q}(\UU^{m}_{m-p,\ldots,m}
)\otimes\alt_{\pp}\Big)\otimes\alt_{m}=\hfill\\
\hskip2cm{}=\ind^{\Sm }_{\vrule height4mm width0pt\hskip-13mm
\StS{m{-}(\pp)}{\pp}}
\hskip-2mm
\Big(\Hbm^{q}(\UU^{m}_{m-p,\ldots,m}
)\otimes\alt_{\pp}\otimes\Res^{\Sm}_{\vrule height4mm width0pt\hskip-6mm
\StS{m{-}(\pp)}{\pp}}
\hskip-10mm
\alt_{m}\quad\Big)
}$$
explique le changement de $\alt_{\pp}$ en $\alt_{m{-}(\pp)}$ dans l'énoncé \ref{theo-approche-i-acyclique}.
 L'égalité $p+q=i+(\mmo)$
découle quant à elle de \ref{iso-tordu}-$(\diamond\diamond)$.

\medskip
(\ref{theo-suite-spectrale-basique-b}) On remarque que dans la mesure où $\Hc(\Fm(\Mg_{\geq0}))=0$, le morphisme de liaison dans la suite longue de cohomologie à support compact \ref{suite-longue-de-base} pour
$\Zg:=\Mg_{\geq0}\times\Delta_{a-1+b}(\Mg_{\geq0}^{a-1}\times \Mg_{>0}^{b})$, à savoir
$$
\Hc(\Delta_{\leq a-1+b}(\Mg_{\geq0}\times\Delta_{a-1+b}(\Mg_{\geq0}^{a-1}\times \Mg_{>0}^{b})))
\to
\Hc(\Delta_{a+b}(\Mg_{\geq0}^{a}\times \Mg_{>0}^{b}))[1]\,,
$$
est un isomorphisme. Comme d'autre part, la projection sur les dernières coordonnées
$$
\pi_{a-1+b}:\Delta_{\leq a-1+b}(\Mg_{\geq0}\times\Delta_{a-1+b}(\Mg_{\geq0}^{a-1}\times \Mg_{>0}^{b}))\onto\Delta_{a-1+b}(\Mg_{\geq0}^{a-1}\times \Mg_{>0}^{b})\,.$$
est un revêtement trivial à $(a-1+b)$ nappes paramétrées par
les applications de $\FFF(a+b,a-1+b)$,
on a un isomorphisme canonique
$$
\Delta_{a-1+b}(\Mg_{\geq0}^{a-1}\times \Mg_{>0}^{b})^{a-1+b}[-1]
\simeq
\Delta_{a+b}(\Mg_{\geq0}^{a}\times \Mg_{>0}^{b})^{a+b}\,,$$
et par induction
$$
\smash{\Hc(\Fb(\Mg_{>0}))^{{(a+b-1)!/ (b-1)!}}}[-a]
\simeq
\Hc(\Delta_{a+b}(\Mg_{\geq0}^{a}\times \Mg_{>0}^{b}))
$$
où \relax{$\displaystyle{(a+b-1)!/(b-1)!}=|\FFF(b,a{+}b)|$}. 

\noindent Par conséquent,
$$\mathalign{\Hc(\UU^{m}_{m-p,\ldots,m})&=&\Hc(\Delta_{m}(\Mg_{\geq0}^{m{-}(\pp)}\times\Mg_{>0}^{\pp}))\hfill\\\noalign{\kern4pt}
&=&
\relax{\bigoplusnl_{\FFF(\pp,m)}}\Hc(\Fg_{\pp}(\Mg_{>0}))[-(m{-}(\pp))]\,,}
$$
et par dualité, l'isomorphisme annoncé dans 
(\ref{theo-suite-spectrale-basique-b})
$$\mathrigid2mu
\Hbm(\UU^{m}_{m-p,\ldots,m})\hf{\varXi^{m}_{\pp}}{\simeq}{0.7cm}
\relax{\bigoplusnl_{\FFF(\pp,m)}}
\Hbm(\Fg_{\pp}(\Mg_{>0}))[-(m{-}(\pp))\,\dMg]\,.\eqno(\dagger)$$

\medskip\binspace2\relspace3
(\ref{theo-suite-spectrale-basique-c})
Comme on a $\Hbm^{i}(\Fg_{\pp}(\Mg_{>0}))=0$ pour tout $i\not\in\iii[0,(\pp)\,\dMg]$, puisque
$\Hc^{i}(\Fg_{\pp}(\Mg_{>0}))=0$ si $i<\pp$ (\ref{pol-poincare}-(\ref{pol-poincare-a})), on comprend par $(\dagger)$ que pour si  $\Hbm^{q}(\UU^{m}_{m-p,\ldots,m})\not=0$, alors
$q\in\iii[(m{-}(\pp))\,\dMg,m\,\dMg]\,,$
d'où (\ref{theo-suite-spectrale-basique-c}). A partir de là,
si l'on fixe $i\in\NN$, la majoration (\ref{theo-suite-spectrale-basique-d}) résulte de l'égalité $ q=i+m{-}(\pp)$.
\enddemo

\subsubsectionline{Termes non nuls de la suite spectrale basique.}
La figure suivante illustre les assertions \ref{theo-suite-spectrale-basique}-(\ref{theo-suite-spectrale-basique-c},\ref{theo-suite-spectrale-basique-d}). Les termes $\IEs(\U^{m})_{r}^{p,q}$ non nuls pour $r\geq1$ sont dans la région hachurée,  et, pour chaque $i\in\NN$ fixé, ceux qui contribuent à $\Hbm^{i}(\Fm(\Mg))$ sont dans la région à la fois hachurée et grisée.
$$\preskip1em
\hss\includegraphics{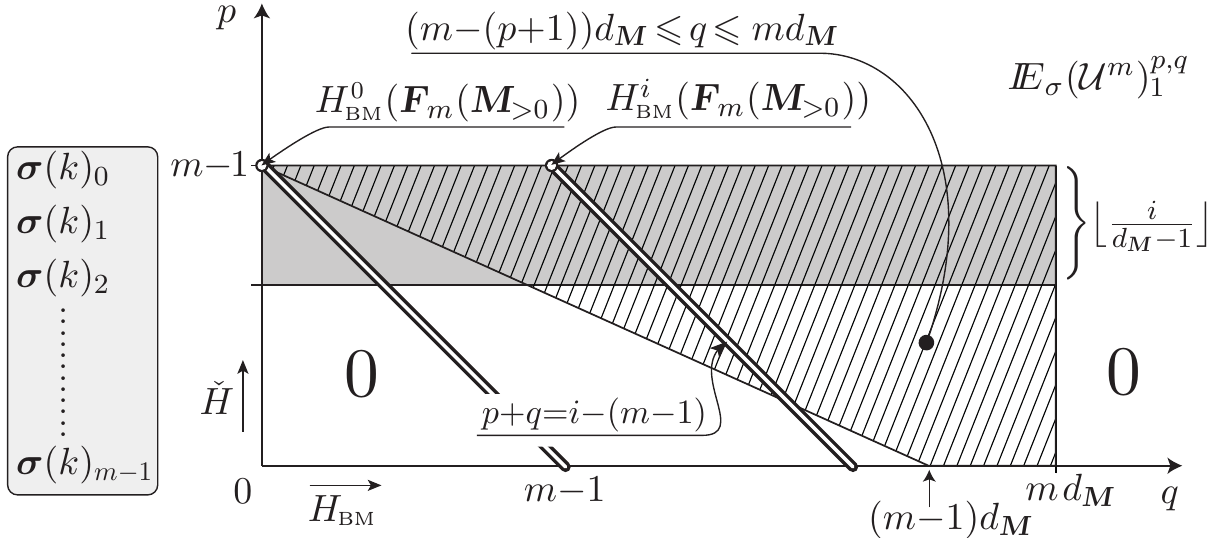}\hss$$
Le petit tableau à gauche rappelle que l'action de $\StS{m{-}(\pp)}{\pp}$ par image-inverse sur les termes de la ligne $p$ est tordue par le caractère $\alt_{m{-}(\pp)}$.

\subsection{Naturalité de la suite spectrale basique pour $\Hbm(\Fm(\Mg))$}\label{naturalite-suite-spectrale-basique}
\subsubsectionnumber
On\label{naturalite-suite-spectrale-basique-preambule} rappelle que l'on note 
$p_m:\Fmm (\Mg)\to\Fg_{m}(\Mg)$ 
la projection sur les $m$ premières coordonnées. 
Nous supposerons que $\Mg$ est une pseudovariété \emph{orientée} de dimension finie $\dMg$ (\footnote{Si dans les précédentes sections cette hypothèse n'était pas vraiment indispensable, maintenant elle l'est car autrement nous ne voyons pas comment donner un sens au morphisme d'image-inverse en cohomologie de Borel-Moore.})
 de sorte que le morphisme d'intégration sur les fibres $p_{m!}:\Hc(\Fmm)\to\Hc(\Fm)$ est défini, ainsi que son dual, le morphisme d'image-inverse $p_{m}^{*}:\Hbm(\Fm)\to\Hbm(\Fmm)$ (\cf \ref{def-action-image-inverse}-(\ref{def-action-image-inverse-a})).
\medskip

Dans cette partie nous préparons le terrain pour l'étude ultérieure en \ref{stabilite-des-familles-i-acycliques} de certaines propriétés asymptotiques des familles des morphismes d'image-inverse en cohomologie de Borel-Moore
$$\bigset\Hbm^{i}(\Fm(\Mg))\hf{ p_{m}^{*}}{}{0.7cm}\Hbm^{i}(\Fmm (\Mg))/_{m\in\NN}\,.$$
Notre but est la construction d'un morphisme de suites spectrales (\cf\ref{morphisme-suistes-spectrales-basiques})
$$
\IEs(q_m^{*}):(\IEs(\U^{m})_r,d_r)\to(\IEs(\U^{\mm})_r,d_r)\,,$$
induisant par passage à la limite, les bi-gradués des morphismes
$p_{m}^{*}$
pour les filtrations du théorème \ref{theo-suite-spectrale-basique}, but qui sera atteint dans le théorème \ref{theo-suite-spectrale-basique-relative}.

\subsubsectionline{Compatibilité des opérations d'intégration sur les fibres.}
En accord\label{compatibilite} avec les notations de \ref{symetries-complexes-cech}, où nous avons introduit l'ouvert 
$$\UU^{\mm}:=\Fmm(\Mg_{\geq0})\moins\Fmm(\Mg_0)\hskip1.8cm$$
et son recouvrement $\U^{\mm}:=\set\UU^{\mm}_{1},\ldots,\UU^{\mm}_{\mm}/$,
on note maintenant
$$\mathalign{\WW^{\mm}&:=&\big(\Fm(\Mg_{\geq0})\times\Mg_{\geq0}\big)\moins
(\Fm(\Mg_0)\times\Mg_{0})\hfill\\\noalign{\kern3pt}
&=&
\big(\UU^{m}\times\Mg_{\geq0}\big)
\cup
\big(\Fm(\Mg_{\geq0})\times\Mg_{>0})\hfill
}
$$
et son recouvrement $\W^{\mm}:=\set\WW^{\mm}_{1},\ldots,\WW^{\mm}_{\mm}/$,
avec
$$\begin{cases}
\WW^{\mm}_{i}:=\UU^{m}_{i}\times\Mg_{\geq0}\,,\hbox{ si $i\leq m$,\quad}\\\noalign{\kern2pt}
\WW^{\mm}_{\mm}:=\Fm(\Mg_{\geq0})\times\Mg_{>0}\,.
\end{cases}
\postskip0pt$$
On a donc:
$$\U^{\mm}=\W^{\mm}\Cap\UU^{\mm}:=\set \WW^{\mm}_{i}\cap \UU^{\mm}=\UU^{\mm}_{i}/\,.$$

Le diagramme suivant, où
les flèches verticales désignent les morphismes de prolongement par zéro, est un morphisme de suites exactes longues:
$$
\xymatrix@C=5mm{
\ar[r]&\Hc(\UU^{\mm})\ar[d]_{\iota_{U!}}\ar[r]&\Hc(\Fmm(\Mg_{\geq0}))\ar[d]\ar[r]&\Hc(\Fmm(\Mg_{0})\ar[d]^{\iota_{\Fg!}}\ar[r]^(0.75){c_{\mm}})&\\
\ar[r]&\Hc(\WW^{\mm})\ar[r]&\Hc(\Fm(\Mg_{\geq0})\times\Mg_{\geq0})\ar[r]&\Hc(\Fm(\Mg_{0})\times\Mg_{0})\ar[r]^(0.85){c'_{\mm}}&\\
}$$
et comme la colonne centrale est nulle, il en résulte le diagramme commutatif où les flèches horizontales sont des isomorphismes:
$$
\xymatrixc{@C=8mm@R=10mm}{
\Hc(\Fmm(\Mg_{0})\ar[d]^{\iota_{\Fg!}}\ar[r]^(0.53){c_{\mm}}_(0.53){\simeq})
&\Hc(\UU^{\mm})\ar[d]^{\iota_{U!}}[1]\\
\Hc(\Fm(\Mg_{0})\times\Mg_{0})\ar[r]^(0.55){c'_{\mm}}_(0.55){\simeq})&\Hc(\WW^{\mm})[1]\\
}\eqno(\dagger)$$

D'autre part, pour $m\geq1$, l'ouvert $\WW^{\mm}\dans\Fm(\Mg_{\geq0})\times\Mg_{\geq0}$ est réunion de deux ouverts $\Hc(\_)$-acycliques
$$
\WW^{\mm}=\WW^{\mm}_{\UU}\cup\WW^{\mm}_{\Fg}
\text{\quad avec\quad }
\begin{cases}
\WW^{\mm}_{\UU}:=\big(\UU^{m}\times\Mg_{\geq0}\big)\\
\WW^{\mm}_{\Fg}:=\big(\Fm(\Mg_{\geq0})\times\Mg_{>0}\big)\end{cases}
\eqno(\dagger\dagger)$$

Le morphisme de liaison dans la suite de Mayer-Vietoris pour la cohomologie à support compact relative à $\set\WW^{\mm}_{\UU},\WW^{\mm}_{\Fg}/$ est donc un isomorphisme:
$$\xymatrix@C=8mm@R=5mm{
\Hc(W^{\mm})\ar[r]^(0.41){c_{\rm MV}}_(0.41){\simeq}&\Hc(\UU^{m}\times\Mg_{>0})[1]\,.}$$

\begin{prop}Soit\label{prop-compatibilite-integrations}
$\Mg$ une pseudovariété \emph{orientée} de dimension $\dMg$. Pour tout $m\geq1$, le
 diagramme suivant où les flèches horizontales sont des isomorphismes, est commutatif au signe près.
$$\preskip4pt\postskip4pt
\xymatrix@C=8mm@R=9mm{
\Hc(\Fmm(\Mg_{0})\ar[d]^{\iota_{\Fg!}}\ar[r]^(0.53){c_{\mm}}_(0.53){\simeq})\xylbl[rd]{\rm(I)}&\Hc(\UU^{\mm})\ar[d]^{\iota_{U!}}[1]\\
\Hc(\Fm(\Mg_{0}))\otimes\Hc(\Mg_{0})
\xylbl[rrd]{\rm(II)}
\ar[d]^{\id\otimes\int_{\Mg_0}}
\ar[r]^(0.6){c'_{\mm}}_(0.6){\simeq})&\Hc(\WW^{\mm})[1]\ar[r]^(0.4){c_{\rm MV}}_(0.4){\simeq}&
\Hc(\UU^{m})\otimes\Hc(\Mg_{>0})[2]
\ar[d]_{\id\otimes\int_{\Mg_{>0}}}\\
\Hc(\Fm(\Mg_0))[-\dMg]\ar[rr]^{c_m[-\dMg]}_{\simeq}&&\Hc(\UU^{m})[1][-\dMg]
}$$
Plus précisément, {\rm(I)} et commutatif et dans {\rm(II)} on a
$$\int_{\Mg_{>0}}(c_{\rm MV}\circ c'_{\mm})(\alpha\otimes\beta)=(-1)^{|\alpha|+|\beta|}\;
c_m\Big(\int_{\Mg_0}\alpha\otimes\beta\Big)
\eqno(\diamond)$$
\end{prop}
\demo\halfdisplayskips   La commutativité de {(I)} a été justifiée dans $(\dagger)$. Pour le sous-diagramme {(II)}, on utilise le fait que les cofaisceaux $\fs\Omegac(\_)$ sont flasques. Quitte à prendre des recouvrements assez fins, un cocycle 
$\omega\in\Omegac(\Fm(\Mg_{0}))\otimes \Omegac(\Mg_{0})$ s'exprime comme la restriction $\sum_{i}\alpha_i\otimes\beta_i$ de 
$$\tilde\omega=\sumnl_{i}\tilde\alpha_i\otimes\tilde\beta_i
\,,\quad\hbox{ avec }\tilde\alpha_i\in\Omegac(\Fm(\Mg_{\geq0}))\text{ et }
\tilde\beta_i\in\Omegac(\Mg_{\geq 0})\,.
$$
On a donc
$
c'_{\mm}\big(\sumnl_{i}\alpha_i\otimes\beta_i\big)=
\sumnl_{i}d\tilde\alpha_i\otimes\tilde\beta_i
+
\sumnl_{i}(-1)^{|\alpha_i|}\tilde\alpha_i\otimes d\tilde\beta_i\,,$
avec 
$
d\tilde\alpha_i\otimes\tilde\beta_i\in\Omegac(\WW^{\mm}_{\UU})$ et
$\tilde\alpha_i\otimes d\tilde\beta_i\in\Omegac(\WW^{\mm}_{\UU})\,.
$
Et alors, par le morphisme de liaison de Mayer-Vietoris $c_{\rm MV}$ relatif à $\set\WW^{\mm}_{\UU},\WW^{\mm}_{\Fg}/$, on a
$$(c_{\rm MV}\circ c'_{\mm})
\big(\sumnl_{i}\alpha_i{\otimes}\beta_i\big)=
\,\sumnl_{i}(-1)^{|\alpha_i|}
d\tilde\alpha_i{\otimes} d\tilde\beta_i\in\Omegac(\UU){\otimes}\Omegac(\Mg_{>0})\,,
$$
de sorte que 
$$\int_{\Mg_{\geq0}}(c_{\rm MV}\circ c'_{\mm})
\big(\sumnl_{i}\alpha_i\otimes\beta_i\big)=
\sumnl_{i}(-1)^{|\alpha_i|+|\beta_i|}\,
d\tilde\alpha_i\int_{\Mg_0} \beta_i\,,
$$
puisque, avec les conventions en cours,
$
\int_{\Mg_{>0}}d\tilde\beta=(-1)^{|\beta|}\int_{\Mg_{0}}\beta\,.$

On conclut par le fait que si $\alpha\in\Z_{\rmc}(\Fm(\Mg_0))$, on a $c_{m}(\alpha)=d\tilde\alpha$.
\enddemo

\begin{rema}[pour le cas $m=0$]L'énoncé de \ref{prop-compatibilite-integrations} doit être modifié, puisqu'alors les morphismes $c_{\rm MV}$ et $c_{m}$ sont nuls. En effet, dans ce cas $\Fm(\Mg_{\geq0})=\Fm(\Mg_0)=\set\pt/$, et alors $\UU^{m}=\emptyset$. Donc, $c_{m}=0$ et le premier terme de l'égalité $(\dagger\dagger)$ est vide, ce qui entraîne la nullité de $c_{\rm MV}$. 
On a aussi $\iota_{\Fg}=\id$, $\iota_{\UU}=\id$ et $c_{\mm}=c'_{\mm}$. Le diagramme de la proposition devient
$$\preskip2pt
\xymatrixc{@C=8mm@R=8mm}{
\Hc(\Mg_{0})
\ar[d]^{\int_{\Mg_0}}
\ar[r]^(0.45){c'_{1}}_(0.45){\simeq})&\Hc(\Mg_{>0})[1]
\ar[dl]+<38pt,2pt>^{\int_{\Mg_{>0}}}\\
\Hc(\pt)[-\dMg]
}
\raise-8pt\hbox{où\quad $\displaystyle
\int_{\Mg_{>0}}c'_{1}(\beta)=(-1)^{|\beta|}
\int_{\Mg_0}\beta$.}
\postskip2em$$ 
\end{rema}

\subsubsectionline{Renormalisation des morphismes de liaison.}Sur\label{renormalisation} le sous-diagram\-me suivant de la proposition \ref{prop-compatibilite-integrations}
$$\def\courbe{\ar`r(11.4,0)|(0.7){\mboxit3{q_{m!}}}`[rdd][rdd]}
\def\Courbe{\ar`r[rr]|(0.7){\mboxit3{q_{m!}}}
`[rdd]
 [rdd]
}
\def\mboxit#1#2{\boxit1{\boxit{#1}{$\scriptstyle#2$}}}
\mathrigid1mu
\initxy
0;<9mm,0mm>:<0mm,10mm>::
\xymatrix@C=9mm@R=10mm{
\Hc(\Fmm(\Mg_{0})\ar[dd]^{\int_{\Mg_0}}\ar[r]|(0.51){\vrule height8pt width0pt\mboxit3{\gamma_{\mm}}}&\Hc(\UU^{\mm})[1]\ar[d]^{\iota_{U!}}
\courbe&
\\
&\Hc(\WW^{\mm})[1]\ar[r]^(0.4){c_{\rm MV}}_(0.4){\simeq}&
\Hc(\UU^{m})\otimes\Hc(\Mg_{>0})[2]
\ar[d]_{\id\otimes\int_{\Mg_{>0}}}\\
\Hc(\Fm(\Mg_0))[-\dMg]\ar[rr]|{\mboxit2{\gamma_m[-\dMg]}}&&\Hc(\UU^{m})[1][-\dMg]
}\endxy$$
on a noté, pour $m\geq1$,
$$q_{m!}:\Hc(\UU^{\mm})\to\Hc(\UU^{m})[-\dMg]\,,\quad
\relax{q_{m!}:=\Big(\id\otimes\int_{\Mg_{>0}}\Big)\circ c_{\rm MV}\circ\iota_{\UU!}}\,,$$
et nous avons remplacé $c_m:\Hc(\Fm(\Mg_0))\to\Hc(\UU^{m})[1]$ par\glossary{${\gamma_{m}(\omega)={(-1)^{|\omega|m+ {\textstyle {m(\mmo)\over2}}d_\Mg}\,c_m(\omega)}}$:renormalisation du morphisme de liaison}
$$\gamma_m:\Hc(\Fm(\Mg_0))\to\Hc(\UU^{\mm})[1]\,,\quad
\gamma_{m}(\omega)={(-1)^{|\omega|m+ {\textstyle {m(\mmo)\over2}}d_\Mg}\,c_m(\omega)}\,.$$

\begin{prop}Soit\label{prop-renormalisation} $\Mg$ une pseudovariété orientée de dimension $\dMg$. 
\mynobreak\begin{enumerate}\nobreak
\item Pour\label{renormalisation-a} $m\geq1$, le diagramme suivant de
morphismes compatibles aux actions des groupes symétriques par image-directe, est un diagramme commutatif,
$$\xymatrix@C=1cm@R=10mm{
\Hc(\Fmm(\Mg_0))\ar[d]_{p_{m!}}^{[-\dMg]}\ar[r]^(0.53){\gamma_{\mm}}_(.53){\simeq}
&\Hc(\UU^{\mm})\ar[d]^{q_{m!}}_{[-\dMg]}[1]\\
\Hc(\Fm(\Mg_0))\ar[r]^(0.53){\gamma_{m}}_(.53){\simeq}&\Hc(\UU^{m})[1]\,.
}$$
\item Le\label{renormalisation-b} diagramme suivant
obtenu en dualisant et concaténant les diagrammes de {\rm (\ref{renormalisation-a})} et où l'on a noté $\Fg_{i}$ pour $\Fg_{i}(\Mg_0)$, est un diagramme commutatif,
$$\hss\mathrigid0mu
\let\Hc\Hbm
\xymatrix@C=4mm@R=10mm{
\Hc(\Fg_{1})\ar[r]_(.47){p_{1}^{*}}&\Hc(\Fg_{2})\ar[r]_(0.62){p_{2}^{*}}&
\cdots\ar[r]
&\Hc(\Fm)\ar[r]_(0.45){p_{m}^{*}}&\Hc(\Fmm)\ar[r]&\\
\Hc(\UU^{1})\ar[r]_(.43){q_{1}^{*}}\ar[u]^{\gamma_{1}^{*}}_{[0]}&
\Hc(\UU^{2})[1]\ar[r]_(0.68){q_{2}^{*}}\ar[u]^{\gamma_{2}^{*}}_{[-1]}&
\cdots\ar[r]&\Hc(\UU^{m})[\mmo]\ar[r]_(0.5){q_{m}^{*}}\ar[u]^{\gamma_{m}^{*}}_{[1{-}m]}&
\Hc(\UU^{\mm})[m]\ar[r]\ar[u]^{\gamma_{\mm}^{*}}_{[-m]}&
}\hss$$
Les morphismes y sont compatibles aux actions des groupes symétriques lorsque l'on tord $\Hbm(\UU^{m})$ par $\alt_{m}$ {\rm(\cf\ref{iso-tordu}-($\diamond\diamond$))}.
Les suites horizontales sont alors des $\FI$-modules (cf. \ref{familles-denombrables}) et la famille $\set \gamma_{m}^{*}[1-m]/$ est un isomorphisme de $\FI$-modules.
\end{enumerate}
\end{prop}
\demo (\ref{renormalisation-a}) Notons $Q(m)={m(\mmo)\over2}\dMg$. Compte tenu de l'égalité $(\diamond)$ de \ref{prop-compatibilite-integrations}, on a
$$\preskip1pt
\mathalign{(q_{m!}\circ\gamma_{\mm})(\omega)&=&
(-1)^{|\omega|m+Q(\mm)}\,
c_{m}(p_{m!}(\omega))\hfill\\
&=&(-1)^{(|\omega|-\dMg)m+Q(m)}c_{m}(p_{m!}(\omega))=
(\gamma_{m}\circ p_{m!})(\omega)
}
$$
La compatibilité par rapport aux actions des groupes symétriques est claire.

(\ref{renormalisation-b}) On verra dans \ref{familles-denombrables} que la condition pour que les suites horizontales soient des $\FI$-modules est que pour tous $m< n$, on ait
$$(p_{n-1}^{*}\circ\cdots\circ p_{m}^{*})(\Hbm(\Fm(\Mg_0)))
\dans\smash{\Hbm(\Fn(\Mg_0))^{\1_{m}\times\S_{n-m}}}$$
et de même pour la deuxième ligne. C'est évident pour la première ligne, et donc aussi pour la seconde qui en est isomorphe.
\enddemo

\subsubsectionline{L'image-inverse sur les suites spectrales basiques.}
Nous allons\label{morphisme-suistes-spectrales-basiques} relever le morphisme $q_{m}^{*}:\Hbm(\UU^{m})\to\Hbm(\UU^{\mm})[1]
$, pour $m\geq 1$, de \ref{prop-renormalisation} en un morphisme
 $\IE(q^{*}_{m}):\IE(\U^{m})\to\IE(\U^{\mm})[1,0]$ de suites spectrales basiques, de sorte que nous aurons le diagramme commutatif
$$\preskip2pt\postskip3pt
\xymatrix@R=6mm@C=1.75cm{
\Hbm(\UU^{m})
\ar[r]|(0.44){\ q_{m}^{*}\ }&\Hbm(\UU^{\mm})[1]\\
\IE(\U^{m})\ar@{=>}[u]\ar[r]|(0.43){\ \IE( q_{m}^{*})\ }&\IE(\U^{\mm})[1,0]\ar@{=>}[u]
}$$
répondant ainsi à la principale motivation de cette section \ref{naturalite-suite-spectrale-basique} (\cf\ref{naturalite-suite-spectrale-basique-preambule}).

\smallskip

Les morphismes $q_{m}^{*}$ sont les duaux des morphisme $q_{m!}$
qui on été définis comme la composée de trois morphismes
$$
\xymatrixc{@R=6mm}{
&
\Hc(\WW^{\mm})\ar[r]^(0.42){c_{\rm MV}}&
\Hc(\UU^{m}\times\Mg_{>0})\ar[rd]+<-8mm,3mm>^(0.65){\quad p_{m!}:=\int_{\Mg>0}}[1]&
\\
\Hc(\UU^{\mm})\ar[ru]^{\iota_{\UU!}}\ar[rrr]|{\ q_{m!}\ }
&&&\Hc(\UU^{m})[-\dMg]
}$$
Dans ce qui suit, nous montrerons que le dual de chacun de ces morphismes admet un relèvement spectral.
Rappelons maintenant quelques notations.

\begin{itemize}\itemsep8pt\mou
\item L'ouvert $\UU^{\mm}\dans\Fmm(\Mg_{\geq0})$ est muni du recouvrement (\ref{symetries-complexes-cech})
$$\U^{\mm}=\set \UU^{\mm}_{1},\ldots,\UU^{\mm}_{\mm}/
\text{\quad où\quad }\big(\cl x\in\UU^{\mm}_{i}\Leftrightarrow x_i\not\in\Mg_{0}\big).$$

\item L'ouvert $\WW^{\mm}\dans\Fm(\Mg_{\geq0})\times\Mg_{\geq0}$ est muni du recouvrement (\ref{compatibilite})
$$\mathrigid2mu
\hss\W^{\mm}:=\set\WW^{\mm}_{1},\cdots,\WW^{\mm}_{\mm}/
\text{\ où\  }
\begin{cases}\noalign{\kern-4pt}
\WW^{\mm}_{i}:=\UU^{m}_{i}\times\Mg_{\geq0}\,,\hbox{ si $i\leq m$,\quad}\\
\WW^{\mm}_{\mm}:=\Fm(\Mg_{\geq0})\times\Mg_{>0}\,.\\\noalign{\kern-2pt}
\end{cases}
\hss$$
On décompose $\WW^{\mm}$ en réunion de deux ouverts $\Hc(\_)$-acycliques
$$\postskip-1ex
\WW^{\mm}=\WW^{\mm}_{\UU}\cup\WW^{\mm}_{\Fg}
\text{ avec }
\begin{cases}\noalign{\kern-4pt}
\WW^{\mm}_{\UU}:=\big(\UU^{m}\times\Mg_{\geq0}\big)\\
\WW^{\mm}_{\Fg}:=\big(\Fm(\Mg_{\geq0})\times\Mg_{>0}\big)\mrlap{\,.}\\\noalign{\kern-2pt}
\end{cases}
\eqno(\dagger\dagger)$$
\end{itemize}

\soustitreline{Relèvement spectral de $\iota_{\UU}^{*}$.}Comme le recouvrement $\U^{\mm}$ est la trace sur $\UU^{\mm}$ du recouvrement $\W^{\mm}$, \idest
$$\U^{\mm}=\W^{\mm}\Cap\UU^{\mm}:=\set \WW^{\mm}_{i}\cap \UU^{\mm}=\UU^{\mm}_{i}/\,,$$
la naturalité de bicomplexes de cochaînes de \vCech-Borel-Moore
vis-à-vis des restrictions ouvertes, induit un morphisme de bicomplexes de degré (0,0)
$$(\iota_{\UU}^{*})_{\bullet,*}:\vC_{<}(\W^{\mm})^{\bullet;*}_{\bm}\to\vC_{<}(\U^{\mm})^{\bullet;*}_{\bm}
$$
et donc un morphisme de suites spectrales
$$(\iota^{*}_{\UU,r}):(\IE(\W^{\mm})_r,d_r)\to(\IE(\U^{\mm})_r,d_r)\,.$$

\def\theotitlestyle{\sl\bfseries\boldmath}
\begin{prop*}[A (relèvement de $i_{\UU}^{*}$)]Le morphisme
$$\preskip0.5ex
\xymatrix@C=1.5cm{\IE(\W^{\mm})_{1}^{p,q}=\vC^{p}(\W;\H_{\bm}^{q}(\_))
\ar[r]^{(\iota_{\UU,1}^{*})_{p,q}}&
\IE(\U^{\mm})_{1}^{p,q}=\vC^{p}(\U;\H_{\bm}^{q}(\_))\,,
}
$$
en fonction de $p\geq0$ et de 
$\mathrigid2mu
1\leq i_0<\cdots<i_p\leq m+1$, prend les valeurs suivantes.

\smallskip
\noindent$\triangleright$ Si $p=0$, on a $(\iota_{\UU,1}^{*})_{0,q}=0$ puisque les $\WW^{\mm}_{i}$ sont $\Hbm(\_)$-acycliques.
$$
\xymatrix@R=5mm@C=20mm{
\big(\IE(\W^{\mm})^{0,q}_{1}\big)_{i_0}\ar[r]_{(\iota_{\UU,1}^{*})_{0,q}=0}\aregal[d]&
\big(\IE(\U^{\mm})^{0,q}_{1}\big)_{i_0}\aregal[d]\\
\Hbm^{q}(\WW^{\mm}_{i_0})\mrlap{=0}&
\Hbm^{q}(\UU^{\mm}_{i_0})
}
$$

\smallskip
\noindent$\triangleright$ Si $p>0$ et $i_p<m+1$, on a aussi $(\iota_{\UU,1}^{*})_{p,q}=0$ puisque 
$$\WW^{\mm}_{i_0,\ldots,i_p}=\UU^{m}_{i_0,\ldots,i_p}\times\Mg_{\geq0}
\postskip0pt$$
est $\Hbm(\_)$-acyclique.
$$\preskip0pt
\xymatrix@R=5mm@C=20mm{
\big(\IE(\W^{\mm})^{p,q}_{1}\big)_{i_0,\ldots,i_p}\ar[r]_{(\iota_{\UU,1}^{*})_{p,q}=0}\aregal[d]&
\big(\IE(\U^{\mm})^{p,q}_{1}\big)_{i_0,\ldots,i_p}\aregal[d]\\
\Hbm^{q}(\UU^{m}_{i_0,\ldots,i_p}\times\Mg_{\geq0})\mrlap{=0}&
\Hbm^{q}(\UU^{\mm}_{i_0,\ldots,i_p})
}
$$

\smallskip
\noindent$\triangleright$ Si $p>0$ et $i_p=m+1$, on a
$$\WW^{\mm}_{i_0,\ldots,i_p=m+1}=\UU^{m}_{i_0,\ldots,i_{p-1}}\times\Mg_{>0}\,,
$$
et $\iota_{\UU,1}^{*}$ s'identifie à la restriction de 
$\big(\UU^{m}_{i_0,\ldots,i_{p-1}}\times\Mg_{>0}\big)$ à 
$\UU^{\mm}_{i_0,\ldots,i_{p-1},m+1}$
$$\preskip3pt
\xymatrix@R=5mm@C=20mm{
\big(\IE(\W^{\mm})^{p,q}_{1}\big)_{i_0,\ldots,i_p=m+1}\ar[r]_{(\iota_{\UU,1}^{*})_{p,q}}\aregal[d]&
\big(\IE(\U^{\mm})^{p,q}_{1}\big)_{i_0,\ldots,i_p=m+1}\aregal[d]\\
\Hbm^{q}(\UU^{m}_{i_0,\ldots,i_{p-1}}\times\Mg_{>0})\ar[r]^{\hbox{\scriptsize restriction}}&
\Hbm^{q}(\UU^{\mm}_{i_0,\ldots,i_{p-1},m+1})
}
$$
\end{prop*}

\soustitreline{Relèvement spectral de $c_{\rm MV}^{*}$.}On commence par un scholie sur les bicomplexes $\vCm^{\bullet,*}(\U;\fsOmega^{*})$ de cochaînes de \vCech\ ordonnées, relatives à un recouvrement $\U=\set U_1,\ldots, U_{\mm}/$ et
à valeurs dans le complexe de faisceaux de cochaînes de Borel-moore $(\fsOmega^{*}_{\bm},d_*)$ (\ref{ss-coh-supp-comp-recouvrement}), ce que nous notons
$$\vC^{\bullet,*}(\U):=\vCm^{\bullet,*}(\U;\fsOmega^{*}_{\bm})\,.$$

Soient $\U':=\set U_1,\ldots,U_m/$ et $\U'':=\set U_{\mm}/$ dont on remarquera tout de suite que $\vC^{p}(\U')=0$ si $p\geq m$, et 
$\vC^{p}(\U'')=0$ si $p\geq1$. 

Définissons les restrictions de cochaînes de \vCech\ $\rho',\rho''$:
$$
\mathalign{
\mllap{\rho:}\vC^{\bullet}(\U)&\hf{}{}{7mm}&\vC^{\bullet}(\U')\oplus\vC^{\bullet}(\U'')\\
\omega&\hfto{}{}{7mm}&(\rho'(\omega),\rho''(\omega))
}\postskip0pt$$
par
$$\preskip-0.5ex\begin{cases}\noalign{\kern-4pt}
\rho'(\omega)_{i_0,\ldots,i_p}=\omega_{i_0,\ldots,i_p}\in\Omega^{*}(U_{i_0,\ldots,i_p})\,,
\hbox{\quad où $1\leq i_0,\ldots,i_p\leq m$,}\\
\rho''(\omega)_{0}=\omega_{\mm}\in\Omega^{*}(U_{\mm})\,.\\\noalign{\kern-2pt}
\end{cases}
$$

Posons ensuite $\U''':=\set U_1\cap U_{\mm},\ldots,U_m\cap U_{\mm}/$, et définissons 
$$\mathalign
{
\mllap{c:}\vC^{\bullet}(\U''')
&\hf{}{}{7mm}&
\vC^{\bullet+1}(\U)
}
\postskip0pt\eqno(\diamond)$$
par
$$
c(\omega)_{i_0,\ldots,i_p}=\begin{cases}\noalign{\kern-4pt}
\omega_{i_0,\ldots,i_{p-1}}\,,\hbox{ si $i_p=m+1$,}\\
0\,,\hbox{ sinon.}\\\noalign{\kern-2pt}
\end{cases}$$

\def\theotitlestyle{\sl\bfseries}
\begin{prop*}
La suite 
$$
0\to\vC^{\bullet}(\U''')\too^{c}\vC^{\bullet+1}(\U)\hf{\rho}{(\rho',\rho'')}{1cm}
\vC^{\bullet+1}(\U')\oplus\vC^{\bullet+1}(\U'')\to0
\eqno(\ddagger)$$
est une suite exacte courte de morphismes de bicomplexes.
Les morphismes qu'elle induit sur les cohomologies des complexes simples associés coïncident alors avec les morphismes de la suite exacte longue de Mayer-Vietoris pour le recouvrement 
$\cup\,\U=(\cup\,\U')\cup(\cup\,\U'')$, soit
$$\hss\mathrigid2mu
\to\Hbm^{*}(\cup\,\U''')\too^{c_{\rm MV}}\Hbm^{*+1}(\cup\,\U)
\too^{\rho}\Hbm^{*+1}(\cup\,\U')\oplus\Hbm^{*+1}(\cup\,\U'')\to\Hbm^{*+1}(\cup\,\U''')\to
\hss$$

En particulier, le morphisme de liaison $c_{\rm MV}$ est l'aboutissement du morphisme des suites spectrales
$(c_r):(\IE(\U'''),d_r)\to(\IE(\U''')[1,0],d_r)$
 induit par le morphisme de bicomplexes $(\diamond)$ et dont la restriction aux termes $\IE_1$ est
$$
\xymatrix@R5mm
{
\IE(\U''')^{p,q}_{1}\mrlap{=\smash{\displaystyle\bigoplusnl_{i_0<\cdots<i_p}}\Hbm^{q}(U_{i_{0},\ldots,i_{p}}\cap U_{\mm})\ar[d]_{(c_1)_{p,q}}}&\kern4cm\\
\IE(\U)^{\pp,q}_{1}\mrlap{=\smash{\displaystyle\bigoplusnl_{i_0<\cdots<i_p}}
\Hbm^{q}(U_{i_{0},\ldots,i_{\pp}})}&
}
\postskip1.em$$
où, si $\omega\in(\IE(\U''')^{p,q}_{1})_{i_0,\ldots,i_p}$, on a
$$
((c_1)_{p,q}(\omega))_{k_0,\ldots,k_{\pp}}=\begin{cases}\noalign{\kern-4pt}
\omega\,,\hbox{ si $(k_0,\ldots,k_{\pp})=(i_0,\ldots,i_p,m+1)$}\\
0\,,\hbox{sinon.}\\\noalign{\kern-2pt}
\end{cases}$$
\end{prop*}
\demo \parskip1ex\mou
Le fait que $\rho$ est un morphisme de bicomplexes surjectif de bidegré $(0,0)$ est immédiat. 
Son noyau $\ker^{\bullet}(\rho)$ est le sous-bicomplexe de $\vC^{\bullet}(\U)$  vérifiant
$$
\ker^{0}(\rho)=0\text{\quad et\quad}
(\forall p>0)\big(\omega\in\ker^{p}(\rho)\Leftrightarrow
\omega_{i_0,\ldots,i_p}=0\,,\hbox{ si $i_p\leq m$}\big)\,,
$$
où l'on reconnaît l'image de $c:\vC^{\bullet-1}(\U''')\to
\vC^{\bullet}(\U)$. La suite $(\dagger)$ est donc bien une suite exacte courte de bicomplexes.

Notons 
$$
\let\quad\, U:=\cup\,\U\,,\qquad
U':=\cup\,\U'\,,\qquad
U'':=\cup\,\U''\,,\qquad
U''':=\cup\,\U'''\,.$$
La suite courte de Mayer-Vietoris de bicomplexes
$$0\to\vC^{\bullet}(\U)\too^{\alpha}
\vC^{\bullet}(\U\Cap U')
\oplus
\vC^{\bullet}(\U\Cap U'')
\too^{\beta}
\vC^{\bullet}(\U\Cap U''')
\to0
$$
est exacte puisque le faisceau $\fsOmega^{*}_{\bm}$ est flasque.
Comme on a $\U\Cap U'\cont\U'$ et $\U\Cap U''\cont\U''$, on dispose de morphismes de restriction de cochaînes de \vCech\ $r',r''$ et du diagramme commutatif
$$\xymatrix@C=4mm@R=10mm
{
&\mllap{0\to}\vC^{\bullet}(\U)\ar[r]^(.25){\alpha}\ar[d]^(0.45){\id}&
\vC^{\bullet}(\U\Cap U')
\oplus
\vC^{\bullet}(\U\Cap U'')
\ar[r]^(.65){\beta}\ar@<8ex>[d]_(0.45){r''}\ar@<-6.3ex>[d]^(0.45){r'}&
\vC^{\bullet}(\U\Cap U''')
\mrlap{\to0}\\
\mllap{0\to}
\vC^{\bullet-1}(\U''')\ar[r]^(.57){c}&
\vC^{\bullet}(\U)\ar[r]^(.3){\rho}&
\vC^{\bullet}(\U')\quad\oplus\quad\vC^{\bullet}(\U'')\mrlap{\to0}
}$$
où $r'$ et $r''$ induisent des quasi-isomorphismes au niveau des complexes simples associés. 

Il en résulté un isomorphisme canonique en cohomologie
$$\xi_{*}:h^{*}\tot\big(\vC^{\bullet,*}(\U\Cap U''')\big)
\to
h^{*}\tot\big(\vC^{\bullet,*}(\U''')
\big)
$$
rendant commutatif le diagramme
$$\xymatrix@R=10mm{
h^{*-1}\tot\big(\vC^{\bullet,*}(\U\Cap U''')\big)\ar[d]^{\xi_{*-1}}\ar[r]^(.55){\gamma}\xylbl[dr]{\hbox{$\bigoplus$}}&
h^{*}\tot\big(\vC^{\bullet,*}(\U)\big)\ar[d]^{\id}\\
h^{*-1}\tot\big(\vC^{\bullet,*}(\U'''))\ar[r]^(.52){c}&
h^{*}\tot\big(\vC^{\bullet,*}(\U)\big)\\
}$$
où $\gamma$ correspond ou morphisme de liaison $c_{\rm MV}:\Hbm(U''')\to\Hbm(U)[1]$ de la suite longue de Mayer-Vietoris pour le recouvrement $U= U'\cup U''$. A partir de là la suite de la proposition est claire.
\enddemo

Cette proposition, appliquée au recouvrement $\W^{\mm}$, fournit le relèvement spectral du morphisme de liaison $c_{\rm MV}^{*}:\Hbm(\UU^{m}\times\Mg_{>0})\to\Hbm(\WW^{\mm})[1]$ 

\def\theotitlestyle{\sl\bfseries\boldmath}
\begin{prop*}[B (relèvement de $c_{\rm MV}^{*}$)]
On munit $\WW^{\mm}$ du recouvrement $\W^{\mm}$, et $\UU^{m}\times\Mg_{>0}$ du recouvrement $\U^{m}\times\Mg_{>0}$ et l'on note
$(\IE(\W^{\mm})_{r},d_r)$ et $(\IE(\U^{m}\times\Mg_{>0})_{r},d_r)$ les suites spectrales correspondantes.
Alors, il existe un morphisme de suites spectrales
$$(c_{\rm MV,r}^{*}):
(\IE(\U^{m}\times\Mg_{>0})_{r},d_r)
\to
(\IE(\W^{\mm})_{r},d_r))[1,0]
$$
convergeant vers le bigradué du morphisme de liaison
$$c_{\rm MV}^{*}:\Hbm(\U^{m}\times\Mg_{>0})\to
\Hbm(\WW^{\mm})[1]\,.$$
L'action de $(c_{\rm MV,r}^{*})$ sur $\IE_1$ est donnée par le plongement canonique
$$\def\ind#1#2{_{\hbox to#1{\hss$\scriptstyle#2$\hss}}}
\def\mbx#1#2{\hbox to#1{$#2$\hss}}
\xymatrix@R=7mm
{
\IE(\U^{m}\times\Mg_{>0})^{p,q}_{1}\ar[d]_{(c_{\rm MV,1})_{p,q}}
\ar@{=}[r]&
\smash{\displaystyle\bigoplus\ind{2cm}{1\leq i_0<\cdots<i_p\leq m}}
\mbx{2cm}{\Hbm^{q}(\UU^{m}_{i_{0},\ldots,i_{p}}\times\Mg_{>0} )}\ar@{=}@<8ex>[d]\\
\IE(\W^{\mm})^{\pp,q}_{1}\ar@{<-^)}[r]&
\displaystyle\smash{\bigoplus\ind{2cm}{1\leq i_0<\cdots<i_p\leq m}}
\mbx{2cm}{\Hbm^{q}(\WW^{\mm}_{i_{0},\ldots,i_{p},m+1})}&
}
\postskip1.ex$$
\end{prop*}

\soustitreline{Relèvement spectral de $\big(\int_{\Mg>0}\big){}^{*}$.}On commence par rappeler que le morphisme d'intégration sur les fibres
a un sens déjà au niveau des complexes de cochaînes à support compact:
$$p_{m!}=\int_{\Mg>0}:\Omegac(\_)\otimes\Omegac(\Mg_{>0})\to\Omegac(\_)[-\dMg{-}1]$$
et c'est un morphisme compatible à l'opération de cobord. Par dualité et naturalité vis-à-vis des inclusions ouvertes, on obtient le morphisme de bicomplexes de cochaînes de \vCech-Borel-Moore de degré $[0,0]$
$$(p_m^{*})_{\bullet,*}:\vC(\U^{m})^{\bullet,*}_{\bm}\to
\vC(\U^{m}\times\Mg_{>0})^{\bullet,*}_{\bm}\,.\eqno(\ddagger)$$
Le théorème suivant et immédiat.

\def\theotitlestyle{\sl\bfseries\boldmath}
\begin{prop*}[C (relèvement de $p_m^{*}$)]Par passage aux complexes simples associés, le morphisme de bicomplexes $(\ddagger)$ induit le morphisme image-inverse
$$p_m^{*}:\Hbm(\U^{m})\to\Hbm(\U^{m}\times\Mg_{>0})\,.$$
Il induit également un morphisme de suites spectrales
$$(p_{m,r}^{*}):(\IE(\U^{m})_r,d_r)\to(\IE(\U^{m}\times\Mg_{>0})_r,d_r)\,.$$
L'action de $(p_{m,r}^{*})$ sur $\IE_1$ est  donnée par les morphismes image-inverse:
$$\def\ind#1#2{_{\hbox to#1{\hss$\scriptstyle#2$\hss}}}
\def\mbx#1#2{\hbox to#1{$#2$\hss}}
\xymatrix@R=7mm
{
\IE(\U^{m})^{p,q}_{1}\ar[d]_{(p_{m,1}^{*})_{p,q}}
\ar@{=}[r]&
\smash{\displaystyle\bigoplus\ind{2cm}{1\leq i_0<\cdots<i_p\leq m}}
\mbx{2cm}{\Hbm^{q}(\UU^{m}_{i_{0},\ldots,i_{p}})}\ar@<7ex>[d]^{p_m^{*}}\\
\IE(\U^{m}\times\Mg_{>0})^{p,q}_{1}\ar@{=}[r]&
\displaystyle\smashtop{\bigoplus\ind{2cm}{1\leq i_0<\cdots<i_p\leq m}}
\mbx{2cm}{\Hbm^{q}(\UU^{m}_{i_{0},\ldots,i_{p}}\times\Mg_{>0} )}&
}
\postskip1.em$$
\end{prop*}

\soustitreline{Conclusion.}La composition des morphismes des suites spectrales des propositions A\,B\,C donne le morphisme de suites spectrales annoncé dans \ref{morphisme-suistes-spectrales-basiques}
$$\IE(q_m^{*}):(\IE(\U^{m})_r,d_r)\to(\IE(\U^{\mm})_r,d_r)[1,0]\,.$$
Le diagramme commutatif suivant traque sa valeur sur les termes $\IE_1$
$$
\def\ind#1#2{_{\hbox to#1{\hss$\scriptstyle#2$\hss}}}
\def\mbx#1#2{\hbox to2.3cm{$#2$\hss}}
\initxy
0;<9mm,0mm>:<0mm,10mm>::\xymatrix
{
\IE(\U^{m})^{p,q}_{1}\ar[d]_{(p_{m,1}^{*})_{p,q}}
\ar@{=}[r]&
\smash{\displaystyle\bigoplus\ind{2cm}{1\leq i_0<\cdots<i_p\leq m}}
\mbx{2cm}{\Hbm^{q}(\UU^{m}_{i_{0},\ldots,i_{p}})}\ar@<7ex>[d]^{p_m^{*}}
\hskip-0mm\ar`(9.2,0)
   `[ddd]|{\vrule height4mm depth3mm width0pt p_m^{*}}
    [ddd]+<83pt,0pt>
          \\
\IE(\U^{m}\times\Mg_{>0})^{p,q}_{1}\ar[d]_{(c_{\rm MV,1})_{p,q}}
\ar@{=}[r]&
\smash{\displaystyle\bigoplus\ind{2cm}{1\leq i_0<\cdots<i_p\leq m}}
\mbx{2cm}{\Hbm^{q}(\UU^{m}_{i_{0},\ldots,i_{p}}\times\Mg_{>0} )}\ar@{=}@<7ex>[d]\\
\IE(\W^{\mm})^{\pp,q}_{1}\ar@{<-^)}[r]\ar[d]_{(\iota_{\UU,1}^{*})_{p,q}}&
\smash{\displaystyle\bigoplus\ind{2cm}{1\leq i_0<\cdots<i_p\leq m}}
\mbx{2cm}{\Hbm^{q}(\WW^{\mm}_{i_{0},\ldots,i_{p},m+1})}\ar@<7ex>[d]^{\iota_{\UU}^{*}}\\
\IE(\U^{\mm})^{\pp,q}_{1}\ar@{<-^)}[r]&
\smash{\displaystyle\bigoplus\ind{2cm}{1\leq i_0<\cdots<i_p\leq m}}
\mbx{2cm}{\Hbm^{q}(\UU^{\mm}_{i_{0},\ldots,i_{p},m+1})}
}
\endxy
\postskip1.25em$$
Toutes ces observations 
conduisent à l'énoncé suivant.\killline

\begin{prop}Soit\label{prop-naturalite-suite-spectrale-basique} $\Mg$ une pseudovariété orientée de dimension $\dMg$.
Il existe un morphisme de suites spectrales basiques
$$\IE(q_m^{*}):(\IE(\U^{m})_r,d_r)\to(\IE(\U^{\mm})_r,d_r)[1,0]\,,$$
qui est compatible aux actions des groupes symétriques
et qui rend le diagramme suivant commutatif {\rm (\footnote{On rappelle que
$\gamma_m^{*}$ est la renormalisation donnée dans \ref{renormalisation} du morphisme de liaison
$c_m^{*}:\Hbm(\UU^{m})\to\Hbm(\Fm(\Mg_{0}))[1{-}m]$. On a $\gamma_{m}(\omega)={(-1)^{|\omega|m+ {\textstyle {m(\mmo)\over2}}d_\Mg}\,c_m(\omega)}$.}).}
$$\preskip5pt\postskip5pt
\xymatrix@C=1.75cm@R=10mm{
\Hbm(\UU^{m})
\ar[r]|(0.44){\ q_{m}^{*}\ }&\Hbm(\UU^{\mm})[1]
\\
\IE(\U^{m})\ar@{=>}[u]\ar[r]|(0.43){\ \IE( q_{m}^{*})\ }&\IE(\U^{\mm})[1,0]\ar@{=>}[u]
}$$
De plus, l'action de $\IE(q_{m}^{*})$ sur les termes $\IE(\_)_1$ est donnée par
$$\preskip5pt\postskip3pt
\def\ind#1#2{_{\hbox to#1{\hss$\scriptstyle#2$\hss}}}
\def\mbx#1#2{\hbox to1.2cm{$#2$\hss}}
\quad\quad
\xymatrixc{@C1.5cm@R=10mm}
{
\IE(\U^{m})^{p,q}_{1}
\ar@{=}[d]
\ar[r]^{\IE(q_m^{*})_{1}^{p,q}}
&\IE(\U^{\mm})^{\pp,q}_{1}\ar@{<-_)}[d]\\
\smash{\displaystyle\bigoplus\ind{1cm}{1\leq i_0<\cdots<i_p\leq m\ }}
{\Hbm^{q}(\UU^{m}_{i_{0},\ldots,i_{p}})}\ar[r]^(0.45){\bigoplus p_m^{*}}
&
\smash{\displaystyle\bigoplus\ind{1cm}{1\leq i_0<\cdots<i_p\leq m}}
{\Hbm^{q}(\UU^{\mm}_{i_{0},\ldots,i_{p},m+1})}\\
}
$$
\end{prop}
\subsection{Relèvement spectral de $p_m^{*}:\Hbm^{i}(\Fm(\Mg))\to\Hbm^{i}(\Fmm(\Mg))$}\label{suivi-spectrale-basique}Si nous combinons la proposition précédente \ref{prop-naturalite-suite-spectrale-basique} aux théorèmes \ref{prop-renormalisation} et \ref{theo-suite-spectrale-basique}, nous obtenons les suites de représentations de groupes symétriques
$$\mathrigid0mu
\def\gamm#1#2{\ar[u]^{\gamma_{#1\ }^{*}}_{\ [#2]}|{\xysimeq}}
\let\Hc\Hbm
\xymatrix@C=4mm@R=7mm{
\Hc(\Fg_{1})\ar[r]_(.47){p_{1}^{*}}&\Hc(\Fg_{2})\ar[r]&
\cdots\ar[r]
&\Hc(\Fm)\ar[r]_(0.45){p_{m}^{*}}&\Hc(\Fmm)\ar[r]&\\
\Hc(\UU^{1})\ar[r]_(.43){q_{1}^{*}}\gamm{1}{0}&
\Hc(\UU^{2})[1]\ar[r]\gamm{2}{-1}&
\cdots\ar[r]&\Hc(\UU^{m})[\mmo]\ar[r]_(0.5){q_{m}^{*}}\gamm{m}{1{-}m}&
\Hc(\UU^{\mm})[m]\ar[r]\gamm{\mm}{-m}&\\
\IEs(\U^1)\ar[r]_(0.4){\vrule height11pt width0pt
\IEs(q_1^{*})}\ar@{=>}[u]&
\IEs(\U^2)[1,0]\ar[r]\ar@{=>}[u]&
\cdots\ar[r]&
\IEs(\U^m)[m{-}1,0]\ar[r]_{\vrule height11pt width0pt
\IEs(q_m^{*})}\ar@{=>}[u]&
\IEs(\U^{\mm})[m,0]\ar[r]\ar@{=>}[u]&
}$$
où la ligne centrale doit être tordue par les caractères $\alt_{m}$ correspondants.
L'énoncé suivant est alors corollaire de \ref{prop-naturalite-suite-spectrale-basique} et \ref{theo-suite-spectrale-basique}.\killline

\begin{theo}Soit\label{theo-suite-spectrale-basique-relative} $\Mg$ une pseudovariété orientée de dimension $\dMg$.
\mynobreak\begin{enumerate}\nobreak\itemsep0pt\parskip0pt
\item Le\label{theo-suite-spectrale-basique-relative-a} morphisme de suites spectrales basiques
$$\IE(q_m^{*}):(\IE(\U^{m})_r,d_r)\to(\IE(\U^{\mm})_r,d_r)[1,0]\,,$$
définit, pour chaque $i\in\NN$, chaque couple $(p,q)$ tel que $q=i+(m{-}(\pp))$ et chaque $r\geq1$, une suite indexée par $m\in\NN$
$$\IEs(q_m^{*})_{r}^{p,q}:\IEs(\U^{m})^{p,q}_{r}\to\IEs(\U^{\mm})^{\pp,q}_{r}\,,$$
qui est un $\FI$-module. La somme de ces suites pour chaque $r\geq1$
converge vers le bi-gradué du morphisme de $\Sm $-modules
$$p_m^{*}:\Hbm^{i}(\Fm(\Mg))\to\Hbm^{i}(\Fmm(\Mg))\postdisplaypenalty10000$$
pour les filtrations régulières induites par les isomorphismes 
$\gamma_{m}^{*}$ et $\gamma_{\mm}^{*}$.

\item Sur\label{theo-suite-spectrale-basique-relative-b} la page $\IEs(\U^{*})_1$ les $\FI$-modules en question, sont naturellement isomorphes aux $\FI$-modules définis par les morphismes
$$\IEs(q_m^{*})_{1}^{p,q}=
\smashbot{\ind^{\Sm }_{
\StS{m{-}(\pp)}{\pp}}
}(p_{m}^{*})$$
qui rendent commutatif le diagramme suivant.
$$\def\mbx#1#2{\hbox to#1{\hss$#2$}{}}
\xymatrix@C=6mm{
\mbx{2.4cm}{\IEs(\U^{m})_{1}^{p,q}}=
\ind^{\Sm }_{\vrule height4mm depth5pt width0pt
\hskip-4mm
\StS{m{-}(\pp)}{\pp}}
\hskip-10mm
\alt\otimes\Hbm^{q}(\UU^{m}_{m-p,\ldots,m})
\ar@{=>}[r]
\ar[d]^(.58){\ind (p_{m}^{*})}&
\Hbm^{i}(\Fm(\Mg))\ar[d]^{p_{m}^{*}}\\
\mbx{2.75cm}{\IEs(\U^{\mm})_{1}^{\pp,q}}=
\ind^{\Smm }_{\vrule height4mm width0pt\hskip-4mm
\StS{m{-}(\pp)}{p+2}}
\hskip-10mm
\alt\otimes\Hbm^{q}(\UU^{\mm}_{m-p,\ldots,m+1})
\ar@<2cm>@{<-}[u]@+<0pt,-8pt>^(.63){\IEs(q_m^{*})_{1}^{p,q}}
\ar@{=>}[r]&
\Hbm^{i}(\Fmm(\Mg))
}$$
où $\alt$ indique que l'action de $\S_{m{-}(\pp)}\times\1$ est tordue $\alt_{m{-}(\pp)}$.
\item Modulo\label{theo-suite-spectrale-basique-relative-c} les isomorphismes $\varXi$ de {\rm\ref{theo-suite-spectrale-basique}-(\ref{theo-suite-spectrale-basique-b})}, on a le diagramme commutatif d'espaces vectoriels
$$\def\ind#1#2{_{\hbox to#1{\hss$\scriptstyle#2$\hss}}}
\def\mbx#1#2{\hbox to1.2cm{$#2$\hss}}
\xymatrixc{@R6mm@C1.5cm}
{
{\Hbm^{q}(\UU^{m}_{m-p,\ldots,m})}
\ar[r]^(0.45){p_m^{*}}
\ar[d]_(0.48){\varXi_{\pp}^{m}\ }|(0.45){\xysimeq}
&
{\Hbm^{q}(\UU^{\mm}_{m-p,\ldots,m+1})}
\\
\displaystyle\bigoplus\ind{1cm}{f\in\FFF(\pp,m)}
\Hbm^{Q}(\Fg_{\pgoth(f)}(\Mg_{>0}))
\ar[r]^{\bigoplus p_{\pp}^{*}}&
\displaystyle\bigoplus\ind{1cm}{f^{\bullet}\in\FFF^{\bullet}(p{+}2,\mm)}
\Hbm^{Q}(\Fg_{\pgoth(f^{\bullet})}(\Mg_{>0}))
\ar@{^(->}[u]+<0pt,-9pt>_(0.6){(\varXi_{p+2}^{\mm})^{-1}}
}
\eqno\,\vtop{\kern0.27cm\hbox to0pt{\hss$(\diamond)$}}
$$
où nous avons noté
\def\varitemizeseps{\itemsep4pt\parskip0pt\topsep4pt}
\begin{itemize}\def\mlbox#1#2{\hbox to#1{\hss$#2$}}\halfdisplayskips
\item $Q:=q-(m{-}(\pp))\,\dMg=i-(m{-}(\pp))\,(\dMg{-}1)$. 
\item $\FFF(\pp,m):={}$ensemble des applications $f:\iii[1,m]$ 
telles que
$$\let\Big\big
\mlbox{5cm}{\Big(\hbox{$x<f(x)$, si $x<m-p$}\Big)}\text{\quad et\quad}\Big(\hbox{$f(x)=x$, autrement}\Big)\,.$$
\item $\FFF^{\bullet}(p{+}2,\mm):={}$ensemble des applications $f:\iii[1,\mm]$ telles que\glossary{$\FFF^{\bullet}(p{+}2,\mm)$:ensemble des applications $f:\iii[1,\mm]$ strictement croissantes sur $\iij[1,m-p]$ qui fixent $\iii[m{-}p,\mm]$} 
$$\let\Big\big
\mlbox{5cm}{\Big(\hbox{$x<f(x)<\mm$, si $x<m-p$}\Big)}\text{\quad et\quad}\Big(\hbox{$f(x)=x$, autrement}\Big)\,.$$
\item $\pgoth(f)$ est la partition de $\iii[1,m]$ définie par les fibres de $f$ et la fonction $f^{\bullet}:\iii[1,\mm]$ est le prolongement de $f$ qui vérifie $f^{\bullet}(\mm)=\mm$.
\end{itemize}

\end{enumerate}\end{theo}
\demo \displayskips65/100
\parskip0.5ex\mou
(\ref{theo-suite-spectrale-basique-relative-a},\ref{theo-suite-spectrale-basique-relative-b})
L'identification des morphismes $\IEs(q_m^{*})_{1}^{p,q}$ comme morphismes induits résulte de \ref{theo-suite-spectrale-basique} qui montre que l'on a
$$\IEs(\U^{m})_{1}^{p,q}=
\ind^{\Sm }_{\vrule height4mm depth5pt width0pt
\hskip-4mm
\StS{m{-}(\pp)}{\pp}}
\hskip-10mm
\alt_{m{-}(\pp)}\otimes\Hbm^{q}(\UU^{m}_{m-p,\ldots,m})\,,
\postskip0pt$$
et comme $\IEs(q_m^{*})_{1}^{p,q}$
est un morphisme de $\Sm$-modules (\ref{prop-naturalite-suite-spectrale-basique}), il est déterminé par sa restriction à $\alt_{m{-}(\pp)}\otimes\Hbm^{q}(\UU^{m}_{m-p,\ldots,m})$ qui n'est autre que $\id\otimes p_{m}^{*}$ où $\id$ désigne l'identité sur $\alt_{m{-}(\pp)}$. Le fait que l'on obtient ainsi un $\FI$-module est alors immédiat. A partir de là, on conclut grâce à la compatibilité des différentielles $d_{r}$ avec les actions des groupes symétriques et au fait que la catégorie des $\FI$-modules est abélienne (\ref{familles-de-representations}).

\medskip
(\ref{theo-suite-spectrale-basique-relative-c})
Compte tenu de \ref{theo-suite-spectrale-basique}-(\ref{theo-suite-spectrale-basique-b}), nous avons seulement à justifier la dernière ligne $(\diamond)$ du diagramme, somme directe des duaux des morphismes d'intégration sur les fibres
$$p_{m!}:
\Hc(\UU^{\mm}_{m-p,\ldots,\mm})\to
\Hc(\UU^{m}_{m-p,\ldots,m})[-d_{\Mgg}]
$$
où, en notant $\Fg_{p+2}:=\Fg_{p+2}(\Mgg)$ pour simplifier, 
$$\UU^{\mm}_{m-p,\ldots,\mm}=\Delta_{\mm}(\Mgge^{m{-}(\pp)}\times\Fg_{p+2})\,.$$

\smallskip
\noindent{\sl Avertissement. Dans un souci d'allégement de notations, nous allons omettre d'écrire `$\Hc(\_)$' autour des termes des diagrammes qui vont suivre.}

\smallskip
Pour tous $a,b\in\NN$, le morphisme d'intégration sur la dernière coordonnée
$$p_{!}:
\Delta_{a+b+1}(\Mgge^{a}\times\Fg_{b+1})\to
\Delta_{a+b}(\Mgge^{a}\times\Fg_{b})
\eqno(\dagger)$$
est défini dans \ref{def-action-image-inverse}-(\ref{def-action-image-inverse-a})  à travers l'inclusion ouverte
$$\xymatrix@R7mm{
\Delta_{a+b+1}(\Mgge^{a}\times\Fg_{b+1})\ar[r]^(0.475){\iota}\ar[rd]_{p_{!}}&
\Delta_{a+b}(\Mgge^{a}\times\Fg_{b})\times\Mgg\\
&\Delta_{a+b}(\Mgge^{a}\times\Fg_{b})[-d_{\Mgg}]
\ar@{<-}[u]+<25pt,-8pt>_(0.80){\int_{\Mgg}}
}$$
et lorsque $a>0$, on est conduit à considérer le diagramme suivant
$$\xymatrix@R=7mm@C=5mm{
\Delta_{a+b+1}(\Mgge^{a}\times\Fg_{b+1})\ar[r]_{\iota}\ar[d]_{\jmath}&
\Delta_{a+b}(\Mgge^{a}\times\Fg_{b})\times\Mgg\ar[d]_{\jmath\times\id}
\ar[r]_(0.76){\ \smashbot{\int_{\Mgg}}\ }&\\
\Mgge\times\Delta_{a+b}(\Mgge^{a-1}\times\Fg_{b+1})\ar[d]_{\rho}\ar[r]^(0.46){\iota'}&
\big[\Mgge\times\Delta_{a+b-1}(\Mgge^{a-1}\times\Fg_{b})\big]\times\Mgg\ar[d]_{\rho\times\id}\\
\smash{\displaystyle\coprod_{\FFF(a+b,a+b+1)}}\hskip-15pt
\Delta_{a+b}(\Mgge^{a-1}\times\Fg_{b+1})\ar@{..>}[r]^(0.48){\iota''}
\ar[d]_(0.6){c}^(0.6){[+1]}&
\hskip-0.6cm
\smash{\displaystyle\coprod_{\FFF(a+b-1,a+b)}}\hskip-15pt
\Delta_{a+b-1}(\Mgge^{a-1}\times\Fg_{b})\times\Mgg
\ar[d]_(0.6){c\times\id}^(0.6){[+1]}\ar[r]^(0.85){\int_{\Mgg}}&
\\&&
}$$
où les colonnes sont les suites exactes longues de cohomologie à support compact habituelles dans nos théorèmes de scindage (\ref{theo-scindage}).
Les termes de la deuxième ligne sont nuls puisque $\Mgge$ apparaît en facteur, les flèches de liaison $c$ sont donc des isomorphismes et la
question que nous cherchons à comprendre concerne l'interprétation de $p_{!}$ dans $(\dagger)$ en termes de la dernière ligne du diagramme. Or, 
la flèche en pointillé n'est pas tout à fait bien définie. En effet, les composantes du fermé
$\hskip-15pt\coprod\limits_{\FFF(a+b,a+b+1)}\hskip-15pt
\Delta_{a+b}(\Mgge^{a-1}\times\Fg_{b+1})$
sont bien les traces des composantes de 
$\hskip-15pt
\coprod\limits_{\FFF(a+b-1,a+b)}\hskip-15pt
\Delta_{a+b-1}(\Mgge^{a-1}\times\Fg_{b})\times\Mgg$
sur l'ouvert $\Mgge\times\Delta_{a+b}(\Mgge^{a-1}\times\Fg_{b+1})$, que l'on va noter $\mathbb U$ dans la suite, à l'exception près de celle indexée par la fonction $f_0:\iii[a+b+1]$ définie par $f_0(1)=a+b+1$, composante, par ailleurs, clairement fermée dans l'ouvert $\Delta_{a+b}(\Mgge^{a}\times\Fg_{b})\times\Mgg$ que l'on va noter $\VV$ dans la suite.

On a donc l'inclusion fermée
$$(\Mgge\times\Delta_{a+b}(\Mgge^{a-1}\times\Fg_{b+1}))_{f_0}\dans \mathbb U\cap\VV$$
et tout cocycle $\omega$ d'Alexander-Spanier (faisceau $c$-mou) à support compact de cette composante se prolonge en une cochaîne $\varpi$ à support compact de $\mathbb U\cap\VV$. Le morphisme de liaison donne alors
$c(\omega)=d\varpi$, ce qui est un cocycle à support compact dans $\Delta_{a+b+1}(\Mgge^{a}\times\Fg_{b+1})$. Or, 
$$p_{!}(c(\omega))=\int_{\Mgg}\iota(c(\omega))=\int_{\Mgg}d\varpi=0\,,$$ 
puisque $\varpi\in\Z_{\rmc}(\VV)$ et que $\VV$ est le domaine de définition de $\int_{\Mgg}$.

Ces remarques prouvent la commutativité du diagramme
$$\xymatrix@R10mm{
\smash{\displaystyle\coprod_{\FFF(a+b,a+b+1)}}\hskip-15pt
\Delta_{a+b}(\Mgge^{a-1}\times\Fg_{b+1})[-1]
\ar[r]^(0.48){\iota''}
\ar[d]^(0.5){\ c}|
{\xysimeq}&
\hskip-0.6cm
\smash{\displaystyle\coprod_{\FFF(a+b-1,a+b)}}\hskip-15pt
\Delta_{a+b-1}(\Mgge^{a-1}\times\Fg_{b})[-d_{\Mgg}][-1]
\ar[d]^(0.5){\ c}|
{\xysimeq}
\\
\Delta_{a+b+1}(\Mgge^{a}\times\Fg_{b+1})\ar[r]^(0.48){p_{!}}&
\Delta_{a+b}(\Mgge^{a}\times\Fg_{b})[-d_{\Mgg}]\\
}$$
où maintenant $\iota''$ est bien définie. Sa valeur, 
qui dépend des composantes, vaut
$$\preskip0.5ex\postskip1.75ex
\iota''=\begin{cases}\noalign{\kern-3pt}
0\text{ sur }(\Mgge\times\Delta_{a+b}(\Mgge^{a-1}\times\Fg_{b+1}))_{f_0}\,,\\\noalign{\kern4pt}
p_{!}:\Delta_{a+b}(\Mgge^{a-1}\times\Fg_{b+1})\to \Delta_{a+b-1}(\Mgge^{a-1}\times\Fg_{b})\,,\text{\ autrement.}\\\noalign{\kern-2pt}
\end{cases}$$

L'itération de ces idées conduit au diagramme commutatif
$$\xymatrix@R10mm{
\smash{\displaystyle
\coprod_{\FFF(b+1,a+b+1)}}
\Fg_{b+1}[-a]
\ar[r]^(0.5){\iota''}
\ar@<1cm>[d]|
(0.5){\xysimeq}&
\smash{\displaystyle\coprod_{\FFF(b,a+b)}}
\Fg_{b}[-d_{\Mgg}][-a]
\ar@<0.5cm>[d]|
(0.5){\xysimeq}
\\
\Delta_{a+b+1}(\Mgge^{a}\times\Fg_{b+1})\ar[r]^(0.45){p_!}&
\Delta_{a+b}(\Mgge^{a}\times\Fg_{b})[-d_{\Mgg}]\\
}$$
où $\iota''=0$ sur les composantes indexées par les fonctions $f:\iii[1,a+b+1]$ telles que $|f^{-1}(a+b+1)|>1$, et c'est l'intégration sur les fibres autrement.

\nobreak En dualisant, on obtient la description de la ligne $(\diamond)$ du diagramme du théorème pour la cohomologie de Borel-Moore.
\enddemo

\comment
Mais avant cela, il nous faut rappeler et adapter des notions classiques de la théorie de bicomplexes de \vCech-deRham, pour le contexte présent des espaces localement compacts et de la cohomologie de Borel-Moore.

\subsubsectionline{Sur la suite spectrale d'un cône de bicomplexes.}
Les constructions des bicomplexes $\vC(\U^{m})^{\bullet,*}_{?}$ des paragraphes précédents se localisent naturellement. 
En effet, pour tout ouvert $W\dans\Mg$, notons 
$$\U^{m}_{W}:=\set\UU^{m}_{1}\cap W,\ldots,\UU^{m}_{m}\cap W/\,,$$ le recouvrement induit sur l'ouvert $\UU^{m}_{W}:=\UU^{m}\cap W$.
Une inclusion ouverte $\iota_{1,2}:W_1\dans W_2$ induit alors par prolongement par zéro un morphisme (injectif) de bicomplexes
$$\iota_{1,2!}:\vC(\U^{m}_{W_1})^{\bullet,*}_{\rm c}\hook
\vC(\U^{m}_{W_2})^{\bullet,*}_{\rm c}\,,$$
et la correspondance $W\fonct\vC(\U^{m}_{W})^{\bullet,*}_{\rm c}$ est fonctorielle sur la catégorie des ouverts de $\Mg$.

Si $W_1,W_2$ sont deux ouverts, la construction classique des suites de Mayer-Vietoris donne lieu à une suite exacte courte de bicomplexes
$$0\to
\vC(\U^{m}_{W_{12}})^{\bullet,*}_{\rmc}\to
\vC(\U^{m}_{W_1})^{\bullet,*}_{\rmc}\oplus
\vC(\U^{m}_{W_2})^{\bullet,*}_{\rmc}\to
\vC(\U^{m}_{W_1\cup W_2})^{\bullet,*}_{\rmc}\to0\,,
$$
et par dualité, à \expression{la suite exacte courte de Mayer-Vietoris pour les bicomplexes de cochaînes de \vCech-Borel-Moore}
$$\let\rmc\bm
0\to
\vC(\U^{m}_{W_1\cup W_2})^{\bullet,*}_{\rmc}\too^{\rho}
\vC(\U^{m}_{W_1})^{\bullet,*}_{\rmc}\oplus
\vC(\U^{m}_{W_2})^{\bullet,*}_{\rmc}\to
\vC(\U^{m}_{W_{12}})^{\bullet,*}_{\rmc}\to0\,,
\eqno(\diamond)$$
où $\rho$ est la somme directe des morphismes de restriction de $W_1\cup W_1$ vers chaque $W_i$.

\soustitreline{Cône d'un morphisme de complexes.}Commençons par quelques rappels généraux sur les cônes de morphismes de complexes.
Si $f_*:\AAA^{*}\to\BBB^{*}$ est un morphisme de complexes de cochaînes, le \expression{cône de $f$} est le complexe
$$\cone{f}:=(\BBB\oplus\AAA[1],d_f)\quad\hbox{avec}\quad
d_{f}(b,a):=(db-fa,-da)\,.\eqno(\ddagger)$$
Il participe dans la suite longue
$$\def\to{\mathop{\longrightarrow}}
\to\AAA\to^{f}\BBB\to^{\beta}\cone f\,\smash{\to^{\alpha}_{\scriptscriptstyle[+1]}}\,\AAA[1]\to^{f[1]}\BBB[1]\to^{\beta[1]}\cone f[1]\to^{\alpha[1]}
\eqno(\ddagger\ddagger)$$
où l'on a noté
$$\begin{cases}\noalign{\kern-4pt}
\beta:\BBB\to\cone f\,,&\beta(b):=(b,0)\\
\alpha:\cone f\to\AAA[1]\,,&\alpha(b,a):=(a)\\
\noalign{\kern-3pt}\end{cases}
\postskip1em$$
et aussi dans le diagramme commutatif de morphismes de complexes
$$\preskip0.6ex\postskip0.6ex
\xymatrixc{@R=4mm}{
0\ar[r]&
\AAA\ar[r]^{f}&
\BBB\ar[r]^{g}&
\CCC\ar[r]&0\\
&
\AAA\aregal[u]\ar[r]^{f}&
\BBB\aregal[u]\ar[r]&
\cone{f}\ar[u]_{\ \tilde g}\ar[r]&
}\eqno(\ddagger{\ddagger}\ddagger)$$
où $g:\BBB\onto\CCC$ est le conoyau de $f$ et où $\tilde g(b,a)=g(b)$.

\'Etant donné un diagramme commutatif de morphismes de complexes
$$\preskip1ex\xymatrixc{@R=5mm}{
\AAA'\ar[r]^{f'}&\BBB'\\
\AAA\ar[u]^{\delta_{\AAA}}\ar[r]^{f'}&\BBB\ar[u]_{\delta_{\BBB}}\\
}$$
l'application $(\delta_\BBB\oplus\delta_\AAA[1]):\cone{f_1}\to\cone{f_2}$
est un morphisme de complexes et rend le diagramme suivant commutatif
$$\def\ligne#1#2#3#4#5{\AAA_{#1}#2\ar[r]^{f_{#1}}&\BBB_{#1}#3\ar[r]^{\beta}&\cone{f_{#1}}#4\ar[r]^{\alpha_{#1}}&\AAA_{#1}#5[1]
}
\xymatrixc{@R=7mm@C=1cm}{
\ligne2{}{}{}{}
\\
\ligne1
{\ar[u]_{\delta_\AAA}}
{\ar[u]_{\delta_\BBB}}
{\ar[u]_{(\delta_\BBB,\delta_\AAA[1])}}
{\ar[u]_{\delta_\AAA[1]}}
\\
}$$

\noindent Le lemme suivant est élémentaire bien connu.\killline
\let\bold\relax\noendpoint\begin{lemm}\label{lemme-cone}
\begin{enumerate}
\item Toutes les\label{lemme-cone-a} constructions $(\ddagger's)$ dépendent fonctoriellement de $f:\AAA\to\BBB$.
\item Par\label{lemme-cone-b} passage en cohomologie, le diagramme $(\ddagger\ddagger)$ devient un isomorphisme de suites exactes longues
$$
\xymatrixc{@C=6mm@R=4mm}{
\cdots\ar[r]&
h^{i}(\AAA)\ar[r]^{f}&
h^{i}(\BBB)\ar[r]^{g}&
h^{i}(\CCC)\ar[r]&
h^{i+1}(\AAA)\ar[r]^{f[1]}&
h^{i+1}(\BBB)\ar[r]&\cdots\\
\cdots\ar[r]&
h^{i}(\AAA)\aregal[u]\ar[r]^{f}&
h^{i}(\BBB)\aregal[u]\ar[r]^(0.43){\beta}&
h^{i}(\cone f)\ar[u]_{\tilde g}^{\simeq}\ar[r]^{\alpha}&
h^{i+1}(\AAA)\aregal[u]\ar[r]^{f[1]}&
h^{i+1}(\BBB)\aregal[u]\ar[r]&\cdots\\
}
$$
\end{enumerate}
\end{lemm}

\soustitre{Cône (horizontal) d'un morphisme de bicomplexes}
Un bicomplexe se voit comme un complexe de complexes ``horizontaux'' ou de complexes ``verticaux''. Dans le premier cas, un morphisme de bicomplexes $f_{\bullet,*}:\AAA^{\bullet,*}\to\BBB^{\bullet,*}$ est une suite de morphismes de complexes (horizontaux) $\set f_{\ell,*}:\AAA^{\ell,*}\to\BBB^{\ell,*}/_{\ell}$. Si on applique les constructions fonctorielles $(\ddagger's)$, on obtient le diagramme commutatif de morphismes de complexes:
$$\def\ligne#1#2#3#4#5{\AAA^{#1}#2\ar[r]^{f_{#1}}&\BBB^{#1}#3
\ar[r]^{\beta_{#1}}&\cone{f_{#1}}#4\ar[r]^{\alpha_{#1}[1]}&\AAA^{#1}[1]#5\\}
\xymatrixc{@R=5mm@C=1cm}{
\ligne{\ell+1}{}{}{}{}
\ligne{\ell}{\ar[u]^{\delta_{l-1}}}{\ar[u]_{\delta_{l-1}[1]}}{\ar[u]_{\delta_{l-1}}}{\ar[u]_{\delta_{l-1}}}
\ligne{\ellmo}{\ar[u]^{\delta_{l-1}}}{\ar[u]_{\delta_{l-1}[1]}}{\ar[u]_{\delta_{l-1}}}{\ar[u]_{\delta_{l-1}}}
}
$$
où chaque est un bicomplexe. La troisième colonne est le bicomplexe \expression{cône (horizontal) de $f_{\bullet,*}$} noté $\coneh f$.
On obtient ainsi l'analogue de la suite longue $(\ddagger\ddagger)$
pour les bicomplexes
$$\def\to{\mathop{\longrightarrow}}\let\cone\coneh
\to\AAA\to^{f}\BBB\to^{\beta}\cone f\,\smash{\to^{\alpha}_{\scriptscriptstyle[0,+1]}}\,\AAA[0,1]\to^{f[0,1]}\BBB[0,1]\to^{\beta[0,1]}\cone f[0,1]\to^{\alpha[1]}
\eqno(\ddagger\ddagger)$$
de même que l'analogue de $(\ddagger{\ddagger}\ddagger)$.

\begin{prop}Pour tout morphisme de bicomplexes $f:\AAA\to\BBB$ on a le diagramme commutatif de morphismes de bicomplexes
$$\preskip0.6ex\postskip0.6ex
\xymatrixc{@R=4mm}{
0\ar[r]&
\AAA\ar[r]^{f}&
\BBB\ar[r]^{g}&
\CCC\ar[r]&0\\
&
\AAA\aregal[u]\ar[r]^{f}&
\BBB\aregal[u]\ar[r]^(.45){\beta}&
\coneh{f}\ar[u]_{\ \tilde g}\ar[r]^{\alpha}&\AAA[0,1]
}$$
où $g:\BBB\onto\CCC$ est le conoyau de $f$ et où $\tilde g(b,a)=g(b)$. Lorsque $\AAA$ et $\BBB$ sont bornés à gauche, le morphisme $\tilde g$
induit un quasi-isomorphisme entre complexes simplexes simples associés
$$\tilde g:\tot^{*}(\coneh{f})\to\tot^{*}(\CCC)$$
\end{prop}
 de complexes déterminés par 
 morphisme injectif de bicomplexes, les constructions $(\ddagger)$ et $(\ddagger)$ s'appliquent ligne à ligne (\ref{lemme-cone}-(\ref{lemme-cone-a})).

On note alors $\coneh f$ le bicomplexe obtenu de

Maintenant, si $0\to\AAA^{\bullet,*}\hf{f_{\bullet,*}}{}{0.8cm}
\BBB^{\bullet,*}\hf{g_{\bullet,*}}{}{0.8cm}
\CCC^{\bullet,*}\to0
$
est une suite exacte de bicomplexes du premier quadrant, on a pour chaque ligne $\ell\in\NN$ un morphisme de suites de complexes
$$\xymatrix{
0\ar[r]&
\AAA^{\ell,*}\ar[r]^{f_{\ell,*}}&
\BBB^{\ell,*}\ar[r]^{g_{\ell,*}}&
\CCC^{\ell,*}\ar[r]&0\\
&
\AAA^{\ell,*}\aregal[u]\ar[r]^{f_{\ell,*}}&
\BBB^{\ell,*}\aregal[u]\ar[r]&
\cone{f_{\ell,*}}\ar[u]_{\tilde g_{\ell,*}}\ar[r]&
}$$

est un morphisme de bicomplexes du premier quadrant de conoyau le bicomplexe $\CCC^{\bullet,*}$
, on note 
note $\cone{\alpha_{i,*}}$ le cône de $\alpha_{i,*}$. On a donc des triangles exacts

a le triangle exact
$$\AAA^{i,*}\hf{\alpha_{i,*}}{}{1cm}\BBB^{i,*}
\hf{}{}{1cm}\cone{\alpha^{i,*}}$$

 on note
$\coneh{\alpha_{\bullet,*}}$ le bicomplexe dont les lignes sont les cônes

Dans\label{eilenberg-moore} l'article célèbre d'Eilenberg-Moore (\cite{e-m}), on considère 
pour tout morphisme de complexes différentiels \emph{filtrés}
$f:\A\to\B$, 
le triangle 
$$\xymatrix{
\A\ar[r]^{f}&
\B\ar[r]^(0.42){\beta}&
\cone f\ar[r]^(0.5){\alpha}&\A[1]
}$$
où $\cone f$\glossary{${\cone f:=(\B\oplus\A[1],d_f)}$:cône d'un morphisme $f:\A\to\B$
 de modules différentiels gradués filtrés, avec $d_{f}(b,a):=(db-f(a),-d(a))$} désigne le \expression{cône de $f$}

est le cône de $f$, avec $(b,a)=db-f(a)$

(\cite{e-m} (7.12), p. 23) (\footnote{Voir aussi l'exercice 5.4.4, p. 135, du livre de Weibel \cite{weibel} dans sa version de 1997.}) affirme que pour tout 
\endcomment

\section{Stabilité des familles
de représentations}\label{stabilite}
\glossarytitle{Stabilité des familles
de représentations}

\subsection{Motivations}Le\label{stabilite-caracteres} comportement lorsque $m\mapsto\infty$ des invariants associés aux espaces $\Fm(\Xg)$, notamment leurs groupes d'homotopie, leurs nombres de Betti, les caractères de leurs cohomologies en tant que $\S _{m}$-modules, ont fait l'objet de très nombreux travaux.  Un des résultats pionniers dans ces questions est dû à V.I.~Arnold (\emph{circa} 1970) et concerne la détermination de la cohomologie des groupes de tresses d'E.~Artin (1925) à un moment où l'on savait, grâce aux travaux de Fadel-Fox-Neuwirth (\cite{FN,FN2}, 1962), que les espaces $\Fm(\CC)$ et $\BFm(\CC):=\Fm(\CC)/\S _{m}$ étaient asphériques, de groupes fondamentaux $\Pi_{1}\Fm(\CC)=P_m$: \expression{le groupe des tresses colorées à $m$ brins}, et $\Pi_{1}\BFm(\CC)=\Bg_m$: \expression{le groupe des tresses non colorées à $m$ brins}, et donc que les cohomologies des espaces $\Fm(\CC )$ et $\BFm(\CC )$ étaient canoniquement isomorphes à celles des groupes de tresses. C'est ainsi que Arnold procède dans \cite{arnold,arnold2} où il énonce son \expression{théorème de stabilité} , théorème qui établit que, pour chaque $i\in\NN$ fixé, le morphisme naturel 
$\Hr^{i}(\BFm(\CC);\ZZ)\to\Hr^{i}(\BF_{\mm}(\CC);\ZZ)$ est bijectif pour $m\geq 2i-2$. C'est un résultat remarquable que ne révèle pas la connaissance explicite du polynôme de Poincaré de $\Fm(\Xg)$ (voir thm.~\ref{theo-conf-symetriques}). Parallèlement, le même phénomène de stabilité pour l'espace $\Fm(\Xg)$ ne se produit pour ainsi dire jamais, déjà si nous appliquons la formule de la remarque \ref{rema-Fg-orientable} pour $\Xg=\CC$, nous voyons que
$\Pr(\Fm(\CC))(T)=(1+T)(1+2T)\cdots(1+(\mmo)T)$
d'où $\Betti^{1}(\Fm(\CC))= {m\choose2}$ (polynomiale en $m$ tout de même !).

Il faut attendre les années 2011 pour comprendre le lien entre ces deux  phénomènes grâce aux travaux de T.~Church et B.~Farb (\cite{chu,chu-far}).
En s'intéressant, non pas à la famille $\set\Betti^{i}(\BFm(\CC))/_{m}$, mais à la famille des représentations $\set\Sm\rep\Hr^{i}(\Fm(\CC))/_{m}$, Church et Farb se livrent à des calculs explicites en s'appuyant sur les recherches de Lehrer (\cite{Le}, 1987) qui donnaient déjà la \expression{polynomialité} des caractères de cette famille \comment
(\footnote{Autrement dit, l'existence un polynôme
$P(\CC;i)(X_1,\ldots,X_{r})$, où $X_i$ est le nombre des cycles de longueur $i$ qui interviennent dans la décomposition d'une permutation en produit de cycles disjoints, tel que l'on a  
$$\chi(\Fm(\CC);i)(\alpha)=P(\CC;i)(X_1(\alpha),\ldots,X_{r}(\alpha))\,,\quad\forall\alpha\in\Sm\,,\ \forall m\gg0\,.$$})
\endcomment (\ref{caracteres-polynomiaux}), et sont conduits à conjecturer que les multiplicités des composantes irréductibles des représentations $\set\S _{m}\rep\Hr^{i}(\Fm(\Xg))/_{m}$ sont stationnaires (dans un sens qui sera précisé dans \ref{repr-stables}), de sorte que lorsque $\Xg$ est en plus de type fini, la polynomialité des caractères en découle moyennant un résultat classique de Macdonald (\cf\ref{macdonald-poly}, 1995). 
Church et Farb développent alors la théorie de \expression{$\FI$-modules} (\ref{familles-de-representations}) et remarquent que le phénomène de multiplicités stationnaires pouvait s'expliquer comme conséquence de ce qu'ils ont appelé la \expression{stabilité} (\ref{rang-stabilite}) de la famille de morphismes $\set p_m^{*}:\Hr^{i}(\Fm(\Xg))\to\Hr^{i}(\Fmm(\Xg))/_{m}$, ce qui s'avéra être effectivement le cas.\killline

\begin{theo*}[(Church \cite{chu}, 2012)]Soit $\Xg$ une variété topologique, connexe et orientable. Pour $i\in\NN$ fixé, la famille $\set p_m^{*}:\Hr^{i}(\Fm(\Xg))\to\Hr^{i}(\Fmm(\Xg))/_{m}$ est stable pour $m\geq 2i$ si $\dim\Xg\geq 3$, et pour $m\geq 4i$ si $\dim \Xg=2$.
Les familles des caractères et des nombres de Betti correspondantes sont  polynomiales
et la famille $\set \Betti^{i}({\Fm(\Xg)/\Sm })/_{m}$ est constante, sur les mêmes rangs.\killline
\end{theo*}

Le théorème de stabilité d'Arnold apparaît de lors comme le
fait que la famille des sous-représentations triviales $\set\Hr^{i}(\Fm(\CC))^{\S _{m}}/_{m}$ est stationnaire.

\medskip
Les paragraphes qui suivent sont destinés à rappeler les bases de la théorie de $\FI$-modules (\cite{cef}) qui a servi à prouver le théorème de Church. Nous introduirons ensuite certains foncteurs d'induction dans la catégorie des $\FI$-modules qui vont nous permettre de généraliser le théorème de Church aux familles $\set\Sm \rep\Delta_{?m-a}\Xg^{m}/_{m}$ où $\Xg$ est une pseudovariété orientable. On procédera dans un premier temps en supposant que $\Xg$ est $i$-acyclique (thm. \ref{theo-stabilite-BM-i-acyclique}), puis, dans un deuxième temps, sans cette hypothèse (thm. \ref{theo-stabilite-BM-pseudo}) mais à l'aide d'un nouvel instrument: les suites spectrales basiques (\ref{theo-suite-spectrale-basique}). 
\comment On remarquera en passant que $\Delta_{\leq m-a}\Xg^{m}$ est généralement singulier, même lorsque $\Xg$ ne l'est pas.
\endcomment

\subsection{Catégorie des $\FI$-modules}\label{familles-de-representations}Suivant \cite{cef}, on note $\FI$\glossary{$\FI$:catégorie des ensembles {\bfit F\/}inis et des applications {\bfit I\/}njectives} la catégorie dont les objets sont les ensembles finis, et dont les morphismes sont les applications injectives. 
Si $\AAA$ est un anneau, on note
 $\ModFI \AAA$\glossary
{$\ModFI \AAA$:catégorie des $\FI$-modules} 
la catégorie dont les objets, les \expression{$\AAA[\FI]$-modules}\index{FI-module@$\FI$-module} (et même simplement $\FI$-modules lorsque $\AAA$ est sous-entendu), sont les foncteurs covariants $\V:\FI\fonct\Mod(\AAA)$ et dont les morphismes sont les transformations naturelles entre ces foncteurs. Si $\V$ et $\W$ sont des $\FI$-modules, on note $\Hom_{\FI}(\V,\W)$ l'ensemble des morphismes de $\V$ vers $\W$.

\smallskip\noindent{\bf Conventions et notations. }Dans ce qui suit, on identifie pour tous $b\geq0$
$$\S _{m}=\Fix\nolimits_{\S _{m+b}}\iii[m+1,m+b]\text{\quad et\quad}
\II_{m}\times\S _{b}=\Fix\nolimits_{\S _{m+b}}\iii[1,m]\,,$$
ce qui donne le sens à l'inclusion $\S _{m}\times\S _{b}\dans\S _{m+b}$.

Dans une notation $V_{m}\boxtimes \AAA_{b}$, on désigne par $V_{m}$ un $\S _{m}$-module et par $\AAA_{b}$ la représentation triviale de $\S _{b}$ dans $\A$. Le produit tensoriel en question est alors le $\S _{m}\times\S _{b}$-module défini par $(\alpha,\beta)\cdot(v\otimes w):=(\alpha\cdot v,\beta\cdot w)$.

\subsubsectionline{Une équivalence de catégories.}La\label{familles-denombrables} restriction d'un $\FI$-module à la sous-catégorie pleine des intervalles de la forme $\iii[1,m]\dans\NN$ est une équivalence de catégories entre
$\ModFI \AAA$ et la catégorie des familles dénombrables $\V=\set \phi_{m}:V_{m}\to V_{\mm}/_{m}$\glossary{${\V=\set \phi_{m}:V_{m}\to V_{\mm}/_{m}}$:famille dénombrable représentant un $\FI$-module}, où les $V_{m}$ sont des $\AAA[\S _{m}]$-modules et où les $\phi_{m}$ sont des morphismes de $\AAA[\S _{m}]$-modules dont les composées
$$\phi_{m+b,m}:=\phi_{m+b}\circ\cdots\circ\phi_{\mm}\circ\phi_{m}\,,\quad\forall b\in\NN\,,\postskip0pt
\eqno(\ddagger)$$
vérifient
$$\phi_{m+b,m}(V_{m})\dans (V_{m+b})^{\II_{m}\times\S _{b}}\,,\quad \forall m,b\in\NN\,.\eqno(\diamond)$$

Les morphismes $\phi_{m+b,m}$ sont les \expression{les morphismes de transition de $\V$}\index{morphisme de transition}\glossary{${\phi_{m+b,m}:V_{m}\to V_{m+b}}$:morphisme de transition de  $\V=\set \phi_{m}:V_{m}\to V_{\mm}/$}, ils se factorisent suivant le diagramme
$$\preskip4pt
\xymatrixc{@C=1.cm@R=4mm}{
V_{m}\ar[rd]+<-50pt,2pt>_(0.6){\iota}\ar[rr]|{\ \phi_{m+b,m}\ }&&V_{m+b}\ar@{<-}[ld]+<50pt,2pt>^(0.55){\ \varPhi_{m+b,m}}
\\
&\ind^{\S _{m+b}}_{\S _{m}\times\S _{b}}(V_{m}\boxtimes \AAA_{b})}
\eqno(\diamond\diamond)$$
où $\iota:V_{m}\to V_{m}\boxtimes \AAA $ est le plongement  $v\mapsto v\otimes \1_{\AAA}$ et où $\varPhi_{m+b,m}$ est le morphisme de $\AAA[\S _{m+b}]$-modules induit~(\footnote{La condition $(\diamond)$ qui porte sur les $\II_{m}\times\S _{b}$-invariants est vide pour $b=1$ et n'est donc pas transitive, il suffit par contre qu'elle soit vérifiée pour $b=2$ et tout $m\in\NN$.}).

Dans cette équivalence, un morphisme $f:\V\to \W$ de $\FI$-modules se voit comme une famille $\set f_{m}:V_{m}\to W_{m}/_{m}$ de morphismes de $\AAA[\S _{m}]$-modules vérifiant $\phi_{m}\circ f_{m}=f_{\mm}\circ\phi_{m}$. Les familles $\set\ker(f_{m})/_{m}$ et $\set\coker(f_{m})/_{m}$
$$\xymatrix@C=1.2cm@R=6mm{
\ker(f_m)\ar@{-->}[d]\arinto[r]
&V_m\ar[d]_{\phi_m}\ar[r]^{f_{m}}\xylbl[rd]{$\bigoplus$}&W_m\ar[d]^{\phi_m}\aronto[r]&\coker(f_m)\ar@{-->}[d]\\
\ker(f_{\mm})\arinto[r]&V_{\mm}\ar[r]_{f_{\mm}}&W_{\mm}\aronto[r]&\coker(f_{\mm})
}$$
 munies des morphismes de transition induits sont des $\FI$-modules et font de $\ModFI \AAA$ une catégorie abélienne (\footnote{Dans \cite{ss}, Sam et Snowden donnent une équivalence de catégories explicite entre $\ModFI \AAA$ et une catégorie de modules sur un  anneau.}).
 
\subsubsectionline{Troncatures de $\FI$-modules.}Soit\label{troncature} $q\in\NN$. Pour tout $\FI$-module $\V=\set \phi_m:V_m\to V_m+1/$, on notera 
$\V_{\geq q}=\set W_m\to W_{\mm}/$\glossary{$\V_{\geq q}$:tronqué de $\V$ qui préserve les $V_m$ pour $m\geq q$ et annule les autres.} le sous-$\FI$-module de $\V$ 
avec $W_m:=0$ si $m\leq{q-1}$, et $W_m:=V_m$ autrement.

\noindent On pose ensuite $\V_{\leq q-1}:=\V/\V_{\geq q}$\glossary{$\V_{\leq q}$:tronqué de $\V$ qui préserve les $V_m$ pour $m\leq q$ et annule les autres.}, d'où la suite exacte courte de $\FI$-modules
$$\xymatrix@C=0.5cm@R=7mm{
\vrule depth6pt width0pt\V_{\geq q}\arhook[d]&\ar[r]&0\ar[d]\ar[r]&0\ar[r]\ar[d]&V_{q}\ar[r]\aregal[d]&V_{q+1}\ar[r]\aregal[d]&V_{q+2}\ar[r]\aregal[d]&\\
\V\aronto[d]&\ar[r]&V_{q-2}\aregal[d]\ar[r]&V_{q-1}\aregal[d]\ar[r]&V_{q}\ar[d]\ar[r]&V_{q+1}\ar[d]\ar[r]&V_{q+2}\ar[d]\ar[r]&\\
\V_{\leq q-1}&\ar[r]&V_{q-2}\ar[r]&V_{q-1}\ar[r]&0\ar[r]&0\ar[r]&0\ar[r]&
}$$
On notera $\ModFI\A_{\geq q}$ (resp. $\ModFI\A_{\leq q}$) la sous catégorie pleine de $\ModFI\A$ dont les objets sont les $\FI$-modules tronqués $\V_{\geq q}$ (resp. $\V_{\leq q}$).

Les correspondances $\V\fonct \V_{\geq q}$ et $\V\fonct \V_{\leq q}$ sont clairement fonctorielles covariantes et exactes. 

\subsubsectionline{Les $\FI$-modules $\M(\ag)$.}Pour\label{M(a)} $0<a\in\NN$, on note \glossary{${\ag:=\iii[1,a]}$: notation pour l'ensemble fini représenté par l'intervalle $\iii[1,a]$}
$\ag:=\iii[1,a]\,,$ et ensuite $\M(\ag):={}$le foncteur covariant représenté par $\ag$\glossary{$\M(\ag)$: le foncteur représentable $(\_)\fonct\AAA[\Mor_{\FI}(\ag,\_)]$}, soit
$$\xymatrix@R=1mm@C=0.cm{\mllap{\M(\ag):}\FI\ar@{~>}[r]&**[r]\Mod(\AAA)\\
S\ar@{|~>}[r]&**[r]\M(\ag)_{S}\mrlap{:=\AAA[\Mor_{\FI}(\ag,S)]}
}
$$
La famille $\set \phi_m:\M(\ag)_{m}\to\M(\ag)_{\mm}/_{m}$ est alors caractérisée par
$$\begin{casesalign}
\hbox{si $b<0$,\quad}&\M(\ag)_{a+b}=0\hfill\\
\hbox{si $b\geq0$,\quad}\hfill&\M(\ag)_{a+b}=\ind^{\S _{a+b}}_{\S _{a}\times\S _{b}}\AAA[\S _{a}]\boxtimes \AAA_{b}\,,
\end{casesalign}$$
et le morphisme de transition $\phi_{a+b}:\M(\ag)_{a+b}\to \M(\ag)_{a+b+1}$ est celui induit à partir de l'identification de $\S _{a}\times\S _{b}$-modules
$k[\S _{a}]\boxtimes\A_{b}=k[\S _{a}]\boxtimes\A_{b+1}$,
$$\xymatrix@C=1.cm@R=7mm{
k[\S _{a}]\boxtimes\A_{b}\aregal[r]\arinto@<3.5pt>[d]&
k[\S _{a}]\boxtimes\A_{b+1}\arinto@<-1.5pt>[d]\\
\ind^{\S _{a+b}}_{\S _{a}\times\S _{b}}\AAA[\S _{a}]\boxtimes \AAA_{b}\ar@{..>}[r]^(0.46){\phi_{a+b}}&
\ind^{\S _{a+b+1}}_{\S _{a}\times\S _{b+1}}\AAA[\S _{a}]\boxtimes \AAA_{b+1}}$$
où les flèches verticales désignent les inclusions canoniques.

\subsubsectionnumber Par le lemme de Yoneda, pour tout $\FI$-module $\V$
 l'application
$$Y:\Hom_{\FI}(\M(\ag),\V)\to V_{a}\,,\quad Y(f)= f(\id_{\ag})\,,$$
est bijective. Si $v\in V_{a}$, on note $\ev_{\ag}(v):=Y^{-1}(v)$, il s'agit du morphisme de foncteurs dont les valeurs 
$\ev_{\ag}(v)_{a+b}:\M(\ag)_{a+b}\to V_{a+b}$ sont
$$\begin{casesalign}
\hbox{si $b<0$,\quad}&\ev_{\ag}(v)_{a+b}=0\hfill\\
\hbox{si $b=0$,\quad}\hfill&\ev_{\ag}(v)_{a}:\AAA[\S _{a}]\to V_{a},\ \ev_{\ag}(v)_{a}(\alpha)=\alpha\cdot v
\hfill\\
\hbox{si $b>0$,\quad}\hfill&\ev_{\ag}(v)_{a+b}\hbox{ est la composée:}
\hfill\\
\end{casesalign}$$
$$\preskip0.ex\ind_{\S _{a}\times\S _{b}}^{\S _{a+b}}\AAA[\S _{a}]\boxtimes \AAA_{b}\hf{\ind(\ev_{\ag}(v)_{a})\ }{}{2cm}
\ind_{\S _{a}\times\S _{b}}^{\S _{a+b}} V_{a}\boxtimes \AAA_{b}
\hf{\varPhi_{a+b,a}\ }{}{1.2cm}
V_{a+b}$$

\subsubsectionline{Sous-$\FI$-modules et $\FI$-modules quotients.}On\label{sous-FI-modules} appelle \expression{sous-$\FI$-module}\index{sous-FI-module@sous-$\FI$-module}  d'un $\FI$-module $\V$ tout
sous-foncteur $\W\dans\V$. En termes des familles dénombrables $\V=\set \phi_m:V_{m}\to V_{\mm}/_{m}$ (\ref{familles-denombrables}), l'inclusion $\W\dans\V$ équivaut à la donnée d'une famille de sous-$\S _{m}$-modules $W_{m}\dans V_{m}$ vérifiant $\phi_{m}(W_{m})\dans W_{\mm}$, auquel cas $\W:=\set \phi_{m}\rest_{W_{m}}:W_{m}\to W_{\mm}/$.  

L'intersection d'une famille de sous-$\FI$-modules est  un sous-$\FI$-module.

On définit dualement la notion de \expression{$\FI$-module quotient}.

\subsubsectionline{Systèmes générateurs de $\FI$-modules.}\label{generateurs-FI-modules}Soit $\V$ un $\FI$-module.
Pour tout sous-ensemble $\Sigma\dans\coprod_{m} V_{m}$, on appelle \expression{sous-$\FI$-module engendré par $\Sigma$} et on le note $\<\Sigma>$\glossary{$\<\Sigma>\dans V$:sous-$\FI$-module de $V$ engendré par $\Sigma\dans\coprod_{m} V_{m}$}, l'intersection de la famille des sous-$\FI$-modules $\W$ de $\V$ tels que $\Sigma\dans\coprod_{m} W_{m}$.

Pour $a\in\NN$, notons $\Sigma_{a}:=\Sigma\cap V_{a}$. D'après \ref{M(a)}, chaque $s\in\Sigma_{a}$ détermine un et un unique morphisme de $\FI$-modules 
$$\halfdisplayskips\ev_{\ag}(s):\M(\ag)\to \V\,,\quad\ev_{\ag}(s)(\id_{\ag})=s\,.$$
 Le lemme suivant est alors immédiat.\killline

\begin{lemm*}[(\cite{cef} 2.3.2)]Le sous-$\FI$-module $\<\Sigma>\dans \V$ est l'image du morphisme
$$\prodnl_{0<a\in\NN,s\in \Sigma_{a}}\ev_{\ag}(s):\bigoplusnl_{0<a;s\in \Sigma_{a}}\M(\ag)\to \V\,,\qquad $$
\smallskip\end{lemm*}

\subsubsectionline{$\FI$-modules de type fini.}Soit\label{type-fini} $\V$ un $\FI$-module.
\begin{enumerate}
\item\leavevmode\label{def-degre-<=d}$\V$ est dit \expression{engendré en degrés $\leq d$} si $V=\<\Sigma>$ avec $\Sigma\dans \coprod_{a\leq d} V_{a}$.
 
\item\leavevmode\label{def-type-fini}$\V$ est dit \expression{de type fini} 
si $\V=\<\Sigma>$ avec $\Sigma$ fini.

\end{enumerate}

\noendpoint\begin{prop}[(\cite{cef} Finitude et noethérianité)]\label{finitude-FI-modules}
\def\varlistskips{\topsep0pt\itemsep2pt}
\begin{enumerate}
\item Un\label{finitude-FI-modules-a} $\FI$-module $\V$ est de type fini si et seulement si, il admet une surjection $\bigoplusnl_{i}\M(\ag_i)\onto \V$ pour une certaine famille finie $\set \ag_i/$

\item Un\label{finitude-FI-modules-b} quotient d'un $\FI$-module de type fini est de type fini.
\item Si\label{finitude-FI-modules-c} l'anneau $\AAA$ des coefficients est noethérien, tout sous-$\FI$-module d'un $\FI$-module de type fini est encore de type fini.~\rm(\footnote{Prouvé dans \cite{cef} en supposant  $\QQ\dans \AAA$, et dans \cite{cefn} en général.})
\end{enumerate}
\end{prop}

\subsection{Caractères polynomiaux et stabilité des $\FI$-modules}\label{car-pol-rep-satb}

\smallskip
{\miseengarde\noindent
 L'anneau $\AAA$ est un corps $k$ de caractéristique nulle.\endgraf\vskip-0.5\baselineskip\vskip0pt}

\subsubsectionline{Caract\`{e}res (éventuellement) polynomiaux.}Si\label{caracteres-polynomiaux}
$\W:=\set\S _{m}\rep W_{m}/_{m}$\glossary{${\W:=\set\S _{m}\rep W_{m}/_{m}}$:famille de représentations de dimensions finies des groupes symétriques} est une famille de représentations de dimensions finies,  on note $\chi(\W):=\set\chi_{\S _{m}}(W_{m})/_{m}$\glossary{${\chi(\W):=\set\chi_{\Sm}(W_{m})/_{m}}$:famille des caractères associée à $\W$} la famille de leurs caractères. 
\comment
Ici, on ne suppose pas que $\W$ provient d'un $\FI$-module et il n'y a donc à priori aucun lien entre les $W_{m}$. \endcomment
Pour chaque $m$, il existe un polynôme $P_{m}\in k[\XX_1,\ldots,\XX_{m}]$ tel que pour $\alpha\in\S _{m}$, on a
$$\chi_{\S _{m}}(W_{m})(\alpha)=P_{m}(\XX_1,\ldots,\XX_{m})(\alpha)\,,
$$
où $\XX_i(\alpha)$ est le nombre des cycles de longueur $i$ dans la décomposition de $\alpha$ en produit de cycles disjoints. 
La famille $\chi(\W)$ peut ainsi être décrite par des familles $\set P_{m}/_{m}$ de polynômes de l'algèbre $k[\cl {\XX\,}]=k[\XX_1,\XX_2,\ldots]$. 

\medskip\noindent{\bf Définition. }La famille de représentations $\W$ est dite \expression{à caractère éventuellement polynomial} s'il existe $N,r\in\NN$ et
$P\in k[\XX_1,\ldots,\XX_r]$ tels que 
$$
\chi_{\S _{m}}(W_{m})(\alpha)=P(\XX_1,\ldots,\XX_r)(\alpha)\,,\quad\forall m\geq N\,,\ \forall\alpha\in\S _{m}\,.\eqno(\ast)$$
On dit alors que $\chi(\W)$ \expression{est polynomiale pour $m\geq N$}.

\noendpoint\begin{remas}\label{remas-caractere-poly}
\def\varlistskips{\topsep0pt\itemsep2pt\mou}
\begin{enumerate}\item
Un caractère du groupe $\Sm$ admet plusieurs écritures dans $k[\XX_1,\ldots,\XX_{m}]$ (\cf\ref{redondance}), par contre, lorsque la condition $(\ast)$ est satisfaite, le polynôme $P$ est  unique (\cf ($^{\hbox{\scriptsize\ref{unicité-poly}}}$), p.~\pageref{unicité-poly}).

\item Notons\label{suite-de-dimensions} $\dim_k(\W):=\set\dim_k W_{m}/_{m}$\glossary{${\dim_k(\W):=\set\dim_k(W_{m})/_{m}}$: famille des dimensions associée à $\W$}. Lorsque $\chi(\W)$ est éventuellement polynomiale, on a, pour $m$ assez grand,
$$\halfdisplayskips
\dim_{k}W_{m}=P(\XX_1,\XX_2,\ldots,\XX_r)(\II_{m})=P(m,0,\ldots,0)$$
et la suite $\dim_{k}(\W)$ est éventuellement polynomiale en $m$. \end{enumerate}\end{remas}

\comment
\killline L'un des intérêts principaux de la finitude des $\FI$-modules réside dans le résultat suivant, corollaire du \emph{théorème de stabilité des $\FI$-modules de type fini}, que nous rappellerons dans \ref{FI-TF=>rep-stat} (voir aussi \ref{827y823=>824}).\killline

\begin{theo}[(\cite{cef} Polynomialité de caractères)]\label{FI-TF=>chi-pol}Un $\FI$-module de type fini $\V$ est à caractère éventuellement polynomial.
\end{theo}
\endcomment

\subsubsectionline{Représentations irréductibles $V(\lambda)_{m}$.}Sous\label{car(k)=0} l'hypothèse en cours de $\car(k)=0$, les représentations des groupes symétriques sont semi-simples et \emph{définies sur $\QQ$}. Les décompositions en composantes irréductibles sont alors ``indépendantes'' de $k\cont\QQ$.

\subsubsection{Reparamétrisation des représentions irréductibles}Pour\label{reparametrisation} $0<a\in\NN$, les représentations irréductibles de $\S _{a}$ sur $k$ sont paramétrées par les décompositions $\lambda\vdash a$. On note $V_{\lambda}$ la représentation irréductible de $\S _{a}$ correspondante à $\lambda\vdash a$.

\smallskip{\bf Définition. }Soit $\lambda=(\lambda_1\geq\cdots\geq\lambda_{\ell})\vdash a$. Pour tout $m\geq |\lambda|+\lambda_1$, on note $\lambda[m]\vdash m$ la décomposition $(m-|\lambda|,\lambda_1,\ldots,\lambda_\ell)$
$$
\yght=3mm\ygwd=3mm
\let\cd\cdot
\def\m#1{\def\donne{#1}\def\vx{\vdots}\ifx\donne\vx\rlap{\kern2.6cm$\vdots$}\else\rlap{\kern2.4cm$(#1)$}\fi}
\null\qquad
\raise-0pt\hbox{$\tableauc{}{\cd&\cd&\cd&\cd\\
\cd\smash{\raise-1pt\hbox{\normalsize\mllap{\lambda=\left\{\vrule height0.50cm width0pt\right.\kern1pt\ \mrlap{\kern1.05cm\left.\vrule height0.50cm width0pt\right\}}}}}
&\cd\\\cd}$}\qquad
\raise-0pt\hbox{$\longmapsto$}\quad
\lambda[m]:=\left\{\tableauc{\cellcolor{0 0 0 0.10}}{
\m{m-|\lambda|}&&&&\noborder\cdots&&\\\cellcolor{0 0 0 0}
\cd\m{\lambda_1}&\cd&\cd&\cd\\
\cd\relax\m{\vdots}&\cd\\
\cd\m{\lambda_\ell}}
\hskip1.5cm\right\}\,,
$$

Pour $m\geq|\lambda|+\lambda_1$, on note $V(\lambda)_{m}$\glossary{$V(\lambda)_{m}$: la représentation irréductible de $\S _{m}$ correspondante à $\lambda[m]$} la représentation irréductible de $\S _{m}$ correspondante à $\lambda[m]$, on pose donc:
$$\displayboxit{V(\lambda)_{m}:=V_{\lambda[m]}}$$
Les notations $V(0)_{m}$ et $V(1)_{m}=V(\Box)_{m}$ désignent ainsi respectivement la re\-pré\-sentation {\mathrigid0mu 
triviale et la re\-pré\-sentation standard de dimension $
\mmo$ de~$\S _{m}$.}

\noendpoint\begin{remas}\label{remas-lambda[m]}
\begin{enumerate}
\item\leavevmode\label{remas-lambda[m]-a}Une décomposition $\mu\vdash m$ s'écrit d'une et d'une {\sl unique\/} manière sous la forme $\mu=\lambda[m]$. Si $\mu=(\mu_1\geq\mu_2\geq\cdots\geq \mu_{\ell})$, on a 
$\lambda=(\mu_2\geq\cdots\geq \mu_{\ell})$ et $\lambda\vdash a:=m-\lambda_1$ (y compris $a=0$).

\item\leavevmode\label{remas-lambda[m]-b}Dans la notation 
`$\lambda[m]$', le nombre `$m$', qui vérifie $m\geq |\lambda|+\lambda_1$, indique la taille finale du diagramme. On a 
$$\ell(\lambda[m])=\ell(\lambda)+1\,,
$$ 
et si $m\leq n$, on a
$\preskip-1ex
(\lambda[m])[n]\not=\lambda[n]\,.$

\item\leavevmode\label{remas-lambda[m]-c}Si $\lambda[m]\in\Y_{\ell}(m)$, on a
$\displaystyle|\lambda|\leq m-\Big\lceil{m\over \ell}\Big\rceil\,.$
~(\footnote{On désigne par $\lceil x\rceil\in\ZZ$ la partie entière par excès de $x\in\RR$, \idest $\lceil x\rceil-1<x\leq\lceil x\rceil $.
Dans (\ref{remas-lambda[m]-c}), le terme $\lceil m/\ell\rceil$ représente le plus petit nombre de colonnes d'un diagramme de Young de taille $m$ possédant $\ell$ lignes, il faudrait donc convenir que $\lceil 0/0\rceil=0$.
})
\glossary{$\lceil x\rceil\in\ZZ$:partie entière par excès de $x\in\RR$, \idest $\lceil x\rceil-1<x\leq\lceil x\rceil $. Le nombre $\lceil m/\ell\rceil$ est le plus petit nombre de colonnes d'un diagramme de Young de taille $m$ possédant $\ell$ lignes, on convient que $\lceil 0/0\rceil=0$}

\end{enumerate}
\end{remas}

\subsubsectionnumber
Les familles $\set V(\lambda)_{m}/_{m}$ apparaissaient déjà dans le livre de Macdonald \cite{mac-livre} comme exemple de famille de représentations à caractère polynomial. L'auteur y établit par un calcul explicite le fait suivant (ex. I.7.14, p.~122).
\killline

\def\theotitlestyle{\sl}\begin{prop*}[(\cite{mac-livre})]\label{macdonald-poly}Pour $\lambda\vdash a$, 
il existe  $P_{\lambda}\in k[\cl X]$ tel que, pour tout
$m\geq |\lambda|+\lambda_1$, on a
$$P_{\lambda}(\alpha)=\chi_{\S _{m}}(V(\lambda)_{m})(\alpha)\,,\mrlap{\quad\forall\alpha\in\S _{m}.}$$
\end{prop*}

\subsection{Monotonie et stabilité des $\FI$-modules}\label{stabilite-et-monotonie}\label{monotonie-et-stabilite}

\subsubsectionline{Stabilité.}Suivant \cite{cef}\label{rang-stabilite}\label{repr-stables}, 
un $\FI$-module $\V=\set \phi_{m}:V_{m}\to V_{\mm}/$ est dit
\expression{(uniformément) stationnaire}, s'il existe $N\in\NN$ tel que les conditions suivantes sont satisfaites pour $m\geq N$.
\def\varlistskips{\topsep2pt\itemsep2pt\mou}
\begin{enumerate}
\item{\sl Injectivité. }\label{repr-stables-a}Les applications $\phi_{m}:V_{m}\to V_{\mm}$ sont injectives.
\item{\sl Surjectivité. }\label{repr-stables-b}$V_{\mm}$ est engendré en tant que $\S _{\mm}$-module par $\im(\phi_{m})$.
\item{\sl Multiplicités. }\label{repr-stables-c}(\footnote{L'adverbe ``uniformément'' est utilisé dans \cite{cef} pour distinguer du cas
où la condition de stabilité de multiplicités (\ref{repr-stables-c}) est demandée séparément pour chaque $\lambda$. }) Dans la décomposition en représentations irréductibles
$$
V_{m}=\bigoplusnl_{\lambda}c(\lambda)_m\,V(\lambda)_{m}\,,$$
les multiplicités $c(\lambda)_m$\glossary{$c(\lambda)_m$:multiplicité de $V(\lambda)_{m}$ dans $V_m$} ne dépendent pas de~$m$. 
\comment En particulier, pour $n\geq N$, on a
$$V_{n}=\bigoplusnl_{\lambda\vdash m\leq N-\lambda_1}\mu_{\lambda,N}\,V(\lambda)_{n}\,.$$\endcomment
\end{enumerate}

\medskip
Lorsque ces conditions sont satisfaites, on dit que $\V$ est un $\FI$-module \expression{stable pour $m\geq N$}.
Le \expression{rang de stabilité de $\V$}, noté $\rks(\V)$\glossary{$\rks(\V)$:rang de stabilité de $\V$}, est  la borne inférieure de tels $N$, et si $\V$ n'est pas stationnaire, on pose $\rks(\V)=+\infty$.

\smallskip\noindent{\bf Fait: }la famille de caractères $\chi(\V)$ est polynomiale pour $m\geq\rkms(\V)$ (\ref{macdonald-poly}).
\subsubsectionline{Monotonie.}Suivant\label{monotonie} la terminologie de Church (\cite{chu}), un $\FI$-module $\V=\set \phi_{m}:V_{m}\to V_{\mm}/_{m}$ est dit \expression{monotone pour $m\geq N$} lorsque pour tout $m\geq N$ et tout $W_m\dans V_m$ isomorphe à $V(\lambda)_{m}^{\oplus k}$ pour certain $\lambda$ et $k$, le sous-$\Smm$-module $\Smm\cdot\phi_{m}(W)\dans V_{\mm}$ contient un sous-module
isomorphe à $V(\lambda)_{\mm}^{\oplus k}$. On notera $\rkm(\V)$\glossary{$\rkm(\V)$:rang de monotonie de $\V$} le plus petit des ces $N$.

\smallskip\noindent{\bf Fait: }les morphismes de transition de $\V$ sont injectifs pour \hbox{$m\geq\rkm(\V)$.}

\begin{lemm}[et définition]
On appellera\label{lemme-rkms}\glossary{$\rkms(\V){:=}\sup\set \rkm(\V),\rks(\V)/$:rang (de monotonie et stabilité) de $\V$} \expression{rang de  monotonie et stabilité} d'un $\FI$-module $\V$, le nombre
$$\displayboxit{\rkms(\V):=\sup\set \rkm(\V),\rks(\V)/}$$
Pour tout $\FI$-module $\V$ et tout $s\in\NN$,
on a $\rkms(\V)\leq\sup\set\rkms(\V_{\geq s}),s/$, et
 pour $\rkms(\V_{\geq s})$, seul deux cas sont possibles
\begin{liste-roman}
\item\leavevmode\label{lemme-rkms-i}
$\rkms(\V_{\geq s})=\sup\set \rkms(\V),s/$\,, ou bien
\item\leavevmode\label{lemme-rkms-ii}
 $\rkms(\V_{\geq s})=0\,,$ auquel cas 
$\V_{\geq\rkms (\V)}=0$ et
$s\geq \rkms (\V)$.
\end{liste-roman}
\end{lemm}
\demo Notons $r:=\rkms (\V)$. L'assertion (\ref{lemme-rkms-i}) est claire et immédiate lorsque $\V_{\geq r}\not=0$ et aussi lorsque $\V_{\geq r}=0$ et que $s< r$. Le cas restant est celui où $\V_{\geq r}=0$ et que $s\geq r$, cas auquel $\V_{\geq s}=0$ et donc  $\rkms(\V_{\geq s})=0.$\enddemo


\begin{prop}\def\MMMS{{\MMM\SSS}}Dans\label{prop-generalites-monotonie-stabilite} les énoncés qui suivent, `monotone' est un raccourci pour `monotone pour $m\geq N$', et de même pour `monotone-stable'.
\def\varlistskips{\topsep0pt\itemsep2pt\mou}
\begin{enumerate}
\item Un\label{prop-generalites-monotonie-stabilite-a} sous-$\FI$-module d'un $\FI$-module monotone est monotone.

\item Tout quotient\label{prop-generalites-monotonie-stabilite-b} monotone $\Q$ d'un $\FI$-module monotone-stable $\V$ est mo\-no\-tone-stable ainsi que le noyau de la surjection canonique $\nu:\V\onto\Q$.

\item Soit\label{prop-generalites-monotonie-stabilite-c} 
$\W\to\V\mathop{\to}\limits\Q\to0$ une suite exacte de $\FI$-modules où $\W$ est monotone-stable. Alors, si $\V$ est monotone (resp. monotone-stable), $\Q$ l'est aussi.

\item Soit\label{prop-generalites-monotonie-stabilite-d} 
$0\to\K\to\V\to\Q\to0$ une suite exacte de $\FI$-modules. Alors
Si $\K$ et $\Q$ sont monotones (resp. monotones-stables), $\V$ l'est aussi.

\item Soit\label{prop-generalites-monotonie-stabilite-e} 
$\cdots\to\V_{-1}\too^{d_{-1}}\V_{0}\too^{d_0}\V_{1}\to\cdots$ un complexe de $\FI$-modules. Alors
\begin{enumerate}
\item Si $\V_{-1}$ est monotone-stable et $\V_0$ est monotone, $H^{0}(\V_*)$ est monotone.
\item Si $\V_{-1}$ et $\V_0$ sont monotones-stables et $\V_1$ est monotone, $H^{0}(\V_*)$ est monotone-stable.
\end{enumerate}
\end{enumerate}
\end{prop}
\demo (\ref{prop-generalites-monotonie-stabilite-a}) Immédiat d'après la définition de monotonie. 

\parskip0pt\mou

\smallskip
(\ref{prop-generalites-monotonie-stabilite-b}) Le noyau $\K:=\ker(\nu)$ est monotone d'après (\ref{prop-generalites-monotonie-stabilite-a}). 
La monotonie assure la condition d'\emph{injectivité} 
\ref{repr-stables}-(\ref{repr-stables-a})
pour $\K$, $\V$ et $\Q$, mais aussi le fait que les multiplicités $ c(\lambda)_{m}$ pour $K_m$, $V_m$ et $Q_m$ sont non décroissantes pour $m\geq N$. Il s'ensuit que si $\V$ est monotone-stable, les familles de multiplicités sont constantes et la condition de \emph{multiplicités} \ref{repr-stables}-(\ref{repr-stables-c}) est  satisfaite. Enfin, la condition de \emph{surjectivité} \ref{repr-stables}-(\ref{repr-stables-b}) est une conséquence logique des propriétés de \emph{monotonie} et de \emph{multiplicités constantes}.

\smallskip
(\ref{prop-generalites-monotonie-stabilite-c}) Pour la monotonie, il suffit, pour $\lambda$ et $n\geq N$ donnés, de se restreindre au sous-$\FI$-module $\Q'\dans\Q$ engendré par $V(\lambda)_{n}^{\oplus k}\dans Q_n$. On note  $\V'$ le sous-$\FI$-module de $\V$ défini par$V'_{m}=\pi_{m}^{-1}(Q'_m)$, d'où la suite exacte $\W\to\V'\to\Q'\to0$ avec $\V'$ monotone d'après (\ref{prop-generalites-monotonie-stabilite-a}).
On considère ensuite le diagramme des morphismes de transition restreint aux composantes $\lambda$-isotypiques en $n$ 
$$\xymatrix@R=5mm{
\strut\K\arhook[d]&V(\lambda)_{n}^{\oplus l}\ar[r]^(.55){\phi_n}\strut\arhook[d]&K''_{n}\mrlap{\dans K_{n{+}1}}
\strut\arhook@<3.em>[d]\\
\V'\aronto[d]_(0.45){\nu}&V(\lambda)_{n}^{\oplus k}\oplus V(\lambda)_{n}^{\oplus l}
\aronto[d]\ar[r]^(.70){\phi_n}&V''_{n}\mrlap{\dans V'_{n{+}1}}\aronto@<3.em>[d]\\
\Q'&V(\lambda)_{n}^{\oplus k}\ar[r]^(0.6){\phi_n}&Q''_{n}\mrlap{\dans Q'_{n{+}1}}\\
}$$
La composante $\lambda$-isotypique de $Q''_{n}:=\S_{n{+}1}\cdot\phi_{n}(\V(\lambda)_{n}^{\oplus k})\dans Q'_{n+1}$ est celle de $V''_{n}:=\S_{n{+}1}\cdot\phi_{n}(\V(\lambda)_{n}^{\oplus k+l})$ modulo 
celle de $K''_{n}:=\S_{n{+}1}\cdot\phi_{n}(\V(\lambda)_{n}^{\oplus l})$. Or, cette dernière est de multiplicité exactement $l$ puisque $\K$ est supposé monotone-stable, tandis que celle en $V''_{n}$ est de multiplicité $\geq k{+}l$ puisque $\V'$ est supposée  (seulement) monotone. 
Par conséquent $\Q$ vérifie la condition de monotonie pour $\lambda$ en $n$. Cette conclusion étant valable pour tout $n\geq N$ et tout $\lambda$, le $\FI$-module $\Q$ est bien monotone. 

Maintenant, si l'on suppose $\V$ monotone \emph{et} stable, $\Q$ est monotone d'après ce qui précède et, donc, monotone-stable par (\ref{prop-generalites-monotonie-stabilite-b}).

\smallskip
(\ref{prop-generalites-monotonie-stabilite-d}). Même type de raisonnements que pour (\ref{prop-generalites-monotonie-stabilite-c}).

\smallskip
(\ref{prop-generalites-monotonie-stabilite-e}). (i) Si $\V_0$ est monotone, $\ker(d_0)$ est monotone par (\ref{prop-generalites-monotonie-stabilite-a}), $\im(d_{-1})$ est monotone-stable par 
(\ref{prop-generalites-monotonie-stabilite-b}), et on applique  
(\ref{prop-generalites-monotonie-stabilite-c}) à $\im(d_{-1})\hook\ker(d_0)\onto H^{0}(\V_*)$. (ii) Si $\V_0$ est monotone-stable et que $\V_1$ est monotone, $\ker(d_0)$ est monotone-stable par (\ref{prop-generalites-monotonie-stabilite-a}) et $H^{0}(\V_*)$ est monotone-stable à nouveau par (\ref{prop-generalites-monotonie-stabilite-c}).
\enddemo


\bigskip
La proposition suivante est une reformulation des résultats de Church concernant les concepts de monotonie et de stabilité (\loccit\ prop. 2.5, p.~475). 

\noendpoint\begin{coro}\label{coro-monotonie-stabilite}
\def\varlistskips{\topsep0pt\itemsep0pt\mou}
\begin{enumerate}
\mynobreak\item La\label{coro-monotonie-stabilite-a}
sous-catégorie pleine $\ModFIsc{\rkm\leq N} k$ des $\FI$-modules 
monotones pour $m\geq N$,  est stable par extensions dans $\ModFI k$.

\item La\label{coro-monotonie-stabilite-b}
sous-catégorie pleine $\ModFIsc{\rkms\leq N} k$ des $\FI$-modules 
monotones et stables pour $m\geq N$,  est une sous-catégorie abélienne
et stable par extensions dans $\ModFI k$. En particulier,
\def\varlistskips{\topsep0pt\itemsep2pt\mou}
\begin{enumerate}
\item Si\label{coro-monotonie-stabilite-bi} 
$(\V_{*},d_{*})$ est un complexe de $\ModFIsc{\rkms\leq N} k$,
ses $\FI$-modules de cohomologie $h^{i}(\V_{*},d_{*})$ vérifient aussi
$\rkms(h^{i}(\V_{*},d_{*}))\leq N\,.$
\item Dans\label{coro-monotonie-stabilite-bii} le groupe $K_0(\ModFI k)$, un $\FI$-module $\V$ appartient au sous-groupe $K_{0}(\ModFIsc{\rkms\leq N} k)$, si et seulement si $\rkms(\V)\leq N$.
\end{enumerate}
\end{enumerate}
\end{coro}

\def\Remarque{Commentaire}\begin{rema}Un\label{comment-poids-degre} théorème important de la théorie de Church et Farb établit l'équivalence pour un $\FI$-module entre le fait d'être de \emph{type fini} et être \emph{éventuellement monotone et stable}
 et (donc) à caractère  polynomial (\cite{cef}, thm.~1.13), mais ce type d'affirmation, de nature uniquement qualitative, ne renseigne pas sur les valeurs des rangs concernés. Les notions de \expression{poids} et de \expression{degré de stabilité} d'un $\FI$-module ont été introduites dans le but de préciser cette question (\cf\ref{rkms=rke}). Les sections \ref{poids}--\ref{deg-stab} qui vont suivre rappellent ces notions.\end{rema}

\comment 
Ensuite, les propositions \ref{critere-stabrk} et \ref{theo-Ind-lambda} indiquent
comment $\poids_{\infty}(\V)$ participe à la majoration de $\rkms(\V)$. La proposition suivante, qui intervient aussi dans la preuve de \ref{theo-Ind-lambda}, donne un majorant du poids d'une représentation induite par le foncteur d'induction $\ind_{\Glambda }^{\Sm}(\_)$ de la section \ref{operateurs-inductions}.

Elle nous sera particulièrement utile dans la preuve du théorème \ref{theo-Ind-lambda}.\endcomment

\subsection{Poids d'un $\FI$-module}\label{poids}
\subsubsectionline{Règles de branchement pour les représentations de $\Sm$.}Elles\label{regles-branchement} expliquent les manipulations à effectuer sur le diagramme de Young $\mu$ correspondant à une représentation simple $V_\mu$ d'un groupe symétrique, pour décrire les composantes irréductibles de ses induits et ses restrictions. Le lemme suivant rappelle ces règles (\cite{cef}, lemma 3.2.3). (\footnote{Une\label{note-LR} bonne référence pour ces questions est \cite{fulhar}, exercices 4.42-45, pp.~57-59. Les règles de branchement découlent plus généralement de la règle de Littlewood-Richardson\glossary{Règles de Littlewood-Richardson, règles de Pieri, règles de branchement$$:} \loccit\ appendice A, eq. (A.8) p.~456.})

\begin{lemm}Le\label{regles-de-branchement} corps $k$ est de caractéristique nulle.
\begin{enumerate}
\item \relax{\rm [Règle de Pieri]} Soit\label{regles-de-branchement-a} $V_\nu$ la représentation irréductible de $\Sa$ correspondante à un diagramme de Young $\nu\vdash a$. On a:
$$\ind_{\Sb \times\Sa }^{\S_{a+b}} k_{b}\boxtimes V_\nu\sim
\ind_{\Sa \times\Sb }^{\S_{a+b}} V_\nu\boxtimes k_{b} =\bigoplusnl_{\mu\vdash a+b} V_{\mu}\,,$$
où les diagrammes $\mu\vdash a+b$ sont ceux qui s'obtiennent en rajoutant une boite sur $b$ colonnes distinctes du diagramme $\nu$. 
\item Soit\label{regles-de-branchement-b}
$V_\mu$ la représentation irréductible de $\S_{a+b}$ correspondante à un diagramme de Young $\mu\vdash a+b$. On a:
$$
\big(\Res_{\Sa \times\Sb }^{\S_{a+b}} V_\mu\big)_{\Sb}=\bigoplusnl_{\nu\vdash a} V_{\nu}\,,$$
où les diagrammes $\nu\vdash a$ sont ceux qui s'obtiennent en enlevant une boite sur $b$ colonnes distinctes du diagramme $\mu$. 
\end{enumerate}
\end{lemm}

\subsubsection{Définitions}\label{def-poids}
\def\varitemizeseps{\itemsep0pt\parskip0pt\topsep1pt}\begin{itemize}
\mynobreak\nobreak\item Soit $m>0$. Le \expression{poids} d'un $\Sm $-module $W$, noté $\poids(W)$, est le plus grand des $|\lambda|$ tels que $V_{\lambda[m]}$ est facteur irréductible de $W$. Si $W=0$, on pose $\poids(W):=0$. 
On a toujours
$\poids(W)< m$ (\ref{remas-lambda[m]}(\ref{remas-lambda[m]-c})).

\item Le \expression{poids}
d'un $\FI$-module $\V=\set V_m/$, noté $\poids(\V)$\glossary{$\poids (\V)$:poids de $\V$}, est la borne supérieure de l'ensemble $\set\poids (V_m)/$ dans $\NN\cup\set+\infty/$.

\comment
\item Le \expression{poids à l'$\infty$}
d'un $\FI$-module $\V=\set V_m/$, noté $\poidsinf(\V)$\glossary{${\poids_{\infty} (\V):=\lim_{m\mapsto+\infty}\poids V_m}$:poids à l'infini de $\V$}, est la limite dans $\NN\cup\set+\infty/$, si elle existe, de $\poids(V_m)$ lorsque $m\mapsto+\infty$.
Si $\V$ est un $\FI$-module avec rang de stabilité fini, on a 
$\poids_{\infty}(\V)\leq\poids(\V)\leq\rks(\V)\,.$

avec égalités seulement si $\poids(\V)=0$, \idest si $V_m\dans V(0)_m$, pour tout $m$. \endcomment
\end{itemize}

\bigskip

\begingroup\def\b{\ell-\uell}\def\Sb{\S_{\b}}
\begin{prop}Pour\label{poids-ind} $\ell,a\in\NN$, soit
$\lambda=(\lambda_1\geq\cdots\geq\lambda_{\ell})\vdash \ell+a$ et notons  $\uell:=\#\set i\mid\lambda_i>1/$. 
Pour tout $\S_{\ell}$-module $W$, on a
$$
\displayboxit{\poids\big(\Ind_{\Glambda }^{\S_{\ell+a}}W\big)
\leq \poids(W)+\uell+a}$$
\end{prop}
\demo\halfdisplayskips
On a 
$$\smash{
\let\Big\big
\Plambda
\ \unlhd\ 
\Big(S_{\lambda_1}\times\cdots\times\StS{\lambda_\uell}{\b}\Big)
\ \unlhd_{(1)}\ 
N_{\S_{|\lambda|}}\Plambda =\Slambda\,,}
\eqno(\dagger)$$
d'où une surjection de $\S_{|\lambda|}$-modules
$$
M:=\ind_{S_{\lambda_1}\times\cdots\times\StS{\lambda_\uell}{\b}}^{\S_{|\lambda|}}W
\onto
\ind_{\Glambda }^{\S_{|\lambda|}}W\,.
\eqno
(\ddagger)$$
(C'est même un isomorphisme lorsque $(\lambda_1,\ldots,\lambda_{\uell})$ est strictement décroissante, car, dans ce cas, l'inclusion $\unlhd_{(1)}$ dans  $(\dagger)$ est une égalité.)
On remarque alors que le terme de gauche de $(\ddagger)$ n'est autre que le $\S_{\ell+a}$-module
$$M=\ind^{\S_{\lambda_1+\cdots+\lambda_{\uell}+\b}}_{
S_{\lambda_1}\times\cdots\times\StS{\lambda_\uell}{\b}}k_{\lambda_1}\boxtimes\cdots\boxtimes
k_{\lambda_\uell}\boxtimes\Res^{\S_{\ell}}_{\Sb}W\,,$$
puisque dans l'action de $S_{\lambda_1}\times\cdots\times\StS{\lambda_\uell}{\b}$ à travers $\Glambda $, les $\S_{\lambda_i}$ agissent trivialement sur $W$. L'itération de la règle de Pieri \ref{regles-de-branchement}-(\ref{regles-de-branchement-a}) conduit alors à la majoration
$$\poids(M)\leq\lambda_1+\cdots+\lambda_\uell+\poids(\Res^{\S_{\ell}}_{\1_{\uell}\times\Sb}W)\,,$$
avec $\lambda_1+\cdots+\lambda_\uell=\uell+a$ et où
$$\poids(\Res^{\S_{\ell}}_{\1_{\uell}\times\Sb}W)\leq\poids(W)\,.\eqno(\diamond)$$
En effet, un cas particulier de \ref{regles-de-branchement}-(\ref{regles-de-branchement-b}) dit que $\poids(\Res^{\S_{1+n}}_{\1_{1}\times\S_{n}} M)\leq\poids(M)\,,$
pour tout $n\in\NN$ et tout $\S_{1+n}$-module $M$, ce qui conduit
à $(\diamond)$ par itération. On a donc montré que $\poids(M)\leq\poids(W)+\uell+a$ et la proposition résulte puisque $\poids(M)$ majore le poids des quotients de $M$, en particulier $(\ddagger)$.
\enddemo

\begin{rema}\halfdisplayskips
Dans\label{rema-poids-ind} \ref{poids-ind}, l'égalité $\poids\big(\ind_{\Glambda }^{\S_{|\lambda|}}W\big)
=\poids(W)+\uell+a$ peut être atteinte. 
En effet, soit $W=V(\nu)_{\ell}$.
Par la règle \ref{regles-de-branchement}-(\ref{regles-de-branchement-b}), $V(\nu)_{\ell-\uell}$ est facteur de $\Res^{\S_{\ell}}_{\S_{\ell-\uell}}V(\nu)_{\ell}$
si et seulement si 
$$\ell-\uell-|\nu|\geq\nu_1\,,\eqno(\ast)$$
auquel cas  $\poids\big(\Res^{\S_{\ell}}_{\S_{\ell-\uell}}V(\nu)_{\ell}\big)=|\nu|$ est le plus grand poids possible. 

Notons 
$M:=\ind_{S_{\lambda_1}\times\cdots\times\StS{\lambda_\uell}{\ell-\uell}}^{\S_{|\lambda|}}V(\nu)_{\ell-\uell}$. Par ce qui précède et
par itération de la règle de Pieri \ref{regles-de-branchement}-(\ref{regles-de-branchement-a}), on aura
$\poids(M)=|\nu|+\uell+a$, si et seulement si
$\ell-\uell-|\nu|\geq\sup\set\nu_1,\lambda_1/\,.$
Or, cette condition est vérifiée si $\lambda=(\nu_1,\nu_{1}-1,\nu_{1}-2,\ldots,1,\ldots,1)$, auquel cas $\uell=\nu_1-1$ et $(*)$ est vérifiée dès que $\ell\geq |\nu|+2\nu_1-1$. Enfin, pour un tel $\lambda$, on a bien $\S_{\lambda}=\PPP_{\lambda}$ et donc $\ind_{G_{\lambda}}^{\S_{|\lambda|}}V(\nu)_{\ell}=M$.
\end{rema}
\endgroup

\subsubsectionline{Scindage d'un $\FI$-module par le poids.}Soit\label{scindage-par-le-poids}
$$W_m=\bigoplusnl_{|\mu|<m} c(\mu)_m V(\mu)_m$$
la décomposition en facteurs irréductibles d'un $\Sm$-module $W_m$.

Pour $t\in\NN$, on pose
$$W_m^{(\geq t)}:=\bigoplusnl_{|\mu|\geq t} c(\mu)_m V(\mu)_m\,,$$
et \emph{mutatis mutandis} pour $W_m^{(> t)}$, $W_m^{(\leq t)}$, $W_m^{(< t)}$ et $W_m^{(=t)}$, aussi noté $W_m^{(t)}$.

\begin{prop}[et définitions]Soit\label{prop-scindage-par-le-poids} $\W:=\set\phi_m:W_m\to W_\mm/_m$ un $\FI$-module.
\begin{enumerate}
\item On\label{prop-scindage-par-le-poids-a} a 
$\phi_m(W_m^{> t})\dans W_\mm^{> t}$ pour tout $t\in\NN$. On définit alors
\begin{enumerate}
\item le\label{prop-scindage-par-le-poids-a-i} sous-$\FI$-module $\W^{(>t)}\dans\W$ déterminé par les restrictions $\phi_{m}^{(>t)}$ des morphismes de transition $\phi_m$ aux sous-modules $W_m^{(>t)}$, soit\glossary{$\W^{(> t)}$, $\W^{(\leq t)}$, $\W^{(t)}$:scindés par le poids de $\W$}
$$\W^{(> t)}:=\bigset\phi^{(>t)}_m:W_m^{(> t)}\to W_\mm^{(> t)}/_m\,\text{\sl;}$$
\item le\label{prop-scindage-par-le-poids-a-ii} $\FI$-module quotient $\W/\W^{(>t)}$
$$\W^{(\leq t)}:=\bigset\phi^{(\leq t)}_m:W_m^{(\leq t)}\to W_\mm^{(\leq t)}/_m
:=\smashtop{\W\over\W^{(>t)}}\,\text{\sl;}$$
\item le\label{prop-scindage-par-le-poids-a-iii} $\FI$-module $\W^{(t)}$ sous-quotient de $\W$ \expression{de poids unique $t$}
$$
\W^{(t)}:=\bigset\phi^{(t)}_m:W_m^{(t)}\to W_\mm^{(t)}/_m
:={\W^{\geq t}\over\W^{>t}}\,.$$ 
\end{enumerate}

\noindent 
(Dans tous les cas, on aura remarqué l'abus de notation pour les morphismes de transition induits\dots)

\item Un\label{prop-scindage-par-le-poids-b} $\FI$-module $\W$
tel que $\poids(\W)<+\infty$, est extension successive de ses $\FI$-modules à poids unique $\W^{(t)}$, $t\leq\poids(\W)$.

\item On\label{prop-scindage-par-le-poids-c} a
\comment
$$
\rkms(\W^{(>t)})\leq\rkms(\W)\,,\quad \rkms(\W^{(\leq t)})\leq\rkms(\W)
$$\endcomment
$
\rkms(\W)=\sup\bigset \rkms(\W^{(> t)}),\rkms(\W^{(\leq t)})/
$,  $\forall t\in\NN$. En particulier
$$
\rkms(\W)=\sup\bigset \rkms(\W^{(t)})\mid t\in\NN/\,.
$$
Et de même avec $\rkm$ à la place de $\rkms$.
\end{enumerate}
\end{prop}
\demo (\ref{prop-scindage-par-le-poids-a}) Si
$V(\mu)_m$ est un facteur irréductible de $W_m$, le sous-$\Smm$-module engendré par $\phi_m(V(\mu)_{m})$ est un quotient de $\ind_{\smashtop\Sm}^{\Smm}V(\mu)_m$. Or, la règle de Pieri \ref{regles-de-branchement}-(\ref{regles-de-branchement-a}) stipule que les composantes irréductibles de $\ind_{\smashtop\Sm}^{\Smm}V(\mu)_m$ sont les $V(\nu)_\mm$ dont le diagramme de Young est obtenu en rajoutant une boite sur une colonne du diagramme de $V(\mu)_m$. Ces composantes sont donc de multiplicité $1$ et de poids $|\nu|\geq|\mu|$, (avec égalité $|\nu|=|\mu|$ si et seulement si $\mu=\nu$). Ceci établit l'inclusion $\phi_m(W_m^{> t})\dans W_\mm^{> t}$. La suite de (\ref{prop-scindage-par-le-poids-a}) est alors immédiate.

(\ref{prop-scindage-par-le-poids-b}) Lorsque $t:=\poids(\W)<+\infty$, on a $\W^{(\geq t)}\not=0$ et $\W^{(> t)}=0$ auquel cas $\W^{(\geq t)}=\W^{(t)}$. On donc la suite exacte courte de $\FI$-modules
$$0\to\W^{(t)}\hook\W\to\W^{(<t)}=\W/\W^{\geq t}\to0$$
où, par construction, $\poids(\W^{(<t)})<t$. Un raisonnement par induction sur le poids termine la preuve.

(\ref{prop-scindage-par-le-poids-c}) On a
$\rkms(\W)=\sup\bigset \rkms(\W^{(> t)}),\rkms(\W^{(\leq t)})/$
par la définition même de $\rkms$. L'égalité \smash{$\rkms(\W)=\sup\bigset \rkms(\W^{(t)})\mid t\in\NN/$} s'ensuit par (\ref{prop-scindage-par-le-poids-b}) et par la stabilité par extensions  de $\ModFIsc{\rkms\leq N} k$ (\ref{coro-monotonie-stabilite}-\ref{coro-monotonie-stabilite-b}).
La preuve du cas $\rkm$ à la place de $\rkms$ est essentiellement la même moyennant \ref{coro-monotonie-stabilite}-(\ref{coro-monotonie-stabilite-a}).
\enddemo

\subsectionline{Les $\FI$-modules $\M_{a}^{H}(W)$.}
Pour\label{M-W} tout $0<a\in\NN$ et tout sous-groupe $\Hg\dans\S _{a}$, on définit le foncteur\glossary{${\M_{a}^{\Hg}:\Mod(k[H])\to\ModFI k}$:}
$$\M_{a}^{\Hg}:\Mod(k[H])\to\ModFI k_{\geq a}\postskip-0.5ex$$
par
$$\M_{a}^{\Hg}(W)_{a+b}=\begin{casesalign}
0\,,\hbox{ si $b<0$,}\hfill\\
\ind_{H\times\S _{b}}^{\S _{a+b}}W\boxtimes k_{b}\,,\hbox{ si $b\geq0$.}
\end{casesalign}$$

Pour tous $0\leq b_1\leq b_2$, l'égalité  $W\boxtimes k_{b_1}=W\boxtimes k_{b_2}$ induit le morphisme de $\Hg\times\S _{b_1}$-modules de 
\smashbot{$W\boxtimes k_{b_1}\to\ind_{\Hg\times\S _{b_2}}^{\S _{a+b_2}}W\boxtimes k_{b_2}$}, d'où le morphisme (canonique) de transition
$$\preskip4pt
\halfsmash{\phi_{a+b_{2},a+b_{1}}:
\ind_{\Hg\times\S _{b_1}}^{\S _{a+b_1}}W\boxtimes k_{b_1}
\to
\ind_{\Hg\times\S _{b_2}}^{\S _{a+b_2}}W\boxtimes k_{b_2}
}$$
dont l'image est clairement invariante sous l'action de $\1_{a+b_1}\times\S _{b_2-b_1} $. 
La famille $\M_{a}^{\Hg}(W):\set\phi_{m+1,m}:\M_{a}^{\Hg}(W)_{m}\to\M_{a}^{\Hg}(W)_{\mm}/$ définit donc bien un $\FI$-module. 

L'action du foncteur $\M_{a}^{\Hg}$ sur les morphismes suit le même principe d'induction et ne sera pas détaillée. 

\subsubsectionnumber Les $\FI$-modules $\M_{a}^{\S _{a}}(k[\S _{a}])$ et $\M(\ag)$  de \ref{M(a)} sont les mêmes.
Dans la suite on notera
$\M_a:=\M_{a}^{\S _{a}}$\glossary{${\M_a:=\M_{a}^{\S _{a}}}$:\emph{alias de notation}}.

\noendpoint\begin{prop}[(\cite{cef})]\label{prop-M(W)}
\begin{enumerate}\itemsep4pt\parskip0pt
\item\leavevmode\label{prop-M(W)-a}Le foncteur $\M_{a}^{\Hg}:\Mod(k[H])\to\ModFI k_{\geq a}$ est additif, exact et  fidèle. 
Il est aussi l'adjoint à gauche du foncteur \expression{d'évaluation en $\ag$}, \idest pour tout $\Hg$-module $W$ et tout $\FI$-module $\V$, on a:
$$\Hom_{\FI}(\M_{a}(W),\V)=\Hom_{\Hg}(W,V_{a})\,.$$
En particulier,  le $\FI$-module $\M_{a}(W)$ est un objet projectif de $\ModFI k$.
\item\leavevmode\label{prop-M(W)-b}Le foncteur $\M_{a}:=\M_{a}^{\S _{m}}$ est pleinement fidèle. 

\item\leavevmode\label{prop-M(W)-c}$\M_{a}^{\Hg}(W)$ est un $\FI$-module de type fini si et seulement si, $\dim_k W<\infty$.
\item\leavevmode\label{prop-M(W)-d}{\rm(\footnote{\cite{cef}, prop. 3.2.4.})} Pour tout $\Sa $-module $W$, on a
$\comment
\poids_{\infty}(\M_{a}(W))=\endcomment
\poids(\M_{a}(W))=a\,.$
\item\leavevmode\label{prop-M(W)-e}{\rm(\footnote{Hemmer \cite{hemmer} thm. 2.4, Church \cite{chu} 
\Spar5 
thm. 2.8, p.~494.})} Pour tout $\S _{a}$-module $W$, le $\FI$-module $\M_{a}^{H}(W)$ est éventuellement monotone et stationnaire. On a
$$\rkm(\M_{a}^{H}(W))\leq a\text{\quad et\quad}
\rks(\M_{a}^{H}(W))\leq2a\,.$$
 En particulier, si $\dim_{k} W<+\infty$, la famille de caractères $\chi(\M_{a}^{H}(W))$ est polynomiale à partir de $m=2a$.
\comment
Plus précisément, si
 $\M_{a}^{H}(W)_{2m}=\displaystyle\bigoplusnl_{\lambda\vdash n<2m} \mu_{\lambda}\,V(\lambda)_{2m}$ est la décomposition en composantes $\S _{2m}$-irréductibles, on a
$$\M_{a}^{H}(W)_{m+b}=\displaystyle\bigoplusnl_{\lambda\vdash n<2m} \mu_{\lambda}\,V(\lambda)_{m+b}\,,\quad\forall b\geq m\,.$$ 
\endcomment
\end{enumerate}
\end{prop}
\def\Demonstration{Indications}
\demo
(\ref{prop-M(W)-a},\ref{prop-M(W)-b},\ref{prop-M(W)-c}) sont immédiates. (\ref{prop-M(W)-d},\ref{prop-M(W)-e}) résultent d'une étude fine des règles de branchement \ref{regles-de-branchement} par Church (\loccit). Lorsque 
$W=V_\lambda$ on a même $\rks(\M_{a}(V_\lambda))=|\lambda|+\lambda_1$.
La dernière partie de (\ref{prop-M(W)-e}) découle du calcul de Macdonald \ref{macdonald-poly}.
\enddemo
\subsectionline{Catégorie des $\FB$-modules.}L'article\label{strategie-stabilite} \cite{cef} s'intéresse également aux représentations de la sous-catégorie pleine $\FB\dans\FI$\glossary{$\FB$:catégorie des ensembles {\bfit F\/}inis et des applications {\bfit B\/}ijectives} des ensembles finis et leurs bijections. Un $k[\FB]$-module est alors, par définition, un foncteur covariant $\R:\FB\to\Vec(k)$, ce qui équivaut à la donnée d'une famille  de représentations de groupes finis $\R:=\set\rho_{m}:\S _{m}\mapsto\Gl_{k}(W_{m})/_{m}$. On a  donc
$$\ModFB k=\prodnl_{m\in\NN}\Mod(k[\S _{m}])\,.$$

Dans $\ModFB k$, la notion de $\FB$-module de type fini est  inintéressante car équivalente à la donnée d'une famille de représentions $\set W_{m}/_{m}$ de dimensions finies presque toutes nulles. Par contre, la notion de \expression{stabilité des multiplicités des représentations $V(\lambda)_{m}$} (\ref{repr-stables}-(\ref{repr-stables-c})) garde tout son intérêt. 

\medskip\noindent{\bf Définition. }On dit que un $\FB$-module est \expression{(éventuellement) stationnaire} lorsque la condition \ref{repr-stables}-(\ref{repr-stables-c}) est satisfaite pour un certain $N\in\NN$.

\subsectionline{Le $\FI$-modules $\M(\lambda)$ et $\V(\lambda)$.}
On\label{M-V-lambda} introduit un certain quotient $\V(\lambda):=\set\phi_m\colon V(\lambda)_m\to V(\lambda)_{\mm}/_m$
du $\FI$-module $\M_{|\lambda|+\lambda_1}(V(\lambda)_{|\lambda|+\lambda_1})$ (\ref{M-W})  qui rassemble les représentations $V(\lambda)_m$. La proposition \ref{prop-V-lambda} donne ses principales propriétés, notamment pour la description des catégories $\ModFIsc{\rkms\leq s} k$
et $\ModFIsc{\rkm\leq s} k$.

\begin{defi}
\'Etant donnée\label{defi-V-lambda}
une décomposition
$\lambda$, notons 
le $\FI$-module $\M_{|\lambda|}(V_\lambda)$ de \ref{M-W}\glossary{$\M(\lambda)\colon=\M_{|\lambda|}(V_\lambda)$:} plus simplement par
$$\M(\lambda):=\set\phi_m:M(\lambda)_m\to M(\lambda)_{\mm}/_m\,.
$$
D'après les règles de Pieri \ref{regles-de-branchement}-(\ref{regles-de-branchement-a}),
pour $m\geq |\lambda|$, la multiplicité $c(\mu)_m$ de $V(\mu)_m$ dans la décomposition en facteurs irréductibles
$$
M(\lambda)_m:=\ind_{\StS{a}{m-a}}^{\Sm}V(\lambda)_a\boxtimes k=
\bigoplusnl_{\mu}c(\mu)_mV(\mu)_m\,,
$$
est le nombre de fois que l'on peut obtenir le diagramme  $\mu[m]$ en rajoutant une boite sur $m{-}|\lambda|$ colonnes distinctes de $\lambda$. On a 
$$c(\mu)_m>0\Rightarrow
\big(|\lambda|\leq|\mu|\leq|\lambda|+\lambda_1\big)\,,$$ 

{\parskip0pt\noindent}\hmodeHabillage{\vrule\ $\halfdisplayskips\xymatrix@R=4mm@C=5mm{
V(\lambda)_m\strut\arinto[d]_{\iota}\\
M(\lambda)_m\ar[r]^(0.43){\phi_m}&
M(\lambda)_{\mm}\aronto[d]_(0.45){\nu}\\
&V(\lambda)_{\mm}
}$
}{2}{-10pt}et le $\FI$-module $\M(\lambda)^{(|\lambda|)}$ (le sous-quotient de $\M(\lambda)$ à poids unique $|\lambda|$ (\ref{prop-scindage-par-le-poids}-(\ref{prop-scindage-par-le-poids-a}))), est un quotient de $\M(\lambda)$. D'autre part, 
on a $c(\lambda)_m=1$ pour $m\geq|\lambda|+\lambda_1$, puisque
$|\mu|=|\lambda|$ équivaut à $\mu=\lambda$. Ceci implique que la composée $\nu\circ\phi_m\circ \iota$ ci-contre,
où~$\nu$ est la surjection canonique, est une \relax{injection} dont l'image engendre $V(\lambda)_{\mm}$. 
{\parskip0pt\noindent}On note alors\glossary{${V(\lambda)_{\mm}:=
{\M(\lambda)^{(\geq|\lambda|)}/\M(\lambda)^{(>|\lambda|)}
}}$:on a $\V(\lambda)=\set\phi_m\colon V(\lambda)_m\to V(\lambda)_{\mm}/_m$}
$$\displaystyle\preskip0pt
\mathrigid4mu
\relax{\V(\lambda)=\bigset \phi_m:V(\lambda)_m\to V(\lambda)_{\mm}/:=
\smashbot{\M(\lambda)^{(\geq|\lambda|)}\over\M(\lambda)^{(>|\lambda|)}
}\,\cdot}$$
\endHabillage
\end{defi}

\noendpoint\begin{prop}\label{prop-V-lambda} 
\def\varlistskips{\itemsep4pt}
\begin{enumerate}
\mynobreak\nobreak\item Pour\label{prop-V-lambda-a} toute décomposition $\lambda$, on a 
$$
\poids(\V(\lambda))=|\lambda|\text{\quad et\quad}
\rkm(\V(\lambda))=0\text { et }\rks(\V(\lambda))=|\lambda|+\lambda_1\,.$$

\comment
\begin{enumerate}
\item On a\label{prop-V-lambda-a-i}
$
\poids(\V(\lambda))=|\lambda|$,
 $\rkm(\V(\lambda))=\rks(\V(\lambda))=|\lambda|+\lambda_1$.

\item\leavevmode\label{prop-V-lambda-a-ii}Le $\FI$-module $\V(\lambda)$ est de type fini et, pour tout $a\geq |\lambda|+\lambda_1$, le morphisme naturel 
$\preskip0pt\postskip0pt
\M_{a}(V(\lambda)_{a})\onto \V(\lambda)_{\geq a}$
est surjectif.
\end{enumerate}
\endcomment
\item\leavevmode\label{prop-V-lambda-b}Soit $\W=\set\phi_{m}: W_{m}\to W_{\mm}/_{m}$ un $\FI$-module à poids unique $t$. Pour $n\in\NN$, si $W_n=\bigoplusnl_{|\mu|=t} c(\mu)_n V(\mu)_n$
est la décomposition en composantes irréductibles, on note
$$\UV(\W,n):=\bigoplusnl_{|\mu|=t} c(\mu)_n \V(\mu)_{\geq n}\,.$$

\begin{enumerate}
\item Pour\label{prop-V-lambda-b-i} tout $\lambda$ et tout $n\geq|\lambda|+\lambda_1$, l'application canonique
$$\Hom_{\FI}(\V(\lambda)_{\geq n},\W)\to\Hom_{\S_n}(V(\lambda)_{n},W_{n})$$
est bijective. 
Dans la suite, on notera $\varUpsilon(n):\UV(\W,n)\to\W$
le morphisme de $\FI$-modules correspondant à l'identité en degré $n$.

\item  $\rkms(\W)\leq s$\label{prop-V-lambda-b-ii}, si et seulement si,  $\varUpsilon(s):\UV(\W,s)\to\W_{\geq s}$ est bijectif.

\penalty-500\item 
$\mathrigid1mu \rkm(\W)\leq s$\label{prop-V-lambda-b-iii}, si et seulement si, 
$\mathrigid1mu \varUpsilon(n):\UV(\W,n)\to\W$
est injectif 
$\mathrigid1mu \forall n\geq s$.

\noindent En particulier, si $\rkm(\W)\leq s$ et si
$W_m\sim V(\lambda)_{m}$ pour tout $m\geq s$, on a $\W_{\geq s}\simeq\V(\lambda)_{\geq s}$.

\end{enumerate}
\item Si\label{prop-V-lambda-c} $\W$ est tel que $s:=\rkms(\W)<+\infty$, le $\FI$-module $\W_{\geq s}$ est extension de $\FI$-modules de la forme $\V(\mu)_{\geq s}$ où $|\mu|+\mu_1\leq s$.

\item Un\label{prop-V-lambda-d} $\FB$-module $\W=\set W_m/_{m}$ est stable pour $m\geq N$, si et seulement si, il existe un $\FI$-module
$\Z=\set\phi_{m}: Z_{m}\to Z_{\mm}/_{m}$ monotone et stable pour $m\geq N$,  tel que $W_m=Z_m$ pour tout $m\geq N$.
\end{enumerate}
\end{prop}
\demo
 (\ref{prop-V-lambda-a}) est implicite dans la définition \ref{defi-V-lambda}

\parskip2pt\mou

\smallskip
(\ref{prop-V-lambda-b})
Pour $m\in\NN$, si nous considérons $V(\mu)_{\mm}$ comme $\Sm $-module, la règle de Pieri \ref{regles-de-branchement}-(\ref{regles-de-branchement-a}), nous dit qu'il contient une et une seule composante irréductible de poids $|\mu|$, à savoir $V(\mu)_m$ et que sa multiplicité est égale à $1$. Il s'ensuit que la restriction $\phi(\mu)_m$ du morphisme de transition de $\W$,
$$\phi_m:\Big(W_m=\bigoplusnl_{|\mu|=t}c(\mu)_mV(\mu)_m\Big)\to \Big(W_\mm=
\bigoplusnl_{|\mu|=t}c(\mu)_{\mm}V(\mu)_{\mm}\Big)\,,$$
à $c(\mu)_mV(\mu)_m$ est à valeurs dans $c(\mu)_{\mm}V(\mu)_{\mm}$. Les sous-familles 
$$\W(\mu):=\set\phi(\mu)_m:c(\mu)_mV(\mu)_m\to 
c(\mu)_{\mm}V(\mu)_{\mm}/_m$$ sont donc des sous-$\FI$-modules de $\W$ et nous avons la décomposition
$$\W=\bigoplusnl_{|\mu|=t}\W(\mu)\,.$$

Chaque morphisme $\phi(\mu)_m$ se factorise suivant le diagramme
$$\xymatrix@C=1cm{
c(\mu)_mV(\mu)_m\ar[r]^{\phi(\mu)_{m}}\arinto[d]_{\iota(\mu)_{m}}&
c(\mu)_{\mm}V(\mu)_{\mm}\\
c(\mu)_m\ind_{\Sm}^{\Smm} V(\mu)_{m}\aronto[r]_(.55){\pi(\mu)_m}\ar[ru]|{\ \vrule depth3pt width0pt\Phi(\mu)_m\ }&c(\mu)_mV(\mu)_{\mm}\ar@{-->}[u]_{\Psi(\W)(\mu)_{\mm}}
}
$$
où $\phi(\mu)_{m}=\Phi(\W)(\mu)_{m}\circ\iota(\mu)_m$ est la factorisation standard (\ref{familles-denombrables}-($\diamond\diamond$)), où
$\pi(\mu)_m:\ind_{\Sm}^{\Smm}V(\mu)_{m}\onto V(\mu)_{\mm}$ est la surjection canonique sur l'unique composante irréductible de
$\ind_{\Sm}^{\Smm}V(\mu)_{m}$ de poids $|\mu|$, et où
$\Psi(\mu)_{m}$ est le morphisme induit par $\Phi(\mu)_m$.

En additionnant, on obtient la factorisation de $\phi_{m}=\sum_{\mu}\phi(\mu)_m$:
$$\halfdisplayskips
\xymatrix@C=1.2cm{
W_m\arinto[r]^(0.4){\iota_m}&\ind_{\Sm}^{\Smm}W_m\aronto[r]^(0.4){\pi_m}&
\big[\ind_{\Sm}^{\Smm}W_m\big]_{t}\ar[r]^(0.6){\Psi(\W)_{\mm}}&
W_{\mm}}\eqno(\ast)$$
où on a noté  $[M]_{t}$ la somme des composantes irréductibles de poids $t$ d'un $\Smm$-module $M$.

La même démarche pour le $\FI$-module $\U:=\UV(\W,n)$ conduit au même type de factorisation $(\ast)$ (à ceci près que $\Psi(\U)_{\mm}$ un isomorphisme). La naturalité de cette factorisation permet de construire, à partir de la donnée d'un morphisme $f_{n}:U_n\to W_n$, une famille (unique) de morphismes $\set f_m:U_m\to W_m/_{m\geq n}$ qui constitue un morphisme de $\FI$-modules $f:\U\to\W$.

En effet, fixons $f_n:U_n\to W_n$ et supposons avoir défini $f_{n+1},\ldots,f_{z}$ de manière compatible aux morphismes de transition de $\U$ et $\W$. On a alors le  morphisme de factorisations $(\ast)$
$$\halfdisplayskips\def\m{z}\def\mm{z{+}1}\def\Sm{\S_z}
\xymatrixc{@C=1.5cm}{
U_{\m}\xylbl[dr]{(I)}\ar[d]|{\strut f_\m}\arinto[r]^(0.4){\iota_{\m}}&\ind_{\Sm}^{\Smm}U_{\m}\ar@{..>}[d]^{\ind_{\Sm}^{\Smm} f_z}\aronto[r]^(0.4){\pi_{\m}}\xylbl[dr]{\hdecale{-5mm}{(II)}}&
\big[\ind_{\Sm}^{\Smm}U_{\m}\big]{\vrule depth3.5pt width0pt}_{t}\ar@{..>}[d]^{[\ind_{\Sm}^{\Smm} f_z]{\vrule depth2.5pt width0pt}_{t}}\ar[r]_(0.6){\sim}^(0.6){\Psi(\U)_{\mm}}
&
U_{\mm}\ar@{-->}[d]|{\strut f_{\mm}}\\
W_{\m}\arinto[r]^(0.4){\iota_{\m}}&\ind_{\Sm}^{\Smm}W_{\m}\aronto[r]^(0.4){\pi_{\m}}&
\big[\ind_{\Sm}^{\Smm}W_{\m}\big]{\vrule depth3.5pt width0pt}_{t}\ar[r]^(0.6){\Psi(\W)_{\mm}}&
W_{\mm}
}\eqno\,\mllap{(\ast\ast)}$$
où les flèches en pointillé sont induites par $f_z$ et sont uniques à rendre commutatifs les sous-diagrammes (I) et (II). Enfin, comme $\Psi(\U)_{z+1}$ est bijectif, l'existence et unicité de $f_{z+1}$ sont claires, ce qui termine la preuve de (\ref{prop-V-lambda-b-i}).

On remarque alors que lorsque $f_z$ est un isomorphisme, $f_{z+1}$ est un isomorphisme, si et seulement si, $\Psi(\W)_{z+1}$ l'est. Or cette dernière condition est une condition de monotonie et stabilité sur $\W$, de sorte que si $s:=\rkms(\W)$, le morphisme $\varUpsilon(s):\UV(\W,s)\to\W_{\geq s}$ est un isomorphisme, d'où (\ref{prop-V-lambda-b-ii}).

De manière analogue, lorsque $f_z$ est une injection, on note $W'_{z}:=\im(f_z)$, et l'on fixe une décomposition $W_z=W'_{z}\oplus N_z$. La dernière ligne de $(\ast\ast)$ se décompose alors en somme directe de deux lignes, ce qui permet de voir que $f_{z+1}$ est injective, si et seulement si, la restriction de $\Psi(\W)_{z+1}$ à
\smash{$\big[\ind_{\Sm}^{\Smm}W'_{z}\big]{\vrule depth3.5pt width0pt}_{t}$} l'est. Or, ceci est très précisément la condition de monotonie en $z$ pour $\W$, et l'assertion (\ref{prop-V-lambda-b-iii}) s'ensuit.

\smallskip
(\ref{prop-V-lambda-c}) Résulte de (\ref{prop-V-lambda-b}) et de \ref{prop-scindage-par-le-poids}-(\ref{prop-scindage-par-le-poids-b},\ref{prop-scindage-par-le-poids-c}).

(\ref{prop-V-lambda-d}) Il suffit de prendre $\Z:=\UV(\W,N)$.
\enddemo

\subsection{Monotonie, stabilité et co-invariants}\label{deg-stab}
\subsubsectionline{Co-invariants.}Dans\label{co-invariants}
\cite{cef}~(\Spar3.1), on introduit, pour
 $0\leq t\leq m$, le foncteur des \expression{$\S_{\mt}$-co-invariants}  (pour l'inclusion $\St \times\S_{\mt}\dans\Sm$)\glossary{${\varPhi_t:=(\_)_{\S_{\mt} }:\Mod(\Sm)\fonct\Mod(\St )}$:foncteur de $\S_{\mt}$-co-invariants}
$$\mathalign{(\_)_{\S_{\mt} }:&\Mod(\Sm)&\fonct&\Mod(\St )\hfill\\
&W&\fonct& W_{\S_{\mt}}:=\relax{k\otimes_{k[\S_{\mt}]}W\,.}}$$

Pour tout $\FI$-module $\W=\set\phi_m: W_m\to W_{\mm}/$ et tout $m\geq t$, notons $\nu_{m}:W_m\onto (W_{m})_{\S_{\mt}}$ la surjection $w\mapsto 1\otimes w$. La composée
$$\preskip0pt
\xymatrix@C=10mm@R=7mm{
\mllap{\hbox to1cm{\rightarrowfill}}W_m\ar[dr]+<-22pt,2pt>|{\vrule depth4pt width0pt\nu_\mm\circ\,\phi_m}
\ar[r]^(0.45){\phi_m}&W_{\mm}\aronto[d]^{\nu_{\mm}}\mrlap{\hbox to1cm{\rightarrowfill}}\\
\relax\phantom{\hskip7mm(W_m)_{\S_{\mt}}}
&(W_{\mm})\mrlap{_{\S_{\mmt}}}
}
$$
se factorise à travers $\nu_m:W_m\onto(W_{m})_{\S_{m-t}}$ en un unique morphisme de $\St $-modules
$(\phi_{m}):
(W_m)_{\S_{\mt}}\to (W_{\mm})_{\S_{\mm{-}t}}$ rendant les diagrammes suivants commutatifs
$$\preskip-1ex
\xymatrix@C=10mm@R=7mm{
\mllap{\hbox to1cm{\rightarrowfill}}W_m\aronto[d]_{\nu_m}\ar[r]^(0.45){\phi_m}&W_{\mm}\aronto[d]^{\nu_{\mm}}\mrlap{\hbox to1cm{\rightarrowfill}}\\
\hskip7mm(W_m)_{\S_{\mt}}\ar[r]^(0.57){(\phi_m)}&(W_{\mm})\mrlap{_{\S_{\mmt}}}
}
$$

\subsubsectionline{Le foncteur $\varPhi_{t}$.}
Pour\label{foncteur-Phi-t} $t\in\NN$, on note dans \cite{cef} (déf. 3.1.2)\glossary{$\varPhi_{t}(\W)$:$k[\St ][T]$-module gradué associé au $\FI$-module $\W$} 
$$\varPhi_{t}:\ModFI k\fonct\Mod^{\NN}(k[\St ][T])\,,$$
le foncteur qui fait correspondre à $\W=\set \phi_m:W_m\to W_{\mm}/_{m}$ le $\St $-module positivement gradué $\varPhi_{t}(\W)^{\ast}$ défini par
$$\varPhi_{t}(\W)^{m-t}:=(W_{m})_{\S_{\mt}}\,,\quad\forall m\geq t\,.$$
muni de l'action de degré $+1$ de $T$ qui vaut $(\phi_{m})$ sur $(W_{m})_{\S_{\mt}}$.

\subsubsubsectionline{A propos des notations.}Le foncteur de co-invariants $(\_)_{\S_{m-t}}:\Mod(\Sm)\fonct\Mod(\St)$ sera parfois aussi noté $\varPhi_{t}(\_):\Mod(\Sm)\fonct\Mod(\St)$.
En particulier, on pourra écrire 
$\varPhi_{t}(\W)^{m-t}=\varPhi_{t}(W_m)$ et $(\phi_m)=\varPhi_{t}(\phi_m)$.

\begin{lemm}
Soit\label{lemme-FB-stable} 
$\V(\lambda)=\set\phi_m:V(\lambda)_{m}\to V(\lambda)_{\mm}/$ le $\FI$-module de \ref{M-V-lambda}.
Pour tout $m\geq t$, le morphisme de $\St $-modules
$$\varPhi_t(\phi_m):\varPhi_{t} (V(\lambda)_{m})\to \varPhi_{t} (V(\lambda)_{\mm})$$
est injectif, et il est bijectif si $t<|\lambda|$ ou si $m\geq t+\lambda_1$.

\smallskip
En particulier, si $\W=\bigoplus_{\lambda} c(\lambda)\V(\lambda)
=\set\phi_m:W_{m}\to W_{\mm}/$. Le morphisme $\varPhi_{t}(\phi_m)$ est injectif pour tout $m\geq t$.
\end{lemm}

\demo 
Pour $m\geq|\lambda|+\lambda_1$ numérotons le diagramme de Young de $\lambda[m]$ de haut en bas et de gauche à droite comme dans le tableau:
$$\ygwd=0.9em \yght\ygwd
\tau_{\lambda}(m):=\tableauc{}{
1&5&8&10&\noborder\cdots&\cdots&m\\
2&6&9&11\\
3&7\\
4}$$
Notons (suivant \cite{fulhar} \Spar4.1) $L_{\lambda}(m)$,
 $P_{\lambda}$ et $Q_{\lambda}$ les sous-groupes de $\S _{m}$ qui laissent respectivement stables, la première ligne, les autres lignes et les colonnes du diagramme sous-jacent
à $\tau_{\lambda}(m)$. Notons ensuite, dans $k[\S_{m}]$,
$$
\Yl_{\lambda}(m):=\sumnl_{\alpha\in L_{\lambda}(m)}\,,\quad
\Ya_{\lambda}:=\sumnl_{\alpha\in P_{\lambda}}\,,\quad
\Yb_{\lambda}:=\sumnl_{\alpha\in Q_{\lambda}}\sgn(\alpha)\alpha\,.
$$
Le \expression{symétriseur de Young} associé à $\tau_{\lambda}(m)$ est l'\'{e}l\'{e}ment de $k[\S _m]$\glossary{${\Yc_{\lambda}(m):=\Yl_{\lambda}(m)\cdot\Ya_{\lambda}\cdot \Yb_{\lambda}}$:symétriseur de Young associé à $\lambda_{[m]}$}
$$\Yc_{\lambda}(m):=\Yl_{\lambda}(m)\cdot \Ya_{\lambda}\cdot \Yb_{\lambda}\,,$$
et le sous-$\S _m$-module à gauche 
$k[\S _m]\cdot \Yc_{\lambda}(m)\dans k[\S _m]$ est 
isomorphe à la représentation irréductible $V(\lambda)_m$. 
On a aussi (\loccit\ ex. 4.4)
$$k[\S _m]\cdot \Yc_{\lambda}(m)=k[\S _m]\cdot \Yb_{\lambda}
\cdot \Ya_{\lambda}\cdot \Yl_{\lambda}(m)
\,.$$

\displayskips8/10

Soient maintenant $|\lambda|+\lambda_1\leq m\leq n$. D'après le choix du schéma de numérotations des tableaux, il est clair que par le plongement $\S _m\dans\S _n$ on a
$$
\Yl_{\lambda}(m)\cdot \Yl_{\lambda}(n)=(m-|\lambda|)!\;\Yl_{\lambda}(n)\,.$$
La multiplication à droite par $\Yl_{\lambda}(n)$ définit par conséquent un morphisme de $\S _m$-modules à gauche injectif:
$$\phi_{n,m}:k[\S _m]\cdot \Yb_{\lambda}\cdot \Ya_{\lambda}
\cdot \Yl_{\lambda}(m)
\hfhook{\vrule depth2pt width0pt\smash{(\_)\cdot \Yl_{\lambda}(n)}}{}{2cm}
k[\S _n]\cdot \Yb_{\lambda}\cdot \Ya_{\lambda}
\cdot \Yl_{\lambda}(n)\,.\eqno(\diamond)$$
On remarque ensuite que toute permutation $\alpha\in\II_{m}\times\S _{n-m}$ fixe l'image de $\phi_{n,m}$ puisque, d'une part, $\alpha$ commute à $k[\S _{m}]$ donc à $k[\S _m]\cdot \Yb_{\lambda}\cdot \Ya_{\lambda}$ et, d'autre part, on a $\alpha\cdot \Yl_{\lambda}(n)=\Yl_{\lambda}(n)$ parce que $\II_{m}\times\S _{n-m}\dans\Yl_{\lambda}(n)$.

Si nous notons maintenant $\phi_m=\phi_{m+1,m}$, 
l'injection $\phi_{n,m}$ est multiple de la composée $\phi_{n-1}\circ\cdots\circ\phi_{m}$ et la famille $\W(\lambda)$
dont les termes $W_m$ sont nuls pour $m<|\lambda|+\lambda_1$ et qui, pour $m\geq|\lambda|+\lambda_1$, coïncide avec
$$\bigset\phi_m:k[\S _m]\cdot \Yc_{\lambda}(m)\to k[\S _{m+1}]\cdot \Yc_{\lambda}(m+1)/_{m}$$
est un $\FI$-module canoniquement isomorphe à $\V(\lambda)$ d'après \ref{prop-V-lambda}-(\ref{prop-V-lambda-b-iii}).

\smallskip

Ce qui précède justifie l'égalité
$$
\varPhi_{t}(V(\lambda)_{m})=k\lquo{\S_{m-t}}{\Sm }\cdot\Yl_{\lambda}(m)\cdot\Ya_{\lambda}\cdot\Yb_{\lambda}\,,\eqno\forall  m\geq t\,,$$
et $\varPhi_t(\phi_m)$ sera injective si et seulement si le morphisme induit par la multiplication à droite par $\Yl_{\lambda}(\mm)$:
$$k\lquo{\S_{m-t}}{\Sm }\cdot\Yl_{\lambda}(m)
\to k\lquo{\S_{\mm-t}}{\Smm }\cdot\Yl_{\lambda}(\mm)
$$
est injectif. Or, ceci revient à montrer que 
l'application naturelle
$$
\biquo{\S_{m-t}}{\Sm }{L_{\lambda}(m)}
\too
\biquo{\S_{\mm-t}}{\Smm }{L_{\lambda}(\mm)}
$$
est injective, ce qui résulte d'une analyse élémentaire.

\smallskip
Maintenant, du fait de l'injectivité de $\varPhi_t(\phi_m)$, sa bijectivité résulte du décompte des composantes irréductibles de $\varPhi_{t}(V(\lambda)_\mm)$ en tant que $\St$-module. La règle de branchement \ref{regles-de-branchement}-(\ref{regles-de-branchement-b}) nous dit que ce nombre est le nombre des manières d'enlever $\mm{-}t$ boites de colonnes différentes de $\lambda[\mm]$. Il s'ensuit que si $\mm{-}t>\lambda_1$, on est obligé d'enlever la dernière boite de la première ligne de $\lambda[\mm]$ et l'on retrouve alors le nombre des composantes irréductibles de $\varPhi_{t}(\V(\lambda)_{m})$ par la même règle de branchement.

\smallskip Le assertion concernant $\W=\bigoplus_{\lambda} c(\lambda)\V(\lambda)$ est conséquence immédiate du cas où la somme et les multiplicités sont finies, ce qui est clair d'après la première assertion.
\enddemo
\begin{defis}Le\label{degre-stabilite} \expression{degré de stabilité (resp. d'injectivité) en $t\in\NN$} d'un $\FI$-module $\W$, noté $\degstab_{t}(\W)$\glossary{$\degstab_t(\W)$:degré de stabilité en $t\in\NN$ d'un $\FI$-module $\W$}
(resp. $\deginj_{t}(\W)$\glossary{$\deginj_t(\W)$:degré d'injectivité en $t\in\NN$ d'un $\FI$-module $\W$}), est le plus petit $d\in\NN\cup\set+\infty/$ tel que l'application $T:\varPhi_{t}(\W)^{n}\to\varPhi_{t}(\W)^{n+1}$ est bijective (resp. injective) pour tout $n\geq d$. 
\end{defis}

\begin{prop}
On fixe $t\in\NN$.\label{lemme-co-invariants}
\begin{enumerate}\itemsep4pt\mou
\item \label{lemme-co-invariants-a}Le foncteur $\varPhi_{t}:\ModFI k\fonct\Mod^{\NN}(k[\St ][T])$ est covariant et exact.
\item Soit\label{lemme-co-invariants-b} $\lambda$ une partition. Pour tous $m\geq t$, notons $\III_{t}(V(\lambda)_{m})$ la famille des partitions $\nu\vdash t$ telles que $V_\nu$ est un facteur irréductible de $\Phi_{t}(V(\lambda)_{m})$, répétées autant de fois que leurs multiplicités. 
\comment 
Alors $\nu\in\III_{t}(V(\lambda)_{m})$ si et seulement si, la décomposition $\nu$ est obtenue de $\lambda[m]$ en en enlevant une boite sur $\mt$ colonnes distinctes. En particulier, 
\endcomment
Alors, on a des inclusions
$$\III_{t}(V(\lambda)_{m})\dans\III_{t}(V(\lambda)_{\mm})\,,\eqno\forall m\geq t\,,$$
et ces inclusions sont des égalités si et seulement si $m\geq t+\lambda_1$.
\item Soit\label{lemme-co-invariants-c} $\lambda$ une partition et soit $\V(\lambda)$ le $\FI$-module de \ref{prop-V-lambda}. Alors, 
\def\varlistskips{\labelsep3pt\labelwidth1cm\advance\leftmargin2pt \itemsep4pt}
\begin{enumerate}\preskip0.75ex
\item\leavevmode Pour\label{lemme-co-invariants-c0} tout $t\in\NN$,
$$\preskip-1.2em\deginj_t(\V(\lambda)_{\geq s})=0\,,\eqno\forall s\in\NN\,.$$
\item\leavevmode\label{lemme-co-invariants-ci}$t<|\lambda|\Leftrightarrow\varPhi_{t}(\V(\lambda))=0$. Dans ces cas,
$$\degstab_t(\V(\lambda)_{\geq s})=0\,,\eqno\forall s\in\NN\,.$$
\item\leavevmode\label{lemme-co-invariants-cii}$t=|\lambda|\Leftrightarrow\varPhi_{t}(\V(\lambda))^{\mt}=V_{\lambda}$
et $T=\id$. Dans ces cas, 
$$\degstab_{|\lambda|}(\V(\lambda)_{\geq s})=\sup\set\lambda_1,s-|\lambda|/\,,\eqno\forall s\in\NN\,.$$
\item Si\label{lemme-co-invariants-ciii}
$t\geq|\lambda|$ et $m\geq t+\lambda_1\,,$
on a $\varPhi_t(\V(\lambda))^{\mt}\simeq
\relax{\ind_{\vrule height0pt width0pt
\hskip-0ex\StS{|\lambda|}{t-|\lambda|}}^{\St }}
V_\lambda\boxtimes k\,,
$
et 
$T:\varPhi_t(\V(\lambda))^{\mt}\to\varPhi_t(\V(\lambda))^{\mmt}$ est bijective. Dans ces cas, 
$$\degstab_{t}(\V(\lambda)_{\geq s})=\sup\set\lambda_1,s-t/\,,\eqno\forall s\in\NN\,.$$
\end{enumerate}
\end{enumerate}
\end{prop}
\demo (\ref{lemme-co-invariants-a}) est immédiat. 

(\ref{lemme-co-invariants-b}) et (\ref{lemme-co-invariants-c0}) par l'injectivité de $(\phi_m)$ du lemme \ref{lemme-FB-stable}. 
\comment Résulte aussi d'appliquer la règle \ref{regles-de-branchement}-(\ref{regles-de-branchement-b}). En effet, pour tout $m\geq|\lambda|+\lambda_1$, si la décomposition $\nu\vdash t$ provient du tableau $\lambda[m]$ en lui enlevant les dernières boites des colonnes $1\leq c_1<\cdots< c_{\mt}\leq
m{-}|\lambda|$, elle provient aussi du tableau $\lambda[\mm]$ en lui enlevant, en plus, la boite isolée sur la dernière colonne $c_{\mm-|\lambda|}$. Cette remarque donne une inclusion canonique $\III_{t}(\V(\lambda)_{m})\dans\III_{t}(\V(\lambda)_{\mm})$, et cette inclusion est une égalité, si et seulement si, pour enlever $\mmt$ boites sur des colonnes différentes de $\lambda[\mm]$, on est obligé d'enlever une boite sur la première ligne de ce tableau, or, ceci est clairement équivalent à demander que $\mmt>\lambda_1$.
\endcomment

(\ref{lemme-co-invariants-ci})-(\ref{lemme-co-invariants-ciii}) c'est le lemme 3.2.7 de \cite{cef}, qui résulte également du lemme \ref{lemme-FB-stable} et de la règle de branchement \ref{regles-de-branchement}-(\ref{regles-de-branchement-b}). Nous omettons les détails.
\enddemo

\comment
\subsubsectionline{Co-invariants pour $\SL$.}Dans ce qui précède le choix du sous-groupe $\1_{t}\times\S_{\mt}$ n'est pas essentiel pour caractériser le degré de stabilité d'un $\FI$-module et une approche plus intrinsèque est parfois plus commode. On en aura notamment besoin dans la démonstration du théorème \ref{theo-Ind-lambda}.

Pour\label{critere-stabilite} 
$L\dans\iii[1,m]$, on identifie le groupe $\SL $\glossary{$\SL $:on a $L\dans\iii[1,m]$ et $\SL :=\Fix_{\Sm}(\iii[1,m]\mmoins L)$} des bijections de $L$, au sous-groupe de $\Sm $ qui fixe  $i\not\in L$. 
Le centralisateur de $\SL $ dans $\Sm$ est le sous-groupe $\S_{\iii[1,m]\mmoins L}$ qui fixe $i\in L$. Tout comme dans \ref{co-invariants}, le foncteur des $\SL$-co-invariants est le foncteur
$$\mathalign{(\_)_{\SL  }:&\Mod(\Sm)&\fonct&\Mod(\S_{\smashtop{\iii[1,m]\mmoins L}})\hfill\\
&W&\fonct& W_{\SL}:=W\otimes_{k[\SL ]}k}$$
Le morphisme canonique $\nu:W\to W_{\SL}$, $w\mapsto w\otimes 1$ est une transformation naturelle de foncteurs $\nu:\id\to(\_)_{\SL}$.

\smallskip

On note $L':=L\sqcup\set m+1/\dans\iii[1,m+1]$\glossary{$L'$:pour $L\dans\iii[1,m]$, on pose $L':=L\coprod\set m+1/$}. Pour tout $\FI$-module $\W$, le morphisme de transition
$\phi_m: W_m\to W_{\mm}$ passe aux co-invariants où il définit le morphisme de $\S_{\iii[1,m]\mmoins L}$-modules 
$$\def\mt{L}\def\mmt{L'}
(\phi_{m}):
(W_m)_{\S_{\mt}}\to (W_{\mm})_{\S_{\mm{-}t}}$$ qui rend commutatif le diagramme
$$\preskip2pt
\def\mt{L}\def\mmt{L'}
\xymatrix@C=15mm@R=6mm{
\mllap{\hbox to1cm{\rightarrowfill}}W_m\aronto[d]_{\nu_m}\ar[r]^(0.45){\phi_m}&W_{\mm}\aronto[d]^{\nu_{\mm}}\mrlap{\hbox to1cm{\rightarrowfill}}\\
(W_m)_{\S_{\mt}}\ar[r]^(0.45){(\phi_m)}&(W_{\mm})_{\S_{\mmt}}
}
$$

\begin{prop}\color{red}Pour\label{deg-stab/L} tout $\FI$-module $\W=\set\phi_m:W_{m}\to W_{\mm}/$, notons $L(\W)$ l'ensemble des parties finies $L\dans\NN$ telles que $(\phi_{m})_{\SL}$ est un isomorphisme pour tout $m$ tel que
$L\dans\iii[1,m]$. Alors, $\degstab(\W)$ est la borne inférieure des cardinaux des $L\in L(\W)$.
\end{prop}
\def\Demonstration{Indication}\demo
\`A l'aide d'un automorphisme intérieur convenable de $\Sm$, on se ramène au cas où $\SL=\1_{t}\times\S_{m-t}$ et où
$(\phi_m)_{\SL}=(\phi_{m})_{\S_{\mt}}$.
\enddemo
\endcomment

\subsubsectionline{Suite du commentaire \ref{comment-poids-degre}.}
L'intérêt\label{comment-poids-degre-2} du poids et du degré de stabilité d'un $\FI$-module $\W$ apparaît dans la proposition 3.3.3 de \cite{cef} qui établit la majoration $\rks(\W)\leq\poids(\W)+\degstab(\W)$ (\footnote{Le\label{note-deg-stab} \expression{degré de stabilité} de $\W$ y est défini par\glossary{${\degstab(\W):=\sup\nolimits_{t\in\NN}\set\degstab_t(\W)/}$:degré de stabilité d'un $\FI$-module $\W$}
$\degstab(\W):=\sup_{t\in\NN}\set\degstab_t(\W)/
$. Il ne sera pas utilisé dans ce travail.}). 
Un outil important dans notre travail est le foncteur d'induction 
$\Ind_{\lambda}:\ModFI k\fonct\ModFI k$ qui sera introduit dans la section \ref{foncteurs-decalages-inductifs}. Ce foncteur est un `recollement' des foncteurs $\Ind_{\smashtop{G_{\lambda\bxx m}}}^{\Sm}$ de \ref{theo-caracteres}-(\ref{theo-caracteres-b}) et lorsqu'il est appliqué à un $\FI$-module $\W$, il est crucial de comprendre le rapport entre les rangs de $\W$ et ceux de $\Ind_\lambda(\W)$.
En ce sens, la proposition 3.3.3 \loccit\ était à priori intéressante  dans la mesure où nous avons un meilleur contrôle de la perturbation de $\poids(\W)$ et de $\degstab(\W)$ lors de ces inductions (\cf\ref{poids-ind}). Cela n'a pourtant pas suffi et nous avons été emmenés à caractériser autrement les rangs de $\W$. La suite décrit cette nouvelle approche. Elle conduira au théorème \ref{theo-Ind-lambda} qui donne les majorations
$\rkms(\Ind_\lambda(\W))\leq\rkms(\W)+2|\ulambda|$
et 
$\rkm(\Ind_\lambda(\W))\leq\rkm(\W)+|\ulambda|$.

\subsubsection{Rangs étendus}Les\label{defi-rang-etendu} \expression{rangs étendus} de $\W=\set \phi_m:W_m\to W_{\mm}/_{m}$ sont les nombres\glossary{${\rkme(\W)\,,\ \rkmse(\W)}$:rangs étendus du $\FI$-module $\W$}
$$\def\text#1{\hbox to2.4cm{\hss#1\hss}}
\left\{\vrule height18pt width0pt\right.\!\!\mathalign{
\rkmse(\W)&:=&
\inf \bigset 
s\in\NN \mid \varPhi_{t}(\phi_m) \text{ est bijective } \forall m\geq s\,,\forall t\leq\poids(\W)\leq m/
\hfill\\\noalign{\kern3pt}
\rkme(\W)&:=&
\inf\bigset 
s\in\NN \mid \varPhi_{t}(\phi_m) \text{ est injective } \forall m\geq s\,,\forall t\leq m/
\hfill}$$
$$\postskip-2em\includegraphics{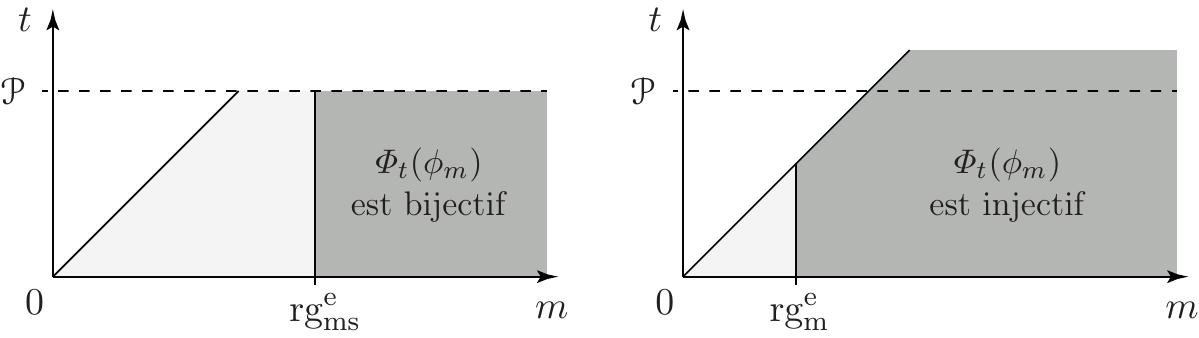}$$

\begin{prop}Pour tout\label{rkms=rke} $\FI$-module $\W$,
on a
$$
\displayboxit{\rkm(\W)=\rkme(\W)
\text{\quad et \quad}
\rkms(\W)=\rkmse(\W)}
$$
\end{prop}

\noindent La démonstration résultera du lemme suivant.

\noendpoint\begin{lemm*}\label{prop-rang-etendu}
\def\varlistskips{\topsep0pt\mou\itemsep4pt\mou}
\begin{enumerate}
\mynobreak
\item Soit\label{prop-rang-etendu-a} $\rke\in\set\rkme,\rkmse/$. 
\begin{enumerate}\halfdisplayskips
\item Si\label{prop-rang-etendu-ai}
$
\xymatrix@C=3mm{
0\to\X\ar[r]&\W\ar[r]&\Y\to0
}
$
est une suite exacte de $\FI$-modules,
on a
$$
\rke(\W)\leq\sup\set\rke(\X),\rke(\Y)/\,.
$$
\item Si\label{prop-rang-etendu-aii} $\set \W^{(t)}/$ est la famille des sous-quotients à poids unique de $\W$, on a 
$$
\rke(\W)=\smashtop{\sup\set \rke(\W^{(t)})/}\,.
\postskip7pt
$$
\end{enumerate}
\item Si\label{prop-rang-etendu-b} $\W$ est à poids unique {\rm (\ref{prop-scindage-par-le-poids}-(\ref{prop-scindage-par-le-poids-a-iii}))}, on a
$$\halfdisplayskips
\rkm(\W)=\rkme(\W)
\text{\quad et\quad }
\rkms(\W)=\rkmse(\W)
\,.$$
\end{enumerate}
\end{lemm*}
\def\Demonstration{Preuve du lemme}\demo 
\displayskips7/10
\parskip1ex \mou
(\ref{prop-rang-etendu-ai}) Pour chaque $t\in\NN$, la suite de $k[\S_t][T]$-modules gradués 
\relax{$
\xymatrix@C=3mm{
0\to\varPhi_t(\X)\ar[r]&\varPhi_t(\W)\ar[r]&\varPhi_t(\Y)\to0
}
$}
est exacte (\ref{lemme-co-invariants}-(\ref{lemme-co-invariants-a})), auquel cas
$\degstab_t(\W)\leq\sup\set\degstab_t(\X),\degstab_t(\Y)/\,,$
et de même pour $\deginj_{t}$.

(\ref{prop-rang-etendu-aii}) On suppose $\poids(\W)<+\infty$. Montrons que pour le plus petit entier $a$ tel que $\W^{(a)}\not=0$, on a
$$
\rke(\W)=\sup\set \rke(\W^{(>a)}),\rke(\W^{(a)})/\,.\eqnnb\label{prop-rang-etendu-ast}$$
En effet, (\ref{prop-rang-etendu-ai}) donne déjà l'inégalité `$\leq$'. Ensuite, $\rke(\W)\geq\rke(W^{(a)})$ vient de ce que si $t\leq a$, on a $\varPhi_{t}(W^{>a})=0$ et alors $\varPhi_{t}(W)\sim\varPhi_{t}(W^{a})$. Il nous reste
$$
\rke(\W)\geq\rke(\W^{(>a)})\,.\eqnnb\label{prop-rang-etendu-star}$$

Le cas où $\rke=\rkme$ est clair puisque $\varPhi_{t}(\W^{(>a)})$ est un sous-$k[\St][T]$-module de \smashtop{$\varPhi_{t}(\W^{(>a)})$}. 
Le cas $\rke=\rkmse$, résulte de se placer dans la situation où $a<t\leq\poids\leq \rkmse(\W)\leq m$ et de considérer
le diagramme de colonnes exactes
$$\halfdisplayskips\preskip-1ex
\xymatrix@R=4mm@C=1.5cm{
\strut \varPhi_{t}(W^{(>a)}_{m})\arinto[d]\ar[r]^{(\phi_{m}^{>a})}&\strut \varPhi_{t}(W^{(>a)}_{\mm})\arinto[d]\\
\varPhi_{t}(W_{m})\aronto[d]\ar@{>->>}[r]^{(\phi_m)}_{\simeq}&\varPhi_{t}(W_{\mm})\aronto[d]\\
\varPhi_{t}(W^{(a)}_{m})\ar[r]^(.48){(\phi_m^{a})}&\varPhi_{t}(W^{(a)}_{\mm})
}$$
où $(\phi_m^{a})$ est bijective, car injective d'après \ref{lemme-co-invariants}-(\ref{lemme-co-invariants-c0}) et surjective d'après les autres flèches. On en déduit la bijectivité de $(\phi_{m}^{>a})$ d'où \ref{prop-rang-etendu-star} et donc \ref{prop-rang-etendu-ast}.

\`A partir de là, l'assertion (\ref{prop-rang-etendu-aii}) résulte par un argument inductif qui affirme que $\W^{(>a)}$ la vérifie.

(\ref{prop-rang-etendu-b}) pour $\rkm(\W)$.  On a tout de suite que
$$\varPhi_{t}(\phi_m)\text{ est injectif,}\quad\forall t<\poids(\W)\,,\ \forall m\,,\eqnnb\label{prop-rang-etendu-eq-0}$$ 
tout simplement parce que 
$\varPhi_{t}(\W)=0$ pour tout $t<\poids(\W)$ (\ref{lemme-co-invariants}-(\ref{lemme-co-invariants-ci})).

Pour étudier les cas où $t\geq\poids(\W)$, on rappelle que d'après
\ref{prop-V-lambda}-(\ref{prop-V-lambda-b-iii}), 
$\rkm(\W)$ est le plus petit des $s\in\NN$ tel que 
$$\varUpsilon(n):\V(\W,n)\to\W
\text{\quad est injectif\quad }\forall n\geq s\,.\eqnnb\label{prop-rang-etendu-eq-2}$$
Si nous appliquons, pour $\poids(\W)\leq t\leq n$,
le foncteur $\varPhi_t(\_)$ au diagramme
$$\preskip0pt\xymatrixc{@C=1.5cm@R=5mm}{
\V(\W,n)\ariinto[d]_{\varUpsilon(n)}&V(\W,n)_n\aregal[d]_{\varUpsilon(n)_{n}}\ar[r]^{\phi_{n',n}}&V(\W,n)_{n'}\strut\ariinto[d]^{\varUpsilon(n)_{n'}}\\
\W&W_n\ar[r]^{\phi^{W}_{n',n}}&W_{n'}
}\postskip0pt$$
on obtient
$$\xymatrixc{@C=0.8cm@R=5mm}{
\varPhi_{t}(\V(\W,n))\ar[d]_{\varPhi_{t}(\varUpsilon(n))}&
\varPhi_{t}(V(\W,n)_n)
\aregal[d]_{\varPhi_t(\varUpsilon(n)_{n})}\arinto[r]&\varPhi_t(V(\W,n)_{n'})\ariinto[d]^{\varPhi_t(\varUpsilon(n)_{n'})}\\
\varPhi_{t}(\W)&\varPhi_{t}(W_n)\ar[r]^{\varPhi_{t}(\phi^{W}_{n',n})}&\varPhi_t(W_{n'})
}\eqnnb\label{prop-rang-etendu-eq-3}$$
où clairement $\varPhi_t(\phi^{W}_{n',n})$ est injective. Par conséquent, et compte tenu de \ref{prop-rang-etendu-eq-0}, on conclut que l'on a
$$\rkme(\W)\leq \rkm(\W)\,.\label{prop-rang-etendu-eq-4}$$
Maintenant, si cette inégalité était stricte, il y aurait un certain $m>\poids(\W)$ tel que $\phi_m:W_m\to W_{\mm}$ n'est pas monotone alors que $\varPhi_{t}(\phi_m)$ serait injective pour tout $t\leq\poids(\W)$, ce qui est clairement impossible déjà pour
$t=\poids(\W)$. L'inégalité \ref{prop-rang-etendu-eq-3} est donc une égalité.

(\ref{prop-rang-etendu-b}) pour $\rkms(\W)$. On procède pareillement en remarquant cette fois que dans \ref{prop-rang-etendu-eq-3} le morphisme $\varPhi_{t}(\phi_{n',n})$ est bijectif pour $n\geq\rkms(\W)$ et $t=\poids(\W)$ auquel cas, $\varPhi_t(\phi^{W}_{n',n})$ est bijective si et seulement si, $\varUpsilon(n):\V(\W,n)\to\W_{\geq n}$ l'est. On conclut ensuite en appliquant
\ref{prop-V-lambda}-(\ref{prop-V-lambda-b-ii}). 
\enddemo

\def\Demonstration{Preuve de la proposition \ref{rkms=rke}}
\demo Le cas non trivial est celui où $\poids:=\poids(\W)<+\infty$. Dans ce cas, si $\set \W^{(t)}/$ est la famille des sous-quotients à poids unique de $\W$, on a vu dans (\ref{prop-rang-etendu-aii}) du lemme et dans \ref{prop-scindage-par-le-poids}-(\ref{prop-scindage-par-le-poids-c}) que l'on a
$$\rke(\W)=\sup\nolimits_{t\leq\poids}\set\rke \W^{(t)}/
\text{\quad et\quad}
\rk(\W)=\sup\nolimits_{t\leq\poids}\set\rk \W^{(t)}/$$
avec $\rk\in\set\rkm,\rkms/$. La proposition résulte alors de 
l'assertion (\ref{prop-rang-etendu-b}) du lemme qui établit les égalités
$
\rkm(\W^{(t)})=\rkme(\W^{t})$
et
$\rkms(\W^{t})=\rkmse(\W^{t})
\,.$
\comment
Lorsque $\W$ est à poids unique, l'assertion est celle du lemme \ref{prop-rang-etendu}-(\ref{prop-rang-etendu-b}) et nous n'avons rien à prouver. Dans le cas général, on note $t=\poids$ et l'on considère la suite exacte de $\FI$-modules
$$0\to\W^{(t)}\to\W\to\W^{(<t)}\to0\,,$$
dont on sait (\ref{prop-scindage-par-le-poids}) que
$$
\rk(\W)=\sup\bigset \rk(\W^{(t)}),\rkms(\W^{(<t)})/
$$
avec $\rk\in\set\rkm,\rkms/$.

\medskip
\noindent{\slshape (A) L'inégalité $\rkmse(\W)\leq\rkms(\W)$. }

\smallskip
\hskip0.8em
L'application itérée de \ref{prop-rang-etendu}-(\ref{prop-rang-etendu-a}) mène à la majoration
$$\rkmse(\W)\leq\supnl_{t\in\NN}\bigset\rkmse(\W^{(t)})/\,,\eqno(\ast)$$
où $\W^{(t)}$ est le sous-quotient de $\W$ à poids unique $t$ de \ref{prop-scindage-par-le-poids}-(\ref{prop-scindage-par-le-poids-a-iii}). Or,  on a $\rkmse(\W^{(t)})=t+\degstab_{t}(\W^{(t)})=\rkms(\W^{(t)})$ par  \ref{lemme-co-invariants}-(\ref{lemme-co-invariants-cii}), et l'inégalité $(\ast)$ se complète en 
$$\supnl_{t\in\NN}\bigset\rkmse(\W^{(t)})/=
\supnl_{t\in\NN}\bigset\rkms(\W^{(t)})/=\rkms(\W)\,,$$
d'après \ref{prop-scindage-par-le-poids}-(\ref{prop-scindage-par-le-poids-c}), ce qui prouve (A).

\medskip
\noindent{\slshape (B) L'inégalité $\rkms(\W)\leq\rkmse(\W)$. }

\smallskip
\hskip0.8em
Soit $\W=\set\phi_m:W_{m}\to W_{\mm}/$ où 
$W_m=\bigoplusnl_{|\mu|\leq \poids} c(\mu)_m V(\mu)_m$
est la décomposition de $W_m$ en composantes irréductibles.
Nous allons prouver~(B) en raisonnant par induction sur le cardinal de l'ensemble $$\Pi(\W):=\set|\mu|\mid\exists m\in\NN \text{ t.q. } c(\mu)_{m}\not=0/\,.$$

Lorsque $|\Pi(\W)|=1$, (B) est l'assertion \ref{prop-rang-etendu}-(\ref{prop-rang-etendu-b}). Pour le cas général, on utilise le scindage par le poids de \ref{prop-scindage-par-le-poids} en $t=\inf(\Pi(\W))$:
$$0\to\W^{>t}\hook\W\onto\W^{(t)}\to0$$

\noindent D'après (\ref{lemme-co-invariants}-(\ref{lemme-co-invariants-c})-(i,ii)), les foncteurs $\Phi_{t}(\_)^{\mt}$ annulent les termes de $\W^{(>t)}$ et rendent isomorphes ceux de $\W$ et $\W^{(t)}$ pour tout $m\geq t$. On a donc
$$\rkms(\W^{(t)})=\rkmse(\W^{(t)})\leq \rkmse(\W)\,,\postskip4pt
\eqno(\dagger)$$
et aussi
$$\postskip4pt\W^{(t)}_{\geq \rkmse(\W)}=\bigoplusnl_{|\mu|=t}c(\mu)_{\rkmse(\W)}\V(\mu)_{\geq \rkmse(\W)}\,,\postdisplaypenalty10000\eqno(\diamond)$$
d'après \ref{prop-V-lambda}-(\ref{prop-V-lambda-b-ii}).

\smallskip
Ensuite, pour tout $t'\in\NN$ tel que $t<t'\leq\poids$, les morphismes
$$(\phi_{m}):\Phi_{t'}(W_m)^{\mt'}\to\Phi_{t'}(W_\mm)^{\mmt'}$$
sont bijectifs si $m\geq\rkmse(\W)\geq t'+\degstab_{t'}(\W)$ de sorte que 
les morphismes 
$$(\tilda\phi_{m}):\Phi_{t'}(W^{(t)}_m)^{\mt'}\to\Phi_{t'}(W^{(t)}_{\mm})^{\mmt'}$$
induits par passage au quotient sont à priori surjectifs pour ces mêmes valeurs de $m$. Or, compte tenu de $(\diamond)$ et par \ref{lemme-co-invariants}-(\ref{lemme-co-invariants-b}), le nombre des composantes irréductibles avec multiplicités des $\S_{t'}$-modules \smashtop{$\Phi_{t'}(W^{(t)}_{m})^{\mt'}$} ne décroît pas avec $m$, et ces surjections sont alors nécessairement des bijections. L'exactitude de $\Phi_{t'}$ permet ensuite de conclure que $t'+\degstab_{t'}(\W^{(>t)})\leq \rkmse(\W)$, et ce pour tout $t'\leq\poids$. Par conséquent,
$$\rkmse(\W^{(>t)})\leq \rkmse(\W)\,,$$
et $\W^{(>t)}$ vérifie (B) par hypothèse de récurrence puisque
$|\Pi(\W^{(>t)}|<|\Pi(\W|$, on a donc
$$\rkms(\W^{(>t)})\leq\rkmse(\W^{(>t)})\leq \rkmse(\W)\,.\eqno(\dagger\dagger)$$

\smallskip
Les inégalités $(\dagger)$ et $(\dagger\dagger)$ et le
fait que pour $s:=\rkmse(\W)$ la catégorie $\ModFIsc{\rkms\leq s} k$ est stable par extensions \ref{coro-monotonie-stabilite}-(\ref{coro-monotonie-stabilite-b}), se conjuguent pour justifier que l'on a
$$\rkms(\W)\leq\rkmse(\W)\,,$$
ce qui termine la preuve de (B).
\endcomment
\vskip-1.5em\enddemo

\comment
\subsubsectionline{Suite du commentaire \ref{comment-poids-degre-2}.}Concernant\label{meilleure-borne} la majoration de la proposition 3.3.3 de \cite{cef}
$\rks(\W)\leq\poids(\W)+\degstab(\W)$ (\cf note ($^{\ref{note-deg-stab}}$)), la proposition \ref{rkms=rke} montre que l'on a
$$\rks(\W)\leq\rkms(\W)=\rkmse(\W)\leq\poids(\W)+\degstab(\W)\,.
$$

Lorsque  $\W=\V(\lambda)$ toutes les inégalités sont des égalités, mais c'est loin d'être le cas en général. En effet, si $\W=\V(\lambda)\oplus\V(\mu)$, on a 
$$\poids(\W)=\sup(|\lambda|,|\mu|)\,\quad\degstab(\W)=\sup(\lambda_1,\mu_1)\,,$$
et donc
$\rkms(\W)=\sup(|\lambda|+\lambda_1,|\mu|+\mu_1)=\rkmse(\W)$ (\ref{prop-rang-etendu}). 
Or, le $\sup$ d'une somme est généralement strictement plus petit 
que la somme des $\sup$'s. 
\endcomment

\comment Par conséquent $\rkmse(\W)$ approche en général mieux $\rkms(\W)$ que ne le fait $\poids(\W)+\degstab(\W)$.\endcomment

\comment
\def\Remarque{Preuve de \ref{FI-TF=>chi-pol}}
\begin{rema}Ce\label{827y823=>824} théorème est corollaire de \ref{prop-M(W)}-(\ref{prop-M(W)-c}). En effet, un $\FI$-module de type fini est quotient d'une somme finie $\V=\bigoplus_i\M(\ag_i)$, chaque $\M(\ag_i)$ étant stationnaire de dimensions finies, les quotients de $\V$ sont forcément stationnaires de dimensions finies, donc à caractères éventuellement polynomiaux. On remarquera que cette justification de donne pas de renseignement utile sur le rangs de stabilité et monotonie.

Un $\FI$-module de type fini est donc toujours stationnaire, et la réciproque, que nous n'allons pas utiliser, est également vraie. 

Le théorème suivant qui résume les dernières remarques est un théorème central dans la théorie des $\FI$-modules.\end{rema}

\mynobreak\nobreak\begin{theo}[(\cite{cef} 1.13, Finitude et stabilité)]\label{FI-TF=>rep-stat}Un $\FI$-module 
{\comment(sur un corps de caractéristique nulle)\endcomment} est de type fini si et seulement si, il est stationnaire et de dimensions finies. 

Plus précisément, pour chaque $\FI$-module de type fini $\V$, il existe une famille de décompositions $\set\lambda_1\vdash a_1,\ldots,\lambda_r\vdash a_r/$ telle que pour $m\gg0$, on a 
$$V_{m}\simeq 
c(\lambda_1)V(\lambda_1)_{m}\oplus\cdots\oplus c(\lambda_r)V(\lambda_r)_{m}\,, $$
où les multiplicités $c(\lambda_i)$ sont indépendantes de $m$.

\nobreak La famille des caractères $\chi(\V)$ est éventuellement polynomiale de même (donc) que la famille des dimensions $\dim_k (\V)$ (\ref{remas-caractere-poly}-(\ref{suite-de-dimensions})).
\end{theo}

\subsection{Majoration du rang d'un $\FI$-module}
\subsubsectionline{Rang d'un $\FI$-module.}Le \expression{rang (de monotonie et stabilité) d'un $\FI$-module $\V$} a déjà été introduit dans \ref{monotonie}, il est par définition\glossary{${\rkms(\V):=\sup\bigset\rkm(\V),\rks(\V)/}$:rang (de monotonie et stabilité) de $\V$}
$$\rkms(\V):=\sup\bigset\rkm(\V),\rks(\V)/\,.$$
Le mot \expression{stabilité} sera utilisé un peu abusivement comme raccourci de \expression{monotonie et stabilité}, le contexte permettant de lever l'ambiguïté.

\begin{prop}Soit\label{critere-stabrk} $\W:\set\phi_m:W_{m}\to W_{\mm}/$ un $\FI$-module de poids $\poids:=\poids(\W)$ fini. Pour $s\geq\poids{}$ donné,  les conditions suivantes sont équivalentes.
\def\varlistskips{\topsep4pt\mou\itemsep4pt\mou}
\begin{enumerate}
\item Pour\label{critere-stabrk-a} tout $m\geq s$
et tout $L\dans\iii[1,m]$ tel que $|L|\geq m-\poids$,
les morphismes  de transition
$$
\phi_{m}:
(W_m)_{\SL}
\to
(W_{\mm})_{\SLp}
$$
sont des isomorphismes de $\S_{\iii[1,m]\moins L}$-modules.
\item Si\label{critere-stabrk-b} $m\geq s\geq\poids\geq t\geq 0$, les morphismes de transition induits
$$(\phi_{m}):
(W_m)_{\S_{m-t}}
\to
(W_{\mm})_{\S_{m+1-t}}
$$
sont des isomorphismes de $\St $-modules.
\item $\rkms(\W)\leq s$\label{critere-stabrk-c}.
\end{enumerate}
\end{prop}

\demo L'équivalence (\ref{critere-stabrk-a})$\Leftrightarrow$(\ref{critere-stabrk-b}) est élémentaire. La preuve de 
 (\ref{critere-stabrk-b})$\Leftrightarrow$(\ref{critere-stabrk-c}) utilise 
la proposition \ref{lemme-co-invariants}.

\noindent (\ref{critere-stabrk-b})$\Rightarrow$(\ref{critere-stabrk-c}). Pour $m\in\NN$, décomposons $W_m$ en facteurs irréductibles
$$W_m=\sumnl_{|\mu|\leq \poids} c(\mu)_m V(\mu)_m\,.$$
Pour chaque $0\leq t$ donné, les sous-$\Sm $-modules
$$W^{(\geq t)}_{m}:=\sumnl_{\poids\geq |\mu|\geq t} c(\mu)_m V(\mu)_m\ \dans\ W_m\,,$$
constituent un sous-$\FI$-module $\W^{(\geq t)}\dans\W$ puisque les composantes irréductibles des sous-modules $\Smm\cdot \phi_{m}(V(\mu)_{m})\dans W_{\mm}$, pour $|\mu|\geq t$, sont de poids minorés par $|\mu|$ d'après la règle de branchement pour les représentations des groupes symétriques qui dit que lors des inductions de $\Sm$ à $\Smm$ les poids ne baissent pas  (\cf p.~\pageref{regles-de-branchement}).
Nous avons alors la suite exacte courte de $\FI$-modules
$$
\postskip1em
\def\tt{t{+}1}
\def\sumnl_#1{\halfsmash{\displaystyle\sum_{#1}}}
\xymatrixc{@C=7mm@R=6mm}{
\W^{(>t)}\arinto[d]&\Big(\ar[r]^(0.22){\phi_{\mmo}}&\sumnl_{|\mu|> t}c(\mu)_m V(\mu)_{m}\ar[r]^(0.45){\phi_{m}}\arinto[d]&\sumnl_{|\mu|> t}c(\mu)_{\mm} V(\mu)_{\mm}\arinto[d]\ar[r]^(0.77){\phi_{\mm}}&\Big)\\
\W^{(\geq t)}\aronto[d]^(.42){\pi_{t}}&\Big(\ar[r]^(0.22){\phi_{\mmo}}&\sumnl_{|\mu|\geq t}c(\mu)_m V(\mu)_{m}\ar[r]^(0.45){\phi_{m}}\aronto[d]^{\pi_{t,m}}
&\sumnl_{|\mu|\geq t}c(\mu)_{\mm} V(\mu)_{\mm}\aronto[d]^{\pi_{t,\mm}}
\ar[r]^(0.77){\phi_{\mm}}&\Big)\\
\W^{(t)}&\Big(\ar[r]^(.22){\tilde\phi_{\mmo}}&\sumnl_{|\mu|=t}c(\mu)_m V(\mu)_{m}\ar[r]^(0.45){\tilde\phi_m}&
\sumnl_{|\mu|=t}c(\mu)_{\mm} V(\mu)_{\mm}\ar[r]^(0.75){\tilde\phi_{\mm}}&\Big)
}\eqno(\dagger)$$
où $\tilde\phi_m$ est le morphisme de transition induit par $\phi_m$ et où 
$\pi_{t}:\W^{(\geq t)}\onto\W^{(t)}$ est la surjection canonique sur le quotient de $\FI$-modules $\W^{(t)}:=\W^{\geq t}/\W^{>t}$.  

Ceci étant, la règle de branchement montre que
dans le $\FI$-module quotient $\W^{(t)}$
on a des inclusions $$\tilde\phi_{m}(c(\mu)_{m}V(\mu)_{m})\dans c(\mu)_{\mm}V(\mu)_{\mm}\,.\eqno(\ddagger)$$ 

Supposons maintenant que $\W^{(\geq t)}$ vérifie l'assertion (\ref{critere-stabrk-b}). En appliquant les foncteurs  $(\_)_{\S_{*{-}t}}$ aux colonnes de $(\dagger)$, 
les assertions \ref{lemme-co-invariants}-(i,ii) montrent que la première ligne est annulée et les deux autres sont isomorphes, auquel cas
l'hypothèse (\ref{critere-stabrk-b}) exprime le fait que les morphismes de $\S_t$-modules
$$
\Big(\sumnl_{|\mu|=t} c(\mu)_m V(\mu)_m\Big)_{\S_{m-t}}
\hf{\tilde\phi_{m}}{}{1.2cm}
\Big(\sumnl_{|\mu|=t} c(\mu)_{\mm} V(\mu)_{\mm}\Big)_{\S_{m+1-t}}\postdisplaypenalty10000$$
sont bijectifs pour tout $m\geq s$. 
On en déduit l'existence d'un isomorphisme de $\FI$-modules tronqués
$$\W^{(t)}_{\geq s}\simeq\sumnl_{|\mu|=t}c(\mu)_{s}\V(\mu)_{\geq s}\,,$$
d'où le fait que 
$\rkms(\W^{(t)}_{\geq s})= s\,.$
Mieux encore, \ref{lemme-co-invariants}-(iii) garantit que $\W^{(t)}_{\geq s}$ vérifie (\ref{critere-stabrk-b}) pour $t'\geq t$, de même, donc, que  $\W^{(\geq t{+}1)}_{\geq s}$ d'après l'exactitude des foncteurs de co-invariants.

\smallskip
Un raisonnant par récurrence croissante à partir de $t=0$, prouve alors que les sous-$\FI$-modules $\W^{(\geq t)}_{\geq s}$ vérifient (\ref{critere-stabrk-b}), et (donc) que 
$\rkms(\W^{(t)}_{\geq s})= s$, et ceci, quel que soit $t\in\NN$. 

Enfin, le fait que $\ModFIsc{\rkms\leq s} k$ est une sous-catégorie abélienne stable par extensions de $\ModFI k$ intervient pour conclure que dans les suites exactes
$$\xymatrix{
\mllap{0\to}\W^{(>t)}\arinto[r]&
\W^{(\geq t)}\aronto[r]^{\pi}&
\W^{(t)}\mrlap{\to0}
}$$
on a $\rkms(\W^{(t)})\geq s$, si et seulement si,
$\rkms(\W^{(t+1)})\geq s$. Un raisonnement par récurrence décroissante sur $t$, vu que $W^{(>\poids)}=0$, permet alors de conclure que
$$\rkms(\W)=\rkms(\W^{(\geq 0)})\leq s\,.$$


\medskip\noindent(\ref{critere-stabrk-c})$\Rightarrow$(\ref{critere-stabrk-b}) Pour $s\geq\rkms(\W)$ on a  $s\geq\poids$, et pour tout $m\geq s$ on a
$$W_m=\sumnl_{|\mu|\leq \poids} c(\mu)\, V(\mu)_m\,,$$
où les multiplicités $c(\lambda)$ sont indépendantes de $m$.
Nous devons alors montrer que  le morphisme
$$\Big(\sumnl_{|\mu|\leq \poids}c(\mu) V(\mu)_{m}\Big)_{\S_{m-t}}
\hf{(\phi_{m})}{}{1.2cm}
\Big(\sumnl_{|\mu|\leq \poids}c(\mu) V(\mu)_{\mm}\Big)_{\S_{m+1-t}}\eqno(\dagger\dagger)$$
est un isomorphisme pour tout $0\leq t\leq\poids$. 

\smallskip
Soit $t_0:=\inf\set |\mu|\mid c(\mu)\not=0/$. En procédant comme dans (\ref{critere-stabrk-b})$\Rightarrow$(\ref{critere-stabrk-c}), on considère la surjection de $\FI$-modules
$$\pi_{t_0}:\W\onto\W^{(t_0)}:=
\bigset \tilde\phi_{m}:W_m^{(t_0)}\to W_{\mm}^{(t_0)}/_{m}\,,$$
à laquelle on applique les foncteurs $(\_)_{\S_{\ast{-}t}}$.
On note 
$$(\tilde\phi_{m}):=(W_{m}^{(t_0)})_{\S_{\mt}}\to
(W_{\mm}^{(t_0)})_{\S_{\mm{-}t}}$$
le morphisme de transition induit.

Pour $0\leq t <t_0$, le morphisme $(\tilde\phi_{m})$ est l'identité entre les modules nuls. Pour $t=t_0$ et $m\geq s$, le morphisme de $\S_{t_0}$-modules
$$\sumnl_{|\mu|= t_0}c(\mu) V_\mu
\hf{(\tilde\phi_{m})}{}{1.3cm}
\sumnl_{|\mu|= t_0}c(\mu) V_\mu\eqno(\ddagger\ddagger)$$
est surjectif puisque $m$ est dans le rang de monotonie de $\W$. Il est alors bijectif aussi puisque les multiplicités $c(\mu)$ dans $(\ddagger\ddagger)$ coïncident. 
Il s'ensuit que \smashbot{$\W_{\geq s}^{(t_0)}$} vérifie  (\ref{critere-stabrk-b}), et $\rkms(\W^{(t_0)})\geq s$ d'après ((\ref{critere-stabrk-b})$\Rightarrow$(\ref{critere-stabrk-c})). 
Par conséquent, $\rkms(\ker(\pi_{t_0}))\leq s$,
puisque noyau d'un morphisme entre deux $\FI$-modules monotones et stables pour $m\geq s$  \ref{coro-monotonie-stabilite}-(\ref{coro-monotonie-stabilite-b}). 
Or, comme l'ensemble des poids des facteurs irréductibles des termes de $\ker(\pi_{t_0})$ est strictement inférieur à celui de $\W$, supposé de poids fini, 
on peut raisonner par récurrence sur le cardinal de ces ensembles et considérer que les analogues des morphismes $(\dagger\dagger)$ pour $\ker(\pi_{t_{0}})$ sont des isomorphismes, ce qui permet aisément de conclure.
\enddemo

\subsubsectionnumber L'implication \ref{critere-stabrk}-(\ref{critere-stabrk-c})$\Rightarrow$(\ref{critere-stabrk-a}) donne une condition sur $m$ garantissant aux morphismes de transition
$$\phi_{m}:
(W_m)_{\SL }
\to
(W_{\mm})_{\SLp }
$$
d'être bijectifs pour tout $L\dans\iii[1,m]$ tel que $\poids(\W)\geq m-|L|$. Dans la preuve du théorème \ref{theo-Ind-lambda}-(\ref{theo-Ind-lambda-d}), nous aurons besoin d'une condition sur $m$ mais pour des valeurs de $m-|L|$ arbitrairement grands. La proposition suivante traite de cette situation.

\begin{prop}Soit\label{critere-c=>a} $\W:\set\phi_m:W_{m}\to W_{\mm}/$ un $\FI$-module de rang de monotonie et stabilité $s:=\rkms(\W)<+\infty$. On fixe $t\in\NN$, et l'on considère
la famille de morphismes pour $m\geq t$,
$$\phi_{t,m}:
(W_m)_{\S_{m-t}}
\to
(W_{\mm})_{\S_{m+1-t}}
$$
\begin{enumerate}
\item Si\label{critere-c=>a-a} $t\leq \poids(\W)$, le morphisme $\phi_{t,m}$ est bijectif pour  $m\geq s$.
\item Si\label{critere-c=>a-b} $t>\poidsinf(\W)$, le morphisme $\phi_{t,m}$ est bijectif pour $m\geq\sup(t,s)$.
\end{enumerate}
\end{prop}
\demo (\ref{critere-c=>a-a}) C'est \ref{critere-stabrk}-(\ref{critere-stabrk-b}). (\ref{critere-c=>a-b}) C'est \ref{lemme-co-invariants}-(iii).
\enddemo
\endcomment

\section{Foncteurs d'induction dans $\ModFI k$}
Nous\label{foncteurs-decalages-inductifs} introduisons  certains foncteurs d'\expression{induction} dans $\ModFI k$ qui \expression{recollent} naturellement les foncteurs $\Ind_{\Glambda }^{\,\S _{\smash{|\lambda|}}}$, $\Ibf _{m}^{m+a}$ et $\Thetabf_{m}^{m+a}$ des
sections \ref{operateurs-inductions} et \ref{inductions-iterees}. Nous étudions ensuite la perturbation du rang de monotonie et stabilité des $\FI$-modules sous l'influence de ces foncteurs (\cf\ref{theo-Ind-lambda}-(\ref{theo-Ind-lambda-d})).

\subsection{Les foncteurs d'induction $\Ind_{\lambda}$ et $\Ind_{\rho,\lambda}$}\label{Ind-lambda-rho}

\begin{lemm}[et définition]Pour\label{a-2a}
$\lambda=(1^{\XX_1},2^{\XX_2},\ldots,m^{\XX_{m}})\vdash m\,$ et $n\geq m$\glossary{$\lambda\bxx{n}$:si $\lambda=(1^{\XX_1},2^{\XX_2},\ldots)\vdash m$, on pose $\lambda\bxx{n}=(1^{\XX_1+n-m},2^{\XX_2},\ldots)$}, on définit $\lambda\bxx n:=(1^{(n-m){+}\XX_1},2^{\XX_2},\ldots,m^{\XX_{m}})\vdash n\,.$
Alors, pour tous $a\in\NN$ et $n\geq m\in\NN$,
l'application
$$
\let\cd\cdot
\yght=1.5ex\ygwd=\yght
\def\paren#1#2#3#4{\llap{\smash{\raise#1\hbox{$ #3\left\{\vrule height#2mm width0pt\right.\,$}#4}\ }}
\def\parend#1#2#3#4{\rlap{\ \smash{#4\raise#1\hbox{$\,\left.\vrule height#2mm width0pt\right\}#3$}}}}
\setbox1=\hbox{\hskip1.2cm\vtop{\kern-1.7em\hbox{\tableau{}{\cd&\cd&\cd&\cd&\cd\\\cd\paren{-4pt}6{m-a}{}&\cd&\cd\\\cd&\cd\\\cd}}}}
\setbox2=\hbox{\hskip1.2cm\vtop{\kern-1.7em\hbox{\tableau{}{\cd&\cd&\cd&\cd&\cd\\
\cd\paren{-4pt}6{m-a}{}&\cd&\cd\\\cd&\cd\\
\cd\parend{0pt}{10.05}{n-a}{\kern1.0cm}\\\cellcolor{0 0 0 0.10}\\
\vdots\paren{0pt}4{n-m}{\kern0.5pt}\\ }{\kern3pt}}}}
\xymatrix@R=0mm{
\Y_{m-a}(m)\ar[r]&\Y_{n-a}(n)\\
\lambda\ar@{|->}[r]&\lambda\bxx{n}\\
\copy1\ar@{|->}[r]+<-40pt,0pt>&\copy2}
$$
est bijective dès que $m\geq2a$. On a donc
$|\Y_{a}(2a)|=|\Y_{m-a}(m)|\,.$
\end{lemm}
\demo En effet, les décompositions de $\lambda\vdash m$
vérifient l'inégalité 
$$m=\XX_1+2\XX_2+\cdots,m\XX_m\geq  \XX_1+2(\XX_2+\cdots+\XX_m)=
 \XX_1+2(m-a-\XX_1)\,,$$
dont on déduit que $\XX_1\geq m-2a$. Il en résulte que si $m\geq2a$ toute décompositions $\lambda\vdash m$ comporte au moins $m-2a$ singletons, le fait de les enlever ou de les rajouter établit la bijection entre $\Y_{m-a}(m)$ et $\Y_{a}(2a)$.
\enddemo


\noendpoint\begin{remas}[et notations]\label{rema-lambda/m}
\begin{enumerate}
\mynobreak\nobreak\item Dans\label{rema-lambda/m-a} la notation 
`$\lambda\bxx n$', le nombre `$n$', qui vérifie $n\geq|\lambda|$,
indique la taille finale du diagramme. 
On a (comparer à \ref{remas-lambda[m]}-(\ref{remas-lambda[m]-b})):
$$\ell(\lambda\bxx n)=\ell(\lambda)+(n-|\lambda|)\,,
$$
et si $n\leq n'$
$$(\lambda\bxx n)\bxx {n'}=\lambda\bxx {n'}\,.$$
\item \'Etant\label{rema-lambda/m-b} donné $\lambda=(1^{\XX_1},2^{\XX_2},3^{\XX_3},\ldots)\in\Y_{\ell}(\ell+a)$, on pose\glossary{$\ulambda$:Si $\lambda=(1^{\XX_1},2^{\XX_2},3^{\XX_3},\ldots)$, on pose $\ulambda=(1^{0},2^{\XX_2},3^{\XX_3},\ldots)$} 
$$\uell:=\ell-\XX_1\text{\quad et\quad }
\ulambda=(1^{0},2^{\XX_2},3^{\XX_3},\ldots)\in\Y_{\uell}(\uell+a)\,.$$
Les données suivantes sont associées au diagramme $\lambda\bxx m$\label{notations-lambda-m}:
$$
\setbox110=\hbox{$(\lambda)$}
\setbox112=\hbox{$(\lambda\bxx m)$}
\setbox114=\hbox{$(\ulambda)$}
\hskip0.5cm\tableauc{}{
&\cdot&\cdot&\cdot&\cdot&\cdot\\
&\cdot&\cdot&\cdot\\
\lparcell{-7pt}{\copy110\ }{1.08cm}{2mm}
\rparcell{-7pt}{\ell}{1.05cm}{1.9cm}
&\cdot\\
\\
\rparcell{-2pt}{\XX_1}{5mm}{2mm}\\
\lparcell{3pt}{\copy112\ }{1.65cm}{1.3cm}
\\\cellcolor{0 0 0 0.10}
\\
\rparcell{-7pt}{m-|\lambda|}{7mm}{2mm}\\
\vdots\\
}
\hskip1.5cm
\vcenter{
\mhbox{
\hbox to2em{$(\ulambda)$\hss}\tableauc{}{
&\cdot&\cdot&\cdot&\cdot&\cdot\\
\rparcell{-2pt}{\hbox{\normalsize$\,\uell:=\ell-\XX_1$}}{5mm}{1.9cm}
&\cdot&\cdot&\cdot\\
&\cdot
}}
\vbox to 2.2cm{
\bigskip
\hbox{\qquad$|\lambda|=\ell+a$}
\hbox{\qquad$|\ulambda|=|\lambda|-\XX_1\leq 2a$}
\smallskip
\hbox{\qquad\yght=1.5ex\ygwd=1.5ex
$\#\set\tableauc{}{\cdot}/=a$}
\vss}
}$$
On a en particulier: $\lambda\bxx m= \ulambda\bxx m$, $\forall m\geq|\lambda|$, ce qui nous conduit à étendre la portée des notations en posant
$$\lambda\bxx m:= \ulambda\bxx m\,,\quad\forall m\geq|\ulambda|=|\lambda|-\XX_1\,.$$

\end{enumerate}
\end{remas}
\subsubsectionline{Le foncteur $\Ind_{\lambda}$  sur $\ModFB k$.}Pour\label{Ind-lambda-FB}
$\lambda\in\Y_{\ell}(\ell+a)$ donné, 
$$\relax{\Ind_{\lambda}:\ModFB k\fonct\ModFB k_{\geq |\ulambda|}}
$$
est le foncteur qui fait correspondre à un $\S_{m-a}$-module $V_{m-a}$, le $\Sm$-module
$$W_{m}:=\begin{cases}\noalign{\kern-4pt}
0\,,\hbox{\quad si $m< |\ulambda|\,,$}\\
\Ind_{\myG_{\ulambda\bxx{m}}}^{\Sm}V_{m-a}\,,\hbox{\quad si $m\geq |\ulambda|$.}\\\noalign{\kern-4pt}
\end{cases}$$
et qui fait correspondre à un morphisme de $\FI$-modules $f:\V\to\W$, la famille de morphismes:
$$\Ind_{\lambda}(f):=\Bigset\Ind_{\myG_{\ulambda\bxx{m}}}^{\Sm}(f_{m{-}a}):
\Ind_{\myG_{\ulambda\bxx{m}}}^{\Sm}V_{m-a}\to
\Ind_{\myG_{\ulambda\bxx{m}}}^{\Sm}W_{m-a}/{\vrule depth5pt width0pt}_{m\geq|\ulambda|}\,.
$$

La proposition suivante étend la définition de $\Ind_{\lambda}$ à la catégorie des $\FI$-modules.
Le théorème \ref{theo-Ind-lambda} établit ensuite le fait que le foncteur est exact et respecte la finitude des $\FI$-modules. Le théorème décrit également son influence sur le rang de monotonie et stabilité des $\FI$-modules.

\noendpoint\begin{prop}[et définition de $\Ind_{\lambda}$ sur $\ModFI k$]\\
\vrule height12pt width0pt 
Soit\label{prop-induction}\label{Ind-lambda}  $\lambda\in\Y_{\ell}(\ell+a)$ et
soient $0\leq |\ulambda|\leq m\leq n$.
\def\varlistskips{\topsep4pt\itemsep4pt}
\begin{enumerate}
\mynobreak\nobreak\item\leavevmode\label{prop-induction-a}L'inclusion
$\Sm \dans\S_n$ induit des inclusions (voir \ref{nota-young}, \ref{lemme-young})
$$
\PPP_{\lambda\bxx{m}}\dans \PPP_{\lambda\bxx{n}}\,,\quad
\S_{\lambda\bxx{m}}\dans \S_{\lambda\bxx{n}}\,,\quad
\myG_{\lambda\bxx{m}}\dans \myG_{\lambda\bxx{n}}\,,
$$
\item\leavevmode\label{prop-induction-b}Une application $k$-linéaire
$\phi_{n,m}:V_{m-a}\to V_{n-a}$
 d'un $\S_{m-a}$-module $V_{m-a}$ vers un $\S_{n-a}$-module $V_{n-a}$,  compatible à l'action de $\S_{m-a}$ et telle que 
$\im(\phi_{n,m})\dans V_{n-a}^{\II_{m-a}\times\S_{n-m}}$, admet un unique prolongement en une application $k$-linéaire
compatible à l'action de $\Sm$
$$\Ind_{\lambda}(\phi_{n,m}):
\Ind_{\myG_{\lambda\bxx{m}}}^{\Sm}V_{m-a}
\to
\Ind_{\myG_{\lambda\bxx{n}}}^{\S_{n}}V_{n-a}
\postskip0pt
$$
De plus,
$$\im(\Ind_{\lambda}(\phi_{n,m}))\dans \big(\Ind_{\myG_{\lambda\bxx{n}}}^{\S_{n}}V_{n-a}\big){}^{\II_{m}\times\S_{n-m}}\,.$$
\item\leavevmode\label{prop-induction-c}Dans  {\rm (\ref{prop-induction-b})}, si l'on suppose en plus que
$k[\S_{n-a}]\cdot\im(\phi_{n,m})=V_{n-a}$,
on a
$$k[\S_{n}]\cdot\im(\Ind_{\lambda}\phi_{n,m})=\Ind_{\lambda}V_{n-a}\,.$$
\item\leavevmode\label{prop-induction-d}Pour tout $\FI$-module $\V=\set\phi_m: V_m\to V_\mm/$, le $\FB$-module $\Ind_{\lambda}(\V)$ muni de la famille de morphismes $\Ind_{\lambda}(\phi_{m})$ est un $\FI$-module qui sera pareillement noté. De plus, si $f:\V\to\W$ est un morphisme de $\FI$-modules, les morphismes des familles $\Ind_{\lambda}(\phi_\V)$ et $\Ind_{\lambda}(\phi_\W)$ commutent à ceux de la famille $\Ind_\lambda(f)$ de \ref{Ind-lambda-FB}. Les correspondances $\V\fonct\Ind_{\lambda}(\V)$, $f\fonct\Ind_{\lambda}(f)$ définissent un foncteur
$$\displayboxit{\Ind_{\lambda}:\ModFI k\fonct\ModFI k_{\geq |\ulambda|}}\eqno(\Ind_{\lambda})$$
\end{enumerate}\end{prop}

\demo 
(\ref{prop-induction-a}) Résulte de remarquer que
$\PPP_{\lambda\bxx m}=\Pulambda \times\II_{m-|\ulambda|}$, ce qui implique
que si $\XX_r$ est le nombre des $\lambda_i=r$, on a  
$\Glambda =\S_{\XX_{\lambda_1}}\times\cdots\times\StS{\XX_2}{\XX_1}$ alors
$$\myG_{\lambda\bxx{m}}=\S_{\XX_{\lambda_1}}\times\cdots\times\StS{\XX_2}{\XX_1+m-|\lambda|}\,.$$

(\ref{prop-induction-b}) Compte tenu de (\ref{prop-induction-a}), le prolongement $\Ind_{\lambda}(\phi)$ annoncé est bien défini et est  unique. L'image de $\Ind_{\lambda}\phi$ est alors le sous-$\Sm$-module de $\Ind_{\myG_{\lambda\bxx{n}}}^{\S_{n}}V_{n-a}$ engendré par $\im(\phi)$, et le fait que $\im(\Ind_{\lambda}\phi)$ est
invariant sous l'action du groupe $\II_{m}\times\S_{n-m}$ résulte de ce que ce groupe commute à l'action de $\Sm\dans\S_n$, qu'il est contenu dans $\S_{\lambda\bxx n}$, et que, 
par la surjection $\S_{\lambda\bxx n}\onto \myG_{\lambda\bxx n}$, il est en bijection avec $\II_{m-a}\times\S_{n-m}$ qui, lui, fixe  $\im(\phi)$ par hypothèse.

(\ref{prop-induction-c},\ref{prop-induction-d}) Clairs.
\enddemo

\begin{theo}Soit\label{theo-Ind-lambda}
$\lambda\in\Y_{\ell}(\ell+a)$. Notons $\uell:=\ell(\ulambda)$.
\begin{enumerate}\itemsep4pt
\item Le\label{theo-Ind-lambda-a} foncteur
$\relax{\Ind_{\lambda}:\ModFI k\fonct\ModFI k_{\geq |\ulambda|}}
$
 de \ref{prop-induction}-(\ref{prop-induction-d}), 
est covariant, additif et exact. Il est fidèle sur la sous-catégorie $\ModFI k{\!}_{\geq\,\uell}$.
\item  Si\label{theo-Ind-lambda-b} $\V$ est (de type fini) engendré en degrés $\leq d$, le $\FI$-module $\Ind_{\lambda}(\V)$ est (de type fini) engendré en degrés $\leq \sup(d+a,|\ulambda|)$.
\item On\label{theo-Ind-lambda-c} a 
$\poids(\Ind_{\lambda}\V)\leq\poids(\V)+|\ulambda|\,.$
\item \leavevmode\relax{On a }\label{theo-Ind-lambda-d}
$
\relax{\rkms(\Ind_{\lambda}\V)\leq 
\rkms(\V)+2|\ulambda|
\text{ et }
\rkm(\Ind_{\lambda}\V)\leq 
\rkm(\V)+a
}$.
\end{enumerate}
\end{theo}
\demo 
(\ref{theo-Ind-lambda-a}) Clair. 
(\ref{theo-Ind-lambda-b}) Un $\FI$-module $\V=\set \phi_m:V_m\to V_{\mm}/_{m}$ est engendré en degrés $\leq d$, si et seulement si, 
$k[\S_{d+n} ]\cdot \im(\phi_{d+n,d})=V_{d+n}$, pour tout $d\geq 0$. L'assertion \ref{prop-induction}-(\ref{prop-induction-c})  montre alors que $\Ind_{\lambda}\V$ est bien engendré en degrés $\leq \sup(\uell,d)+a$. La deuxième partie de (\ref{theo-Ind-lambda-b}) résulte de même, puisque 
$\V$ est de type fini, si et seulement si, il est engendré en degrés $\leq d$ (pour $d$ assez grand), et $\dim_{k}V_m<+\infty$ pour tout $m\leq d$. (\ref{theo-Ind-lambda-c}) est établi dans \ref{poids-ind}. 

\medskip
\noindent (\ref{theo-Ind-lambda-d}) pour $\rkms(\Ind_{\lambda}\V)$. 

Si $\rkms(\V)=+\infty$ on n'a rien à prouver, on suppose donc $\rkms(\V)<+\infty$.

\smallskip
\noindent\relax{\slshape ($A$) Réduction au cas $\V:=\V(\mu)$.}

\nobreak \noindent  Soit $\V$ un $\FI$-module et notons pour simplifier
$$\poids:=\poids(\V)\text{\quad et\quad }s:=\rkms(\V)\,.
$$
Par \ref{prop-V-lambda}-(\ref{prop-V-lambda-c}), $\V_{\geq s}$ est extension de $\FI$-modules $\V(\mu)_{\geq s}$ où $|\mu|+\mu_1\leq s$. Maintenant, 
pout toute
extension de deux tels $\FI$-modules
$$
\xymatrix@C=3mm{
0\to\V(\mu)_{\geq s}\ar[r]&\W\ar[r]&\V(\mu')_{\geq s}\to0\,,
}
\postskip-0.2ex
$$
la suite
$$
0\to\Ind_{\lambda}(\V(\mu)_{\geq s})\to\Ind_{\lambda}\W\to\Ind_{\lambda}(\V(\mu')_{\geq s})\to0
$$
est exacte et on aura
$$\rkms(\Ind_{\lambda}\W)\leq s +2|\ulambda|\,,
\eqnnb\label{theo-star}
$$
en appliquant \ref{coro-monotonie-stabilite}-(\ref{coro-monotonie-stabilite-b}) si les extrémités vérifient cette majoration.

Or, on a $\Ind_{\lambda}(\V(\mu)_{\geq s})=\Ind_{\lambda}(\V(\mu))_{\geq s+a}$ et 
$$\preskip1.2ex\mathalign{
\rkms(\Ind_{\lambda}(\V(\mu))_{\geq s+a})&\leq&\sup\set \rkms(\Ind_{\lambda}(\V(\mu)),s+a/\hfill\\\noalign{\kern2pt}
&&\quad{}=\sup\set|\mu|+\mu_1+2|\ulambda|,s+a/\leq s+2|\ulambda|\,,
}$$
et de même pour le terme en $\mu'$, en supposant (\ref{theo-Ind-lambda-d}) vérifiée pour les extrémités.

En itérant l'idée, la majoration \ref{theo-star} est vérifiée pour toute extension finie de $\FI$-modules $\V(\mu)_{\geq s}$, donc par  $\V_{\geq s}$. On considère alors la suite exacte
$$
0\to\Ind_{\lambda}(\V_{\geq s})\to\Ind_{\lambda}\V\to\Ind_{\lambda}(\V_{< s})\to0\,,
$$
où $\Ind_{\lambda}(\V_{< s})_m=0$ pour tout $m\geq s+a$. On a donc bien
$$\rkms(\Ind_{\lambda}\V)=\rkms(\Ind_{\lambda}(\V_{\geq s}))\leq \rkms(\V)+2|\ulambda|\,.$$

\medskip
\noindent{\slshape ($B$) Le cas $\V=\V(\mu)$.}

\nobreak \noindent Nous devons prouver la majoration:
$$\rkms(\Ind_{\lambda}(\V(\mu)))\leq 
|\mu|+\mu_1+2|\ulambda|\,,
$$ 
ce pour quoi, il suffira de prouver que l'on a
$$\preskip3pt\sup_{0\leq t\leq|\mu|+|\ulambda|}
\bigset t+\degstab_{t}\big(\Ind_{\lambda}(\V(\mu))\big)/\leq
|\mu|+\mu_1+2|\ulambda|\,,\eqnnb\label{theo-diamond}$$
par \ref{rkms=rke}-(\ref{prop-rang-etendu-b}) et puisque $\poids(\Ind_{\lambda}(\V(\mu)))\leq|\mu|+|\ulambda|$ par (\ref{theo-Ind-lambda-c}).

\smallskip
Dans ces questions il sera important de savoir à partir de quel $m\in\NN$ le poids de $\Ind_{\lambda}(\V(\mu))$ est atteint. Les idées de la remarque \ref{rema-poids-ind} montrent que pour avoir 
\smashbot{$\poids(\Ind_{\myG_{\ulambda\bxx{m}}}^{\Sm}V(\mu)_{m-a})=|\mu|+|\ulambda|$}, il faut $\mu[m{-}a]_1\geq\uell+\sup\set\lambda_1,\mu_1/$ 
$$\hskip-3cm\let\cd\relax\binmuskip0mu 
\def\m#1#2{\def\donne{#1}\def\vx{\vdots}\ifx\donne\vx\rlap{\kern3.6cm$\vdots$}\else\rlap{\kern3.4cm\smash{$(#1)#2$}}\fi}
\mu[m-a]:=
\tableauc{\cellcolor{0 0 0 0.10}}{
\m{m-a-|\mu|}{\geq\uell+\sup\set\lambda_1,\mu_1/}&&&&\noborder\cdots&&&&&\\\cellcolor{0 0 0 0}
\cd\m{\mu_1}{}&\cd&\cd&\cd\\
\cd\relax\m{\vdots}{}&\cd\\
\cd\m{\vdots}{}}
\postskip2pt$$
autrement dit, il faudra que
$$\preskip3pt\displayboxit{m\geq|\mu|+|\ulambda|+\sup\set\mu_1,\lambda_1/}\eqnnb\label{theo-dagger}$$
ce que nous supposons désormais.

\smallskip
\noindent La preuve de \ref{theo-diamond} demande plusieurs étapes.

\def\ma{m{-}a}\def\mma{\mm{-}a}
\def\Sma{\S_{\ma}}

\smallskip\noindent{\slshape ($B_1$) Réduction de $\Ind_{\myG_{\ulambda\bxx{m}}}^{\Sm}V(\mu)_{m-a}$ à $\Ind_{\1_{|\ulambda|}\times\S_{m{-}|\ulambda|}}^{\Sm}V(\mu)_{m-a}$.}

\smallskip

\noindent
D'après \ref{prop-induction}, on a
$\S_{\lambda\bxx{m}}=\S_{\ulambda}\times\S_{m-|\ulambda|}$, d'où
$$\myG_{\ulambda\bxx{m}}=\myG_{\ulambda}\times\S_{m-|\ulambda|}\dans\S_{m-a}\,. $$
On considère alors le morphisme surjectif de $\Sm$-modules
$$
\xymatrixc{@R3mm}{
\Ind_{\1_{|\ulambda|}\times\S_{m-|\ulambda|}}^{\Sm}V_{m-a}\aronto[r]^(0.55){\rho_m}\aregal[d]&
\Ind_{\myG_{\ulambda\bxx{m}}}^{\Sm}V_{m-a}\aregal[d]\\
k[\Sm]\varOtimes{0em}{7pt}_{\S_{m-|\ulambda|}}V_{m-a}
\aronto[r]^(0.5){\rho_m}&
k[\Sm]\varOtimes{1em}{7pt}_{\myG_{|\ulambda|}\times\S_{m-\ulambda}}V_{m-a}
}\eqnnb\label{theo-rho}$$

On remarque ensuite que la surjection \ref{theo-rho} établit un isomorphisme de $\Sm$-modules à gauche entre 
$$
\xymatrix{\Big(k[\Sm]\varOtimes{4pt}{7pt}_{\S_{m-|\ulambda|}}V_{m-a}\Big){\vrule height7pt width0pt}^{\myG_{\ulambda}}\ar[r]_(0.55){\simeq}&
k[\Sm]\varOtimes{1em}{7pt}_{\myG_{|\ulambda|}\times\S_{m-\ulambda}}V_{m-a}}
\eqnnb\label{theo-clef}$$
où le groupe $\myG_{\ulambda}$ agit à droite sur le terme de gauche par 
$$(P\otimes v)\cdot\alpha:=P\alpha\otimes\alpha^{-1} (v)\,.$$

Maintenant, si $\phi_{\ma}:V(\mu)_{\ma}\to V(\mu)_{\mma}$ est un morphisme de transition de $\V(\mu)$, le morphisme induit
$$\xymatrix@C=1.8cm{
k[\Sm]\varOtimes{1ex}{7pt}_{\S_{m-|\ulambda|}}V_{m-a}
\ar[r]^(.47){\id\otimes\phi_{\ma}}&
k[\Smm]\varOtimes{1em}{7pt}_{\S_{\mm-|\ulambda|}}V_{\mm-a}
}
$$
est un morphisme de $\Sm\times\myG_{\ulambda}$-bimodules.
Il suffira par conséquent que le morphisme de transition de co-invariants
$$\xymatrix@C=1.8cm{
\mllap{\varPsi_{m,t}:}\Big(k[\Sm]\varOtimes{1ex}{7pt}_{\S_{m-|\ulambda|}}V_{m-a}\Big)_{\S_{\mt}}\hskip-1.5em
\ar[r]^(.45){(\id\otimes\phi_{\ma})}&
\Big(k[\Smm]\varOtimes{1em}{7pt}_{\S_{\mm-|\ulambda|}}V_{\mm-a}\Big)_{\S_{\mmt}}\hskip-1.5em
}
\eqnnb\label{theo-varPsi-m-t}$$
soit un isomorphisme de $\St $-modules (à gauche) pour qu'il soit aussi un isomorphisme de $\myG_{\ulambda}$-modules à droite et donc pour qu'il induise aussi un isomorphisme sur le sous-$\St $-module des $\myG_{\ulambda}$-invariants à droite, autrement dit, pour que le morphisme de transition induit
$$
\xymatrix@C=1.8cm{
\Big(\Ind_{\myG_{\ulambda\bxx{m}}}^{\Sm}V(\mu)_{m-a}\Big)_{\S_{\mt}}
\hskip-1.5em\ar[r]^(.45){(\ind_{\lambda}\phi_{\ma})}&
\Big(\Ind_{\myG_{\ulambda\bxx{\mm}}}^{\Smm}V(\mu)_{\mm-a}\Big)_{\S_{\mmt}}
\hskip-1.5em}$$
soit un isomorphisme (\footnote{La réciproque n'est à priori pas vraie, et nos estimations pour les degrés de stabilité ne seront donc pas toujours optimales, elles le seront lorsque l'action de $\myG_{\ulambda}$ sur $V(\mu)_m$ est triviale, par exemple si la suite $\ulambda=(\lambda_1,\lambda_2,\ldots)$ est strictement décroissante.}).

\medskip\noindent{\slshape ($B_2$) \'Equivalences dans la description du morphisme $\varPsi_{m,t}$ dans \ref{theo-varPsi-m-t}.}

\nobreak
\noindent Pour tout $m\in\NN$ et tout $L\dans\iii[1,m]$, identifions le groupe $\SL $\glossary{$\SL $:on a $L\dans\iii[1,m]$ et $\SL :=\Fix_{\Sm}(\iii[1,m]\mmoins L)$} des bijections de $L$, au sous-groupe de $\Sm $ qui fixe  $i\not\in L$. 
Le centralisateur de $\SL $ dans $\Sm$ est le sous-groupe $\S_{\iii[1,m]\mmoins L}$ qui fixe $i\in L$. Tout comme dans \ref{co-invariants}, le foncteur des $\SL$-co-invariants est le foncteur
$$\mathalign{(\_)_{\SL  }:&\Mod(\Sm)&\fonct&\Mod(\S_{\smashtop{\iii[1,m]\mmoins L}})\hfill\\
&W&\fonct& W_{\SL}:=k\otimes_{k[\SL ]}W\,.}$$
On note ensuite $L':=L\sqcup\set m+1/\dans\iii[1,m+1]$\glossary{$L'$:pour $L\dans\iii[1,m]$, on pose $L':=L\coprod\set m+1/$}. Pour tout $\FI$-module $\W$, le morphisme de transition
$\phi_m: W_m\to W_{\mm}$ passe aux co-invariants où il définit le morphisme de $\S_{\iii[1,m]\mmoins L}$-modules 
$$\def\mt{L}\def\mmt{L'}
(\phi_{m}):
(W_m)_{\S_{\mt}}\to (W_{\mm})_{\S_{\mmt}}\,.$$ 

\smallskip
Maintenant, pour $t\leq m$ fixés,
notons 
$$L:=\iii[t+1,m]\text{\quad et\quad}A:=\iii[|\ulambda|+1,m]\,.
$$
L'étude de $\varPsi_{m,t}$ se ramène alors à l'étude du morphisme de $\St \times\myG_{\ulambda}$-bimodules
$$\def\lquo#1#2{\Otimes_{#1}k[#2]}
\xymatrix@C=1.cm{
\displaystyle
k\lquo{\SL}{\Sm }
\Otimes_{\SA}
V(\mu)_{m-a}
\ar[r]^(0.46){\varPsi_{m,t}}&
\displaystyle
k\lquo{\SLp}{\S_{m+1}}
\Otimes_{\SAp }
V(\mu)_{\mma}\,.}\hskip-0.2cm
\eqnnb\label{theo-StGl}$$

\smallskip
\noindent\begingroup\sl Lemme 1. On\label{theo-Ind-lambda-Lemme-1} a des identifications de $\S_A$-modules à droite
$$
k\otimes_{\SL}k[\Sm]=k\lquo{\SL}{\Sm }=
\bigoplusnl_{
\cl \alpha\in\big(\biquosb {\SA }{\Sm}{\SL}\big)\quad}\hskip-0.2cm
 k\lquo{\S_{A\cap\alpha L}}{\SA }\,,
$$
d'où un isomorphisme canonique de foncteurs sur $\Mod(k[\S_{\ma}])$
$$(\ind_{\SA}^{\Sm}(\_))_{\SL}\ \simeq\ \
\bigoplusnl_{\cl \alpha\in\big(\biquosb {\SA }{\Sm }{\SL }\big)}
(\_)_{\S_{A\cap\alpha L}}\,.$$\endgroup

\medskip\noindent{\sl Indication. }Les représentations en question ont les même caractères.
\hfill$\boxminus$

\medskip
Grâce à ce lemme, le morphisme \ref{theo-StGl} se voit comme le morphisme 
$$
\bigoplus_{\cl \alpha\in\big(\biquosb {\SA }{\Sm }{\SL }\big)\quad}\hskip-0.5cm
(V(\mu)_{\ma})_{\S_{A\cap\alpha L}}
\ \hf{\varPsi_{m,t}}{}{1.cm}\hskip-0.5cm
\bigoplus_{\cl \alpha\in\big(\biquosb {\SAp }{\S_{\mm}}{\SLp }\big)\quad}\hskip-1cm
(V(\mu)_{\mma})_{\S_{A'\cap\alpha L'}}
$$
induit terme à terme
 par le morphisme de transition $\phi_{\ma}$ via l'inclusion naturelle
$$\let\big\Big
\big(\raise2pt\hbox{$\biquo {\SA }{\Sm }{\SL }$}\big)\ \dans\ 
\big(\raise2pt\hbox{$\biquo {\SAp }{\S_{\mm}}{\SLp }$}\big)\,.\eqnnb\label{theo-inclusion-2-classes}$$

\medskip
\noindent\begingroup\sl Lemme 2. L'inclusion \ref{theo-inclusion-2-classes} est une égalité lorsque
$$\relax{m\geq |\ulambda|+t}\eqnnb\label{theo-lemme-1}$$
auquel cas, on a un isomorphisme de foncteurs sur $\Mod(k[\S_{\ma}])$
$$
(\ind_{\S_{m-|\ulambda|}}^{\Sm}(\_))_{\S_{\mt}}\quad\simeq\quad
\bigoplusnl_{\cl \alpha\in\big(\biquosb {\SA }{\Sm }{\SL }\big)}
(\_)_{\S_{|A\cap\alpha L|}}\,.$$
De plus, on a
$$\uell\ \leq\ \big((\ma)-|A\cap\alpha L|\big)\ \leq\
\uell+t\,.\eqnnb\label{theo-nouveau-t}$$
et $\alpha\in\Sm$ est tel que $A\cap\alpha L$ est l'intervalle des derniers $|A\cap\alpha L|$ éléments de l'intervalle $\iii[1,m]$.\endgroup

\medskip\noindent{\sl Preuve du lemme 2. }Pour toute double classe $\cl\alpha\in\big(\txtquo {\SAp }{\Smm }{\SLp }\big)$, on a
$$|A'\cap\alpha L'|\geq |L'|+|A'|-(\mm)=\big(m{-}(|\ulambda|{+}t)\big)+1\,,
$$
et lorsque \ref{theo-lemme-1} est vérifié, il existe $\alpha\in\cl\alpha$ 
tel que $L'\cap\alpha A'$ est l'intervalle des derniers $|L'\cap\alpha A'|$ éléments de $\iii[1,\mm]$, en particulier $\alpha\in\Sm$.
Dans ce cas, on a $L'\cap\alpha A'=(L\cap\alpha A)\sqcup\set\mm/$ et l'inclusion \ref{theo-inclusion-2-classes} est bien une égalité. La fin de la preuve est immédiate d'après le lemme 1.\hfill$\boxminus$

\medskip\noindent{\slshape ($B_3$) Condition d'isomorphie du morphisme $\varPsi_{m,t}$ dans \ref{theo-varPsi-m-t}.}

\nobreak\noindent Nous sommes maintenant en mesure de prouver l'inégalité \ref{theo-diamond}. D'après le lemme 2, l'application $\varPsi_{m,t}$ sera bijective lorsque $m\geq|\ulambda|+t$ et que
les
$$
(V(\mu)_{\ma})_{\S_{|A\cap\alpha L|}}
\ \hf{(\phi_{\ma})}{}{1.5cm}
(V(\mu)_{\mma})_{\S_{|A'\cap\alpha L'|}}
\eqno(\diamond)$$
sont bijectifs. 
Or, si nous notons $ t'':=(\ma)-|A\cap\alpha L|$, on sait par \ref{lemme-co-invariants}-(\ref{lemme-co-invariants-c})  que l'application $(\diamond)$ est bijective dès que
$\ma\geq t''+\mu_1\,,$
et comme $t''\leq \uell+t$ d'après \ref{theo-nouveau-t}, 
et que $t\leq|\mu|+|\ulambda|$ dans \ref{theo-diamond}, on voit qu'il suffit que 
$$m\geq a+\uell+|\mu|+|\ulambda|+\mu_1=|\mu|+\mu_1+2|\ulambda|
$$ 
(condition qui garantit aussi \ref{theo-dagger}) pour que $\varPsi_{m,t}$ soit bijective. L'inégalité \ref{theo-diamond} est donc bien satisfaite et ceci termine la preuve de la partie de (\ref{theo-Ind-lambda-d}) qui concerne $\rkms(\Ind_{\lambda}\V)$.

\medskip
\noindent (\ref{theo-Ind-lambda-d}) pour $\rkm(\Ind_{\lambda}\V)$, on suit la même démarche que pour $\rkms$. 

\smallskip
\noindent\relax{\slshape ($A$) Réduction au cas $\V$ est à poids unique $t$. }(\footnote{Cette étape n'est pas indispensable, mais elle simplifie quelque peu l'étape \slshape ($B$).})

\noindent Soit $p$ est le plus petit poids tel que $\V^{(p)}\not=0$. On a la suite exacte courte
$$0\to\V^{(>p)}\to\V\to\V^{(p)}\to0\,,$$
d'où la suite exacte d'extensions
$$
0\to\Ind_{\lambda}(\V^{(>p)})\to\Ind_{\lambda}(\V)\to\Ind_{\lambda}(\V^{(p)})\to0\,.
$$
Maintenant, si (\ref{theo-Ind-lambda-d}) satisfaite par les extrémités, on aura 
$$\mathalign{\rkm(\Ind_{\lambda}(\V))&\leq&\sup\set 
\rkm(\Ind_{\lambda}(\V^{(>p)})),\rkm(\Ind_{\lambda}(\V^{(p)}))/\hfill\\\noalign{\kern2pt}
&&\quad{}=\sup\set
\rkm(\V^{(>p)})+a,\rkm(\V^{(p)})+a/
=
\rkm(\V)+a\,.\hfill
}
$$
par \ref{coro-monotonie-stabilite}-(\ref{coro-monotonie-stabilite-a}) et d'après \ref{prop-scindage-par-le-poids}-(\ref{prop-scindage-par-le-poids-c}), et $\Ind_{\lambda}(\V)$ vérifie (\ref{theo-Ind-lambda-d}).

L'itération de cette idée montre que si
(\ref{theo-Ind-lambda-d}) est satisfaite par les $\FI$-modules à poids unique, elle sera satisfaite par les quotients $\V^{(\leq t)}$ pour tout $t\in\NN$, donc par $\V$, qu'elle soit de poids fini ou non.

\medskip
\noindent{\slshape ($B$) Le cas où $\V$ est à poids unique.}

\noindent Notons $\phi^{\Ig}_{m}$ le morphisme de transition  à l'ordre $m$ de $\Ind_{\lambda}(\V)$ et utilisons l'égalité
$\rkm=\rkme$ de \ref{rkms=rke}-(\ref{prop-rang-etendu-b}). On a donc
$$\rkm(\Ind_{\lambda}(\V))=\ind\bigset
s\in\NN\mid \varPhi_{t}(\phi^{\Ig}_{m})\text{ est injective }\forall m\geq s\,,\forall t\leq m/\,.$$

L'étude de l'injectivité de $\varPhi_{t}(\phi_{m}^{\Ig})$ suit les mêmes étapes $(B_i)$ que pour le cas $\rkms$, mais en plus simple puisque seul l'injectivité nous intéresse. De ce fait, on ne rencontre aucune contrainte sur $m$ en $(B_1)$ et $(B_2)$, mais si en $(B_3)$ où il faudra que $\varPhi_{t}(\phi_{m-a})$ soit injective, soit que
$m{-}a\geq\rkm(\V)$, grâce, une fois de plus, à l'égalité $\rkm=\rkme$ de \ref{rkms=rke}-(\ref{prop-rang-etendu-b}).\enddemo

\subsubsection{Le foncteur $\Ind_{\rho,\lambda}$}
Lors\label{Ind-alt-lambda} de l'étude de
rangs de monotonie et stabilité des termes des suites spectrales basiques de \ref{rang-stabilite-termes-basiques}, nous aurons besoin d'une version légèrement plus élaborée des foncteurs d'induction $\Ind_{\lambda}$.

\'Etant donnés,
$\lambda\in\Y_{\ell}(\ell+a)$ 
\emph{et une représentation $\rho:\HHH_{\ulambda}\to\Gl_{k}(\alt)$ d'un sous-groupe $\HHH_{\ulambda}\dans\S_{\ulambda}$}, on définit le foncteur 
$$\displayboxit{\Ind_{\rho,\lambda}:\ModFI k\fonct\ModFI k_{\geq |\ulambda|}}\eqno(\Ind_{\rho,\lambda})$$
qui fait correspondre à un $\S_{m-a}$-module $V_{m-a}$, le $\Sm$-module \comment
(\footnote{La notation $\Ind_{\HHH_{\ulambda}\times\myG_{\ulambda\bxx{m}}}^{\Sm}(\_)\boxtimes(\_)$ est,
conformément à la convention de \ref{operateurs-inductions},
un raccourci pour $\ind_{\HHH_{\ulambda}\times\S_{\ulambda\bxx{m}}}^{\Sm}(\_)\boxtimes(\_)$ où le sous-groupe
$\PPP_{\ulambda\bxx{m}}\dans\S_{\ulambda\bxx{m}}$ agit trivialement.})
\endcomment
$$W_{m}:=\begin{cases}\noalign{\kern-2pt}
0\,,\hbox{\quad si $m< |\ulambda|\,,$}\\
\ind_{\HHH_{\ulambda}\times\S_{m-|\ulambda|}}^{\Sm}\alt\otimes V_{m-a}\,,\hbox{\quad si $m\geq |\ulambda|$.}\\\noalign{\kern-3pt}
\end{cases}$$
Les analogues de la proposition \ref{prop-induction} et du théorème \ref{theo-Ind-lambda} sont encore vérifiés. En particulier, on contrôle toujours le rang de monotonie et stabilité.

\begin{theo}Soit\label{theo-Ind-rho-lambda}
$\lambda\in\Y_{\ell}(\ell+a)$ et $\rho:\HHH_{\ulambda}\to\Gl_{k}(\alt)$ une représentation d'un sous-groupe $\HHH_{\ulambda}\dans\S_{|\ulambda|}$. Alors, 
$\poids(\Ind_{\rho,\lambda}\V)\leq\poids(\V)+|\ulambda|$ et
$$
\displayboxit{
\rkms(\Ind_{\rho,\lambda}\V)\leq 
\rkms(\V)+2|\ulambda|\,,\quad
\rkm(\Ind_{\rho,\lambda}\V)\leq 
\rkm(\V)+a
}$$
\end{theo}
\def\Demonstration{Indication}
\demo On a une surjection naturelle de $\Sm$-modules
$$\ind_{\HHH_{\ulambda}\times\S_{m-|\ulambda|}}^{\Sm}\alt\boxtimes \big(\Res^{\S_{m-a}}_{\Sm-|\ulambda|}V_{m-a}\big)\onto
\ind_{\HHH_{\ulambda}\times\S_{m-|\ulambda|}}^{\Sm}\alt\otimes V_{m-a}$$
et le règles de Littlewood-Richardson (\cf note ($^{\ref{note-LR}}$), p.~\pageref{note-LR}) appliquées au $\Sm$-module de gauche confirment le majorant $\P(\V)+|\ulambda|$. La preuve de \ref{theo-Ind-lambda} s'applique alors telle quelle grâce à la réduction $B_1$ qui fait abstraction de la représentation $\alt$.
\enddemo
\subsubsection{Rangs des produits tensoriels d'un espace  gradué}\'Etant\label{FI-tensor} donné un $k$-espace vectoriel positivement gradué
$$\AAA:=\bigoplusnl_{i\in\NN}\AAA^{i}\in\Vec^{\NN}(k)\text{\quad avec\quad}\AAA^0=k\,,$$
on munit $\AAA^{\otimes m}:=\smash{\underparenthesis{\AAA\otimes_k\cdots\otimes_{k}\AAA}{\vrule depth6pt width0pt}_{m}}$ de l'action de $\alpha\in\Sm $
$$\alpha\cdot (a_1\otimes\cdots\otimes a_m):=
a_{\alpha^{-1}(1)}\otimes\cdots\otimes a_{\alpha^{-1}(m)}\,.
$$
L'application
$\phi_m:\AAA^{\otimes m}\to \AAA^{\otimes m+1}$,
$\omega\mapsto \omega\otimes \1_k$, 
où $\1_{k}\in\AAA^0$, est alors un morphisme de $\Sm $-modules positivement gradués, et 
$$\halfdisplayskips
\AAA^{\otimes}:=\set \phi_m:\AAA^{\otimes m}\to\AAA^{\otimes m+1}/_m\eqno(\A^\otimes)$$ est un  $\FI$-module.
\comment qui est de type fini (si et) seulement si $\dim\A=1$, car autrement la famille $\set\dim_k\AAA^{\otimes m}=(\dim_{k}\AAA)^{m}/_{m}$ n'est pas polynomiale (\ref{remas-caractere-poly}-(\ref{suite-de-dimensions})).
\endcomment
Pour chaque $i\in\NN$, les restrictions des $\phi_m$ aux composantes homogènes de degré $i$ donnent le sous-$\FI$-module de $(\AAA^{\otimes}){}^{i}\dans\AAA^{\otimes}$
$$\halfdisplayskips
(\AAA^{\otimes}){}^{i}:=\bigset\phi_m:(\AAA^{\otimes m})^{i}\to (\AAA^{\otimes m+1})^{i}/_{m}\,.\eqno(\AAA^{\otimes}){}^{i}$$

\begin{prop}[(\cite{cef})]Soit\label{AAA-tensor-stable} $\AAA$ un $k$-espace vectoriel positivement gradué avec $\AAA^{0}=k$. Alors, pour tout $i\in\NN$, on a 
$$
\poids((\AAA^{\otimes})^{i})\leq i\,,\quad
\rkm((\AAA^{\otimes})^{i})\leq i
\text{\ \ et\ \ }
\rkms((\AAA^{\otimes})^{i})\leq 2i\,.$$
\end{prop}
\demo
Pour tout $m\geq i$, on a un isomorphisme
$$
\halfdisplayskips
\def\underparenthesisnoalign{\kern0pt}
(\AAA^{\otimes m})^{i}=\bigoplus_{\lambda\in\Y_{i}(2i)}
\ind_{\myG_{\ulambda}\times\S_{m-\uell}}^{\Sm}
\underparenthesis{\A^{\lambda_1-1}\otimes\cdots\otimes\A^{\lambda_\uell-1}}_{\uell}
\boxtimes \underparenthesis{k\otimes\cdots\otimes k}_{m-\uell}\eqno(\ast)$$
Où $\myG_{\ulambda}\times\S_{m-\uell}$ est sous-groupe de $\S_{\uell}\times\S_{m-\uell}\dans\Sm$. Les règles de Littlewood-Richardson (\ref{regles-branchement}) montrent aussitôt que le poids de chaque induit est majoré par $\uell=\ell(\ulambda)\leq i$. Chaque terme de la somme $(\ast)$ n'est autre que le $\FI$-module 
$\M_{\smashbot{\uell}}^{\smashbot{\myG_{\ulambda}}}(\A^{\lambda_1-1}\otimes\cdots\otimes\A^{\lambda_\uell-1})$ de \ref{M-W} dont le rang de 
monotonie est majoré par $\uell$ et celui de stabilité par $2\uell$ d'après \ref{prop-M(W)}-(\ref{prop-M(W)-e}).
\enddemo

\def\Remarque{Commentaire}
\begin{rema}Si $W$ est une représentation de \relax{$\myG_{\ulambda}$}, la famille de foncteurs $\set\smashtop{\ind_{\myG_{\ulambda}\times\S_{m-\uell}}^{\Sm}}W\boxtimes(\_)/_{m}$ se recolle en un foncteur $\smash{\Ind_{\myG_{\ulambda}}W\boxtimes(\_)}$ défini sur $\ModFI k$. Il est alors possible d'adapter la démarche des théorèmes \ref{theo-Ind-lambda} et \ref{theo-Ind-rho-lambda} 
pour démontrer que l'on a 
$$\halfdisplayskips
\poids(\Ind_{\myG_{\ulambda}}W\boxtimes\V)\leq\poids(\V)+\uell
\text{\quad et\quad}
\pmathalign{\rkms(\Ind_{\myG_{\ulambda}}W\boxtimes\V)&\leq&\rkms(\V)+2\uell\,,\hfill\\
\rkm(\Ind_{\myG_{\ulambda}}W\boxtimes\V)&\leq&\rkm(\V)+\uell\,.\hfill}$$
Le cas $W:=\A^{\lambda_1-1}\otimes\cdots\otimes\A^{\lambda_\uell-1}$ et  $\V:=\V(0)$ aurait alors fournit une nouvelle preuve de \ref{AAA-tensor-stable}.
\end{rema}
\subsection{Les foncteurs d'induction $\Ibf^{a}$, $\Ibf(\agoth)$ et $\Thetabf^{a}$}Dans\label{foncteur-I-a} la section \ref{operateurs-inductions}, nous avons introduit le foncteur $\Ibf_{\ell}^{\ell+a}$ qui associe à une représentation de $\S_\ell$ une représentation de $\S_{\ell+a}$ en induisant suivant tous les diagrammes de Young $\lambda\in\Y_{\ell}(\ell+a)$:
$$\halfdisplayskips
\relax{\Ibf ^{\ell+a}_{\ell}:\Mod(k[\S_{\ell}])\fonct\Mod(k[\S_{\ell+a}])\,,\quad
\Ibf ^{\ell+a}_{\ell}:=\sumnl_{\lambda\in\Y_{\ell}(\ell+a)}\Ind^{\S _{\ell+a}}_{\Glambda }\,.
}\eqno(\Ibf )$$ 

Ce foncteur apparaît naturellement dans l'étude des espaces $\Delta_{\ell}\Xg^{\ell+a}$. Il intervient notamment dans la proposition \ref{coro-car-l-m}, où il donne l'égalité
$$\chic(\Delta_{\ell}\Xg^{\ell+a};i)=\Ibf_{\ell}^{\ell+a}(\chic(\Fl(\Xg);i))\,.$$

Nous allons maintenant recoller la famille des foncteurs $\set\Ibf_{m-a}^{m}/_{m}$ en un unique foncteur 
$\Ibf^{a}:\ModFI k\fonct\ModFI k$ dans le but d'étudier la stabilité (\cf\ref{strategie-stabilite}) des $k[\FB]$-modules 
$
\set\Hc^{i}(\Delta_{m-a}\Xg^{m})/_{m}$ et 
$\set\Hr^{i}(\Delta_{m-a}\Xg^{m})/_{m}\,.
$
\comment
Cela établira à priori, grâce au théorème \ref{FI-TF=>rep-stat}, que ces familles de représentations ont des caractères éventuellement polynomiaux. Bien évidemment, nous cherchons aussi à donner des bornes pour les rangs de monotonie et stabilité, ce qui présuppose que nous puissions contrôler les perturbations de rangs
que le foncteur $\Ibf^{a}$ provoque.
\endcomment

\medskip

Dans la section \ref{Ind-lambda}, nous avons associé à chaque $\lambda\in\Y_{\ell}(\ell+a)$ un foncteur 
$$\Ind_{\lambda}:\ModFI k\fonct\ModFI k_{\geq|\ulambda|}\,.$$
Sa définition est telle que
$\Ind_{\lambda}=\Ind_{\ulambda}$
où $|\ulambda|\leq 2a$\comment (\ref{rema-lambda/m}-(\ref{rema-lambda/m-b}))\endcomment, de sorte que la correspondance
$$\preskip-1.5ex\Big(\bigcup\nolimits_{1\leq i\leq a}\Y_{i}(i+a)\Big)\ni\lambda\mapsto \Ind_{\lambda}
\postskip4pt$$
est surjective sur l'ensemble de tels foncteurs. Mieux encore, le lemme \ref{a-2a} montre que cette correspondance devient  bijective pour peu que l'on tronque les foncteurs:
$\Y_{a}(2a)\ni\lambda
\gooddownstackrel{0pt}\longleftrightarrow\simeq
\Ind_{\lambda}(\_)_{\geq 2a}\,.$
On pose alors pour $a\in\NN$:\label{Ibfa}
$$\displayboxit{\Ibf ^{a}:=\sumnl_{\lambda\in\Y_{a}(2a)}\Ind_{\lambda}(\_)_{\geq 2a}:\ModFI k\fonct\ModFI k_{\geq 2a}
}\postdisplaypenalty10000
\eqno(\Ibf^{a} )$$ 
qui est un foncteur exact, fidèle sur $\ModFI k{}_{\geq a}$.

\subsubsectionnumber 
Nous pouvons maintenant facilement introduire l'analogue des inductions composées $\Ibf(\sigma)$ et $\Thetabf_{\ell}^{\ell+a}$ de \ref{inductions-iterees}.

\smallskip

Pour $a\in\NN$, notons $\Agoth(a)$ l'ensemble des suites d'entiers
$\agoth=(a_1,a_2,\ldots,a_r)$, avec $a_i>0$ et tels que $a=\sum_ia_i$.
On pose  pour chaque $\agoth\in\Agoth(a)$:
$$\Ibf(\agoth):=\Ibf^{a_1}\circ\Ibf^{a_2}\circ\cdots\circ\Ibf^{a_r}:
\ModFI k\to\ModFI k_{\geq 2a}\,.
$$
C'est un foncteur exact. La valeur de troncature de l'image, fixée à $2a$, est plus large que nécessaire, mais elle est ainsi indépendante de $\agoth$. De même, le noyau de $\Ig(\agoth)$ dépend beaucoup de la suite $\agoth$. La figure suivante illustre ces remarques, en grisé la zone de fidélité d'un foncteur composé pour $a=20$.
%
\begin{figure}[h!]
$$\hss\preskip0pt\includegraphics[scale=1]{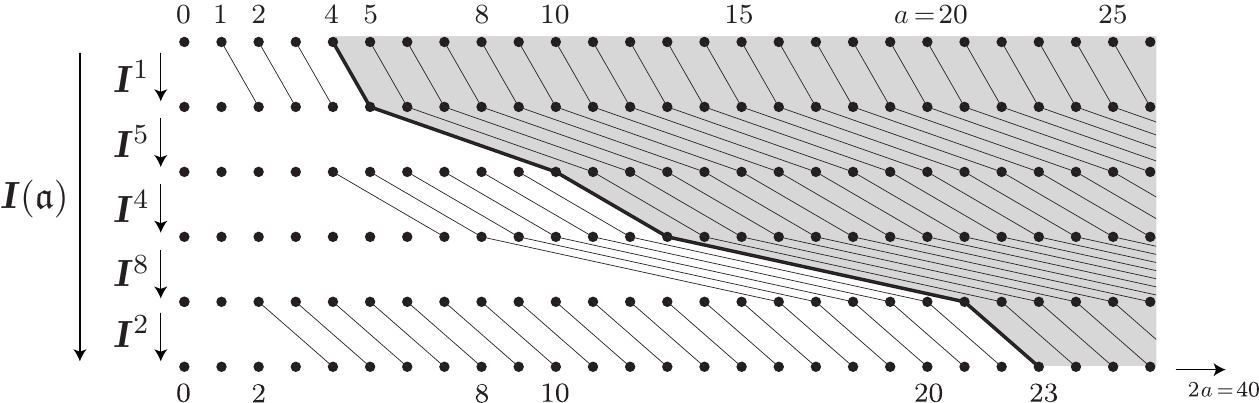}\hss$$
\mycaption{\small\slshape Foncteur d'inductions composées $\Ig(\agoth)=\Ig(2,8,4,5,1)$}
\end{figure}

\noindent On pose\label{Thetabfa} ensuite
$$\displayboxit{\mathrigid3mu
\Thetabf^{a}:=(-1)^{a}\sumnl_{\agoth\in\Agoth}(-1)^{|\agoth|}\Ibf(\agoth)
\,:\,\ModFI k\to K(\ModFI k_{\geq 2a})}\eqno(\Thetabf^{a})$$
où la présence de signes indique que l'on a affaire à des $\FI$-modules virtuels.

\medskip
\noindent Le corollaire suivant du théorème \ref{theo-Ind-lambda} est immédiat.\killline

\begin{theo}\let\Rg\Lg
Soient\label{theo-Ind-FI} $a\in\NN$ et $\agoth\in\Agoth(a)$ et désignons par $\Rg^{a}$ l'un des trois foncteurs $\Ibf^{a},\Ibf(\agoth),\Thetabf^{a}$.
\def\varlistskips{\topsep0pt\itemsep2pt}
\begin{enumerate}
\item Le\label{theo-Ind-FI-a} foncteur
$\mathrigid2mu
\Rg^{a}:\ModFI k\fonct\ModFI k_{\geq 2a}$
est covariant, additif, exact.
\item  Si\label{theo-Ind-FI-b} $\V$ est (de type fini) engendré en degrés $\leq d$, le $\FI$-module $\Rg^{a}(\V)$ est (de type fini) engendré en degrés $\leq \sup(d+a,2a)$.
\item \leavevmode\rlap{On a}\label{theo-Ind-FI-c}\hfill
$
\displayboxit{
\rkms(\Rg^{a}\V)\leq \rkms(\V)+4\,a\,,\quad
\rkm(\Rg^{a}\V)\leq \rkm(\V)+a
}$
\hfill\null
\end{enumerate}
\end{theo}

\section{Stabilité des familles $\set\Sm\rep\Hbm^{i}(\Delta_{?m-a}\Xg^{m})/_{m}$}
Dans cette section on\label{stabilite-des-familles} généralise le théorème de stabilité de Church (\cite{chu}) sur les familles 
$
\set\Fm(\Mg)/_{m}
$
où $\Mg$ est une variété différentielle orientée,
au cas des des familles
$\set\Delta_{?m-a}\Mg^{m}/_{m}$ où $\Mg$ est une \emph{pseudovariété} orientable, et où $a\in\NN$ est quelconque.

\subsection{Structure de $\FI$-module de $\set\Hbm^{i}(\Delta_{?m-a}\Mg^{m})/_{m}$}Dans\label{structure-FI-modules} les sections précédentes (\cf\ref{coh-BM}, \ref{renormalisation}), nous avons muni les familles $\Hbm^{i}(\Xg^{m})$ et $\Hbm^{i}(\Delta_{m-a}\Xg^{m})$ de structures de $\FI$-modules induites par les projections $p_m:\Xg^{m+1}\to\Xg^{m}$. On étend maintenant ces structures aux familles $\set\Hbm^{i}(\Delta_{?m-a}\Mg^{m})/_{m}$ pour toute pseudovariété orientée $\Mg$. Pour cela on considère le diagramme commutatif de décompositions ouvertes-fermées:
$$\halfdisplayskips
\def\tt{\vrule height 0pt depth6pt width0pt}
\def\uu{\vrule height 10pt depth0pt width0pt}
\def\o{{\rm ouvert}}
\def\f{{\rm ferm\acute e}}
\xymatrix@R=5mm{
\tt\Delta_{m+1-a}\Mg^{m+1}\arinto[r]_(0.48){\o}\arinto[d]^{j}&
\Delta_{\leq m+1-a}\Mg^{m+1}\ar@{<-^{)}}[r]_{\f}\ar@{<-^{)}}[d]^{i}&
\Delta_{\leq m-a}\Mg^{m+1}\ar@{<-^{)}}[d]^{i}\\
\Delta_{m-a}\Mg^{m}\times\Mg\arinto[r]_(0.48){\o}\ar[d]^{p_m}&
\uu\Delta_{\leq m-a}\Mg^{m}\times\Mg\ar@{<-^{)}}[r]_{\f}\ar[d]^{p_m}&
\uu\Delta_{\leq m-a-1}\Mg^{m}\times\Mg\ar[d]^{p_m}\\
\Delta_{m-a}\Mg^{m}\arinto[r]^(0.48){\o}&
\Delta_{\leq m-a}\Mg^{m}\ar@{<-^{)}}[r]^{\f}&
\Delta_{\leq m-a-1}\Mg^{m}
}$$
où $j$ est une inclusion ouverte et les inclusions $i$ sont fermées entre des espaces de même dimension cohomologique $(m{-}a)\,\dMg$ et mêmes orientations. On en déduit par les règles de fonctorialité de la cohomologie de Borel-Moore, le diagramme commutatif de suites exactes longues:
$$\halfdisplayskips
\def\tt{\vrule height 0pt depth0pt width0pt}
\def\uu{\vrule height 0pt depth0pt width0pt}
\def\o{{\rm ouvert}}
\def\f{{\rm ferm\acute e}}
\def\arleft{\ar@{<-}[r]}
\def\HH#1/{\Hbm^{*}(#1)}
\def\HHm#1/{\Hbm^{*-\dMg}(#1)}
\mathrigid0mu
\hss
\xymatrix@R=5mm@C=2mm{
\tt
\HH\Delta_{m+1-a}\Mg^{m+1}/
\arleft\ar@{<-}[d]^{(j_{!})\dual}&
\HH\Delta_{\leq m+1-a}\Mg^{m+1}/
\arleft\ar@{<-}[d]^{(i^{*})\dual}&
\HHm\Delta_{\leq m-a}\Mg^{m+1}/
\ar@{<-}[d]^{(i^{*})\dual}\arleft&\\
\HH\Delta_{m-a}\Mg^{m}\times\Mg/
\arleft\ar@{<-}[d]^{p_m^{*}}&
\uu
\HH\Delta_{\leq m-a}\Mg^{m}\times\Mg/
\arleft\ar@{<-}[d]^{p_m^{*}}&
\uu
\HHm\Delta_{\leq m-a-1}\Mg^{m}\times\Mg/
\ar@{<-}[d]^{p_m^{*}}\arleft&\\
\HH\Delta_{m-a}\Mg^{m}/\arleft&
\HH\Delta_{\leq m-a}\Mg^{m}/\arleft&
\HHm\Delta_{\leq m-a-1}\Mg^{m}/\arleft&
}\hss$$
où les composées des flèches verticales sont compatibles à l'action des groupes symétriques et vérifient la condition pour définir des $\FI$-modules. 

Ces observations constituent l'essentiel de la preuve des assertions (\ref{prop-structure-FI-modules-a},\ref{prop-structure-FI-modules-b}) de proposition suivante.\killline

\begin{prop}[et définitions]Soit\label{prop-structure-FI-modules} $\Mg$ une pseudovariété orientée de dimension $\dMg$, et soient $0<a\leq m$.
\def\varlistskips{\topsep2pt\itemsep2pt}
\begin{enumerate}
\item Dans\label{prop-structure-FI-modules-a} le diagramme précédent,
les morphismes verticaux définissent les $\FI$-modules 
$\set\Hbm^{i}(\Delta_{?m-a}\Mg^{m})/_{m}$ et les horizontaux 
définissent alors la suite exacte longue de $\FI$-modules
$$\hss\mathrigid0mu
\to
\bigset\Hbm^{i-\dMg}(\Delta_{\leq m-a-1}\Mg^{m})/
\to
\bigset\Hbm^{i}(\Delta_{\leq m-a}\Mg^{m})/
\to
\bigset\Hbm^{i}(\Delta_{m-a}\Mg^{m})/
\to
\hss$$

\item Les\label{prop-structure-FI-modules-b} complexes fondamentaux en cohomologie de Borel-Moore de la famille $\set \Delta_{\leq m-a}\Mg^{m}/_{m}$, obtenus en dualisant les complexes \ref{complexe-fondamental}, 
 s'organisent naturellement en \expression{complexe fondamental de $\FI$-modules} 
$$\def\ell{m-a}
\medmuskip=0mu
0\to
\set\Hbm ^{i}(\Delta_{\leq\ell}\Mg^{m})/_m
\to
\Hr(\ell)
\to
\Hr(\ellmo)
\to\cdots\to
\Hr(1)
\to
0
$$
où $\Hr(m{-}a{-}b):=\bigset\Hbm ^{i-b\,(\dMg{-}1)}(\Delta_{m-a-b}\Mg^{m})/_{m}$ et $i\geq b(\dMg{-}1)$.

\noindent Lorsque $\Mg$ est $i$-acyclique, le complexe fondamental est exact et l'on a l'égalité suivante dans $K(\ModFI k)$,
$$\bigset\Hbm^{i}(\Delta_{\leq m-a}\Mg^{m})/=
\hskip-1em\sum_{0\leq b<m-a}\hskip-1em
(-1)^{b}\big(\bigset\Hbm^{i-b\,(d_\Mg-1)}(\Delta_{m-a-b}\Mg^{m})/_{m-a-b}\big)\,.
$$

\item Le\label{prop-structure-FI-modules-c} foncteur 
$\Ibf^{a}:\ModFI k\fonct\ModFI k_{\geq 2a}$ identifie les $\FI$-modules
$$\mathalign{
\hfill\bigset\Hbm^{i}(\Delta_{m-a}\Mg^{m})/_{m\geq2a}&=&
\Ibf ^{a}\big(\bigset\Hbm^{i}(\Fg_{m-a}(\Mg))/_{m-a}\big)\,.}
$$

\penalty-5000
\item Si\label{prop-structure-FI-modules-d} $\Xg$ est $i$-acyclique, les foncteurs $\Thetabf^{a},\Ibf^{a}:\ModFI k\fonct\ModFI k_{\geq 2a}$ donnent les identifications suivantes  dans $K(\ModFI k)$,
$$\def\sep{\noalign{\kern4pt}}\def\ssep{\noalign{\kern4pt}}
\dimen10=1.2cm\dimen11=5cm
\mathalign{\noalign{\kern6pt}
\hfill\bigset\Hbm^{i}(\Fm(\Xg))/_{m\geq2\tilde \imath}&=&\sumnl_{0\leq a<m}\Thetabf^{a}\big(\bigset\Hbm^{i-a\,(d_{\Xg}{-}1)}(\Xg^{m-a})/_{m-a}\big)\hfill\\\ssep
\bigset\Hbm^{i}(\Delta_{m-a}\Xg^{m})/_{m\geq2a+2\tilde \imath}&=&\\\sep
&&\hskip-\dimen11\hskip2.5\dimen10
\mathrigid0mu
\mllap{=\ }\Ibf ^{a}\Big(
\sumnl_{0\leq b<m-a}
\hskip-0.5ex
\Thetabf^{b}\big(\bigset\Hbm^{i-b\,(d_{\Xg}{-}1)}(\Xg^{m-(a+b)}/_{m-(a+b)}\big)\Big)\hfill\\\ssep
\bigset\Hbm^{i}(\Delta_{\leq m-a}\Xg^{m})/_{m\geq2a+2\tilde \imath}&=&\\\sep
&&\hskip-\dimen11\hskip\dimen10
\mathrigid0mu
\mllap{=}\sum_{0\leq b<m-a}\hskip-0.5ex(-1)^{b}\ 
\Ibf ^{a+b}
\Big(\hfill\\\noalign{\kern-8pt}
&&\hskip-\dimen11\hskip2\dimen10
\sum_{0\leq c<m-a-b}
\hskip-1em
\Thetabf^{c}\big(\bigset\Hbm^{i-(b+c)(d_{\Xg}{-}1)}(\Xg^{m-(a+b+c)}/_{m-(a+b+c)}\big)\ \ \Big)\hfill\\\noalign{\kern4pt}
}$$
avec $\tilde \imath=\lfloor i/(d_{\Xg}{-}1)\rfloor$.
 
\end{enumerate}\end{prop}
\def\Demonstration{Indications}\demo (\ref{prop-structure-FI-modules-c})
résulte
de \ref{coro-car-l-m} et de la définition de $\Ibf^{a}$ (\ref{foncteur-I-a}), et (\ref{prop-structure-FI-modules-d}) résulte de
la définition de de $\Thetabf^{a}$ (\loccit) et des égalités
 \ref{theo-caracteres-devissage}-(\ref{theo-caracteres-devissage-a})-(i,ii,iii), le tout modulo les équivalences $\Hc^{i}(\Delta_{\leq r}\Xg^{s})\dual=
\Hbm^{r\,d_{\Xg}-i}(\Delta_{\leq r}\Xg^{s}) $.
\enddemo

\subsection{Les familles de représentations $\set\Hbm^{i}(\Delta_{?m-a}\Xg^{m})/_{m}$}


Nous\label{stabilite-des-familles-i-acycliques} supposons dans un premier temps que $\Xg$ est une pseudovariété orientable \emph{$i$-acyclique}. Plus tard, dans \ref{stabilite-des-familles-pseudovarietes}, nous nous affranchirons de l'hypothèse de $i$-acyclicité.

\subsubsectionline{Le cas de la stabilité de la cohomologie à support compact.}Ce cas est inintéressant si $\Xg$ est $i$-acyclique puisque, d'après \ref{=l-m-poly-univ} et \ref{<=l-m-poly-univ}, les polynômes de Poincaré pour $\Hc(\Delta_{?m-a}\Xg^{m})$ sont de valuation $m-a$, ce qui implique que
$$\Hc^{i}(\Delta_{?m-a}\Xg^{m})=0\,,\quad\forall m>a+i\,. $$

\subsubsectionline{Le cas de la stabilité de la cohomologie de Borel-Moore.}\label{rappel-coh-BM}


Lorsque\label{rem-connexite} $\Mg$ est de type fini, la condition de stabilité de représentations impose sur $\Mg$ qu'elle soit connexe et que $\dim_{\RR}\Mg\geq2$, car, autrement, les $\FI$-modules $\Hbm^{i}(\Delta_{m-a}\Mg^{m})$ ne sont pas de type fini.
En effet, si $r:=|\pi_0(\Mg)|$, on a $|\pi_0(\Mg^{m})|=r^{m}$ et $|\pi_0(\Fm(\Mg))|=(r^{m})\usp{m}$ ce qui exclut toute possibilité de finitude pour le $\FI$-module $\set\Hbm^{i}(\Delta_{m-a}\Mg^{m})/_{m}$ pour peu que $r>1$. D'autre part, si $\Mg$ est connexe et $\dim_{\RR}(\Mg)=1$, c'est une courbe illimitée avec un nombre fini $f$ de points multiples et alors $|\pi_{0}(\Fm(\Mg))|\geq(m-f)!$ d'où encore une obstruction à la finitude. 

{\sl On suppose que  $\Xg$ est une pseudovariété $i$-acyclique connexe orientable.}

\medskip
\noindent{$\bullet$ {\bf(A)} \sl Le cas de $\set \Delta_{\leq m}(\Xg^{m})/_{m}$.--- }On s'intéresse au $\FI$-module 
$$\A(\Xg;i):=\set\Hbm^{i}(\Xg^{m})/_{m}\,.$$
La proposition \ref{AAA-tensor-stable} s'applique, car $\Hbm^{0}(\Xg;k)=k$ (\cf\ref{comm-coh-BM}-(\ref{comm-coh-BM-Hbm0})), et donne
$$\rkms(\A(\Xg;i))\leq 2i\,, \text{\quad et\quad }\rkm(\A(\Xg;i))\leq i\,.\eqno(\rm A)$$

\medskip
\noindent{$\bullet$ {\bf(B)} \sl Le cas de $\set \Delta_{m}(\Xg^{m})/_{m}$.--- } (\footnote{Le cas particulier où
$\Xg$ est une variété topologique est
étudié par Church dans \cite{chu}, (\cf\ref{stabilite-caracteres}); il a suscité ce chapitre de généralisations.})
On s'intéresse au $\FI$-module 
$$\B(\Xg;i):=\set\Hbm^{i}(\Fm(\Xg)/_{m}\,.$$
L'égalité \ref{prop-structure-FI-modules}-(\ref{prop-structure-FI-modules-d}) nous conduit à chercher un majorant pour les rangs de stabilité des $\FI$-modules
$$\B_a(\Xg;i):=\Thetabf^{a}\big(\set\Hbm^{i-a\,(d_{\Xg}{-}1)}(\Xg^{m-a})/_{m-a}\big)\eqno(\ast)$$
pour $0\leq a<m$ vérifiant
$$a(d_{\Xg}{-}1)\leq i\,.\eqno(\ddagger)$$
On sait d'après {\bf(A)} que le
$\FI$-module $
\set\Hbm^{i-a(d_{\Xg}{-}1)}(\Xg^{m-a})/_{m-a}$
est stable pour
$m{-}a\geq 2\big(i{-}a\,(d_{\Xg}{-}1)\big)$ et le théorème \ref{theo-Ind-FI}-(\ref{theo-Ind-FI-c}) permet d'évaluer les rangs de monotonie et stabilité de $\B_a(\Xg;i)$. On aboutit alors aux majorations
$$\preskip1.5em
\pmathalign{
\rkms(\B_a(\Xg;i)_{\geq 2\tilde\imath})
&\ \leq\ &
\smash{2\big(i-a(d_{\Xg}{-}1)\big)+4a}
\ \leq\ 
\smashtop{\relspace1
\left\{\vrule height5mm width0pt\right.
\!\mathalign{
4i,&\text{ si $d_{\Xg}=2$,}\\\noalign{\kern2pt}
2i,&\text{ si $d_{\Xg}\geq3$,}
}}\\\noalign{\kern4pt}
\rkm(\B_a(\Xg;i)_{\geq 2\tilde\imath})
&\ \leq\ &
\big(i-a(d_{\Xg}{-}1)\big)+a=i-a(d_{\Xg}{-}2)
\ \leq\ i\,,\hfill}
\postdisplaypenalty10000
\eqno(\rm B)$$
dont la première résulte d'expliciter $2\big(i+(3-d_{\Xg})a\big)$
suivant que $d_{\Xg}=2$ ou $d_{\Xg}\geq3$ et
moyennant la majoration ($\ddagger$).

Il faut encore remarquer que le foncteurs $\Thetabf^{a}$ considérés sont à valeurs dans $\ModFI(k)_{\geq2a}$ et que dans tous les cas on a $2a\leq2i$ de sorte que les majorations (B) sont aussi valables pour $\B(\Xg;i)$.

\medskip
\noindent{$\bullet$ {\bf(C)} \sl Le\label{Cas-C} cas de $\set \Delta_{m-a}(\Xg)/_{m}$.--- }On s'intéresse au $\FI$-module
$$
\C_a(\Xg;i):=\set \Hbm^{i}(\Delta_{m-a}\Xg^{m})/_{m}\,.
$$
On raisonne exactement comme dans (B) à l'aide de
la formule \ref{prop-structure-FI-modules}-(\ref{prop-structure-FI-modules-d}), 
ce qui nous conduit à calculer les rangs de stabilité des $\FI$-modules
$$\Ibf ^{a}\Big(
\sumnl_{0\leq b<m-a}
\hskip-0.5ex
\Thetabf^{b}\big(\bigset\Hbm^{i-b\,(d_{\Xg}{-}1)}(\Xg^{m-(a+b)}/_{m-(a+b)}\big)\Big)\,.$$
On est conduit à trouver le borne supérieures
$$\preskip1em
\pmathalign{
\rkms &\leftrightarrow& \relax{\sup\nolimits_{0\leq b}\bigset2\big(i-b\,(d_{\Xg}{-}1)\big)+4b+4a/}\hfill\\\noalign{\kern4pt}
\rkm &\leftrightarrow& \relax{\sup\nolimits_{0\leq b}\bigset\big(i-b\,(d_{\Xg}{-}1)\big)+b+a/}\hfill}$$
sachant que $b\,(d_{\Xg}{-}1)\leq i$. On trouve les majorations
$$\preskip1.2em
\pmathalign{
\rkms(\C_a(\Xg;i)_{\geq 2a+2\tilde\imath})
&\leq&
\smashtop{\left\{\vrule height5mm width0pt\right.
\!\!\mathalign{
4i+4a,&\text{ si $d_{\Xg}=2$,}\\\noalign{\kern2pt}
2i+4a,&\text{ si $d_{\Xg}\geq3$,}\\
}}\\\noalign{\kern4pt}
\rkm(\C_a(\Xg;i)_{\geq 2a+2\tilde\imath})
&\leq&i+a\,,\hfill
}
\postdisplaypenalty10000
\eqno(\rm C)$$
qui sont également valables pour $\C_{a}(\Xg;i)$ puisque
 $\tilde \imath=\lfloor i/(d_{\Xg}{-}1)\rfloor\leq i$.

\medskip
\noindent{$\bullet$ {\bf(D)} \sl Le\label{Cas-D} cas de $\set \Delta_{\leq m-a}(\Xg)/_{m}$.--- }On s'intéresse au $\FI$-module (\ref{structure-FI-modules})
$$
\D_a(\Xg;i):=\set \Hbm ^{i}(\Delta_{\leq m-a}\Xg^{m})/_{m}
$$
On trouve encore les mêmes majorations des rangs que dans (C). Pour s'en convaincre, on pourrait raisonner à l'aide de
la formule \ref{prop-structure-FI-modules}-(\ref{prop-structure-FI-modules-d}), mais il s'avère plus intéressant d'appliquer {({\bf C})} au complexe fondamental de $\FI$-modules 
pour la cohomologie de Borel-Moore (\ref{prop-structure-FI-modules}-(\ref{prop-structure-FI-modules-b})), 
$$\def\ell{m-a}
\medmuskip=0mu
0\to
\set\Hbm ^{i}(\Delta_{\leq\ell}\Xg^{m})/_m
\to
\Hr(\ell)
\to
\Hr(\ellmo)
\to\cdots\to
\Hr(1)
\to
0
$$
avec $\Hr(m{-}a{-}\ell):=\set\Hbm ^{i-(d_{\Xg}{-}1)\ell}(\Delta_{m-a-\ell}\Xg^{m})/_{m}$ et aussi $i\geq (d_{\Xg}{-}1)\ell$, et
qui est exact dans le cas présent puisque $\Xg$ est $i$-acyclique.

Le $\FI$-module $\set\Hbm ^{i}(\Delta_{\leq\ell}\Xg^{m})/_m$ apparaît ainsi comme somme alternée des $\FI$-modules $\Hr(m{-}a{-}\ell)$  dans le groupe de Grothendieck $K_0(\ModFI k)$, et hérite donc des mêmes rangs de monotonie et stabilité d'après
\ref{coro-monotonie-stabilite}-(\ref{coro-monotonie-stabilite-bii}).

\medskip
\noindent Résumons les conclusions de ces observations sous forme de théorème.\vskip-1.5ex\vskip0pt

\begin{theo}Soit\label{theo-stabilite-BM-i-acyclique} $\Xg$ une pseudovariété $i$-acyclique connexe orientable de dimension $d_{\Xg}\geq2$. Pour  $a,i\in\NN$,
le $\FI$-module
$\set\Sm \rep\Hbm ^{i}(\Delta_{?m-a}\Xg^{m})/_{m}$
est monotone pour $m\geq i+a$ et est monotone
et stable pour $m\geq4i+4a$, si $d_{\Xg}=2$, et pour $m\geq2i+4a$, si $d_{\Xg}\geq3$. Les familles des caractères et des nombres de Betti correspondantes sont (donc) polynomiales sur les mêmes intervalles d'entiers $m$.
\end{theo}

\begin{rema}Dans les arguments qui précèdent, $\set\Hbm ^{i}(\Delta_{?m-a}\Xg^{m})/_{m}$ a une structure à priori de $\FI$-module, celle donnée par les morphismes image-inverse (\ref{structure-FI-modules}). Lorsque nous utilisons le théorème \ref{prop-structure-FI-modules}-(\ref{prop-structure-FI-modules-d}) nous passons de la catégorie des $\FI$-modules à son groupe de Grothendieck $K_0(\ModFI k)$ mais la conclusion se fait de nouveau dans $\ModFI k$. Par exemple, le $\FI$-module $\set\Hbm ^{i}(\Delta_{\leq m-a}\Xg^{m})/_{m}$ $({\dagger})$ se retrouve identifié
à une somme alternée de $\FI$-modules ayant les bons rangs de monotonie et stabilité, mais cette somme alternée provient en fait d'une résolution dans $\ModFI k_{\geq ?}$ par des $\FI$-modules à la fois monotones et stables. Le fait que la sous-catégorie pleine 
\smashbot{$\ModFI k_{\geq?}$} de ce type de $\FI$-modules en soit une sous-catégorie abélienne permet de conclure que \smash{$({\dagger})$} est également monotone et stable pour $m\geq\,?$. C'est d'ailleurs ce que dit la proposition \ref{coro-monotonie-stabilite}-(\ref{coro-monotonie-stabilite-bii}).
\end{rema}

\subsection{Stabilité des familles de représentations $\set\Hbm^{i}(\Delta_{?m-a}\Mg^{m})/_{m}$}\label{stabilite-des-familles-pseudovarietes}

{\sl On rappelle que $\Mg$ désigne une pseudovariété orientée de type fini.}

Nous démontrerons l'analogue du dernier théorème \ref{theo-stabilite-BM-i-acyclique} pour la famille de représentations $\set\Sm\rep\Hbm^{i}(\Fm(\Mg))/_{m}$ où $\Mg$ est une pseudovariété orientée générale, \idest qu'elle soit $i$-acyclique ou non.

\subsubsectionline{Stabilité dans les suites spectrales basiques.}
Le\label{stabilite-termes-basiques} théorème \ref{theo-suite-spectrale-basique-relative} établit que les $\FI$-modules définis par les morphismes
$$\IEs(q_m^{*})_{1}^{p,q}:\IEs(\U^{m})^{p,q}_{1}\to\IEs(\U^{\mm})^{\pp,q}_{1}\,,$$
sont canoniquement isomorphes aux $\FI$-modules définis par les morphismes
$$
\ind^{\Sm }_{\vrule height4mm depth5pt width0pt
\hskip-4mm
\StS{m{-}(\pp)}{\pp}}
\hskip-10mm
\alt\otimes\Hbm^{q}(\UU^{m}_{m-p,\ldots,m})
\hf{\ind p_{m}^{*}}{}{1.3cm}
\ind^{\Smm }_{\vrule height4mm depth5pt width0pt
\hskip-4mm
\StS{m{-}(\pp)}{p+2}}
\hskip-10mm
\alt\otimes\Hbm^{q}(\UU^{\mm}_{m-p,\ldots,m+1})
$$
où $\alt$ indique que l'action de $\S_{m{-}(\pp)}\times\1$ est tordue $\alt_{m{-}(\pp)}$ et lorsque le couple $(p,q)\in\NN^{2}$ est soumis à la contrainte $q=i+(m{-}(\pp))$.

Le théorème (\loccit) décrit aussi, dans son assertion (\ref{theo-suite-spectrale-basique-relative-c}), le diagramme commutatif d'\emph{espaces vectoriels}
$$\preskip0pt\def\ind#1#2{_{\hbox to#1{\hss$\scriptstyle#2$\hss}}}
\def\mbx#1#2{\hbox to1.2cm{$#2$\hss}}
\xymatrixc{@R6mm@C1.5cm}
{
{\Hbm^{q}(\UU^{m}_{m-p,\ldots,m})}
\ar[r]^(0.45){p_m^{*}}
\ar[d]_(0.48){\varXi_{\pp}^{m}\ }|(0.45){\xysimeq}
&
{\Hbm^{q}(\UU^{\mm}_{m-p,\ldots,m+1})}
\\
\displaystyle\bigoplus\ind{1cm}{f\in\FFF(\pp,m)}
\Hbm^{Q}(\Fg_{\pgoth(f)}(\Mg_{>0}))
\ar[r]^{\bigoplus p_{\pp}^{*}}&
\displaystyle\bigoplus\ind{1cm}{f^{\bullet}\in\FFF^{\bullet}(p{+}2,\mm)}
\Hbm^{Q}(\Fg_{\pgoth(f^{\bullet})}(\Mg_{>0}))
\ar@{^(->}[u]+<0pt,-9pt>_(0.6){(\varXi_{p+2}^{\mm})^{-1}}
}
$$
où les isomorphismes $\Xi^{m}_{\pp}$ pourraient nous inciter à transporter l'action des $\StS{m{-}(\pp)}{\pp}$ pour en faire un diagramme commutatif de \emph{représentations}, quelque chose dont nous n'avons pas encore eu besoin. Or, un tel transfert formel de structure sur des sommes directes dont la définition est basée sur un ordre précis des $m-(\pp)$ premières coordonnées ne rend pas transparent le fait que l'on est en présence de sommes directes de $\FI$-modules. C'est pour cette raison que nous introduisons dans les sections qui suivent une nouvelle décomposition de $\Hbm(\UU^{m}_{m-p,\ldots,m})$ qui respecte par construction les symétries et pour laquelle les analogues du diagramme précédent seront à priori compatibles à l'action de $\StS{m{-}(\pp)}{\pp}$.

\subsubsectionline{Tableaux et fonctions de $\FFF(\pp,m)$.}
Nous\label{tableaux<->FFF} encodons une fonction $f\in\FFF(\pp,m)$ sous la forme d'un tableau à $\pp$ lignes et $m$ boites.

L'idée est la suivante. Soit $a>0$. Pour toute suite $\cl x=(x_1,\ldots,x_a)$ d'éléments deux à deux distincts de $\iii[1,a]$, associons la fonction $f_{\cl x}\in\FFF(1,a+1)$ définie pour $i\leq a$,
par la règle suivante

{\sl Pour $i\in\iii[1,a]$, notons $I_{\cl x}(i):=\set j\leq i\mid x_j>x_i/$. Alors
$$f_{\cl x}(x_i)=\begin{cases}\noalign{\kern-4pt}
a+1\,,\text{ si $I_{\cl x}(i)=\emptyset$}\\
x_{\sup I_{\cl x}(i)},\text{ autrement.}
\end{cases}\postskip-0.5ex$$
}
Par exemple:
$$\preskip2pt\binspace0
\mathalign{
\cl x=(1,2,\ldots,a-1,a)&\mapsto&f_{\cl x}(i)=a+1\hfill\\
\cl x=(a,a-1,\ldots,2,1)&\mapsto&f_{\cl x}(i)=\sup\set i+1,a+1/\\
}\postskip2pt$$

Cette correspondance est injective. En effet, si $\cl x\not=\cl y$, il existe un premier $i$ tel que $x_i\not=y_i$. Supposons que l'on ait $x_i< y_i$. Alors, comme il existe $h>i$ avec $y_h=x_i$, on aura 
$f_{\cl y}(x_i)=f_{\cl y}(y_h)=y_{k}\not=a+1$ pour un certain $k\geq i$, tandis que $f_{\cl x}(x_i)=x_{k}=y_k$ pour un certain $k<i$, ou bien $f_{\cl x}(x_i)=a+1$. dans tous les cas, $f_{\cl x}\not= f_{\cl y}$.
$$\halfdisplayskips
\includegraphics[scale=0.9]{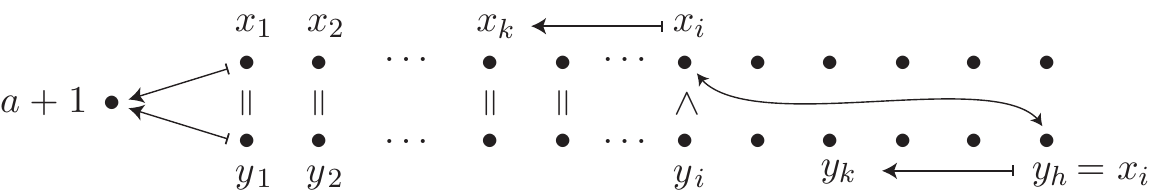}$$
Comme d'autre part $a!=|\FFF(1,a+1)|$, on conclut que la correspondance en question est bijective.
Cette idée est la base de la preuve du lemme suivant.

\begin{lemm}[et notations]L'ensemble\label{lemme-TTT-FFF} $\FFF(\pp,m)$ est en bijection avec l'ensemble 
$\TTT(\pp,m)$\glossary{$\TTT(\pp,m)$:ensemble de tableau paramétrant $\FFF(\pp,m)$} des tableaux à $(\pp)$ lignes et $m$ boites de formes:
$$
\yght=0.42cm\ygwd=1.3cm
\let\ygvdots\Ygvdots
\def\everyygcell{\footnotesize}
\tau:=\left\{\tableauc{}{
m{-}p&\bullet&\bullet\\
m{-}\pp&\bullet\\
m{-}p{+}2&\bullet&\bullet&\bullet&\bullet\\
\vdots&\noborder\rlap{\hskip-1.5ex
$(j)$ \raise5pt\hbox{\rotatebox{180}{$\hfto{}{}{2.cm}$}} $(i)$}\\
m{-}2\\
m{-}1&\bullet\\
m&\bullet&\bullet&\bullet}\right.$$
où les boites marquées par {\rm`$\bullet$'} sont remplies par les entiers de $\iii[1,m{-}(\pp)]$ de manière biunivoque.
La correspondance associe au tableau $\tau$ la fonction $f_{\tau}$
qui, appliquée à $i\in\iii[1,m{-}(\pp)]$, donne le premier entier $j>i$ à  gauche de $i$ sur la ligne de $\tau$ où il figure.

L'ensemble $\FFF^{\bullet}(p{+}2,m{+}1)$ est en bijection avec l'ensemble des tableaux $\TTT^{\bullet}(p{+}2,m{+}1)$\glossary{$\TTT^{\bullet}(\pp,m)$:ensemble de tableau paramétrant $\FFF^{\bullet}(\pp,m)$} dont la ligne $(m{+}1)$ est de longueur $1$.
L'application\glossary{${(\_)^{\bullet}:\TTT(\pp,m)\to\TTT^{\bullet}(p{+}2,m{+}1)}$:application qui rajoute une boite isolée}
$$(\_)^{\bullet}:\TTT(\pp,m)\to\TTT^{\bullet}(p{+}2,m{+}1)$$
qui associe à $\tau$ le tableau $\tau^{\bullet}$ obtenu en rajoutant la boite $(m{+}1)$ est bijective.
\end{lemm}

\subsubsection{Décomposition symétrique de $\Hc(\Delta_{a+b}(\Zg^{a}\times\Xg^{b}))$}\label{decomposition-symétrique-Delta-Za-Xb}
\noindent Avant d'aborder le cas spécifique de $U^{m}_{m-p,\ldots,m}=\Delta_m(\Mg_{\geq0}^{m-(\pp)}\times\Mg_{>0}^{\pp})$,
on se place dans un contexte plus général avec des notations plus légères.

\ssubsubsection{Données et notations}\label{donnees-notas-dec-sym}
\begin{list}{\theenumi)}{\def\theenumi{N-\arabic{enumi}}\usecounter{enumi}\leftmargin2em
\labelsep1ex\labelwidth3ex}
\item Soit $\Xg\dans\Zg$ une inclusion d'espaces $i$-acycliques et 
$\Hc(\Zg)=0\,.$
\item Pour $\Ig\dans\iii[1,a]$, on note $\Zg^{\Ig}$ le produit $\Zg^{|\Ig|}$ où les coordonnées sont indexées par $\Ig$. On notera
$\Fg_{\Ig}:=\Delta_{|\Ig|}(\Zg^{\Ig})$ et 
 $\Fb:=\Fb(\Xg)$. Les notations $\Delta_{?\ell}(\Fg_{\Ig}\times\Fb)$ et
 $\Delta_{?\ell}(\Zg^{\Ig}\times\Fb)$ ont alors le sens habituel de \ref{notas}-(\ref{nota-Delta-}). 
\item Soit $\TTT_{b}(\Ig)$\glossary{$\TTT_{b}(\Ig)$:ensemble des tableaux à $(|\Ig|{+}b)$ boites dont la première colonne est $\iii[1,b]$ et dont $\Ig$ remplit les autres boites}
l'ensemble des tableaux à $(|\Ig|{+}b)$ boites dont la première colonne est $(1,\ldots,b)$ et dont la forme est
$$
\yght=0.35cm\ygwd=0.7cm
\let\ygvdots\Ygvdots
\tau:=\left\{\tableauc{}{
1&\bullet&\cellcolor{0 0 0 0.0}\bullet\\
\cellcolor{0 0 0 0}2&\cellcolor{0 0 0 0.0}\bullet\\
\cellcolor{0 0 0 0}3&\bullet&\bullet&\bullet&\bullet&\bullet&\cellcolor{0 0 0 0.0}\bullet\\
\vdots\\
\cellcolor{0 0 0 0}b{-}2&\cellcolor{0 0 0 0.0}\bullet\\
\cellcolor{0 0 0 0}b{-}1\\
\cellcolor{0 0 0 0}b&\bullet&\bullet&\cellcolor{0 0 0 0.0}\bullet}\right.$$
où les boites {\rm`$\bullet$'} sont remplies avec tous les éléments de $\Ig$. Le groupe
$\Sg_{\Ig}\times\Sb$ agit sur $\TTT_{b}(\Ig)$ par son action sur le contenu des tableaux. Par exemple, si $\tau:=\tableauc{}{1&x&y\\
2&z}$, on a $(1,2)(x,y,z)\cdot\tau=\tableauc{}{2&y&z\\
1&x}=\tableauc{}{1&x\\2&y&z}$

\item Pour\label{tau-goth} $\tau\in\TTT_{b}(\Ig)$, on notera $\pgoth(\tau)$\glossary{$\pgoth(\tau)$:partition de $\iii[1,b]$ déterminée par le tableau $\tau\in\TTT_{b}(\Ig)$} la partition de $\Ig\sqcup\iii[1,b]$ déterminée 
par les lignes de $\tau$. Par exemple, 
$$
\tau:=\tableauc{}{1&x&y&z\\
2\\
3&a}\quad\Longrightarrow\quad
\pgoth(\tau)=\bigset{\set1, x,y,z/},{\set 2/},{\set3, a/}/\,.$$
Alors, conformément à la notation de \ref{nota-pgoth}, on pose 
$$\Fg_{\tau}:=\Fg_{\pgoth(\tau)}\dans\Delta_{b}(\Zg^{\Ig}\times\Fb)\,.$$
La projection $\pi_{b}:\Fg_{\tau}\to\Fb$ est clairement un homéomorphisme.

\comment
\begingroup\leftskip1em\parindent0pt
\everypar{\llap{--\ }}\parskip1pt

On note $\egoth(\tau)\dans\Ig$ l'ensemble des éléments \expression{extrémaux de $\tau$}, \idest ceux figurant à l'extrémité droite des lignes de $\tau$ (cases grisées). 
Par exemple, on a  $\egoth(\tau)=\set w,a/$ pour
$$
\tau:=\tableauc{}{1&u&v&w\\
2\\
3&a}\eqno(\ddagger)$$

Si $w\in\egoth(\tau)$, on note $\tau{-}w\in\TTT_{b}(\Ig{-}w)$ le tableau obtenu en effaçant la case de $\tau$ contenant $w$.
Dans l'exemple $(\ddagger)$,
$$
\tau{-}w:=\tableauc{}{1&u&v\\
2\\
3&a}\hskip0.8cm$$

On note $\pgoth(\tau)$ la partition de $\Ig$ déterminée 
par le regroupement des éléments des lignes de $\tau$ en sous-ensembles de $\Ig$. Dans l'exemple $(\ddagger)$, 
$$\pgoth(\tau)=\bigset{\set u,v,w/},{\set a/}/\,.$$
Alors, conformément à la notation de \ref{nota-pgoth}, on note 
$$\Fg_{\tau}:=\Fg_{\pgoth(\tau)}\dans\Delta_{b}(\Zg^{\Ig}\times\Fb)\,.$$
La projection $\pi_{b}:\Fg_{\tau}\to\Fb$ est clairement un homéomorphisme.
\par\endgroup
\endcomment
\end{list}

\ssubsubsection{Les isomorphismes $\Psi_{b}^{\Ig}$}

\noindent On\label{defi-Psi-b-Ig} définit par induction sur $|\Ig|$ un isomorphisme de $\S_{\Ig}\times\Sb$-modules\glossary{${\Psi_{b}^{\Ig}:\bigoplusnl_{\tau\in\TTT_b(\Ig)}
\Hc(\Fg_{\tau}){\to}\Hc(\Delta_{|\Ig|+b}(\Zg^{\Ig}\times\Fb))}$:isomorphisme de $\S_{\Ig}\times\Sb$-modules}
$$\displayboxit{\Psi_{b}^{\Ig}:\bigoplusnl_{\tau\in\TTT_b(\Ig)}
\Hc(\Fg_{\tau})\hfonto{}{+|\Ig|}{1cm}\Hc(\Delta_{|\Ig|+b}(\Zg^{\Ig}\times\Fb))}$$

On utilisera de manière systématique le complexe fondamental (\ref{theo-complexe-pgoth-exact}) pour l'espace $(\Fg_{\Ig}\times\Fb)$. Lorsque $\Ig\not=\emptyset$, on a $\Hc(\Fg_{\Ig}\times\Fb)=0$, car $\Hc(\Zg)=0$, et ce complexe donne la présentation
de $\S_{|\Ig|}\times\Sb$-modules
$$\def\ab{|\Ig|+b}
\def\Fa{\Fg_{\Ig}}\def\Za{\Zg^{\Ig}}
\mathalign{
\Delta_{\ab-2}(\Fa\times\Fb)\hf{c_2}{+1}{0.7cm}
\Delta_{\ab-1}(\Fa\times\Fb)\hf{c_1}{+1}{0.7cm}
\Delta_{\ab}(\Fa\times\Fb)}
\eqno(\diamond)$$
où la notation `$\Hc(\_)$' est absente pour gagner de la place.

\comment
On utilisera de manière systématique le morphisme de complexes fondamentaux (\ref{theo-complexe-pgoth-exact}) correspondant à l'inclusion ouverte $(\Fg_{\Ig}\times\Fb)\dans(\Zg^{\Ig}\times\Fb)$. Lorsque $\Ig\not=\emptyset$, on a $\Hc(\Fg_{\Ig}\times\Fb)=\Hc(\Zg^{\Ig}\times\Fb)=0$, car $\Hc(\Zg)=0$, et le morphisme en question est un morphisme de résolutions
de $\S_{|\Ig|}\times\Sb$-modules
$$\hss\mathrigid2mu
\def\Fa{\Fg_{\Ig}}\def\Za{\Zg^{\Ig}}
\def\ab{|\Ig|+b}
\hdecale{-0.cm}{\xymatrixc{@C=2mm@R=5mm}{
\mllap{\cdots\to}
\Delta_{\ab-2}(\Fa\times\Fb)[-2]
\ar[d]_{\hdecale{-2.3cm}{c_2}}^{\hdecale{2.25cm}{c_3}}
\ar[r]&
\Delta_{\ab-1}(\Fa\times\Fb)[-1]
\ar[d]_{\hdecale{-2.3cm}{c_1}}
\ar[r]&
\Delta_{\ab}(\Fa\times\Fb)\aregal[d]
\mrlap{\to0}
\\
\mllap{\cdots\to}
\Delta_{\ab-2}(\Za\times\Fb)[-2]
\ar[r]&
\Delta_{\ab-1}(\Za\times\Fb)[-1]
\ar[r]&
\Delta_{\ab}(\Za\times\Fb)
\mrlap{\to0}
}}
\eqno(\diamond)$$
où la notation `$\Hc(\_)$' est absente pour gagner de la place. 
\endcomment

\smallskip
\noindent On rappelle qu'on est sous l'hypothèse $b>0$.

\soustitreline{\slshape Le cas $|\Ig|=0$. }L'ensemble $\TTT_{b}(\emptyset)$ comporte l'unique tableau
\smashbot{$\tau=
\tableauc{}{1\\
2\\
\vdots\\
b}
$}. On a clairement
$\Fg_{\tau}=\Fb=\Delta_{|\Ig|+b}(\Zg^{\Ig}\times\Fb)$ et l'on pose 
$\Psi_{b}^{\emptyset}=\id_{\Fb}\,.$

\soustitreline{\slshape Le cas $|\Ig|=1$. }L'ensemble $\TTT_{b}(\Ig)$ comporte $b$ tableaux
$$
\tableauc{}{1&\bullet\\
2\\
\vdots\\
b}\,,
\qquad
\tableauc{}{1\\
2&\bullet\\
\vdots\\
b}
\,,\qquad\ldots
\qquad
\tableauc{}{1\\
2\\
\vdots\\
b&\bullet}\,.
$$

Le complexe fondamental $(\diamond)$ est réduit au seul isomorphisme
$$c_{1}:\Hc(\Delta_{b}(\Zg\times\Fb))\hf{\simeq}{+1}{0.7cm}
\Hc(\Delta_{1+b}(\Zg\times\Fb))
$$
où le terme de gauche se décompose naturellement en
$$\Hc(\Delta_{b}(\Zg\times\Fb))=
\bigoplusnl_{\tau\in\TTT_{b}(\Ig)}\Hc(\Fg_{\tau})\,.
$$
On définit alors $\Psi_{b}^{\Ig}=c_1$, c'est un isomorphisme de $\S_{\Ig}\times\Sb$-modules
$$
\Psi_{b}^{\Ig}:\bigoplusnl_{\tau\in\TTT_{b}(\Ig)}\Hc(\Fg_{\tau})
\hf{\Sigma\Psi_{\tau}}{+|\Ig|}{1cm}
\Hc(\Delta_{1+b}(\Zg\times\Fb))\,.
$$
où on a noté $\Psi_{\tau}$ la restriction de $\Psi_{b}^{\Ig}$ au facteur $\Hc(\Fg_{\tau})$.

\soustitreline{\slshape Le cas $|\Ig|>1$. }Le complexe fondamental $(\diamond)$ fournit la surjection
$$c_{1}:\Hc(\Delta_{|\Ig|+b-1}(\Fg_{\Ig}\times\Fb))\hfonto{}{+1}{0.7cm}
\Hc(\Delta_{|\Ig|+b}(\Fg_{\Ig}\times\Fb))
$$
où le terme de gauche se décompose en
$$
\relax{\Hc(\Delta_{|\Ig|+b-1}(\Fg_{\Ig}\times\Fb))=
\bigoplusnl_{x\in\Ig\,,1\leq j\leq b}\Hc(\Delta^{(x,j)}_{|\Ig|+b-1}(\Fg_{\Ig}\times\Fb))}
$$
où $\Delta^{(x,j)}_{|\Ig|+b-1}(\Fg_{\Ig}\times\Fb)$ est l'ouvert-fermé de $\Delta_{|\Ig|+b-1}(\Fg_{\Ig}\times\Fb)$ des uplets dont les éléments en coordonnées $x\in\Ig$ et $j\in\iii[1,b]$ coïncident.
Notons $c_{1}(x,j)$ la restriction de $c_1$ au terme en question, \idest
$$
\smash{c_{1}(x,j):\Hc(\Delta^{(x,j)}_{|\Ig|+b-1}(\Fg_{\Ig}\times\Fb))
\hf{}{+1}{0.7cm}
\Hc(\Delta_{|\Ig|+b}(\Fg_{\Ig}\times\Fb))\,.}
$$
En négligeant la coordonnée redondante $x$, on  identifie ensuite
$$
\relax{\mathalign{
\Hc(\Delta^{(x,j)}_{|\Ig|+b-1}(\Fg_{\Ig}\times\Fb))
&=&
\Hc(\Delta_{|\Ig|+b-1}(\Fg_{\Ig- x}\times\Fb))\,,
}}
\postdisplaypenalty10000$$
où le terme de droite est l'image de l'isomorphisme 
$\Sg_{\Ig-x}\times\Sb$-modules
$$\Psi_{b}^{\Ig-x}:\bigoplusnl_{\kern-3mm\vrule height15pt width0pt
\tau\in\TTT_{b}(\Ig{-}x)}\hskip-1cm
\Hc(\Fg_{\tau})\hf{\sim}{+|\Ig{-}x|}{1.2cm}
\Hc(\Delta_{|\Ig|+b-1}(\Fg_{\Ig- x}\times\Fb))
$$
que l'on peut supposer déjà défini (induction).
On a donc les morphismes
%
$$\mathalign{
\bigoplusnl_{\kern-3mm\vrule height15pt width0pt
\tau\in\TTT_{b}(\Ig{-}x)}\hskip-1cm
\Hc(\Fg_{\tau})\hf{\Psi_{b}^{\Ig-x}}{+|\Ig|}{1cm}
\Hc\big(\Delta^{(x,j)}_{[-1]}(\Fg_{\Ig}\times\Fb)\big)\hf{c_1(x,j)}{+1}{1.1cm}
\Hc\big(\Delta_{[0]}(\Fg_{\Ig}\times\Fb)\big)\hfill
}
\postskip4pt
\eqno(\diamond\diamond)$$
où la notation $\Delta_{[-i]}(\_)$ abrège $\Delta_{|\Ig|+b-i}(\_)$.

\medskip

Maintenant, pour $\nu\in\TTT_{b}(\Ig{-}x)$ notons $\nu{+}(x,j)$ le tableau de
$\TTT_{b}(\Ig)$ obtenu en rajoutant la case $x$ à l'extrémité droite de la ligne $j$ de $\nu$. Par exemple,
$$
\text{si\ }\nu=\tableauc{}{1&u&v\\
2\\
3&a}\,,\qquad
\nu{+}(x,1)=\tableauc{}{1&u&v&x\\
2\\
3&a}\text{\quad et\quad }
\nu{+}(x,2)=\tableauc{}{1&u&v\\
2&x\\
3&a}\,.
$$

Alors, pour chaque $(x,j)$ et $\nu\in\TTT_{b}(\Ig{-}x)$, on 
identifie $\Fg_{\nu{+}(x,j)}$ à $\Fg_{\nu}$ par l'oubli de la $x$-ième coordonnée et l'on définit le morphisme (de degré $|\Ig|$)
$$\Psi_{\nu,(x,j)}:\Hc(\Fg_{\nu+(x,j)})\hf{}{+|\Ig|}{0.9cm}\Hc(\Delta_{|\Ig|+b}(\Fg_{\Ig}\times\Fb))\eqno(\ddagger\ddagger)$$ 
comme la restriction de la composée de morphismes $(\diamond\diamond)$ à $\Hc(\Fg_{\nu})$.

\begin{lemm}\label{lemm-Psi-tau} $\Psi_{\nu,(x,j)}$ est déterminé par le tableau $\tau:=\nu{+}(x,j)$.
\end{lemm}
\demo Dans la suite, l'opérateur `$\Hc(\_)$' étant omniprésent, il sera effacé pour gagner de la place dans les diagrammes. De même, la notation
$\Delta_{[-i]}$ sera un raccourci pour  $\Hc(\Delta_{|\Ig|+b-i}(\Fg_{\Ig}\times\Fb))$.

Soient $w,x,y,z\in\Ig$ et $i,j,k,l\in\iii[1,b]$, avec $x\not=y$ et $j\not=k$, et considérons le diagramme suivant.

\penalty-1500
\begin{figure}[h!]
$$\preskip-1ex\let\FIFb\relax
\def\po{{\smash{\Big(}}}
\def\pf{\smash{\Big)}}
\rlap{\hskip5mm\vtop to0pt{\kern-68pt
\def\dots{\ar@{..}}\def\fdots{\ar@{<..}}
\xymatrix@C=10.5cm@R=3.5cm{
\fdots[d]|{\vrule height7pt depth4pt width0pt c_2}\dots[r]&
\fdots[d]|{\vrule height7pt depth4pt width0pt c_2}\\
\dots[r]&}
\vss}}
\def\vfluhook#1#2#3{\hdecale{-1mm}{\rotatebox{90}{$\hfhook{#1}{}{#3}$}}}
\def\L#1{\mllap{(#1)}\quad}
\mathalign{
&(0)&&(1)&&(2)\\\noalign{\kern3mm}
\L0&\Delta_{[0]}\FIFb&\hfhook{\delta}{}{1cm}\po&
\Delta_{[0]}\FIFb&\oplus&
\Delta_{[0]}\FIFb&\pf\\
&\vfluonto{c_1}{+1}{0.5cm}&&
\vflu{c_1(x,j)}{+1}{0.5cm}&&
\vflu{c_1(y,k)}{+1}{0.5cm}&&\\
\L1&
\bigoplus_{(w,i)}\Delta_{[-1]}^{(w,i)}
&\hfonto{\rm proj}{}{1cm}\po&
\Delta_{[-1]}^{(x,j)}\FIFb&\oplus&\Delta_{[-1]}^{(y,k)}\FIFb&\pf\\
\noalign{\kern-7pt}
&\vfld{}{}{0.5cm}&&\vfld{}{}{0.5cm}&&\vfld{}{}{0.5cm}\\
\L2&\Delta_{\leq[-1]}\FIFb&\hf{\rm rest}{}{1cm}\po&
\Delta_{\leq[-1]}^{(x,j)}\FIFb&\oplus&\Delta_{\leq[-1]}^{(y,k)}\FIFb&\pf\\
&\vfld{0}{}{0.5cm}&&\vfld{0}{}{0.5cm}&&\vfld{0}{}{0.5cm}\\
\L3&\Delta_{\leq[-2]}\FIFb&\hf{\rm rest}{}{1cm}\po&
\Delta_{\leq[-2]}^{(x,j)}\FIFb&\oplus&\Delta_{\leq[-2]}^{(y,k)}\FIFb&\pf\\
&\vfluonto{}{}{0.5cm}&&\vfluonto{}{}{0.5cm}&&\vfluonto{}{}{0.5cm}\\
\noalign{\kern-3pt}
\L4&
\bigoplus_{w\not=z,i\not=l}\Delta_{[-2]}^{(w,i)(z,l)}
&\hf{\rm projj}{}{1cm}\po&
\bigoplus_{x\not=z,j\not=l}\Delta_{[-2]}^{(x,j)(z,l)}
&\oplus&
\bigoplus_{y\not=z,k\not=l}\Delta_{[-2]}^{(y,k)(z,l)}\FIFb&\pf\\
\noalign{\kern-7pt}
&\vfluhook{}{}{0.5cm}&&\vfluhook{}{}{0.5cm}&&\vfluhook{}{}{0.5cm}\\
\L5&\Delta_{[-2]}^{(x,j)(y,k)}\FIFb&\hfhook{\delta}{}{1cm}\po&
\Delta_{[-2]}^{(x,j)(y,k)}\FIFb&\oplus&\Delta_{[-2]}^{(y,k)(x,j)}\FIFb&\pf\\
}
\postskip0pt$$
\vskip1em\mycaption{\slshape Diagramme I}\vskip-1ex
\end{figure}

\noindent
Dans le diagramme $I$
\def\varitemizeseps{\leftmargin1em\itemsep1pt\parsep0pt\parskip0pt\topsep3pt}
\begin{itemize}
\mynobreak\nobreak
\item Les lignes (0) et (5) sont les plongements diagonaux.

\item Les lignes (1,2,3) résultent d'appliquer verticalement la suite  de foncteurs
$$
\to
\Hc(\Delta_{[-1]}(\_))
\to
\Hc(\Delta_{\leq[-1]}(\_))
\too^{\rho}
\Hc(\Delta_{\leq[-2]}(\_))
\to\,,
$$ 
scindée en $\rho$ (\ref{coro-scindage-F-pgoth}). 
Les morphismes de liaison entre lignes $(3)[-1]\to(1)$ sont alors injectifs, et ceux des colonnes (1) et (2) sont même bijectifs puisque 
l'on a \smashtop{$\Hc(\Delta_{\leq[-1]}^{(x,j)})=\Hc(\Delta_{\leq[-1]}^{(y,k)})=0,.$}

\item La flèche `$\hfonto{\rm proj}{}{0.8cm}$' en ligne (1) est la surjection de $\Hc(\Delta_{[-1]})$ sur les facteurs indiqués de sa décomposition en somme directe. Tandis que la flèche `$\hf{\rm projj}{}{0.8cm}$' en ligne (4) 
est la projection de $\Hc(\Delta_{[-2]})$, colonne par colonne, sur les sommes de facteurs indiqués, elle n'est pas nécessairement surjective.
\item Les flèches `$\hf{\rm rest}{}{0.8cm}$' désignent les restrictions aux les sous-espaces fermés.
\item La ligne (5) est restriction de la (4) aux espaces indiqués et c'est un plongement diagonal. L'injection $(5)\hook(4)$ est par prolongement par zéro.

\item
Le diagramme est commutatif par construction.
\end{itemize}

Suite à ces remarques, on remplace le sous-diagramme en pointillé par les morphismes $c_2$ de complexes fondamentaux. On obtient alors le diagramme commutatif {\slshape II}.

\begin{figure}[h!]
$$\preskip0.em\let\FIFb\relax
\dimen11=.8cm
\def\po{{\smash{\Big(}}}
\def\pf{\smash{\Big)}}
\def\vfluhook#1#2#3{\hdecale{-1mm}{\rotatebox{90}{$\hfhook{#1}{}{#3}$}}}
\def\L#1{\relax{(#1)\,}}
\mathalign{
\L0&\Delta_{[0]}\FIFb&\hfhook{\delta}{}{\dimen11}\po&
\Delta_{[0]}\FIFb&\oplus&
\Delta_{[0]}\FIFb&\pf\\
&\vfluonto{c_1}{+1}{0.5cm}&&
\vflu{c_1(x,j)}{+1}{0.5cm}&&
\vflu{c_1(y,k)}{+1}{0.5cm}&&\\
\L1&
\bigoplus_{(w,i)}\Delta_{[-1]}^{(w,i)}\hskip3.4ex
&\hfonto{\rm proj}{}{\dimen11}\po&
\Delta_{[-1]}^{(x,j)}\FIFb&\oplus&\Delta_{[-1]}^{(y,k)}\FIFb
&\oplus\;&\Delta_{[-1]}^{(w,l)}\;
&\pf\\
\noalign{\kern-12pt}
&\vflu{c_2}{+1}{0.5cm}&&\vfluonto{c_2(x,j)}{+1}{0.5cm}&&\vfluonto{c_2(y,k)}{+1}{0.5cm}&&\vfluonto{c_2(w,l)}{+1}{0.5cm}\\
\L4&
\bigoplus_{w\not=z,i\not=l}\Delta_{[-2]}^{(w,i)(z,l)}
&\hf{\rm projj}{}{\dimen11}\po&
\bigoplus_{x\not=z,j\not=l}\Delta_{[-2]}^{(x,j)(z,l)}
&\oplus&
\bigoplus_{y\not=z,k\not=l}\Delta_{[-2]}^{(y,k)(z,l)}\FIFb
&\oplus\;&\Delta_{[-2]}^{(w,l)}\;&\pf\\
\noalign{\kern-7pt}
&\vfluhook{}{}{0.5cm}&&\vfluhook{}{}{0.5cm}&&\vfluhook{}{}{0.5cm}\\
\L5&\Delta_{[-2]}^{(x,j)(y,k)}\FIFb&\hfhook{\delta}{}{\dimen11}\po&
\Delta_{[-2]}^{(x,j)(y,k)}\FIFb&\oplus&\Delta_{[-2]}^{(y,k)(x,j)}\FIFb&\pf\\
}
\postskip0pt$$
\vskip1em\mycaption{\slshape Diagramme II}
\end{figure}

\noindent
Maintenant, pour $\tau=\tableauc{}{
{}\vdots\\
j&\dots&\dots&x\\
{}\vdots\\
k&\dots&y\\
{}\vdots
}\in\TTT_{b}(\Ig)\,,$
notons 
$\tau{-}x$, $\tau{-}y$, $\tau{-}xy$ les tableaux obtenus en enlevant les boites $x,y$  comme indiquent les notations.
Le morphisme $\Psi_{\tau-x,(x,j)}$ est, par la définition
$(\ddagger\ddagger)$, la composée
$$
\mathalign{
\Hc(\Fg_{\tau-xy})
\hf{\Psi_{\tau-xy,(y,k)}}{+|\Ig-xy|}{1.5cm}
\Delta_{[-1]}(\Fg_{\Ig-x}\times\Fb)=
\Delta_{[-2]}^{(x,j)}
\hf{c_{1}(x,j)\circ c_{2}(x,j)}{+2}{2.3cm}}
\Delta_{[0]}\,.
$$
Or, l'image de $\Psi_{\tau-xy,(x,j)}$ est contenue dans le facteur
$\Delta_{[-2]}^{(x,j)(y,k)}$ de $\Delta_{[-2]}^{(x,j)}$, facteur qui est le même pour $\Psi_{\tau-xy,(y,k)}$ modulo le plongement diagonal de la ligne (5).
L'égalité
$$\Psi_{\tau-x,(x,j)}=\Psi_{\tau-y,(y,k)}$$
résulte alors de ce que la projection de 
\goodsmash{1}{.9}{$c_2\big(\Delta_{[-2]}^{(x,j)(y,k)}\big)$} sur un facteur de $\Delta_{[-1]}$ de la forme \halfsmashtop{$\Delta_{[-1]}^{(w,l)}$} avec $(w,l)\not\in\set (x,j),(y,k)/$ est nul. Ce qui est clair déjà sur la ligne (4) où $\Delta_{[-2]}^{(x,j)(y,k)}\cap
\Delta_{[-1]}^{(w,l)}=0$  et par commutativité du diagramme.
\enddemo

\medskip
 Le lemme précédent montre que le morphisme $\Psi_{(\tau{-}x),(x,j)}$ de  la définition $(\dagger\dagger)$ est indépendant de l'écriture
$\tau=(\tau{-}x){+}(x,j)$, raison pour laquelle il sera noté simplement\glossary{$\Psi_{\tau}$:restriction de $\Psi_{b}^{\Ig}$ à $\Hc(\Fg_{\tau})$}
$$\Psi_{\tau}
:\Hc(\Fg_{\tau})\hf{}{+|\Ig|}{0.9cm}\Hc(\Delta_{|\Ig|+b}(\Fg_{\Ig}\times\Fb))\,.
$$

\begin{prop}Le\label{prop-Psi-iso} morphisme
$$
\Psi_{b}^{\Ig}:\bigoplusnl_{\tau\in\TTT_{b}(\Ig)}\Hc(\Fg_{\tau})
\hf{\Sigma\Psi_{\tau}}{+|\Ig|}{1cm}
\Hc(\Delta_{|\Ig|+b}(\Fg_{\Ig}\times\Fb))\,.
$$
est un isomorphisme de $\S_{\Ig}\times\Sb$-modules.
\end{prop}
\demo Le fait que c'est un morphisme de $\S_{\Ig}\times\Sb$-modules résulte de ce que, par construction, pour tout $\alpha\in\S_{\Ig}\times\Sb$ le diagramme
$$\mathalign{
\Hc(\Fg_{\tau})&\hf{\Psi_\tau}{+|\Ig|}{0.9cm}&\Hc(\Delta_{|\Ig|+b}(\Fg_{\Ig}\times\Fb))\\
\vfld{\alpha}{}{0.5cm}&&\vfld{\alpha}{}{0.5cm}\\
\Hc(\Fg_{\alpha(\tau)})&\hf{\Psi_{\alpha(\tau)}}{+|\Ig|}{0.9cm}&\Hc(\Delta_{|\Ig|+b}(\Fg_{\Ig}\times\Fb))}$$
est commutatif.
Il est surjectif puisqu'il en est ainsi de
$$\bigoplusnl_{(x,j)}\Delta_{[-1]}^{(x,j)}=\Delta_{[-1]}\hfonto{c_1}{}{0.7cm}\Delta_{[0]}\,,$$
et que \smashtop{$\Psi_{b}^{\Ig{-}\ast}$} est une surjection sur 
chaque \smashtop{$\Delta_{[-1]}^{(\ast,\cdot)}$} par hypothèse de récurrence.
Pour conclure maintenant que l'on a un isomorphisme il suffit de comparer les dimensions. On a
$$\smashbot{\dim\Big(\bigoplusnl_{\tau\in\TTT_{b}(\Ig)}\Hc(\Fg_{\tau})\Big)=|\TTT_{b}(\Ig)|\cdot\dim\Hc(\Fb)}\eqno(d_1)$$ 
et, par \ref{theo-suite-spectrale-basique}(\ref{theo-suite-spectrale-basique-b}),
$$\dim\big(\Hc(\Delta_{|\Ig|+b}(\Fg_{\Ig}\times\Fb))\big)=
|\FFF(b,|\Ig|+b)|\cdot\dim\Hc(\Fb)\,.\eqno(d_2)$$
L'égalité de dimensions dans $(d_1)$ et $(d_2)$ résulte alors
du lemme \ref{lemme-TTT-FFF} qui établit que l'on a $|\FFF(b,|\Ig|+b)|=|\TTT(b,|\Ig|+b)|=|\TTT_b(\Ig)|$.
\enddemo

\medskip

Cette proposition achève la dernière étape de la définition inductive des isomorphismes $\Psi_{b}^{\Ig}$ commencée dans \ref{defi-Psi-b-Ig}.

\subsubsection{Structure de $\StS{m-b}{b}$-module de $\Hc(\Delta_{m}(\Zg^{m-b}\times\Fb(\Xg)))$}
On\label{structure-module-Delta-Za-Xb} applique les considérations précédentes au cas où $\Ig=\iii[1,m{-}b]$ et où l'intervalle
$\iii[1,b]$ est décalé vers $\iii[m{-}b{+}1,m]$. Il est alors avantageux de remplacer la notation de l'ensemble $\TTT_{b}^{\Ig}$ par celle, équivalente, de l'ensemble de tableaux $\TTT(b,m)$ de la section \ref{tableaux<->FFF}.

\ssubsubsectionline{Tableaux normaux.}Chaque\label{tableaux-normaux} orbite de $\StS{m-b}{b}$ dans $\TTT(b,m)$ contient un unique tableau où la longueur des lignes est décroissante et le contenu des boites `$\bullet$' est strictement croissant. 
Par exemple, si $m=50$ et $m{-}b=14$,
$$\mycount=14
\def\ss{\relax\global\advance\mycount by 1
\relax\the\mycount}
\def\everyygcell{\scriptsize}
\yght=3.5mm\ygwd=0.45cm
\let\ygvdots\Ygvdots
\tau=\left\{\tableauc{}{
\ss&1&2&3&4&5&6\\
\ss&7&8&9\\
\ss&10&11&12\\
\ss&13&14\\
\ss\rparcell{0pt}{b}{1.7cm}{31mm}\\
\vdots\\
\cellcolor{0 0 0 0.10}29\\
\vdots\rparcell{0pt}{b{-}(m{-}b)=22}{6mm}{3mm}\\
50}\right.\eqno(\ast)
$$
Un tel tableau sera dit \expression{normal}, leur ensemble est noté $\TTT_{0}(b,m)$\glossary{$\TTT_{0}(b,m)$:ensemble de tableaux \expression{normaux} dans $\TTT(b,m)$}.

\begin{rema}On a indiqué en gris la plus petite colonne de lignes de cardinal $1$ dans les tableaux normaux, son cardinal, $b{-}(m{-}b)$, est atteint lorsque la longueur des lignes de $\tau$ est $\leq2$ et que $m\geq2(m{-}b)$. Il y a alors autant de tableaux normaux que des partitions de $\iii[1,m{-}b]$ et l'application 
$(\_)^{\bullet}:\TTT(b,m)\to\TTT^{\bullet}(b{+}1,m{+}1)$ de \ref{lemme-TTT-FFF} est une bijection sur $\TTT(b{+}1,m{+}1)$.
\end{rema}

 \begin{lemm}[et notation]
Le\label{HHH-tau} stabilisateur $\HHH_{\tau}$ de $\tau$ dans $\StS{m-b}{b}$
\glossary{$\HHH_{\tau}$:le stabilisateur de $\tau$ dans $\StS{m-b}{b}$} est le produit direct
$$\preskip-1ex
\HHH_{\tau}=\HHH_{\utau}\times\S_{m{-}|\utau|}\,,$$
où $\utau$ désigne le sous-tableau des lignes de longueur $>1$
de $\tau$, et où $\HHH_{\utau}$ désigne le stabilisateur de $\utau$ dans $\S_{|\utau|}\cap\big(\StS{m-b}{b}\big)$. On a 
$$\pmathalign{
\HHH_{\utau}\dans
\S_{m{-}b}\times\S_{\iii[\mkern2mu m{-}b{+}1,m-|\utau|\mkern2mu]}\,,
\hfill\\\noalign{\kern2pt}
\1_{|\tau|}\times\S_{m{-}|\utau|}\dans\1_{m{-}b}\times\S_{b}\,.\hfill
}$$
La projection $p_2:\HHH_{\utau}\to\S_{\iii[\mkern2mu m{-}b{+}1,m-|\utau|\mkern2mu]}$ est un isomorphisme sur son image.
\end{lemm}

\def\Remarque{Commentaire}\begin{rema}Dans \ref{tableaux-normaux}-$(\ast)$, on a $\HHH_{\utau}\dans\S_{\iii[1,14]}\times\S_{\iii[15,18]}$. C'est le groupe cyclique d'ordre $2$ engendré par l'involution
{\small$\arraycolsep2pt\begin{pmatrix}
16&7&8&9\cr
17&10&11&12
\end{pmatrix}$}.

\end{rema}

\begin{prop}
Pour\label{prop-iso-sym} toute inclusion d'espaces $i$-acycliques $\Xg\dans\Zg$ avec $\Hc(\Zg)=0$, il existe un isomorphisme canonique de $\StS{m-b}{b}$-modules
$$
\relax{\Psi_{b}^{m}:\bigoplusnl_{\kern-0mm\vrule height13pt width0pt
\tau\in\TTT(b,m)}\hskip-0.cm
\ind^{\StS{m-b}{b}}_{\HHH_{\tau}}\Hc(\Fg_{\tau}(\Xg))\hf{\simeq}{+(m{-}b)}{1.2cm}\Hc(\Delta_{m}(\Zg^{m-b}\times\Fb(\Xg)))}
$$
\end{prop}
\demo Corollaire de \ref{prop-Psi-iso}.
\enddemo

\subsubsection{Structure de $\StS{m{-}(\pp)}{\pp}$-module de $\Hbm(U^{m}_{m-p,\ldots,m})$
}
On\label{action-StS->FFF} reprend maintenant le sujet de \ref{stabilite-des-familles-pseudovarietes} qui concerne la donnée d'une pseudovariété orientée $\Mg$ de type fini de dimension $d_{\Mg}$.
On appliquera la proposition \ref{prop-iso-sym} au cas où
$\Zg:=\Mg_{\geq0}$, où $\Xg:=\Mg_{>0}$ et où $b=\pp$.
Pour chaque tableau $\tau\in\TTT(\pp,m)$, on a \glossary{${\Fg_{\tau}:=\set\cl x\in
\Mgge^{m{-}(\pp)}\times\Fg_{\pp}\mid x_i=x_{f_{\tau}(i)}/}$:}
$$\halfdisplayskips
\Fg_{\tau}(\Mg_{>0}):=\bigset\cl x\in
\Mgge^{m{-}(\pp)}\times\Fg_{\pp}(\Mgg)\mid x_i=x_{f_{\tau}(i)}/\,.$$

\noendpoint\begin{prop}\label{prop-tableau-normaux-ss-basiques}
\begin{enumerate}\mynobreak\nobreak
\item Il\label{prop-tableau-normaux-ss-basiques-a} existe un isomorphisme canonique de $\StS{m{-}(\pp)}{\pp}$-modules gradués, de degré $-(m{-}(\pp))d_{\Mg}$,
$$\halfdisplayskips
\mathalign{
\alt\otimes\Hbm(\UU^{m}_{m-p,\ldots,m})&\cong&
\bigoplus_{\hskip-4mm
\tau\in\TTT(\pp,m)
\hskip-4mm}\ind^{\StS{m{-}(\pp)}{\pp}}_{\HHH_{\tau}}
\alt^{d_{\Mg}}\otimes\Hbm(\Fg_{\tau}(\Mg_{>0}))\,,}
\postdisplaypenalty10000$$
où $\alt$ est la signature de $\S_{m{-}(\pp)}$.
On a $|\TTT(\pp,m)|\leq|\Pgoth(m{-}(\pp))|$ avec égalité pour tout $m\geq2(m{-}(\pp))$.
\item Le\label{prop-tableau-normaux-ss-basiques-b} $\FI$-module défini 
dans \ref{theo-suite-spectrale-basique-relative} par les morphismes des suites spectrales 
$$\IEs(q_m^{*})_{1}^{p,q}:\IEs(\U^{m})^{p,q}_{1}\to\IEs(\U^{\mm})^{\pp,q}_{1}\,,$$
où $\mathrigid1mu
q=i+(m{-}(\pp))$, est canoniquement isomorphe au $\FI$-module défini par les morphismes
$$
\bigoplus_{\tau\in\TTT(\pp,m)}\left(\vcenter{
\setbox0=\hbox{\ $\displaystyle
\ind^{\Sm}_{\HHH_{\utau}\times\S_{m-|\utau|}}
\alt^{d_{\Mg}}\otimes\Hbm^{Q}(\Fg_{\pp}(\Mgg))
$\ }
\setbox1=\hbox{\ $\displaystyle
\ind^{\Smm}_{\HHH_{\utau}\times\S_{\mm-|\utau|}}
\alt^{d_{\Mg}}\otimes\Hbm^{Q}(\Fg_{p+2}(\Mgg))
$\ }
\dimen0=\wd0\ifdim\dimen0<\wd1\dimen0=\wd1\fi
\advance\dimen0 by1ex
\hbox to\dimen0{\hss\box0\hss}
\hbox to\dimen0{\hss$\vfld{}{\ind p_{\pp}^{*}}{5mm}$\hss}
\hbox to\dimen0{\hss\box1\hss}
}
\right)
$$
où $\alt$ est la restriction à $\HHH_{\utau}$ de la 
 signature de $\S_{m{-}(\pp)}$, où
$\Fg_{\pp}$ et $\Fg_{p{+}2}$ réfèrent respectivement à $\Fg_{\pgoth(\tau)}$ et $\Fg_{\pgoth(\tau^{\bullet})}$ {\rm(\cf\ref{donnees-notas-dec-sym}-\ref{tau-goth} et \ref{lemme-TTT-FFF})}, et où
$$
Q:=q-(m{-}(\pp))\,\dMg=i-(m{-}(\pp))\,(\dMg{-}1)\,.$$
\end{enumerate}
\end{prop}
\def\Demonstration{Indication}
\demo  (\ref{prop-tableau-normaux-ss-basiques-a}) En dualisant \ref{prop-iso-sym} 
et en incorporant les caractères signature nécessaires (\ref{def-action-image-inverse}-(\ref{def-action-image-inverse-b})),
on a l'isomorphisme de $\StS{m-(\pp)}{\pp}$-modules
$$
\mathalign{
(\Psi_{b}^{m})\dual:
\alt^{d_{\Mg}{+}1}_{m-(\pp)}\otimes\alt^{d_{\Mg}{+}1}_{\pp}
\otimes\Hbm(\UU^{m}_{m-p,\ldots,m})
&\hf{\simeq}{-(m{-}(\pp))d_{\Mg}}{1.2cm}\hfill\\
\noalign{\kern6pt}
&\hskip-4cm
\bigoplusnl_{\tau\in\TTT(b,m)}
\ind^{\StS{m-b}{b}}_{\HHH_{\tau}}\alt^{d_{\Mg}+1}_{\pp}\otimes
\Hbm(\Fg_{\tau}(\Xg))\,,}
$$
dont on conclut en simplifiant les caractères redondants (\cf note ($^{\ref{note-ind-res}}$)). L'assertion 
 (\ref{prop-tableau-normaux-ss-basiques-b}) est alors une simple reformulation de \ref{theo-suite-spectrale-basique-relative}-{\ref{theo-suite-spectrale-basique-relative-b}}.
\enddemo

\subsubsection{Monotonie et stabilité dans les suites spectrales basiques}
Pour chaque $\tau\in\TTT(\pp,m)$, on\label{rang-stabilite-termes-basiques} reconnaît dans \ref{prop-tableau-normaux-ss-basiques}-(\ref{prop-tableau-normaux-ss-basiques-b}) le $\FI$-module
$$\Ind_{\rho,\tau}\big(\bigset
\Hbm^{Q}(\Fg_{\pp}(\Mgg))/_{\pp}\big)\in
\ModFI k_{\geq 2(m{-}(\pp))}\,.$$
où $\rho$ désigne l'action de $\HHH_{\utau}$ sur $\alt_{m{-}(\pp)}$ et 
$\Ind_{\rho,\tau}$ est le foncteur d'induction introduit dans \ref{Ind-alt-lambda}. 

\begin{theo}Soit\label{theo-stabilite-BM-pseudo} $\Mg$ une pseudovariété connexe orientée avec $\dMg{\geq}2$.
\mynobreak\begin{enumerate}\nobreak
\item Pour\label{theo-rang-pseudovariete-a} $Q:=i-(m{-}(\pp))\,(\dMg{-}1)$ et $\tau\in\TTT(\pp,m)$, on a
$$\preskip1.3em\let\Big\big
\let\strut\relax\pmathalign{\rkms\Big(
\Ind_{\rho,\tau}\big(\set\Hbm^{Q}(\Fg_{\pp}(\Mgg))/_{\pp}\big)\Big)&\leq&
\smashtop{\begin{cases}\noalign{\kern-4pt
}4i\,,\text{ si $\dMg=2$}\\
2i\,,\text{ si $\dMg\geq3$.}\\\noalign{\kern-3pt}
\end{cases}}\\\noalign{\kern2pt}
\rkm\Big(
\Ind_{\rho,\tau}\big(\set\Hbm^{Q}(\Fg_{\pp}(\Mgg))/_{\pp}\big)\Big)&\leq&i\,.\hfill}
$$

\item \tolerance7000
Pour  $a,i\in\NN$ fixés, le\label{theo-rang-pseudovariete-b} $\FI$-module 
$\set\Hbm ^{i}(\Delta_{?m-a}\Mg^{m})/_{m}$
est monotone et stable pour $m\geq4i+4a$, si $\dMg=2$, et pour $m\geq2i+4a$, si $\dMg\geq3$. Les familles des caractères et des nombres de Betti correspondantes sont (donc) polynomiales sur les mêmes intervalles d'entiers $m$.
\end{enumerate}
\end{theo}
\demo\halfdisplayskips
 (\ref{theo-rang-pseudovariete-a}) On applique le théorème \ref{theo-stabilite-BM-i-acyclique} à l'espace $i$-acyclique $\Mgg$. 
Comme $\dim(\Mgg)\geq3$, la famille $\set\Hbm^{Q}(\Fg_{\pp}(\Mgg))/_{\pp}$ est monotone pour $(\pp)\geq Q$ et est monotone 
 et stable pour $(\pp)\geq 2Q$. Le théorème \ref{theo-Ind-rho-lambda} garantit alors que le $\FI$-module induit par $\Ind_{\rho,\tau}$ est monotone pour tout
$$
m\geq Q+\big(m{-}(\pp)\big)=i-\big(m{-}(\pp)\big)(\dMg{-}2)\,,
$$
et est monotone 
 et stable pour tout
$$m\geq 2Q+4\big(m{-}(\pp)\big)=2i+\big(m{-}(\pp)\big)\big(4{-}2(\dMg{-}1)\big)\,,\eqno(\ast)$$
et comme on dispose de l'inégalité 
\ref{theo-suite-spectrale-basique}-(\ref{theo-suite-spectrale-basique-d}): 
$i\geq (m{-}(\pp))(d_{\Mg}{-}1)$,  le dernier terme de $(\ast)$ est majoré par $4i$ si $\dMg=2$, et par $2i$ si $\dMg\geq 3$.

\medskip
(\ref{theo-rang-pseudovariete-b}$_1$) Le cas de la famille $\bigset 
p_{m}^{*}:\Hbm^{i}(\Fm(\Mg))\to\Hbm^{i}(\Fmm(\Mg))/_{m}$.

Fixons $i\in\NN$. D'après \ref{theo-suite-spectrale-basique} les termes $\IEs(\U^{m})^{p,q}_{r}$ de la suite spectrale basique 
$(\IEs(\U^{m})_{r},d_r)$ qui convergent vers $\Hbm^{i}(\Fm(\Mg))$ sont ceux pour lesquels on a $q=i+(m{-}(\pp))$. Si maintenant on fixe $q$, la différence $(m{-}p)$ est constante et si l'on augmente $m$ et $p$ simultanément, ces contraintes décrivent un facteur direct du $\FI$-module défini par la famille de morphismes de suites spectrales basiques  (\ref{theo-suite-spectrale-basique-relative})
$$\bigset\IEs(q_m^{*})_{r}^{p,q}:\IEs(\U^{m})^{p,q}_{r}\to\IEs(\U^{\mm})^{\pp,q}_{r}/{\vrule depth3pt width0pt}_{m}\,,$$
facteur que nous allons noter 
$$\EEE^{i,q}_{r}:=
\bigset\IEs(q_m^{*})_{r}^{p,q}:\IEs(\U^{m})^{p,q}_{r}\to\IEs(\U^{\mm})^{\pp,q}_{r}\mid 
q=i+(m{-}(\pp))/{\vrule depth3pt width0pt}_{m}\,.$$
Les différentielles des suites spectrales basiques définissent alors un complexe de $\FI$-modules
$$
(\EEE_r(q),d_r):=\big(\cdots\too
\EEE_{r}^{i-1,q+r-1}
\too^{d_r}
\EEE_{r}^{i,q}
\too^{d_r}
\EEE_{r}^{i+1,q-r+1}
\too
\cdots\big)\,.
\eqno(\ast)$$
Nous avons des isomorphismes canoniques
$$h^{i}(\EEE_r(q),d_r)\simeq\EEE_{r+1}^{i,q}
\text{\quad et \quad}\EEE_{r}^{i,q}=\EEE_{r+1}^{i,q}\,,\ \forall r> q{+}1\,,$$
et donc la condition de convergence (\ref{theo-suite-spectrale-basique-relative}-(\ref{theo-suite-spectrale-basique-relative-a}))
$$
\Big(\bigoplusnl_{q\in\NN}\EEE^{i,q}_{q+2}\Big)
\Longrightarrow
\bigset 
p_{m}^{*}:\Hbm^{i}(\Fm(\Mg))\to\Hbm^{i}(\Fmm(\Mg))/_{m}\,.\eqno(\diamond)$$

\smallskip
\noindent\begingroup\sl Lemme. 
Pour tout $r\geq1$, on a
$$
\pmathalign{
\rkm(\EEE_{r}^{i{+}1,q})&=&\rkms(\EEE_{r}^{i,q})=0\,,\hfill&\quad\forall i<0\,,\hfill\\\noalign{\kern1pt}
\rkms(\EEE_{r}^{i,q})&\leq&\sup\set i\,\epsilon,1/\,,\hfill&\quad\forall i\geq0\,,\hfill\\\noalign{\kern1pt}
\rkm(\EEE_{r}^{i,q})&\leq&\sup\set(i{-}1)\epsilon ,1/\,,\hfill&\quad\forall i\geq1\,,\hfill
}
$$ 
avec $\epsilon=4$ si $d_\Mg=2$, et $\epsilon=2$ si $d_\Mg\geq3$.
\endgroup
$$\def\strut{\vrule height12pt depth5pt width0pt}
\def\bb#1/{\hbox to1cm{\hss$#1$\hss}}
\begin{array}{|c||c|c|c|c|c||c|}
\hline
\strut
\bb \EEE_r/&
\bb\EEE_r^{-1}/&
\bb\EEE_r^{0}/&
\bb\EEE_r^{1}/&
\bb\EEE_r^{2}/&
\bb\EEE_r^{3}/&
\bb\EEE_r^{i}/\\
\hline
\strut\rkms\leq&
0&1&\epsilon&2\epsilon&3\epsilon&i\epsilon\\
\hline
\strut\rkm\leq&
0&0&1&\epsilon&2\epsilon&(i{-}1)\epsilon\\
\hline
\end{array}$$

\normaldisplayskips

\begingroup
\parskip0.5ex
\medskip\noindent{\sl Preuve. } $\bullet$ {\slshape Le cas $r=1$. }On a par définition et \ref{prop-tableau-normaux-ss-basiques}-(\ref{prop-tableau-normaux-ss-basiques-b})
$$(\EEE_{1}^{i,q})_m=\bigoplusnl_{\tau\in\TTT(\pp,m)}\Ind_{\rho,\tau}
\Hbm^{Q}(\Fpp(\Mg_{>0}))\,.\eqno(\ddagger)$$
L'assertion (\ref{theo-rang-pseudovariete-a}) donne,
pour chaque $\tau\in\TTT(\pp,m)$ et tout $i\geq0$
$$\pmathalign{
\rkms(\Ind_{\rho,\tau}(\Hbm^{Q}(\Fpp(\Mg_{>0})))&\leq&i\epsilon \hfill\\
\rkm(\Ind_{\rho,\tau}(\Hbm^{Q}(\Fpp(\Mg_{>0})))&\leq&i
}$$
La majoration $\rkm(\EEE_{1}^{i,q})\leq i$ s'ensuit  puisque $(\_)^{\bullet}:\TTT(\pp,m)\to\TTT(p{+}2,\mm)$ est toujours injective, mais pour avoir 
\smash{$\rkms(\EEE_{1}^{i,q})\leq \epsilon i$}, il faut s'assurer que si l'on a $m\geq \epsilon i$ et \smashtop{$\Hbm^{Q}(\Fpp(\Mg_{>0}))\not=0$}, l'application $(\_)^{\bullet}$ est bijective. Raisonnons par l'absurde, si $(\_)^{\bullet}$ n'est pas bijective, on a $(\pp)<m/2$, et si de plus
$m\geq \epsilon i$, on a\label{bullet-non-bijective}
$$
(d_{\Mg}{-}1)(m-(\pp))>\smash{{d_\Mg{-}1\over 2}m\geq {d_\Mg{-}1\over 2}\epsilon i\geq i\,,}$$
donc $Q=i-(d_{\Mg}{-}1)(m-\pp)<0$ et alors $\Hbm^{Q}(\Fpp(\Mg_{>0}))=0$. 

Le tableau qui suit illustre ces conclusions dont on remarquera qu'elles sont compatibles à celles du lemme pour $\EEE_{1}$ puisque $i\leq\epsilon(i{-}1)$ dès que $2\leq i$.
$$
\def\strut{\vrule height12pt depth5pt width0pt}
\def\bb#1/{\hbox to1cm{\hss$#1$\hss}}
\begin{array}{|c||c|c|c|c|c||c|}
\hline
\strut
\bb \EEE_1/&
\bb\EEE_1^{-1}/&
\bb\EEE_1^{0}/&
\bb\EEE_1^{1}/&
\bb\EEE_1^{2}/&
\bb\EEE_1^{3}/&
\bb\EEE_1^{i}/\\
\hline
\strut\rkms&
0&0&\epsilon&2\epsilon&3\epsilon&i\epsilon\\
\hline
\strut\rkm&
0&0&1&2&3&i\epsilon\\
\hline
\end{array}\eqno(\ddagger\ddagger)$$

\smallskip
{\slshape$\bullet$ Le cas général. }Si nous appliquons au tableau $(\ddagger\ddagger)$ les règles suivantes, établies dans \ref{prop-generalites-monotonie-stabilite}-(\ref{prop-generalites-monotonie-stabilite-e}),
$$\pmathalign{
\rkm\EEE_{r+1}^{i}&\leq&\sup\set
\rkms\EEE_{r}^{i-1},\rkm\EEE_{r}^{i}/\hfill\\
\rkms\EEE_{r+1}^{i}&\leq&\sup\set
\rkms\EEE_{r}^{i-1},\rkms\EEE_{r}^{i},\rkm\EEE_{r}^{i+1}/}
$$
on obtient aussitôt le tableau du lemme pour $\EEE_{2}$, et comme ce tableau est laissé stable par ces règles le lemme résulte pour tout $r\geq1$.
 \hfill$\boxminus$
\endgroup

\smallskip
Ceci étant, on rappelle que la convergence $(\diamond)$ aboutit en fait sur le bi-gradué du $\FI$-module $\H^{i}:=\set 
p_{m}^{*}:\Hbm^{i}(\Fm(\Mg))\to\Hbm^{i}(\Fmm(\Mg))/_{m}$ relativement à la filtration par le degré $p$ des cochaînes de \v Cech (\cf\ref{comm-coh-BM-suite-spectrale}). Or, comme nous sommes soumis à la condition
$\binspace1
Q=i-(m-(\pp))(d_{\Mg}{-}1)>0$, la filtration en question possède un nombre fini de termes ($\leq$\smash{${i\over d_{\Mg}{-}1}+1$}), et
les majorations des rangs de $\EEE_{q+2}^{i,q}$ sont aussi valables pour $\H^{i}$ d'après \ref{coro-monotonie-stabilite}.

Dans le cas particulier où $i=0$, ces raisonnements donnent seulement la majoration $\rkms\H^{0}\leq1$ (c'est ce qui arrive lorsque $\Hbm^{0}(\Mg)=0$). Lorsque $\Hbm^{0}(\Mg)=k$, un calcul direct élémentaire montre que l'on a $\H^{0}=\V(0)$ et donc que $\rkms\H^{0}=0$.

\medskip
(\ref{theo-rang-pseudovariete-b}$_{2}$) Le cas des familles $\set\Hbm ^{i}(\Delta_{m-a}\Mg^{m})/_{m}$ résulte de \ref{theo-rang-pseudovariete-a}) exactement comme dans le cas ({\bf C}) des espaces $i$-acycliques 
(p.~\pageref{Cas-C}).

\medskip
(\ref{theo-rang-pseudovariete-b}$_{3}$)
Pour le cas général des familles $\set\Hbm^{i}(\Delta_{\leq m-a}\Mg^{m})/_{m}$, on ne peut pas faire appel au complexe fondamental des $\FI$-modules, comme dans le cas ({\bf D}) des espaces $i$-acycliques (p.~\pageref{Cas-D}), puisque ce complexe n'est plus exact. On fait plutôt appel aux suites longues de $\FI$-modules
de \ref{prop-structure-FI-modules}
$$
\to
\set\Hbm^{i-\dMg}(\Delta_{\leq m-a-1}\Mg^{m})/
\to
\set\Hbm^{i}(\Delta_{\leq m-a}\Mg^{m})/
\to
\set\Hbm^{i}(\Delta_{m-a}\Mg^{m})/
\to
$$
où l'on peut supposer que les termes de droite et gauche sont monotones et stables pour $m\geq 2i+4a$ ou $m\geq4i+4a$ suivant le cas. On en déduit (\ref{coro-monotonie-stabilite}-(\ref{coro-monotonie-stabilite-b})) le même rang de monotonie et stabilité pour $\set\Hbm^{i}(\Delta_{\leq m-a}\Mg^{m})/$ puisque extension d'un noyau et conoyau de $\FI$-modules monotones est stables sur le même rang des entiers $m$.
\cqfd
\enddemo

\def\Remarque{Commentaire}
\begin{rema}
Au delà du fait que $\Xg$ est maintenant une pseudovariété, la démarche  de Church dans la démonstration de son théorème de stabilité (\cite{chu}, \cf\ref{stabilite-caracteres}) est très différente de la nôtre. Elle utilise de manière essentielle les résultats de Totaro (\cite{tota}) sur la de suite spectrale de Leray associée au plongement $\Fm(\Xg)\hook\Xg^{m}$, valable parce que $\Xg$ est lisse.
\end{rema}

\section{Calcul du caractère de $\S _{m}$-module de $\Hc(\Fm (\Xg))$}
\glossarytitle{Calcul du caractère de $\S _{m}$-module de $\Hc(\Fm (\Xg))$}
On\label{caractere} étend les résultats de Macdonald (\cite{mac}) concernant le caractère du $\Sm$-module $\Hc(\Xg^{m})$ au
cas du $\Sm$-module  $\Hc(\Fm (\Xg))$ lorsque $\Xg$ est $i$-acyclique.

\def\Remarque{\miseengarde Avertissement}
\begin{rema*}{\miseengarde Dans cette section et la suivante,  $\car(k)=0$ et
la donnée d'un
espace $i$-acyclique $\Xg$ présuppose que $\dim\Hc(\Xg)<\infty$.}
\killline
\end{rema*}

\subsection{Série de caractères de $\S _{m}$}
\subsubsectionline{Série de caractères d'un $\S _{m}$-module gradué.}\label{serie-de-traces}Nous appellerons $\S _{m}$-module gradué, la donnée d'un $k$-espace vectoriel gradué $\smashbot{V:=\bigoplus_{i\in\ZZ}V^{i}}$, tel que chaque $V^i$ est un $\S _{m}$-module. 

Lorsque, de plus, chaque $V^{i}$ de dimension finie sur $k$, la \expression{série de caractères de $V$}, est l'application $\chi_{V}:\S _{m}\to k[[T]]$ définie par (\cf\cite{mac} (2.3))
\glossary{${\chi_{V}:\S _{m}\to k[[T]]}$:série de caractères d'un $\S _{m}$-module gradué $V$}
$$\S _{m}\ni\alpha\mapsto
\chi_{V}(\alpha,T):=\sumnl_{i\in\ZZ}(-1)^{i}\tr(\alpha\actson V^{i})\, T^{i}\,.$$
Le lemme suivant, pendant du lemme \ref{poincare-elem}, est élémentaire.\killline

\nopoint\begin{lemm}\label{serie-elem}
\begin{enumerate}
\item\leavevmode\label{serie-elem-a}Si $(\dots\to V_{i-1}\to V_{i}\to V_{i+1}\to\cdots)$ est un complexe borné de $\S _{m}$-modules gradués, on a 
$\sumnl_{i\in\ZZ}(-1)^{i}\,\chi_{V^{i}}(\alpha,T)=0\,.$
\item\leavevmode\label{serie-elem-b}Pour $r\in\ZZ$, on a
$\chi_{V[-r]}(\alpha,T)=\chi_{V}(\alpha,T)\cdot(-T)^{r}\,.$
\item\leavevmode\label{serie-elem-c}$\chi_{V_1\otimes_k V_2}(\alpha_1\otimes\alpha_2,T)=\chi_{V_1}(\alpha_1,T)\cdot\chi_{V_2}(\alpha_2,T)$.
\item\leavevmode\label{serie-elem-d}$\P(\mathop{\rm Homgr}_{k}(V_1,V_2))(T)=\P(V_1)(1/T)\cdot\P(V_2)(T)$
\end{enumerate}
\end{lemm}

\subsubsectionline{Série de caractères d'un $\S _{m}$-espace topologique.}Si\label{def-serie-caracteres} $\Mg$ est un espace topologique de type fini, muni d'une action de $\S _{m}$, les \expression{séries de caractères de $\Hc(\Mg)$ et de $\Hr(\Mg)$} seront respectivement notées\glossary{${\chic(\Mg)(\_,T):\S _{m}\to k[[T]]}$:série de caractères de $\S _{m}$-module de $\Hc(\Mg)$}
$$\chic(\Mg)(\alpha,T):=\chi_{\Hc(\Mg,k)}(\alpha,T)\,,\quad\hbox{et}\quad
\chi(\Mg)(\alpha,T):=\chi_{\Hr(\Mg,k)}(\alpha,T)\,.$$
On a donc vis-à-vis des définitions de \ref{char-i}, l'égalité
$$\chic(\Mg)(\_,T)
=\sumnl_{i\in\NN}(-1)^{i}\, \chi_{\Hc^{i}(\Mg)}(\_)\, T^{i}
=\sumnl_{i\in\NN}(-1)^{i}\, \chic(\Mg;i)(\_)\, T^{i}\,.$$

\begin{rema}Il\label{rema-serie-caracteres} convient de retenir pour la suite les expression suivantes.
Pour tout espace topologique $\Xg$ et tout $m\geq0$, on a :
$$\mathalign{
\hfill{\Pc(\Xg^{m})(-T)\over T^m}&=&\Big({\chic(\Xg)(\1,T)\over T}\Big)^m\,,\hfill\\\noalign{\kern5pt}
{\Pc(\Fm (\Xg))(-T)\over T^{m}}&=&\Big({\chic(\Xg)(\1,T)\over T}\Big)\usp m\,.}$$\vskip-1.5em\vskip-1ex
\end{rema}

\subsection{Séries de caractères de $\Fg^{\qgoth}(\Xg)$}L'exactitude des complexes fondamentaux associés aux espaces $\Delta_{\leq\ell}(\Fg^{\qgoth})$ du théorème \ref{theo-complexe-pgoth-exact} s'avère particulièrement adaptée à la détermination des séries de caractères. La proposition suivante est un ingrédient important dans le calcul des caractères.

\begin{theo}\label{caracterisation-trace}Soit $\Xg$ un espace $i$-acyclique. Pour toute partition $\qgoth$ de $\iii[1,m]$, tout $\ell\leq m$ et tout $\alpha\in\S ^{\qgoth}$ (\ref{Sqgoth}), on a
$$\halfdisplayskips{\chic\big(\Delta_{\leq\ell}\Fg^{\qgoth}(\Xg)\big)(\alpha,T)\over T^{\ell}} =
\sumnl_{\aa\geq 0}
{\chic\big(\Delta_{\ell-\aa}\Fg^{\qgoth}(\Xg)\big)(\alpha,T)\over T^{\ell-\aa}}\,.
\postskip0pt$$
\end{theo}
\demo Corollaire immédiat de \ref{theo-complexe-pgoth-exact} et du lemme \ref{serie-elem}.
\enddemo

\subsection{Séries de traces pour un $m$-cycle}
L'action de $\S _{m}$ sur $\iii[1,m]$ induit une action sur l'ensemble $\Pgoth(m)$ des partitions de $\iii[1,m]$ qui sera sous-entendue dans la suite. 
Pour $m>0$ donné, on note $\sigma_{m}\in\S _{m}$\glossary{$\sigma_{m}\in\S _{m}$:la permutation cyclique $(1,2,\ldots,m)$}, ou simplement $\sigma$ s'il est superflu de préciser l'entier $m$, la permutation cyclique $(1,2,\ldots,m)$. On notera $\Cg_{m}:=\langle\sigma_{m}\rangle$\glossary{${\Cg_{m}:=\langle\sigma_{m}\rangle}$:sous-groupe de $\S _{m}$ engendré par $\sigma_{m}$} le sous-groupe de $\S _{m}$ engendré par $\sigma_{m}$.\killline

\begin{lemm}\label{partitions-Cm-fixes}Les points fixes de l'ensemble $\Pgoth(m)$ sous l'action de $\Cg_{m}$ sont les partitions en orbites de $\iii[1,m]$ sous l'action des différents sous-groupes de $\Cg_{m}$. On a 
$$\halfdisplayskips\begin{casesalign}
\Pgoth_{d}(m)^{\Cg_{m}}=\bigset{\iii[1,m]/{\langle\sigma^{m/d}_{m}\rangle}}/\,,&\hbox{\ si $d\div m$,}\\
\Pgoth_{d}(m)^{\Cg_{m}}=\emptyset\,,\hfill&\hbox{\ sinon.}\hfill
\end{casesalign}
$$ 
En particulier, $\big|\Pgoth(m)^{\Cg_{m}}\big|=\big|\bigset d\in\NN\hbox{\ t.q.\ } d\div m/\big|\,.$
\end{lemm}
\demo L'application $\Cg_{m}\to\iii[1,m]$, $\alpha\mapsto\alpha(1)$ est un isomorphisme de $\Cg_{m}$-espaces lorsque l'on munit $\Cg_{m}$ de son action par multiplications à gauche. Les partitions $\Cg_{m}$ stables par multiplication à gauche sont aussi stables par multiplication à droite puisque $\Cg_{m}$ est abélien. Or, les partitions de $\Cg_{m}$ stables par multiplication à droite sont les ensemble des classes d'équivalence à gauche de $\Cg_{m}$ des sous-groupes $H\dans\Cg$, autrement dit, ce sont les partitions de $\Cg_{m}$ en $H$-orbites. Enfin, comme $\Cg_{m}$ est cyclique engendré par $\sigma$, ses sous-groupes sont les $\langle\sigma^{d}\rangle$ avec $d\div m$.
\enddemo

\begin{theo}\label{theo-trace-s}Soit $\Xg$ un espace $i$-acyclique. Pour tout $m\geq1$, la série de traces de l'action de $\sigma_{m}$ sur $\Hc(\Fm (\Xg))$ vérifie
$$
\displayboxit{{\chic(\Fm (\Xg)(\sigma_{m},T))\over T^{m}}=\sum_{d\div m}\mu\Big({m\over d}\Big)
{\chic(\Xg)(\1,T^{d})\over T^{d}}}\,,$$
où $\mu(\_)$\glossary{$\mu(\_)$:fonction de Möbius} est la fonction de Möbius. 
\end{theo}
\demo
Le théorème \ref{caracterisation-trace} appliqué à la partition $\qgoth=(1,\ldots,1)$ et pour $\ell=m$ donne l'égalité de séries de traces
$${\chic(\Xg^{m})(\sigma,T)\over T^m}=
\sum_{\aa\geq0}
{\chic(\Delta_{m-\aa}\Xg^{m})(\sigma,T)\over T^{m-\aa}}
\eqno(\ast)$$

D'après \ref{connexes}, on a la décomposition 
$$\Hc(\Delta_{m-\aa}\Xg^{m})=\bigoplusnl_{\pgoth\in\Pgoth_{m-\aa}(m)}\Hc(\Fg_{\pgoth})\,,$$ 
et la trace de l'action de $\sigma$ sur $\Hc(\Delta_{m-\aa}\Xg^{m})$ se lit sur les termes $\Hc(\Fg_{\pgoth})$ tels que $\sigma\cdot\pgoth=\pgoth$. Or, on a vu dans le lemme \ref{partitions-Cm-fixes} qu'il n'y a de telles partitions que lorsque $d:=m-\aa$ est un diviseur de $m$, et dans ces cas, il y a une \emph{et une seule} partition telle, à savoir $\pgoth:=\Cg_{m}\cdot\langle\sigma_{m}^{d}\rangle$. L'action de $\sigma_{m}$ sur $\Fg_{\pgoth}$ coïncide avec l'action de $\sigma_{d}$ sur $\Fg_{d}$, on a donc l'égalité de séries de traces
$$
\begin{casesalign}
\chic(\Delta_{d}\Xg^{m})(\sigma_{m},T)&=&\chic(\Fg_d(\Xg))(\sigma_d,T)\,,&
\hbox{ si $d\div m\,,$}\\
\chic(\Delta_{d}\Xg^{m})(\sigma_{m},T)&=&0\,,\hfill&\hbox{ si $d\notdiv m\,,$}
\end{casesalign}$$
dont on déduit l'expression du deuxième membre de $(\ast)$ suivante:
$$\sum_{\aa\geq0}
{\chic(\Delta_{m-\aa}\Xg^{m})(\sigma_{m},T)\over T^{m-\aa}}
=\sum_{d\div m}{\chic(\Fg_{d}(\Xg))(\sigma_{d},T)\over T^{d}}\,.
\eqno(\dagger)$$

Maintenant, si $p_i:\Xg^{m}\to\Xg$ est la projection $\cl x\mapsto x_i$, on a $\sigma^{*}\circ p_{\sigma i}^{*}=p_{i}^{*}$ et l'action de $\sigma$ sur
les tenseurs simples de $\Hc^{k}(\Xg)^{\otimes m}$ se fait par permutation signée   (\cf\cite{mac}):
$$\sigma^{*}(\omega_1\otimes\omega_2\otimes\cdots\otimes\omega_{m})=(-1)^{k(\mmo)}
(\omega_2\otimes\omega_3\otimes\cdots\otimes\omega_{1})\,.
$$
Comme l'opérateur $\sigma^{*}$ agit sur une base de tenseurs simples de $\Hc(\Xg)^{\otimes m}$, sa trace le lit sur ceux de tels tenseurs fixés par $\sigma$, donc de la forme $\omega\otimes\cdots\otimes\omega$ pour $\omega$ appartenant à une base de $\Hc(\Xg)$. On a alors l'égalité
$$\left\{\mathalign{
\tr(\sigma^{*}\actson\Hc^{km}\Xg^{m})&=&(-1)^{k(\mmo)}\dim\Hc^{k}(\Xg)\\
\tr(\sigma^{*}\actson\Hc^{j}\Xg^{m})\hfill&=&0\,,\hbox{\ si $m\notdiv j$,}\hfill}\right.$$
et, par définition de la série des traces \ref{serie-de-traces}, 
$$\mathalign{\chic(\Xg^{m})(\sigma_{m},T)&=&\sum_{k\geq 0}(-1)^{km}\cdot(-1)^{k(\mmo)}\dim\Hc^{k}(\Xg)\cdot T^{km}\hfill\\
&=&\sum_{k\geq 0}(-1)^{k}\dim\Hc^{k}(\Xg)\cdot (T^{m})^{k}=\chic(\Xg)(\1,T^{m})\,.\hfill}\eqno(\ddagger)$$

En reportant les égalités $(\dagger)$ et $(\ddagger)$ dans la formule ($\ast$), on obtient:
$${\chic(\Xg)(\1,T^{m})\over T^{m}}=\sum_{d\div m}{\chic(\Fg_{d}(\Xg))(\sigma_{d},T)\over T^{d}}\,,\quad\forall m\geq1\,,$$
et la formule d'inversion de Möbius donne l'égalité recherchée:
$${\chic(\Fm (\Xg)(\sigma_{m},T))\over T^{m}}=\sum_{d\div m}\mu\Big({m\over d}\Big)
{\chic(\Xg)(\1,T^{d})\over T^{d}}\,.$$
\vskip-2em\enddemo

\begin{rema}[pour $\Xg^{m}$]\label{theo-trace-s-Xm}D'après l'égalité $(\ddagger)$ de la preuve précédente qui donne
$${\chic(\Xg^{m})(\sigma_{m},T)}={\chic(\Xg)(\1,T^{m})}\,,$$
l'analogue du théorème \ref{theo-trace-s} pour l'espace $\Xg^{m}$ s'énonce par la même formule où la sommation est restreinte
au seul terme d'indice $d:=m$. On retrouve alors les résultats de Macdonald (\cite{mac}, formule (3.2)).
\end{rema}

\subsection{Séries des traces pour une puissance d'un $m$-cycle}
\subsubsection{Partitions en orbites $\qgoth_r$}On\label{intro-trace-puissance} s'intéresse maintenant au cas où $m=dr$ avec $1\leq d\leq m$ et à la trace de l'opérateur $\sigma_{dr}^{r}$ agissant sur $\Fg_{dr}(\Xg)$. Notons $\qgoth_{r}$ la partition de $\iii[1,dr]$ et $\langle\sigma_{dr}^{r}\rangle$-orbites. Pour l'étude de $\chic(\Fg_{dr}(\Xg))(\sigma_{dr}^r,T)$ nous allons utiliser le théorème \ref{caracterisation-trace} avec $\qgoth:=\qgoth_r$ et $\ell:=dr$. Dans ce cas, on a  (\ref{serie-elem}-(\ref{serie-elem-c}))
$${\chic\big(\Fg^{\qgoth_r}(\Xg)\big)(\sigma_{dr}^{r},T)\over T^{dr}}=\Big({\chic(\Fg_d(\Xg))(\sigma_d,T)\over T^{d}}\Big)^{r}\,,\postskip0pt$$
d'où l'égalité:
$$\Big({\chic(\Fg_d(\Xg))(\sigma_d,T)\over T^{d}}\Big)^{r}=
\sum_{\aa\geq 0}
{\chic\big(\Delta_{dr-\aa}\Fg^{\qgoth_r}(\Xg)\big)(\sigma_{dr}^{r},T)\over T^{dr-\aa}}\,,$$
qui nous emmène à l'étude de l'espace $\Delta_{dr-\aa}\Fg^{\qgoth_r}(\Xg)$ sous l'action de $\sigma_{m}^{r}$.

\begin{prop}\label{prop-transversalite-qr}
Soit $\qgoth_{r}$ la partition de $\iii[1,dr]$ en $\langle\sigma_{dr}^{r}\rangle$-orbites et soit $\pgoth=\set I_1,\ldots,I_\ell/\trans\qgoth_{r}$ telle que 
$\sigma_{dr}^{r}\cdot\pgoth=\pgoth$.
Alors,  $\ell=dr'$, pour un certain $r'$ vérifiant $1\leq r'\leq r$, et l'on a une identification d'espaces munis d'actions
$$(\sigma_{dr}^{r}\actson\Fg_{\pgoth}(\Xg))\simeq(\sigma_{dr'}^{r'}\actson\Fg_{dr'}(\Xg))\,.
$$
Le cardinal de l'ensemble de telles partitions $\pgoth$ est 
$$\Big|\big(\qgoth^{\trans}\cap\Pgoth_{dr'}(dr)\big){}^{\sigma^{r}_{dr}}\Big|=\Parties r {r'}\,d^{(r-r')}\,.$$
\end{prop}
\demo Commençons par remarquer qu'une partition $\pgoth\in\Pgoth(dr)$ vérifie 
$\sigma^{r}\cdot\pgoth=\pgoth$, si et seulement si, la relation $\relp$ est \expression{$\langle\sigma^{r}\rangle$-équivariante}, c'est-à-dire:
$$i\relp j\Rightarrow
\sigma^{ra}(i)\relp\sigma^{ra}(j)\,,\quad\forall a\in\NN.$$

 L'illustration ci-dessous représente, sur la partie gauche, la partition $\qgoth_r$. Chaque colonne étant l'une des parties de $\qgoth_r$, est remplie d'indices $i$ de $\iii[1,dr]$ pour lesquels les coordonnées $x_i$ de $\cl x\in\Fg^{\qgoth}$ sont deux à deux distinctes. Le cycle à l'intérieur des colonnes rappelle l'action simplement transitive du groupe $\langle\sigma^{r}\rangle$ sur chacune d'elles.
 
\medskip
Sur la partie droite, on représente, à l'aide de flèches, les identifications définies par l'équivalence $\relp$. Il faut remarquer le parallélisme des flèches dans le sens vertical qui reflète la $\sigma^{r}$-équivariance de $\relp$.
$$\includegraphics{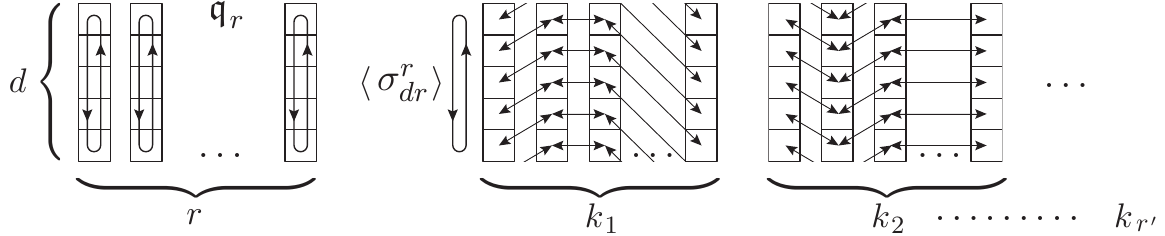}$$
Les regroupements $k_1,k_2,\ldots,k_{r'}$ sont tels que $r=\sum_i k_i$. On voit donc bien qu'au final on a $\Fg_{\pgoth}\simeq\Fg_{dr'}$ et que l'action $\sigma_{dr}^{r}\actson\Fg_{\pgoth}$ s'identifie à $\sigma_{dr'}^{r'}\actson\Fg_{dr'}$.

La partie de droite de la figure permet aussi de comprendre que pour chaque partition de l'ensemble $\qgoth_{r}$ en $r'$ parties non vides de cardinaux $k_i$, les différentes manières de définir les partitions $\pgoth\in(\qgoth^{\trans})^{\sigma_{dr}^{r}}$
 pour les regroupements $k_i$ sont en bijection avec les différentes manières de dessiner les flèches, soit au nombre $d^{k_1-1}d^{k_2-1}\cdots d^{k_{r'}-1}=d^{r-r'}$, nombre qui seul dépend de $r'$ et pas de la suite $(k_1,\ldots,k_{r'})$. La proposition en découle.
\enddemo

\begin{theo}\label{theo-trace-sr}Soit $\Xg$ un espace $i$-acyclique. Pour tous $d,r\geq1\in\NN$, la série de traces de l'action de $\sigma_{dr}^{r}$ sur $\Hc(\Fg_{dr}(\Xg))$ vérifie
l'égalité :
$${\chic(\Fg_{dr}(\Xg))(\sigma_{dr}^{r},T)\over (dT^{d})^r}=\Big({\chic(\Fg_d(\Xg))(\sigma_{d},T)\over dT^{d}}\Big)\usp r\,,\postskip0pt$$
soit
$$\preskip1ex\displayboxit{{\chic(\Fg_{dr}(\Xg))(\sigma_{dr}^{r},T)\over (dT^{d})^{r}}=
\Big(\sum_{e\div d}\mobius de
{\chic(\Xg)(\1,T^{e})\over dT^{e}}\Big)\usp r}$$
où $\mu$ est la fonction de Möbius et $(\_)\usp r$ est la factorielle décroissante de \ref{factorielles}.\end{theo}
\demo Dans la formule de l'introduction \ref{intro-trace-puissance}:
$$\Big({\chic(\Fg_d(\Xg))(\sigma_d,T)\over T^{d}}\Big)^{r}=
\sum_{\aa\geq 0}
{\chic\big(\Delta_{dr-\aa}\Fg^{\qgoth_r}(\Xg)\big)(\sigma_{dr}^{r},T)\over T^{dr-\aa}}\,,$$
nous pouvons maintenant préciser le membre de droite. Par \ref{prop-transversalite-qr}, on a
$$\begin{cases}
\chic\big(\Delta_{\ell}\Fg^{\qgoth_r}(\Xg)\big)(\sigma_{dr}^{r},T)=0\,,\hbox{\quad si $d\notdiv \ell\,$,}\\[5pt]
\chic\big(\Delta_{dr'}\Fg^{\qgoth_r}(\Xg)\big)(\sigma_{dr}^{r},T)=
\smashbot{\Parties r {r'}\,d^{(r-r')}}\,\chic\big(\Fg_{dr'}(\Xg)\big)(\sigma_{dr'}^{r'},T)\,.
\end{cases}$$
En substituant les expressions, on obtient, pour tout $r\geq1$,
$$\Big({\chic(\Fg_d(\Xg))(\sigma_d,T)\over (dT)^{d}}\Big)^{r}=
\sum_{r'\geq 1}\,
\Parties r {r'}\,{\chic\big(\Fg_{dr'}(\Xg)\big)(\sigma_{dr'}^{r'},T)\over (dT^{d})^{r'}}\,,$$
soi, en écriture vectorielle,
$$\let\ds\displaystyle\def\tt{\vrule height 10pt depth 10pt width0pt}
\left(\begin{array}{c}
\ds\Big({\chic(\Fg_d)(\sigma_{d},T)\over dT^{d}}\Big){\vrule height10pt width0pt}^1
\cr\vdots\cr
\ds\Big({\chic(\Fg_d)(\sigma_{d},T)\over dT^{d}}\Big){\vrule height10pt width0pt}^r
\end{array}\right)
=\left(\!\Parties ij\!\right)
\left(\begin{array}{c}
\tt\ds{\chic(\Fg_{d\cdot1})(\sigma_{d\cdot 1}^{1},T)\over (dT^{d})^{1}}
\cr\vdots\cr
\tt\ds{\chic(\Fg_{d\cdot r})(\sigma_{dr}^{r},T)\over (dT^{d})^{r}}
\end{array}\right)\,.
$$

La matrice $\big(\sbParties ij\big)$ est la matrice des nombres de Stirling de deuxième espèce $(\uStirling ij)$ (\ref{Stirling-Cardinaux}). Son inverse est la matrice des nombres de Stirling de première espèce signés $(\ustirling ij)$. Cette matrice est aussi la matrice de changement d'un système de puissances $(X^1,\ldots,X^r)$ vers un système de factorielles décroissantes $(X\usp1,\ldots,X\usp r)$ (\ref{rema-initiales}), par conséquent:
$${\chic(\Fg_{dr}(\Xg))(\sigma_{dr}^{r},T)\over (dT^{d})^r}=\Big({\chic(\Fg_d(\Xg))(\sigma_{d},T)\over dT^{d}}\Big)\usp r\,.$$

Le théorème découle alors de l'égalité du théorème \ref{theo-trace-s} qui donne:
$$
\relax{{\chic(\Fg_{d}(\Xg)(\sigma_d,T))\over dT^{d}}=\sum_{e\div d}\mobius de
{\chic(\Xg)(\1,T^{e})\over dT^{e}}}\,.\postskip0pt$$
\vskip-2.5em
\enddemo

\begin{rema}[pour $\Xg^{m}$]\label{theo-trace-sr-Xm}La remarque de l'introduction \ref{intro-trace-puissance} est la même pour l'espace $\Xg^m$, mais en plus simple. On a clairement (\cf\ref{theo-trace-s-Xm}, mais aussi \cite{mac} (eq.~4.5))
$${\chic(\Xg^{dr})(\sigma_{dr}^{r},T)}
=\chic(\Xg)(\1,T^{d})^{r}\,.
$$
Aussi, l'analogue du dernier théorème \ref{theo-trace-sr} pour l'espace $\Xg^{m}$ est donné par la même formule où la factorielle décroissante $(\_)\usp r$ est remplacée par la puissance $(\_)^r$ et où la sommation est restreinte au terme d'indice $e:=d$. 
\end{rema}

\subsection{Séries des traces pour une permutation générale}
\subsubsection{Permutations immiscibles}
\`A\label{permutations-immiscibles} une permutation $\alpha\in\S _{m}$, on associe la partition $\pgoth_{\alpha}$  de $\iii[1,m]$ en $\langle\alpha\rangle$-orbites. Les cardinaux des $\langle\alpha\rangle$-orbites définissent une décomposition de $m$ que l'on note $\lambda(\alpha)=(1^{\XX_1(\alpha)},2^{\XX_2(\alpha)},\ldots,m^{\XX_{m}(\alpha)})$, ce qui signifie le fait qu'il y a exactement  $\XX_i(\alpha)$ orbites de cardinal $i$. 

\begin{defi*}\label{def-immiscible}\'Etant données deux parties non vides $I$ et $J$  de $\iii[1,m]$ et deux permutations $\alpha\in\S _{I}$ et $\beta\in\S _{J}$ (voir \ref{SgI}), on dit qu'elles sont \expression{immiscibles}, si l'on a $\sum_{i}\XX_{i}(\alpha)\XX_{i}(\beta)=0$, autrement dit, si les orbites de $\langle\alpha\rangle$ dans $I$ ont toutes des cardinaux différents de celles de $\langle\beta\rangle$ dans $J$.
\end{defi*}

\begin{prop}\label{prop-immiscible}Soit $\Xg$ un espace $i$-acyclique. Pour toute partition $\qgoth:=\set J_1,\ldots,J_r/\in\Pgoth(m)$, et toute famille de permutations $\set\alpha_{i}\in\S _{J_{i}}/_{i=1,\ldots,r}$ deux à deux immiscibles, on a
$$\chic(\Fm (\Xg))(\alpha,T)=\prodnl_{i=1}^{r}\chic(\Fg_{|J_i|}(\Xg))(\alpha_i,T)$$
où $\alpha$ désigne le recollement des $\alpha_i$ et où,
dans les terme de droite, nous avons identifié le support $J_i$ de $\alpha_i$ avec l'intervalle $\iii[1,|J_i|]$.
\end{prop}
\demo Grâce au complexe fondamental de $\Fg^{\qgoth}(\Xg)$, pour $\ell=m$ et pour $\alpha$ le recollement des $\alpha_i$, on a l'égalité (\ref{caracterisation-trace})
$${\chic(\Fg^{\qgoth})(\alpha,T)\over T^{m}}=
\sumnl_{\aa\geq0}{\chic(\Delta_{m-\aa}\Fg^{\qgoth}(\Xg))(\alpha,T)\over T^{m-\aa}}\,,$$
où, dans le premier membre, on a clairement
$$\chic(\Fg^{\qgoth})(\alpha,T)=\prodnl_{i=1}^{r}\chic(\Fg_{|J_i|}(\Xg))(\alpha_i,T)\,.$$
Dans le second membre, pour $\aa=0$, on a bien
$$\chic(\Delta_{m}\Fg^{\qgoth}(\Xg))(\alpha,T)=\chic(\Fm (\Xg))(\alpha,T)$$
ce qui nous emmène à montrer que pour tout $\aa\geq1$, on a 
$$\chic(\Delta_{m-\aa}\Fg^{\qgoth}(\Xg))(\alpha,T)=0\eqno(\ast)$$
Or, d'après \ref{prop-transversalite}, on a la décomposition en sous-espaces ouverts
$$\Delta_{\ell}\Fg^{\qgoth}(\Xg):=\coprod\nolimits_{
\let\stacklhss\relax
\scriptstack{\pgoth\in\Pgoth_{\ell}(m),\,
\pgoth\trans\qgoth}}\Fg_{\pgoth}(\Xg)\,,$$
et la série de traces de $\alpha$ pour $\Delta_{\ell}\Fg^{\qgoth}$ est concentrée sur les espaces $\Fg_{\pgoth}$ correspondants aux partitions $\pgoth$ vérifiant, de plus, $\alpha\cdot\pgoth=\pgoth$, donc, telles que
$$(i\relp j)\Longrightarrow (\alpha^{r}i\relp \alpha^{r}j)\,.$$

Maintenant, s'il existait une telle partition pour $\ell<m$, on aurait un certain couple $(i\not=j)$ vérifiant $(i\relp j)$, auquel cas, par transversalité, on aurait $(i\not\relq j)$ et il existerait des parties $J_{a}\not= J_{b}$ de $\qgoth$ telles que $i\in J_a$ et $j\in J_{b}$. En particulier, les orbites $\langle\alpha\rangle\cdot i=\langle\alpha_a\rangle\cdot i$ et $\langle\alpha\rangle\cdot j=\langle\alpha_{b}\rangle\cdot j$ seraient de cardinaux différents puisque $\alpha_a$ et $\alpha_{b}$ sont immiscibles. Si $\epsilon:=\mathop{\rm ord}(\alpha_a)<\mathop{\rm ord}(\alpha_{b})$, on aurait 
$$\preskip0pt(i\relp j)\text{ et\ }(i=\alpha_{a}^{\epsilon}\,i\relp \alpha_{b}^{\epsilon}\,j)\,, \text{\quad et donc\quad}
(j\relp \alpha_{b}^{\epsilon}\,j)\,,$$
avec $\alpha_{b}^{\epsilon}\,j\not= j$, et ceci contredirait la transversalité de $\pgoth$ et $\qgoth$ puisque $(j\relq \alpha^{\epsilon}\,j)$.

\nobreak L'égalité $(\ast)$ se trouve ainsi justifiée et la proposition résulte.
\enddemo

\medskip\noindent On peut maintenant énoncer le théorème principal de cette section.\vskip-1ex\vskip0pt

\begin{theo}Soit\label{theo-trace-gen} $\Xg$ un espace $i$-acyclique. \'Etant donné une permutation $\alpha\in\S _{m}$, soit $\lambda(\alpha)=(1^{\XX_1},2^{\XX_{2}},\ldots,m^{\XX_{m}})$ la décomposition de $m$ déterminée par la partition de $\iii[1,m]$ en $\langle\alpha\rangle$-orbites. Alors, 
$$\displayboxit{{\chic(\Fm (\Xg))(\alpha,T)\over T^{m}}=\prod_{d=1}^{m}d^{\,\XX_d}
\Big(\sum_{e\div d}\mobius de
{\chic(\Xg)(\1,T^{e})\over dT^{e}}\Big)\usp {\XX_d}}
$$
où $\mu$ est la fonction de Möbius et $(\_)\usp r$ est la factorielle décroissante de \ref{factorielles}.
\end{theo}
\demo 
Pour chaque $d=1,2,\ldots,m$, notons $J_d$ la réunion des $\langle\alpha\rangle$-orbites de $\iii[1,m]$ qui sont de cardinal $d$. Notons $\qgoth$ la partition de $\iii[1,m]$ définie par les parties $J_d$ non vides, \idest telles que $\XX_d\not=0$.

Pour chaque $J\in\qgoth$, notons $\alpha_{J}$ la restriction de $\alpha$ à $J$.
On a $\alpha\in\S ^{\qgoth}$, et ses différentes composantes $\alpha_J$ sont deux à deux immiscibles. On peut donc appliquer la proposition \ref{prop-immiscible} à la partition $\qgoth$ et à la permutation $\alpha$. On obtient l'égalité
$${\chic(\Fm (\Xg))(\alpha,T)\over T^{m}}=
\prodnl_{J\in\qgoth}{\chic(\Fg_{|J|}(\Xg))(\alpha_J,T)\over T^{|J|}}\,.\eqno(\ast)$$

Maintenant, comme une partie $J_d\in\qgoth$ est la réunion de $\XX_d$ orbites de cardinal $d$, l'action de $\alpha$ sur $J_d$ est équivalente à l'action de $\sigma_{d\XX_{d}}^{\XX_d}$ sur $\iii[1,d\XX_d]$ et chaque terme du produit dans $(\ast)$ est, d'après \ref{theo-trace-sr},
$${\chic(\Fg_{d\XX_d}(\Xg))(\sigma_{d\XX_d}^{\XX_d},T)\over T^{d\XX_d}}=
d^{\XX_d}\Big({\chic(\Fg_d(\Xg))(\sigma_d,T)\over dT^{d}}\Big)\usp {\XX_d}\,.
$$
Le théorème résulte alors de remarquer que dans le produit de la formule à démontrer, les termes d'indices $d$ tels que $\XX_d=0$ sont tous égaux à $1$.
\enddemo

\begin{rema}[pour $\Xg^{m}$]\label{theo-trace-gen-Xm}Tout comme dans les remarques 
\ref{theo-trace-s-Xm} et \ref{theo-trace-sr-Xm}, l'analogue du dernier théorème \ref{theo-trace-gen} pour l'espace $\Xg^{m}$ est donné par la même formule où $(\_)\usp r$ est remplacée par $(\_)^r$ et où la sommation est restreinte au seul terme d'indice $e:=d$ (\cf aussi \cite{mac} (eq.~4.5)). 
On a donc pour tout $\alpha\in\S _{m}$ avec $\lambda(\alpha)=(1^{\XX_1},2^{\XX_2},\ldots,m^{\XX_{m}})$:
$$\chic(\Xg^{m})(\alpha,T)=\prodnl_{d=1}^{m}\,d^{\XX_d}\,\chic(\Xg)(\1,T^{d})^{\XX_{d}}\,.
$$\end{rema}

\subsection{Comparaison entre $\chic(\Fm (\Xg))$ et $\chic(\Fm (\Xg\mmoins a))$}
Pour $\alpha\in\S _{m}$, notons $\lambda(\alpha)=(1^{\XX_1(\alpha)},2^{\XX_{2}(\alpha)},\ldots\ldots)\vdash m$
son type. Le coefficient $\XX_1(\alpha)$ est le cardinal de l'ensemble des $\alpha$-orbites réduites à un point, autrement dit, de l'ensemble des $i\in\iii[1,m]$ fixés par $\alpha$.

\subsubsectionline{Le cas ou $\XX_1(\alpha)=0$.}
On\label{moins-un-point-comparaison} a vu dans le corollaire de \ref{moins-un-point} que si $\Xg$ est $i$-acyclique, on dispose d'une suite exacte courte de $\Sm$-modules gradués
$$0\to\Hc\big(\Fm ^{\bullet}(\Xg)\big)[-1]\to\Hc(\Fm (\oo\Xg))\to
\Hc(\Fm (\Xg))\to0\eqno(\ast)$$
où $\oo\Xg:=\Xg\mmoins\set \bullet/$ pour un certain $\bullet\in\Xg$, et où $\Fm ^{\bullet}(\Xg)$ est la réunion disjointe des  $\Fm ^{x_i=\bullet}(\Xg):=\set (x_1,\ldots,x_{m})\in\Fm  \mid x_i=\bullet/$ où $i=1,\ldots,m$. 

Or, par l'égalité évidente $\alpha(\Fm ^{x_i=\bullet})=\Fm ^{x_{\alpha(i)}=\bullet}$, on voit que si \hbox{$\XX_1(\alpha)=0$}, on a $\alpha(\Fm ^{x_i=\bullet})\not=(\Fm ^{x_i=\bullet})$, et donc que $\chic(\Fm ^{\bullet})(\alpha,T)=0$. On conclut alors, par l'exactitude de $(\ast)$, que $\chic(\Fm (\Xg\mmoins1))(\alpha,T)=\chic(\Fm (\Xg))(\alpha,T)$, d'où le lemme suivant.

\begin{lemm}\label{lemme-moins-a-points}Soit $\Xg$ un espace $i$-acyclique. Pout tout $m>0$ et toute permutation $\alpha\in\S _{m}$ telle que $\XX_1(\alpha)=0$, on a
$$\chic(\Fm (\Xg\mmoins a))(\alpha,T)=\chic(\Fm (\Xg))(\alpha,T)$$
pour tout $a\in\NN$.
\end{lemm}

\begin{rema}Ce lemme peut aussi être vérifié en appliquant le théorème \ref{theo-trace-gen}. Dans ce cas, il suffit de voir que pour tout $d>1$ on a
$$\sum_{e\div d}\mobius de
{\chic(\Xg\mmoins1)(\1,T^{e})\over T^{e}}=\sum_{e\div d}\mobius de
{\chic(\Xg)(\1,T^{e})\over T^{e}}\,.$$
Or, d'après \ref{poincare-espaces}-(\ref{poincare-espaces-c}), on a
$\chic(\Xg\mmoins1)(\1,T)=\chic(\Xg)(\1,T)-T$, et alors
$$\mathalign{\sum_{e\div d}\mobius de
{\chic(\Xg\mmoins1)(\1,T^{e})\over T^{e}}&=&
\sum_{e\div d}\mobius de
\Big({\chic(\Xg)(\1,T^{e})\over T^{e}}-1\Big)\,.
}$$
On conclut par le fait bien connu sur la fonction de Möbius que affirme que l'on a $\sum_{e\div d}\mu(d/e)=0$ lorsque  $d>1$.
\end{rema}

\subsubsectionline{Le cas général.}Pour $\alpha\in\S _{m}$, notons $I:=\set i\in\iii[1,m]\mid\alpha(i)=i/$ et $J:=\iii[1,m]\mmoins I$. Les permutations $\id_{I}\in\S _{I}$ et $\alpha\rest _{J}\in\S _J$ sont clairement immiscibles et on peut appliquer \ref{prop-immiscible}. On a
$$\halfdisplayskips\chic(\Fm (\Xg))(\alpha,T)=\chic(\Fg_{|I|}(\Xg))(\1,T)\cdot \chic(\Fg_{|J|}(\Xg))(\alpha\rest_{J},T)\,,$$
d'où la proposition suivante.

\begin{prop}\label{prop-moins-a-points}Soit $\Xg$ un espace $i$-acyclique. On a
$$ \let\ds\displaystyle
{\chic(\Fm (\Xg))(\alpha,T)
\vrule depth5pt width0pt\over \vrule height18pt width0pt
\quad\Big(\ds{\chic(\Xg)(\1,T)\over T}\Big)\usp{\XX_1(\alpha)}\quad}
=
{\chic(\Fm (\Xg\mmoins a))(\alpha,T)
\vrule depth5pt width0pt\over \vrule height18pt width0pt
\quad \Big(\ds{\chic(\Xg)(\1,T)\over T}-a\Big)\usp{\XX_1(\alpha)}\quad}
$$
\end{prop}

\subsection{Comparaison entre $\chic(\Fba(\Xg))$ et $\chic(\Fb(\Xg\mmoins a))\cdot \chic(\Fa (\Xg))$}\label{comparaison}Dans la section \ref{Leray}, nous allons nous intéresser de plus près aux projections $\pi_a:\Fba (\Xg)\to\Fa (\Xg)$ de fibre $\Fg_{b}(\Xg\mmoins a)$. Le théorème \ref{degen} établit que lorsque $\Xg$ est, de plus, $i$-acyclique et localement connexe, la suite spectrale de Leray associée à $\pia $ est dégénérée ($d_r=0$, pour $r\geq2$), auquel cas on a un isomorphisme
$$\Hc(\Fba (\Xg))\sim
\Hc(\Fg_{b}(\Xg\mmoins a))\otimes\Hc(\Fa (\Xg))\,.$$
Notons $\S _{a}$ (resp. $\S _{b}$) le sous-groupe des permutations $\alpha\in\S _{b+a}$ telles que
$\alpha(i)=i$ pour tout $i\leq b$ (resp. $b<i$).
Le groupe $\S _{a}$ agit sur la base $\Fa (\Xg)$, et le groupe $\S _{b}$ sur la fibre $\Fg_{b}(\Xg\mmoins a)$. 

La proposition suivante montre qu'il n'est généralement pas vrai que le caractère de $\S _{b}\times\S _{a}$ sur $\Hc(\Fba )$ soit le produit des deux autres caractères, \idest \emph{on n'a pas} pour tous $(\alpha,\beta)\in\S _{b}\times\S _{a}$ l'égalité:
$$\chic(\Fba (\Xg))((\beta,\alpha),T)=\chic(\Fg_{b}(\Xg\mmoins a))(\beta,T)\cdot\chic(\Fa (\Xg))(\alpha,T)\,.\eqno(\ddagger)$$

\begin{prop}\label{compare-b+a}Soit $\Xg$ un espace $i$-acyclique. 
\mynobreak\begin{enumerate}
\nobreak\item\leavevmode\label{compare-b+a-a}Pour $d>1$ et $r>r'\geq1$, et pour $\pi:\Fg_{dr}(\Xg)\to\Fg_{dr'}(\Xg)$ la projection sur la réunion de $r'$ orbites de $\langle \sigma_{dr}^{r}\rangle$. On a
$$\mathrigid1mu
\mathalign{
\hfill\chic(\Fg_{dr}(\Xg))(\sigma_{dr}^{r})=(dT^{d})^{r}Q_{d}(\Xg)\usp r\hfill\\\noalign{\kern4pt}
\chic(\Fg_{b}(\Xg\mmoins dr'))(\sigma_{dr''}^{r''})\cdot\chic(\Fg_{dr'}(\Xg))(\sigma_{dr'}^{r'})=(dT^{d})^{r}
Q_{d}(\Xg)\usp {r''}\cdot Q_{d}(\Xg)\usp {r'}
}$$
où $r'':=r-r'$ et 
$Q_{d}(\Xg):=\chic(\Fg_d(\Xg))(\sigma_d,T)/(dT^{d})$. 

\nobreak
Dans cette situation, l'égalité {\rm($\ddagger$)} n'est donc jamais vérifiée.
\item\leavevmode\label{compare-b+a-b}Pour $a,b\geq1$. On a pour tout $\alpha\in\S _a$,
$$\mathrigid1mu
\mathalign{
\hfill\chic(\Fba (\Xg))(\1_{b}\times\alpha,T)=
T^{b+a}{Q_1}\usp{b+\XX_1(\alpha)}\cdot R
\hfill\\\noalign{\kern4pt}
\chic(\Fg_{b}(\Xg\mmoins a))(\1_{b},T)\cdot\chic(\Fa (\Xg))(\alpha,T)=
T^{b+a}(Q_1-a)\usp{b}\cdot
{Q_1}\usp{\XX_1(\alpha)}\cdot R
}$$
où $Q_1:=({\chic(\Xg)(\1,T)/T})$ et 
$R:=\chic(\Fba (\Xg))(\alpha',T)/T^{b+a}$ avec $\alpha'$ tel que 
$\XX_1(\alpha')=0$ et $\XX_{i}(\alpha')=\XX_{i}(\alpha)$ pour les autres indices $i$.\\
Dans cette situation, la formule {\rm($\ddagger$)} est vérifiée si et seulement si
$\alpha=\1_a$.
\end{enumerate}
\end{prop}
\demo
(\ref{compare-b+a-a}) La fibre de $\pi$ étant $\Fg_{dr''}(\Xg\mmoins dr')$, on a par \ref{theo-trace-sr}
$$\begin{casesalign}
\hfill\chic(\Fg_{dr}(\Xg))(\sigma_{dr}^{r},T)&=&(dT^{d})^{r}Q_{d}(\Xg)\usp r\hfill\\\noalign{\kern4pt}
\hfill\chic(\Fg_{dr'}(\Xg))(\sigma_{dr'}^{r'},T)&=&(dT^{d})^{r'}Q_{d}(\Xg)\usp {r'}\hfill\\\noalign{\kern4pt}
\chic(\Fg_{dr''}(\Xg\mmoins dr'))(\sigma_{dr''}^{r''},T)&=&(dT^{d})^{r''}Q_{d}(\Xg\mmoins dr')\usp {r''}\,.\\
\end{casesalign}$$
où, dans la troisième égalité on a 
$Q_{d}(\Xg\mmoins dr')=Q_{d}(\Xg)$, d'après \ref{lemme-moins-a-points}.

\smallskip
(\ref{compare-b+a-b}) Les égalités résultent aussitôt de \ref{theo-trace-gen}. La formule ($\ddagger$) est vérifiée si et seulement si, 
$$Q_1{}\usp{b+\XX_1(\alpha)}
=(Q_1-a)\usp{b}\cdot
{Q_1}\usp{\XX_1(\alpha)}$$
et la conclusion résulte de ce que $Q_1{}\usp{b+\XX_1(\alpha)}={Q_1}\usp{\XX_1(\alpha)}\cdot (Q_1-\XX_1(\alpha))\usp{b}$.
\enddemo

\begin{prop}\label{comparaison-ok}Soit $\Xg$ un espace $i$-acyclique. Si $\alpha$ et $\beta$ sont immiscibles et si $\XX_1(\beta)=0$, l'égalité $(\ddagger)$ est vérifiée.
\end{prop}
\demo Conséquence du théorème  \ref{theo-trace-gen} et du lemme \ref{lemme-moins-a-points}. 
\enddemo
\section{Quotients d'espaces de configuration généralisés}
\glossarytitle{Nombres de Betti des quotients des espaces de configuration}
La\label{quotients} possibilité d'un algorithme de calcul des caractères de $\S _{m}$-modules $\Hc^{i}(\Delta_{\ell}\Xg^{m})$ donnée par le théorème \ref{theo-caracteres} et évoquée dans la remarque \ref{algo-caracteres}, ouvre la porte à la détermination des polynômes de Poincaré de la cohomologie des espaces d'orbites $(\Delta_{?\ell}\Xg^{m})/\Hg$ où $\Hg$ est un sous-groupe de $\S _{m}$.

\subsection{Polynômes de Poincaré de $(\Delta_{?\ell}\Xg^{m})/\Hg$}L'énoncé bien connu suivant rappelle le lien entre le polynôme de Poincaré des quotients et les séries de caractères.

\begin{prop}\label{Poincare-invariants}Soit $\Zg$ un espace topologique de type fini muni d'une action de $\S _{m}$. Pour tout sous-groupe $\Hg\dans\S _{m}$, on a
$$\mathalign{\Pc(\Zg/\Hg)(-T)&=&\smash{1\over |\Hg|}\sumnl_{h\in\Hg}\chic(\Zg)(h,T)\hfill\\\noalign{\kern4pt}
&=&{1\over m!}\sumnl_{g\in\S _{m}}\chic(\Zg)(g,T)\cdot \chi(\ind_{\Hg}^{\S _{m}}\1)(g)\hfill\\\noalign{\kern4pt}
&=&{1\over m!}
\sumnl_{\lambda}h_{\lambda}\,\chic(\Zg)(g_{\lambda},T)\cdot \chi(\ind_{\Hg}^{\S _{m}}\1)(g_{\lambda})\hfill
}
$$
Dans la dernière formule, la sommation est indexée par les décompositions $\lambda=\set 1^{\XX_1},2^{\XX_2},\ldots,m^{\XX_{m}}/$ de $m$, le nombre $h_{\lambda}$ est le cardinal de l'ensemble $\pi(\lambda)$ des permutations dont la décomposition en produit de cycles disjoints est de type $\lambda$, l'élément $g_{\lambda}$ est alors un représentant quelconque de $\pi(\lambda)$. 
\end{prop}
\demo Un théorème classique de Grothendieck (\cite{gro})\footnote{{\it Loc.cit. }Théorème 5.3.1 et corollaire de la Proposition 5.2.3.} donne l'équivalence 
$\Hc^{i}(\Zg/H)=\Hc^{i}(\Zg)^{H}$, pour tout $i\in\NN$. On a donc
$$\mathalign{\Pc(\Zg/\Hg)(-T)&=&\sum_{i\in\NN}\Big(
{1\over |\Hg|}\sumnl_{h\in\Hg}\tr(h\actson\Hc^{i}(\Zg))\Big) (-T)^{i}\\
&=&
{1\over |\Hg|}\sumnl_{h\in\Hg}\chic(\Zg)(h,T)\,.\hfill
}$$
Les autres égalités sont classiques (\cf\cite{mac} eq. (6.1) et (6.2)).\enddemo

\def\Remarque{Commentaire}\begin{rema}Lorsque l'espace $\Xg$ est $i$-acyclique, cette proposition et l'explicitation des séries de caractères des $\S _{m}$-modules $\Hc(\Fm (\Xg))$ du théorème \ref{theo-trace-gen}, donnent des formules très explicites des polynômes de Poincaré
pour les espaces de configurations cycliques et non-ordonnées associés à $\Fm (\Xg)$, ce qui constitue le sujet des deux sections suivantes.
\end{rema}

\subsection{Espaces de configurations cycliques $\Cg\Fm (\Xg)$}Soit $\Cg_{m}$ le sous-groupe de $\Sm$ engendré par le cycle $\sigma_{m}:=(1,\ldots,m)$. L'\expression{espace de configurations cycliques de $\Xg$}, noté $\Cg\Fm (\Xg)$, est l'espace\glossary{${\Cg\Fm (\Xg):=\Fm (\Xg)/\Cg_{m}}$:espace de configurations cycliques de $\Xg$}
$$
\Cg\Fm (\Xg):=\Fm (\Xg)/\Cg_{m}\,.$$

\begin{theo}\label{theo-conf-cycliques}Soit $\Xg$ un espace $i$-acyclique.
Pour tout $m\in\NN$, on a
$$\displayboxit{{\Pc(\Cg\Fm )(-T)}=\relax{T^{m}\over m}\sum_{d\div m}\phi(d)\,
d^{m/d}
\Big(\sum_{e\div d}\mobius de
{\Pc(\Xg)(-T^{e})\over dT^{e}}\Big)\usp {m/d}}
\,,$$
où $\phi(\_)$\glossary{$\phi(\_)$:fonction indicatrice d'Euler} est la fonction indicatrice d'Euler, $\mu(\_)$ est la fonction de Möbius et $(\_)\usp r$ est la factorielle décroissante de \ref{factorielles}.
\end{theo}
\demo Par \ref{Poincare-invariants}, on a
$${\Pc(\Cg\Fm )(-T)\over T^{m}}=\smash{1\over m}\sumnl_{dr= m}\phi(d)\,{\chic(\Fm )(\sigma^{r}_{m},T)\over T^{m}}\,.
$$
Par \ref{theo-trace-sr}, on a alors
$${\Pc(\Cg\Fm )(-T)\over T^{m}}=\smash{1\over m}\sumnl_{dr=m}\phi(d)\,
d^{r}
\Big({\chic(\Fg_d(\Xg))(\sigma_d,T)\over dT^{d}}\Big)\usp r\,,
$$
et la proposition résulte de l'égalité \ref{theo-trace-s}:
$$
\relax{{\chic(\Fg_{d}(\Xg)(\sigma_d,T))\over T^{d}}=\sumnl_{e\div d}\mobius de
{\Pc(\Xg)(-T^{e})\over T^{e}}}\,.$$
\vskip-2em
\enddemo

\begin{rema}Un cas particulier de la proposition précédente est lorsque $m$ est un nombre premier $p$. Dans ce cas, on a simplement
$$\mathrigid1mu
\mathalign{{\Pc(\Cg\Fg_{p})(-T)\over T^{p}}&=&\smash{1\over p}\Big(
{\chic(\Fg_{p})(\1,T)\over T^{p}}+(p-1){\chic(\Fg_{p})(\sigma_p,T)\over T^{p}}\Big)\hfill\\\noalign{\kern4pt}
&&\hskip-1cm{}=\mathrigid1mu
{1\over p}\bigg(\!\!\!\Big({\chic(\Xg)(\1,T)\over T}\Big)\usp p +(p-1)\Big({\chic(\Xg)(\1,T^{p})\over T^p}-{\chic(\Xg)(\1,T)\over T}\Big)\!\!\!\bigg)
}
$$
Où l'on remarquera que si nous notons $Q(T):=\chic(\Xg)(\1,T)/T$, la dernière expression entre les grandes parenthèses est
$$Q(T)\usp p+(p-1)\big(Q(T^{p})-Q(T)\big)\,,$$
qui est nulle modulo $p$ puisque l'on a
$$Q(T)\usp p\equiv_{p}Q(T)^{p}-Q(T)\equiv_{p} Q(T^p)-Q(T)\,.$$
\end{rema}

\begin{rema}En suivant la méthode décrite dans $\cite{mac}$, la proposition résulte également comme application de la troisième égalité de \ref{Poincare-invariants}. Dans ce cas, $\chi(\ind_{\Cg_{m}}^{\S _{m}})(g_{\lambda})\not=0$, si et seulement si, $g_{\lambda}$ est conjuguée à $\sigma_{m}^{m/d}$ pour un certain $d\div m$, auquel cas 
$$\def\\{\,,\quad}
\goodsmash{0.80}{0.80}{\lambda=\bigset {{d}{}^{m/d}}/\\
h_{\lambda}={m!\over d^{m/d}}\\
g_{\lambda}:=\sigma_{m}^{m/d}\\
\chi(\ind_{\Cg_{m}}^{\S _{m}})(g_{\lambda})=\phi\big({d}\big)\,{d^{m/d}\over m}\,,}
\postskip1ex$$
d'où
$$\preskip0pt\Pc(\Cg\Fm (\Xg))(-T)={1\over m}\sumnl_{d\div m}\phi(d)\,\chic(\Fm (\Xg))(\sigma_{m}^{m/d},T)\,,$$
et la conclusion suit par application de \ref{theo-trace-gen} à la décomposition $\lambda$ en question, ce qui donne
$$\preskip0pt\chic(\Fm (\Xg))(\sigma_{m}^{m/d},T)=T^{m}\ d^{m/d}\ 
\Big({\chic(\Fg_{d}(\Xg))(\sigma_{d},T)\over d T ^{d}}
\Big)\usp{m/d}\,.$$
\end{rema}

\begin{rema}[pour $\Xg^{m}$]\label{theo-prod-cycliques}En accord avec la remarque \ref{theo-trace-gen-Xm}, l'analogue du théorème \ref{theo-conf-cycliques} pour le produit cyclique $\Xg[m]:=\Xg^{m}/\Cg_{m}$~(\footnote{Notation de \cite{mac} p.~568.}) est donné par la même formule où la factorielle décroissante $(\_)\usp r$ est remplacée par la puissance $(\_)^r$ et où la deuxième sommation est restreinte au seul terme d'indice $e:=d$. On a donc :
$${\Pc(\Xg[m])(T)}=\relax{1\over m}\sumnl_{d\div m}\phi(d)\,
{\Pc(\Xg)((-1)^{d+1}T^{d})}{}\sp {m/d}
\,,$$
qui est très précisément la formule (8.4) de \cite{mac}.
\end{rema}

\subsection{Espaces de configurations non-ordonnées $\BFm (\Xg)$}On\label{BFm} appelle \expression{espace de configurations non ordonnées de $\Xg$}, noté $\BFm (\Xg)$, l'espace\glossary{${\BFm (\Xg):=\Fm (\Xg)/\Sm}$:espace de configurations non ordonnées de $\Xg$}
$$\BFm (\Xg):=\Fm (\Xg)/\Sm\,.$$
\begin{theo}Soit\label{theo-conf-symetriques} $\Xg$ un espace $i$-acyclique.
Pour tout $m\in\NN$, on a
$$\mathrigid1mu
\hss\displayboxit{{\Pc(\BFm (\Xg))(-T)\over T^{m}}=\decale-5pt{1\over m!}
\mkern-25mu
\vscalesum{1.75}{-1.3pt}_{\vrule height0pt width0pt\lambda:=(1^{\XX_1},\ldots,m^{\XX_{m}})\vdash m}
\mkern-30mu 
\decale3pt{h_{\lambda}}
\prodnl_{d=1}^{m}d^{\XX_d}
\bigg(\sum_{e\div d}\mobius de
{\Pc(\Xg)(-T^{e})\over dT^{e}}\bigg)\usp {\XX_{d}}}
\hss$$
où $\mu(\_)$ est la fonction de Möbius, $(\_)\usp r$ est la factorielle décroissante  (\ref{factorielles}), et $h_{\lambda}$ est le cardinal de l'ensemble des permutations de $\S _{m}$ dont la décomposition en cycles disjoints est de type $\lambda:=(1^{\XX_1},\ldots,m^{\XX_{m}})\vdash m$, soit:
$$h_{\lambda}={m!\over \XX_{1}!\,\XX_{2}!\,\ldots\XX_{m}!
\ (1!)^{\XX_1}(2!)^{\XX_2}\cdots(m!)^{\XX_{m}}}\,.$$
\end{theo}
\demo Corollaire immédiat de la troisième formule de \ref{Poincare-invariants} modulo le calcul de séries de traces de \ref{theo-trace-gen}. La valeur de $h_{\lambda}$ est classique et bien connue (\cf\cite{and} thm.~13.2, p. 215).
\enddemo

\begin{rema}[pour $\Xg^{m}$]\label{theo-prod-symetriques}Comme pour les autres remarques concernant $\Xg^{m}$, l'analogue du dernier théorème \ref{theo-conf-symetriques} pour le produit symétrique $\Xg(m):=\Xg^{m}/\S _{m}$~(\footnote{Notation de \cite{mac} p.~568.}) est donné par la même formule où la factorielle décroissante $(\_)\usp r$ est remplacée par la puissance $(\_)^r$ et où la deuxième sommation est restreinte au seul terme d'indice $e:=d$. On a donc :
$$\mathrigid1mu
\hss
\displayboxit{{\Pc(\Xg(m))(T)}=\decale-5pt{1\over m!}
\mkern-25mu
\vscalesum{1.75}{-1.3pt}_{\vrule height0pt width0pt\lambda:=(1^{\XX_1},\ldots,m^{\XX_{m}})\vdash m}
\mkern-30mu 
\decale3pt{h_{\lambda}}
\prodnl_{d=1}^{m}
{\Pc(\Xg)((-1)^{d+1}T^{d})}\sp {\XX_{d}}
}\hss$$
qui est très précisément la formule (8.3) de \cite{mac}.

Concernant cette formule, on rappelle que Macdonald en a donné une très belle fonction génératrice (\loccit\ eq (8.5)): le $k$-ième coefficient de $\Pc(\Xg(m))$ est le coefficient en $x^{k}t^{m}$ du développement en série entière de 
$$\def\terme#1#2#3{(1#1x^{#2}t)^{\beta_{#3}}}
{
\terme +{}1
\terme +33\cdots
\over
(1-t)^{\beta_0}
\terme -22
\terme -44
\cdots
}$$
où $\beta_{k}$ est le $k$-ième coefficient de $\Pc(\Xg)$.

Nous ne connaissons pas de résultat semblable pour $\Pc(\BFm (\Xg))$, la difficulté essentielle, par rapport à l'approche de Macdonald, réside dans les factorielles décroissantes.
\end{rema}

\def\Remarque{Commentaire}\begin{rema}Conformément\label{stabilite-Betti} au théorème  \ref{theo-stabilite-BM-i-acyclique}
de monotonie et stabilité des familles
\smash{$\set\Sm\rep\Hbm^{i}(\Fm(\Mg))/_{m}$}, lorsque $\Mg$ est une pseudovariété connexe orientable de dimension $\dMg\geq2$ et pour chaque $i\in\NN$ fixé, la famille de  polynômes de Poincaré
$$\bigset{\Pc(\BFm (\Mg))(1/T)\cdot T^{m\dMg}\ (\mathop{\rm mod}\ T^{i+1})}/_{m}\eqno(\ast)$$
est stationnaire. C'est un phénomène qui n'est pas apparent à la lecture de la formule  \ref{theo-conf-symetriques}. Cela nous a intrigué et nous a conduit à la rédaction de la section suivante où nous déterminons assez précisément le rang de stabilité de la famille de nombres $\set\Bettibm^{i}(\BFm(\Mg))/_{m}$. On y procède en deux temps. D'abord, lorsque $\Mg$ est $i$-acyclique (\ref{theo-stabilite-i-acyclique-betti-BFm}), en étudiant le coefficient de $T^{i}$ dans développement de la série $(1/m!)\sum_{\alpha\in\Sm}\chibm(\Fm (\Mg))(\alpha,T)$ d'après la formule  \ref{theo-trace-gen}. Ensuite,
lorsque $\Mg$ est général, à l'aide des suites spectrales basiques. Ce faisant, nous obtenons la proposition \ref{prop-stabilite-betti-dim3-BFm} qui est une importante
amélioration  de l'estimation de la plage de stabilité de la famille de nombres $(\ast)$ puisque l'on passe de $\set m\geq4i/$ à $\set m\geq2i/$ si $\dMg=2$,
et de $\set m\geq2i/$ à $\set m\geq i/$ si $\dMg\geq3$. Cette dernière proposition généralise les théorèmes A et B de Randal-Williams \cite{R-W}, mais aussi le corollaire 3 de Church (\cite{chu},~p.~470)  valables seulement lorsque $\Mg$ est une variété topologique.\footnote{Il convient cependant de signaler que pour $\Mg$ une variété topologique de dimension $d_{\Mg}=2$, Church améliore encore la plage de stabilité à $\set m>i/$ ce que nous n'avons pas établi, mais que nous n'excluons pas comme accessible avec nos méthodes.}
\end{rema}

\subsection{Sur le rang de stabilité de la famille $\set\Bettibm^{i}(\BFm(\Mg))/_{m}$}\label{stabilite-betti-BFm}
\subsubsectionline{Rang de stabilité des caractères polynomiaux.}Dans\label{rang-stabilite-caracteres} \ref{caracteres-polynomiaux}, nous avons rappelé la notion de polynomialité d'une famille de caractères $\set\chi_{\Sm}(W_m)/_{m}$. Nous y avons introduit la fonction $\XX_i:\Sm\to\NN$ qui fait correspondre à $\alpha\in\Sm$ le nombre $\XX_i(\alpha)$ de cycles de longueur $i$ dans la décomposition de $\alpha\in\S _{m}$ comme produit de cycles disjoints. 

\subsubsectionline{Polynomialité des fonctions centrales.}Pour $m\geq0$, on note $\QQ_{\rm cl}(\Sm)$\glossary{$\QQ_{\rm class}(\Sm)$:$\QQ$-algèbre des fonctions $f:\Sm\to\QQ$ \emph{centrales}, \idest 
constantes sur les classes de conjugaison.} la $\QQ$-algèbre des \emph{fonctions $f:\Sm\to\QQ$ dites centrales}, \idest constantes sur chaque classe de conjugaison de $\Sm$.
On note ensuite
$$\rho_m:k[\cl\XX]:=k[\XX_1,\XX_2,\ldots]\to\QQ_{\rm cl}(\Sm)$$
l'homomorphisme d'algèbres défini et faisant correspondre à $\XX_i$ la fonction centrale $\alpha\mapsto\XX_i(\alpha)$. On a donc
$$\rho_m(\XX_1^{d_{1}}\XX_2^{d_{2}}\cdots\XX_\ell^{d_{\ell}})(\alpha):=\prodnl_{k}\XX_k(\alpha)^{d_k}\,,\quad\forall d_k\in\NN\,.$$

\begin{prop}\label{redondance}L'homomorphisme de $\QQ$-algèbres 
$$\rho_m:k[\cl \XX  ]\to\QQ_{\rm cl}(\Sm)$$
est surjectif et son noyau contient les polynômes:
$$\mathrigid2mu
\big(\XX _1+2\XX _2+\cdots+m\XX _m-m\big)
\text{\quad and \quad}
\big(\XX _i(\XX _i-1)\cdots\big(\XX _i-\big\lfloor m/i\big\rfloor\big)\big)\,.
$$
En particulier, les caractères de $\SSS_m$ peuvent être réprésentés par des polynômes à coefficients rationnels en les variables $\XX _1,\ldots, \XX _{m-1}$. 
\end{prop}
\demo Pour $i\in\iii[1,m]$ et $j\in\iii[0,\lfloor m/i\big\rfloor]$
soit le polynome de $\ZZ[\Zsf]$:
$$\Zsf_{i,j}:=
{\Zsf(\Zsf-1)\cdots\widehat{(\Zsf-j)}\cdots\big(\Zsf-\big\lfloor m/i\big\rfloor\big)}
$$
et notons 
$$\Qsf_{i,j}:=\Zsf_{i,j}(\XX_i)/\Zsf_{i,j}(j)\in k[\cl\XX]\,.$$
On a clairement l'égalité
$$\rho_m(\Qsf_{i,j})(\alpha)=\delta_{\XX_{i}(\alpha),j}\hbox{ (delta de Kronecker)\hskip-2cm}
$$
En particulier, si $\mu=(1^{a_1},\ldots,m^{a_m})\vdash m$, d'image par $\rho_m$ du polynôme
$$\Qsf_\mu:=\prodnl_{i=1}^{m} \Qsf_{i,a_i}$$
est la fonction caractéristique de l'ensemble des permutations de type $\mu$. 
\enddemo

\subsubsectionline{La graduation $\ddeg$ de $k[\cl\XX]$.}Il s'agit de la graduation des polynômes dans l'algèbre  $k[\cl\XX]:=k[\XX_1,\XX_2,\ldots]$ qui résulte de poser $\ddeg\XX_i=i$\glossary{$\ddeg$:graduation de $k[\XX_1,\XX_2,\ldots]$ telle que $\ddeg\XX_k=k$}.

\begin{prop}Pour\label{prop-rang-stabilite-caracteres} 
$P\in k[\XX_1,\ldots,\XX_{\ell}]$, la famille de nombres
$$\Bigset{1\over m!}\sumnl_{\alpha\in\Sm}P\big(\XX_1(\alpha),\ldots,\XX_{\ell}(\alpha)\big)/{\vrule depth6pt width0pt}_{m}$$
est constante pour $m\geq\ddeg(P)$.
\end{prop}
\demo Il suffit de le prouver pour les monômes $P:=\XX_1^{d_{1}}\XX_2^{d_{2}}\cdots\XX_\ell^{d_{\ell}}$, où $\ddeg(P)>0$, vu que pour $P=1$ la proposition est évidente.

\medskip\noindent
{\slshape On introduit les ensembles $\C_{i}(m)$, pour $1\leq i\leq m$.}

On pose $\C_{1}(m):=\iii[1,m]$, et, pour $i>1$,
on pose $\C_{i}(m):={}$le sous-ensemble de $\Sm$ des $i$-cycles.
On fait agir $\Sm$ sur $\C_{1}(m)$ par son action naturelle et sur
$\C_{i}(m)$ par conjugaison, $\alpha\cdot \gamma:=\alpha\gamma \alpha^{-1}$. On note ensuite $\E_{i}(m)$ le $k[\Sm]$-module engendré par les éléments de $\C_{i}(m)$, et l'on note 
$$
\E_P(m):=
\E_{1}(m)^{\otimes d_1}\otimes\cdots\otimes
\E_{\ell}(m)^{\otimes d_\ell}\,.$$
muni de sa structure de $\Sm$-module produit tensoriel de représentations.

\medskip
\noindent{\slshape$(A)$ Rang de stabilité de $\dim_{k}(\smash{\E_{P}(m)^{\Sm})}$.}

\nobreak
L'espace $\E_{P}(m)$ a une base paramétrée par l'ensemble
de $\ddeg(P)$-uplets
$$\C_{P}(m):=\bigset(\gamma_{1,1},\ldots,\gamma_{1,d_1},
\gamma_{2,1},\ldots,\gamma_{2,d_2},
\ldots,
\gamma_{\ell,1},\ldots,\gamma_{\ell,d_\ell}
)/\,,$$
où $\gamma_{i,j}\in\C_i(m)$, de sorte que
le sous-espace $\E_{P}(m)^{\Sm}$ admet une base paramétrée par l'ensemble $\C_{P}(m)/\Sm$
des orbites de l'action $\Sm$ sur $\C_{P}(m)$.

L'inclusion canonique $\Sm\dans\Smm$ induit des inclusions 
$\C_{P}(m)\dans \C_{P}(\mm)$, et donc $\E_P(m)\dans\E_{P}(\mm)$,
compatibles aux actions des groupes symétriques. On vérifie alors aisément que les applications induites
$$\C_{P}(m)/{\Sm}\dans \C_{P}(\mm)/{\Smm}\text{\quad et\quad}
\E_{P}(m)^{\Sm}\dans \E_{P}(\mm)^{\Smm}\,,\eqno(\ast)$$
sont injectives.

Maintenant, pour $i>1$ et $\gamma\in\C_{i}(m)$, notons $|\gamma|$ le \expression{support de $\gamma$}, \ie l'ensemble des $x\in\iii[1,m]$ tels que $\gamma(x)\not= x$, et pour $i=1$ et $\gamma\in\C_1(m)$, notons plus simplement $|\gamma|=\set\gamma/$. Avec ces conventions, on définit le support de $\cl\gamma=(\gamma_{i,j})\in\C_{P}(m)$ par 
$\preskip0pt|\cl\gamma|:=\bigcup\nolimits_{i,j}|\gamma_{i,j}|\,,$
où la réunion est paramétrée par les couples $(i,j)$ tels que $d_{i}\geq 1$ et $1\leq j\leq d_{i}$.

Le support de $\cl\gamma\in\C_{P}(m)$ vérifie  $\#(|\cl\gamma|)\leq\ddeg(P)$. Lorsque $m\geq\ddeg(P)$, on peut trouver, pour $\cl\gamma$ donné, un élément $\alpha\in\Smm$ tel que $|\alpha\cdot\cl\gamma|\dans\iii[1,\ddeg(P)]$. On vérifie alors que les injections $(\ast)$ sont bijectives, et donc
$$\dim_{k}\big(\E_{P}(m)^{\Sm}\big)=\#(\C_{P}(m)/{\Sm})
\eqno(\diamond)$$
est indépendant de $m\geq\ddeg(P)$.

\medskip\goodbreak

\medskip\noindent{\slshape$(B)$ Formule de traces.}

\nobreak Notons $\YY_{i}:\Sm\to k$ la fonction centrale définie par
$$\YY_i(\alpha):=\tr(\alpha\Rep \E_{i}(m))\,.\postskip0pt$$
On a 
$$\dim_{k}\big(\E_{P}(m)^{\Sm}\big)={1\over m!}\sumnl_{\alpha\in\Sm} 
\YY_1(\alpha)^{d_1}\cdots\YY_\ell(\alpha)^{d_\ell}\,.\eqno(\diamond\diamond)$$
et la sommation de droite est bien indépendant de $m\geq\ddeg(P)$, d'après $(\diamond)$.

\medskip 

\noindent Nous devons obtenir la même conclusion en remplaçant les $\YY$ par des $\XX$.

\begin{itemize}
\item Pour $i=1$, on a $\YY_1(\alpha)=\tr(\alpha\Rep\E_1(m))=\C_{1}(m)^{\langle\alpha\rangle}=\XX_1(\alpha)$.

\item Pour $i\geq 2$, on a 
$$\YY_i(\alpha)=\tr(\alpha\Rep\E_i(m))=\C_{i}(m)^{\langle\alpha\rangle}=\#\set \gamma\in\C_{i}(m)\mid \alpha\gamma=\gamma\alpha/\,.$$

\smallskip
La condition $\alpha\gamma=\gamma\alpha$ implique que $\alpha$ laisse stable l'ensemble $|\gamma|$ qui va donc se décomposer en réunion disjointe de cycles de $\alpha$, tous de la même longueur puisque $\gamma$ agit transitivement sur son support. Il s'ensuit que les éléments de $\C_{i}(m)^{\langle\alpha\rangle}$ sont les $i$-cycles de $\Sm$ obtenus, pour chaque $d\div i$, en concaténant $i/d$ cycles distincts de $\alpha$ de longueur~$d$. On en déduit l'égalité
$$\YY_i(\alpha)=\phi(i)\,\XX_i( \alpha)+\sumnl_{d\div i,\ d<i}
\phi(d)\,\frac{d^{(i/d)}}{i}\,\XX_{d}(\alpha)^{\usp {i/d}}\,,\eqno(\ddagger)$$
où $\phi$ est l'indicatrice d'Euler.
Il convient de remarquer que dans le membre de droite tous les termes sont de degré $\ddeg$ majoré par $i$, et que ceux après le signe somme ne concernent que des variables $\XX_j$ avec $j<i$ avec un degré \emph{polynôme} strictement supérieur à $1$.
\end{itemize}

\medskip\noindent{\slshape$(C)$ Conclusion. }La partie $(A)$ nous dit que la proposition est vraie pour les monômes $\YY_{1}^{d_1}\cdots\YY_{\ell}^{d_\ell}$. Mais, si nous remplaçons chaque $\YY_i$ par son expression en termes des $\XX_i$ de $(\ddagger)$, nous avons une somme de termes de la forme
$$\YY_{1}^{d_1}\cdots\YY_{\ell}^{d_\ell}=
\ast\, \XX_{1}^{d_1}\cdots\XX_{\ell}^{d_\ell}
+\sum \ast\, \XX_{1}^{r_1}\cdots\XX_{\ell}^{r_\ell}
$$
avec  $\ddeg(\XX_{1}^{r_1}\cdots\XX_{\ell}^{r_\ell})\leq\ddeg(P)$, mais où le degré polynôme de $\XX_{1}^{r_1}\cdots\XX_{\ell}^{r_\ell}$ est strictement plus grand que celui de $P$.
Il est donc envisageable de raisonner par récurrence inverse sur le degré \emph{polynôme} des monômes $\XX_{1}^{d_1}\cdots\XX_{\ell}^{d_\ell}$, ce qui nous fait aboutir au monôme $X_1^{r}$, pour lequel la conclusion résulte de l'égalité $\XX_1=\YY_1$.
\enddemo

\subsubsectionline{La série de caractères $\chibm(\Yg)(\_,T)$.}De manière analogue à \ref{def-serie-caracteres}, nous définissons la série des caractères pour la cohomologie de Borel-Moore d'une pseudovariété $\Yg$ orientable (\ref{def-cohomologie-bm}) et munie d'une action de $\S _{m}$ (\ref{ad}-\ref{ad3}), par la formule 
$$\halfsmash{\chibm(\Yg)(\alpha,T)
=\sumnl_{i\in\ZZ}(-1)^{i}\, \chi_{\Hbm^{i}(\Yg)}(\alpha)\, T^{i}
\,.}$$

Nous avons montré dans \ref{def-action-image-inverse}-(\ref{def-action-image-inverse-b}) que si $\Mg$ est une pseudovariété orientable de type fini de dimension $\dMg$, on a une identification de $\Sm$-modules
$$\halfdisplayskips\relax{\Hbm ^{i}(\BFm(\Mg))=
\sgn(\_)^{\dMg}\otimes\Hc^{m\dMg-i}(\BFm(\Mg);k)^{\vee}\,,}
\postskip0pt$$ 
d'où
$$\halfdisplayskips
\def\Yg{\BFm(\Mg)}
\mathalign{\chibm(\Yg)(\alpha,T)
&=&\big(\sgn(\alpha)^{\dMg}\,(-T)^{m\dMg}\big)\,
\chic(\Yg)(\alpha^{-1},T^{-1})\,\hfill
}\,,
$$
avec 
$\mathalign{
\sgn(\alpha)^{\dMg}\,
(-T)^{m\dMg}=
\textstyle\prod_d (-T^{d})^{d_\Mg\XX_d(\alpha)}
}$ d'après un calcul élémentaire. 

\smallskip
Maintenant, si $\Mg$ est en plus $i$-acyclique, le théorème \ref{theo-trace-gen} et des manipulations simples  (\ref{rema-serie-caracteres}), nous conduisent à la formule explicite
$$\mathrigid2mu
\def\U{-T}
\displayboxit{
\chibm(\Fm (\Mg))(\alpha,T)=
\prod_{d=1}^{m}
\Big({(-T^{d})^{\dMg}\over d^{-1}\,T^{d}}\Big){\vrule height 9pt width0pt}^{\XX_d} 
\Big(\sum_{e\div d}\mobius de
{\Pc(\Mg)(-T^{-e})\over d\,T^{-e}}\Big)\usp {\XX_d}
}
$$
où `$\XX_d$' est un raccourci pour `$\XX_d(\alpha)$'. 


\subsectionline{Stabilité de $\set\Bettibm^{i}(\BFm(\Mg))/_{m}$, cas $i$-acyclique.} 
La proposition \ref{prop-stabilite-chi-bg} qui suit étudie la formule de $\chibm(\Fm (\Mg))(\alpha,T)$ indépendamment de toute théorie de représentation ou topologie. La formule dépend alors seulement de la donnée de $m$ et $\dMg\in\NN$ et de la suite des nombres de Betti compacts $\cg=(\cg_0,\ldots,\cg_\dMg)\dans\ZZ^{\dMg+1}$, on notera 
$$\chi(\cg)_{m}:=\chibm(\Fm (\Mg))(\alpha,T)\,.$$

\smallskip\noindent{\slshape\bfseries\boldmath Polynomialité de la famille $\chideg{i}(\cg)$. }\'Evaluée en $\alpha\in\Sm$, la formule de $\chi(\cg)_{m}$ donne le polynôme de $\ZZ[T,1/T]$ noté 
$$\chi(\cg)_{m}(\alpha)=\sumnl_{i\in\ZZ}\chideg{i}(\cg)_{m}(\alpha)\,T^{i}\,,
$$
d'où les familles de fonctions centrales 
$\chideg{i}(\cg)=\set \chideg{i}(\cg)_{m}:\Sm\to\ZZ/_{m}\,.$

\smallskip
On dira d'une famille de fonctions centrales
$\fgg=\set f_{m}:\Sm\to\ZZ/_{m}$
qu'elle est \expression{(éventuellement) polynomiale} s'il existe un polynôme $Q_{\fgg}\in\ZZ[\XX_1,\XX_2,\ldots]$ tel que $f_{m}(\alpha)=Q_{\fgg}(\alpha)$ pour tout $m$ (assez grand) et tout $\alpha\in\Sm$. Un tel $Q_{\fgg}$ est unique (\footnote{En\label{unicité-poly} effet, si $Q(\XX_1,\ldots,\XX_\ell)\in\ZZ[\cl \XX]$ est non nul, il existe $a_i\in\NN$ tels que $Q(a_1,\ldots,a_\ell)\not=0$. Il est alors facile de fabriquer des permutations $\alpha_{m}\in\Sm$, pour tout $m\geq\sum_{i} a_i\,i$, telles que $a_i=\XX_{i}(\alpha_{m})$, auquel cas $Q(\alpha_m)=Q(a_1,\ldots,a_\ell)$ et $Q$ n'est pas  éventuellement nulle.}) et le \expression{degré polynôme de $\fgg$} est alors celui de $Q_{\fgg}$. Si $f_m=Q_{\fgg}$ pour tout $m\geq N$, on dira 
que \expression{$\fgg$ est polynomiale sur $\set m\geq N/$}.

\def\vartitle{Critère de polynomialité}\begin{var}Soit\label{critere-poly} $\fgg(\XX_{i_1},\ldots,\XX_{i_\ell})$
une expression telle que pour tout $m$, la fonction
$f_{m}(\alpha):=\fgg(\XX_{i_1}(\alpha),\ldots,\XX_{i_\ell}(\alpha))$
soit bien définie sur $\Sm$ à valeurs dans $\ZZ$. Alors, pour que
la famille de fonctions centrales $\fgg=\set f_{m}:\Sm\to\ZZ/$ soit polynomiale sur $\set m\geq N/$, il est nécessaire que $\fgg(n_1,\ldots,n_\ell)$ soit un polynôme en les entiers $n_j\geq0$ tels que $\sum_j n_j\cdot i_{j}\geq N$.
\end{var}
\demo
Pour $\ng:=(n_{1},\ldots,n_{\ell})$, notons $|\ng|=\sum_j n_{j}\cdot i_{j}$ et choisissons $\sigma(\ng)\in\S_{|\ng|}$ vérifiant $\XX_{i_j}(\sigma(\ng))=n_j$, pour $1\leq j\leq \ell$. Alors, si $Q\in\ZZ[\XX_1,\XX_2,\ldots]$ est un polynôme tel que $f_m=Q$ sur tous les $\Sm$ tels que $m\geq N$, on a 
$$\fgg(n_1,\ldots,n_{\ell})=
f_{|\ng|}(\sigma(\ng))=Q(\sigma(\ng))=Q'(n_1,\ldots,n_{\ell})$$
où $Q'\in\ZZ[\YY_1,\ldots,\YY_\ell]$ est obtenu de $Q$ en posant $\YY_j:=\XX_{i_j}$ et en annulant les autres variables.
\enddemo

\smallskip
\noindent{\slshape\bfseries Factorielle décroissante modifiée. }On utilisera une extension de la notion de factorielle décroissante (\ref{factorielles}). 
Pour $a,b$ dans un anneau et $\XX\in\NN$, on note\glossary{${a\uspp{0}[b]=1\,,\quad\text{et}\quad
a\uspp \XX [b]:=a(a-h)(a-2h)\cdots(a-(\XX{-}1)h)}$:factorielle décroissante étendue}
$$a\uspp{0}[b]=1\,,\quad\text{et}\quad
a\uspp \XX [b]:=a(a-b)(a-2b)\cdots(a-(\XX{-}1)b)\,,
$$
On a alors: 
$
b^{\XX}a\uspp \XX[c]=(ba)\uspp \XX [bc]
$ et $
a\uspp {\XX+1} [b]=a\uspp {\XX} [b](a-\XX\, b)
\,.$

\comment
\medskip
D'autre part
$$
\Betti_{i}(\BFm(\Mg))={(-1)^{i}\over m!}\sum_{\alpha\in\Sm}{1\over i!}{\partial^{i}\over\partial T^{i}}\Big(\chibm(\Fm (\Mg))(\alpha,T)\Big)\vrule height15pt depth12pt{\vrule depth10pt width0pt}_{\,T=0}$$
Des formules qui, modulo le théorème de stabilité de caractères polynomiaux \ref{prop-rang-stabilite-caracteres}, réduisent la détermination du rang de stabilité de la famille $\set\Betti_{i}(\BFm(\Mg))/_{m}$, à celle du degré $\ddeg$ du terme en $T^{i}$ de la série
$$Q(\cl\XX,T):=\prod_{d=1}^{m}
\goodsmash{0.8}{0.5}{\Big({-}d
(-T)^{d(d_\Mg-1)}\Big){\vrule height 9pt width0pt}^{\XX_d} 
\Big({-}{1\over d}\sum_{e\div d}\mobius de
{\Pc(\Mg)(-T^{-e})\over -T^{-e}}\Big)\usp {\XX_d}}
\,.$$
\endcomment

\begin{prop}
\'Etant\label{prop-stabilite-chi-bg} donnés $\dMg\in\NN$ et 
$\cg=\set\cg_{0},\ldots,\cg_{\dMg}/\dans\ZZ$, on note
$P_{\cg}(T)=\cg_0+\cg_1T+\cdots+\cg_{\dMg}T^{\dMg}$ et 
$$\chi(\cg)_{m}(\cl\XX,T):=
\prod_{d=1}^{m}\Big({(-T^{d})^{\dMg}\over d^{-1}\,T^{d}}\Big){\vrule height 12pt width0pt}^{\XX_d} 
\Big(\sum_{e\div d}\mobius de
{P_{\cg}(-T^{-e})\over d\,T^{-e}}\Big){\vrule height 12pt width0pt}\usp {\XX_d}\,.$$
Alors, pour tout $i\in\NN$, la famille $\chideg{i}(\cg)$ est éventuellement polynomiale
 si et seulement si, $\cg_{\dMg}\in\set0,1/$. Dans ce cas, on a les propriétés suivantes de $\chi(\cg)$ suivant les valeurs de $\dMg$.
$${\def\minipage#1{\vtop{\hsize8cm\parindent0pt#1}}
\def\tt#1 #2 {\vrule height#1pt depth#2pt width0pt}
\def\ligne#1 #2 #3{\tt#1 #2 \tt 10 4 $\bullet$ #3}
\def\hline{\noalign{\hrule}}
\vbox{
\halign{\vrule\tt15 6 \quad\hfil$#$\hfil\quad\vrule&\quad\minipage{#}\hfil\quad\vrule\crcr
\noalign{\hrule}
\tt0 8 \dMg&$\hfil\chideg i(\cg_0,\ldots,\cg_{\dMg})$
\cr\hline
0&
\ligne0 0 {$\ddeg\chideg {-i}(\cg_0{=}1)\leq2i\,,$\hfill$\forall i\geq0$}\\
\ligne0 8 {$\chideg {-i}(\cg_0{=}0)_{m}=0\,,$\hfill$\forall m\,,\forall i\geq0$}
\cr\hline
1&
\ligne15 0 {$\chideg {0}(\cg_0{=}0,\cg_1)=(\cg_1)\,\uspp{\XX_1}[-1]\prod_{d>1}0\,\uspp{\XX_d}[-d]$\,,}\\
\ligne0 8 {$\chideg {i}(\cg_0{=}0,\cg_1)_{m}=0\,,$\hfill$\forall m\,,\forall i>0$}
\cr\hline
2&
\ligne0 0 {$\ddeg\chideg {i}(\cg_0,\cg_1,\cg_2{=}1)\leq2i\,,$\hfill$\forall i\geq0$}\\
\ligne0 8 {$\chideg {i}(\cg_0,\cg_1,\cg_2{=}0)_{m}=0\,,$\hfill$\forall m>i\geq0$}
\cr\hline
\geq3&
\ligne0 0 {$\ddeg\chideg {i}(\cg_0,\ldots,\cg_{\dMg}{=}1)\leq i\,,$\hfill$\forall i\geq0$}\\
\ligne0 11 {$\ddeg\chideg {i}(\cg_0,\ldots,\cg_{\dMg}{=}0)\leq {i\over \dMg{-}1}\,,$\hfill$\forall i\geq0$}
\cr
\hline}}}$$
\end{prop}
\demo \let\displayboxit\relax
Pour $d\geq 1$, notons
$$
\leftbrace{-0.5mm}{8.mm}\mathalign{
\displaystyle\Q_{d}(\cg)(\XX_d,T)&:=&
\Big({(-T^{d})^{\dMg}\over d^{-1}\,T^{d}}\Big){\vrule height 12pt width0pt}^{\XX_d} 
&
\Big(\sumnl_{e\div d}\mobius de\,
{P_{\cg}(-T^{-e})\over d\,T^{-e}}\Big)
{\vrule height 12pt width0pt}\usp {\XX_d}\,,\\
\noalign{\kern2pt}
\displaystyle
\hfill\R(\cg)_{d}(T)&:=&\Big({(-T^{d})^{\dMg}\over d^{-1}\,T^{d}}\Big) 
\hfill&
\Big(\sumnl_{e\div d}\mobius de\,
{P_{\cg}(-T^{-e})\over d\,T^{-e}}\Big)\,,\hfill
}
\eqnnb\label{Eq-Rd}
$$
de sorte que 
(\footnote{On omet d'indiquer `$\cg$' lorsqu'il est sous-entendu superflu.})
$$
\Q_d=(\R_{d})\uspp{\XX_d}[
(-1)^{d_\Mg}dT^{d(d_\Mg{-}1)}]\,.
\eqnnb\label{Eq-Rd-rec}
$$
Lorsque l'on évalue $\Q_{d}$ en $\alpha\in\Sm$, on a
$\XX_d(\alpha)\in\NN$ et le développent de la factorielle décroissante  génère une série 
$$\Q_{d}(\XX_d(\alpha),T)=\sumnl_{i\in\ZZ}\Qdeg{i}_{d}(\XX_d(\alpha))\,T^{i}$$
où \smash{$\Qdeg{i}_{d}$} définit la famille fonctions centrales 
\smashtop{$\set\Sm\ni\alpha\mapsto\Qdeg{i}_{d}(\XX_d(\alpha))\in\ZZ/_m$.}

\smallskip
Le terme $\R_{d}(\cg)$ dans \ref{Eq-Rd} vaut, suivant les cas $d=1$ ou $d>1$,
$$
\hskip-3.9mm
\leftbrace{2.5mm}{8mm}
\mathalign
{
\hfill\R_{1}(\cg)(T)&=&
\hskip-3mm\sum_{0\leq j\leq \dMg}
\cg_{j}
\,(-T)^{(\dMg{-}j)}\,,\hfill\\
\hfill\R_{d}(\cg)(T)&=&
\hskip-3mm
\sum_{
\scriptstack{
\vrule height7pt width0pt
0\leq j\leq \dMg\mrlap{\scriptstyle,\ j\not=1\,,\ e\div d}}}
\mobius de\ 
\cg_{j}(-1)^{d_\Mg-j}
\,T^{d(\dMg-1)-e(j-1)}\,.\hskip-7mm
}
\raiseeqnnb{15pt}\label{Eq-Rd-dev}
$$
Il s'agit de polynômes, en $T^{\pm}$, en $T$ si $\dMg>0$ et en $1/T$ si $\dMg=0$. Le coefficient du terme constant est $\cg_{\dMg}$, sauf si $\dMg=1$ et $d>1$ auquel cas il est nul. Ensuite, le premier terme non constant, de la forme $? T^{\epsilon\not=0}$, 
varie suivant les situations. Le calcul explicite à partir des formules \ref{Eq-Rd-dev} donne
\begin{itemize}
\item \noindent\rlap{Si $d=1$,}\hfill $(\forall\dMg)
\quad\hskip-2pt ? T^{\epsilon}=-\cg_{\dMg-1}\,T^{\pm}\,.$
\hfill\null
\item \leavevmode\rlap{Si $d>1$,}\hfill
$\preskip-1.5ex
\def\ss#1{\hbox to1em{\mathsurround0pt\hss$#1$}}
\let\strut\relax\left\{\ \mathalign{
\dMg=0,\ &?T^{\epsilon}&=&
\ss{}\cg_0\, T^{-d(p-1)/p}\hfill\\
\dMg=1,\hfill&?T^{\epsilon}&=&
\ss{-}\cg_0\,\mu(d)\, T\hfill\\
\dMg=2,\hfill&?T^{\epsilon}&=&
\ss{}\cg_2\, T^{d(p-1)/p}\hfill\\
\dMg\geq 3,\hfill&?T^{\epsilon}&=&
\ss{-}\cg_{\dMg-1}\, T^{d}\hfill\\
}\right.
$\hskip3mm\hfill\null

\end{itemize}
d'où les expressions
$$\def\ss#1{\hbox to2cm{#1,\hss}}
\let\strut\relax\left\{\ \mathalign{
\forall\dMg\ &\Q_{1}(\cg)&=&
\halfsmashtop{\big(?T^{\pm}+\cg_{\dMg}\big)\uspp{\XX_1}[(-1)^{\dMg}T^{(d_{\Mg}{-}1)}]}\hfill\\
\dMg=0,\ &\Q_{d}(\cg)&=&
\relax{\big(?T^{-d(p-1)/p}+\cg_0\big)\uspp{\XX_d}[d\,T^{-d}]}\hfill\\
\dMg=1,\hfill&\Q_{d}(\cg)&=&
\big(?T\big)\uspp{\XX_d}[-d]\hfill\\
\dMg=2,\hfill&\Q_{d}(\cg)&=&
\big(?T^{d(p-1)/p}+\cg_{2}\big)\uspp{\XX_d}[d\,T^{d}]\hfill\\
\dMg\geq 3,\hfill&\Q_{d}(\cg)&=&
\big(?T^{d}+\cg_{\dMg}\big)\uspp{\XX_d}[(-1)^{\dMg}d\,T^{d(\dMg{-}1)}]\hfill\\
}\right.\postdisplaypenalty10000
\eqnnb\label{Eq-Qds}$$
où $p$ désigne le plus petit facteur premier de $d$.

\smallskip
\noindent{\slshape$\bullet$ Conditions nécessaires à la polynomialité. }En dehors du cas $\dMg=1$ qui pose de problèmes particuliers, dans tous les autres cas, le critère de polynomialité \ref{critere-poly} appliqué aux expressions \ref{Eq-Qds} donne aussitôt la condition
$$\Big(\Qdeg i_{d}(\cg)\vcenter{\hsize3cm\footnotesize\advance\baselineskip-1pt\centering est éventuellement polynomiale\par}\Big)\Longrightarrow \Big(\cg_{\dMg}\in\set0,1/\Big)$$
En effet, c'est déjà clairement la condition pour que
$\Qdeg 0_{d}(\cg)=\cg_{\dMg}^{\XX_d}$ soit polynomial. Ensuite, 
par exemple dans les cas $\dMg\geq3$, le terme qui suit est
$\Qdeg d_{d}=?\cg_{\dMg}^{\XX_d{-}1}T^{d}$ où la même condition apparaît, etc\dots

\smallskip
\noindent{\slshape$\bullet$ Les cas $d_\Mg=1$. }On a $\cg=(\cg_0,\cg_1)$ et
$$\leftbrace{0.mm}{6mm}\mathalign{
\Q_1(\XX_1,T)&=&\big({-\cg_0}\,T+\cg_1\big)\uspp{\XX_1}[-1]\hfill\\\noalign{\kern-2pt}
\Q_d(\XX_d,T)&=&\Big({-\cg_0}\sumnl_{e\div d}\mobius de\, T^{e}\Big){\vrule height10pt width0pt}\uspp{\XX_d}[-d]\,,\text{ si $d>1$,}
}
$$
Le critère de polynomialité appliqué à $\Q_1$ nous emmène à regarder les produits de $n$ facteurs
$$\big({-\cg_0}\,T+\cg_1\big)\big({-\cg_0}\,T+\cg_1+1\big)\cdots\big({-\cg_0}\,T+\cg_1+(n-1)\big)\,.$$
Pour chaque $i>0$, le coefficient en $T^{i}$ du produit
est
$(-\cg_0)^{i}\big(\sum x_1\cdots x_{n-i}\big)$, où la somme concerne les $(n{-}i)$-uplets d'éléments deux à deux distincts de l'ensemble $\binspace0
\set \cg_1,\cg_1+1,\ldots,\cg_1+(n-1)/$. Une telle somme est donnée par un polynôme en $n$ de degré $\binspace0n-i+1$, le coefficient en question n'est pas polynomial, sauf évidemment si $\cg_0=0$ auquel cas on a
$$\leftbrace{-0.5mm}{6mm}\mathalign{
\Qdeg{0}_{1}(0,\cg_1)&=&\cg_1\uspp{\XX_1}[-1]
\hfill\\
\Qdeg{0}_{d}(0,\cg_1)&=&0\uspp{\XX_d}[-d]\hfill\\
}
\text{\quad et\quad }\Qdeg{i}_{d}(0,\cg_1)=0\,,\ \forall i\geq 1\,,\forall d\geq 1\,.
$$
Par conséquent, 
$$\leftbrace{-0.5mm}{6mm}\mathalign{\chideg{0}(0,\cg_1)&=&\cg_1\uspp{\XX_1}[-1]
\Big(\prodnl_{d>1}0\,\uspp{\XX_d}[-1]\Big)\\
\chideg{i}(0,\cg_1)&=&0\,,\quad \forall i>0\,.\hfill}$$


\smallskip
\noindent{\slshape$\bullet$ Les cas $d_\Mg\in\set0,2/$. }Pour $\cg$ respectivement $(\cg_0)$ et $(\cg_0,\cg_1,\cg_2)$ on a
$$\mathalign{
\hfill\Q_{d}(\cg_0)&=&\Big(\sumnl_{e\div d}\mobius de {\cg_0\over T^{d-e}}\Big)\uspp{\XX_d}[d/T^{d}]\,,\hfill\\\noalign{\kern-2pt}
\Q_{d}(\cg_0,\cg_1,\cg_2)&=&\Big(\sumnl_{e\div d}\mobius de 
\big({\cg_0\, T^{d+e}}+\cg_1\,T^{d}+{\cg_2\,T^{d-e}}\big)
\Big)\uspp{\XX_d}[d\,T^{d}]\,,
}$$
où l'on constate aussitôt l'égalité
$$\Q_{d}(\cg_0)(\XX_d,1/T)=\Q_{d}(0,0,\cg_0)(\XX_d,T)\,,\eqnnb\label{Eq-dMg-0<=>2}$$
qui réduit le cas $\dMg=0$ à un cas particulier de $\dMg=2$.

Pour $\dMg=2$, on a vu dans \ref{Eq-Qds} que l'on a les expressions
$$\leftbrace{-0.5mm}{6mm}\mathalign{
\Q_{1}(\cg)&=&\big(?\,T+\cg_2\big)\uspp{\XX_1}[T]\hfill\\\noalign{\kern-2pt}
\Q_d(\cg)&=&\big(?\,T^{d((p-1)/p)}+\cg_2\big)\uspp{\XX_d}[dT^{d}]\,,\text{ si $d>1$.}
}
\eqnnb\label{Qd-dM=2-b2=1}
$$
où seuls les cas $\cg_{2}\in\set0,1/$ nous intéressent.

\begin{itemize}
\item Pour $\cg_2=1$, on a  $\Qdeg{0}_{d}=1$. Ensuite, les propriétés élémentaires des factorielles décroissantes dans les expressions \ref{Qd-dM=2-b2=1} justifient les relations 
$$
\let\quad\relax
\mathalign{
\Q_1(\XX_1,T)-\Q_1(\XX_1{-}1,T)&=&
\Q_1(\XX_1{-}1,T)\,(?-(\XX_1{-}1))\, T\hfill\\
\Q_{d}(\XX_d,T)-\Q_{d}(\XX_d{-}1,T)&=&
\Q_{d}(\XX_d{-}1,T)\,(?\,T^{d((p-1)/p)}-?(\XX_d{-}1)\,T^{d})\,,}
$$
dont on déduit, pour tout $i\geq 1$,
$$\def\sss#1 {\hbox to1.5em{\hfill(#1)}}
\mathalign{
\sss i \ddeg\Qdeg{i}_{1}&\leq&2&+&\ddeg\Qdeg{i{-}1}_{1}\hfill\\
\sss ii \ddeg\Qdeg{i}_{d}&\leq&d&+&\sup\Bigset
\ddeg\Qdeg{i-d((p-1)/p)}_{d},\ 
d+\ddeg\Qdeg{i-d}_{d}/\,.}
$$
On peut alors prouver par induction sur $i$ que l'on a
$$
\mathalign{
\ddeg\Qdeg{i}_{d}(\cg_0,\cg_1,\cg_2{=}1)&\leq& 2i\,,&\quad\forall i\geq0\,,\forall d\geq1\,.
}\eqnnb\label{Eq-dMg=2-b2=1}
$$
En effet, pour $i=0$ l'assertion est claire puisque $\Qdeg 0_{d}=1$. Ensuite, si l'on suppose $\ddeg\Qdeg{i{-}1}_{1}\leq2(i{-}1)$, l'égalité (i) donne
$\ddeg\Qdeg i_{1}\leq 2i$. Pour $\Qdeg i_{d}$, on remarque que 
$2i-2d(p-1)/p\leq 2i-d\,,
$ puisque $1/2\leq(p-1)/p$, et, à l'aide de (ii), on tire
$$\ddeg\Qdeg i_{d}\leq d+\sup\bigset 2i-d,\ d +2(i-d)/=2i\,,$$
ce qui prouve \ref{Eq-dMg=2-b2=1}. Par conséquent,
$$
\displayboxit{
\mathalign{\ddeg\chideg{i}(\cg_0,\cg_1,\cg_2{=}1)&\leq& 2i\,,\quad\forall i\geq0\,,
}}$$
puisque
$\chideg{i}(\cg)=\sum_{i_1+i_2+\cdots=i}\Qdeg{i_1}(\cg)_{1}\,\Qdeg{i_2}(\cg)_{2}\cdots$. 
On en déduit
$$
\displayboxit{
\mathalign{
\ddeg\chideg{-i}(\cg_0{=}1)&\leq& 2i\,,\quad\forall i\geq0\,,}}$$
grâce à l'égalité \ref{Eq-dMg-0<=>2}.

\item Lorsque $\cg_2=0$, on a
$$\binspace0
\leftbrace{-0.5mm}{13mm}\mathalign{
\Q_1(\XX_1,T)&=&\big(\cg_{0}\,T^{2}+\cg_1\,T\big)\uspp{\XX_1}[T]=
(-T)^{\XX_1}(-\cg_0\,T-\cg_1)\uspp{\XX_1}[-1]\hfill\\\noalign{\kern4pt}
\Q_d(\XX_d,T)&=&\Big((-T^{d})\sumnl_{e\div d}\mobius de (-\cg_0)\,T^{e}\Big)\uspp{\XX_d}[dT^{d}]\,,\text{ si $d>1$,}\hfill\\
&=&\big(-T^{d}\big)^{\XX_{d}}
\Big(\sumnl_{e\div d}\mobius de (-\cg_0)\,T^{e}\Big)\uspp{\XX_d}[-d]
\hfill}
$$
ce qui montre que pour tout $m$, la fonction centrale
$$\chi(\cg_0,\cg_1,\cg_2{=}0):\Sm\to\ZZ[T]\postskip0pt$$
 coïncide avec
$$(-T)^{m}\chi(\cg_0,-\cg_1)
:\Sm\to\ZZ[T]\,.$$
En particulier, 
$$\displayboxit{\chideg{i}(\cg_0,\cg_1,\cg_2{=}0)_{m}=0\,,\quad\forall m> i\geq0\,,}$$
alors que $\chi(\cg_0{=}0)=\prod_{d\geq1}0\uspp{\XX_d}[d/T^{d}]$ et donc
$$\displayboxit{\chideg {-i}(\cg_0{=}0)_{m}=0
\,,\quad\forall m\geq1\,,\forall i\geq0\,.}$$
\end{itemize}

\smallskip
\noindent{\slshape$\bullet$ Les cas $d_\Mg\geq3$. }On suit la même démarche que pour le cas $\dMg=2$ en se restreignant (donc) aux cas où $\cg_{\dMg}\in\set0,1/$.
\begin{itemize}
\item Pour $\cg_{\dMg}=1$ et tout $d\geq1$, on a $\Qdeg{0}_{d}=1$, et grâce à \ref{Qd-dM=2-b2=1}, la relation
$$
\let\quad\relax
\mathalign{
\Q_{d}(\XX_d,T)-\Q_{d}(\XX_d{-}1,T)&=&
\Q_{d}(\XX_d{-}1,T)\,(?\,T^{d}-?(\XX_d{-}1)\,T^{d(d_{\Mg}{-}1)})\,,}
$$
dont on déduit, pour tout $i\geq 1$, la majoration
$$
\def\sss#1 {\hbox to1.5em{\hfill(#1)}}
\mathalign{
\ddeg\Qdeg{i}_{d}&\leq &d&+&\sup\Bigset
\ddeg\Qdeg{i-d}_{d},\ 
d+\ddeg\Qdeg{i-d(d_{\Mg}-1)}_{d}/\,,}
\eqnnb\label{Eq-dMg=3-b2=1}$$
où $d(\dMg{-}1)\geq2d$, ce qui permet de prouver par induction sur $i$
$$
\mathalign{
\ddeg\Qdeg{i}_{d}(\cg_{d_\Mg}{=}1)&\leq& i\,,&\quad\forall i\geq0\,,\forall d\geq1\,.
}
$$
et donc
$$
\displayboxit{
\mathalign{\ddeg\chideg{i}(\cg_{\dMg}{=}1)&\leq& i\,,\quad\forall i\geq0\,,
}}$$

\item Lorsque $\cg_{\dMg}=0$, la majoration \ref{Eq-dMg=3-b2=1} se simplifie en
$$
\def\sss#1 {\hbox to1.5em{\hfill(#1)}}
\mathalign{
\ddeg\Qdeg{i}_{d}&\leq \sup\Bigset
\ddeg\Qdeg{i-d}_{d},\ 
d+\ddeg\Qdeg{i-d(d_{\Mg}-1)}_{d}/\,,}
$$
ce qui permet de prouver par induction sur $i$
$$
\mathalign{
\ddeg\Qdeg{i}_{d}(\cg_{\dMg}{=}0)&\leq& 
\smashbot{{i\over\dMg{-}1}}\,,&\quad\forall i\geq0\,,\forall d\geq1\,.
}
$$
et donc
$$\preskip0pt
\displayboxit{
\mathalign{\ddeg\chideg{i}(\cg_{\dMg}{=}0)&\leq& 
\smashtop{{i\over\dMg{-}1}}\,,\quad\forall i\geq0\,,
}}$$
ce qui termine la preuve de la proposition.\enddemobox
\end{itemize}
\enddemo

\begin{coro}
Soit\label{theo-stabilite-i-acyclique-betti-BFm} $\Mg$ une pseudovariété $i$-acyclique de type fini telle que $\dim\Hc^{\dMg}(\Mg;\QQ)\leq1$. 
Pour $i\in\NN$, la famille $\set\Bettibm^{i}(\BFm(\Mg;\QQ))/_{m}$ est constante pour $m\geq i$. On a aussi,
\begin{enumerate}
\item si\label{theo-stabilite-i-acyclique-betti-BFm-a} $\dMg=1$, constance sur $\set m\geq 1/$;
\item si\label{theo-stabilite-i-acyclique-betti-BFm-b} $\dMg=2$ et $\cg_2=0$, nullité sur $\set m>i/$;
\item si\label{theo-stabilite-i-acyclique-betti-BFm-c} $\dMg\geq3$ et $\cg_{\dMg}=0$, constance sur $\bigset m\geq {i\over\dMg{-}1}/$.
\end{enumerate}
\end{coro}
\demo On applique les estimations de \ref{prop-stabilite-chi-bg} avec $\cg_\dMg\in\set0,1/$. Pour $\dMg=1$, seul le cas $i=0$ et $\cg_0{=}0$ est à considérer. On a pour tout $m\geq1$
$$\mathalign{{1\over m!}\sum_{\alpha\in\Sm}\chideg {0}(\cg_0{=}0,\cg_1)(\alpha)&=&{1\over m!}\sum_{\alpha\in\Sm}(\cg_1)\,\uspp{\XX_1}[-1]\prodnl_{d>1}0\,\uspp{\XX_d}[-d](\alpha)\\\noalign{\kern4pt}
&=&\let\strut\relax\left\{\mathalign{1\text{ si $\cg_1=1$}\\0\text{ si $\cg_1=0$}}\right.\hfill}$$

Les estimations de \ref{prop-stabilite-chi-bg} pour le degré $\ddeg$ des familles de fonctions centrales polynomiales $\chideg i(\cg)$ jointes à la proposition \ref{prop-rang-stabilite-caracteres} justifient les autres assertions à l'exception du cas $\dMg=2$. Dans ce cas, le procédé donne la constance seulement sur $\set m\geq 2i/$. Pour aller plus loin, on remarque que $\chibm\minideg i(\Mg)$ seul dépend des nombre de Betti compacts de $\Mg$. Le cas $\Mg=\CC$ est connu depuis Arnold (\cite{arnold2}) qui montre que la cohomologie rationnelle de $\Bg\Fm(\CC)$ est nulle en degrés $i\geq2$ et est de dimension constante $1$ sur $\set m\geq i/$ pour $i\in\set0,1/$, ce qui règle le cas $i$-acyclique avec $\dMg=2$ et $\cg_1=0$. Le cas $\cg_1>0$ résulte de considérer la suite exacte courte de $\Sm$-modules de \ref{moins-un-point-comparaison}
$$0\to\Hc\big(\Fm ^{\bullet}(\Mg)\big)[-1]\to\Hc(\Fm (\oo\Mg))\to
\Hc(\Fm (\Mg))\to0\eqno(\ast)$$
où $\oo\Mg:=\Mg\mmoins\set \bullet/$ pour un certain $\bullet\in\Mg$, et où $\Fm ^{\bullet}(\Mg)$ est la réunion disjointe des ouverts $\Fm ^{j=\bullet}(\Mg):=\set \cl x\in\Fm  \mid x_i=\bullet/$ où $j=1,\ldots,m$. 
On remarquera ici que $\Hbm^{k}(\Mg)$ et $\Hbm^{k}(\oo\Mg)$ seul diffèrent en degré $d_{\Mg}{-}1$ où 
$$\dim\Hbm^{d_{\Mg}{-}1}(\oo\Mg)=\dim\Hbm^{d_{\Mg}{-}1}(\Mg){+}1\,.$$
En dualisant $(\ast)$ on obtient:
$$0\to
\Hbm^{i}\big(\Fm(\Mg)\big)\to
\Hbm^{i}\big(\Fm(\oo\Mg)\big)\to
\ind_{\Smmo}^{\Sm}\Hbm^{i-(d_{\Mg}{-}1)}\big(\Fmmo(\oo\Mg)\big)\to
0$$
et donc
$$\Betti^{i}(\BFm(\oo\Mg))=\Betti^{i}(\BFm(\Mg))+\Betti^{i-(d_{\Mg}{-}1)}(\BFmm(\oo\Mg))\,.$$

Lorsque $d_{\Mg}=2$, ces remarques montrent que 
$\set\Betti^i(\BFm(\CC\mmoins \cg_1))/_m$ est constante pour $m\geq i$, alors que $\Betti^{1}(\CC\mmoins\cg_1)=\cg_1$. Ceci termine la justification du cas $\dMg=2$, et achève la preuve du corollaire.\enddemo

\subsection{Stabilité de $\set\Bettibm^{i}(\BFm(\Mg))/_{m}$, cas général}\label{stabilite-Betti-BFm-cas-general}
On étend le corollaire \ref{theo-stabilite-i-acyclique-betti-BFm} au cas des pseudovariétés générales $\Mg$ orientables. Nous procédons comme dans la preuve du théorème \ref{theo-stabilite-BM-pseudo}, à l'aide de la suite spectrale basique
$$\halfdisplayskips\mathalign{
\IEs(\U^{m})_{1}^{p,q}&=&
\bigoplus_{\hskip-4mm
\tau\in\TTT(\pp,m)
\hskip-4mm}
\ind^{\Sm}_{\HHH_{\utau}\times\S_{m-|\utau|}}
\alt^{d_{\Mg}}\otimes\Hbm^{Q}(\Fg_{\tau}(\Mg_{>0}))\ \Rightarrow\ 
\Hbm^{i}(\Fm(\Mg))}
\postdisplaypenalty10000$$
avec $Q=i-(d_{\Mg}{-}1)(m{-}(\pp))$.

La compatibilité entre les morphismes $p_{m}^{*}:\Hbm^{i}(\Fm(\Mg))\to\Hbm^{i}(\Fmm(\Mg))$ et ceux des suites spectrales basiques (\ref{prop-tableau-normaux-ss-basiques}-(\ref{prop-tableau-normaux-ss-basiques-b}))
$$\IEs(q_m^{*})_{1}^{p,q}:\IEs(\U^{m})^{p,q}_{1}\to\IEs(\U^{\mm})^{\pp,q}_{1}\,,\eqno(\dagger)$$
y compris avec les actions des groupes symétriques, ramènent la question de la stabilité de la famille $\set\Bettibm^{i}(\BFm(\Mg))/_{m}$, à celle de la stabilité des morphismes de suites spectrales de co-invariants (\ref{co-invariants}) induits par $(\dagger)$, \idest
$$\varPhi_{0}(\IEs(q_m^{*})_{1}^{p,q}):
\varPhi_0(\IEs(\U^{m})^{p,q}_{1})\to\varPhi_0(\IEs(\U^{\mm})^{\pp,q}_{1})\,.$$
On est alors emmené à trouver le rang de stabilité des familles de morphismes 
$$
\bigoplus_{\hskip-5mm
\tau\in\TTT(\pp,m)
\hskip-5mm}
\Hbm^{Q}(\Fg_{\tau}(\Mg_{>0}))_{\HHH_{\utau}\times\S_{m-|\utau|}}
\hf{\sum p_{\tau}^{*}}{}{1cm}
\bigoplus_{\hskip-7mm
\tau\in\TTT(p+2,\mm)
\hskip-7mm}
\Hbm^{Q}(\Fg_{\tau}(\Mg_{>0}))_{\HHH_{\utau}\times\S_{\mm-|\utau|}}
\eqno(\diamond)$$
définis par les morphismes de $\HHH_{\utau}$-modules
$$
\varPhi_{\uell}(p_{\tau}^{*}):
\varPhi_{\uell}(\Hbm^{Q}(\Fg_{\tau}(\Mg_{>0}))
\to
\varPhi_{\uell}(\Hbm^{Q}(\Fg_{\tau^{\bullet}}(\Mg_{>0})))
\eqno(\diamond\diamond)$$ 
où $p_{\tau}:\Fg_{\tau^{\bullet}}(\Mg_{>0})\to\Fg_{\tau}(\Mg_{>0})$ est la projection sur les $\pp$ premières coordonnées, et où $\uell$ est tel que $m{-}|\utau|=(\pp){-}\uell$,
conformément à la notation \ref{rema-lambda/m}-(\ref{rema-lambda/m-b}). L'étude de cette question conduit à l'énoncé suivant.

\begin{prop}Soit\label{prop-stabilite-betti-dim3-BFm} $\Mg$ une pseudovariété connexe orientée. Pour tout $i\in\NN $ fixé, le morphisme 
$$\varPhi_0(p_{m}^{*}):\varPhi_0(\Hbm^{i}(\Fm(\Mg)))\to
\varPhi_0(\Hbm^{i}(\Fmm(\Mg)))$$
 est bijectif pour $m\geq 2i$ si $d_{\Mg}=2$, et pour $m\geq i$ si $d_{\Mg}\geq 3$.
\end{prop}
\demo Soit $\epsilon:=2$ si $d_{\Mg}=2$ et $\epsilon:=1$ si $d_{\Mg}\geq3$.
D'après les remarques préliminaires, la proposition résultera de montrer que les morphismes $\sum p^{*}_{\tau}$ dans $(\diamond)$ sont des isomorphismes pour $m\geq \epsilon i$, ce qui découle clairement de ce que pour tout $m\geq \epsilon i$,
\def\varlistskips{\topsep0cm\itemsep0.5pt\mou\parskip0pt\mou}
\begin{enumerate}
\item le\label{dem-stabilite-betti-dim3-BFm-a} morphisme $\varPhi_{\uell}(p^{*}_{\tau})$
dans $(\diamond\diamond)$ est un isomorphisme;
\item l'application\label{dem-stabilite-betti-dim3-BFm-b} $(\_)^{\bullet}:\TTT(\pp,m)\to\TTT(p{+2},\mm)$ est bijective pour $Q\geq 0$.
\end{enumerate}

\smallskip
(\ref{dem-stabilite-betti-dim3-BFm-a}) Commençons la généralisation suivante du corollaire \ref{theo-stabilite-i-acyclique-betti-BFm}.

\noindent{\slshape Lemme. Soit $\Xg$ est une pseudovariété $i$-acyclique connexe orientée de dimension $d_{\Xg}\geq2$. Pour tout $t\geq 0$, le morphisme induit
$$\varPhi_{t}(p_{m}^{*}):
\varPhi_{t}(\Hbm^{i}(\Fm(\Xg)))
\to
\varPhi_{t}(\Hbm^{i}(\Fmm(\Xg)))$$
est bijectif pour $m\geq i+t$.}

\begingroup
\parskip0.5ex
\medskip\noindent{\sl Preuve. }On raisonne par induction sur $t\in\NN$ puis, pour chaque $t$ par induction sur $i\in\NN$.
Le cas $t=0$ est réglé pour tout $i\in\NN$, par le corollaire \ref{theo-stabilite-i-acyclique-betti-BFm} modulo le théorème \ref{theo-stabilite-BM-i-acyclique} qui donne le rang de monotonie $m\geq i$. Supposons donc $t>0$ et le lemme vérifié pour $(t{-}1)$.
En dualisant la suite exacte du théorème de scindage \ref{theo-scindage}-(\ref{theo-scindage-a}) et en appliquant le foncteur $\varPhi_{t}(\_)$ (possible puisque $t>0$), on a la suite exacte courte
$$\mathalign{0\to
\bigoplus_{a+b=i}\Hbm^{a}(\Xg)\otimes&
\varPhi_{t-1}(\Hbm^{b}(\Fmmo))\to
\varPhi_{t}(\Hbm^{i}(\Fm))\hfill\\\noalign{\kern-8pt}
&\hskip2mmm{}\to
\varPhi_{t-1}(\ind^{\S_{m-1}}_{\S_{m-2}}\Hbm^{i-(d_{\Xg}-1)}(\Fmmo))
\to
0}
\eqno(\ddagger)$$
qui montre que la bijectivité de $\varPhi_{t}(p_{m}^{*})$  au niveau du terme central est subordonnée par celle des termes extrêmes. 

Lorsque $i=0$, le terme de droite de $(\ddagger)$ est nul puisque $d_{\Xg}>1$ et le lemme résulte des terme de gauche $\varPhi_{t-1}(\Hbm^{b}(\Fg_{m-1}(\Xg)))$ par l'hypothèse inductive.

Lorsque $i>0$, le terme de droite n'est pas forcément nul.
On remarque alors que la preuve du lemme~1~(p.~\pageref{theo-Ind-lambda-Lemme-1}) dans la démonstration du théorème \ref{theo-Ind-lambda}, donnait déjà l'identification
$$\varPhi_{t-1}(\ind^{\S_{m-1}}_{\S_{m-2}}(\_))
\ \simeq\ \
\bigoplusnl_{\cl \alpha\in\big(\biquosb {\S_{m-2} }{\S_{m-1} }{\S_{m-t} }\big)}
(\_)_{\S_{\iii[2,m-1]\cap\alpha \iii[t,m-1]}}\,,
\postskip-0.2ex$$
où 
$$\big|\iii[2,m-1]\cap\alpha \iii[t,m-1]\big|=(m{-}1){-}t'\,,\quad\text{avec $t'\in\iii[t-1,t]$.}$$ On y montrait aussi que l'application naturelle
$$(\_)^{\bullet}_{m}:
\Big(\biquo {\S_{m-2} }{\S_{m-1} }{\S_{m-t} }\Big)
\too
\Big(\biquo {\S_{m-1} }{\S_{m} }{\S_{m-(t-1)} }\Big)
\,,$$
où $\alpha^{\bullet}=\alpha$ sur $\iii[1,m{-}1]$ (donc $\alpha^{\bullet}(m)=m$), est bijective dès que
$m\geq t+1\,.$

De ces remarques, suit que sur les termes de droite de $(\ddagger)$, l'application
$$
\varPhi_{t-1}(\ind^{\S_{m-1}}_{\S_{m-2}}(\Hbm^{i-(d_{\Xg}-1)}(\Fmmo)))\to
\varPhi_{t-1}(\ind^{\S_{m}}_{\S_{m-1}}(\Hbm^{i-(d_{\Xg}-1)}(\Fm)))\eqno(\ast)$$
sera bijective dès lors que 
$$m\geq i+t\,.\eqno(\diamond)$$
En effet, pour que $(\ast)$ soit bijective, il suffit de vérifier que
\def\varlistskips{\itemsep0pt\mou\parsep0pt\parskip0pt\topsep2pt
\def\theenumi{\roman{enumi}}
\def\labelenumi{{\rm\theenumi)}}}
\begin{enumerate}
\item l'application\label{prop-stabilite-betti-dim3-BFm-lemme-a} $(\_)^{\bullet}_{m}$ est bijective, donc que $m\geq t+1$;

\item l'application\label{prop-stabilite-betti-dim3-BFm-lemme-b} $\varPhi_{t'}(\Hbm^{i-(d_{\Xg}-1)}(\Fmmo))\to
\varPhi_{t'}(\Hbm^{i-(d_{\Xg}-1)}(\Fm))$ est bijective pour $t'\in\iii[t-1,t]$, donc, que $m{-}1\geq i{-}(d_{\Xg}{-}1){+}t'$, par hypothèses inductives puisque $i{-}(d_{\Xg}{-}1)<i$ et $t'\leq t$;
\end{enumerate}
et comme $i\geq 1$, la condition $(\diamond)$ suffit pour avoir 
(\ref{prop-stabilite-betti-dim3-BFm-lemme-a},\ref{prop-stabilite-betti-dim3-BFm-lemme-b}), ce qui termine la preuve du lemme.
 \hfill$\boxminus$
\endgroup

\smallskip
Revenons à  la preuve de l'assertion (\ref{dem-stabilite-betti-dim3-BFm-a}). Lorsque $d_{\Mg}\geq 1$, le lemme s'applique aux morphismes $(\diamond\diamond)$
$$
\varPhi_{\uell}(p_{\tau}^{*}):
\varPhi_{\uell}(\Hbm^{Q}(\Fg_{\tau}(\Mg_{>0})))
\to
\varPhi_{\uell}(\Hbm^{Q}(\Fg_{\tau^{\bullet}}(\Mg_{>0})))
$$ 
et montre que
$\varPhi_{\uell}(p_{\tau}^{*})$ est bijective sous les conditions
$$\bigpar (p{+}1)\geq (Q{+}\uell)
/\text{\quad ou bien\quad }\bigpar Q< 0/
\,.\eqno(\ddagger)$$

En substituant (\cf notations  \ref{rema-lambda/m}-(\ref{rema-lambda/m-b}))
$$\binspace1\pmathalign{
Q&\leftrightarrow&i-(d_{\Mg}-1)(m-(\pp))\,,\\
\uell&\leftrightarrow&|\utau|-(m-(\pp))\,,\hfill\\
m{-}(\pp)&\leftrightarrow&|\utau|-\uell\,,\hfill}$$ 
la première condition dans $(\ddagger)$ s'écrit aussi
$$
\bigpar m\geq i-(d_{\Mg}{-}2)|\utau| + (d_{\Mg}{-}1)\uell/
\,,\eqno(\ddagger\ddagger)$$
et comme $|\utau|\geq2\uell$, on voit aisément que si $d_{\Mg}=2$, l'intervalle $\set m\geq 2i/$ est conforme à  $(\ddagger\ddagger)$, tandis que si $d_{\Mg}\geq3$, c'est le cas de $\set m\geq i/$.

\smallskip
(\ref{dem-stabilite-betti-dim3-BFm-b}) On raisonne par l'absurde. 
Nous avons vu dans la preuve du théorème \ref{theo-stabilite-BM-pseudo}-(\ref{theo-rang-pseudovariete-b}) (p.~\pageref{bullet-non-bijective}) que l'application $(\_)^{\bullet}$ n'est pas bijective si et seulement si, $2(\pp)<m$. Dans ce cas, si $m\geq\epsilon i$, on a
$$\preskip1exQ=i{-}(m{-}(\pp))
<i-(d_\Mg{-}1){m\over 2}\leq 
\bigpar{(2/\epsilon){+}1{-}d_{\Mg}}/{m\over 2}\leq 0\,,$$
puisque $d_{\Mg}\geq 2$, ce qui est contraire à la condition $Q\geq 0$.\enddemo

\subsection{Stabilité de $\set\Bettibm^{i}(\Delta_{?m-a}(\Mg^{m}){/}\Sm)/_{m}$, cas général}\label{stabilite-Betti-Delta-cas-general}Dans 
cette dernière section nous esquissons la démarche à suivre pour démontrer la généralisation suivante de la proposition
\ref{prop-stabilite-betti-dim3-BFm}.

\begin{theo}Soit $\Mg$ une pseudovariété connexe orientable, et soient $a,i\in\NN$. 
\begin{enumerate}
\item La famille
$\bigset
\Bettibm^{i}(\Delta_{m-a}(\Mg^{m}){/}\Sm)
/{}_{m}$
est constante pour 
$$\pmathalign{
m&\geq&  2i&+&2a\,,&\text{\quad si $\dMg=2$,}\\
m&\geq&  i&+&2a\,,\hfill &\text{\quad si $\dMg\geq3$.}}
$$
\item La famille
$\bigset
\Bettibm^{i}(\Delta_{\leq m-a}(\Mg^{m}){/}\Sm)
/{}_{m}$
est constante pour 
$$\pmathalign{
m&\geq&  2i&+&4a\,,&\text{\quad si $\dMg=2$,}\\
m&\geq&  i&+&\dMg a\,,\hfill &\text{\quad si $\dMg\geq3$.}}
$$
\end{enumerate}
\end{theo}
\def\Demonstration{Indications}
\demo Des modifications simples dans la preuve de la proposition \ref{prop-stabilite-betti-dim3-BFm} conduisent à l'énoncé suivant:

\noindent{\slshape Proposition.
Soit $\Mg$ une pseudovariété connexe orientée. Pour tous $t,i\in\NN $ fixés, le morphisme 
$$\varPhi_t(p_{m}^{*}):
\varPhi_t(\Hbm^{i}(\Fm(\Mg)))\to
\varPhi_t(\Hbm^{i}(\Fmm(\Mg)))$$
 est bijectif pour $m\geq \epsilon i+t$, où $\epsilon=2$ si $d_{\Mg}=2$, et $\epsilon=1$ si $d_{\Mg}\geq 3$.
}

\smallskip
\noindent{\slshape$\bullet$ Les cas $\bigset
\Bettibm^{i}(\Delta_{m-a}(\Mg^{m}){/}\Sm)
/{}_{m}$. }Par la proposition \ref{coro-car-l-m}, on a 
$$\Delta_{m-a}(\Mg^{m})/\Sm
\cong\coprod\nolimits_{\lambda\in\Y_{m-a}(m)}\Fg_{m-a}(\Mg)/\S_{\lambda}\,.$$
où $\S_{\lambda}\cont\1_{a}\times \S_{m-2a}$, de sorte 
que la stabilité recherchée est majorée par celle des familles
$\set\Betti^{i}(\varPhi_a\Hbm(\Fg_{m-a}(\Mg))/$, donc par
$\set m-a\geq \epsilon i+a/ $ d'après la {\slshape Proposition}.

\smallskip
\noindent{\slshape$\bullet$ Les cas $\bigset
\Bettibm^{i}(\Delta_{\leq m-a}(\Mg^{m}){/}\Sm)
/{}_{m}$. }
Par induction sur $a\in\NN$. Le cas $a=0$ est le lemme 4.2 de Church (\cite{chu}, p.~492). Dans le cas général on considère le complexe exact de $\Sm$-modules de cohomologie de Borel-Moore:
$$
\mathalign{
\to
\Hbm^{i-1}(\Delta_{\leq m-a})
\to&
\Hbm^{i-1}(\Delta_{m-a})
\to\hfill\\
&\hskip0.5cm{}\to\Hbm^{i-\dMg}(\Delta_{\leq m-(a+1)})\to\hfill\\
&\hskip2cm\to
\Hbm^{i}(\Delta_{\leq m-a})
\to
\Hbm^{i}(\Delta_{m-a})
\to}
$$
et le théorème résulte du lemme des cinq.
\enddemo

\section{Suites spectrales de Leray}
\glossarytitle{Suites spectrales de Leray}
Cette\label{Leray} section est consacrée à l'étude de  la dégénérescence des suites spectrales de Leray associées aux fibrations $\pi_a=\GaDelta_{?\ell}\Xg^{m}\to\Fa (\Xg)$ (en particulier $\pi_a=\Fba (\Xg)\to\Fa (\Xg)$),
 lorsque $\Xg$ est $i$-acyclique et localement connexe.

\def\Remarque{\miseengarde Avertissement}
\begin{rema*}{\miseengarde Le corps  $k$ est de caractéristique quelconque et l'espace $i$-acyclique $\Xg$ n'est pas supposé de type fini.}
\end{rema*}

\setbox111=\hbox{\ \vrule\ $\xymatrix@C=3mm@R=4mm{\Yg\ar[rd]_{\pi_{\Yg}}\ar[rr]|{\ f\ }&&\Zg\ar[dl]^{\pi_{\Zg}}\\
&\Bg}$}
\subsection{Cohomologie à support $\pi$-propre}
\subsubsectionline{Catégorie d'espaces au-dessus de $\Bg$.}\hmodeHabillage{\hdecalage{-0.68\hsize}{\copy111}}{2}{-10pt}Soit\label{c-mous} $\Bg$ un espace topologique.  On rappelle qu'on appelle \expression{espace de base $\Bg$, ou au-dessus de $\Bg$,} la donnée d'une application continue $\pi_{\Yg}:\Yg\to\Bg$, et \expression{morphisme d'espaces au-dessus de $\Bg$\/} de $\pi_{\Yg}:\Yg\to\Bg$ vers $\pi_{\Zg}:\Zg\to\Bg$ la donnée une application continue $f:\Yg\to\Zg$ rendant commutatif le diagramme ci-contre.
\endHabillage

Pour $\pi_{\Yg}:\Yg\to\Bg$ donné, la \expression{cohomologie à support $\pi_{\Yg}$-propre}, notée $\Hpi(\pi_\Yg)$ où $\Hpi(\Yg;\Bg)$ où même $\Hpi(\Yg)$, suivant les contextes\glossary{$\Hpi(\pi_\Yg)$, $\Hpi(\Yg;\Bg)$, $\Hpi(\Yg)$:cohomologie à support $\pi$-propre pour la fibration  $\pi_a:\Yg\to\Bg$}, comme la cohomologie du complexe 
$$
(\C_{\cv}(\Yg;\Bg),d_{\cv}):=\IR\Gamma(\Bg,\IR\pi_{\Yg!}\,\fs k_{\Yg})\,.$$
Les morphismes naturels de complexes
$$
\IR\Gammac(\Bg,\IR\pi_{\Yg!}\,\fs k_{\Yg})\to
\IR\Gamma(\Bg,\IR\pi_{\Yg!}\,\fs k_{\Yg})\to
\IR\Gamma(\Bg,\IR\pi_{\Yg*}\,\fs k_{\Yg})$$
induisent des morphismes naturels en cohomologie:
$$\Hc(\Yg)\to\Hpi(\Yg;\Bg)\to\Hr(\Yg)\,.\eqno(\ddagger\ddagger)$$

\subsubsectionnumber On rappelle (\cf\ref{resolution-c-molle}) que dans la mesure où la résolution du faisceau constant par les cochaînes d'Alexander-Spanier  $\fs k_{\Yg}\to(\AS^{\bullet}(\Xg;k),d_{\bullet})$ est une résolution $\varPhi$-molle pour toute famille paracompactifiante $\varPhi$ de $\Bg$, on dispose du candidat canonique des \expression{cochaînes à support $\pi$-propre}:
$$(\C_{\cv}(\Yg;\Bg),d_{\cv}):=\Gamma(\Bg,\pi_{\Yg!}(\AS^{\bullet}(\Xg;k),d_{\bullet}))\,,$$
puisque $\pi_{\Yg!}$ conserve la propriété d'être $c$-mou (\footnote{\label{ref-c-mou-c-mou}\cf \cite{KS} proposition~2.5.7-(ii), p.~105.}) et qu'un $c$-mou sur la pseudovariété dénombrable à l'infini $\Bg$ est $\Gamma(\Bg,\_)$-acyclique (\footnote{\label{ref-c-mou-c*-injectif}\cf \cite{KS}~proposition~2.5.10, p.~106.}).
Des raisons qui expliquent aussi que $\Hc(\Yg)$ est calculée par le complexe des \expression{cochaînes à support compact}
$$(\C_{\rmc}(\Yg;\Bg),d_{\rmc}):=\Gammac(\Bg,\pi_{\Yg!}(\AS^{\bullet}(\Xg;k),d_{\bullet}))\,.$$
Dans la suite on notera 
$$\halfdisplayskips\Z_{\cv}(\Yg;\Bg):=\ker(d_{\cv})\text{\quad et\quad}\B_{\cv}(\Yg;\Bg):=\im(d_{\cv}).$$
et de manière analogue pour le complexe $(\C_{\rmc}(\Yg;\Bg),d_{\rmc})$.

\subsubsectionnumber Lorsque $\Bg=\set\pt/$, on a $\Hpi(\_;\Bg)=\Hc(\_)$. La cohomologie à support $\pi$-propre étend la cohomologie à support compact et tout comme elle, si $j:U\hook\Yg$ est une inclusion ouverte, complémentaire de l'inclusion fermée $i:F\hook\Yg$, on munit $U$ et $F$ des structures d'espaces basés sur $\Bg$ en composant les inclusions avec $\pi_{\Yg}:\Yg\to\Bg$, et l'on dispose alors des morphismes
\def\varitemizeseps{\itemsep1pt\parskip0pt\topsep3pt}\begin{itemize}
\item de prolongement par zéro $j_{!}:\Hpi(U;\Bg)\to\Hpi(\Yg;\Bg)$
\item de restriction à un fermé $i^{*}:\Hpi(\Yg;\Bg)\to\Hpi(F;\Bg)$
\end{itemize}
et d'une suite exacte longue de cohomologie à support $\pi$-propre
$$
\mathrigid1mu
\cdots\to\Hpi^{i-1}(F;\Bg)\to\Hpi^{i}(U;\Bg)\to\Hpi^{i}(\Yg;\Bg)\to\Hpi^{i}(F;\Bg)\to\Hpi^{i+1}(U;\Bg)\to\cdots$$
qui est aussi un complexe de $\Hr(\Bg)$-modules.

\smallskip
En effet, une décomposition en parties respectivement ouverte et fermée $\Yg=U\sqcup F$ donne lieu à la suite exacte courte 
$0\to j_{!}\,\fs k_{U}\to
\fs k_{\Yg}\to
i_{!}\,\fs k_{F}\to0$ de $\Mod(\fsk_{\Yg})$,
et donc au triangle exact de $D^{+}_{k}(\Bg)$\glossary{$D^{+}_{k}(\Yg)$:catégorie dérivée des complexes de faisceaux de $k$-espaces vectoriels sur $\Yg$
bornés inférieurement}:
$$\IR\pi_{\Yg!}\,j_{!}\,\fs k_{U}\to
\IR\pi_{\Yg!}\,\fs k_{\Yg}\to
\IR\pi_{\Yg!}\,i_{!}\,\fs k_{F}\to\,,$$
où nous avons 
$\IR\pi_{\Yg!}\,j_{!}=\IR\pi_{U!}$ et $\IR\pi_{\Yg!}\,i_{!}=\IR\pi_{F!}$, puisque $j_!$ et $i_{!}$ sont exacts et transforment $c$-mous en $c$-mous ($^{\ref{ref-c-mou-c-mou}}$).

Enfin, dans le modèle des cochaînes d'Alexander-Spanier, les morphisme naturels 
$\pi_{\Yg!}j_{!}\,\AS^{i}_U\to \pi_{\Yg!}\,\AS^{i}_\Yg$ et 
$\pi_{\Yg!}\,\AS^{i}_\Yg\to \pi_{\Yg!}\,\AS^{i}_F$ correspondent clairement aux opérations de prolongement par zéro et de restriction habituels.

\subsubsectionnumber Si $f:\Zg\to\Yg$ est une application continue et \emph{propre} entre deux espaces au-dessus de $\Bg$, le morphisme d'adjonction $\fs k_{\Yg}\to (\IR f_{*}\circ f^{-1})\,\fs k_{\Yg}$ donne lieu au morphisme de complexes de faisceaux
$$\IR\pi_{\Yg!}\,\fs k_{\Yg}\to (\IR\pi_{\Yg!}\circ\IR f_{*})\, \fs k_{\Zg}
= (\IR\pi_{\Yg!}\circ\IR f_{!})\, \fs k_{\Zg}=
\IR\pi_{\Zg!}\,\fs k_{\Zg}\,,
$$
d'où le morphisme \expression{image-inverse} pour la cohomologie à support $\pi$-propre $f^{*}:\Hpi(\Yg;\Bg)\to\Hpi(\Zg;\Bg)$. La naturalité des morphismes $(\ddagger\ddagger)$ de \ref{c-mous}, donne alors lieu au diagramme  commutatif de morphismes image-inverse
$$\xymatrix@R=5mm{
\Hc(\Yg)\ar[d]^{f^{*}}\ar[r]&\Hpi(\Yg;\Bg)\ar[d]^{f^{*}}\ar[r]&\Hr(\Yg)\ar[d]^{f^{*}}\\
\Hc(\Zg)\ar[r]&\Hpi(\Zg;\Bg)\ar[r]&\Hr(\Zg)\\
}$$

 L'analogue pour la cohomologie $\Hpi(\_;\Bg)$ de la propriété caractéristique des espaces $i$-acycliques, le théorème \ref{caracterisation}-(\ref{caracterisation-clef}), est également vérifiée.\killline

\begin{prop}\label{caracterisation-basee}Soit $\Xg$ un espace $i$-acyclique. \'Etant donnés des espaces basés $\pi_\Zg:\Zg\to\Bg$ et $\pi_\Yg:\Yg\to\Bg$, soit $\pi_{\Xg\times\Yg}:\Xg\times\Yg\to\Bg$, $(x,y)\mapsto\pi_{\Yg}(y)$ et notons $p_2:\Xg\times\Yg\to\Yg$, $(x,m)\mapsto m$. Soit $f:\Zg\to\Xg\times\Yg$ un morphisme d'espaces au-dessus de $\Bg$, notons $f_2:=p_2\circ f$. Soit $j:V\hook\Yg$ un plongement ouvert, notons $f'_{2}:f^{-1}(V)\to V$ la restriction de $f_{2}$. On a le diagramme commutatif d'espaces au-dessus de $\Bg$ suivants,
$$\xymatrix@R=6mm{
f_{2}^{-1}(V)\ar[rd]|{f'_2}\ar@{^(->}[r]&\Zg\ar[rd]|{f_2}\ar[r]|(0.4){\,f\,}&\Xg\times\Yg\ar[d]^{p_2}\\
&V\ar@{^{(}->}[r]|(0.45){\,j\,}&\Yg\mrlap{\,.}
}$$
Alors, si $f$ et $f'_2$ sont propres et si  $j_{!}:\Hpi(V;\Bg)\to\Hpi(\Yg;\Bg)$ est surjectif, on a
$$\preskip-1ex\big(f^{*}:\Hpi(\Xg\times\Yg;\Bg)\to\Hpi(\Zg;\Bg)\big)=0\,.
$$
\end{prop}
\demo
C'est presque littéralement la même preuve que \ref{caracterisation}-(\ref{caracterisation-clef}).

Soient $p_1:\Xg\times\Yg\to\Xg$ la projection canonique et $f_{1}:=p_{1}\circ f$. 
Par Künneth, on a $\Hpi(\Xg\times\Yg;\Bg)=\Hc(\Xg)\otimes\Hpi(\Yg;\Bg)$ et il suffit de montrer que pour tous cocycles $\omega\in\Z_{\rmc}(\Xg)$ et $\varpi\in\Z_{\cv}(\Yg;\Bg)$, le cocycle
$$f^{*}(\omega\otimes\varpi)=f^{*}(p_1^{*}\omega\cup p_2^{*}\varphi)=
f_{1}^{*}\omega\cup f_{2}^{*}\varpi\in\Z_{\cv}(\Zg;\Bg)\eqno(\ast)$$
est la différentielle d'une cochaîne à support $\pi$-propre de $\Zg$.
Or, il existe par hypothèse $\tau\in\Z_{\cv}(V;\Bg)$ qui représente la classe de $\varpi$ dans $\Hpi(\Yg;\Bg)$. On peut donc remplacer dans $(\ast)$
$f_{2}^{*}\varpi$ par $f_{2}^{*}\tau$ qui est à support $\pi$-propre dans $f_{2}^{-1}(V)$ donc dans $\Zg$ par le prolongement par zéro $j_{!}$.
D'autre part, comme $\Xg$ est $i$-acyclique, on a $f_{1}^{*}\omega=f_{1}^{*}d\alpha$ pour une certaine cochaîne $\alpha$ de $\Xg$ (à support non nécessairement compact), et alors
$$f_{1}^{*}\omega\cup f_{2}^{*}\tau=d(f_1^{*}\alpha\cup f_{2}^{*}\tau)\,,$$ 
où $f_1^{*}\alpha\cup f_{2}^{*}\tau$ est à support $\pi$-propre de $\Zg$ puisqu'il en est ainsi de $f_{2}^{*}\tau$. On a donc bien $f^{*}(\omega\otimes\varpi)=0$ dans $\Hpi(\Zg;\Bg)$.
\enddemo

\subsection{Localisation du théorème de scindage}\label{localisation}
\subsubsectionline{Naturalité de $\Hpi(\_;\Bg)$ relative aux ouverts de $\Bg$.}\label{naturalite-Hcv}
Soit $\Bg'$ un ouvert de $\Bg$. \'Etant donné $\pi_\Yg:\Yg\to\Bg$, on note $\Yg':=\pi_\Yg^{-1}(\Bg')$ et $\pi_{\Yg'}:\Yg'\to\Bg'$ la restriction de $\pi_\Yg$. Nous avons ainsi un diagramme cartésien de plongements ouverts:
$$\preskip-1.ex
\xymatrix@R=4mm@C=1cm{
\Yg'\ar[d]_{\pi_{\Yg'}}\xylbl[rd]{$\Box$}\arinto[r]&\Yg\ar[d]^{\pi_{\Yg}}\\
\Bg'\arinto[r]&\Bg
}$$
et un morphisme naturel de restriction de complexes
$$\IR\Gamma(\Bg,\IR\pi_{\Yg!}\,\fs k_{\Yg})\to\IR\Gamma(\Bg',\IR\pi_{\Yg!}\,\fs k_{\Yg})=\IR\Gamma(\Bg',\IR\pi_{\Yg'!}\,\fs k_{\Yg'})
\,.$$
donnant lieu à un morphisme naturel de restriction:
$$\Hpi(\Yg;\Bg)\to\Hpi(\Yg';\Bg')\,.$$
Nous aurons besoin du résultat classique suivant (\footnote{Conséquence immédiate de la proposition 2.5.2 p.103, \cite{KS}.}).\killline

\begin{prop}Soient\label{germe-compact} $\Yg$ et $\Bg$ des espaces localement compacts et 
soit $\pi:\Yg\to\Bg$ une application continue. Pour tout $b\in\Bg$, le morphisme canonique
$$\limind_{\Bg'\ni b}\Hpi(\pi^{-1}(\Bg');\Bg')\to\Hc(\pi^{-1}(b))$$
est un isomorphisme.
\end{prop}

\subsubsectionline{Notations.}\label{notas-localisation}Pour la suite de cette section nous fixons un espace  $\Xg$ et un entier  $0\leq a\in\NN$. On note $\Fa:=\Fa (\Xg)$, si $a>0$, et $\Fg_0:=\set\pt/$. 

Pour tous $\aa\leq\ell\leq m\in\NN$, nous avons introduit (\ref{notas}-(\ref{nota-FG}))
les espaces 
$$\GaDelta_{?\ell}\Xg^{m}:=\Delta_{?\ell}\Xg^{m}\cap (\Xg^{m-a}\times\Fa)\,,$$
que nous notons aussi\glossary{{$\GaDelta_{?\ell,m}$}:notation abrégée de $\GaDelta_{?\ell}\Xg^{m}:=\Delta_{?\ell}\Xg^{m}\cap (\Xg^{m-a}\times\Fa)$}
$\GaDelta_{?\ell,m}:=\GaDelta_{?\ell}\Xg^{m}$
lorsque l'on aura besoin d'une écriture plus compacte. Nous y avons aussi introduit 
la projection  $\pi_a:\Xg^{m}\to\Xg^{a}$ sur les dernières $a$ coordonnées, et les espaces de base $\Fa$
$$
\pia :\GaDelta_{?\ell,m}\to\Fa
\text{\quad et\quad}
\pia :\Fm \to\Fa
\,.$$
Maintenant, pour tout ouvert $\Ug\dans\Fa$, nous notons
$$
\Ug\GaDelta_{?\ell,m}:= \pia ^{-1}\Ug\cap\GaDelta_{?\ell,m}
\text{\quad et\quad}
\Ug\Fm := \pia ^{-1}\Ug\cap\Fm\,, 
$$
et considérons les produits fibrés:
$$\hskip-2em
\xymatrix@C=6mm@R=6mm{
\Ug\GaDelta_{?\ell,m}\ar@{^{(}->}[r]\ar[d]_{\pia}\xylbl[rd]{$\Box$}
&\GaDelta_{?\ell,m}\ar[d]^{\pia}\\
\Ug\ar@{^{(}->}[r]&\Fa\\
}
\quad
\xymatrix@C=8mm@R=7.5mm{
\Ug\Fm \ar@{^{(}->}[r]\ar[d]_{\pia}\xylbl[rd]{$\Box$}
&\Fm \ar[d]^{\pia}\\
\Ug\ar@{^{(}->}[r]&\Fa\\
}
\quad
\xymatrix@C=6mm@R=7.5mm{
\Xg\times\Ug\Fm \ar@{^{(}->}[r]\ar[d]_{\pia}\xylbl[rd]{$\Box$}
&\Xg\times\Fm \ar[d]^{\pia\mrlap{\mkern55mu\scriptstyle\pia}}\mrlap{\ni\smash{\vtop{\setbox11=\hbox{$(x,\cl y)$}\copy11
\hbox to\wd11{\hss$\vfld{}{}{7.68mm}$\hss}}}}\\
\Ug\ar@{^{(}->}[r]&\Fa\mrlap{\mkern10mu \ni\mkern10mu \pia(\cl y)}\\
}\qquad
$$
Les décompositions ouvertes-fermées 
$$\mathalign{
\hfill
\Ug\GaDelta_{\leq\ell,m}&=&\Ug\GaDelta_{\ell,m}&\sqcup&\Ug\GaDelta_{\leq\ellmo,m}\hfill\\\noalign{\kern4pt}
\Xg\times\Ug\Fm &=&\Ug\Fg_{1+m}&\sqcup& \Delta_{m}(\Xg\times\Ug\Fm )}
$$ sont alors des décompositions d'espaces basés sur $\Ug$. 

\smallskip
On peut maintenant énoncer l'analogue du théorème \ref{theo-scindage} dans le contexte d'espaces basés sur $\Ug$.
Il est important de remarquer que si dans \ref{theo-scindage} la $i$-acyclicité était une condition suffisante, maintenant elle apparaît comme une condition nécessaire et suffisante.

\begin{theo}[de scindage local]\label{theo-scindage-basee}Soient $1\leq a\leq\ell\leq m\in\NN$. On fixe un ouvert $\Ug\dans\Fa$.
Les espaces dans cet énoncé sont des sous-espaces de $\pi_a:\Xg^{m-a}\times\Ug\to\Ug$. On note $\Hpi(\_)$ la cohomologie à support $\pia$-propre.
Alors, $\Xg$ est $i$-acyclique, si et seulement si, il vérifie les assertions suivantes.
\mynobreak\begin{enumerate}\nobreak
\item\leavevmode\label{theo-scindage-basee-a}Pour $m\geq a\geq1$ et un ouvert $\Ug\dans\Fa$, le morphisme de restriction
$$\halfdisplayskips\Hpi(\Xg\times\Ug\Fm )\to
\Hpi(\Delta_{\leq m}(\Xg\times\Ug\Fm ))$$
est nul et la suite
$$\halfdisplayskips0\to
\Hpi(\Ug\Fm )[-1]^{m}\to
\Hpi(\Ug\Fmm )\to
\Hpi(\Xg\times\Ug\Fm )\to0\,,
$$
extraite de la suite longue de cohomologie à support $\pi$-propre, est exacte.
\item\leavevmode\label{theo-scindage-basee-b}Pour $m\geq\ell\geq a\geq1$  et un ouvert $\Ug\dans\Fa$, le morphisme de restriction
$$\Hpi(\Ug\GaDelta_{\leq\ell,m})\to\Hpi(\Ug\GaDelta_{\leq\ellmo,m})
\postskip-2pt$$
est nul et la suite 
$$\halfdisplayskips0\to
\Hpi(\Ug\GaDelta_{\leq \ellmo,m})[-1]\to
\Hpi(\Ug\GaDelta_{\ell,m})\to
\Hpi(\Ug\GaDelta_{\leq \ell,m})\to
0\,,
$$
extraite de la suite longue de cohomologie à support  $\pi$-propre, est exacte.
\end{enumerate}
\end{theo}
\demo On suppose que $\Xg$ est $i$-acyclique. La preuve de 
(\ref{theo-scindage-basee-a}) et (\ref{theo-scindage-basee-b}) est alors la même que dans \ref{theo-scindage} modulo la version de la propriété fondamentale des espaces $i$-acycliques\ pour les espaces basés sur $\Ug$ de la proposition \ref{caracterisation-basee}. 

Pour l'assertion  (\ref{theo-scindage-basee-a}), on est conduit au diagramme
$$\xymatrix@R=5mm{
\Delta_{\leq m}(\Xg\times\Ug\Fm )\ar@{^(->}[r]\ar[rd]|(0.55){\ f_2\ }&
\Xg\times \Ug\Fm \ar[d]^{p_2}\\
&\Ug\Fm 
}$$
où l'application $f_2$ est propre car revêtement trivial d'ordre $m$.

Pour (\ref{theo-scindage-basee-b}), on raisonne par induction sur $\ell$. Sa valeur la plus petite est $\ell:=a$, auquel cas $\GaDelta_{\leq\ellmo,m}=\emptyset$ et l'assertion est claire quel que soit $m$. Ensuite, en supposant l'assertion établie pour $\ellmo\geq a$, 
on est conduit à considérer le diagramme commutatif suivant  où $m\geq \ell$
$$\xymatrix@R=5mm{
f_{2}^{-1}(\Ug\GaDelta_{\ellmo,\mmo})\ar@{^(->}[r]\ar[rd]|{\ f'_2\ }&\Ug\GaDelta_{\leq\ellmo,m}\ar@{^(->}[r]\ar[rd]|{\ f_2\ }&
\Xg\times \Ug\GaDelta_{\leq\ellmo,\mmo}\ar[d]^{p_2}\\
&\Ug\GaDelta_{\ellmo,\mmo}\ar@{^(->}[r]&\Ug\GaDelta_{\leq\ellmo,\mmo}
}$$
où $f'_2$ est propre car restriction (par l'image) de l'application propre de même nom $f'_2:f_2^{-1}(\Delta_{\ellmo}^{\mmo})\to\Delta_{\ellmo}^{\mmo}$ de la preuve de \ref{theo-scindage}-(\ref{theo-scindage-b}).

\smallskip
Réciproquement, si $m=\ell=a=1$ et si $\Ug=\Xg$, le morphisme 
$\Hpi(\Xg\times\Xg)\to\Hr(\Delta_1(\Xg\times\Xg))$ est nul. Or, ce morphisme s'identifie au cup-produit
$$\cup:\Hc(\Xg)\otimes\Hr(\Xg)\to\Hr(\Xg)\postskip2pt$$
et $\Xg$ est bien $i$-acyclique.
\enddemo

\subsubsectionnumber
Le corollaire suivant est un résultat technique qui sera très utile pour montrer que les faisceaux des cohomologie à support $\pi$-propre de la fibration $\pia :\Fm \to\Fa$ sont constants sur les composantes connexes de $\Fa$. Ceci participera ensuite de manière décisive dans l'étude des suites spectrales de Leray correspondantes. (Voir la remarque \ref{trivialite-cohomologique}.)

\begin{prop}\label{coro-critique}Soit $\Xg$ un espace $i$-acyclique. Soient $a\leq m\in\NN$ et $\pi_{m,a}:\Fm (\Xg)\to\Fa(\Xg)$, $(x_1,\ldots,x_{m})\mapsto(x_{m-a+1},\ldots,x_{m})$.
Pour\label{coro-critique-a} tout $\cl x\in\Fa$, le morphisme de restriction
$$\halfdisplayskips
\rho_{m,a}:\Hpi(\Fm (\Xg))\to\Hc(\pi_{m,a} ^{-1}\cl x)\postskip0pt$$
est surjectif.
\end{prop}
\demo\halfdisplayskips
 On raisonne par induction sur $m\geq a$. Lorsque $m=a$, on a $\pi_a=\id_{\Fa}$ et donc $\Hpi(\Fm )=\Hr(\Fa)$, et cette cohomologie se surjecte clairement sur $\Hc(\cl x)=\Hr^{0}(\cl x)=k$.

Dans le cas général, on considère décomposition en parties respectivement ouverte et fermée 
$$\Xg\times\Fm=\Fmm\ \sqcup\ \Delta_{m}(\Xg\times\Fm)$$
et la décomposition qu'elle induit sur le fermé $\Xg\times\pi_{m,a}^{-1}\cl x$, à savoir
$$\Xg\times\pi_{m,a}^{-1}\cl x=\pi_{\mm,a}^{-1}\cl x\ \sqcup\ \Delta_{m}(\Xg\times\pi_{m,a}^{-1}\cl x)\,.$$
On considère alors le morphisme de restrictions des suites de cohomologie à supports $\pi$-propres
$$\halfdisplayskips
\cxymatrix{@C=4mm@R6mm}{
0\ar[r]&
\Hpi(\Fm )[-1]^{m}\ar[r]\aronto[d]_{(\rho_{m,a})^{m}}&
\Hpi(\Fmm )\ar[r]\ar[d]_{\rho_{m+1,a}}&
\Hpi(\Xg{\times}\Fm )\ar[r]\aronto[d]_{\xi}&
0\\
\ar[r]&
\Hc(\pi_{m,a}^{-1}\cl x)[-1]^{m}\ar[r]&
\Hc(\pi_{m+1,a}^{-1}\cl x)\ar[r]&
\Hc(\Xg{\times}\pi_{m,a}^{-1}\cl x)\ar[r]&
\\
}\eqno(\ast)$$
où la première ligne est exacte d'après \ref{theo-scindage-basee}-(\ref{theo-scindage-basee-a}) (la seconde aussi, mais on n'en aura pas besoin). 
La colonne de gauche est surjective puisque $\rho_{m,a}$ l'est par hypothèse de récurrence, et la colonne de droite l'est puisque le morphisme de restriction $\xi:\Hpi(\Xg\times\Fm \to\Hc(\Xg\times \pi_{m,a}^{-1}\cl x)$ s'identifie trivialement à $\id\otimes\rho_{m,a}:\Hc(\Xg)\otimes\Hpi(\Fm )\to
\Hc(\Xg)\otimes\Hc(\pi_{m,a}^{-1}\cl x)$. Une chasse au diagramme élémentaire montre alors la surjectivité de $\rho_{\mm,a}$.
\enddemo

\subsubsection{Polynôme de Poincaré de la cohomologie à support $\pi$-propre}Pour tout ouvert $\Ug\dans\Fa$ et tout
 $\pi:\Zg\to\Ug$ tel que $\dim_{k}\Hpi(\Zg;\Ug)<\infty$, on note, de manière analogue à \ref{Poincare-r-c},\glossary{$\Pcv(\Zg;\Ug)(T)$:polynôme de Poincaré pour la cohomologie à support $\pi$-propre d'un espace $\Zg$ basé sur $\Ug\dans\Fa (\Xg)$}
$$\Pcv(\Zg;\Ug)(T):=\sumnl_{i\in\ZZ}\dim_{k}(\Hpi^i(\Zg;\Ug))\,T^{i}\,.$$

Le corollaire suivant de \ref{theo-scindage-basee}-(\ref{theo-scindage-basee-a}) est un ingrédient essentiel dans la preuve de dégénérescence de la suite spectrale de Leray du théorème \ref{degen}. (Comparer à \ref{PcUFm}.)

\begin{coro}\label{PcvUFm}Soit $\Xg$ un espace $i$-acyclique tel que $\dim_k\Hc(\Xg)<\infty$ et soit $\Ug\dans\Fa(\Xg)$ un ouvert tel que $\dim_k\Hr(\Ug)<\infty$. Alors, pour tout $b>0$, 
$$
\Pcv(\Ug\Fba (\Xg);\Ug)=\Pr(\Ug)\cdot\Pc(\Fg_{b}(\Xg\mmoins a))\,.
$$
\end{coro}
\demo Par \ref{theo-scindage-basee}-(\ref{theo-scindage-basee-a}), on a l'égalité
$$\Pcv(\Ug\Fm ;\Ug)=\Pcv(\Ug\Fg_{\mmo};\Ug)\cdot(\Pc(\Xg)+(\mmo)T)$$
que l'on va itérer jusqu'à ce que $\mmo=a$, auquel cas $\Pcv(\Ug\Fa;\Ug)=\Pr(\Ug)$ puisque $\Ug\Fg_a=\Ug$ et que $\pi_a=\id_{\Ug}$. Le produit résiduel 
$$(\Pc(\Xg)+aT)\cdots(\Pc(\Xg)+(\mmo)T)$$ étant alors clairement égal à $\Pc(\Fg_{m-a}(\Xg\mmoins\aa))$.
\enddemo

\subsection{Constance des faisceaux de cohomologie à support $\pi$-propre}\label{constance-de-faisceau}
Pour $a\leq\ell\leq m$, reprenons les espaces $\Fm $ et $\GaDelta_{?\ell,m}$  de base $\Fa$ et les décompositions ouvertes-fermées de base $\Fa$
$$\mathalign{
\hfill\Xg\times\Fm &=&\Fg_{1+m}\sqcup \Delta_{m}(\Xg\times\Fm )\hfill&\mrlap{\hskip0.9cm(\diamond_1)}\\
\hfill\GaDelta_{\leq\ell,m}&=&\GaDelta_{\ell,m}\sqcup\GaDelta_{\leq\ellmo,m}\hfill&\mrlap{\hskip0.9cm(\diamond_2)}}
$$

\comment
Notons $\fs k_{\GaDelta_{?\ell,m}}$ le faisceau constant sur $\GaDelta_{?\ell,m}$. Soit $c:\GaDelta_{?\ell,m}\to\bullet$ l'application constante.
La théorie générale de faisceaux nous donne l'égalité
$$\IR(c_{?}\circ\pi_{a!})(\_)=
\IR\Gamma_{?}(\Fa;\IR \pi_{a!}(\_))\,,$$
où $?={}!,*$. Les groupes de cohomologie de ce complexe évalué en $\fs k_{\GaDelta_{?\ell,m}}$ sont respectivement $\Hc(\GaDelta_{?\ell,m})$ et $\Hpi(\GaDelta_{?\ell,m})$.
Avant d'aborder le problème de la convergence des suites spectrales de Leray, nous étudions dans cette section les faisceaux $\IR^{i}\pi_{a!}(\_)$ pour les espaces de configuration généralisés associés à une variété topologique $\Xg$ supposée $i$-acyclique.
\endcomment

\subsubsection{Les faisceaux $\cHpi^{i}(\_)$ sur $\Fa$}\label{cHcv}\smallskip\noindent{\sl Le cas $\GaDelta_{?\ell,m}$.}
Reprenons l'application $\pia :\GaDelta_{?\ell,m}\to\Fa$ et introduisons les faisceaux\glossary{${\cHpi^{i}(\Yg):=\IR^{i}\pi_{a!}(\fs k_{\Yg})}$:faisceau de cohomologie à support $\pi$-propre pour l'application $\pia :\Yg\to\Fa$}:
$$\displayboxit{\cHpi^{i}(\GaDelta_{?\ell,m}):=\IR^{i}\pi_{a!}(\fs k_{\GaDelta_{?\ell,m}})
\,,\quad\forall i\in\ZZ}$$
dont les fibres sont (\ref{germe-compact})
$$
\displayboxit{\cHpi^{i}(\GaDelta_{?\ell,m})_{\cl x}=\Hc^{i}(\pia ^{-1}(\cl x))\,,\quad\forall\cl x\in\Fa(\Xg)}
$$
D'autre part, le foncteur $\IR\pi_{a!}$ appliqué au
le triangle exact de $D^{+}_{k}(\GaDelta_{\leq\ell,m})$
$$\preskip0.2ex
j_{!}\fs k_{\GaDelta_{\ell,m}}\too^{j_!}\fs k_{\GaDelta_{\leq\ell,m}}\too^{\rho} i_{!}\fs k_{\GaDelta_{\leq\ellmo,m}}\to
$$
associé à la décomposition $(\diamond_2)$, où $j_!$ est le prolongement par zéro et $\rho$ la restriction, donne la suite exacte longue de faisceaux localement constants:
$$
\too\cHpi^{i}(\GaDelta_{\ell,m})\too^{j_{!i}}
\cHpi^{i}(\GaDelta_{\leq\ell,m})\too^{\rho_i}
\cHpi^{i}(\GaDelta_{\leq\ellmo,m})\too^{+1}\,.
\eqno(\ddagger\ddagger)$$

\smallskip\noindent{\sl Le cas $\Fm (\Xg)$.}  La même démarche sur la décomposition  $(\diamond_1)$, nous conduit à la suite longue de faisceaux localement constants sur $\Fa$:
$$
\too_{+1}
\cHpi^{i}(\Fg_{1+m})\too^{j_{!i}}
\cHpi^{i}(\Xg{\times}\Fm )\too^{\rho_i}
\cHpi^{i}(\Delta_{m}(\Xg{\times}\Fm) )\too_{+1}\,.
\eqno(\ddagger)$$

\begin{theo}[de scindage]\label{autre-scindage}Soit $\Xg$ un espace $i$-acyclique.
\mynobreak\begin{enumerate}\nobreak
\item\leavevmode\label{autre-scindage-a}Pour $a\leq m\in\NN$, les morphismes de faisceaux de la suite $(\ddagger)$
$$\cHpi^{i}(\Xg{\times}\Fm )\too^{\rho_i}
\cHpi^{i}(\Delta_{m}(\Xg{\times}\Fm) )\,,\quad\forall i\in\ZZ\,,$$
sont  nuls et l'on a une suite exacte courte de faisceaux
$$
0\to
\cHpi(\Delta_{m}\Xg{\times}\Fm )[-1]\to
\cHpi(\Fg_{1+m})\to
\cHpi(\Xg{\times}\Fm )
\to
0
\,.
$$

\item\leavevmode\label{autre-scindage-b}
Pour $0\leq a\leq\ellmo$ et $\ell\leq m\in\NN$, les morphismes de faisceaux dans $(\ddagger\ddagger)$
$$
\cHpi^{i}(\GaDelta_{\leq\ell,m})\too^{\rho_i}
\cHpi^{i}(\GaDelta_{\leq\ellmo,m})\,,\quad\forall i\in\ZZ\,,
$$
sont nuls et l'on a la suite exacte courte de faisceaux
$$
0\to
\cHpi(\GaDelta_{\leq\ellmo,m})[-1]\to
\cHpi(\GaDelta_{\ell,m})\to
\cHpi(\GaDelta_{\leq\ell,m})
\to
0
\,.
$$
\item Si\label{autre-scindage-c} $\Xg$ est, de plus, localement connexe, les faisceaux  dans {\rm(\ref{autre-scindage-a},\ref{autre-scindage-b})} sont  cons\-tants sur les composantes connexes de $\Fa(\Xg)$ \rm (\cf \ref{comm-independance-Hc-X-a} et \ref{fibres-pia}).
\end{enumerate}
\end{theo}
\demo\halfdisplayskips
 La preuve de {\rm(\ref{autre-scindage-a})} et {\rm(\ref{autre-scindage-b})} est essentiellement la même que celle de \ref{theo-scindage}. On indique brièvement les modifications à faire pour prouver (\ref{autre-scindage-b}). La preuve pour (\ref{autre-scindage-a}) suit exactement la même démarche.

\noindent {\sl Preuve de \rm (\ref{autre-scindage-b}).} On reprend les notations de \ref{notas-localisation}. Pour un ouvert $\Ug\dans\Fa$, on pose  $\Ug\GaDelta_{?\ell,m}:=\pia ^{-1}(\Ug)$\glossary{${\Ug\Yg}$:$\pia :\Yg\to\Fa$ et $\Ug\dans\Fa$, on pose $\Ug\Yg:=\pia ^{-1}\Ug$}. On a les suites exactes courtes (\ref{theo-scindage-basee}-(\ref{theo-scindage-basee-b})):
$$0\to
\Hpi(\Ug\GaDelta_{\leq \ellmo,m})[-1]\to
\Hpi(\Ug\GaDelta_{\ell,m})\to
\Hpi(\Ug\GaDelta_{\leq \ell,m})\to
0\,.\eqno(\ast)
$$
Cela étant, la naturalité de la cohomologie à supports propres relative à la base (\ref{naturalite-Hcv}) permet de dire que si
 $\pia :\Mg\to\Fa$ est un espace topologique basé sur $\Fa$, et que l'on note $\Ug\Mg:=\pia ^{-1}(\Ug)$, la correspondance
$$\Ug\fonct\smash{\Hpi(\Ug\Mg)}$$
définit le préfaisceau $\cHHcv(\Mg)$ de $k$-espaces vectoriels sur $\Fa$. Cette idée appliquée aux suites $(\ast)$
donne la suite exacte de préfaisceaux sur $\Fa$:
$$\let\cHpi\cHHcv
0\to
\cHpi(\GaDelta_{\leq\ellmo,m})[-1]\to
\cHpi(\GaDelta_{\ell,m})\to
\cHpi(\GaDelta_{\leq\ell,m})
\to
0
\,,
$$
et comme le faisceau engendré par $\cHHcv(\_)$ est précisément le faisceau $\cHpi(\_)$, on obtient l'exactitude de la suite courte des faisceaux sur $\Fa$:
$$
0\to
\cHpi(\GaDelta_{\leq\ellmo,m})[-1]\too
\cHpi(\GaDelta_{\ell,m})\too
\cHpi(\GaDelta_{\leq\ell,m})
\to
0
\,,
$$
et donc aussi l'annulation du morphisme $\rho_i$ dans la suite longue $(\ddagger\ddagger)$.

\medskip\noindent{\small\rightskip0em\leftskip3em
\sl Commentaire. On remarquera que grâce à l'exactitude des ces suites courtes, l'assertion: {\rm ``\,le faisceau  
\smashbot{$\cHpi(\GaDelta_{\leq\ell,m})$} est constant sur une composante connexe $\Cg$\glossary{$\Cg$:composante connexe de $\Fa$} de $\Fa$\,''} sera conséquence de ce que
le faisceaux $\cHpi(\GaDelta_{\leq\ellmo,m})$ et $\cHpi(\GaDelta_{\ell,m})$ la vérifient\comment
(\footnote{Le quotient de deux faisceaux constants et un faisceau constant.})\endcomment. Or, comme la plus petite valeur possible de $\ell$ vérifie $a=\ellmo$ et que $\cHpi(\GaDelta_{\leq a,m})=\cHpi(\GaDelta_{a,m})$, une preuve de l'assertion par récurrence découle de prouver seulement que les faisceaux $\cHpi(\GaDelta_{\ell,m})$ sont constants sur $\Cg$ et ce, pour tout $\ell\geq a$ (et pas seulement $\ell>a$). \par}

\medskip
\noindent {\sl Preuve de {\rm(\ref{autre-scindage-c})}}. 
Nous commençons par montrer que les faisceaux $\cHpi(\Fm )$ sont constants sur une composante connexe  $\Cg$ de $\Fa$. 

Nous procédons par induction sur $m\geq a$, le cas $m=a$ étant évident. 
Supposons maintenant que $\cHpi(\Fm )$ est constant sur $\Cg$. 
Comme $\Delta_{m}(\Xg\times\Fm )$ est isomorphe au produit $\iii[1,m]\times\Fm $, le faisceau $\cHpi(\Delta_{m}\Xg{\times}\Fm )$ est isomorphe à $\cHpi(\Fm )^{m}$ et il est donc constant sur  $\Cg$. D'autre part, on a  $\cHpi(\Xg{\times}\Fm )\simeq{ \Hc(\Xg)}\otimes_{k}\cHpi(\Fm )$ et $\cHpi(\Xg{\times}\Fm )$ est aussi constant sur $\Cg$. Par conséquent, dans la suite exacte courte de (\ref{autre-scindage-a})
$$
0\to
\cHpi(\Delta_{m}\Xg{\times}\Fm )[-1]\to
\cHpi(\Fg_{1+m})\to
\cHpi(\Xg{\times}\Fm )
\to
0
\,,\eqno(**)
$$
les faisceaux de droite et de gauche sont constants. Nous allons en déduire que celui du centre l'est également.

\smallskip\noindent$\bullet$ {\sl Critère de trivialité d'extensions de faisceaux constants.\label{critere} {\rm (\footnote{Voir aussi la remarque \ref{rema-Fg-orientable}.})}
Dans une suite exacte courte de faisceaux sur une composante connexe $\Cg$ de $\Fa$
$$\0\to\fs A\to\cL\to\fs B\to0$$
où $\fs A$ et $\fs B$ sont constants, le faisceau $\cL$ est constant, si et seulement si la suite est scindée, donc si et seulement si, l'application naturelle 
$$\Hom_{\fs k_{\Fa}}(\fs B,\cL)\to\Hom_{\fs k_{\Fa}}(\fs B,\fs B)$$ est surjective, et comme $\fs B$ est constant sur un espace localement connexe, ceci équivaut à la surjectivité de }
$$
\Gamma(\Cg,\cL)\onto\Gamma(\Cg,\fs B)=B_{z}\,,\quad\forall z\in\Cg\,.$$

Dans le cas de la suite de faisceaux $(**)$, si $\cl x\in\Fa$ et si $\Cg$ est la composante connexe de $\Fa$ contenant $\cl x$, on a
$$\mathalign{
\Gamma(\Cg,\cHpi(\Xg{\times}\Fm ))&=&\cHpi(\Xg{\times}\Fm )_{\cl x}=\Hc(\Xg)\otimes\cHpi(\Fm )_{\cl x}\\\noalign{\kern4pt}
&=&\Hc(\Xg)\otimes\Hc(\pia ^{-1}\cl x)\,,\hfill}$$
d'après \ref{germe-compact}.
D'autre part, le morphisme $\Gamma(\Cg,\_)\to(\_)_{\cl x}$ \og germe en $\cl x$ 
\fg, donne sur le \emph{\sl préfaisceau} $\cHHcv(\Xg{\times}\Fm )$
$$\cxymatrix{@R=5mm@C=1.2cm}{
\Gamma(\Cg,\cHHcv(\Xg{\times}\Fm ))\ar@<-2.ex>@{=}[d]\ar[r]&\cHHcv(\Xg{\times}\Fm )_{\cl x}\ar@<-2.5ex>@{=}[d]\mrlap{=\Gamma(\Cg,\cHpi(\Xg{\times}\Fm ))}\\
\Hc(\Xg)\otimes\Hpi(\Cg\Fm )\ar@{-->>}[r]^{\id\otimes\rho_{\Cg}}&\Hc(\Xg)\otimes\Hc(\pia ^{-1}\cl x)&\kern1.8cm
}\eqno\vtop{\kern-25pt\hbox{$(\dagger)$}}$$
où l'on voit apparaître la restriction
$\rho_{\Cg}:\Hpi(\Cg\Fm )\to\Hc(\pia ^{-1}\cl x)$
qui est surjective comme conséquence de la surjectivité de $\rho_{m,a}:\Hpi(\Fm )\to\Hc(\pia ^{-1}\cl x)$, d'après \ref{coro-critique}, et du fait que, 
dans la décomposition 
$$\Hpi(\Fm )=\bigoplusnl_{\Cg'\in\Pi_{0}\Fa}\Hpi(\Cg'\Fm )\,,$$ la restriction
$\Hpi(\Cg'\Fm )\to\Hc(\pia ^{-1}\cl x)$ est nulle si $\Cg'\not\ni\cl x$.
Le morphisme $(\dagger)$ est par conséquent surjectif. 
Comme d'autre part une section globale d'un préfaisceau détermine une section globale du faisceau associé, le diagramme naturel induit par les prolongements par zéro
$$\preskip4pt\postdisplaypenalty=10000\xymatrix@R=5mm@C1.5cm{
\Gamma(\Cg,\cHHcv(\Fg_{1+m}))\ar@<-2.ex>[d]\ar@{->>}[r]^(0.45){j_{!}}&\Gamma(\Cg;\cHHcv(\Xg{\times}\Fm ))\ar@<-2.5ex>@{->>}[d]^{\;(\dagger)}\\
\Gamma(\Cg,\cHpi(\Fmm ))\ar[r]^(0.45){j_{!}}&\Gamma(\Cg;\cHpi(\Xg{\times}\Fm ))
}$$
permet de conclure à la surjectivité du morphisme
$$\Gamma(\Cg,\cHpi(\Fg_{1+m}))\ontoo\Gamma(\Cg,\cHpi(\Xg{\times}\Fm ))\,,$$
et le critère de trivialité d'extensions de faisceaux constants s'applique, et le faisceau $\cHpi(\Fg_{1+m})$ est bien constant sur $\Cg$.

\medskip 
\noindent {\sl Preuve de la deuxième partie de \rm (\ref{autre-scindage-c}).} D'après le commentaire à la fin de la preuve de (\ref{autre-scindage-b}), il suffit de prouver que les faisceaux
$\cHpi(\GaDelta_{\ell,m})$
sont constants sur les composantes connexes de $\Fa$.

On a la décomposition en parties ouvertes (\cf\ref{connexes})
$$\halfdisplayskips\Delta_{\ell}^{m}=\coprod\nolimits_{\pgoth\in\Pgoth_{\ell}(m)}\Fg_{\pgoth}\,.$$
L'inclusion $\Fg_{\pgoth}\dans\GaDelta_{\ell,m}$ a lieu, si et seulement si, la partition $\pgoth$ de $\iii[1,m]$ décompose l'intervalle $\iii[m-a+1,m]$ en singletons. Appelons une 
telle partition \og $\pi_a$-adaptée\fg. Si $\pgoth$ n'est pas $\pi_a$-adaptée,
on a $\Fg_{\pgoth}\cap\GaDelta_{\ell,m}=\emptyset$. 

Si $\pgoth$ est $\pi_a$-adaptée, la projection $\pia :\Fg_{\pgoth}\to\Fa$ (sur les $a$ derniers termes) est trivialement équivalente à des projections $\pia :\Fg_\ell\to\Fa$ déjà traitées. Le faisceau $\cHpi(\Fg_{\pgoth})$ est donc constant sur les composantes connexes de $\Fa$.
On conclut ensuite grâce à la décomposition évidente
$$\halfdisplayskips\cHpi(\GaDelta_{\ell,m})=\bigoplusnl_{\pgoth}\cHpi(\Fg_{\pgoth})\,,$$
où $\pgoth$ décrit l'ensemble des partitions $\pgoth\in\Pgoth_{\ell}(m)$ qui sont $\pi_a$-adaptées.
\enddemo

\begin{rema}[et notation]L'assertion\label{ameliore-independance-Hc-X-a} \ref{autre-scindage}-(\ref{autre-scindage-c}) est la généralisation de la trivialité monodromique annoncée dans le commentaire \ref{comm-independance-Hc-X-a} par le fait qu'elle montre pour $\pia:\Fba(\Xg)\to\Fa(\Xg)$ l'une identification canonique des fibres $\cHpi(\Fba(\Xg))_{\cl x}=\Hc(\Fb(\Xg\mmoins \cl x))$ lorsque $\cl x$ parcourt une composante connexe $\Cg$ de $\Fa(\Xg)$. Il est intéressant d'observer que la connexité par arcs n'est pas concernée et que le phénomène dépasse de ce fait le cadre classique de la monodromie. L'assertion permet de comprendre le sens à donner à la notation $\Hc(\Fm(\Xg\mmoins a))$ où $a,m\in\NN$ (\cf\ref{encore-independance-Hc-X-a}).
\end{rema}

\subsection{Dégénérescence des suites spectrales de Leray}\label{degenerescence}

\subsubsectionline{Critère élémentaire de dégénérescence.}Nous donnons ici un critère de dégénérescence de suites spectrales par comparaison à une suite spectrale dégénérée. Le critère est très élémentaire, mais il sera utilisé à plusieurs reprises, ce qui justifie que nous l'énoncions séparément.\label{critere-degen}

Comme il est d'usage, l'expression \expression{\og la suite spectrale $(\IE_r,d_r)$ est dégénérée\fg} sera synonyme de \expressiong{$d_r=0$, pour tout $r\geq2$}. 

\begin{prop}\label{prop-critere-degen}Soit $\bigset\varphi_r:(\IE_{r},d_r)\to(\IE'_{r},d'_r)\mid{r\in\NN}/$ un morphisme de suites spectrales.
\begin{enumerate}\itemsep3pt
\item\leavevmode\label{prop-critere-degen-a}Si $(\IE_{r},d_r)$ est dégénérée et $\varphi_2$ est surjectif, alors $(\IE'_{r},d'_r)$ est dégénérée.
\item\leavevmode\label{prop-critere-degen-b}Si $(\IE'_{r},d'_r)$ est dégénérée et $\varphi_2$ est injective, alors $(\IE_{r},d_r)$ est dégénérée.
\end{enumerate}
Dans les deux cas on a $\varphi_r=\varphi_2$, pour tout $r\geq2$.
\end{prop}
\demo (\ref{prop-critere-degen-a}) On montre  par induction sur $\mathrigid2mu
r\geq2$, que $d'_r=0$ et que $\mathrigid2mu
\varphi_{r}=\varphi_2$. Lorsque $r:=2$, le complexe $(\IE'_2,d'_2)$ est quotient de $(\IE_2,0)$ par hypothèse, donc $d'_2=0$. Maintenant, si nous avons montré que $d'_r=0$, on aura $\varphi_{r+1}=\varphi_r$ et alors $\varphi_{r+1}=\varphi_2$ puisque par hypothèse inductive $\varphi_{r}=\varphi_2$, le morphisme de complexes $\varphi_{r+1}$ est alors surjectif et $d'_{r+1}=0 $. L'assertion (\ref{prop-critere-degen-b}) se démontre par un raisonnement dual.
\enddemo

\comment
\subsubsectionnumber 
Nos espaces topologiques $\Xg$ sont des pseudovariétés et en tant que telles, ils admettent des filtrations croissantes $\Xg=\mathop{\bigcup{\uparrow}}_{n\in\NN} U_{n}$ par des parties ouvertes $U_{n}$ qui sont cohomologiquement de type fini (\cf note ($^{\ref{tf}}$)). Les morphismes naturels suivants, respectivement de restriction et de prolongement par zéro,
$$\Hr(\Xg)\to\limproj_{n\in\NN}\Hr(U_{n})\text{\quad et\quad }
\limind_{n\in\NN} \Hc(U_{n})\to\Hc(\Xg)$$
sont alors des isomorphismes.

Maintenant, supposons donnée d'une fibration localement triviale $\pi:\Mg\to\Bg$ de fibre $F$ telle que les faisceaux $\cHpi(\Mg)$ sont constants donne lieu
\endcomment

\subsubsectionline{Les suites spectrales de Leray.}\label{suites-spectrales}
Le théorème \ref{autre-scindage} établit que pour un espace $i$-acyclique $\Xg$ localement connexe, les faisceaux
$$\halfdisplayskips\cHpi^{i}(\GaDelta_{?\ell,m}):=\IR^{i}\pi_{a!}(\fs k_{\GaDelta_{?\ell,m}})\,,\quad\forall i\in\ZZ\,,$$
sont constants sur les composantes connexes de $\Fa(\Xg)$.

Notons $c:\Fa\to\set\pt/$\glossary{${c:\Yg\to\set\pt/}$:application constante} l'application constante. Comme nous l'avons déjà rappelé dans \ref{c-mous}, la théorie des foncteurs dérivés donne des identifications
$$\IR\big(c_{\phi}\circ\pi_{a!}\big)(\fs k_{\GaDelta_{?\ell,m}})=
\IR c_{\phi}\big(\IR \pi_{a!}(\fs k_{\GaDelta_{?\ell,m}})\big)\,,\eqno
\hbox to1ex{\hss$\text{où }\phi\in\set {}\rm *,!/\,,$}$$
dont on tire les suites spectrales de Leray 
$$(\IE_{r}(\GaDelta_{?\ell,m})_{\varphi},d_r)\,,\eqno
\hbox to1ex{\hss$\text{où }\varphi\in\set \pi!,\rmc/\,.$}$$
Ces suites se décomposent en sommes directes de suites spectrales
$$\IE_{r}(\GaDelta_{?\ell,m})_{\varphi}=
\bigoplusnl_{\Cg\in\Pi_0\Fa}\IE_r(\Cg\GaDelta_{?\ell,m})_{\varphi}\,,
$$ 
où $\Cg$ est une composante connexe de $\Fa$.
Les termes de $\IE_2(\Cg\GaDelta_{?\ell,m})_{\varphi}$ sont
$$\IE_2(\Cg\GaDelta_{?\ell,m})_{\varphi}
:=\Hr_{\phi}(\Cg,\cHpi(\GaDelta_{?\ell,m}))\,,\qquad\eqno
\hbox to1em{\hss$\text{où }\Hr_{\phi}\in\set\Hr_{\pi!},\Hc/\,,$}$$
 et comme le faisceau $\cHpi(\GaDelta_{?\ell,m})$ est constant sur $\Cg$, on a respectivement
$$\left\{\mathalign{
\IE_{2}^{p,q}(\Cg\GaDelta_{?\ell,m})_{\rmc}:=\Hc^{p}(\Cg)\otimes\Hc^{q}(F_\Cg)&\Longrightarrow&
\Hc^{p+q}\big(\GaDelta_{?\ell,m}\big)\hfill\\\noalign{\kern4pt}
\smashbot{\IE_{2}^{p,q}(\Cg\GaDelta_{?\ell,m})_{\cv }:=\Hr^{p}(\Cg)\otimes\Hc^{q}(F_\Cg)}&\Longrightarrow&
\smashbot{\Hpi^{p+q}\big(\GaDelta_{?\ell,m}\big)}\,.\hfill\\
}\right.\quad
\postdisplaypenalty10000
\eqno\,\mllap{(\IE_2(\Cg))}$$
où $F_{\Cg}$ désigne la fibre de $\pia:\GaDelta_{?\ell,m}\to\Fa $ au-dessus de $\Cg$.

\begin{rema}[sur les notations]Concernant l'expression \expression{$F_{\Cg}$ est la fibre de $\pia$}\label{encore-independance-Hc-X-a}, on rappelle que bien que ces fibres ne sont généralement pas homéomorphes, nous les notons depuis \ref{PcUFm} par la notation $\Fb(\Xg\mmoins a)$. C'est une notation commode mais qui demande à être réinterprétée en fonction du contexte. 
Par exemple, les termes $\IE_2(\Fba(\Xg))_{\varphi}$ sont notés de manière succincte
$$\IE_2(\Ug\Fba(\Xg))_{\varphi}=\Hr_{\phi}(\Ug)\otimes\Hc(\Fb(\Xg\mmoins a))\,.\eqno(\ast)$$
Lorsque $\Ug\dans\Fa(\Xg)$ est connexe, la constance de $\cHpi(\Fba(\Xg))$ au dessus de $\Ug$ enlève toute ambiguïté à la notation $\Hc(\Fb(\Xg\mmoins a))$ (\ref{ameliore-independance-Hc-X-a}). Autrement, l'égalité ($\ast$) sous-entend de décomposer $\Ug$ en ses composantes connexes $\Cg$ et d'interpréter $\IE_2(\Ug\Fba(\Xg))_{\varphi}$ comme la somme $\bigoplus_{\Cg}\IE_2(\Cg\Fba(\Xg))_{\varphi}$ des formules $(\IE_2(\Cg))$ ci-dessus.

L'intérêt pratique de la notation ($\ast$) est que comme  $\Pc(\Fb(\Xg\mmoins a))$ est intrinsèque d'après \ref{PcUFm}, elle donne une description immédiate du polynôme de Poincaré de termes $\IE_2(\Ug\Fba(\Xg))_{\varphi}$.

\end{rema}
\subsubsectionline{Le cas où $\dim_k(\Hc(\Xg))<+\infty$.}\label{suites-spectrales-type-fini}
Lorsque $\Xg$ est $i$-acyclique et que $\dim\Hc(\Xg)<+\infty$, on a $\Hc(\Fm (\Xg\mmoins n))<+\infty$, pour tous $m,n\in\NN$, d'après \ref{prop-finitude-Hc}-(\ref{prop-finitude-Hc-b}). Les remarques \ref{PcUFm} et \ref{PcvUFm} s'appliquent alors au cas de la  fibration  $\pi_a:\Fba (\Xg)\to\Fa(\Xg)$ de fibre $\Fb(\Xg\mmoins\aa)$, de sorte que, pour tout ouvert $\Ug$ de $\Fa(\Xg)$, on a
\begin{itemize}\vskip1ex
\item[$(A)$] si $\Hc(\Ug)<+\infty$, alors 
$\Pc(\Ug\Fba(\Xg))=\Pc(\Ug)\cdot\Pc(\Fb(\Xg\mmoins a))$;
\item[$(B)$] si $\Hr(\Ug)<+\infty$, alors 
$\Pcv(\Ug\Fba(\Xg))=\P(\Ug)\cdot\Pc(\Fb(\Xg\mmoins a))$.
\end{itemize}

\smallskip
\noindent On en déduit, respectivement dans chaque cas, que:
\begin{itemize}\vskip1ex\halfdisplayskips
\item[$(A)$] si $\Hc(\Ug)<+\infty$, alors 
$$\dim_{k}\IE_2(\Ug\Fba)_{\rmc}=\Hc(\Ug\Fba)=\dim_{k}\IE_\infty(\Ug\Fba)_{\rmc}\,;$$
\item[$(B)$] si $\Hr(\Ug)<+\infty$, alors  
$$\dim_{k}\IE_2(\Ug\Fba)_{\cv}=\Hpi(\Ug\Fba)=\dim_{k}\IE_\infty(\Ug\Fba)_{\cv}\,;$$
\end{itemize}
et on conclut que les suites  $(\IE_r(\Ug\Fba)_{\rmc},d_r)$ et $(\IE_r(\Ug\Fba)_{\cv },d_r)$ sont dégénérées, car si jamais on avait $d_r\not=0$ pour un certain $r$, on aurait forcément 
$$\halfdisplayskips\dim(\IE_{2})\geq\dim(\IE_{r})>\dim(\IE_{r+1})\geq\dim(\IE_{\infty})\,,$$
ce qui n'est pas le cas. 

\begin{rema}Dans cette approche, la finitude joue un rôle essentiel mais nous verrons qu'elle n'est pas nécessaire. La proposition suivante s'affranchi des hypothèses de finitude sur $\Ug$, mais demande encore celle sur $\Xg$. Ce sera le théorème \ref{degen} qui aura l'énoncé le plus général, \idest sans aucune hypothèse de finitude sur $\Xg$.
\end{rema}

\begin{prop}\label{prop-ss-degen-type-fini}Soit $\Xg$ un espace $i$-acyclique, localement connexe et tel que $\dim_{k}\Hc(\Xg)<+\infty$. Alors, pour tout ouvert $\Ug\dans\Fa(\Xg)$,
les suites spectrales $(\IE_r(\Ug\Fba(\Xg))_{\rmc},d_r)$ et $(\IE_r(\Ug\Fba(\Xg))_{\cv },d_r)$ sont dégénérées.\label{degen-ouvert-type-fini}
\end{prop}

\demo Comme $\Fa(\Xg)$ localement connexe, il suffit de considérer le cas où $\Ug$ est connexe. Nous pouvons alors fixer une famille croissante $\UUU:=\set{\uparrow}\Ug_\mgoth/_{\mgoth\in\NN}$ d'ouverts connexes et de type fini qui recouvre $\Fa(\Xg)$. Pour chaque $\Ug_{\mgoth}\in\UUU$, on considère le produit fibré
$$
\xymatrixc{@R=6mm@C=1cm}{
\kern2cm&\mllap{\Ug_{\mgoth}\times_{\Ug}\Ug\Fba={}}\Ug_{\mgoth}\Fba\xylbl[rd]{$\boxtimes$}
\ar[d]_{\pia}\ar@{^{(}->}[r]|(0.53){\, j_{\mgoth}\,}&
\Ug\Fba \ar[d]^{\pia}\\
&\Ug_{\mgoth}\ar@{^{(}->}[r]|{\,j_{\mgoth}\,}&\Ug
}
$$

\sldash{La suite spectrale $\IE(\Ug\Fba)_{\rmc}$}Le cas où $\Hc(\Ug)<+\infty$ est celui déjà considéré dans \ref{suites-spectrales-type-fini}-$(A)$ où  $(\IE_r(\Fa)_{\rmc},d_r)$ est bien dégénérée.
Dans le cas général on a recours au recouvrement $\UUU$. Les morphismes canoniques
$$j_{\mgoth !}\,\IR \pi_{a!}\,\fs k_{\Ug_{\mgoth}\Fba}=\IR \pi_{a!}\,j_{\mgoth!}\,\fs k_{\Ug_{\mgoth}\Fba}\to \fs k_{\Ug\Fba}$$
induisent un morphisme de ``prolongement par zéro'' de suites spectrales
$$j_{\mgoth!}:(\IE_{r}(\Ug_{\mgoth}),d_r)\to(\IE_{r}(\Ug,d_r)\,,\eqno(\IE_{r})$$
qui s'identifie pour $r=2$, au morphisme de complexes
$$\Hc(\Ug_{\mgoth})\otimes\Hc(\Fb(\Xg\mmoins a))\hf{j_{\mgoth !}\otimes\id}{}{1.1cm}
\Hc(\Ug)\otimes\Hc(\Fb(\Xg\mmoins a))\,.\eqno(\IE_{2})$$
Maintenant, le fait que $\limind_{\mgoth} j_{\mgoth !}:\limind_{\mgoth}\Hc(\Ug_{\mgoth})=\Hc(\Ug)$ implique que le morphisme de complexes
$$\limind_{\mgoth}
j_{\mgoth!}\otimes\id:\limind_{\mgoth}(\IE_{2}(\Ug_{\mgoth}),d_2)\to(\IE_{2}(\Ug,d_2)$$
est bijectif, et alors, comme les suites $(\IE(\Ug_{\mgoth}),d_r)$ sont dégénérées, 
le critère de dégénérescence \ref{prop-critere-degen} s'applique
et   $(\IE_{r}(\Ug),d_r)$ est bien dégénérée.

\smallskip
\sldash{La suite spectrale $\IE(\Ug\Fba )_{\cv}$}Le cas où $\Hr(\Ug)<+\infty$ est celui déjà considéré dans \ref{suites-spectrales-type-fini}-$(B)$ où  $(\IE_r(\Fa)_{\cv},d_r)$ est bien dégénérée. Dans le cas général on a recours au recouvrement $\UUU$. 
Les morphismes canoniques
$$\IR \pi_{a!}\,\fs k_{\Ug\Fba}\to
\IR j_{\mgoth*}\,j_{\mgoth}^{-1}\,  \IR \pi_{a!}\,\fs k_{\Ug\Fba}
=
\IR j_{\mgoth*}\,  \IR \pi_{a!}\,\fs k_{\Ug_{\mgoth}\Fba}
$$
donnent les morphismes ``de restriction'' de complexes
$$\IR\Gamma(\Ug;\IR \pi_{a!}\,\fs k_{\Ug\Fba})\to
\IR\Gamma(\Ug_{\mgoth}; \IR \pi_{a!}\,\fs k_{\Ug_{\mgoth}\Fba})
$$
d'où les morphismes   de suites spectrales
$$(\IE_r(\Ug\Fba )_{\cv},d_r)\to(\IE_r(\Ug_{\mgoth}\Fba)_{\cv},d_r)\,,$$
et donc le morphisme de suites spectrales
$$(\IE_r(\Ug\Fba )_{\cv},d_r)\to\limproj_{\mgoth}(\IE_r(\Ug_{\mgoth}\Fba)_{\cv},d_r)\,.\eqno(\IE_r)$$
Or, comme l'application en homologie $\limind_{\mgoth}H_*(\Ug_{\mgoth},k)\to H_*(\Ug,k)$ est bijective, on a par dualité, que 
$\Hr(\Ug)=\limproj_{\mgoth}\Hr(\Ug_\mgoth)$ l'est également. On en déduit les identifications suivantes pour les termes $\IE_2$,
$$\mathrigid0mu
\mathalign{
\Ug\IE_2(\Fba )_{\cv}&&{}=\Hr(\Ug )\otimes\Hc(F)\hfill\\
&&\quad{}=\big(\limproj_{\mgoth}\Hr(\Ug_{\mgoth})\big)\otimes\Hc(F)=\limproj_{\mgoth}\big(\Hr(\Ug_{\mgoth}))\otimes\Hc(F)\big)\\
&&\qquad {}=\limproj_{\mgoth}\IE_2(\Ug_{\mgoth}\Fba)_{\cv}\,,\hfill}
\eqno(\IE_2)$$
où on a noté $\Hc(F):=\Hc(\Fg_{b}(\Xg\mmoins a))$ et où l'égalité de la ligne centrale est justifiée par le fait que $\dim_{k}\Hc(F)<+\infty$. Ces identifications sont compatibles aux différentielles $d_2$. Les conclusions de $(B)$ s'appliquent aux suites $(\IE_r(\Ug_{\mgoth}\Fba)_{\cv},d_r)$,  et alors $d_2=0$ dans $\IE_2(\Ug\Fba )_{\cv}$. Le critère de dégénérescence \ref{prop-critere-degen} s'applique et  $(\IE_r(\Ug\Fba )_{\cv},d_r)$ est bien dégénérée.
\enddemo

\comment
\begin{rema}\label{degen-ouvert-type-fini}Les mêmes raisonnements montrent pour la fibration $\pia:\Fba(\Xg)\to\Fa(\Xg)$, que si pour un ouvert $\Ug\dans\Fa(\Xg)$, on note $\Ug\Fba(\Xg):=\pia^{-1}(\Ug)$, les suites spectrales pour les foncteurs $\Hc(\_)$ et $\Hpi(\_)$ associées à la fibration localement triviale $\pia:\Ug\Fba(\Xg)\to\Ug$, sont également dégénérées.
\end{rema}
\endcomment
\subsubsectionline{Le cas général.}\label{suites-spectrales-general}Nous nous affranchissons maintenant de l'hypothèse de finitude pour $\Hc(\Xg)$ et montrons la dégénérescence des suites spectrales en question en nous appuyant sur la proposition précédente \ref{prop-ss-degen-type-fini}.

\begin{theo}
\label{degen}Soient $a\leq\ell\leq m\in\NN$. Soit $\Xg$ un espace topologique $i$-acyclique et localement connexe. Pour tout ouvert $\Ug\dans\Fa(\Xg)$, les applications $\pi_a:\Ug\GaDelta_{?\ell,m}\to\Ug$  donnent lieu à des suites spectrales de Leray pour les cohomologies $\Hc(\_)$ et $\Hpi(\_)$
notées respectivement
$$
(\IE_r(\Ug\GaDelta_{?\ell,m})_{\rmc},d_r)
\text{\quad et\quad }
(\IE_r(\Ug\GaDelta_{?\ell,m})_{\cv},d_r)
\eqno(\IE_r)\,.$$
telles que
$$\begin{cases}\noalign{\kern-4pt}
\bigoplus_{\Cg\in\Pi_{0}\Ug}\IE_2^{p,q}(\Cg\GaDelta_{?\ell,m})_{\rmc}=\Hc^{p}(\Cg)\otimes\Hc^{q}(F_{\Cg})\Rightarrow \Hc^{p+q}(\GaDelta_{?\ell,m})\\\noalign{\kern2pt}
\bigoplus_{\Cg\in\Pi_{0}\Ug}\IE_2^{p,q}(\Cg\GaDelta_{?\ell,m})_{\cv}=\Hr^{p}(\Cg)\otimes\Hc^{q}(F_{\Cg})\Rightarrow \Hpi^{p+q}(\GaDelta_{?\ell,m})\\\noalign{\kern-1pt}
\end{cases}\postskip7pt$$
où $\Cg$ désigne une composante connexe de $\Ug$ et où
$F_{\Cg}$ désigne une fibre quelconque de $\pi_a$ au-dessus de $\Cg$. De plus, les suites spectrales $(\IE_{r})$ sont dégénérées.
\end{theo}
\demo Compte tenu des résultats qui précèdent, il ne nous reste qu'à justifier la dégénérescence des suites spectrales, ce pour quoi nous allons nous limiter aux ouverts $\Ug$ connexes. Les faisceaux $\cHpi(\GaDelta_{?\ell,m})$ sont donc constants sur $\Ug$.

On fixe une famille croissante $\VVV:=\set{\uparrow}\Vg_\ngoth/_{\ngoth\in\NN}$ d'ouverts $\Vg_{\ngoth}$ (pas forcément connexes) de type fini qui recouvre $\Xg$.

\smallskip\bfsldash{Le cas des  fibrations $\pi_a:\Ug\Fba(\Xg) \to\Ug\dans\Fa(\Xg)$}

\sldash{Dégénérescence de $(\Ug\IE_r(\Fba)_{\rmc,\cv},d_r)$}Compte tenu de la proposition \ref{prop-ss-degen-type-fini}, nous avons juste à vérifier le cas où $\dim_{k}(\Hc(\Xg))=+\infty$. 

Pour chaque $\mgoth\in\NN$, notons $\Ug_{\mgoth}:=\Fa(\Vg_{\mgoth})\cap\Ug$. 
 La famille $\UUU:=\set \Ug_{\mgoth}/_{\mgoth\in\NN}$ est clairement un recouvrement ouvert croissant de $\Ug$. On remarquera aussi que l'on des inclusions ouvertes 
$$\halfdisplayskips\Ug_{\mgoth}\dans\Fa(\Vg_{\mgoth})\dans\Fa(\Vg_{\ngoth})\dans\Fa(\Xg)$$
pour tous $\mgoth\leq\ngoth$.

Considérons ensuite les diagrammes commutatifs où $\ngoth\geq\mgoth$
$$
\xymatrixc{@R=6mm@C=1.cm}{
\Fba(\Vg_{\ngoth})\ar[d]_{\pia}\xylbl[rd]{(I)}
\ar@{<-^{)}}[r]
&
\Ug_{\mgoth}\Fba (\Vg_{\ngoth})\ar[d]_{\pia}\ar@{^{(}->}[r]|{\;j_{\mgoth,\ngoth}\;}
\xylbl[rd]{(II)}&
\Ug_{\mgoth}\Fba (\Xg)\ar[d]_{\pia}
\ar@{^{(}->}[r]|(0.55){\;j_{\mgoth}\;}\xylbl[rd]{(III)}&
\Ug\Fba (\Xg)\ar[d]^{\pia}\\
\Fa(\Vg_{\ngoth})
\ar@{<-^{)}}[r]
&
\Ug_{\mgoth}\ar@{=}[r]&\Ug_{\mgoth}\ar@{^{(}->}[r]|(0.5){\;j_{\mgoth}\;}&
\Ug
}$$

Les sous-diagrammes (I) et (III) sont des produits fibrés et les cohomologies à support compact des fibres de $\pia$ y sont (constantes car $\Ug$ est connexe) respectivement $\Fb(\Vg_{\ngoth}\mmoins a)$ et $\Fb(\Xg\mmoins a)$.

Dans le sous-diagramme (II), la base est constante et seul les fibres changent. La famille croissante de plongements ouverts $\set j_{\mgoth,\ngoth}/_{\ngoth}$  recouvre $\Ug_{\mgoth}\Fba(\Xg)$ et
les morphismes  de complexes
$$
\limind_{\ngoth\in\NN}\IR c_{\phi}\, \IR\pi_{a!}\, j_{\mgoth,\ngoth!}\,\fs k_{\Ug_{\mgoth}\Fba(\Vg_{\ngoth})}
\to
\IR c_{\phi}\, \IR\pi_{a!}\, \fs k_{\Ug_{\mgoth}\Fba(\Xg)}
$$
où $c:\Ug_{\mgoth}\to\set\pt/$ est l'application constante et $\phi\in\set *,!/$,
induisent alors les morphismes de suites spectrales
$$\limind_{\ngoth\in\NN}\big(\IE_{r}(\Ug_{\mgoth}\Fba(\Vg_{\ngoth}))_{\varphi},d_r\big)
\to 
\big(\IE_{r}(\Ug_{\mgoth}\Fba(\Xg))_{\varphi},d_r\big)
\eqno(\IE_r)$$
où $\varphi\in\set \cv,\rmc /$, qui s'identifient pour $r=2$, au morphisme
$$\mathrigid0mu
\limind_{\ngoth\in\NN}\big(\Hr_{\phi}(\Ug_{\mgoth})\otimes\Hc(\Fb(\Vg_{\ngoth}\mmoins \aa))\hf{\id\otimes j_{\ngoth!}}{}{1cm}
\Hr_{\phi}(\Ug_{\mgoth})\otimes\Hc(\Fb(\Xg\mmoins \aa))\big)\,,\eqno(\IE_2)$$
où $\Hr_{\phi}\in\set H,\Hc/$ et où
$j_{\ngoth}:\Hc(\Vg_{\ngoth})\to\Hc(\Xg)$ est le prolongement par zéro.
Or, comme la famille $\set \Fa(\Vg_{\ngoth})/_{\ngoth\geq\mgoth\in\NN} $ recouvre $\Fb(\Xg\mmoins\aa)$, la limite inductive des morphismes $(\id\otimes j_{\ngoth!})$ est un isomorphisme. 

Enfin, les suites spectrales $(\IE_{r}(\Ug_{\mgoth}\Fba(\Vg_{\ngoth}))_{\varphi},d_r)$ sont dégénérées d'après \ref{degen-ouvert-type-fini} puisque $\dim\Hc(\Vg_{\ngoth})<+\infty$. Le critère de dégénérescence \ref{prop-critere-degen} s'applique et les suite spectrales $(\IE_r(\Ug_{\mgoth}\Fba(\Xg) )_{\varphi},d_r)$ sont bien dégénérées.

\medskip
Revenons maintenant au sous-diagramme (III). Ici, les fibres sont fixes et égales à $\Fb(\Xg\mmoins a)$ mais la base change. Aussi, nous avons deux situations différentes à considérer.

\smallskip
\sldash{La suite spectrale $(\IE_r(\Ug\Fba(\Xg))_{\rmc},d_r)$} Les morphismes  de complexes des prolongements par zéro
$$
\limind_{\mgoth}\IR c_{!}\, \IR\pi_{a!}\, j_{\mgoth!}\,\fs k_{\Ug_{\mgoth}\Fba(\Xg)}
\to
\IR c_{!}\, \IR\pi_{a!}\, \fs k_{\Ug\Fba(\Xg)}\,,
$$
induisent un morphisme de suites spectrales
$$
\limind_{\mgoth}(\IE_{r}(\Ug_{\mgoth}\Fba(\Xg))_{c},d_r)
\to
(\IE_{r}(\Ug\Fba(\Xg))_{c},d_r)\,,
\eqno(\IE_r)$$
qui s'identifie pour $r=2$, au morphisme
$$\mathrigid1mu
\limind_{\mgoth}\big(\Hc(\Ug_{\mgoth})\otimes\Hc(\Fb(\Xg\mmoins \aa))\big)
\hf{\limind_{\mgoth} j_{\mgoth!}\otimes\id}{}{1.8cm}
\Hc(\Ug)\otimes\Hc(\Fb(\Xg\mmoins \aa))\,,\eqno(\IE_2)$$
clairement bijectif puisque $\set\Ug_{\mgoth}/$ recouvre $\Ug$. Ce fait, et la dégénérescence déjà établie des suites spectrales $(\IE_{r}(\Ug_{\mgoth}\Fba(\Xg))_{c},d_r)$ font que le critère
de dégénérescence \ref{prop-critere-degen} s'applique et la suite spectrale $(\IE_r(\Ug\Fba(\Xg) )_{c},d_r)$ est dégénérée.

\smallskip
\sldash{La suite spectrale $(\IE_r(\Ug\Fba(\Xg))_{\cv},d_r)$}On commence par rappeler le morphisme canonique de restriction de complexes dans $D^{+}(\Ug)$
$$\IR\pi_{a!}\, \fs k_{\Ug\Fba(\Xg)}\to
\IR\pi_{a!}\, \IR j_{\mgoth*}\,\fs k_{\Ug_{\mgoth}\Fba(\Xg)}=
\IR j_{\mgoth*}\,\IR\pi_{a!}\, \fs k_{\Ug_{\mgoth}\Fba(\Xg)}\,,
$$
où l'égalité est justifiée puisque (III) est un produit fibré. On en déduit le morphisme de complexes
$$\halfdisplayskips
\IR \Gamma(\Ug;\IR\pi_{a!}\, \fs k_{\Ug\Fba(\Xg)})
\to
\IR \Gamma(\Ug_\mgoth;\IR\pi_{a!}\,\fs k_{\Ug_{\mgoth}\Fba(\Xg)})\,,
$$
qui donne lieu au morphisme de suites spectrales
$$
(\IE_{r}(\Ug\Fba(\Xg))_{\cv},d_r)
\to
\limproj_{\mgoth}(\IE_{r}(\Ug_{\mgoth}\Fba(\Xg))_{\cv},d_r)\,,
\eqno(\IE_{r})$$
qui s'identifie pour $r=2$, au morphisme
$$\mathrigid1mu
\Hr(\Ug)\otimes\Hc(\Fb(\Xg\mmoins \aa))
\hf{\limproj_{\mgoth} j_{\mgoth}^{*}\otimes\id}{}{2cm}
\limproj_{\mgoth}\big(\Hr(\Ug_{\mgoth})\otimes\Hc(\Fb(\Xg\mmoins \aa))\big)\,.\eqno(\IE_2)$$
Or, comme $\set\Ug_{\mgoth}/$ recouvre $\Ug$, on a $\Hr(\Ug)=\limproj_{\mgoth}\Hr(\Ug_{\mgoth})$ et ($\IE_2$)  est injectif. Nous savons d'autre part, que les suites spectrales $(\IE_{r}(\Ug_{\mgoth}\Fba(\Xg))_{\cv},d_r)$ sont dégénérées. Le critère
de dégénérescence \ref{prop-critere-degen} s'applique et la suite spectrale $(\IE_{r}(\Ug\Fba(\Xg))_{\cv},d_r)$ est dégénérée.

\comment\color{blue}
\hrule depth1cm
plongements ouverts définis par les produits fibrés $\Fa(\Vg_{\mgoth})\times_{\Fa(\Xg)}\Fba(\Xg)$
$$\xymatrixc{@R=6mm}{
\kern2cm&\mllap{\Fa(\Vg_{\mgoth})\times_{\Fa(\Xg)}\Fba(\Xg)={}}\pia^{-1}\Fa (\Vg_{\mgoth})\ar[d]_{\pia\mrlap{\kern1.3cm\boxtimes}}\ar@{^{(}->}[r]^(0.53){j_{\mgoth}}&\Fba (\Xg)\ar[d]^{\pia}\\
&\Fa(\Vg_{\mgoth})\ar@{^{(}->}[r]^{j_{\mgoth}}&\Fa(\Xg)
}\eqno(\VVV\VVV)$$
Ces plongements donnent lieu à des morphismes naturels de restriction 
$$
\IR\pi_{a!}\, \fs k_{\Fba(\Xg)}
\to
\IR\pi_{a!}\, \IR j_{\mgoth*}\,j_{\mgoth}^{-1}\fs k_{\Fba(\Xg)}
=\IR j_{\mgoth*}\,\IR\pi_{a!}\,j_{\mgoth}^{-1}\fs k_{\Fba(\Xg)}
$$
et aux morphismes de complexes
$$\IR\Gamma\big(\Fa(\Xg); \IR\pi_{a!}\, \fs k_{\Fba(\Xg)}\big)
\to
\IR\Gamma\big(\Fa(\Vg_{\mgoth}); \IR\pi_{a!}\,j_{\mgoth}^{-1}\fs k_{\Fba(\Xg)}\big)\,,\eqno(\diamond\diamond)$$
induisant des morphisme ``de restriction'' de suites spectrales
$$
(\IE_{r}(\Fa(\Xg))_{\cv},d_r)
\to 
(\IE_{r}(\Fa(\Vg_{\mgoth}))'_{\cv},d_r)
$$
où nous avons noté $(\IE_{r}(\Fa(\Vg_{\mgoth}))'_{\cv},d_r)$ la suite spectrale associée au terme de droite de $(\diamond\diamond)$. Ce morphisme se lit pour $r=2$ comme le morphisme
$$\Hr(\Fa(\Xg))\otimes\Hc(\Fb(\Xg\mmoins \aa))\hf{j_{\mgoth}^{*}\otimes \id}{}{1.2cm}
\Hr(\Fa(\Vg_{\mgoth}))\otimes\Hc(\Fb(\Xg\mmoins \aa))\,.\eqno(**)$$
où $j_{\mgoth}^{*}:\Hr(\Fa(\Xg))\to\Hr(\Fa(\Vg_{\mgoth}))$ est le morphisme de restriction. Or, comme $\Hr(\Fa(\Xg))=\limproj_{\Vg_{\mgoth}\in\VVV}\Hr(\Fg(\Vg_{\mgoth}))$, on comprend que la dégénérescence de $(\IE_{r}(\Fa(\Xg))_{\cv},d_r)$ résultera de celles de $(\IE_{r}(\Fa(\Vg_{\mgoth}))'_{\cv},d_r)$, pour tout $\mgoth\in\NN$. Ce que nous montrons à continuation.

Pour $\mgoth\in\NN$ donné, la fibration $\pia:\pia^{-1}\Fa(\Vg_{\mgoth})\to \Fa(\Vg_{\mgoth})$ dans $(\VVV\VVV)$ est à fibres de la forme $\Fb(\Xg\mmoins\aa)$ dont la cohomologie à support compact peut ne pas être de dimension finie. L'idée alors est de faire intervenir les restrictions des fibrations $\pia:\Fba(\Vg_{\ngoth})\to \Fa(\Vg_{\ngoth})$, pour $\ngoth\geq\mgoth$, à l'ouvert
$\Fa(\Vg_{\mgoth})\dans\Fa(\Vg_{\ngoth})$. Autrement dit, on fait intervenir les produits fibrés
$\Fa(\Vg_{\mgoth})\times_{\Fa(\Vg_{\ngoth})}\Fba(\Vg_{\ngoth})$:
$$\def\Xg{\Vg_{\ngoth}}
\xymatrixc{@R=6mm}{
\kern2cm&\mllap{\Fa(\Vg_{\mgoth})\times_{\Fa(\Xg)}\Fba(\Xg)={}}\pia^{-1}\Fa (\Vg_{\mgoth})\ar[d]_{\pia\mrlap{\kern1.3cm\boxtimes}}\ar@{^{(}->}[r]^(0.53){j_{\mgoth}}&\Fba (\Xg)\ar[d]^{\pia}\\
&\Fa(\Vg_{\mgoth})\ar@{^{(}->}[r]^{j_{\mgoth}}&\Fa(\Xg)
}$$
On considère alors les plongement ouverts 
$$\xymatrixc{@C=5mm}{\Fa(\Vg_{\mgoth})\times_{\Fa(\Vg_{\ngoth})}\Fba(\Vg_{\ngoth})\ar[rr]^{j_{\ngoth}}\ar[rd]_{\pia}
&&
\Fa(\Vg_{\mgoth})\times_{\Fa(\Xg)}\Fba(\Xg)\ar[ld]^{\pia}\\
&\Fa(\Vg_{\mgoth})}
\eqno(\diamond{\diamond}\diamond)$$
qui induisent un morphisme ``de prolongement par zéro'' de suites spectrales
$$
(\IE_{r}(\Fa(\Vg_{\mgoth}))''_{\cv},d_r)
\too^{j_{\ngoth !}} 
(\IE_{r}(\Fa(\Vg_{\mgoth}))'_{\cv},d_r)
$$
où nous avons noté $(\IE_{r}(\Fa(\Vg_{\mgoth}))''_{\cv},d_r)$ la suite spectrale associée au terme de gauche de $(\diamond{\diamond}\diamond)$. Ce morphisme se lit pour $r=2$ comme le morphisme
$$\Hr(\Fa(\Vg_{\mgoth}))\otimes\Hc(\Fb(\Vg_{\ngoth}\mmoins \aa))\hf{\id\otimes j_{\ngoth !}}{}{1.2cm}
\Hr(\Fa(\Vg_{\mgoth}))\otimes\Hc(\Fb(\Xg\mmoins \aa))\,.\eqno(*{*}*)$$
où $j_{\ngoth !}:\Hr(\Fb(\Vg_{\ngoth}\mmoins\aa))\to\Hr(\Fb(\Xg\mmoins\aa))$ est le prolongement par zéro.

Or, comme d'après la remarque \ref{degen-ouvert-type-fini}
la suite spectrale  $(\IE_{r}(\Fa(\Vg_{\mgoth}))''_{\cv},d_r)$ est dégénérée, et que 
$\limind_{\Vg_{\ngoth}\in\VVV}\Hr(\Fb(\Vg_{\ngoth}\mmoins\aa))\to\Hr(\Fb(\Xg\mmoins\aa))$ est un isomorphisme, la différentielle $d_2$ dans  $\IE_{2}(\Fa(\Vg_{\mgoth}))'_{\cv}$ est nulle et
$\IE_{3}(\Fa(\Vg_{\mgoth}))'_{\cv}=\IE_{2}(\Fa(\Vg_{\mgoth}))'_{\cv}$. On montre ainsi, une fois de plus par récurrence que $d_r=0$ dans $(\IE_{r}(\Fa(\Vg_{\mgoth}))'_{\cv},d_r)$ et cette suite est bien dégénérée.
\endcomment

Ceci termine la preuve du cas des fibrations $\pia:\Ug\Fba(\Xg)\to\Ug\dans\Fa(\Xg)$.

\smallskip
\bfsldash{Le cas des fibrations $\pi_a:\Ug\GaDelta_{\ell}\Xg^{m}\to\Ug\dans\Fa(\Xg)$}

{\parskip1pt\mou\nobreak\leavevmode}Résulte du cas précédent, en raison de la décomposition ouverte (\footnote{\Cf fin de la démonstration du théorème \ref{autre-scindage}.})
$$\halfdisplayskips\GaDelta_{\ell}\Xg^{m}=\coprod\nolimits_{\pgoth}\Fg_{\pgoth}$$
où $\pgoth$ est une partition  $\pi_a$-adaptée de $\Pgoth_{\ell}(m)$, et où $\Fg_{\pgoth}\simeq\Fg_{\ell}$.

\smallskip
\bfsldash{Le cas des fibrations $\pi_a:\Ug\GaDelta_{\leq\ell}\Xg^{m}\to\Ug\dans\Fa(\Xg)$}

{\parskip1pt\mou\nobreak\leavevmode}Notons $j:\Ug\GaDelta_{\ell,m}\to\Ug\GaDelta_{\leq\ell,m}$ l'inclusion ouverte.
Le morphisme de prolongement par zéro dans $D^{+}(\Ug)$
$$j_{!}\,\IR\pi_{a!}\,\fs k_{\Ug\GaDelta_{\ell,m}}=\IR\pi_{a!}\,j_{!}\,\fs k_{\Ug\GaDelta_{\ell,m}}\to
\IR\pi_{a!}\,\fs k_{\Ug\GaDelta_{\leq\ell,m}}$$
donne lieu aux morphismes de complexes
$$\IR \Gamma_{\phi}\big(\Ug;\IR\pi_{a!}\,\fs k_{\Ug\GaDelta_{\ell,m}}\big)\to
\IR \Gamma_{\phi}\big(\Ug;\IR\pi_{a!}\,\fs k_{\Ug\GaDelta_{\leq\ell,m}}\big)\,,$$
avec $\Gamma_{\phi}\in\set \Gamma,\Gammac/$, qui induisent des morphismes de suites spectrales 
$$\big(\IE_{r}(\Ug\GaDelta_{\ell,m})_{\varphi},d_r\big)
\to 
\big(\IE_{r}(\Ug\GaDelta_{\leq\ell,m})_{\varphi},d_r\big)\,,
\eqno(\IE_r)$$
avec $\varphi\in\set \cv,\rmc/$, qui s'identifient pour $r:=2$, aux morphismes
$$\Hr_{\phi}(\Ug)\times\Hc(F_{\ell})
\hf{\id\otimes j_{!}}{}{1cm}
\Hr_{\phi}(\Ug)\times\Hc(F_{\leq\ell})
\eqno(\IE_2)$$
avec $\let\Gamma\Hr\Gamma_{\phi}\in\set \Gamma,\Gammac/$, et où $F_{?\ell}$ désigne la fibre de $\pia:\Ug\GaDelta_{?\ell,m}\to\Ug$, et  
$$j_{!}:\Hc(F_{\ell})\to\Hc(F_{\leq\ell})\eqno(\dagger)$$ est le prolongement par zéro.
Or, le théorème \ref{theo-scindage-basee}-(\ref{theo-scindage-basee-b}) établit la surjectivité du prolongement par zéro
$$
\Hpi(\Vg\GaDelta_{\ell,m})\too
\Hpi(\Vg\GaDelta_{\leq \ell,m})
\eqno(\ddagger)$$
pour tout $\Vg\dans\Fa$. Si maintenant $\Vg$ parcours une base de voisinages connexes d'un point, la proposition \ref{germe-compact} 
assure la surjectivité de $(\dagger)$, et donc aussi celle des morphismes ($\IE_2$). Comme nous avons déjà établi la dégénérescence de $(\IE_{r}(\Ug\GaDelta_{\leq\ell,m})_{\varphi},d_r)$, le critère
 \ref{prop-critere-degen} s'applique et $(\IE_{r}(\Ug\GaDelta_{\leq\ell,m})_{\varphi},d_r)$ est aussi dégénérée.
Ceci termine la preuve du théorème.
\enddemo

\subsubsection{Sur l'action de $\Sb \times \Sa $ sur la suite spectrale de Leray}
Reprenons la discussion de \ref{comparaison}. Pour $a,b\in\NN$, réalisons $\Sb \times \Sa $ comme le sous-groupe de $\S_{b+a}$ des permutations laissant stables les sous-intervalles $\iii[1,b]$ et $\iii[b+1,b+a]$ et faisons-le agir sur $\Fba $ par son action naturelle à travers $\S_{b+a}$.
Faisons-le ensuite agir sur $\Fa $, à travers de la projection sur $\Sa $. La projection (\ref{notas}-(\ref{nota-pi-a}))
$$
\pi_a:\Fba (\Xg)\to\Fa (\Xg)\,,
$$
est alors $\Sb \times \Sa $-équivariante. 

\soustitre{Action de $\Sb $ sur $\IE_{2}^{p,q}(\Cg)$}Notons par $\varPhi_h:\Fba \to\Fba $ l'homéomorphismes de l'action de $h\in\Sb $. On a $\pia \circ\varPhi_h=\pia $. L'action de $\varPhi_{h}$ respecte donc les fibres de $\pi_a$ et l'action induite sur $\Hc(\Cg,\cHpi^{i}(\Fa ))$ correspond à l'action de $\Sb $ sur le second facteur du produit tensoriel $\Hc(\Cg)\otimes\Hc(F_{\Cg})$.

\soustitre{Action de $\Sa $ sur $\IE_{2}^{p,q}(\S _{a}\Cg)$}Les exemples de la proposition \ref{compare-b+a}, montrent bien que l'action de $\S _a$
sur $\Hc(\Sa \Cg,\cHpi^{i}(\Fa ))$\glossary{${\Sa \Cg:=\bigcup_{g\in\Sa }g\cdot\Cg}$:saturé d'une partie $\Cg\dans\Fa$ sous l'action de $\Sa $} ne correspond pas à l'action de $\Sa $ sur le premier facteur du produit tensoriel $\Hc(\Sa \Cg)\otimes\Hc(F_{\Cg})$, et ne se voit donc pas sur la suite spectrale de Leray.

\comment
Notons par $\varPhi_g:\Fba \to\Fba $ et $\varphi_g=\Fa \to\Fg_a$ les homéomorphismes de l'action de $g\in\Sa $. 
Fixons une résolution injective $\fs k_{\Fba }\to\I^{\mb}$ et un quasi-isomorphisme $\xi_{g}:\varPhi_{g}^{-1}\I^{\mb}\to\I^{\mb}$. Notons ensuite (\footnote{Remarquer pour la suite que $\varPhi_{g!}=\varPhi_{g*}$ puisque $\varPhi_g$ est propre et que $\varPhi_{g*}$ est exact
puisque $\varPhi_{g*}=\varPhi_{g^{-1}*}^{-1}$. Idem pour $\varphi_g$.})
$$\Xi_g:=\varPhi_{g*}(\xi_g):(\I^{\mb}=\varPhi_{g*}\varPhi_{g}^{-1}\I^{\mb})\to\varPhi_{g*}\I^{\mb}\,.$$

L'action $g$ sur $\Hc^{i}(\Fba )$ est celle induite par l'action $\Gammac(\Fba ;\Xi_{g})$ sur le complexe $\Gammac(\Fba ,\I^{\mb})$, et comme on a par ailleurs
$$\IR\Gammac(\Fba ;\fs k_{\Fba })=\Gammac(\Fa ;\pi_{a!}\I^{\mb})\,,$$ cette action se voit aussi comme celle de $\Gammac(\Fa ;\pi_{a!}\Xi_{g})$, où le morphisme
$$\pi_{a!}\Xi_g:
\pi_{a!}\I^{\mb}
\to\pi_{a!}\varPhi_{g*}\I^{\mb}\,,$$
équivaut, par équivariance de $\pia $, à un morphisme
$$\pi_{a!}\Xi_g:
\pi_{a!}\I^{\mb}
\to
\varphi_{g*}\pi_{a!}\I^{\mb}\,,\eqno(\ast)$$
qui induit, en reprenant les notations de \ref{cHcv}, les morphismes de faisceaux
$$h^{i}(\pi_{a!}\Xi_g):
\cHpi^{i}(\Fba )
\to\varphi_{g*}\cHpi^{i}(\Fba )\,,$$
puisque $\varphi_{g*}$ est exact.

Comme d'après \ref{autre-scindage}, le faisceau $\cHpi^{i}(\Fba )$ est constant sur une composante connexe $\Cg\dans\Fa $, on comprend que l'action que $h^{i}(\pi_{a!}\Xi_g)$ induit sur 
$\Hc(\Sa \Cg,\cHpi^{i}(\Fa ))$\glossary{${\Sa \Cg:=\bigcup_{g\in\Sa }g\cdot\Cg}$:saturé de la composante $\Cg\dans\Fa$ sous l'action de $\Sa $} est très précisément l'action de $\Sa $ sur le premier facteur du produit tensoriel $\Hc(\Sa \Cg)\otimes\Hc(F_{\Cg})$.
\endcomment

\soustitre{Conclusion}Ces remarques montrent que seule l'action $\Sb \times\1_a$ sur $\Hc(\Fba )$ et $\Hpi(\Fba )$ se voit sur les termes $\IE_2$ des suites spectrales de Leray (\cf\ref{degen})  $$\begin{cases}\noalign{\kern-4pt}
\bigoplus_{\Cg}\IE_2^{p,q}(\Cg):=\Hc^{p}(\Cg)\otimes\Hc^{q}(F_{\Cg})\Rightarrow \Hc^{p+q}(\Fba )\,,\\\noalign{\kern4pt}
\bigoplus_{\Cg}\IE_2^{p,q}(\Cg):=\Hr^{p}(\Cg)\otimes\Hc^{q}(F_{\Cg})\Rightarrow \Hpi^{p+q}(\Fba )\,.\\\noalign{\kern-2pt}
\end{cases}$$
Elle se voit comme étant l'action de $\Sb $ sur $\Hc(F_{\Cg})=\Hc(\Fg_{b}(\Xg\mmoins a))$. 

\comment
La proposition suivante est la conséquence immédiate de cette observation.\vskip-1ex\vskip0pt

\begin{prop}\label{SbxSa}Soit $\Xg$ une variété topologique $i$-acyclique. \'Etant donné des sous-groupes $H_a\dans\Sa $ et $H_{b}\dans\Sb $, on a
$$\def\decale#1#2{\kern#1#2\kern-#1}
\mathalign{\Pc\big(\Fba (\Xg)^{H_{b}\times H_a}\big)&=&
\Pc\big(\Fa (\Xg)^{H_a}\big)\hfill
&\decale{-2mm}\cdot&
\Pc\big(\Fg_{b}(\Xg\mmoins a)^{H_{b}}\big)\\
&=&
\Pc\big(\Fa (\Xg\mmoins b)^{H_a}\big)
&\cdot&
\Pc\big(\Fg_{b}(\Xg)^{H_{b}}\big)\,.\hfill
}$$
\end{prop}
\endcomment

\section{Questions diverses sur les espaces $i$-acycliques}
Nous\label{appendice-exemples} rassemblons ici des exemples d'espaces $i$-acycliques ainsi que des contre-exemples à certaines propriétés de l'$i$-acyclicité que l'on aurait souhaité avoir mais qui ne sont pas vérifiées en toute généralité.

\secskip=1.25ex\mou

\comment
\subsection{Exemples élémentaires d'espaces $i$-acycliques}Conséquence de la définition, mais aussi des propriétés élémentaires justifiées dans \ref{prop-acycliques}.
\def\varlistskips{\leftmargin1.3em}\begin{enumerate}\itemindent1em
\labelsep0.5ex\itemsep0pt\mou\sl
\item Tout espace contractile non compact est $i$-acyclique.
\item Tout fibré vectoriel $p:\Eg\to\Bg$ de variétés différentielles
dont la dimension de la base $\Bg$ est strictement majoré par le rang de $\Eg$ est $i$-acyclique.
\item Si $\Xg$ est $i$-acyclique, les espaces de configuration généralisés $\Delta_{?\ell}\Xg^{m}$ sont aussi $i$-acycliques.
\demo (a) et (b) sont évidentes. (c) Si $\Xg$ est $i$-acyclique, $\Xg^{m}$ l'est aussi (\ref{prop-acycliques}-(\ref{prop-acycliques-d})), et donc aussi l'ouvert $\Fm (\Xg)\dans\Xg^{m}$ (\ref{prop-acycliques}-(\ref{prop-acycliques-c})). La décomposition de $\Delta_{\ell}\Xg^{m}$ en réunion disjointe ouverte d'espaces homéomorphes à $\Delta_{\ell}(\Xg)$ (\ref{connexes}) montre qu'il est, lui aussi, $i$-acyclique. 
\enddemo
\end{enumerate}
\endcomment

\subsection{Espaces de configuration généralisés}
\proclaim{Si\label{config-gener} $\Xg$ est $i$-acyclique, l'espace  $\Delta_{?\ell}\Xg^{m}$ est $i$-acyclique pour  $0<\ell\leq m$. 
}
\demo
Si $\Xg$ est $i$-acyclique, $\Xg^{m}$ l'est aussi (\ref{prop-acycliques}-(\ref{prop-acycliques-d})), et donc aussi l'ouvert $\Fm (\Xg)\dans\Xg^{m}$ (\ref{prop-acycliques}-(\ref{prop-acycliques-c})) de même que tout espace $\Delta_{\ell}\Xg^{m}$, car réunion disjointe ouverte d'espaces homéomorphes à $\Fg_{\ell}(\Xg)$ (\ref{connexes}).

\def\Up{U_{\pgoth}}
\def\Kp{K_{\pgoth}}
\def\dUp{\partial U_{\pgoth}}
\def\dKp{\partial K_{\pgoth}}
\def\clUp{\cl U_{\pgoth}}

Pour établir la $i$-acyclicité des espaces $\Delta_{\leq\ell}\Xg^{m}$, nous utilisons le critère \ref{caracterisation}-(\ref{caracterisation-b}) selon lequel un espace $\Yg$ est $i$-acyclique lorsque pour tout compact $K\dans\Yg$, la
restriction $\rho_{K}:\Hc(\Yg)\to\Hr(K)$ est nulle. 
Dans le cas présent, comme en plus le prolongement par zéro $j_{!}:\Hc(\Delta_{\ell}\Xg^{m})\to\Hc(\Delta_{\leq\ell}\Xg^{m})$
est surjectif (\ref{theo-scindage}-(\ref{theo-scindage-b})) et que $\Delta_{\ell}\Xg^{m}$ est réunion disjointe des ouverts
$\Up:=\Fg_{\pgoth}(\Xg)$, avec $\pgoth\in\Pgoth_{\ell}(m)$ (\cf\ref{nota-pgoth}), le critère en question sera validé par la nullité, pour tous $\pgoth$ et $K$,
 des composées 
$$
\Hc(\Up)\too^{j_{!}} \Hc(\Delta_{\leq\ell}\Xg^{m})\too^{\rho_{K}}\Hr(K)\,.\eqno(\diamond)$$

\begingroup\setbox222=\vtop to0pt{\kern-4.5mm\hbox{\qquad 
\includegraphics{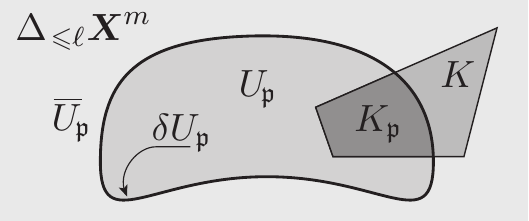}}\vss}
\dimen222=\hsize\advance\dimen222 by -\wd222

\parshape=7
0pt \dimen222
0pt \dimen222
0pt \dimen222
0pt \dimen222
0pt \dimen222
0pt \dimen222
0pt \hsize
Pour\vadjust{\hbox to\hsize{\hfill\copy222}} $\pgoth$ et $K$ donnés, notons $\clUp$ l'adhérence de $\Up$ dans $\Delta_{\leq\ell}\Xg^{m}$, puis  
$$\displaywidth\dimen222\delta \Up:=\clUp\mmoins \Up\text{\quad et\quad }\Kp:=K\cap\clUp\,.\postdisplaypenalty10000$$
Les triplets d'espaces $(\Up\dans\clUp\cont\delta \Up)$ et $(\Fg_{\ell}(\Xg)\dans\Xg^{\ell}\cont\Delta_{\leq\ellmo}^{\ell})$ sont  homéomorphes et l'annulation de $(\diamond)$ résultera du lemme suivant qui généralise quelque peu le théorème \ref{theo-scindage}-(\ref{theo-scindage-b}).

\endgroup

\def\KK{J}
\subproclaim{Lemme. Soit $\Xg$ un espace $i$-acyclique. Pour tout compact $\KK \dans\Xg^{\ell}$, le morphisme de restriction 
$$\halfdisplayskips
\Hc(\Xg^{\ell})\to\Hc((\Delta_{\leq\ellmo}\Xg^{\ell})\cup \KK )\eqno(\ddagger\ddagger_{\KK })$$ est nul. En particulier, le prolongement par zéro $\Hc(\Fg_{\ell}(\Xg)\mmoins \KK )\to\Hc(\Xg^{\ell})$ est surjectif.

Preuve du lemme. On procède comme dans la démonstration de \ref{theo-scindage}. On considère le diagramme commutatif
$$\preskip0.5ex\postskip0.5ex
\xymatrix@R=3mm{
f_2^{-1}(\Fg_{\ellmo}(\Xg))\ar@{^(->}[r]\ar[rd]|{f'_{2}}&\big((\Delta_{\leq\ellmo}\Xg^{\ell})\cup \KK \big)\ar[rd]|{f_2}\ar@{^(->}[r]|(0.6){\,f\,}
&\Xg\times\Xg^{\ellmo}\mrlap{{}=\Xg^{\ell}}\ar[d]^(0.4){\, p_2}\\
&\Fg_{\ellmo}(\Xg)\ar@{^(->}[r]|{\;j\;}&\Xg^{\ellmo}
}$$
où $f$ est l'inclusion (fermée) et $j$ est l'inclusion (ouverte). La restriction $f'_{2}$ de $f_2:=p_2\circ f$ est propre puisque si $L\dans\Fg_{\ellmo}(\Xg)$ est compact, on a
$$f_{2}^{-1}(L)=\big(p_{2}^{-1}(L)\cap \Delta_{\ellmo}\Xg^{\ell}\big)\cup \big(p_{2}^{-1}(L)\cap \KK \big)
$$ 
où $p_{2}^{-1}(L)\cap \Delta_{\ellmo}\Xg^{\ell}$ est compact puisque la restriction de $p_{2}$
à $\Delta_{\ellmo}\Xg^{\ell}$ est un revêtement fini au-dessus de $\Fg_{\ellmo}(\Xg)$
{\rm(\cf\loccit)}. L'annulation de $(\ddagger\ddagger_{\KK })$ résulte alors d'appliquer  \ref{caracterisation}-(\ref{caracterisation-clef}), exactement comme dans \ref{theo-scindage}.\hfill\raise-0.5pt\hbox{\rotatebox{90}{$\boxminus$}}}

\medskip
Dans la situation présente, ce lemme (avec $\KK:=\Kp$) nous dit que le morphisme de prolongement par zéro
$$\halfdisplayskips
\Hc(\Up\mmoins K)\onto\Hc(\clUp)$$
est surjectif. Une classe de cohomologie $[\alpha]\in\Hc(\Up)$ est donc représentée par une cocycle d'Alexander-Spanier $\alpha\in\AS(\Up)$ dont le support $|\alpha|$ est compact dans $\Up\mmoins K$. Comme cette partie est ouverte dans $\Delta_{\leq\ell}\Xg^{m}$, le cocycle $\alpha$ représente aussi l'image $j_{!}[\alpha]$ de $[\alpha]$ dans $\Hc(\Delta_{\leq\ell}\Xg^{m})$. Or, on a $|\alpha|\cap K=\emptyset$ et donc $\rho_{K}(j_{!}[\alpha])=0$. L'annulation des morphismes $(\diamond)$ est ainsi établie et l'espace $\Delta_{\leq\ell}\Xg^{m}$ est bien $i$-acyclique.
\enddemo

\secskip=1.75ex\mou
\subsection{Variétés toriques affines}
\label{quo-tor}
\proclaim{On suppose le corps $k$ de caractéristique nulle.
\def\varlistskips{\leftmargin2.3em}\begin{enumerate}
\labelwidth1em\labelsep0.5ex
\item\leavevmode\label{quo-tor-a}Si $\Gg$ un groupe fini agissant sur un espace $i$-acyclique $\Xg$,
l'espace des orbites $\Xg/\Gg$ est $i$-acyclique.

\item\leavevmode\label{quo-tor-b}Une variété torique affine $U_{\sigma}$ associée à un cône simplicial $\sigma$ est $i$-acyclique.
\end{enumerate}}
\demo (\ref{quo-tor-a}) Résulte aussitôt des égalités $\Hc(\Xg/\Gg)=\Hc(\Xg)^{\Gg}$ et $\let\Hc\Hr\Hc(\Xg/\Gg)=\Hc(\Xg)^{\Gg}$ et de la commutativité du diagramme
$$\xymatrix@R=0.8cm@C=5mm{
\Hc(\Xg/\Gg)\ar[d]\ar@{=}[r]&\Hc(\Xg)^{\Gg}\ar@{^(->}[r]&\Hc(\Xg)\ar[d]^{0}\\
\Hr(\Xg/\Gg)\ar@{=}[r]&\Hr(\Xg)^{\Gg}\ar@{^(->}[r]&\Hr(\Xg)
}
$$

(\ref{quo-tor-b}) On a $U_{\sigma}=\AA_{n}(\CC)^{\Gg}\times(\CC^{*})^{m}$, où $\Gg$ est un groupe abélien fini agissant sur l'espace affine complexe $\AA_{n}(\CC)$
 (\cf\cite{ful}, \Spar2.1, p. 29 et \Spar2.2, p. 34).
\enddemo

\subsection{Groupes de Lie non-compacts}\label{Lie}
\proclaim{Tout groupe de Lie réel connexe non compact est $i$-acyclique.
}
\demo Comme un groupe de Lie $\Gg$ est toujours une variété différentielle orientable, il revient au même de montrer que $\Gg$ est \uacyclique\ (\ref{prop-acycliques}-(\ref{prop-acycliques-b}). Pour cela, on remarque que l'application
$$\mathalign{\Psi:&\Delta_\Gg\times\Gg&\to&\Gg\times\Gg\cr
&(x,x)g&\mapsto&(x,xg)}$$
est bien un homéomorphisme échangeant les plongements fermés 
$$\Delta_{\Gg}\times\set e/\dans\Delta_{\Gg}\times\Gg\,,\text{\quad et\quad }
\Delta_{\Gg}\dans\Gg\times\Gg$$
Le morphisme de restriction à la diagonale s'identifie alors par Künneth au morphisme 
$$\mathalign{
\Hc(\Delta_\Gg)\otimes\Hc(\Gg)&\to&\Hc(\Delta_\Gg)\otimes\Hc(\set e/)\cr
\omega\otimes\varpi&\mapsto&\omega\otimes\varpi\rest_{e}
}$$
clairement nul lorsque $\Gg$ n'est pas compact.
\enddemo

\subsection{Ouverts $i$-acycliques de $\IP_{n}(\RR)$}\label{proj-R}
\proclaim{Si $\car k\not=2$, tout ouvert de $\IP_{n}(\RR)$ est $i$-acyclique.}
\demo Résulte de ce que $\Hr^{i}(\IP_{n}(\RR);k)=0$ pour tout $i<n$.
\enddemo

\subsection{Ouverts $i$-acycliques de $\IP_{n}(\CC)$}\label{complementaires}
\proclaim{Soit $\Ug$ un ouvert de $\IP_{n}(\CC)$ dont le complémentaire contient
une hypersurface complexe $Z$.
\def\varlistskips{\leftmargin2.3em}\begin{enumerate}
\labelwidth1em\labelsep0.5ex
\item\leavevmode\label{complementaires-a}Si $\car(k)=0$, l'espace $\Ug$ est $i$-acyclique.

\item\leavevmode\label{complementaires-b}Si $\car(k)=p>0$, l'espace $\Ug$ est $i$-acyclique si $n=1\mod2$, ou bien si $n=0\mod2$ et si les degrés des composantes irréductibles de $Z$ sont premiers à $p$.
\end{enumerate}
}

\demo Par \ref{prop-acycliques}-(\ref{prop-acycliques-c}), l'assertion résulte aussitôt du cas où l'ouvert est $U:=\IP_{n}(\CC)\mmoins Z$. L'ouvert $U$ est alors une variété complexe \emph{affine} lisse de dimension $n$. D'après un théorème de Hamm (\cf \cite{dim} th. (6.8) p. 26), il a le type d'homotopie d'un CW-complexe fini de dimension réelle $n$. Il s'ensuit que les groupes d'homologie $\Hr_{i}(\Xg;A)$ sont nuls, pour $i>n$ et tout anneau~A, et de même en cohomologie. Comme $U$ est non singulière et orientable, on dispose de la dualité de Poincaré-Lefschetz $H_{i}(U;A)\simeq\Hc^{2n-i}(U;A)$, et donc $\Hc^{i}(U;A)=0$ pour tout $i<n$. Par conséquent, on a à priori
$$\big(\epsilon_{U,i}:\Hc^{i}(U;A)\to\Hr^{i}(U;A)\big)=0\,,\quad\forall i\not=n\,,$$
et la preuve de la proposition se réduit à montrer l'annulation de $\epsilon_{U,n}$.

Notons $j:U\hook\IP_{n}$ l'inclusion ensembliste.

Le morphisme $\epsilon_{U,n}$ se factorise à travers $\Hr^{n}(\IP_{n})$ suivant les morphismes
$$\mathalign{
\Hc^{n}(U;A)&\too^{j_{!}}&\Hr^{n}(\IP_{n};A)&\too^{j^{*}}&\Hr^{n}(U;A)\\
\mrlap{\let\biguparrow\uparrow
\sqrupc{\,\epsilon_{U,n}\,}{5.3cm}}}$$
où $j_{!}$ désigne le prolongement par zéro et $j^{*}$ la restriction. On a aussitôt deux cas à considérer suivant la parité de $n$.

\noindent{\bold $n=1\mod 2$. }On a $\Hr^{n}(\IP_{n};A)=0$ et donc $\epsilon_{U,n}=0$. La proposition est alors prouvée et ce, quel que soit l'anneau de coefficients $A$.

\noindent{\bold $n=0\mod 2$. }On a $\Hr(\IP_{n};A)=A$. Dans ce cas, on a intérêt à reprendre notre convention et considérer que l'anneau de coefficients $A$ est un corps, auquel cas on a l'équivalence
$$(\epsilon_{U,n}=0)\Longleftrightarrow \left(\!\vcenter{\parindent0pt\leftskip0pt plus1fill\rightskip\leftskip
\hsize4.8cm le morphisme de restriction\\ $\rho_{n}:\Hr^{n}(\IP_{n})\to\Hr^{n}(Z)$\\ est injectif}\!\right)\eqno(\ddagger)$$

En effet, comme $\Hr^{n}(\IP_{n})$ est une droite vectorielle, l'annulation de $\epsilon_{U,n}$ équivaut au fait que soit $j_!$, soit $j^{*}$, est nul. Or, la dualité de Poincaré échange ces morphismes et il sont tous les deux nuls ou non nuls. Ainsi, l'annulation de $\epsilon_{U,n}$ équivaux à l'annulation de $j_{!}$ et donc à l'injectivité de $\rho_{n}$, par la suite exacte longue de cohomologies.
$$\cdots\too\Hc^{n}(U)\too^{j_{!}}\Hr^{n}(\IP_{n})\too^{\rho_{n}}\Hr^{n}(Z)\too\cdots
$$

Notons 
$L$ l'opérateur sur $\Hr(\IP_{n})$ et $\Hr(Z)$ de multiplication par la classe fondamentale $\omega\in\Hr^{2}(\IP_{n})$ de $\IP_{1}\dans\IP_{n}$, et considérons le diagramme commutatif suivant
où la deuxième ligne est la suite exacte de cohomologies
$$
\vcenter{\hbox{$\xymatrix@C1cm@R8mm{
&\Hr^{n}(\IP_{n})\mput{0.65cm}{-2.5mm}{\bigoplus}
\ar[d]_{L^{(n-2)/2}}^{\simeq}\ar[r]^{\rho_{n}}&\Hr^{n}(Z)\ar[d]^{L^{(n-2)/2}}\\
\cdots\Hc^{2n-2}(U)\ar[r]^{j_!}&\Hr^{2n-2}(\IP_{n})\ar[r]^{\rho_{2n-2}}&\Hr^{2n-2}(Z)\ar[r]^{c}&
\Hc^{2n-1}(U)\cdots\\
}$}}
\eqno(\ddagger\ddagger)$$

Nous avons maintenant besoin de plus de renseignements sur le morphisme de liaison $c$. Pour cela, on commence par remarquer que l'on peut supposer $Z$ \emph{irréductible}. En effet, si $Z'$ est une composante irréductible de $Z$, l'ouvert $U=\IP_{n}\mmoins Z$ est contenu dans l'ouvert $U'=\IP_{n}\mmoins Z'$ et il est donc $i$-acyclique si $U'$ l'est (\ref{prop-acycliques}-\ref{prop-acycliques-c}). 

Supposons donc que $Z=V(f)$ est irréductible et, compte tenu des hypothèses, que $\deg f$ n'est pas multiple de la caractéristique du corps de coefficients. L'espace vectoriel $\Hr^{2n-2}(Z)$ est alors de dimension $1$ et comme
$$\Hc^{2n-1}(U;\ZZ)=H_1(U;\ZZ)\simeq\ZZ/(\deg f)\eqno(\star)$$
(\cf\cite{dim} chap. 4, prop. 1.3, p. 102), on a
$\Hc^{2n-1}(U)=\Hr_{1}(U)=0\,.$
Ces données reportées sur $(\ddagger\ddagger)$ assurent que $\rho_{n-2}$ est surjective, et donc injective aussi. L'injectivité de $\rho_{n}$ découle alors  de la commutativité du même diagramme.
\enddemo

\begin{rema}\parskip2pt\mou
L'hypothèse dans la proposition \ref{complementaires} concernant la caractéristique du corps de coefficients de la cohomologie est optimale. En effet, si $Z=V(f)\dans\IP_2(\CC)$ avec $f$ irréductible homogène de degré $p$, on a pour $U:=\IP_2(\CC)\mmoins V(f)$ d'après ($\star$):
$$\halfdisplayskips\Hc^{3}(U;\IF_{p})\simeq\Hr_{1}(U;\IF_p)\simeq\IF_p\,,$$
auquel cas, la suite exacte longue de cohomologie devient
$$\halfdisplayskips\Hr^{2}(\IP_{2};\IF_p)\too^{\rho_2}
(\Hr^{2}(V(f);\IF_p)=\IF_p)\too^{c}
(\Hr^{3}(U;\IF_p)=\IF_p)\to0\,,$$
le morphisme $c$ est bijectif, et alors nécessairement $\rho_{2}=0$.

On en déduit la surjectivité de $j_!:\Hc^{2}(U;\IF_p)\to\Hc^{2}(\IP_2;\IF_p)$ et, par dualité,
l'injectivité de $j^{*}:\Hr^{2}(\IP_2;\IF_p)\to\Hr^{2}(U;\IF_p)$. \`A partir de là, la non nullité de
$\epsilon_{U;2}=j^{*}\circ j_!$ découle de ce que $\Hr^{2}(\IP_2(\CC);\IF_p)=\IF_p$. L'ouvert $U$ n'est donc pas $i$-acyclique pour la cohomologie à  coefficients dans $\IF_p$.
\end{rema}

\begin{rema}Dans le cas où la proposition \ref{complementaires} intéresse la cohomologie à coefficients dans un corps de caractéristique nulle, le fait que $\Hr_{!}(U)$ soit concentrée en dimension moitié $n$, peut aussi être justifié en invoquant le théorème de comparaison de Grothendieck entre la cohomologie des formes différentielles holomorphes (dont le degré est à priori majoré par $n$) et la cohomologie du faisceau constant.

\addhabille-1\Habillage{\vrule\ \hbox{$\xymatrix@C1.5cm@R0.55cm{
\Hc^{n-2}(U)\ar[d]_{j_!}\\
\Hr^{n-2}(\IP_{n})\ar[d]_{\rho_{n-2}}\ar[r]^{L}_{\sim}&
\Hr^{n}(\IP_{n})\ar[d]^{\rho_{n}}\\
\Hr^{n-2}(Z)\ar[r]^{L}&\Hr^{n}(Z)
}$}}{1}{-8pt}Dans le même ordre d'idées, si l'on se restreint au cas où l'hypersurface est supposée non-singulière, l'injectivité de $\rho_{n}$ (et donc la $i$-acyclicité de $U$) admet aussi une justification via le théorème de Lefschetz vache. En effet, notons $L$ l'opérateur de multiplication par la classe génératrice $\omega\in\Hr^{2}(\IP_{n}(\CC))$ et considérons, pour $n$ pair, le diagramme commutatif ci-contre.
\endHabillage

L'injectivité de $\rho_{n}:\Hr^{n}(\IP_{n})\to\Hr^{n}(Z)$ découle alors par une chasse au diagramme élémentaire de ce que:
\begin{itemize}\itemsep0pt
\item $L:\Hr^{n-2}(\IP_{n})\to\Hr^{n}(\IP_{n})$ est trivialement bijectif,

\item $L:\Hr^{(n-1)-1}(Z)\to\Hr^{(n-1)+1}(Z)$ est bijectif d'après le théorème vache de Lefschetz.

\item $\rho_{n-2}:\Hr^{n-2}(\IP_{n})\to\Hr^{n-2}(Z)$ est injectif puisque $\Hc^{n-2}(U)=0$ en raison de précisément de l'affinité de $U$. 
\end{itemize}
\end{rema}

\subsection{Ouverts non $i$-acycliques de $\IP_{n}(\CC)$}\label{complementaires2}
\proclaim{Un ouvert de $\IP_{n}(\CC)$ dont le complémentaire 
est contenu dans un fermé $\Fg$ algébrique complexe (ou réel) et tel que $\dim_{\RR}(\Fg)<n$, n'est pas $i$-acyclique.}
\demo Un tel ouvert contient le complémentaire $\Ug$ d'un fermé algébrique $\Yg$ de petite dimension. Il suffira donc, d'après \ref{prop-acycliques}-(\ref{prop-acycliques-c}), de ne considérer que ce cas.
Des suites longues de cohomologie, on retire les sous-suites
$$\mathalign{
\Hr^{n-1}(\Yg)&\to&\Hc^{n}(\Ug)&\to&\Hr^{n}(\IP_{n})&\to&\Hr^{n}(\Yg)\mrlap{=0}\\
\mllap{0=}\Hr^{n}_{\Yg}(\IP_{n})&\to&\Hr^{n}(\IP_{n})&\to&\Hr^{n}(\Ug)&\to&\Hr_{\Yg}^{n+1}(\IP_{n})
}
$$
où $\Hr^{n}(\Yg)=0$ puisque $\dim_{\RR}(\Yg)<n$, et $\Hr_{\Zg}^{n}(\IP_{n})=\Hr^{n}(\Yg)^{\vee}=
0$  (\footnote{Par le formalisme de la dualité de Poincaré-Grothendieck-Verdier, cela résulte de ce que 
$\ID^{\bullet} \IR c_{\Yg!}c_{\Yg}^{!}\fs\ZZ_{\IP}=
\IR c_{\Yg*}c_{\Yg}^{-1}\fs\ID_{\IP}^{\bullet}
=
\IR c_{\Yg*}c_{\Yg}^{-1}\fs\ZZ_{\IP}[2n]
$ et du fait que $\fs\ID_{\IP}^{\bullet}=\fs\ZZ_{\IP}[2n]$.
}).
\enddemo

\subsection{Courbes algébriques $i$-acycliques}
\label{courbes-iacycliques}
\proclaim{Une courbe algébrique complexe irréductible $\Cg$ qui est
$i$-acyclique est rationnellement lisse. Si de plus $\Cg$ est affine, elle est homéomorphe à un ouvert algébrique de $\CC$.}

\demo Pour chaque $x\in\Xg$, et tout voisinage ouvert $V\ni x$, 
le morphisme $\alpha_x:\Hr_{x}(\Cg)\to\Hr(\Cg)$ se factorise à travers $\Hc(\Cg)$ de sorte que si $\Cg$ est supposée $i$-acyclique, on a $\alpha_{x}=0$ d'où la suite exacte courte
$$0\to\Hr(\Cg)\to\Hr(\Cg\mmoins\set x/)\to\Hr_{x}(\Cg)[1]\to0\,.$$
Il s'ensuit que lorsque $\Cg$ est complexe irréductible, on a $\Hr^{1}_{x}(\Cg)=0$ $(\ddagger)$ puisque $\Cg$ et $\Cg\mmoins\set x/$ sont connexes. D'autre part, la singularité $\Sigma$ de $\Cg$ est une partie fermée de dimension $0$ donc discrète, et $x$ admet des voisinages ouverts coniques $V=\hat c(\LL(x,\Cg))$, où
$\LL(x,\Cg)$, le \expression{lien de $x$ dans $\Cg$}, est une réunion disjointe de cercles $\sqcup_{i=1}^{r} \cSS^{1}$.
On a donc d'après $(\ddagger)$, l'égalité
$$0=\Hr_{x}^{1}(\Cg)=\Hr_{x}^{1}(\hat c(\LL(x,\Cg)))=\cl\Hr{}^{0}(\LL(x,\Cg))=k^{r-1}\,,$$
où $\cl\Hr$ désigne la cohomologie réduite. Par conséquent, $r=1$ et $x$ est rationnellement lisse dans $\Cg$.

Supposons maintenant que $\Cg$ est en plus affine. Quitte à remplacer les voisinages coniques  $\hat c(\cSS^{1})$ par des disques $\ID^{2}$, nous pouvons supposer que $\Cg$ est une sous-variété fermée et différentiable de $\AA_{n}(\CC)$ à laquelle nous pouvons appliquer la théorie de Morse (\footnote{\Cf le théorème fondamental de structure dans \cite{nic}, thm. 2.7, p. 47.}). Soit donc $h:\Cg\to\RR_{+}$ une fonction de Morse et notons $\Cg_{<r}:=\set x\in\Cg\mid h(x)<r/$, c'est un espace $i$-acyclique car ouvert de $\Cg$ qui l'est par hypothèse. Montrons que $\Cg_{<r}$ est homéomorphe à un ouvert algébrique de $\CC$. 

Pour $r$ petit on a $\Cg_{r}\sim\CC$. Supposons ensuite que pour une certaine valeur critique $r$, l'espace $\Cg_{<r}$ est homéomorphe à $\CC$ privé d'un ensemble fini de disques fermés disjoints homéomorphes à $\ID^{2}$. Le bord $\Cg_{r}=\partial\Cg_{<r}$ est alors une réunion finie de cercles $\cSS^{1}$.

Par la théorie de Morse, $\Cg_{<r+\epsilon}$ est homéomorphe au récolement de  $\Cg_{<r}$ et d'une anse $A_{1,2}=\cl\ID{}^{1}\times\cl\ID{}^{1}$ de noyau de dimension $1$, ou bien d'une anse $A_{2,2}=\cl\ID{}^{2}\times\cl\ID{}^{0}$ de noyau de dimension $2$. Comme $\Cg_{<r+\epsilon}$ est un ouvert d'une variété algébrique complexe, il est orientable, ce qui exclu certains recollements.
$$\preskip0pt\includegraphics{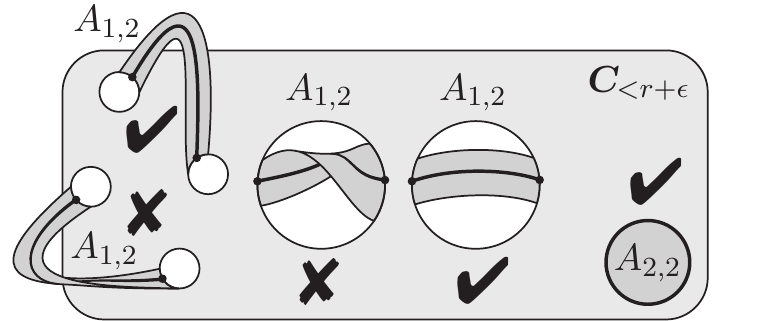}$$
Pour une anse $A_{1,2}$, il y a deux recollements orientables possibles qui résultent du fait que le bord du cœur de $A_{1,2}$, à savoir $\cl\ID{}^{1}=[0,1]$, possède deux composantes connexes $\set 0,1/$. Ces composantes peuvent être collées tantôt sur deux composantes connexes de $\partial\Cg_{<r}$ tantôt sur la même.
Dans le premier cas, $\Cg_{<r+\epsilon}$ serait homéomorphe à un tore $\Tg:=\Sg^{1}\times\Sg^{1}$ privé d'un ensemble fini de points $\Fg$, mais alors on aurait la suite exacte
$$\halfdisplayskips0\to\cl\Hr{}^{0}(\Fg)\to\Hc^{1}(\Tg\mmoins\Fg)\to\Hr^{1}(\Tg)\to0\,,$$
et $\Hr^{1}(\Tg)\hook\Hr^{1}(\Tg\mmoins\Fg)$ par dualité. Donc $\dim\Hi^{1}(\Cg_{r+\epsilon})=\dim\Hr^{1}(\Tg)=2\,$, ce qui contredit la $i$-acyclicité de $\Cg_{<r+\epsilon}$. Dans le deuxième cas, un trou de $\Cg_{r}$ est partiellement comblé par l'anse et $\Cg_{<r+\epsilon}$ est bien comme annoncé.

Pour une anse $A_{2,2}$. On recolle le disque $\ID^{2}$ par son bord $\cSS^{1}$ sur le bord de $\Cg_{<r}$. Par connexité, il n'y a qu'une seule manière de le faire et c'est en comblant l'un des trous de $\Cg_{<r}$. L'espace $\Cg_{<r+\epsilon}$ est donc, de nouveau, homéomorphe à $\CC$ privé d'un nombre fini de points.
\enddemo

\begin{rema}Dans la preuve de \ref{courbes-iacycliques}, le tore $\Tg:=\cSS^{1}\times\cSS^{1}$ privé d'un nombre fini non nul de points est une variété affine complexe non singulière $\Cg$ qui n'est pas $i$-acyclique. Les produits finis $\Pg:=\prod_{i}\Cg_i$ de tels espaces sont des exemples de variété affine complexe non singulière non $i$-acycliques, contrairement à tout produit de la forme $\CC\times\Pg$ (\ref{prop-acycliques}).
\end{rema}

\let\bq\bouquet
\subsection{Bouquet d'espaces $i$-acycliques}\label{bouquet}
\proclaim{Un bouquet d'espaces $i$-acycliques est $i$-acyclique. Et de même en remplaçant $i$-acyclique par (totalement) $\cup$-acyclique.}

\demo \def\pt{x=y}
Soient $\Xg$ et $\Yg$ deux espaces $i$-acycliques. Notons $\Xg\bq\Yg$ le bouquet qui identifie $x\in\Xg$ et $y\in\Yg$. 
$$\includegraphics{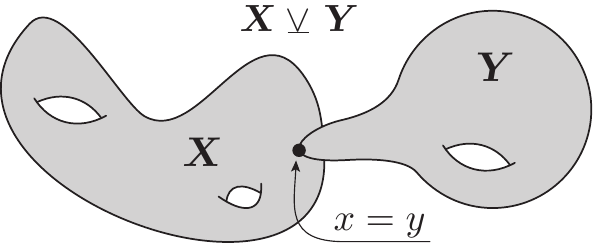}$$
On considère le morphisme suivant de suites exactes de Mayer-Vietoris associées à la décomposition en parties fermées $\Xg\bq\Yg=\Xg\cup\Yg$.
$$\halfdisplayskips\xymatrix@R=5mm
{
\ar[r]&\Hc(\pt)\aregal[d]\ar[r]&
\Hc(\Xg\bq\Yg)\ar[d]^{\epsilon_{\Xg\bq\Yg}}\ar[r]&\Hc(\Xg)\oplus\Hc(\Yg)\ar@<-8mm>[d]_{\epsilon_{\Xg}}\ar@<8mm>[d]^{\epsilon_{\Yg}}\ar[r]&\\
\ar[r]&\Hr(\pt)\ar[r]&\Hr(\Xg\bq\Yg)\ar[r]^(.45){\rho}&\Hr(\Xg)\oplus\Hr(\Yg)\ar[r]&\\
}
$$ 
où $\rho:\Hr^{+}(\Xg\bq\Yg)\to\Hr^{+}(\Xg)\oplus\Hr(\Yg)$ est clairement injectif. La nullité de $\epsilon_{\Xg\bq\Yg}$ résulte alors immédiatement de celles de $\epsilon_{\Xg}$ et $\epsilon_{\Yg}$.

\nobreak Le même raisonnement prouve que $\Xg\bq\Yg$ est (totalement) $\cup$-acyclique si les espaces $\Xg$ et $\Yg$ le sont.
\enddemo

\subsection{Sommes amalgamées d'espaces $i$-acycliques}
\proclaim{La somme amalgamée de variétés $i$-acycliques est $i$-acyclique  si et seulement si au moins l'une des variétés  est orientable.}

\demo\def\QQ{k}\relax
Soient $\Xg$ et $\Yg$ deux variétés topologiques $i$-acycliques de dimension $n$. Notons $\Xg'$ et $\Yg'$ des complémentaires d'un point dans $\Xg$ et $\Yg$ respectivement, ce sont des espaces $i$-acycliques d'après \ref{prop-acycliques}-(\ref{prop-acycliques-c}). La somme amalgamée $\Xg\amalgam \Yg$ est la variété topologique obtenue en recollant homéomorphiquement $\Xg'$ et $\Yg'$ le long du cylindre ouvert bordant $\cSS^{n-1}\times\RR$ noté $\Vg $ bordant respectivement $\Xg'$ et $\Yg'$.
$$\includegraphics{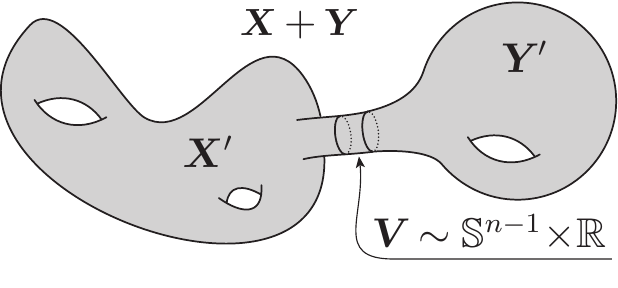}$$
On considère alors les suites exactes longues de Mayer-Vietoris  pour le recouvrement ouvert $\Xg\amalgam \Yg=\Xg'\cup\Yg'$ 
$$
\mvcenter{\xymatrix{
\ar[r]&
\Hc^{i}(\Vg )\ar[r]^(0.35){\alpha_i}\ar[d]|{\ttt85\epsilon_{\Vg,i}}&
\Hc^{i}(\Xg')\oplus\Hc^{i}(\Yg')\ar[r]^(0.58){\beta_i}\ar[d]|{\ttt85\epsilon_{\Xg',i}\oplus
\epsilon_{\Yg',i}}&
\Hc^{i}(\Xg\amalgam \Yg)\ar[d]|{\ttt85\epsilon_{\Xg\amalgam \Yg,i}}\ar[r]&\\
\ar@{<-}[r]&
\Hr^{i}(\Vg )\ar@{<-}[r]^(0.36){\delta_i}&
\Hr^{i}(\Xg')\oplus\Hr^{i}(\Yg')\ar@{<-}[r]^(0.59){\gamma_i}&
\Hr^{i}(\Xg\amalgam \Yg)\ar@{<-}[r]&
}}
\eqno(\star)$$
où 
$$\begin{cases}
\Hc^{i}(\Vg )=\QQ(1)\otimes\big(\QQ(0)\oplus\QQ(n-1)\big)=\QQ(1)\oplus\QQ(n)\,,\\
\Hr^{i}(\Vg )=\QQ(0)\oplus\QQ(n-1)\,.\end{cases}$$

\noindent{\bold Cas $i\not=n-1$. }Le morphisme $\beta_i$ est surjectif et $\gamma_i$ est injectif, soit parce que $i<n$ auquel cas $\Hr^{i-1}(\Vg)=0$, soit parce que $i=n$ auquel cas les cohomologies $\Hr^{n}(\_)$ sont nulles. On considère alors le diagramme
$$\mathalign{
\Hc^{i}(\Xg')\oplus\Hc^{i}(\Yg')&\hfonto{\beta_i}{}{0.7cm}&\Hc^{i}(\Xg\amalgam \Yg)&\too&0\hfill\cr
&&\vfld{}{\epsilon_{\Xg\amalgam \Yg,i}}{0.5cm}\cr
\hfill\Hr^{i-1}(\Vg )&\hf{}{}{0.7cm}&\Hr^{i}(\Xg\amalgam \Yg)&\hfhook{\gamma_i}{}{0.7cm}&\Hr^{i}(\Xg')\oplus\Hr^{i}(\Yg')
}\eqno(\ast)$$
où les lignes sont exactes. On y voit que 
$$\gamma_i\circ\epsilon_{\Xg\amalgam \Yg,i}\circ\beta_i=\epsilon_{\Xg',i}\oplus\epsilon_{\Yg',i}=0\,,$$
et nous concluons, indépendamment de l'orientabilité de $\Xg$ et $\Yg$, que 
$$\epsilon_{\Xg\amalgam \Yg,i}=0\,,\quad\forall i\not=n-1\,.$$

\noindent{\bold Cas $i=n-1$. }La suite de Mayer-Vietoris en degré $n$ est
$$\too\big(\Hc^{n}(\Vg)=\QQ\big)\too^{\alpha_{n}}\Hc^{n}(\Xg')\oplus\Hc^{n}(\Yg')\too\Hc^{n}(\Xg\amalgam \Yg)\to\0\,.\eqno(\ddagger)$$
montre que $\alpha_{n}$ est injective pour peu que $\Xg$ ou $\Yg$ soit orientable. Dans ces cas, $\beta_{n-1}$ sera surjective et le diagramme $(\ast)$ pour $i=n-1$ fournit toujours l'annulation $\epsilon_{\Xg\amalgam \Yg,n-1}=0$ puisque $\gamma_{n-1}$ est injective. Ceci termine la preuve de l'un des sens de la proposition.

Réciproquement, ni $\Xg$ ni $\Yg$ ne sont orientables, ce qui équivaut à dire que ni $\Xg'$ ni $\Yg'$ ne le sont. Nous allons nous intéresser à des suites longues de cohomologies associées à la décomposition ouvert/fermé
$$\Xg'=\Vg\sqcup (\Xg'\mmoins\Vg)$$
où $\Xg'\mmoins\Vg$ est clairement est une variété à bord $\cSS^{n-1}$ et d'intérieur homéomorphe à $\Xg'$, raison pour laquelle on la notera $\cl{\Xg'}$. Enfin, l'adhérence de $\Vg$ dans $\Xg'$ qui sera notée $\cl\Vg$, est une variété à bord dont le bord $\cSS^{n-1}$ est aussi le bord de $\cl{\Xg'}$.
$$\includegraphics{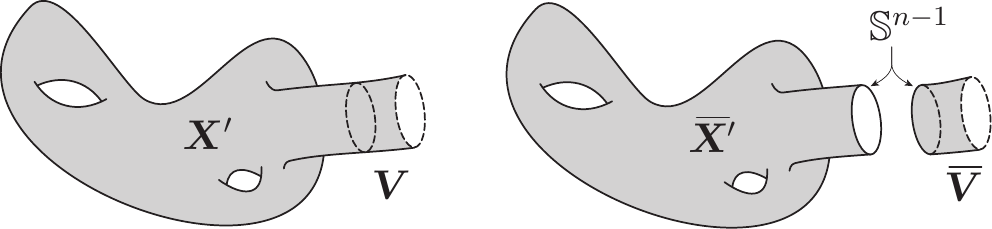}$$

On a alors le diagramme de suites exactes longues :
$$\xymatrix@C0.7cm@R=7mm{
\Hc^{n-1}(\Xg')\ar[r]^{\rm rest}\ar[d]_{\rm rest}&
\Hc^{n-1}(\cl{\Xg'})\xylbl[rd]{(I)}\ar[r]^(0.44){c_{\Xg}}\ar[d]_{\rm rest}&
(\Hc^{n}(\Vg)=k)\ar[r]\ar[d]_{\id}&
(\Hc^{n}(\Xg')=0)\ar@<-3ex>[d]_{\rm rest}\\
(\Hc^{n-1}(\cl{\Vg})=0)\ar[r]&
\Hc^{n-1}(\cSS^{n-1})\ar[r]^(0.44){c_{\Xg}}_(0.44){\simeq}&
(\Hc^{n}(\Vg)=k)\ar[r]&
(\Hc^{n}(\cl\Vg)=0)\\
}$$
On y constate que la classe fondamentale de $\Vg$ se relève bien en un classe $[\Vg_{\Xg}]\in\Hc^{n-1}(\cl{\Xg'})$ et que cette classe se restreint sur la classe fondamentale $[\cSS]$ de $\cSS^{n-1}$, d'où le carré:
$$\preskip1.5ex\xymatrix@R7mm@C=1.5cm{
\mllap{\Hc^{n-1}(\cl{\Xg'})\ni}[\Vg_{\Xg}]\xylbl[rd]{(I)}\ar@{|->}[d]_{\rm rest}\ar@{|->}[r]^(0.53){c_{\Xg}}&[\Vg]\ar@{=}[d]\\
\mllap{\Hr^{n-1}(\cSS^{n-1})\ni}[\cSS]\ar@{|->}[r]^(0.53){c_{\Xg}}&[\Vg]\\
}$$
Maintenant, en remarquant que $\Hc(\cl{\Xg'})\too^{\rm rest}\Hc(\cSS^{n-1})$ se factorise suivant 
$$\Hc(\cl{\Xg'})\too^{\epsilon_{\cl{\Xg'}}}\Hr(\cl{\Xg'})\too^{\rm rest}
\Hr(\cSS^{n-1})\,,$$ on déduit que 
$$\epsilon_{\cl{\Xg'}}([\Vg_{\Xg}])\not=0\,.\eqno(**)$$

Ces observations s'appliquent telles quelles également à $\Yg$ et conduisent aux mêmes conclusions où $\Yg$ vient remplacer $\Xg$.

Considérons à présent la décomposition ouvert/fermé
$$\Xg\amalgam \Yg=\Vg\sqcup (\Xg\amalgam \Yg\mmoins\Vg)=\Vg\sqcup(\cl{\Xg'}\sqcup\cl{\Yg'})$$
et le morphisme de suites exactes de cohomologie
$$
\xymatrix@C1.2cm{
\Hc^{n-1}(\Xg\amalgam \Yg)\mput{5mm}{-5mm}{\bf II}\ar[r]\ar[d]|{\ttt76\epsilon_{\Xg\amalgam \Yg}}&
\Hc^{n-1}(\cl{\Xg'})\oplus\Hc^{n-1}(\cl{\Yg'})\ar[d]|{\ttt76\epsilon_{\cl{\Xg'}}\oplus\epsilon_{\cl{\Yg'}}}\ar[r]^(.59){c_{\Xg} + c_{\Yg}}&
(\Hc^{n}(\Vg)=k)\ar[d]|{\ttt76\epsilon_{\Vg}}\\
\Hr^{n-1}(\Xg\amalgam \Yg)\ar[r]&\Hr^{n-1}(\cl{\Xg'})\oplus\Hr^{n-1}(\cl{\Yg'})\ar[r]^(.59){c_{\Xg} + c_{\Yg}}&
(\Hr^{n}(\Vg)=0)
}$$

La classe $([\Vg_{\Xg}],-[\Vg_{\Yg}])\in\Hc^{n-1}(\cl{\Xg'})\oplus\Hc^{n-1}(\cl{\Yg'})$ appartient clairement au noyau de $c_{\Xg}\oplus c_{\Yg}$ et survit au morphisme $\epsilon_{\cl{\Xg'}}\oplus\epsilon_{\cl{\Yg'}}$ d'après $(**)$. La commutativité de $({\bf II})$ permet alors de conclure que $\epsilon_{\Xg\amalgam \Yg}\not=0$ et donc que l'espace $\Xg\amalgam \Yg$ n'est pas $i$-acyclique.
\enddemo

\subsection{La bouteille de Klein épointée}
Dans\label{U-non-I} l'implication (\uacyclique)$\,\Rightarrow\,$($i$-acyclique) de l'assertion \ref{prop-acycliques}-(\ref{prop-acycliques-b}), l'hypothèse d'orientabilité est indispensable. On donne ici un contre-exemple lorsque cette hypothèse fait défaut.\begingroup\def\x{x}\def\y{y}

\proclaim{ La bouteille de Klein épointée vérifie le théorème de scindage \ref{theo-scindage}. Elle n'est ni $i$-acyclique ni totalement $\cup$-acyclique, mais elle est bien $\cup$-acyclique.}

\smallskip
Notons $\Tg$ le tore $\RR^2/\ZZ^{2}$.
L'anneau de cohomologie $\Hr(\Tg)$ est engendré par les $1$-cocycles $d\x $ et $d\y $.
Si nous notons $C_{\x }:={\RR/\ZZ}\times\set0/$ et $C_{\y }:=\set0/\times(\RR/\ZZ)$, on peut voir que 
$$\halfdisplayskips
\Hr^{1}_{C_{\x }}(\Tg)=\langle d\y \rangle_{k}\dans\Hr^{1}(\Tg)\,,\quad
\Hr^{1}_{C_{\y }}(\Tg)=\langle d\x \rangle_{k}\dans\Hr^{1}(\Tg)\,.
$$

\addhabille1
\Habillage{\includegraphics{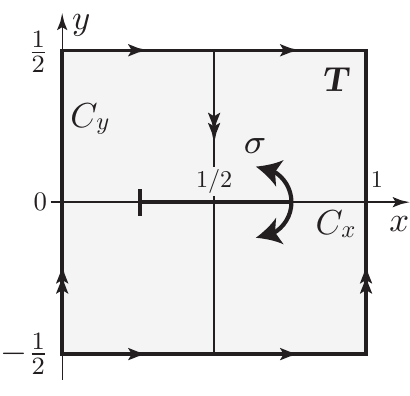}}{0}{-5pt}L'involution $\sigma:(\x ,\y )\mapsto(\x +1/2,-\y)$, induit sur $\Hr(\Tg)$ l'action
$\sigma: d\x \mapsto d\x $ et $\sigma: d\y \mapsto -d\y $ de sorte que l'on a
$$\Hr(\Tg)^{\langle\sigma\rangle}\simeq k(0)\oplus \langle d\x \rangle_{k}(1)\,.$$
Comme la bouteille de Klein est la variété quotient $\Kg:=\Tg/\langle\sigma\rangle$, si
$\nu:\Tg\to \Kg$ 
est la surjection canonique, le morphisme $\nu^{*}:\Hr(\Kg)\to\Hr(\Tg) $ identifie $\Hr(\Kg)$ à $\Hr(\Tg)^{\langle\sigma\rangle}\simeq k(0)\oplus k(1)$. 
\endHabillage

\smallskip
\noindent On pose maintenant $\Xg:=\Kg\mmoins\pt$.

\soustitreline{$\cup$-acyclicité de $\Xg$.}
 L'examen des suites exactes longues
$$\xymatrixc{@C=5mm@R=1mm}{
\ar[r]&\Hr^{i}_{\pt}(\Xg)=k(2)\ar[r]&\Hr^{i}(\Kg)\ar[r]^{\rho_i}&\Hr^{i}(\Xg)\ar[r]&\\
\ar[r]&\Hc^{i}(\Xg)\ar[r]^{\iota_i}&\Hc^{i}(\Kg)\ar[r]&\Hc^{i}(\pt)=k(0)\ar[r]&\\
}\eqno(\diamond)$$
montre que l'on a
$$\mathalign{
\Hr(\Xg)&=&k(0)\oplus 
\big(\langle d\x \rangle_{k}\oplus \langle d\x  \wedge d\y \rangle_{k}\big)(1)=k(0)\oplus k(1)^{2}\,,\\\noalign{\kern2pt}
\Hc(\Xg)&=&\Hr^{+}(\Kg)=\langle d\x \rangle_{k}(1)=k(1)\,.\hfill}$$
et $\Hc(\Xg)\wedge\Hc(\Xg)=0$. L'espace $\Xg$ est donc bien $\cup$-acyclique.

\soustitreline{Non $i$-acyclicité de $\Xg$.}
Comme le cycle $C_\y $ est plongé (via $\nu$) dans $\Xg$, l'espace $\Hc^{1}(\Xg)$ est engendré par $d\x \in\Hr_{C_\y }^{1}(\Xg)\dans\Hc^{1}(\Xg)$ et le morphisme 
$\iota_1:\Hc^{1}(\Xg)\to\Hc^{1}(\Kg)$ est bijectif. D'autre part, il est clair par la première suite  dans $(\diamond)$ que le morphisme $\rho_1:\Hr^{1}(\Kg)\hook\Hr^{1}(\Xg)$ est injectif. On a donc
$$\mathalign{
\hfill\Hc^{1}(\Xg)&\hf{\iota_1}{\sim}{0.7cm}&\Hc^{1}(\Kg)\cr
&&\vegal{0.5cm}\cr
0=\Hr_{\pt}^{1}(\Kg)&\too^{0}&\Hr^{1}(\Kg)&\hfhook{\rho_1}{}{0.7cm}\Hr^{1}(\Xg)\hfonto{}{}{0.7cm}\Hr_{\pt}^{2}(\Kg)=k(2)
}$$
et
$\epsilon_{\Xg,1}=\rho_1\circ\iota_1$ est non nul.
L'espace $\Xg$ n'est donc pas $i$-acyclique.

\soustitreline{Non totale $\cup$-acyclicité de $\Xg$.}
Soit maintenant $\Ug$ un ouvert connexe de $\Xg$. Comme $\Hc(\Xg)=k(1)$, le cup-produit 
$$\cup:\Hc(\Xg)\times\Hc(\Ug)\to\Hc(\Ug)\eqno(\ast)$$ 
est automatiquement nul si $\Ug$ n'est pas orientable car alors $\Hc^{2}(\Ug)=0$. Lorsque, par contre, $\Ug$ est orientable, la dualité de Poincaré s'applique et le cap produit $(\ast)$ est nul, si et seulement si, la restriction $\Hc^{1}(\Xg)\to\Hr(\Ug)$ est nulle, donc si et seulement si $d\x \rest_{\Ug}=0$ dans $\Hr(\Ug)$.

On a deux cas possibles. 
\begin{itemize}\setbox11=\hbox{\includegraphics{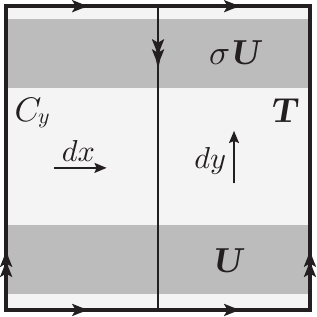}}
\dimen0=\wd11\advance\dimen0 by1em
\dimen1=\hsize \advance\dimen1 by-\dimen0
\parshape 6 
0pt \hsize
0pt \hsize
0pt \hsize
0pt \hsize
0pt \hsize
0pt \dimen1
\item[--]
Cas $C_\y \dans\Ug$. La restriction du revêtement $\nu:\Tg\onto\Kg$ à $\nu^{-1}(\Ug)$ est  triviale à deux nappes $\Ug_1$ et $\Ug_{2}=\sigma\Ug_1$, car $\Ug$ orientable. On peut supposer $C_\y \dans\Ug_1$ et alors $\sigma C_y\cap\Ug_{1}=\emptyset$. On a donc $\Ug\simeq\Ug_1\dans\Tg\mmoins C_\y $, et $\Ug$ est $i$-acyclique puisque ouvert de $\Tg\mmoins C_\y \sim\cSS^{1}\times\RR$, espace $i$-acyclique. Maintenant, comme $C_\y \dans\Ug_1$, l'image de $\Hr^{1}_{C_\y }(\Xg)=\Hc^{1}(\Xg)\to\Hr^{1}(\Ug)$ est l'image de $\Hr^{1}_{C_\y }(\Ug)\to\Hr^{1}(\Ug)$, morphisme qui se factorise naturellement à travers $\epsilon_{\Ug}:\Hc(\Ug)\to\Hr(\Ug)$, nul puisque $\Ug$ $i$-acyclique. Par conséquent, si $C_\y \dans\Ug$ le cup produit $(\ast)$ est nul.

\item[--]\tolerance700
\parshape 4 
0pt \dimen1
0pt \dimen1
0pt \dimen1
0pt \hsize
Cas\vadjust{\vtop to0pt{\kern-2.4cm\hbox to\hsize{\hss \copy11}\vss}} $C_\y \not\dans\Ug$. Il existe des ouverts $\Ug\dans\Xg$ tels que $(\ast)$ n'est pas nul.
En effet, si $\Ug$ est l'ouvert indiqué par la figure ci-contre, on y voit que c'est un cylindre et que $d\x $ ne s'intègre pas en une fonction sur $\Ug$. Par conséquent, $d\x \not=0$ dans $\Hr^{1}(\Ug)$ et l'espace $\Xg$ n'est pas totalement $\cup$-acyclique. 
\end{itemize}

\smallskip\soustitreline{Le théorème de scindage \ref{theo-scindage} pour $\Xg$.}L'espace $\Xg$ vérifie les assertions (\ref{theo-scindage-a}) et (\ref{theo-scindage-b}) de ce théorème, et la raison vient de ce que $\Hc(\Xg)$ est concentrée en degré $1$. En effet, dans de tels cas, on s'aperçoit rapidement que $\Hc(\Delta_{?\ell}\Xg^{m})$ doit être concentré en degré $\ell$, ce que l'on démontre par induction sur $\ell$ et prouve, en passant, les assertions en question.

Plus précisément, dans le cas des espaces $\Fm$, on s'intéresse au morphisme
$$\Hc(\Xg)\otimes\Hc(\Fm)\to
\Hc(\Delta_{\leq m}(\Xg\times\Fm ))\sim\Hc(\Fm)^{m}\eqno(\ddagger)$$
où le terme de gauche est concentré en degré $m+1$ puisque $\Hc(\Xg)=k(1)^{\#}$ et que $\Hc(\Fm)=k(m)^{\#}$, par hypothèse inductive. Le morphisme $(\ddagger)$ est alors nul puisque $\Hc(\Delta_{\leq m}(\Xg\times\Fm ))k(m)^{\#}$, également par hypothèse inductive. On a donc la suite exacte courte
$$0\to
\Hc(\Fm )[-1]^{m}\to
\Hc(\Fmm )\to
\Hc(\Xg\times\Fm)\to0\,,
$$
et elle montre que $\Hc(\Fmm )$ est concentré en degré $m+1$, et l'étape inductive peut être itérée.

Dans le cas des espaces $\Delta_{\leq\ell}\Xg^{m}$, on s'intéresse au morphisme
$$\Hc(\Delta_{\leq\ell}\Xg^{m})\to\Hc(\Delta_{\leq\ellmo}\Xg^{m})\eqno(\ddagger\ddagger)$$
dont la preuve de \ref{theo-scindage} montre que son annulation résulte de l'annulation de 
$$\Hc(\Xg)\otimes\Hc(\Delta_{\leq\ellmo}^{\mmo})\to\Hc(\Delta_{\leq\ellmo}^{m})\,.$$
Cela résulte pour les mêmes raisons de degré que précédemment, sous l'hypothèse d'induction que $\Delta_{\leq\ellmo}^{??}$ est concentré en degré $\ellmo$. \`A partir de là, on a la suite exacte courte
$$0\to
\Hc(\Delta_{\leq \ellmo}\Xg^{m})[-1]\to
\Hc(\Delta_{\ell}\Xg^{m})\to
\Hc(\Delta_{\leq \ell}\Xg^{m})\to
0\,,
$$
où le terme central est concentré en degré $\ell$ car somme directe d'espaces gradués isomorphes à $\Hc(\Fg_{\ell})$ (\cf\ref{connexes}). On conclut que $\Hc(\Delta_{\leq \ell}\Xg^{m})$ est concentré en degré $\ell$, et l'étape inductive peut être itérée.

\endgroup

\subsection{Revêtements non $i$-acycliques à base $i$-acyclique}
\proclaim{Un revêtement fini de base $i$-acyclique peut ne pas être $i$-acyclique.}

Une  manière élémentaire de produire des contre-exemples consiste à prendre une variété orientable $\Yg$ non $i$-acyclique, donc telle que $\Hc(\Yg)\wedge\Hc(\Yg)\not=0$, et à faire agir librement sur $\Yg$ un groupe fini $\Wg$ qui conserve l'orientation et tel que $\Hc(\Yg)^{\Wg}\wedge\Hc(\Yg)^{\Wg}=0$. 

\Habillage{\includegraphics{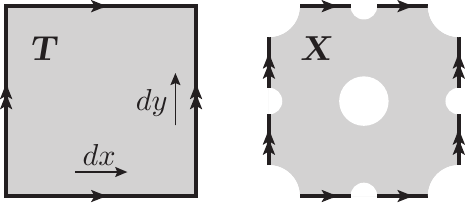}}{0}{-5pt}Le tore $\Tg$ de dimension $2$ est le domaine carré ci-contre sur lequel on a identifié les côtés opposés. Cet espace privé de quatre points, noté $\Xg$, est représenté par  la figure à droite. La suite exacte longue
\endHabillage
\vskip-1.5em$$\mathrigid2mu 
0\to\Hc^{0}(\Tg)\to\Hc^{0}({4})\to\Hc^{1}(\Xg)\to\Hc^{1}(\Tg)\to0\to
\Hc^{2}(\Xg)\to\Hc^{2}(\Tg)\to0\,,
\eqno(\dagger)$$
montre que les $1$-cocycles $dx$ et $dy$ dans $\Tg$ sont représentés par des $1$-cocycles à support compact dans $\Xg$. La classe fondamentale de $\Xg$ se retrouve alors représentée comme le produit de ces deux cocycles. Il s'ensuit que $\Xg$ n'est pas \uacyclique\ (ni $i$-acyclique). La suite $(\dagger)$ révèle aussi l'apparition de trois nouveaux $1$-cocycles dans $\Xg$ provenant de $\Hc^{0}({4})$, nous avons donc
$$\begin{lcases}
\Hc^{0}(\Xg)=\0\cr
\Hc^{1}(\Xg)=\QQ^{5}\cr
\Hc^{2}(\Xg)=\QQ\,.
\end{lcases}$$
Voici une représentation de $1$-cycles de $\Xg$ générant $\Hc^{1}(\Xg)(=H_{1}(\Xg))$:
$$\vcenter{\hbox{\includegraphics{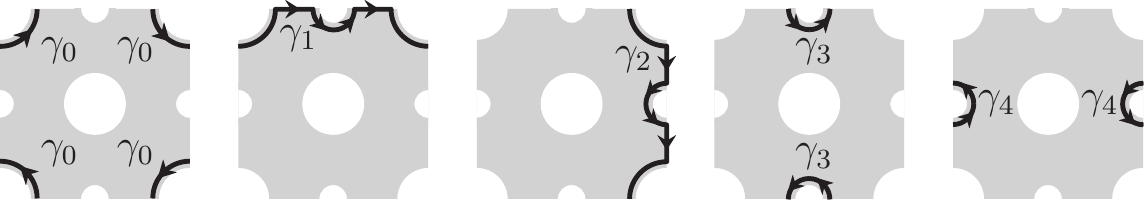}}}\eqno(\ddagger)$$

\medskip
Notons maintenant $\theta:\Tg\to\Tg$ l'isomorphisme d'ordre $4$ qui se voit sur  $(\ast)$ par une rotation horaire d'angle $\pi/4$. Il est clair que $\theta$ agit librement sur $\Xg$ et préserve son orientation. 
L'espace quotient $\Bg:=\Xg/\langle\theta\rangle$ est par conséquent une variété différentiable orientable et $$\Hc(\Bg)=\Hc(\Xg)^{\langle\theta\rangle}\,.$$
Or, l'action de $\theta$ sur les $1$-cycles
se lit facilement sur $(\ddagger)$. On y voit que
$$\theta(\gamma_0)=\gamma_0\,,\quad\begin{cases}
\theta(\gamma_1)=\gamma_2\cr
\theta^{2}(\gamma_1)=\gamma_0-\gamma_1+\gamma_3\cr
\theta^{3}(\gamma_1)=\gamma_0-\gamma_2+\gamma_4\cr
\end{cases}
\text{\quad et\quad}
\begin{cases}\theta(\gamma_3)=\gamma_4\cr
\theta(\gamma_4)=\gamma_3\,,
\end{cases}
$$
de sorte que
$\Hc^{0}(\Bg)=0$, $\Hc^{1}(\Bg)=[\gamma_0]\QQ\oplus[\gamma_3+\gamma_4]\QQ$, et $\Hc^{2}(\Bg)=\QQ$,
où $[\gamma_0]\wedge[\gamma_3+\gamma_4]=0$ puisque l'intersection des cycles sous-jacents est vide. On  conclut que
$\Hc(\Bg)\wedge\Hc(\Bg)=\0$
et $\Bg$ est \uacyclique\ et donc $i$-acyclique. (En fait, on peut vérifier que $\Bg$ est un cylindre ouvert épointé $(\cSS^{1}\times\RR)\mmoins\point$.)

Le revêtement $\Xg\to\Bg$ est bien à base $i$-acyclique alors que $\Xg$ ne l'est pas.

\comment
\begin{prop}\label{Euler}Soit $\Xg$ une variété différentiable compacte et orientable. Notons $\Tg\Xg$ son fibré tangent. Une condition nécessaire pour que $\Tg\Xg$ soit $i$-acyclique  (ou \uacyclique) est que la caractéristique d'Euler-Poincaré $\chi(\Xg)$ soit nulle.
\end{prop}
\demo En effet, si $\omega$ désigne la classe de la section nulle de $\Tg\Xg$, notée $\Xg$ dans la suite, on a bien $\omega\in\Hc(\Tg\Xg)$. On sait que $\int_{\Xg}\omega\rest_\Xg=\chi(\Xg)$. Or, comme $\Xg$ est compacte, la restriction
$\Hc(\Tg\Xg)\too\Hc(\Xg)=\Hr(\Xg)$ se factorise bien en 
$\Hc(\Tg\Xg)\too^{\epsilon_{\Tg\Xg}}\Hr(\Tg\Xg)\too^{\rm rest}\Hr(\Xg)$ et $\omega\rest_{\Xg}=0$ lorsque $\Tg\Xg$ est $i$-acyclique.
\enddemo
\endcomment

\def\TIP{\TTT\IP}
\subsection{Fibrations non $i$-acycliques à fibre et base $i$-acycliques}
\proclaim{Une fibration à fibre et base $i$-acycliques peut ne pas être $i$-acyclique.}

\smallskip
Rappelons que la cohomologie de l'espace projectif complexe $\IP:=\IP_{n}(\CC)$ est isomorphe à l'anneau gradué $\QQ[\zg]/(\zg^{n+1})$ où $\deg(\zg)=2$.

Le fibré tautologique $\TIP$ au-dessus de $\IP$ est le sous-espace de $\IP\times\CC^{n+1}$ des couples $(\D,z)$ où $\D$ est une droite vectorielle de $\CC^{n+1}$ et $z\in\D$. L'application $\pi:\TIP\to\IP$, $(\D,z)\mapsto\D$, est une fibration vectorielle localement triviale de fibres isomorphes à $\CC$. Notons $\IP\dans\TIP$ sa section nulle. La classe de Thom $\Phi_{(\IP,\TIP)}\in\Hc^{2}(\TIP)$ de la section nulle se restreint à la classe génératrice de $\Hr(\IP)$, autrement dit,
$$\Phi_{(\IP,\TIP)}\rest_{\IP}=\zg\postskip0pt$$
En fait, on a 
$$\Hc(\TIP)=\QQ[\Phi]/(\Phi^{n+2})$$
et la restriction à la section nulle s'identifie à l'isomorphisme canonique
$$\QQ[\Phi]/(\Phi^{n+2})\to\QQ[\zg]/(\zg^{n+1})\,,\quad \Phi\mapsto\zg$$

En particulier, la multiplication par $\Phi_{\IP,\TTT}$ n'est pas une opération nulle dans $\Hc^{+}(\TIP)$ pour peu que $n\geq1$. 

\medskip
Maintenant, si nous restreignons le fibré tautologique au complémentaire de deux points $\IP':=\IP\mmoins2$, nous avons la suite exacte longue
$$\0\to\Hc^{2}(\TIP')\too\Hc^{2}(\TIP)=\Phi\cdot\CC\too^{\rho}\Hc^{2}(\CC)^{2}\too\Hc^{3}(\TIP')\to\0$$
où $\rho$ est injective. 
On en déduit que $\Hc^{i}(\TIP)=0$ pour $i\in\set0,1,2/$, puis que $\Hc^{3}(\TIP)=\QQ$, et enfin que 
$$\Hc^{i}(\TIP')=\Hc^{i}(\TIP)\,,\quad\text{pour tout $i\geq4$.}\eqno(\ast)$$
Par conséquent, si $n\geq4$ on a
$$\def\tt{\vrule height12pt depth6pt width0pt}
\begin{tabular}{|c|c|c|c|c|c|c|c|c|c|c|}
\hline
\tt $i$&0&1&2&3&4&5&6&7&8&\dots\\
\hline
\tt $\dim\Hc^{i}(\TIP')$&0&0&0&0&1&0&1&0&1&\dots\cr
\hline
\tt $\dim\Hc^{i}(\IP')$&0&1&1&0&1&0&1&0&1&\dots\cr
\hline
\end{tabular}$$
et si $\omega$ est non nulle de $\Hc^{4}(\TIP')$, on a $\omega\wedge\omega\not=0$ dans $\Hc^{8}(\TIP')$ compte tenu de $(\ast)$ et du fait qu'il en est ainsi pour $\Hc(\TIP)$ (pour $n\geq3$).

\smallskip
Considérons maintenant $\IP'$ pour $n=3$. Dans ce cas, on dispose sur $\IP_3(\CC)$ de l'involution $\theta$ sans point fixes :
$$\theta(\langle a,b,c,d\rangle)=\langle b,- a, d,- e\rangle$$
dont on sait que son action sur $\zg$ est $\theta\zg=-\zg$, de même par ailleurs que l'action sur le générateur de $\Hc^{1}(\IP')$. Il s'ensuit que l'on a
$$\def\tt{\vrule height12pt depth6pt width0pt}
\begin{tabular}{|c|c|c|c|c|c|c|c|c|c|}
\hline
\tt  $i$&0&1&2&3&4&5&6&7&8\\
\hline
\tt $\dim\Hc^{i}(\TIP')$&0&0&0&0&1&0&1&0&1\cr
\hline
\hline
\tt $\dim\Hc^{i}(\IP')$&0&1&1&0&1&0&1&&\cr
\hline
\tt $\dim\Hc^{i}(\IP'/\langle\theta\rangle)$&0&0&0&0&1&0&0&&\cr
\hline
\hline
\tt $\dim\Hr^{i}(\IP')$&1&0&1&0&1&1&0&&\cr
\hline
\tt $\dim\Hr^{i}(\IP'/\langle\theta\rangle)$&0&0&1&0&0&0&0&&\cr
\hline
\end{tabular}$$
Et on voit bien que la variété $(\IP_3(\CC)\mmoins2)/\langle\theta\rangle$ n'est pas \uacyclique.

\section{Rappels sur les nombres de Stirling}
\glossarytitle{Rappels sur les nombres de Stirling}
Pour\label{factorielles-Stirling} les notions de cette section de rappels, nous renvoyons à \cite{knuth} pour plus de détails.

\subsection{Factorielles croissantes et décroissantes}Pour\label{factorielles} $n\in\NN$, les \expression{$n$-ièmes factorielles croissante et décroissante}\glossary{$x\osp{n}{:=}x(x+1)\cdots(x+(n-1))$:$n$-ième factorielle croissante de $x$.}\glossary{$x\usp{n}{:=}x(x-1)\cdots(x-(n-1))$:$n$-ième factorielle décroissante de $x$.} (\footnote{\expression{rising and falling factorials} en anglais.}) d'un élément $x$ d'un anneau, notées respectivement $x\osp n$ et $x\usp n$, sont définies par
$$\def\et{\text{\quad et\quad}}
\begin{cases}
x\osp{0}:=1\,,\et
x\osp{n}:=x(x+1)(x+2)\cdots(x+(n-1))\,,\\
x\usp{0}:=1\,,\et
x\usp{n}:=x(x-1)(x-2)\cdots(x-(n-1))\,.
\end{cases}
$$
On a clairement 
$$\preskip0pt
x\usp{n+1}=x\usp{n}(x-n)\text{\quad et\quad}x\osp{n+1}=x\osp{n}(x+n)\,.$$

\subsection{Nombres de Stirling de première espèce}Dans\label{nombres-stirling} l'anneau de polynômes $\ZZ[X]$, le sous-module $\ZZ\sp n[X]$ des polynômes de degré majoré par $n$ admet les trois bases suivantes
$$
\B:=\set x\sp{0},x\sp{1},\ldots,x\sp{n}/\,,\quad
\Uline{0.5pt}{1pt}\B:=\set x\FD{0},x\FD{1},\ldots,x\usp{n}/\,,\quad
\let\sp\FA
\Oline{2pt}{0.5pt}\B:=\set x\sp{0},x\sp{1},\ldots,x\sp{n}/\,,
$$

Le développement des polynômes $X\usp n$ et $X\osp n$ en somme de monômes donne les coefficients des matrices de passage de $\B$ vers les deux autres bases. On note ces matrices respectivement par $(\ustirling ij)$ et $(\ostirling ij)$. On a donc
$$
X\usp i=\sumnl_{i\geq j\geq 0}\ustirling ij X^j
\text{\quad et\quad }
X\osp i=\sumnl_{i\geq j\geq 0}\ostirling ij X^j\,.$$
Les matrices $(\ustirling ij)$ et $(\ostirling ij)$ sont triangulaires inférieures avec des $1$ sur la diagonale. On voit clairement que $\ostirling ij\geq0$ et que $(-1)^{i-j}\ustirling ij\geq0$.

On étend la définition de $\ustirling ij$ et $\ostirling ij$ à tous les indices $i,j\in\NN$, par la valeur $0$ lorsque $j>i$. Ainsi, les  sommations ci-dessous peuvent être indexées tout simplement indexées par `${j\geq 0}$'.

\begin{defi}\label{def-stirling}Pour $i,j\in\NN$, les entiers naturels $\ostirling ij\in\NN$ sont  \expression{les nombres de Stirling de première espèce (non signés)}\glossary{$\ustirling ij,\ostirling ij$:nombres de Stirling de première espèce resp. signés et non signés}, et les entiers relatifs $\ustirling ij\in\ZZ$ sont \expression{les nombres de Stirling de première espèce (signés)}.
\end{defi}

\begin{lemm}\label{lemme-stirling}Les nombres de Stirling de première espèce vérifient les propriétés suivantes.
\begin{enumerate}
\item\leavevmode\label{lemme-stirling-a}Pour tous $i,j\in\NN$, on a $\ustirling ij=(-1)^{i-j}\,\ostirling ij$. En particulier, si $D$ désigne la matrice diagonale $\diag(1,-1,\ldots,(-1)^i)$, on a
$$(\ustirling ij)=D(\ostirling ij)D^{-1}$$
\item\leavevmode\label{lemme-stirling-b}Pour tout $i\geq0$, on a
$\ustirling ii=1$ et $\ostirling ii=1$.
\item\leavevmode\label{lemme-stirling-c}Pour tout $i\geq1$, on a $\ustirling0i=\ustirling i0=0$ et $\ostirling0i=\ostirling i0=0$.
\item\leavevmode\label{lemme-stirling-d}Pour tous $i,j\geq1$, on a 
$$\halfdisplayskips
\begin{cases}\noalign{\kern-4pt}
\ustirling{i}j=\ustirling {i-1}{j-1}-(i-1)\,\ustirling {i-1}j\,,\\\noalign{\kern1pt}
\ostirling{i}j=\ostirling {i-1}{j-1}+(i-1)\,\ostirling {i-1}j\,.\\\noalign{\kern-4pt}
\end{cases}$$
\item\leavevmode\label{lemme-stirling-e}Pour tout $i\geq1$, on a $\ustirling i1=(-1)^{i-1}\, (i-1)!$ et $\ostirling i1=(i-1)!$.
\end{enumerate}
\end{lemm}
\demo (\ref{lemme-stirling-a}) Évident puisque $X\osp i=(-1)^{i}(-X)\usp i$.
(\ref{lemme-stirling-b},\ref{lemme-stirling-c}) Pour tout $i>0$, les polynômes $X\usp i$ sont clairement de coefficient constant $0$ et de coefficient principal $1$. (\ref{lemme-stirling-d}) Résulte de ce que pour $i\geq1$, on a
$$\mathalign{
X\usp{i}&=&X\usp{i-1}(X-(i-1))=\Big(\sumnl_{j\geq0}\ustirling {i-1}j X^j\Big)(X-(i-1))\hfill\\\noalign{\kern2pt}&=&
\sumnl_{j\geq0}\ustirling {i-1}j X^{j+1}-
\sumnl_{j\geq0}\ustirling {i-1}j (i-1)X^{j}=\sum_{j\geq0}\ustirling ijX\sp{j}
\,.
\hfill}
$$
(\ref{lemme-stirling-e}) Pour $i>0$, on a 
$$\big(\ustirling i1=(X-1)\cdots(X-(i-1))\big)\rest_{X=0}=(-1)(-2)\cdots(-i+1)\,.$$
Les égalités concernant les coefficients $\ostirling ij$ résultent ensuite de (a).
\enddemo

\subsection{Nombres de Stirling de deuxième espèce}Les\label{nombres-Stirling} coefficients des matrices  
$(\uStirling ij):=(\ustirling ij)^{-1}$ et $(\oStirling  ij):=(\ostirling ij)^{-1}$ vérifient:
$$X^i=\sumnl _{j\geq0}\uStirling ijX\usp j=\sumnl _{j\geq0}\oStirling  ijX\osp j\,,\eqno\forall\ i\geq0. 
$$
Les matrices $(\uStirling ij)$ et $(\oStirling ij)$ sont triangulaires inférieures avec des $1$ sur la diagonale. On verra que $\uStirling ij\geq0$ et que $(-1)^{i-j}\oStirling ij\geq0$.

\begin{defi}\label{def-Stirling}Pour $i,j\in\NN$, les entiers naturels $\uStirling ij\in\NN$ sont  \expression{les nombres de Stirling de seconde espèce (non signés)}\glossary{$\oStirling ij,\uStirling ij$:nombres de Stirling de seconde espèce resp. signés et non signés}, et les entiers relatifs $\oStirling  ij\in\ZZ$ sont \expression{les nombres de Stirling de seconde espèce (signés)}.
\end{defi}

\begin{lemm}\label{lemme-Stirling}Nombres de Stirling de deuxième espèce vérifient les propriétés suivantes
\begin{enumerate}
\item\leavevmode\label{lemme-Stirling-a}Pour tous $i,j\in\NN$, on a $\oStirling  ij=(-1)^{i-j}\,\uStirling ij$. En particulier, si $D$ désigne la matrice diagonale $\diag(1,-1,\ldots,(-1)^i)$, on a
$$(\oStirling  ij)=D(\uStirling ij)D^{-1}$$
\item\leavevmode\label{lemme-Stirling-b}Pour tout $i\geq0$, on a
$\uStirling ii=1$ et $\oStirling  ii=1$.
\item\leavevmode\label{lemme-Stirling-c}Pour tout $i\geq1$, on a $\uStirling 0i=\uStirling i0=0$ et $\oStirling 0i=\oStirling  i0=0$.
\item\leavevmode\label{lemme-Stirling-d}Pour tous $i,j\geq1$, on a 
$$
\begin{cases}
\uStirling{i}{j}=\uStirling {i-1}{j-1}+j\,\uStirling {i-1}j\,,\\[3pt]
\oStirling {i}{j}=\oStirling  {i-1}{j-1}-j\,\oStirling  {i-1}j\,.
\end{cases}$$
\item\leavevmode\label{lemme-Stirling-e}Pour tous $i,j\in\NN$, on a $\uStirling ij\geq0$ et $(-1)^{i-j}\oStirling ij\geq0$.
\item\leavevmode\label{lemme-Stirling-f}Pour tout $i\geq0$, on a $\uStirling i1=1$ et $\oStirling  {i}1=(-1)^{i+1}$.
\end{enumerate}
\end{lemm}
\demo (\ref{lemme-Stirling-a}) Évident d'après \ref{lemme-stirling}-(\ref{lemme-stirling-a}).
(\ref{lemme-Stirling-b},\ref{lemme-Stirling-c}) Évidents. (\ref{lemme-Stirling-d}) Résulte de ce que pour $j\geq0$, on a $X\usp {j}X=X\usp{j+1}+jX\usp {j}$ et alors
$$\mathalign{
X\sp{i}&=&X\sp{i-1}X=\Big(\sum_{j\geq0}\uStirling {i-1}j X\usp j\Big)\,X\hfill\\\noalign{\kern2pt}
&=&
\sum_{j\geq0}\uStirling {i-1}j X\usp{j+1}+
\sum_{j\geq0}\uStirling {i-1}j jX\usp{j}=
\sum_{j\geq0}\uStirling {i}j X\usp{j}\,.
\hfill}
$$
(\ref{lemme-Stirling-e}) La positivité de $\uStirling ij$ résulte inductivement de (\ref{lemme-Stirling-d}) à partir de la positivité de $\uStirling 0j$ et $\uStirling 1j$, ce qui a été établi dans (\ref{lemme-Stirling-b}) et (\ref{lemme-Stirling-c}).
(\ref{lemme-Stirling-f}) La question (\ref{lemme-Stirling-b}) fixe le cas $i=1$, pour $i>1$, on a  $\uStirling i1=\uStirling {i-1}0+\uStirling{i-1}1=\uStirling{i-1}1$, d'après (\ref{lemme-Stirling-c}) et (\ref{lemme-Stirling-d}), et par induction, $\uStirling i1=\uStirling 11=1$.
\enddemo

\noendpoint\begin{remas}\label{rema-initiales}\mynobreak\begin{itemize}\nobreak\itemsep0pt\mou
\item {\sl Valeurs initiales des nombres de Stirling. }Les nombres de Stirling ont été indexés par les couples $(i,j)\in\NN^2$. Dans tous les cas, les coefficients de la colonne $(j=0)$ et la ligne $(i=0)$ sont nuls sauf pour $i=j=0$ où ils valent $1$. Tous les autres termes découlent de ces \expression{valeurs initiales} via les quatre règles de récurrence  (\ref{lemme-stirling-d}) des lemmes \ref{lemme-stirling} et \ref{lemme-Stirling}.

\item On remarquera aussi les égalités suivantes qui concernent les sous-matrices de nombres de Stirling d'indices non nuls.
$$\begin{casesalign}
 (\uStirling ij_{i,j\geq1})^{-1}&=&(\ustirling ij_{i,j\geq1})\\ 
 (\oStirling ij_{i,j\geq1})^{-1}&=&(\ostirling ij_{i,j\geq1})
\end{casesalign}$$
\end{itemize}\end{remas}

\subsection{Nombres de Stirling non signés et cardinaux}\label{Stirling}Pour tous $i, j\geq0\in\NN$, on définit:

\nobreak\hangindent1.9cm\hangafter1 
$\bullet$ \smash{$\Cycles ij:={}$}\glossary{$\sbCycles m\ell=\ostirling ij$:
cardinal de l'ensemble des permutations d'un ensemble à $i$ éléments qui sont produits d'exactement $j$ cycles}cardinal de l'ensemble des permutations d'un ensemble à $i$ éléments qui sont produits d'exactement $j$ cycles.

\hangindent1.9cm\hangafter1 
$\bullet$ \smash{$\Parties ij:={}$}\glossary{$\sbParties ij=\uStirling ij$:cardinal de l'ensemble $\Pgoth_j(E)$ de partitions d'un ensemble $E$ à $i$ éléments en $j$ parties non vides}cardinal de l'ensemble $\Pgoth_j(E)$ de partitions d'un ensemble $E$ à $i$ éléments en $j$ parties non vides.

On remarquera l'égalité \goodsmash{1}{1}{$\Parties 00=1$} qui dit qu'il y a une unique partition de l'ensemble vide en 0 parties non vides, et l'égalité \halfsmashtop{$\Cycles 00=1$} qui dit qu'il y a une unique permutation qui soit produit de 0 cycles.

\begin{prop}\label{Stirling-Cardinaux}Pour tous $i\geq j\geq0\in\NN$, on a 
\begin{enumerate}
\item \leavevmode\relax{$\Parties ij=\uStirling ij$}. Lorsque $i\geq1$, on a
$\halfsmash{\Parties ij={1\over j! }\displaystyle\sum_{k=0}^{j}(-1)^{j-k}\binome{j}{k}\,k^{i}\,.}$
\item $\Cycles ij=\ostirling ij=(-1)^{i-j}\ustirling ij$.
\end{enumerate}
\end{prop}
\demo On a bien \smashtop{$1=\Cycles 00=\Parties 00$ et 
$0=\Parties i0=\Cycles i0=\Parties 0j=\Cycles 0j$} pour tous $i,j\not=0$. Les familles des nombres en question ont donc bien les même valeurs initiales que les nombres de Stirling.

Pour $i\geq j\geq1$, on trie les partitions de $\Pgoth_{j}(\iii[1,i])$ en deux parties suivant qu'elles contiennent ou non le singleton $\set i/$. Les cardinaux de ces parties sont respectivement \smashbot{$\Parties{i-1}{j-1}$ et $j\Parties{i-1}{j}$}.
On a donc 
$\Parties ij=\Parties{i-1}{j-1}+j\Parties{i-1}{j}\,,$
ce qui correspond à la récurrence \ref{lemme-Stirling}-(\ref{lemme-Stirling-d}) pour les nombres $\uStirling ij$.

De même, en triant les permutations de $\iii[1,i]$, suivant que $\set i/$ est fixé ou non, on obtient deux classes de cardinaux \smashbot{$\Cycles {i-1}{j-1}$ et $(i-1)\Cycles {i-1}{j}$}, on a donc:
$\Cycles ij=\Cycles{i-1}{j-1}+(i-1)\Cycles{i-1}j\,,$
ce qui correspond à la récurrence \ref{lemme-stirling}-(\ref{lemme-stirling-d}) pour les nombres $\ostirling ij$.

Enfin, on rappelle que la formule dans (a) provient du dénombrement des surjections de $\iii[1,i]\onto\iii[1,j]$ (modulo les permutations de $\iii[1,j]$). Cet ensemble est le complémentaire $F_*$ dans l'ensemble $F$ de toutes les applications de $\iii[1,i]\to \iii[1,j]$ de l'ensemble des applications qui ne sont pas surjectives. Notons $F_t$ le sous-ensemble ses applications de $F$ qui n'atteignent pas la valeur $t\in\iii[1,j]$. Notons 
$F_{t_1,\ldots t_k}:=F_{t_1}\cap\cdots\cap F_{t_k}$.
On a alors 
$$\let\nolimits\relax
\mathalign{|F_*|
&=&\sumnl_{k=1}^{j}\ (-1)^{k-1}
\hskip-7mm
\sumnl_{1\leq t_1<\cdots <t_k\leq j}
\hskip-5mm
|F_{t_1,\ldots t_k}|
&=&\sumnl_{k=1}^{j}(-1)^{k-1}\binome jk (j-k)^{i}}\,,
$$
et la formule découle aussitôt.
\enddemo

\def\refname{Références bibliographiques}
\ifarxiv
\vskip0.5em
\hbox to\hsize{\hss \vrule height4pt depth-3.6pt width2cm \hss}
\else\newpage\fi

\ifarxiv
\vskip1em
\hbox to\hsize{\hss \vrule height4pt depth-3.6pt width2cm \hss}
\vskip-1ex\vskip0pt\else\newpage\fi

\printglossary

\end{document}